\documentclass[a4paper, 11pt]{article}
\usepackage{settings}

\pdfoutput=1

\begin{document}
\pagenumbering{roman}
\thispagestyle{empty}
\begin{center}
    \includegraphics[width=8cm,height=8cm,keepaspectratio]{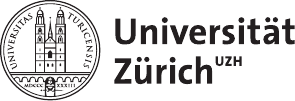}
    \label{fig:uzh}
    
   	\vspace{1cm}
	\Huge
	\textbf{Discretizations of Teichmüller Geodesic Flow and Enumeration of Pseudo-Anosov Diffeomorphisms}

    \vspace{1cm}
    
    \begin{center}
        \begin{tikzpicture}[scale = 1.3]
            \draw[smooth] (0,1) to[out=30,in=150] (2,1) to[out=-30,in=210] (3,1) to[out=30,in=150] (5,1) to[out=-30,in=30] (5,-1) to[out=210,in=-30] (3,-1) to[out=150,in=30] (2,-1) to[out=210,in=-30] (0,-1) to[out=150,in=-150] (0,1);
            \draw[smooth] (0.4,0.1) .. controls (0.8,-0.25) and (1.2,-0.25) .. (1.6,0.1);
            \draw[smooth] (0.5,0) .. controls (0.8,0.2) and (1.2,0.2) .. (1.5,0);
            \draw[smooth] (3.4,0.1) .. controls (3.8,-0.25) and (4.2,-0.25) .. (4.6,0.1);
            \draw[smooth] (3.5,0) .. controls (3.8,0.2) and (4.2,0.2) .. (4.5,0);
            \fill[white] (0,-3) circle (0.1pt);

            \draw[orange, thick] (2.5,0) -- (2.9,0);
            \draw[orange, thick, dash pattern=on 1.2pt off 1.6pt] (2.9,0) -- (3.4,0);
            \draw[teal, thick] (2.5,0) -- (2.1, 0);
            \draw[teal, thick, dash pattern=on 1.2pt off 1.6pt] (2.1, 0) -- (1.6,0);
            
            \draw[noamblue, thick] (2.1, -0.4) -- (2.5,0);
            \draw[noamblue, thick, dash pattern = on 1.2pt off 1.6pt] (2.1, -0.4) -- (1.7,-0.8);
            \draw[purple, thick] (2.5,0) -- (2.9, 0.4);
            \draw[purple, thick, dash pattern = on 1.2pt off 1.6pt] (2.9,0.4) -- (3.3, 0.8);

            \draw[violet, thick] (2.1, 0.4) -- (2.5,0) ;
            \draw[violet, thick, dash pattern = on 1.2pt off 1.6pt] (2.1, 0.4) -- (1.7,0.8);
            \draw[gray, thick, dash pattern = on 1.2pt off 1.6pt] (2.9, -0.4) -- (3.3, -0.8);
            \draw[gray, thick] (2.5,0) -- (2.9, -0.4);
            \fill[darkgray] (2.5,0) circle (2pt);
        \end{tikzpicture}
    \end{center}
	\vspace{-1cm}
	\huge
	Master Thesis in Mathematics
	
	\vspace{1cm}	
	\Large
	\textbf{Noam Mordehai Isaac Szyfer}
	
	\vspace{2cm}
	\normalsize	
	Supervised by \\
	\Large
	Prof. Dr. Corinna Ulcigrai\\
	Dr. Daniele Galli\\
	
	\vspace{2cm}
	\large
	\today
\end{center}
\vspace{\stretch{3}}
\pagebreak

\thispagestyle{empty}
\

\pagebreak

\pagenumbering{roman}
\setcounter{page}{3}
\thispagestyle{plain}
\begin{abstract}
    In this thesis, we study the Teichmüller geodesic flow on the space of translation surfaces by introducing two related discrete-time dynamical systems. First, we discuss the Rauzy--Veech induction, highlighting its connections to interval exchange transformations and continued fraction expansions. While effective for addressing ergodic properties, this method faces challenges in counting closed orbits. Second, we introduce diagonal changes, a discretization better suited for counting and enumeration problems, initially applied to hyperelliptic components—subsets of translation surfaces with additional symmetries. Understanding closed orbits is significant due to the one-to-one correspondence with pseudo-Anosov mapping classes. After detailing this connection, we demonstrate how diagonal changes can produce a complete list of pseudo-Anosov mapping classes ordered by dilatation, and extend the algorithm to general components of strata.
\end{abstract}

\pagebreak
\thispagestyle{empty}
\
\pagebreak

\thispagestyle{plain}
\epigraph{The product of mathematics is clarity and understanding. Not theorems, by themselves. The world does not suffer from an oversupply of clarity and understanding. Mathematics only exists in a living community of mathematicians that spreads understanding and breaths life into ideas both old and new. The real satisfaction from mathematics is in learning from others and sharing with others. All of us have clear understanding of a few things and murky concepts of many more. There is no way to run out of ideas in need of clarification.}{---\textup{William P. Thurston, 2010}}

\pagebreak

\thispagestyle{empty}
\
\pagebreak

\thispagestyle{plain}
\begin{center}
   \textbf{Acknowledgements}
\end{center}
This thesis would not have been possible without the invaluable support I had the privilege to receive from so many people. 

First and foremost, I want to express my deep gratitude to Prof.\ Dr.\ Corinna Ulcigrai, not only for agreeing to supervise my thesis and suggesting this interesting topic, but also for being an incredible lecturer reigniting my enthusiasm for Mathematics when my excitement was fading. I also want to extend my sincere thanks to my co-advisor Dr.\ Daniele Galli both for listening and explaining, patiently pointing me in the right direction when I myself could not see things clear enough.

I am also extremely grateful to the many people sharing their knowledge and passion with me during my journey through Mathematics. Severin, Andres and Nikolai, I will fondly remember your classes forever, thank you for going above and beyond in trying to make the difficult things seem a bit less difficult. Special thanks to Niklas for teaching me just enough Algebra to pass my exam and for being a great friend, and thanks also to Frank, Przemek and especially to Yuriy for explaining the intricacies of dynamical systems to me. I also want to thank my peers, Bettina, Sebastian, Tabea, Flo, Franziska, Adina, Silvan, Ephraim, Alain, Rathes and all the others for sharing both the good times and the struggles. I thank Zouhair for sharing with me his beautiful perspective of Mathematics as a whole.  

I would have never embarked on this journey without Maxime and Vali paving the way, convincing me that returning to University was not only possible but actually a great idea. Many thanks also to Niculin, Mischa and Maxime for reverting me into a musician at least once a year and hence granting me a true break from Mathematics.

The journey towards my degree in Mathematics would not have been possible without the support of my family. I am very grateful to my parents for giving me the opportunity to always pursue my interests. Lastly, words cannot express my gratitude to Simona, whose love and support were the things keeping me sane during the ride.  

\pagebreak

\thispagestyle{empty}
\
\pagebreak

\thispagestyle{plain}
\newgeometry{left=2.5cm,right=2.5cm, top=2.3cm, bottom=3cm}

{\footnotesize\tableofcontents}

\restoregeometry
\pagebreak

\thispagestyle{empty}
\
\pagebreak

\pagenumbering{arabic}

\section{Introduction}\label{sec:introduction}

\thispagestyle{plain}

The mathematical area of Dynamical Systems can be broadly described as the \emph{study of things that evolve in time}. This very general description reveals one particular appeal of this theory: Almost any system one wishes to study can be seen as something that changes over time. If one is interested in the applied side of things, then examples of such systems could be a vessel containing a gas in Physics or Chemistry, a population of a species in Biology or the weather in Meteorology. But also in pure Mathematics there is an abundance of objects that can be seen as evolving in time, for example the iteration of a function or the evolution in time of a solution of some systems of equations. 

It is therefore unsurprising that this subject intersects many other mathematical fields. Dynamical approaches have proven useful in Mathematical Physics, Number Theory and Probability Theory, to name a few. In this thesis, we will examine Dynamical Systems in the context of Geometry. More precisely, we are interested in the study of a special class of surfaces, that is, two-dimensional real Riemannian manifolds. 

The thesis is divided into six parts. In the \emph{introduction}, we give some motivation as to why we choose to study a particular class of surfaces, namely \emph{translation surfaces}, which we will then properly define. In the second part, we will introduce a general principle in the study of translation surfaces, which has been applied with great success to several questions regarding the dynamics on such surfaces. Instead of studying a single, fixed surface, we will look at this surface as a point in a \emph{family of surfaces} sharing the same topology, which includes for example the genus. The dynamics of \enquote{walking on the surface in a straight line} are translated into deforming this surface in a particular way. The specific type of deformation we will be most interested in is known as the \emph{Teichmüller geodesic flow}, which will lie at the center of the remaining four parts. These first two parts of the thesis are inspired by two wonderfully written survey papers by Anton Zorich \cite{zorich2006flat} and Daniel Massart \cite{massart2022short} as well as a lecture given by Corinna Ulcigrai at the University of Zurich in 2022.

In the third part, we will connect translation surfaces to a special kind of maps from an interval to itself, called \emph{interval exchange maps}. We will work out explicitly a connection to continued fraction expansions, which will allow us to see the dynamics on surfaces of higher genus as a generalization of the Euclidean algorithm, and we will see that this algorithm is closely tied to the Teichmüller geodesic flow we have just mentioned. Explicitly, it can be seen as a \emph{discretization} of the flow, which is classically known as \emph{Rauzy--Veech induction}. This discretization is a powerful tool, used very successfully in answering various questions about the dynamical properties of flows on surfaces. However, it seems impossible to use Rauzy--Veech induction in order to find solutions to so-called \emph{counting problems}, where one is interested in the number of closed geodesics up to a certain length. For this reason, an alternative discretization was developed, which is based on \emph{diagonal changes}. We will introduce a special case of this algorithm, restricted to \emph{hyperelliptic} surfaces, in detail in the fourth part, after which in the fifth part of the thesis we will explain how exactly diagonal changes can be used in counting problems. Lastly, in the sixth part, we present a generalization of diagonal changes going beyond hyperelliptic surfaces, which can be seen as a translation and refinement of the work of Sebastien Ferenczi in \cite{ferenczi2015diagonal} to our setting. 

\subsection{Motivation}
We want to study surfaces, or more precisely, 2-dimensional compact orientable real manifolds endowed with a Riemannian metric. However, studying the dynamics related to this very large class of objects turns out to be unfeasible in this generality, so we are compelled to restrict ourselves to a smaller class of objects. One possibility is to only consider \emph{flat} surfaces, in the sense that the Gaussian curvature is zero everywhere. But this may be too restrictive. As will follow from the famed Gauss--Bonnet theorem ({Theorem \ref{thm:gauss_bonnet}}), considering only manifolds without boundary there is only one topological surface that allows the Gaussian curvature to vanish throughout, the torus. 

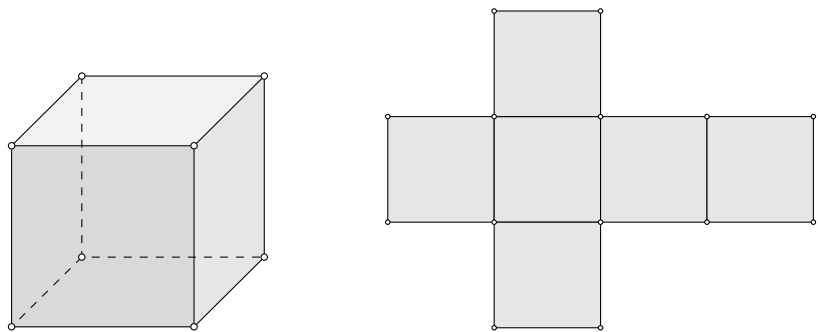
\begin{figure}[ht]
    \centering
    \begin{tikzpicture}[scale=1.2,transform shape]
    \coordinate (A) at (0,0,0);
    \coordinate (B) at (2,0,0);
    \coordinate (C) at (2,2,0);
    \coordinate (D) at (0,2,0);
    \coordinate (E) at (0,0,2);
    \coordinate (F) at (2,0,2);
    \coordinate (G) at (2,2,2);
    \coordinate (H) at (0,2,2);

    \filldraw[gray!30] (E) -- (F) -- (G) -- (H) -- cycle;
    \filldraw[gray!10] (G) -- (H) -- (D) -- (C) -- cycle;
    \filldraw[gray!20] (B) -- (F) -- (G) -- (C) -- cycle;
    
    \draw (B) -- (C) -- (D);
    \draw (E) -- (F) -- (G) -- (H) -- cycle;
    \draw[dashed] (A) -- (E);
    \draw (B) -- (F);
    \draw (C) -- (G);
    \draw (D) -- (H);
    \draw[dashed] (D) -- (A) -- (B);

    \foreach \vertex in {A,B,C,D,E,F,G,H}
        \filldraw[white, draw=black] (\vertex) circle (1pt);
\end{tikzpicture}\qquad\qquad 
\begin{tikzpicture}[scale=.7,transform shape]
    \draw[fill = gray!20] (0,0) rectangle (2,2);
    \draw[fill = gray!20] (2,0) rectangle (4,2);
    \draw[fill = gray!20] (4,0) rectangle (6,2);
    \draw[fill = gray!20] (6,0) rectangle (8,2);
    \draw[fill = gray!20] (2,2) rectangle (4,4);
    \draw[fill = gray!20] (2,-2) rectangle (4,0);
    
    \filldraw[white, draw=black] (0,0) circle (1.2pt);
    \filldraw[white, draw=black] (2,2) circle (1.2pt);
    \filldraw[white, draw=black] (2,0) circle (1.2pt);
    \filldraw[white, draw=black] (4,2) circle (1.2pt);
    \filldraw[white, draw=black] (4,0) circle (1.2pt);
    \filldraw[white, draw=black] (6,2) circle (1.2pt);
    \filldraw[white, draw=black] (6,0) circle (1.2pt);
    \filldraw[white, draw=black] (8,2) circle (1.2pt);
    \filldraw[white, draw=black] (4,4) circle (1.2pt);
    \filldraw[white, draw=black] (2,-2) circle (1.2pt);
    \filldraw[white, draw=black] (2,4) circle (1.2pt);
    \filldraw[white, draw=black] (4,-2) circle (1.2pt);
    \filldraw[white, draw=black] (0,2) circle (1.2pt);
    \filldraw[white, draw=black] (8,0) circle (1.2pt);
\end{tikzpicture}
    \caption{The cube viewed as a flat surface with eight singularities.}
    \label{fig:cube}
\end{figure}

We can remedy this situation by allowing our surfaces to have a certain number of special points, or \emph{singularities}. One way to think of these singular points is to imagine the surface as being made of a malleable material such as aluminum foil, so that we can flatten large parts of the surface and push all the curvature to smaller and smaller regions, so that finally all the curvature is concentrated in a single point, the singularity. A simple example of this is given by a cube as in Figure \ref{fig:cube}, which has 8 singular points which are exactly the vertices of the cube. The cube can be realized by identifying the correct edges in the collection of squares on the right side of Figure \ref{fig:cube}.

Note that the edges might appear singular as well, but this is just a consequence of embedding the cube into the Euclidean space $\R^3$. Indeed, taking a small neighborhood of an interior point of an edge of the cube, we can remove the crease given by the edge very easily just by flattening out the neighborhood. This is not the case for the vertices. Taking any neighborhood of a vertex, we can see on the right of Figure \ref{fig:cube} that such a neighborhood is isometric to a neighborhood of a vertex of a \emph{cone}. Indeed, going around the vertex, one does not travel an angle of $2\pi$ but only $\frac{3\pi}{2}$, before returning to the start point. For this reason, vertices of these types are known as \emph{conical} singularities. 

Considering a manifold such as the cube endowed with a Riemannian metric, natural objects to study are \emph{geodesics}, paths that can be seen a generalization of a straight line in Euclidean space. And since the metric we consider is flat, locally geodesics on the cube are isometric to straight lines. Typical questions about geodesics concern for instance their long-term behavior. If we imagine living on the surface, and we start walking straight ahead, so along a geodesic, will we ever return to our starting point? If not, will we eventually reach every other part of the surface? Will we spend the same time in every part of the surface? In more technical terms, questions about geodesics that do not return are about the \emph{ergodic} properties of the geodesic flow. Other questions revolve about cases where we \emph{do} come back to the starting point. If this is the case, we speak of a \emph{closed} geodesic. Can we find a closed geodesic on any surface which does not go through a singularity? If yes, how many of them (up to a certain length $L$) can we find?

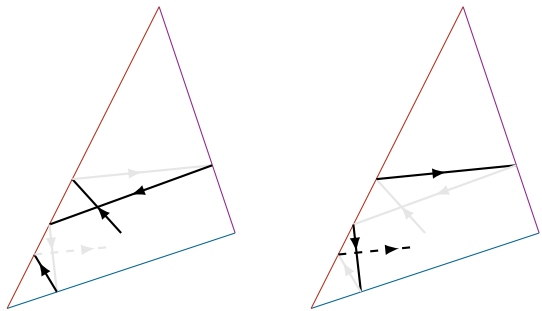
\begin{figure}[ht]
    \centering
    \begin{tikzpicture}[scale = 1]
  \coordinate (A) at (0,0);
  \coordinate (B) at (3,1);
  \coordinate (C) at (2,4);

  \draw[name path = AB, color = MidnightBlue] (A) -- (B);
  \draw[name path = BC,color = Plum] (B) -- (C);
  \draw[name path = CA,color = BrickRed] (C) -- (A);

  \coordinate (D) at (4,0);
  \coordinate (E) at (7,1);
  \coordinate (F) at (6,4);

  \draw[color = MidnightBlue] (D) -- (E);
  \draw[color = Plum] (E) -- (F);
  \draw[color = BrickRed] (F) -- (D);

    \coordinate (start) at (1.5, 1); 
    \coordinate (end) at ($(start) + (-1.8, 2)$); 

    \path[name path=trajectory] (start) -- (end);
    \path[name intersections={of=trajectory and CA, by=reflection}];
    \draw[thick, black, postaction={
    decorate, decoration={ markings, mark=at position 0.5 with{\arrow{latex}}}}] (start) -- (reflection);
    \draw[thick, gray!20,postaction={
    decorate, decoration={ markings, mark=at position 0.5 with{\arrow{latex}}}}] ($(start) + (4,0)$) -- ($(reflection) + (4,0)$);

    \coordinate (start) at (reflection);
    \coordinate (dir) at ($(start) + (2,0.2)$);
    \path[name path=trajectory] (reflection) -- (dir);
    \path[name intersections={of=trajectory and BC, by=reflection}];
    \draw[thick, gray!20, postaction={
    decorate, decoration={ markings, mark=at position 0.5 with{\arrow{latex}}}}] (start) -- (reflection);
    \draw[thick, black, postaction={
    decorate, decoration={ markings, mark=at position 0.5 with{\arrow{latex}}}}] ($(start) + (4,0)$) -- ($(reflection) + (4,0)$);

    \coordinate (start) at (reflection);
    \coordinate (dir) at ($(start) + (-3,-1.1)$);
    \path[name path=trajectory] (reflection) -- (dir);
    \path[name intersections={of=trajectory and CA, by=reflection}];
    \draw[thick, black, postaction={
    decorate, decoration={ markings, mark=at position 0.5 with{\arrow{latex}}}}] (start) -- (reflection);
    \draw[thick, gray!20, postaction={
    decorate, decoration={ markings, mark=at position 0.5 with{\arrow{latex}}}}] ($(start) + (4,0)$) -- ($(reflection) + (4,0)$);
    
     \coordinate (start) at (reflection);
    \coordinate (dir) at ($(start) + (0.1,-0.9)$);
    \path[name path=trajectory] (reflection) -- (dir);
    \path[name intersections={of=trajectory and AB, by=reflection}];
    \draw[thick, gray!20, postaction={
    decorate, decoration={ markings, mark=at position 0.4 with{\arrow{latex}}}}] (start) -- (reflection);
    \draw[thick, black, postaction={
    decorate, decoration={ markings, mark=at position 0.4 with{\arrow{latex}}}}] ($(start) + (4,0)$) -- ($(reflection) + (4,0)$);   

     \coordinate (start) at (reflection);
    \coordinate (dir) at ($(start) + (-0.3,0.5)$);
    \path[name path=trajectory] (reflection) -- (dir);
    \path[name intersections={of=trajectory and CA, by=reflection}];
    \draw[thick, black, postaction={
    decorate, decoration={ markings, mark=at position 0.8 with{\arrow{latex}}}}] (start) -- (reflection);
    \draw[thick, gray!20, postaction={
    decorate, decoration={ markings, mark=at position 0.8 with{\arrow{latex}}}}] ($(start) + (4,0)$) -- ($(reflection) + (4,0)$);  

     \coordinate (start) at (reflection);
    \coordinate (dir) at ($(start) + (1,0.1)$);
    \draw[thick, gray!20, dashed, postaction={
    decorate, decoration={ markings, mark=at position 0.8 with{\arrow{latex}}}}] (start) -- (dir);
    \draw[thick, black, dashed, postaction={
    decorate, decoration={ markings, mark=at position 0.8 with{\arrow{latex}}}}] ($(start) + (4,0)$) -- ($(dir) + (4,0)$);  

\end{tikzpicture}
    \caption{Correspondence between a billiard trajectory in a triangle and a linear flow on a sphere via the Fox--Kershner construction, where edges of the same color are identified. }
    \label{fig:fox_kershner}
\end{figure}

One is easily convinced that answers to these questions for general Riemannian surfaces are difficult to obtain, thus it should not be surprising that not much is known. Perhaps more surprisingly, the case for \emph{flat} surfaces is exactly the same. To give an example, problems of the kind introduced above remain open even in very simple cases, considering a sphere with only three singularities. To see why this is the case, we refer to a seemingly different, but actually equivalent problem, namely billiards in a triangle. A well-known open problem is the following.
\begin{problem}
    Does every obtuse triangle $T$ admit a closed billiard trajectory?
\end{problem}
If the triangle $T$ is acute or has a right angle, the answer is affirmative. Also, if all angles are rational multiples of $\pi$ then we are always able to find a closed billiard trajectory. For obtuse triangles without rational angles, the problem is only partially solved. In \cite{schwartz2009obtuse}, Schwartz shows that if the obtuse angle is of at most 112.5 degrees, then again the answer is yes. Unfortunately, the approach does not work for wider angles, so that the general problem remains open. To translate this billiard problem into a statement on closed geodesics on a flat surface, we may employ the so-called Fox--Kershner construction introduced in \cite{fox1936concerning}. The idea is to consider two copies of the triangle, but only a single billiard trajectory that jumps between the two copies at every reflection. Identifying the boundaries of the triangles, one obtains a path on the resulting topological sphere, which has conical singularities exactly at the vertices of the triangles. This is illustrated in Figure \ref{fig:fox_kershner}.

Another way of thinking about conical singularities is in terms of \emph{parallel transport}. Parallel transport along a null-homotopic curve transports a vector tangent to the surface back to itself. If the curve is not homotopic to a trivial one, i.e., if it goes around a singularity, this must no longer be the case: the returning vector has turned by some angle. In fact, parallel transport around a conical singularity makes the vector turn by exactly the cone angle. If this can happen, i.e., a tangent vector does not necessarily return to itself after performing any parallel transform, we say that the surface has \emph{nontrivial linear holonomy}. In this situation, a geodesic is forced to come back to itself and intersect itself again and again in different directions. Such a thing cannot happen on a flat surface of genus 1, i.e., a torus. The reason is, that any parallel transport brings any tangent vector back to itself, meaning that the torus has \emph{trivial linear holonomy}. It turns out that this additional constraint, \emph{trivial linear holonomy}, makes answering the questions above tractable. So the objects we will consider can be seen as \emph{very flat} surfaces, or \emph{translation surfaces} as we will define them below, which are closed orientable surfaces endowed with a flat metric having a finite number of conical singularities and having trivial linear holonomy. Note that triviality of linear holonomy implies in particular that all cone angles at conical singularities must be integer multiples of $2\pi$.

\subsection{Definitions of Translation Surfaces}

Let us now define translation surfaces, which will be our main objects of study. There are several equivalent way of defining translation surfaces, among them are definitions using the language of Euclidean Geometry, Riemannian Geometry and Complex Analysis. The perspective coming from Euclidean Geometry may be seen as the most natural approach, as it needs very little technical prerequisites. For this reason, we will use this as our main Definition. Nonetheless, both other points of view will be useful throughout this thesis, so we will present these alternative definitions as well, and we will motivate why the definitions are equivalent. The correspondences between the different viewpoint are worked out in great detail, in \cite{athreyamasur2023translationsurfaces} and also in \cite{massart2022short}. 

\begin{definition}[Translation surface]\label{def:translation_surface}
    Let $\mathcal{P} = \{P_1, \ldots, P_n\}$ be a nonempty collection of disjoint polygons $P_i \subseteq \R^2$ (or $P_i \subseteq \C$), such that the collection $\mathcal{E}$ of the edges of the polygon can be grouped into pairs $(e_1, e_2), e_i \in \mathcal{E}$ so that each edge $e$ of $P_j \in \mathcal{P}$ belongs to exactly one pair, and the two sides $(e_1, e_2)$ in each pair are parallel and isometric. We orient each edge $e \in \mathcal{E}$ such that the Euclidean translation taking $e_1$ to $e_2$ preserves the orientation. We now identify the two edges of each pairing, requiring that the polygon of one of the oriented sides lies to the left of this side, while the polygon of the other oriented side lies to the right of this polygon.

    Formally, we now define a \emph{translation surface} $X$ to be
    \begin{equation*}
        X = \mathcal{P}/_\sim,
    \end{equation*}
    where we write $\sim$ for the equivalence relation between edges described above. 
\end{definition}
\begin{remark}
    Even though the collection of polygons $\mathcal{P}$ will not be connected in general, exactly if it consists of more than one polygon, due to the identifications of the sides the resulting quotient $\mathcal{P}/_\sim$ may still be connected. We will always assume that a translation surface is connected, for the simple reason that if it is not, we may study the individual components instead. Let us also mention that the collection $\mathcal{P}$ in Definition \ref{def:translation_surface} is assumed to be finite. Without this assumption, one enters the area of \emph{infinite} translation surfaces, where many tools used here become unavailable. 
\end{remark}

\begin{figure}[ht]
    \centering
    \begin{tikzpicture}[scale = 1]
        \draw[thick, blue] (0,0) -- (3,0) node[midway, below, black] {$b$};
        \draw[thick] (0,0) -- (0,3) node[midway, left] {$a$};
        \draw[thick] (3,0) -- (3,3) node[midway, right] {$a$};
        \draw[thick] (0,3) -- (3,3) node[midway, above] {$b$};
    \end{tikzpicture}\qquad 
\begin{tikzpicture}
  \node (octagon) [regular polygon, regular polygon sides=8, minimum size=4cm, draw, thick] at (0,0) {};
  
  \node[above] at (octagon.side 1) {$a$};
  \node[below] at (octagon.side 5) {$a$};
  \node[yshift = 5pt, xshift = -5pt] at (octagon.side 2) {$b$};
  \node[yshift = -5pt, xshift = 5pt] at (octagon.side 6) {$b$};
  \node[left] at (octagon.side 3) {$c$};
  \node[right] at (octagon.side 7) {$c$};
  \node[yshift = -5pt, xshift = -5pt] at (octagon.side 4) {$d$};
  \node[yshift = 5pt, xshift = 5pt] at (octagon.side 8) {$d$};
\end{tikzpicture}\qquad
\begin{tikzpicture}[scale = 1]
  \node (pentagon1) [regular polygon, regular polygon sides=5, minimum size=3cm, draw, thick] at (0,0) {};
  
  \node (pentagon2) [regular polygon, regular polygon sides=5, minimum size=3cm, draw, thick, rotate=36] at (2.4cm,-0.78cm) {};

    \node[yshift = 5pt, xshift = -3pt] at (pentagon1.side 1) {$a$};
    \node[yshift = 0pt, xshift = -5pt] at (pentagon1.side 2) {$b$};
    \node[below] at (pentagon1.side 3) {$c$};
    \node[yshift = 3pt, xshift = -5pt] at (pentagon1.side 4) {$d$};
    \node[yshift = 5pt, xshift = 3pt] at (pentagon1.side 5) {$e$};

    \node[yshift = -3pt, xshift = 5pt] at (pentagon2.side 1) {$d$};
    \node[yshift = -5pt, xshift = -5pt] at (pentagon2.side 2) {$e$};
    \node[yshift = -5pt, xshift = 3pt] at (pentagon2.side 3) {$a$};
    \node[yshift = 3pt, xshift = 5pt] at (pentagon2.side 4) {$b$};
    \node[above] at (pentagon2.side 5) {$c$};
\end{tikzpicture}
    \caption{Three examples of translation surfaces. The edges with the same label are identified by Euclidean translations. }
    \label{fig:translation_surfaces}
\end{figure}
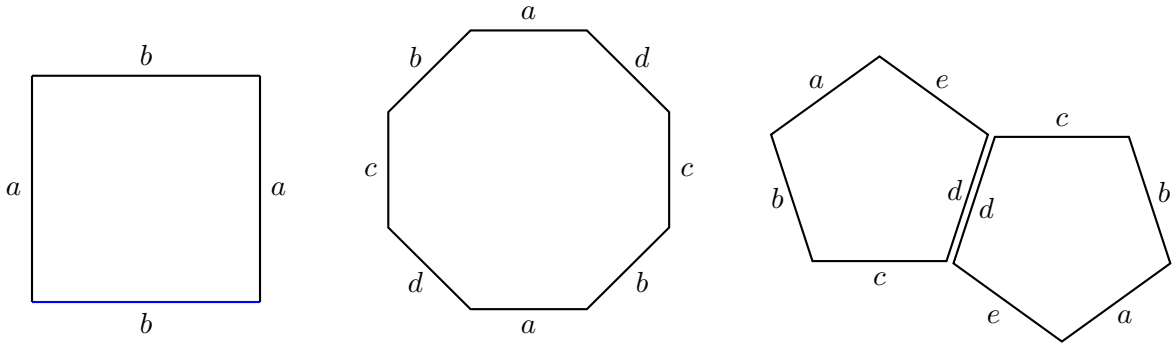

Three different translation surfaces arising from this polygon construction can be seen in Figure \ref{fig:translation_surfaces}. Topologically, the left example corresponds to a surface of genus $\mathbf{g} = 1$, while the middle and right examples correspond to surfaces of genus $\mathbf{g} = 2$.

Let us stress that we distinguish between a translation surface $X = \mathcal{P}/_\sim$ and the translation surface $X' = \mathcal{P}'/_\sim$, where $\mathcal{P}' = \{\e^{\ii\theta}P_1, \ldots \e^{\ii \theta}P_n\}$ is the collection of polygons, where each individual polygon has been rotated by the angle $\theta$. Put differently, each translation surface as constructed above comes with a choice of vertical direction (which is induced by the vertical direction in $\R^2$), and a rotation of (the polygons inducing) the surface yields a \emph{distinct} surface.

\begin{definition}[Singularity set]
    The finite set $\Sigma \subseteq X$ consisting of the images of the vertices of polygons in $\mathcal{P}$ under the identification map is called the \emph{singularity set} of the translation surface $X$.
\end{definition}

For example, one checks that all examples in Figure \ref{fig:translation_surfaces} have only one singularity, i.e., $|\Sigma| = 1$ in all three cases. We also want to mention that allowing singular points is absolutely necessary in order to obtain translation surfaces of genus $\mathbf{g} \geq 2$. To see why this is the case, we must recall the celebrated theorem of Gauss--Bonnet.

\begin{theorem}[Gauss--Bonnet]\label{thm:gauss_bonnet}
    Suppose $X$ is a 2-dimensional Riemannian manifold without boundary. Let $K$ be the Gaussian curvature of $X$. Then,
    \begin{equation*}
        \int_X K \, \dd \operatorname{Area} = 2\pi \chi(X),
    \end{equation*}
    where $\chi(X)$ denotes the Euler characteristic of $X$. If $X$ is orientable, then Euler's formula implies that $\chi(X) = 2 - 2\mathbf{g}$, where $\mathbf{g}$ denotes the genus of $X$. 
\end{theorem}
We see that if the Gaussian curvature vanishes \emph{everywhere}, this forces the genus of $X$ to be 1. In Proposition \ref{prop:genus_formula} below, we will see that we can interpret the singularities $\Sigma$ as Dirac-masses of curvature, i.e., all negative curvature is contained in the finite set $\Sigma$ in a way such that taking the integral with respect to $\dd \operatorname{Area}$ gives back exactly the finite value such that we recover the Gauss--Bonnet theorem.

Since translations are Euclidean isometries, the quotient $X$ inherits a local Euclidean metric from the polygons away from the singularity set. Let us make this more precise in showing that the quotient space is a complex 1-dimensional manifold, which under the usual identification $\R^2 \cong \C$ is equivalent to showing that it is a 2-dimensional real manifold. 

\begin{proposition}[Translation surface is a manifold, \cite{massart2022short}]\label{prop:translation_surface_is_manifold}
    Let $X$ be a translation surface, then it is a 1-dimensional complex manifold, i.e., every point $x \in X$ has a neighborhood $U_x$ in $X$ and a chart $\phi_x\colon U_x \to \C$, such that for any $x,y$ in $X$, if $U_x \cap U_y \neq \emptyset$, then $\phi_y\circ\phi_x^{-1}$ is a homeomorphism between open subsets of $\C$.
\end{proposition}
\begin{proof}
    For points $x \in X$ that are the images of interior points $\Bar{x} \in \operatorname{int}(P_i)$, we can take a small enough neighborhood $\Bar{U}_x$ of $\Bar{x}$ such that $\Bar{U}_x \subseteq \operatorname{int}(P_i)$ and take the identity map as the chart $\phi_x$. 
    
    If $x$ is the image of an interior point of an edge $e \in \mathcal{E}$, then $x$ has exactly two preimages in $P_i$ and $P_j$, where possibly, $i = j$, and $x$ has a neighborhood in $X$ whose preimage in $P_i$ and $P_j$ is given by two half-disks of equal radius and parallel diameters. We can define a chart $\phi_x$ by declaring it to be the identity on one of the half-disks and an appropriate translation on the other half-disk.   

    Lastly, if $x \in \Sigma$ is the image of a vertex, then all preimages of $X$ are vertices of polygons $P_{i_1}, \ldots, P_{i_k} \in \mathcal{P}$. We can see that $x$ has a neighborhood $U_x$ in $X$ whose pre-image in $\mathcal{P}$ is a union of circular sectors with common radius $r$ and straight boundaries $r_{i_j}$ and $s_{i_j}$, so that $s_{i_j}$ is parallel to $r_{i_{j+1}}$ for $j \in [k-1]$ and $s_{i_k}$ is parallel to $r_{i_1}$. Let us denote by $\theta_i$ the angle of sector $i$, i.e., the sector between $r_{i_j}$ and $s_{i_j}$. Since $s_{i_k}$ is parallel to $r_{i_1}$ and identifications are done by translations, the angles $\theta_i$ must sum to an integer multiple of $2\pi$, say $2\xi\pi$ for $\xi \in \N_{\geq 1}$. Let us write
    \begin{equation*}
        \Theta_i = \sum_{j = 1}^{i-1} \theta_j
    \end{equation*}
    for the partial sums of the angles. This allows us to define a chart which \enquote{compresses} each sector by scaling the angles such that the total angle in the image is again exactly $2\pi$, i.e., we obtain a Euclidean neighborhood also for singular points. Explicitly, the chart is defined on the $j$\textsuperscript{th} sector by
    \begin{equation*}
        p_{i_j} + \rho \e^{\ii(\theta + \Theta_j)} \mapsto \rho \e^{\frac{\ii(\theta + \Theta_j)}{\xi}} \quad 0\leq \rho \leq r, 0 \leq \theta \leq \theta_j,
    \end{equation*}
    where $p_{i_j}$ denotes the vertex in $\mathcal{P}$. Note that the image of the chart is given by a ball of radius $r$ around the origin in $\C$. This shows that a translation surface has the structure of a \emph{topological manifold}. By examining the transition maps, we can see that we may endow a translation surface with more structure.
    
    Adapting the convention from \cite{massart2022short}, we call a point $x \in X$ which is the image of an interior point of some $P_i$ of type I, a point $x \in X$ which is the image of two interior points of edges of $P_i, P_j$ of type II, and lastly we say singular points are of type III. 

    If $x$ and $y$ are both of type I and $U_x \cap U_y$ is nonempty, then the transition map $\phi_y \circ \phi_x^{-1}$ is simply the identity. If $x$ and $y$ are either of type I or type II and $U_x \cap U_y$ is nonempty, then the transition map $\phi_y \circ \phi_x^{-1}$ is a translation, i.e., of the form $z \mapsto z + c$. If $x$ and $y$ are both of type III, we may simply choose the neighborhoods $U_x$ and $U_y$ in a way so that $U_x \cap U_y = \emptyset$. Finally, if $x$ is of type III and $y$ is of type $I$ or $II$, then the transition map $\phi_y\circ\phi_x^{-1}$ can be seen as a branch of the $\xi$\textsuperscript{th} root. 
\end{proof}

Note that the transition maps being given by translations away from the singularity set $\Sigma$ justifies the use of the term \emph{translation} surface. Moreover, since the $\xi$\textsuperscript{th} root is holomorphic when restricted appropriately to a branch, translation surfaces come equipped with a \emph{complex} structure as well. 

In fact, the property we have proven in Proposition \ref{prop:translation_surface_is_manifold} can be taken as a definition which is \emph{equivalent} to Definition \ref{def:translation_surface}.

\begin{definition}[Translation surface, 2\textsuperscript{nd} definition]\label{def:translation_surface_2}
    A \emph{translation surface} $X$ is a compact 1-dimensional complex manifold with a finite collection of points $\Sigma \subseteq X$, called the \emph{singularity set}, and an atlas $\mathcal{A} = \{(U_i, \phi_i)\}_{i \in I}$ of $X$ such that the transition maps between any two charts whose domains are disjoint of $\Sigma$ are translations, i.e., of the form $z \mapsto z + c$ for some $c\in \C$, and transition maps between a chart whose domain contains a singularity and a domain who contains no singularity are of the form $z \mapsto z^k$ for some $k \in \N$. 
\end{definition}
The fact that this definition is equivalent to Definition \ref{def:translation_surface} is proven in \cite{athreyamasur2023translationsurfaces}. Morally, to obtain a collection of polygons $\mathcal{P}$ given an atlas as in Definition \ref{def:translation_surface_2}, we join the singularities by enough geodesics so that cutting open the surface along these geodesics yields components that are simply connected and, since geodesics are locally straight lines by the fact that transition maps are translations, these components have straight lines as boundaries and hence are polygons. This construction of course yields again an atlas, where the charts are exactly as described in Proposition \ref{prop:translation_surface_is_manifold}. Even though this atlas might not be \emph{identical} to the initial atlas, but they share a common \emph{maximal atlas}, hence the atlases are the same in the usual sense.

Lastly, we can also define translation surfaces using concepts from Complex Analysis. More precisely, a translation surface can be seen as a complex manifold $X$ equipped with a holomorphic differential $\omega$. Let us recall briefly the most important concepts. 

\begin{definition}[Holomorphic differential]
    For a complex manifold $X$, a \emph{holomorphic differential} $\omega$ assigns to each $x \in X$ a complex linear map $\omega_x$ from the tangent space to $X$ at $x$ to $\C$, with the condition that $\omega_x$ depends \emph{holomorphically} on $x$. More explicitly, for any $z \in T_xX \cong \C$, the map
    \begin{equation*}
        x \mapsto \omega_x(z)
    \end{equation*}
    is holomorphic. 
\end{definition}

We only consider the simple case where our complex manifold $X$ has (complex) dimension 1. The only complex linear maps from $\C$ to itself are maps of the form
\begin{equation*}
    z \mapsto \lambda z,
\end{equation*}
where $\lambda \in \C$. We can think of the identity map $\operatorname{id}\colon w \to w$ as a differential $\dd z$ on $\C$, i.e., $\dd z \colon w \mapsto w$ for any $z \in \C$. Locally (in a chart), we can then write the differential $\omega$ as $f(z)\,\dd z$ for some holomorphic map $f$. So the corresponding linear map in the point $z$ in the chart is given by $w \mapsto f(z) w$.  

We require an additional condition, namely, we need to ensure the compatibility of the charts. If two charts intersect and $\omega$ reads $f(z) \, \dd z$ in the first and $g(z)\, \dd z$ in the second chart, if we write $T(z)$ for the transition map then we require
\begin{equation*}
    f(T(z))\cdot T'(z) = g(z),
\end{equation*}
which we recognize as the usual chain rule. Put differently, if $x \in X$ and $\mathbf{v} \in T_xX$ and we have two different charts $\phi$ and $\psi$ at $x$ such that $\omega$ reads $f(z)\, \dd z$ in $\phi$ and $g(z) \, \dd z$ in $\psi$, then 
\begin{equation*}
    f(\phi(x))\phi'(x)\cdot \mathbf{v}=\omega_x(v) = g(\psi(x))\psi'(x)\cdot \mathbf{v}.
\end{equation*}
Since this is true for all $\mathbf{v} \in T_xX$, it follows that
\begin{equation*}
    f(\phi(x))\phi'(x) = g(\psi(x))\psi'(x).
\end{equation*}
The transition map is given by $T = \phi\circ\psi^{-1}$, so by setting $z = \psi(x)$ we obtain
\begin{equation*}
    f(T(z))\phi'(\psi^{-1}(z)) = g(z)\psi'(\psi^{-1}(z)) \quad \Leftrightarrow \quad f(T(z))\cdot T'(z) = g(z),
\end{equation*}
as above. Not every compact manifold admits a non-zero holomorphic differential, but translation surfaces do, i.e., the holomorphic differential $\dd z$ on the plane descends to a holomorphic differential on the quotient space $X = \mathcal{P}/_\sim$, with zeroes exactly at the vertices. The order of the zeroes is given by the angle of the conical singularities. Indeed, this can be seen as a consequence of the chain rule. If the differential reads $g(z)\,\dd z$ in a chart around a singularity, and it reads $f(z)\,\dd z$ around a regular chart, then the transition map is given by $z \mapsto z^k$ to that
\begin{equation*}
    g(z)\, \dd z = f(z^k)\, \dd z^k = kf(z^k) z^{k-1} \,\dd  z, 
\end{equation*}
hence $g$ has a zero of order $k-1$ at the singularity. 

\begin{definition}[Translation surface, 3\textsuperscript{rd} definition]\label{def:translation_surface_3} A translation surface $(X, \omega)$ is a compact 1-dimensional complex manifold $X$ together with a choice of non-zero holomorphic differential $\omega$. 
\end{definition}

The above discussion essentially shows how to obtain a holomorphic differential given a translation surface in the form of an atlas with the appropriate transition maps. Let us sketch how given a holomorphic differential $\omega$ on a manifold $X$ we can recover an atlas as in Definition \ref{def:translation_surface_2}. Expositions providing further details can be found in \cite{athreyamasur2023translationsurfaces} or in \cite{massart2022short}.

We proceed in two steps. First, we want to obtain charts where $\omega$ reads $\dd z$. Suppose $\omega$ reads $g(z) \, \dd z$ in some chart $(U, \phi)$ at $x_0 \in X$ with $\phi(x_0) = 0$, which we may assume by translating the chart if necessary, and suppose that $g(0) \neq 0$ so that $x_0$ is a regular point of $\omega$. 

Then, we may take a primitive of $\omega$ as a new chart at $x_0$, i.e., we set
\begin{align*}
    \psi \colon V &\to \C, \\
    x &\mapsto \int_{x_0}^x \omega,
\end{align*}
where $V$ is some neighborhood of $x_0$ with $V \subseteq U$. Since $\omega$ does not vanish at $x_0$, the map $\psi$ is locally biholomorphic. Let $G$ be the local primitive of $g$, defined in $V$, such that $G(0) = 0$, which exists since we may assume that $V$ is simply connected. Note that
\begin{equation*}
    \psi(x) = G(\phi(x)),
\end{equation*}
which implies in particular that $G$ \emph{is} the transition map between the charts $\psi$ and $\phi$. Assuming that $\omega$ reads $f(z) \, \dd z$ in the chart $(V, \psi)$, by the chain rule we obtain
\begin{equation*}
    f(G(z))\underbrace{G'(z)}_{=g(z)} = g(z),
\end{equation*}
hence $f$ is constantly equal to 1 and $\omega$ reads $\dd z$ in $(V, \psi)$. 

Secondly, we will show that transitions between two such charts constructed in the first step are translations. Assume we have two different charts at different regular points $x$ and $y$ whose domains intersect. Assume further that $\omega$ reads $\dd z$ in both charts. Let $T$ be the transition map, then we have
\begin{equation*}
    \underbrace{f(T(z))}_{\equiv 1} T'(z) = \underbrace{g(z)}_{\equiv 1},
\end{equation*}
so that $T' \equiv 1$ showing that $T$ is a translation.

We will also mention the related concept of \emph{half}-translation surface, since it will play a role in our discussion of hyperellipticity in section \ref{sec:diagonal_changes}. In the same way as above for translation surfaces, taking different viewpoints leads to distinct but equivalent definitions also for half-translation surfaces. 

\begin{definition}[Half-translation surface]
    Given a collection $\mathcal{P}$ of polygons with edges $\mathcal{E}$ as in Definition \ref{def:translation_surface}, a \emph{half-translation surface} is obtained by identifying two edges of a pairing $(e_1, e_2), e_i \in \mathcal{E}$ where the orientation does not necessarily need to be preserved. Put differently, we may identify two edges not just by translations but by translations and a rotation by $\pi$. 

    Equivalently, we can define a half-translation surface via an atlas $\mathcal{A} = \{(U_i, \phi_i)\}_{i \in I}$ where the transition maps between regular charts are given by \emph{half-translations}, i.e., by maps of the form
    \begin{equation*}
        z \mapsto \pm z + c,
    \end{equation*}
    where $c \in \C$ is a constant. 

    A third equivalent way to define a half-translation surface is given by endowing a compact surface $X$ with a \emph{quadratic} holomorphic differential $q$. In local coordinates, $q$ reads
    \begin{equation*}
        f(z) \, \dd z^2 \coloneqq f(z) \, (\dd z \otimes \dd z),
    \end{equation*}
    where $f$ is holomorphic.
\end{definition}

\begin{figure}[ht]
    \centering
\begin{tikzpicture}[scale=2.3]
    \draw[thick, black, postaction={
    decorate, decoration={ markings, mark=at position 0.5 with{\arrow{latex}}}}] (0,0) -- (0,1);
    \draw[thick, black, postaction={
    decorate, decoration={ markings, mark=at position 0.5 with{\arrow{latex}}}}] (2,0) -- (2,1);
    \draw[thick, black, postaction={
    decorate, decoration={ markings, mark=at position 0.45 with{\arrow{latex}}, mark=at position 0.55 with {\arrow{latex}]}}}] (1,1) -- (2,1);
    \draw[thick, black, postaction={
    decorate, decoration={ markings, mark=at position 0.45 with{\arrow{latex}}, mark=at position 0.55 with {\arrow{latex}]}}}] (1,1) -- (0,1);
    \draw[thick, black, postaction={
    decorate, decoration={ markings, mark=at position 0.47 with{\arrow{to}}, mark=at position 0.53 with {\arrow{to}]}}}] (1,0) -- (2,0);
    \draw[thick, black, postaction={
    decorate, decoration={ markings, mark=at position 0.47 with{\arrow{to}}, mark=at position 0.53 with {\arrow{to}]}}}] (1,0) -- (0,0);
    \foreach \vertex in {(0,0), (0,1), (1,0), (1,1), (2,0), (2,1)}
        \filldraw[white, draw=black] \vertex circle (0.7pt);
\end{tikzpicture}
\qquad\qquad
\begin{tikzpicture}[scale = 2.3]
    \draw[thick, black, postaction={
    decorate, decoration={ markings, mark=at position 0.5 with{\arrow{latex}}}}] (0,0) -- (1,0);
    \draw[thick, black, postaction={
    decorate, decoration={ markings, mark=at position 0.54 with{\arrow{latex}}, mark=at position 0.46 with {\arrow{latex}}}}] (2.2,0) -- (1.2,0);
    \draw[thick, black, postaction={
    decorate, decoration={ markings, mark=at position 0.5 with{\arrow{latex}}}}] (3.4,0) -- (2.4,0);
    \draw[thick, black, postaction={
    decorate, decoration={ markings, mark=at position 0.5 with{\arrow{to}}}}] (0,1) -- (1,1);
    \draw[thick, black, postaction={
    decorate, decoration={ markings, mark=at position 0.5 with{\arrow{to}}}}] (2.2,1) -- (1.2,1);
    \draw[thick, black, postaction={
    decorate, decoration={ markings, mark=at position 0.54 with{\arrow{latex}}, mark=at position 0.46 with {\arrow{latex}}}}] (3.4,1) -- (2.4,1);

    \draw[thick, black] (0,0) -- (0,1);
    \path (0,0.5) node[right] {$a$};
    \draw[thick, black] (1,0) -- (1,1);
    \path (1,0.5) node[left] {$c$};
    \draw[thick, black] (1.2,0) -- (1.2,1);
    \path (1.2,0.5) node[right] {$b$};
    \draw[thick, black] (2.2,0) -- (2.2,1);
    \path (2.2, 0.5) node[left] {$a$};
    \draw[thick, black] (2.4,0) -- (2.4,1);
    \path (2.4,0.5) node[right] {$c$};
    \draw[thick, black] (3.4,0) -- (3.4,1);
    \path (3.4,0.5) node[left] {$b$};

    \foreach \vertex in {(0,0), (1,0), (1.2,0), (2.2,0), (2.4,0), (3.4,0),(0,1), (1,1), (1.2,1), (2.2,1), (2.4,1), (3.4,1)}
        \filldraw[white, draw = black] \vertex circle (0.7pt);
\end{tikzpicture}
    \caption{Two examples of half-translation surfaces. Topologically, the left surface corresponds to a sphere with 4 conical singularities of cone angle $\pi$ each, while the right surface corresponds to a surface of genus $\mathbf{g} = 1$ with two conical singularities of cone angle $\pi$ and one singularity of cone angle $4\pi$.}
    \label{fig:half_translation_surfaces}
\end{figure}
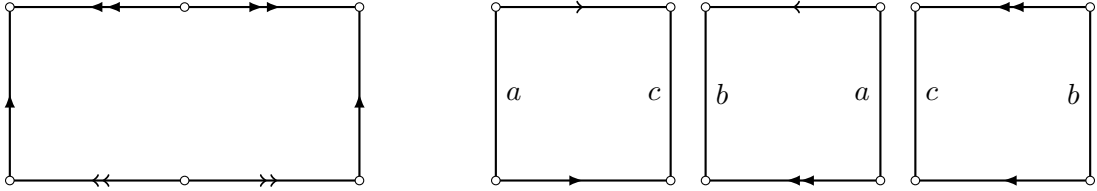

Singularities on half-translation surfaces are conical as well, where the cone angles are integer multiples of $\pi$, i.e., they are of the form $k \pi$ for $k \in \N_{\geq 1}$. Two examples of half-translation surfaces which are not translation surfaces can be seen in Figure \ref{fig:half_translation_surfaces}.

\subsection{Basic Concepts}
First, we want to derive a formula for the genus $\mathbf{g}$ of a translation surface $X$, which depends on the number and type of conical singularities.

\begin{proposition}[Genus formula]\label{prop:genus_formula}
    Let $X$ be a translation surface with genus $\mathbf{g}$ and singularity set $\Sigma = \{p_1, \ldots, p_n\}$ where $p_i$ has cone angle $2\pi k_i$. Then,
    \begin{equation*}
        \sum_{i = 1}^n k_i = 2\mathbf{g} - 2 + n.
    \end{equation*}
\end{proposition}
\begin{proof}
    We'll consider the special case when $X$ is obtained from a single polygon $P$, i.e., $X = P/_\sim$. Suppose that $P$ is a $2d$-gon. We have
    \begin{equation*}
        \sum_{i = 1}^n 2\pi k_i = \pi (2d - 2),
    \end{equation*}
    by the formula for the sum of internal angles of a $2d$-gon. This is equivalent to
    \begin{equation}\label{eq:genus_formula_1}
        \sum_{i = 1}^n k_i = d-1.
    \end{equation}
    By Euler's formula, we have
    \begin{equation}\label{eq:genus_formula_2}
        2 - 2\mathbf{g} = \chi(X) = \#\mathcal{F} - \#\mathcal{E} + \#\mathcal{V} \stackrel{*}{=} 1 - d + n,
    \end{equation}
    where $\#\mathcal{F}, \#\mathcal{E}$ and $\#\mathcal{V}$ denote the cardinalities of the faces, edges and vertices of the polygon and $*$ follows from the fact that we have one face, $d$ edges and $n = 2d$ vertices. Combining equations \eqref{eq:genus_formula_1} and \eqref{eq:genus_formula_2} finishes the proof. 
\end{proof}

\begin{remark}
    One can obtain this formula also as a direct corollary of the Gauss--Bonnet theorem (Theorem \ref{thm:gauss_bonnet}), see for instance \cite{athreyamasur2023translationsurfaces}. 
\end{remark}

We can rearrange the formula from the statement to read
\begin{equation*}
    2\pi(2-2\mathbf{g}) = -\sum_{i = 1}^n (k_i - 1) 2\pi.
\end{equation*}
The left-hand side gives exactly $2\pi \chi(X)$, while the terms $(k_i - 1)2\pi$ in the sum on the right-hand side can be seen as the excess angle at $p_i$, i.e., the amount of additional angle that is present at this singularity preventing it from being a Euclidean point. Recalling the Gauss--Bonnet theorem (Theorem \ref{thm:gauss_bonnet}) this now justifies the interpretation of the singularities $p_i$ as Dirac-masses we mentioned above. 

There are several structures on translation surfaces naturally induced by $\R^2$. First, if $X$ is a translation surface, the area of a subset $E \subseteq X$ is given by the Euclidean Area on the polygons $P_i \in \mathcal{P}$ in the definition of $X$, i.e., we can write
\begin{equation*}
    \operatorname{Area}(E) = \sum_{i = 1}^n \operatorname{Area}(E\cap P_i).
\end{equation*}
Moreover, as we have a Riemannian metric on $X$ induced by the Euclidean metric. More explicitly, at regular points $p \in X \setminus \Sigma$, we can write the metric as
\begin{equation*}
    \dd s^2 = \dd x^2  + \dd y^2 \stackrel{*}{=} \dd r^2 + r^2 d\phi^2,
\end{equation*}
where $*$ denotes a change to polar coordinates. At singular points $ p \in \Sigma$ we have
\begin{equation}\label{eq:conical_metric}
    \dd s^2 = r^k(\dd r^2 + (r\dd \phi)^2),
\end{equation}
where $k \in \N_{\geq 0}$ is such that the cone angle at $p$ is $2\pi(k+1)$. Note that if $k = 0$, then \eqref{eq:conical_metric} reduces again to the Euclidean metric. This is for example the case for the unique conical singularity of the torus. Singularities of this type are sometimes referred to as \emph{fake} singularities and are treated as regular points.

\begin{figure}[ht]
    \centering
    \begin{tikzpicture}
  \node (octagon) [regular polygon, regular polygon sides=8, minimum size=4cm, draw, thick] at (0,0) {};
  
  \node[above] at (octagon.side 1) {$a$};
  \node[below] at (octagon.side 5) {$a$};
  \node[yshift = 5pt, xshift = -5pt] at (octagon.side 2) {$b$};
  \node[yshift = -5pt, xshift = 5pt] at (octagon.side 6) {$b$};
  \node[left] at (octagon.side 3) {$c$};
  \node[right] at (octagon.side 7) {$c$};
  \node[yshift = -5pt, xshift = -5pt] at (octagon.side 4) {$d$};
  \node[yshift = 5pt, xshift = 5pt] at (octagon.side 8) {$d$};

  \draw[thick, line cap = round, purple] ($(octagon.corner 5) + (0.02,0.02)$)  -- ($(octagon.corner 7) - (0.02,0.0)$);

  \draw[thick, line cap = round, purple] ($(octagon.corner 4) + (0.02,0)$) -- ($(octagon.corner 8) - (0.02,0.0)$);

  \draw[thick, line cap = round, purple] ($(octagon.corner 3) + (0.02,0)$) -- ($(octagon.corner 1) - (0.01,0.01)$);

\end{tikzpicture}\qquad\qquad
\begin{tikzpicture}
  \node (octagon) [regular polygon, regular polygon sides=8, minimum size=4cm, draw, thick] at (0,0) {};
  
  \node[above] at (octagon.side 1) {$a$};
  \node[below] at (octagon.side 5) {$a$};
  \node[yshift = 5pt, xshift = -5pt] at (octagon.side 2) {$b$};
  \node[yshift = -5pt, xshift = 5pt] at (octagon.side 6) {$b$};
  \node[left] at (octagon.side 3) {$c$};
  \node[right] at (octagon.side 7) {$c$};
  \node[yshift = -5pt, xshift = -5pt] at (octagon.side 4) {$d$};
  \node[yshift = 5pt, xshift = 5pt] at (octagon.side 8) {$d$};

  \draw[thick, line cap = round, green!10!blue!95!red!70] ($(octagon.corner 5) + (0.02,0.02)$)  -- ($(octagon.corner 2)+(0.5,-0.03)$);
  \draw[thick, line cap = round, green!10!blue!95!red!70] ($(octagon.corner 5) + (0.5,0.02)$)  -- ($(octagon.corner 2)+(1,-0.03)$);
  \draw[thick, line cap = round, green!10!blue!95!red!70] ($(octagon.corner 5) + (1,0.02)$)  -- ($(octagon.corner 1)-(0,0.02)$);
\end{tikzpicture}
    \caption{The unique conical singularity $p$ of the translation surface $X$ induced by the regular octagon has 3 outgoing lines in direction $\frac{\pi}{8}$ as seen on the left. These lines are also saddle connections. An additional saddle connection is depicted on the right.}
    \label{fig:octagon_outgoing}
\end{figure}
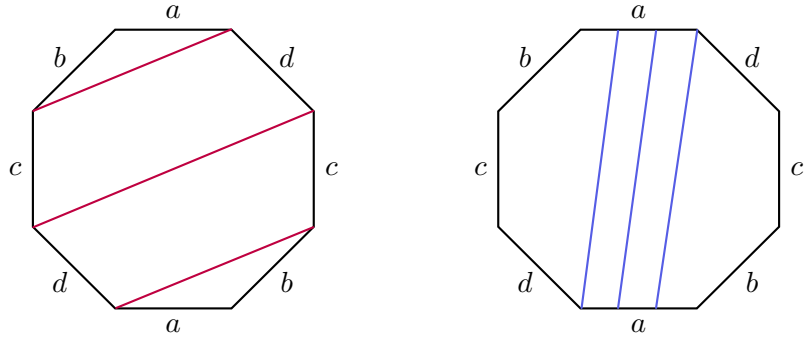

Another structure induced by $\R^2$ is the Euclidean distance, where the length of a segment is computed by pulling back the Euclidean distance in $\R^2$. Lastly, let us stress again that a translation surface comes with a notion of direction $\theta \in S^1$ which is preserved by the transition maps since they are translations. We adapt the convention that $\theta = 0$ is the horizontal direction pointing right and $\theta = \frac{\pi}{2}$ is the vertical direction pointing upwards. Note that at singular points $p \in \Sigma$, any direction $\theta$ may not be unique. If the cone angle at $p$ is $2\pi k$, then there are $k$ different outgoing lines in direction $\theta$ for all $\theta \in S^1$. As an example, consider Figure \ref{fig:octagon_outgoing}, where there are 3 outgoing lines in direction $\frac{\pi}{4}$ at the unique conical singularity $p$ of cone angle $6\pi$. 

From a dynamical point of view, a natural object to study on a surface is the geodesic flow. Intuitively, we are interested in the outcome of a point moving in a straight line on the surface. 

\begin{definition}[Linear flow]
    For all $\theta \in S^1$, we define the \emph{linear flow in direction $\theta$} on the translation surface $X$ as the solution of the ordinary differential equation
    \begin{equation*}
        \begin{cases}
            x'(t) = \cos \theta, \\
            y'(t) = \sin \theta,
        \end{cases}
    \end{equation*}
    which is unique given an initial condition $p = (x_0, y_0) \in X$. We will denote the solution by $\varphi_\R^\theta = (\varphi_t^\theta)_{t \in \R}$. Put differently, the linear flow maps a point $p$ to the point $\varphi_t^\theta(p)$, which is obtained by moving the point in a straight line in direction $\theta$ at unit speed for time $t$.
\end{definition}
Note that the linear flow is no longer well-defined if a singularity is reached at some time $t_0$, since there are multiple outgoing straight lines. Often, we will be interested in the \emph{vertical} linear flow, i.e., the flow associated to the angle $\theta = \frac{\pi}{2}$, which we will denote simply by $\varphi_\R$. Since on a translation surface, geodesics are straight lines in a literal sense, we may use the terms \enquote{linear flow} and \enquote{geodesic flow} interchangeably. 

\begin{definition}[Regular (forward) trajectory]
    A point $p \in X$ has a \emph{regular trajectory} for $\varphi_\R^\theta$, if $(\varphi_t^\theta)_{t\in \R} \cap \Sigma = \emptyset$. The point $p$ is said to have a regular \emph{forward} trajectory for $\varphi_\R^\theta$, if $(\varphi_t^\theta(p))_{t> 0} \cap \Sigma = \emptyset$, i.e., if the \emph{forward} trajectory is well-defined for all $t > 0$.  
\end{definition}
\begin{remark}
    A singular point $p \in \Sigma$ can have a regular forward trajectory. 
\end{remark}

The following proposition shows that a generic point on $X$ is regular.

\begin{proposition}[Generic points are regular]\label{prop:generic_points_regular}
    Almost every point on $X$ is regular, i.e., if we write $\mathcal{N}(X) = \{p \in X \mid p \text{ has a non-regular trajectory}\}$, then
    \begin{equation*}
        \operatorname{Area}(\mathcal{N}(X)) = 0.
    \end{equation*}
    Consequently, $\varphi_\R^\theta$ is defined almost everywhere. 
\end{proposition}
\begin{proof}
    Since there are only finitely many singular points, the set $\mathcal{N}(X)$ can be written as a finite union of lines in $\R^2$ which projects to a countable union of line segments on the surface, which has (2-dimensional) Lebesgue measure zero.
\end{proof}

The following two concepts, \emph{saddle connections} and their \emph{systoles}, which will be important in the sequel. 

\begin{definition}[Saddle connection]\label{def:saddle_connection}
    A \emph{saddle connection} is a geodesic for the metric induced by the Euclidean metric joining two singularities without any singularities in its interior. 
\end{definition}
In the polygonal representation of a translation surface $X$, we see that every edge is given by a saddle connection. In Figure \ref{fig:octagon_outgoing}, more saddle connections are pictured, in particular saddle connections which are \emph{not} given by the sides of the polygon.

To any saddle connection $\gamma$ we have a natural way to associate to it its length $\ell(\gamma)$, namely by integration along the saddle connection. 

\begin{definition}[Systole]\label{def:systole}
    Given a translation surface $X$ its \emph{systole} is given by
    \begin{equation*}
        \operatorname{sys}(X) \coloneqq \min \{\ell(\gamma) \mid \gamma \text{ a saddle connection on }X\},
    \end{equation*}
    that is, the systole of $X$ is the length of the shortest saddle connection. 
\end{definition}

\begin{remark}
    The systole is well-defined in that the quantity in the definition is indeed a minimum, not an infimum. One way to see this is through the concept of \emph{holonomy vectors}, which we introduce in section \ref{sec:diagonal_changes}, see in particular Definition \ref{def:holonomy_vector} and Theorem \ref{thm:holonomy_discrete}.
\end{remark}

We now present \emph{Keane's Theorem}, a result of great importance for the setting of sections \ref{sec:rauzy_veech}, \ref{sec:diagonal_changes} and \ref{sec:counting}. First, the following definition is needed.

\begin{definition}[Minimal]
    We say that a linear flow $\varphi_\R^\theta$ on a translation surface $X$ is \emph{minimal}, if every regular forward trajectory (i.e., starting at any point $p$) is dense in $X$.
\end{definition}

\begin{theorem}[Keane's Theorem]\label{thm:keanes}
    If there exists no saddle connection in direction $\theta$, then the linear flow $\varphi_\R^\theta$ is minimal.
\end{theorem}

A proof can be found in the appendix. Using a similar argument as in the proof of Proposition \ref{prop:generic_points_regular}, we see that there are only countably many saddle connections on any translation surface $X$. Thus, we immediately obtain the following corollary to Theorem \ref{thm:keanes}.

\begin{corollary}
    For any translation surface $X$ and almost all directions $\theta$ with respect to the natural Haar measure on $S^1$, the linear flow $\varphi_\R^\theta$ on $X$ in direction $\theta$ is minimal.
\end{corollary} \pagebreak
\thispagestyle{empty}
\
\pagebreak
\section{A General Principle: Studying Families of Surfaces}\label{sec:general_principle}
\thispagestyle{plain}

In this section, we will describe a general principle in the study of translation surfaces, which has proven to be extremely powerful in answering questions about the geometry or dynamics of some given translation surface $X$. 

We start with an informal exposition. Instead of studying the dynamical properties of the linear flow, say, on the surface $X$ directly, we view the surface as a single point that belongs to a \emph{family} of surfaces and define a new flow which acts on this family by deforming the surface $X$ in a specific way. Perhaps surprisingly, we will be able to obtain information about the original surface $X$ by studying its behavior under this deformation action, which would not be accessible otherwise. 

The rest of this section is dedicated to formalizing this informal description and giving an example in the form of \emph{Masur's Criterion} (Theorem \ref{thm:masurs_criterion}). 

\subsection{Deformation of Surfaces}

Consider a real $2\times 2$ matrix $A \in \operatorname{Mat}_{2 \times 2}(\R)$ with real coefficients given by
\begin{equation*}
    A = 
    \begin{bmatrix}
        a & b \\
        c & d
    \end{bmatrix},
\end{equation*}
and let $X = \mathcal{P}/_\sim$ be a translation surface. The matrix $A$ acts linearly on $\R^2$ by
\begin{equation*}
    \begin{bmatrix}
        x \\
        y
    \end{bmatrix}
    \mapsto
    A \cdot 
    \begin{bmatrix}
        x \\
        y
    \end{bmatrix}
    = 
    \begin{bmatrix}
        ax + by \\
        cx + dy
    \end{bmatrix}.
\end{equation*}
In particular, $A$ acts on the polygons $P_i \in \mathcal{P}$ in a way such that isometric parallel lines stay isometric and parallel. Therefore, the action of $A$ descends to the translation surface $X$.

\begin{definition}
    Let $X$ be a translation surface and $A \in \operatorname{Mat}_{2 \times 2}(\R)$. The matrix $A$ acts on $X$ by
    \begin{equation*}
        A \cdot X \coloneqq (A \cdot \mathcal{P}) /_\sim,
    \end{equation*}
    where $A \cdot \mathcal{P} = \{A\cdot P_1, \ldots, A\cdot P_n\}$ and $\sim$ denotes the equivalence relation from Definition \ref{def:translation_surface}.
\end{definition}

We are mostly interested in matrices that preserve the orientation as well as the natural area of $X$ under the action defined above, i.e., matrices that satisfy
\begin{equation*}
    \operatorname{Area}(A\cdot X)= \operatorname{Area}(X).
\end{equation*}

These matrices are exactly given by the \emph{special linear group} $\operatorname{SL}(2,\R)$, which is defined as
\begin{equation*}
    \operatorname{SL}(2,\R)= \left\{
    A = \begin{bmatrix}
        a & b \\
        c & d
    \end{bmatrix} \mid 
    a,b,c,d \in \R, \det A = 1
    \right\}.
\end{equation*}

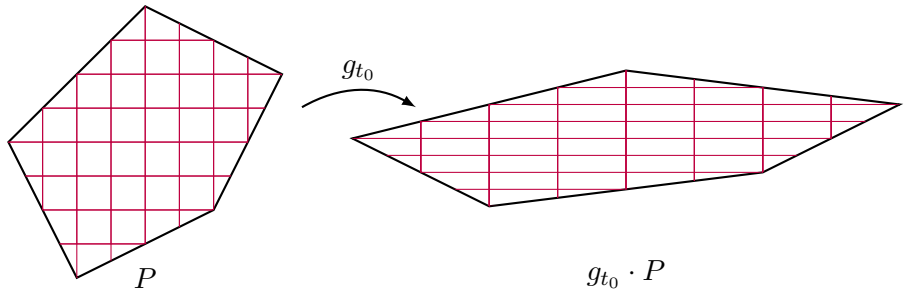
\begin{figure}[ht]
    \centering
    \begin{tikzpicture}[scale=0.9]
        \coordinate (A) at (0,0);
        \coordinate (B) at (2,2);
        \coordinate (C) at (4,1);
        \coordinate (D) at (3,-1);
        \coordinate (E) at (1,-2);
        
        \draw[thick] (A) -- (B) -- (C) -- (D) -- (E) -- cycle;
        \node at (2,-2) {$P$};

        \clip (A) -- (B) -- (C) -- (D) -- (E) -- cycle;
        
        \foreach \i in {-2,-1.5,...,4}
        \foreach \j in {-2,-1.5,...,2}
        {
            \draw[line cap = round, purple] (\i,\j) -- (\i+1,\j);
            \draw[line cap = round, purple] (\i,\j) -- (\i,\j+1);
        }
    \end{tikzpicture}
    \qquad
    \begin{tikzpicture}[scale=0.9]
        \coordinate (A) at (0,0);
        \coordinate (B) at (4,1);
        \coordinate (C) at (8,0.5);
        \coordinate (D) at (6,-0.5);
        \coordinate (E) at (2,-1);
        
        \draw[thick] (A) -- (B) -- (C) -- (D) -- (E) -- cycle;
        \node at (4,-2) {$g_{t_0} \cdot P$};

        \clip (A) -- (B) -- (C) -- (D) -- (E) -- cycle;
        
        \foreach \i in {-2,...,8}
        \foreach \j in {-2,-1.75,...,2}
        {
            \draw[line cap = round, purple] (\i,\j) -- (\i+1,\j);
            \draw[line cap = round, purple] (\i,\j) -- (\i,\j+1);
        }
    \end{tikzpicture}
    
    \begin{tikzpicture}[remember picture, overlay]
        \draw[->, bend left=30, thick, >=latex] (-2,3) to node[above] {$g_{t_0}$} (-0.5,3);
    \end{tikzpicture}
    \caption{The matrix $g_{t_0}$ acting on a polygon for $t_0 > 0$.}
    \label{fig:teichmueller_polygon}
\end{figure}

The subgroup consisting of \emph{diagonal} matrices is of particular interest, i.e., matrices of the form
\begin{equation}\label{eq:teichmueller_matrix}
    g_t = \begin{bmatrix}
        \e^t & 0 \\
        0 & \e^{-t}
    \end{bmatrix},
\end{equation}
for $t \in \R$. On the level of polygons, a matrix $g_t$ acts by expanding either the horizontal or vertical direction and by contracting the vertical or horizontal direction, respectively. The action of $g_{t_0}$ on a polygon for some $t_0 > 0$ is depicted in Figure \ref{fig:teichmueller_polygon}. In case the polygon represents a translation surface, the action of $g_t$ naturally descends to this surface as well.

Associated to the matrices $g_t$ we will define a \emph{flow}, i.e., we will be interested in the behavior of the translation surface $X$ as we let $t \to \pm \infty$. A priori, this would not be very interesting, since it is easy to see that any translation surface $X$ would simply degenerate under the action of $g_t$ for $t \to \pm \infty$, since the polygons $P_i \in \mathcal{P}$ would just become vertical or horizontal lines in the limit. We now introduce an equivalence relation, which will indicate when we regard two translation surfaces as equal.

\begin{definition}[Equivalent translation surfaces]\label{def:translation_surface_equivalence}
    We say that two translation surfaces $X_1 = \mathcal{P}_1/_\sim$ and $X_2 = \mathcal{P}_2/_\sim$ are \emph{equivalent}, if $\mathcal{P}_2$ can be obtained from $\mathcal{P}_1$ by a sequence of \emph{cut-and-paste} operations. These are given by either cutting a polygon $P$ along a straight line connecting two vertices which is contained in $P$ or by forming a larger polygon by combining $P_i$ and $P_j$ along an edge which is identified by $\sim$ by a translation.
\end{definition}

\begin{remark}
    We can formulate the above equivalence of translation surfaces also in terms of Definition \ref{def:translation_surface_2} using the language of Riemannian Geometry or in terms of Definition \ref{def:translation_surface_3} using the language of Complex Analysis. Details can be found in \cite{athreyamasur2023translationsurfaces}.
\end{remark}

\begin{figure}[ht]
    \centering
    \begin{tikzpicture}[scale=0.9]
        \coordinate (A) at (0,0);
        \coordinate (B) at (3,2);
        \coordinate (C) at (5,1);
        \coordinate (D) at (2,-1);
        
        \draw[thick] (A) -- (B) -- (C) -- (D) -- cycle;
        \draw[thick, dashed, line cap = round, green!10!blue!95!red!70] (D) -- (B);

        \node[rotate =70] at ($(B) + (0.08,0.22)$) {\Large \Leftscissors};
    \end{tikzpicture}
    \qquad
    \begin{tikzpicture}
        \coordinate (A) at (0,0);
        \coordinate (B) at (3,2);
        \coordinate (C) at (2,-1);

        \draw[thick] (A) -- (B) -- (C) -- cycle;
        \draw[green!10!blue!95!red!70, thick] (B) edge (C);

        \coordinate (X) at (3,-1);
        \coordinate (Y) at (4,2);
        \coordinate (Z) at (6,1);

        \draw[thick] (X) -- (Y) -- (Z) -- cycle;
        \draw[green!10!blue!95!red!70, thick] (X) edge (Y);
    \end{tikzpicture}
    \begin{tikzpicture}[remember picture, overlay]
        \draw[>=latex, ->, bend left=30, thick] (-7,2.5) to node[above] {cut} (-5,2.5);
        \draw[>=latex, ->, bend left=30, thick] (-6,0.5) to node[below] {paste} (-8,0.5);
    \end{tikzpicture}
    
    \caption{The rectangle and the two triangles induce equivalent translation surfaces.}
    \label{fig:equivalent_translation_surfaces}
\end{figure}
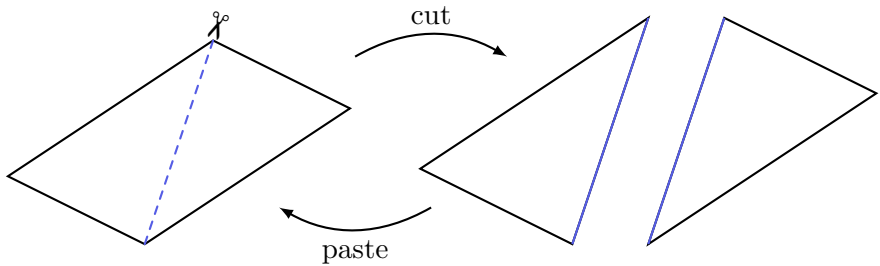

An example of two collections of polygons which induce equivalent translation surfaces can be seen in Figure \ref{fig:equivalent_translation_surfaces}. We cut the rectangle on the left along the diagonal to obtain the two triangles on the right, and conversely we can identify the two triangles along any side to obtain again a rectangle.

\subsection{Moduli Space and Strata}
The equivalence relation defined in the end of the last section allows us to properly define the notion of \emph{family of surfaces}, which we mentioned in the beginning of this chapter. For fixed $\mathbf{g} \in \N_{>0}, n \in \N$ and $k_1, \ldots k_n \in \N$ such that
\begin{equation*}
    \sum_{i = 1}^n k_i = 2 \mathbf{g} - 2 + n,
\end{equation*}
we can find examples of translation surfaces $X$ of genus $\mathbf{g}$ and $n$ singularities of cone angles 
\begin{equation*}
    2\pi k_1, \ldots, 2\pi k_n.
\end{equation*}
This motivates the following definition.

\begin{definition}[Strata and Moduli space]
    Let $\mathbf{g}, n$ and $k_1, \ldots, k_n$ as above. We define the \emph{stratum} $\mathcal{H}(k_1 - 1, \ldots, k_n - 1)$ as the set of all translation surfaces of genus $\mathbf{g}$ with $n$ conical singularities of cone angles $k_1, \ldots, k_n$, modulo the equivalence relation introduced in Definition \ref{def:translation_surface_equivalence}. Further, we will call the set $\mathcal{H}_\mathbf{g}$ of all translation surfaces of genus $\mathbf{g}$ (without specifying the number and type of conical singularities) the \emph{moduli space of translation surfaces of genus $\mathbf{g}$}. 
\end{definition}

\begin{remark}
    Referring to the definition of translation surfaces using holomorphic differentials (Definition \ref{def:translation_surface_3}), the set $\mathcal{H}_\mathbf{g}$ is also known as the \emph{moduli space of holomorphic (or Abelian) differentials} on surfaces of genus $\mathbf{g}$, equipped with the appropriate equivalence relation. 
\end{remark}

As explained in the previous section, the special linear group $\operatorname{SL}(2, \R)$ acts on a stratum $\mathcal{H}(k_1-1, \ldots, k_n-1)$ by
\begin{align*}
    X \mapsto A \cdot X. 
\end{align*}
This is in particular true for matrices of the form $g_t$ as in \eqref{eq:teichmueller_matrix}. Viewing $t$ as a real parameter, we can define a flow on the strata called the \emph{Teichmüller geodesic flow}.

\begin{definition}[Teichmüller geodesic flow]\label{def:teichmueller_geodesic_flow}
    The flow $g_\R = (g_t)_{t \in \R}$ obtained by the action of the 1-parameter diagonal subgroup
    \begin{equation*}
        \mathbf{A} = 
        \left\{
        g_t = 
        \begin{bmatrix}
            \e^t & 0 \\
            0 & \e^{-t}
        \end{bmatrix} \mid t \in \R
        \right\}
    \end{equation*}
    on $\mathcal{H}(k_1, \ldots, k_n)$ is called the \emph{Teichmüller geodesic flow.}
\end{definition}
The terminology is appropriate due to the fact that $(g_t \cdot X)_{t \in \R}$ is a geodesic with respect to the so-called \emph{Teichmüller metric}. This fact is known as \emph{Teichmüller's Theorem}. The definition of the metric and a proof of the theorem can be found in \cite{ahlfors2006lectures}. We can further define a \emph{directional} version $g_t^\theta$ of the Teichmüller geodesic flow, which is simply obtained by rotating the surface by the angle $-\theta$, applying $g_t$ and rotating it back. 

We can endow $\mathcal{H}(k_1-1, \ldots, k_n-1)$ with a measure in a natural way using so-called \emph{period coordinates}, which are local coordinates defined in a neighborhood of some fixed surface $X$. 

\begin{definition}[Period coordinates]\label{def:period_coordinates}
    Let $X = \mathcal{P}/_\sim$ be a translation surface and let $\mathcal{E}_1$ denote the set of edges that are on the left in the grouping in Definition \ref{def:translation_surface}, so that in particular the lengths and directions of \emph{all} edges in $\mathcal{E}$ are fully determined. For each $e_i \in \mathcal{E}_1$, if we write $p_i$ for its initial vertex and $q_i$ for its end vertex, the \emph{period} of $e_i$ is the complex number $c_i$ such that $q_i = p_i + c_i$.

    Given fixed sides $e_1, \ldots, e_n$ which determine $X$, we set
    \begin{equation*}
        \operatorname{Per}(X) = (c_1, \ldots, c_n) \in \C^n.
    \end{equation*}
    We call $\operatorname{Per}$ the \emph{period map}. In a neighborhood $U$ of the surface $X$, we have
    \begin{align*}
        \operatorname{Per} \colon U &\to \C^n, \\
        X & \mapsto (c_1(X), \ldots, c_n(X)).
    \end{align*}    
\end{definition}

Period coordinates allow us to define a measure, commonly known as the \emph{Masur--Veech measure}, as the pullback of the Lebesgue measure on $\C \cong \R^2$.

\begin{definition}[Masur--Veech measure]\label{def:masur_veech_measure}
    We define a measure $\mu$ on $\mathcal{H}(k_1-1, \ldots,k_n-1)$ by
    \begin{equation*}
        \mu(E) = \operatorname{Leb}(\operatorname{Per}(E)).
    \end{equation*}
\end{definition}
The measure defined above is well-defined, i.e., it is independent of the polygonal representations used to define the translation surfaces in the set $E$. Note that clearly, the measure is infinite, since we can scale any translation surface $X \in \mathcal{H}(k_1 - 1, \ldots, k_n - 1)$ by any positive real number $s \in \R_+$. Restricting ourselves to the space $\mathcal{H}^1(k_1 - 1 , \ldots, k_n - 1) \subseteq \mathcal{H}(k_1 - 1, \ldots, k_n - 1)$ of \emph{unit area} translation surfaces, which is a codimension 1 hypersurface, we do in fact obtain a \emph{finite} measure $\mu_1$. Since $\mathcal{H}^1(k_1 - 1 , \ldots, k_n - 1)$ is not compact, this is not obvious. A proof of this fact can be found in \cite{masur1991hausdorff}. We will adopt the convention that the measure $\mu_1$ has been normalized, so that the total space has measure 1.

\begin{lemma}[Masur--Veech measure is $\operatorname{SL}(2, \R)$ invariant]
    The measure $\mu$ is $\operatorname{SL}(2, \R)$-invariant, i.e., for all $A \in \operatorname{SL}(2, \R)$ we have
    \begin{equation*}
        \mu(A \cdot E) = \mu(E),
    \end{equation*}
    for all measurable subsets $E \subseteq \mathcal{H}(k_1 - 1, \ldots, k_n - 1)$.
\end{lemma}
\begin{proof}
    Note that
    \begin{equation*}
        \operatorname{Per}(A \cdot X) = A\cdot \operatorname{Per}(X),
    \end{equation*}
    for every $X\in \mathcal{H}(k_1-1,\ldots, k_n - 1)$ and $A \in \operatorname{SL}(2,\R)$, where $\cdot$ on the left-hand side denotes the action of $\operatorname{SL}(2, \R)$ on $\mathcal{H}(k_1 - 1, \ldots, k_n - 1)$ and on the right-hand side it denotes the usual matrix multiplication. Hence, for any measurable set $E$ we have
    \begin{equation*}
        \mu(A\cdot E) = \operatorname{Leb}\big(\operatorname{Per}(A\cdot E)\big) = \operatorname{Leb}\big(A\cdot \operatorname{Per}(E)\big) = \operatorname{Leb}\big(\operatorname{Per}(E)\big) = \mu(E),
    \end{equation*}
    where we use the fact that $\det A = 1$.
\end{proof}
\begin{remark}
    The induced measure $\mu_1$ enjoys the same property, as the exact same proof applies to $\mu_1$ as well.
\end{remark}

\subsection{Masur's Criterion}
We will now introduce the language necessary to state the result known as \emph{Masur's Criterion}. As mentioned above, This result will in particular provide us with an example of the general principle of answering questions about a single translation surface by exploring its orbit under a flow in its stratum. To start, we will need to introduce a notion of convergence on the strata $\mathcal{H}(k_1 - 1, \ldots, k_n - 1)$. We can give a simple description using the concept of \emph{flat triangulations}.

\begin{definition}[Flat triangulation]\label{def:flat_triangulation}
    A triangulation $\tau = \{T_1, \ldots, T_m\}$ of $X \in \mathcal{H}(k_1- 1, \ldots, k_n - 1)$ is a \emph{flat triangulation}, if the following are satisfied.
    \begin{enumerate}
        \item All vertices of all triangles $T_i \in \tau$ are singularities of $X$.
        \item No triangle $T_i \in \tau$ contains singularities in its interior.
        \item The sides of the triangles are Euclidean geodesics, which do not contain singularities in their interior.
        \item The edges of the triangles are identified by translations.
    \end{enumerate}
\end{definition}

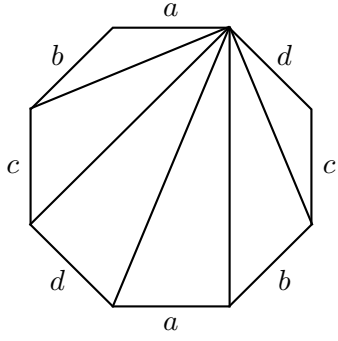
\begin{figure}
    \centering
    \begin{tikzpicture}
  \node (octagon) [regular polygon, regular polygon sides=8, minimum size=4cm, draw, thick] at (0,0) {};
  
  \node[above] at (octagon.side 1) {$a$};
  \node[below] at (octagon.side 5) {$a$};
  \node[yshift = 5pt, xshift = -5pt] at (octagon.side 2) {$b$};
  \node[yshift = -5pt, xshift = 5pt] at (octagon.side 6) {$b$};
  \node[left] at (octagon.side 3) {$c$};
  \node[right] at (octagon.side 7) {$c$};
  \node[yshift = -5pt, xshift = -5pt] at (octagon.side 4) {$d$};
  \node[yshift = 5pt, xshift = 5pt] at (octagon.side 8) {$d$};

    \foreach \i in {3,4,5,6,7}{
        \draw[thick] (octagon.corner 1) -- (octagon.corner \i);
    }
  
\end{tikzpicture}
    \caption{An example of a flat triangulation.}
    \label{fig:flat_triangulation}
\end{figure}

An example of a flat triangulation can be seen on the left-hand side of Figure \ref{fig:flat_triangulation}. Given a collection $\mathcal{P}$ as in the (first) definition of translation surface (Definition \ref{def:translation_surface}), we obtain a flat triangulation by joining enough vertices by Euclidean geodesics so that we are only left with triangles. Note that it is always possible to do so, meaning that we can always find a flat triangulation for any translation surface $X$. A proof of this fact can be found as Lemma 1.1 in \cite{masur1991hausdorff}. Moreover, since the Euler characteristic is well-defined on translation surfaces, the number of triangles $m$ in any triangulation $\tau$ is fixed, depending only on the number of singularities and on the genus of the surface.

\begin{definition}[Convergence in strata]\label{def:convergence_stratum}
    A sequence $(X_n)_{n \in \N} \subseteq \mathcal{H}(k_1 - 1, \ldots, k_n - 1)$ of translation surfaces in some fixed stratum converges to the translation surface $X_\infty \in \mathcal{H}(k_1 - 1, \ldots, k_n - 1)$ and we write $X_n \xrightarrow{n \to \infty} X_\infty$, if each $X_n$ for $n \in \N \cup \{\infty\}$ has a flat triangulation $\tau_n = \{T^n_1, \ldots, T^n_m\}$ in $m$ triangles, such that
    \begin{enumerate}
        \item $T_i^n \xrightarrow{n \to \infty} T_i^\infty$ as Euclidean triangles, i.e., the vertices converge as points in $\R^2$.
        \item The identification patterns, identifying edges of the triangles in $\tau_n$, are eventually constant.
        \item There exists a uniform constant $c>0$ such that $\ell(\gamma) \geq c$ for every saddle connection $\gamma$ in $X_n$, where $\ell(\gamma)$ denotes the length of the saddle connection.
    \end{enumerate}
\end{definition}
The third condition ensures that triangles cannot degenerate, or put differently, that singularities cannot collide. Without this assumption, the sequence of translation surfaces could converge to some surface in a different stratum. Combining this notion of convergence with the Teichmüller geodesic flow (Definition \ref{def:teichmueller_geodesic_flow}) enables us to define the concepts of \emph{divergence} and \emph{recurrence}.

\begin{definition}[Divergent and recurrent]\label{def:divergent_recurrent}
    A translation surface $X \in \mathcal{H}(k_1 - 1, \ldots, k_n - 1)$ is \emph{divergent in direction} $\theta$, if
    \begin{equation*}
        \operatorname{sys}(g_t^\theta \cdot X) \xrightarrow{t \to \infty} 0.
    \end{equation*}
    The surface $X$ is said to be \emph{recurrent in direction} $\theta$, if
    \begin{equation*}
        \operatorname{sys}(g_t^\theta \cdot X) \centernot{\xrightarrow{t\to\infty}} 0.
    \end{equation*}
\end{definition}

We are now ready to state \emph{Masur's Criterion}, which gives a convenient way to check if the linear flow in direction $\theta$ is uniquely ergodic, a property of some \emph{fixed} translation surface $X$, by studying its orbit under the (directional) Teichmüller geodesic flow $g_t^\theta$.

\begin{theorem}[Masur's Criterion]\label{thm:masurs_criterion}
    If a translation surface $X \in \mathcal{H}(k_1 - 1, \ldots, k_n - 1)$ is recurrent in direction $\theta$, then the linear flow in direction $\theta$ on $X$ is \emph{uniquely ergodic}. Equivalently, every regular trajectory on $X$ is \emph{equidistributed}.
\end{theorem}
A proof, based on \cite{ulcigrai2022surfacedynamics} and \cite{monteil2010introduction}, can be found in the appendix. The condition in Masur's Criterion is sufficient, but not necessary. In \cite{cheung2007slow}, Cheung and Eskin prove the following variant of Theorem \ref{thm:masurs_criterion}.

\begin{theorem}
    For any translation surface $X \in \mathcal{H}(k_1 - 1, \ldots, k_n-1)$, there is an $\varepsilon > 0$, depending only on the genus $\mathbf{g}(X)$ and the number of singularities $n$, such that if 
    \begin{equation*}
        \operatorname{sys}(g^\theta_t X) > t^{-\frac{\varepsilon}{2}} \cdot c
    \end{equation*}
    for some constant $c > 0$ and all $t>0$, then the linear flow in direction $\theta$ on $X$ is uniquely ergodic. 
\end{theorem}
In other words, if the divergence in the sense of Definition \ref{def:divergent_recurrent} is slow enough, then the linear flow on $X$ is still uniquely ergodic. 

\sloppy Results of this type provide ample motivation to understand the orbits of $g_t$ in the strata ${\mathcal{H}(k_1 - 1, \ldots, k_n - 1)}$. We will see more perspectives on the interplay between the linear flow on a single surface and the Teichmüller geodesic flow in the following sections. 
 \pagebreak

\section{IETs, Renormalization and Rauzy--Veech Induction}\label{sec:rauzy_veech}
\thispagestyle{plain}
In this chapter we want to introduce interval exchange transformations (usually abbreviated as IETs), give a conceptual description of a broadly used technique called \emph{renormalization} and give a brief description of a classical algorithm in the study of IETs and linear flows on surfaces called the \emph{Rauzy--Veech induction.}

This chapter is heavily inspired by chapter 5 of \cite{zorich2006flat}, so the main focus lies on geometric interpretations as well. Readers interested in a more formulaic approach to the same topic are referred to \cite{yoccoz2006continued}. Moreover, the general definition of an IET and many notational aspects are taken from \cite{viana2006ergodic}, which itself uses notation introduced in \cite{marmi2005cohomological}.

\subsection{IETs and Renormalization}\label{ch:IETs_and_Renormalization}

Before giving the definition of an IET, let us see how these kinds of transformations arise very naturally when studying the linear flow on a surface. As a toy example, we want to start with a torus, a translation surface of genus $\mathbf{g}=1$ with a single conical singularity of cone angle $2\pi$. In fact, we will assume for now that the torus is glued from the unit square. A very powerful approach when studying dynamical systems in general is to choose some section transverse to the flow of interest and consider the endomorphism of this section induced by the first return of the flow. This classical approach was first introduced by Henri Poincaré in \cite{poincare1881memoire}, for which reason such transverse section is often referred to as a \emph{Poincaré section} and the map itself is called the \emph{Poincaré first return map}. It is a remarkable fact, that one may study the first return map, which is a discretization of the continuous-time flow, to obtain results about the flow itself. General results on first return maps can be found in any textbook about Dynamical Systems, we recommend in particular \cite{viana2021differential}.

\begin{figure}[ht]
    \centering
    \begin{tikzpicture}
    \draw (0,0) rectangle (3.5,3.5);
    \draw[thick, color = green!10!blue!95!red!70] (0,0) -- (2,0);
    \draw[thick, noamblue] (0,3.5) -- (2,3.5); 
    \fill[noamblue] (0,3.5) circle (1pt); \fill[noamblue] (2,3.5) circle (1pt);
    \fill[green!10!blue!95!red!70] (0,0) circle (1pt); \fill[green!10!blue!95!red!70] (2,0) circle (1pt);

    \draw[black, dashed, postaction={
    decorate, decoration={ markings, mark=at position 0.8 with{\arrow{latex}}}}] (0,0) -- (1,3.5);
    \draw[black, dashed,postaction={
    decorate, decoration={ markings, mark=at position 0.8 with{\arrow{latex}}}}] (1,0) -- (2,3.5);
    \draw[black, dashed,postaction={
    decorate, decoration={ markings, mark=at position 0.8 with{\arrow{latex}}}}] (1.5,0) -- (2.5,3.5);
    \draw[black, dashed,postaction={
    decorate, decoration={ markings, mark=at position 0.8 with{\arrow{latex}}}}] (2,0) -- (3,3.5);
    \draw[black, dashed,postaction={
    decorate, decoration={ markings, mark=at position 0.8 with{\arrow{latex}}}}] (2.5,0) -- (3.5,3.5);
    \draw[black, dashed,postaction={
    decorate, decoration={ markings, mark=at position 0.5 with{\arrow{latex}}}}] (3,0) -- (3.5,1.75);
    \draw[black, dashed,postaction={
    decorate, decoration={ markings, mark=at position 0.6 with{\arrow{latex}}}}] (0,1.75) -- (0.5,3.5);

    \fill[color = green!10!blue!95!red!70, opacity = 0.2] (0,0) -- (1,0) -- (2,3.5) -- (1,3.5) -- cycle;
    \fill[color = purple, opacity = 0.2] (1,0) -- (1.5,0) -- (2.5,3.5) -- (2,3.5) -- cycle;
    \fill[color = teal, opacity = 0.2] (1.5,0) -- (2, 0) -- (3, 3.5) -- (2.5, 3.5) -- cycle; 
    \fill[color = purple, opacity = 0.2] (2,0) -- (2.5,0) -- (3.5,3.5) -- (3,3.5) -- cycle;
    \fill[color = teal, opacity = 0.2] (2.5,0) -- (3, 0) -- (3.5, 1.75) -- (3.5, 3.5) -- cycle; 
    \fill[color = teal, opacity = 0.2] (0,1.75) -- (0,3.5) -- (0.5, 3.5) -- cycle;
    \fill[color = purple, opacity = 0.2] (0,0) -- (1,3.5) -- (0.5, 3.5) -- (0, 1.75) -- cycle;
    \fill[color = purple, opacity = 0.2] (3,0) -- (3.5, 0) -- (3.5, 1.75) -- cycle;
\end{tikzpicture}
    \caption{A torus glued from the unit square with a transverse section from which the linear flow is emitted. }
    \label{fig:torus_IET_1}
\end{figure}
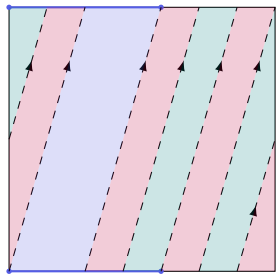

So let us fix a (non-horizontal) direction $\theta$ of the linear flow on the torus and pick a transverse section $S$, say a horizontal segment whose left endpoint coincides with the singularity of the torus. The linear flow emitted from this transverse section will wind around the torus before returning to the section eventually. This is illustrated in Figure \ref{fig:torus_IET_1}.

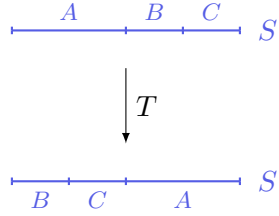
\begin{figure}[ht]
    \centering
    \begin{tikzpicture}
    \draw[color = green!10!blue!95!red!70, thick] (0,2) -- (3,2);
    \node[right, color =green!10!blue!95!red!70] at (3.1,2) {$S$};

    \draw[color =green!10!blue!95!red!70, thick] (0,2) - ++(0,0.05); 
    \draw[color =green!10!blue!95!red!70, thick] (0,2) - ++(0,-0.05); 

    \draw[color =green!10!blue!95!red!70, thick] (1.5,2) - ++(0,0.05); 
    \draw[color =green!10!blue!95!red!70, thick] (1.5,2) - ++(0,-0.05); 

    \draw[color =green!10!blue!95!red!70, thick] (9/4,2) - ++(0,0.05); 
    \draw[color =green!10!blue!95!red!70, thick] (9/4,2) - ++(0,-0.05); 

    \draw[color =green!10!blue!95!red!70, thick] (3,2) - ++(0,0.05); 
    \draw[color =green!10!blue!95!red!70, thick] (3,2) - ++(0,-0.05); 

    \draw[>=latex, ->] (1.5,1.5) to node[right] {$T$} (1.5,0.5);

    \draw[color = green!10!blue!95!red!70, thick] (0,0) -- (3,0);
    \node[right, color =green!10!blue!95!red!70] at (3.1,0) {$S$};

    \draw[color =green!10!blue!95!red!70, thick] (0,0) - ++(0,0.05); 
    \draw[color =green!10!blue!95!red!70, thick] (0,0) - ++(0,-0.05); 

    \draw[color =green!10!blue!95!red!70, thick] (3/4,0) - ++(0,0.05); 
    \draw[color =green!10!blue!95!red!70, thick] (3/4,0) - ++(0,-0.05); 

    \draw[color =green!10!blue!95!red!70, thick] (1.5,0) - ++(0,0.05); 
    \draw[color =green!10!blue!95!red!70, thick] (1.5,0) - ++(0,-0.05); 

    \draw[color =green!10!blue!95!red!70, thick] (3,0) - ++(0,0.05); 
    \draw[color =green!10!blue!95!red!70, thick] (3,0) - ++(0,-0.05); 

    \node[above, color =green!10!blue!95!red!70] at (3/4,2) {\footnotesize $A$};
    \node[above, color =green!10!blue!95!red!70] at (15/8,2) {\footnotesize $B$};
    \node[above, color =green!10!blue!95!red!70] at (21/8,2) {\footnotesize $C$};

    \node[below, color =green!10!blue!95!red!70] at (3/8,0) {\footnotesize $B$};
    \node[below, color =green!10!blue!95!red!70] at (9/8,0) {\footnotesize $C$};
    \node[below, color =green!10!blue!95!red!70] at (9/4,0) {\footnotesize $A$};

\end{tikzpicture}
    \caption{The map $T \colon S \to S$ induced by the linear flow on the torus.}
    \label{fig:torus_IET_2}
\end{figure}

We can see clearly how the induced map $T$ acts on the transverse section $S$. In the example pictured in Figure \ref{fig:torus_IET_1}, the segment gets cut up into three pieces exactly at the points from which the emitted flow eventually hits the singularity or the right endpoint of the segment. The points on each of these three pieces return to the section \enquote{together}, meaning that $T$ is given by rearranging these three intervals. The order of the \emph{intervals} is \emph{exchanged} by the first return \emph{map}. Forgetting about the torus and the linear flow on it and only considering the map from the section to itself, we can illustrate the return map as in Figure \ref{fig:torus_IET_2}.

Nothing in this construction depended on the fact that we are considering the torus obtained by gluing the unit square, thus it readily generalizes to more general translation surfaces of higher genus as well. Moreover, the construction encompasses exactly the idea of an interval exchange transformation.

\begin{definition}[Interval exchange transformation -- IET]
    Let $I \subseteq \R$ be an interval closed on the left and open on the right, and denote by $\mathcal{A}$ some alphabet with $d \geq 2$ symbols. Suppose further that $\{I_\alpha \colon \alpha \in \mathcal{A}\}$ is a partition of $I$ into subintervals, closed on the left and open on the right. An \emph{interval exchange map} (or \emph{IET}) is a bijective map from $I$ to $I$ which is a translation when restricted to each subinterval $I_\alpha$.
\end{definition}
We remark that it will not be crucial whether the endpoints belong to the intervals or not, since taking the measure theoretic perspective we can always dismiss sets of measure zero. The following lemma follows immediately from the definition. Here and in the remainder of the thesis, we use the notation $[d] \coloneqq \{1, 2, \ldots, d\}$ to denote the set of the first $d$ nonzero natural numbers. 
\begin{lemma}\label{lem:IET_characterization}
    Any IET is determined by combinatorial and length data as follows.
    \begin{enumerate}
        \item A pair $\boldsymbol\pi = (\topp, \bott)$ of bijections $\pi_\varepsilon \colon \mathcal{A} \to [d]$ for $\varepsilon \in \{\mathrm{top}, \mathrm{bot}\}$ describing the ordering of the subintervals $I_\alpha$ before and after the map is applied. 

        We will represent $\boldsymbol{\pi}$ as
        \begin{equation*}
            \boldsymbol{\pi} = \begin{pmatrix}
                \alpha_1^\mathrm{top} & \alpha_2^\mathrm{top} & \cdots & \alpha_d^\mathrm{top} \\
                \alpha_1^\mathrm{bot} & \alpha_2^\mathrm{bot} & \cdots & \alpha_d^\mathrm{bot}
            \end{pmatrix},
        \end{equation*}
        where $\alpha_j^\varepsilon = \pi_\varepsilon^{-1}(j)$ for $\varepsilon \in \{\mathrm{top}, \mathrm{bot}\}$ and $j \in [d]$. 
        \item A vector $\boldsymbol{\lambda} = (\lambda_\alpha)_{\alpha \in \mathcal{A}}$ with positive entries, where $\lambda_\alpha$ is the length of the subinterval $I_\alpha$. 
    \end{enumerate}
\end{lemma}

\begin{remark}
    For convenience, we will always assume that the left endpoint of $I$ is 0.
\end{remark}

Let us write the map $T \colon S \to S$ which we have obtained above as the first return to a transverse section on the torus in this new language.

\begin{example}
    For the IET depicted in Figure \ref{fig:torus_IET_2}, we have $\mathcal{A} = \{A, B, C\}$, the bijections $\boldsymbol{\pi}$ are given by
    \begin{align*}
        \topp &\colon \mathcal{A} \to \{1,2,3\}, \quad A \mapsto 1, B \mapsto 2, C \mapsto 3, \\
        \bott &\colon \mathcal{A} \to \{1,2,3\}, \quad A \mapsto 3, B \mapsto 1, C \mapsto 2.
    \end{align*}
    
    Thus, it can be represented by the combinatorial datum
    \begin{equation*}
        \boldsymbol{\pi} = 
        \begin{pmatrix}
            A & B & C \\
            B & C & A
        \end{pmatrix}
    \end{equation*}
    and the length datum
    \begin{equation*}
        \boldsymbol{\lambda} = (\lambda_A, \lambda_B, \lambda_C) = \left(\frac{1}{2}, \frac{1}{4}, \frac{1}{4}\right),
    \end{equation*}
    assuming that the section has length 1.
\end{example}

Note that the characterization using a pair of bijections $\boldsymbol{\pi} = (\topp, \bott)$ appears in some sense more complicated than strictly needed. What we mean is the following. We will call
\begin{equation*}
    p = \bott \circ \topp^{-1} \colon [d] \to [d]
\end{equation*}
the \emph{monodromy invariant} of the pair $\boldsymbol{\pi} = (\topp, \bott)$. Given any $(\boldsymbol{\pi}, \boldsymbol{\lambda})$ as in Lemma \ref{lem:IET_characterization} and any bijection $\phi \colon \hat{\mathcal{A}} \to \mathcal{A}$ we can define
\begin{equation*}
    \hat{\pi}_\varepsilon \coloneqq \pi_\varepsilon \circ \phi, \quad \text{for } \varepsilon \in \{\mathrm{top}, \mathrm{bot}\} \quad \text{ and }\quad  \hat{\lambda}_{\hat{\alpha}} \coloneqq \lambda_{\phi(\hat{\alpha})}, \quad \text{for }\hat{\alpha} \in \hat{\mathcal{A}}.
\end{equation*}
Then, $(\boldsymbol{\pi}, \boldsymbol{\lambda})$ and $(\boldsymbol{\hat\pi}, \boldsymbol{\hat{\lambda}})$ have the same monodromy invariant, and they characterize the same IET. Indeed, we have
\begin{equation*}
    \hat{\pi}_\mathrm{bot} \circ \hat{\pi}_\mathrm{top}^{-1} = \bott \circ \phi \circ \phi^{-1} \circ \topp^{-1} = \bott \circ \topp^{-1} = p.  
\end{equation*}
This just expresses the fact that relabeling the intervals does not change the IET in any meaningful way. So in principle we could just reduce everything to the case where $\mathcal{A} = [d]$ and $\topp = \operatorname{id}$, in which case $\bott$ coincides with the monodromy invariant $p$. However, the description of the algorithm below will be more elegant using the initial notation $\boldsymbol{\pi} = (\topp, \bott)$. 

In this thesis, we use IETs to obtain a convenient language to talk about certain Poincaré first return maps of the linear flow on a translation surface. But the study of the properties of IETs is an interesting topic in and of itself, of which we can only scratch the very surface. Nonetheless, we want to briefly present some results about IETs that have been obtained in the past decades, using, among other techniques, the Rauzy--Veech induction algorithm we will describe below.  

Dynamical systems in general can be roughly classified into the following three classes, see \cite{ulcigrai2014shearing}.
\begin{enumerate}
    \item \emph{Hyperbolic dynamical systems}: Nearby trajectories diverge \emph{exponentially fast.}
    \item \emph{Parabolic dynamical systems}: Nearby trajectories diverge \emph{subexponentially / polynomially}.
    \item \emph{Elliptic dynamical systems}: \emph{No} divergence, or possibly divergence that is \emph{slower than polynomially}. 
\end{enumerate}

It can be shown that IETs are examples of parabolic dynamical systems. 

In some aspects, they are closer to elliptic systems such as a rotation on a circle, than to hyperbolic systems such as geodesic flows on compact manifolds of constant negative curvature. For instance, the following theorem, established by A. Katok in \cite{katok1980interval}, states that IETs are never mixing, which is the case for circle rotations as well.

\begin{theorem}[\cite{katok1980interval}] 
    Let $T \colon I \to I$ be an IET and $\mu$ any Borel probability measure invariant with respect to $T$. Then, the map $T$ considered as an automorphism of the measure space $(I, \mu)$, is not mixing. 
\end{theorem}

We refer readers not familiar with Ergodic Theory to the excellently written textbook \cite{dajani2021first}. Notice however that more recently in \cite{avila2007weak}, it has been shown that a typical IET is weakly mixing. 

\begin{theorem}[\cite{avila2007weak}]
    Let $T \colon I \to I$ be an IET determined by $\mathcal{A} = [d]$ and $\boldsymbol{\pi} = (\operatorname{id}, \bott)$ such that $\bott$ is not a \emph{rotation}, i.e., we do not have $\bott(i+1) = \bott(i)+1 \mod d$. Then, for almost every $\boldsymbol{\lambda} = (\lambda_1, \ldots, \lambda_d) \subseteq \R_{>0}^d$, the IET $T \colon I \to I$ is weakly mixing. 
\end{theorem}

Note that \emph{if} $\bott$ is a rotation, then $T$ is conjugate to a rotation of the circle, hence it cannot be weakly mixing. 

Other properties of IETs put them closer to hyperbolic systems. For instance, typical IETs exhibit sensitive dependence on initial condition, which is the main characteristic of a \emph{chaotic} dynamical system and also a typical aspect in the hyperbolic case. Morally, it is clear that this chaotic behavior only develops rather slowly: The only possibility for two orbits that start close by to diverge, is that at some point, the iterates lie in two different intervals $I_\alpha$ and $I_\beta$ with $\alpha \neq \beta$. Until this happens, the two points will stay at a fixed distance. This is in stark contrast to the situation exhibited by hyperbolic systems, where we always find close points that \enquote{drift apart immediately}. For this reason, such systems are sometimes referred to as \emph{slowly chaotic dynamical systems} (see e.g., \cite{ulcigrai2021slow}). 

There is one property of IETs that lies at the very heart of the study of these maps, which can be seen as a form of \emph{self-similarity}. If we shorten the horizontal segment $S$ to a smaller segment $S'$, then the new segment $S'$ induces a new Poincaré map $T' \colon S' \to S'$ exactly in the same way the map $T\colon S \to S$ was induced. A priori, it could be that $T'$ is much more complicated than $T$ by, for instance, dividing $S'$ in many more intervals. However, the following lemma shows that the IET induced on a subinterval $S' \subseteq S$ decomposes the subinterval into (almost) the same number of segments as the initial IET.

\begin{lemma}[Self-similarity of IETs]\label{lem:self_similarity}
    Let $X$ be a translation surface with singularity set $\Sigma$, $\theta$ a generic direction, i.e., a direction in which the flow is minimal, and let $S$ be a horizontal segment whose left endpoint is some $\sigma \in \Sigma$. Then, the number of discontinuity points $D_T$ of the Poincaré first return map $T \colon S \to S$ is given by
    \begin{equation*}
        |D_T| = \sum_{\sigma \in \Sigma} \mathrm{deg}(\sigma) + \xi,
    \end{equation*}
    where $\operatorname{deg}(\sigma) = k$ if $\sigma$ is a conical singularity of angle $2\pi k$ (or equivalently, if $\sigma$ corresponds to a zero of order $k-1$ of the corresponding Abelian differential) and $\xi \in \{0,1\}$. Moreover, it is always possible to choose the length of the section $S$ in a way such that $\xi = 0$. 
\end{lemma}
\begin{proof}
    An interior point $x \in S$ is a point of discontinuity for $T$ if and only if one of the following happens.
    \begin{enumerate}[1)]
        \item The forward trajectory of $x$ hits a non-singular endpoint of $S$ before returning to $S$.
        \item The forward trajectory of $x$ hits a conical singularity before returning to $S$.
    \end{enumerate}
    Since we assume that the left endpoint of $S$ is a singularity, the first possiblilty above is responsible for $\xi \in \{0,1\}$ points of discontinuity. 

    A conical singularity with cone angle, say, $2\pi k$ has $k$ incoming trajectories. Following these $k$ trajectories backwards in time until we hit $S$ produces $k$ points of discontinuity on $S$. Hence, all the singularities together account for $\sum_{\sigma \in \Sigma} \operatorname{deg}(\sigma)$ points of discontinuity, so we conclude that the formula above holds. 

    For the last assertion, note that we can always choose the length of the segment in a way that either the backwards or forwards flow starting at the right endpoint of $S$ hits some conical singularity before returning to $S$. Call the right endpoint $r$ and suppose the backwards flow hits some $\sigma \in \Sigma$ before it hits $S$. In this case, the forward flow of any point in $S$ would hit $\sigma$ before hitting $r$ so that in this case $\xi = 0$. If the forward flow starting at $r$ hits some $\sigma \in \Sigma$, then the flow starting from any point $s \in S$ would first hit $r$ and not $\sigma$, reducing the number of discontinuities that occur by hitting a singularity by 1, so we have $\xi = 0$ in this case as well. 
\end{proof}

\begin{remark}
    Note that in the proof above, genericity of the direction is needed to make sure that the backwards flow of any singularity does in fact eventually hit $S$. This is guaranteed by Keane's Theorem (Theorem \ref{thm:keanes}).
\end{remark}

\begin{example}[IET on 4 intervals]\label{ex:octagon_IET}
    As an example of Lemma \ref{lem:self_similarity}, let us see how we can obtain an IET on 4 intervals induced by the linear vertical flow on a (tilted) regular octagon. The situation is pictured in Figure \ref{fig:octagon_IET}, where we have chosen a section $S$ such that the forwards flow of the right endpoint of $S$ hits the single conical singularity of cone angle $6\pi$. 

    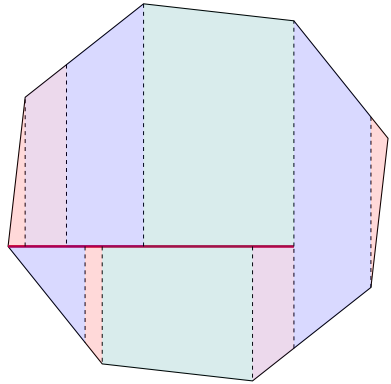
\begin{figure}[ht]
        \centering
        \begin{tikzpicture}[rotate around={16:(0,0)},transform shape, scale =1.3, baseline={(0,0)}]
    \def\radius{2cm} 
    
    \foreach \angle/\coord in {0/{A}, 45/{B}, 90/{C}, 135/{D}, 180/{E}, 225/{F}, 270/{G}, 315/{H}}
    \coordinate (\coord) at (\angle:\radius);
      
    \draw[name path = octagon] (A) -- (B) -- (C) -- (D) -- (E) -- (F) -- (G) -- (H) -- cycle;
    
    \path[name path=section] (E) --+ (-16:4);
    \path[name path=vert_5] (B) --+(-106:4);

    \draw[line width = 1pt,purple,name intersections={of=section and vert_5,by={Int5}}] (E) -- (Int5);

    \path[name path=vert_1] (D) --+ (-106:4);
    \path[name path=vert_2] (F) --+ (74:4);
    \path[name path=vert_3] (C) --+ (-106:4);
    \path[name path=vert_4] (G) --+ (74:4);
    \path[name path=vert_6] (H) --+ (74:4);

    \path[name intersections={of=vert_1 and section, by=Int1}];
    \path[name intersections={of=vert_2 and section, by=Int2}];
    \path[name intersections={of=vert_3 and section, by=Int3}];
    \path[name intersections={of=vert_4 and section, by=Int4}];
    \path[name intersections={of=vert_6 and octagon,by=Int6}];

    \draw[black, dashed, line width=0.3pt, dash pattern=on 1.5pt off 1.5pt] (Int1) -- (D);
    \draw[black, dashed, line width=0.3pt, dash pattern=on 1.5pt off 1.5pt] (Int2) -- (F);    
    \draw[black, dashed, line width=0.3pt, dash pattern=on 1.5pt off 1.5pt] (Int3) -- (C);
    \draw[black, dashed, line width=0.3pt, dash pattern=on 1.5pt off 1.5pt] (Int4) -- (G);
    \draw[black, dashed, line width=0.3pt, dash pattern=on 1.5pt off 1.5pt] (Int5) -- (B);
    \draw[black, dashed, line width=0.3pt, dash pattern=on 1.5pt off 1.5pt] (Int6) -- (H);

    \path[name path=vert_7] (Int5) --+ (-106:2);
    \path[name intersections={of=vert_7 and octagon, by=Int7}];

    \draw[black, dashed, line width=0.3pt, dash pattern=on 1.5pt off 1.5pt] (Int7) -- (Int5);    
    
    \path[name path = helper_line] (Int6) --+ (-157.5:4);
    \path[name intersections={of=helper_line and octagon, name=Int8}];
    \path[name path = vert_8] (Int8-2) --+ (74:3);
    \path[name intersections={of=vert_8 and section, by=Int9}];

    \draw[black, dashed, line width=0.3pt, dash pattern=on 1.5pt off 1.5pt] (Int9) -- (Int8-2);    

    \path[name intersections={of=vert_5 and octagon, name=Int10}];
    \path[name path = helper_line_2] (Int10-2) --+ (-247.5:4);
    \path[name intersections={of=helper_line_2 and octagon, name = Int11}];
    \path[name path = helper_line3] (Int11-1) --+ (-106:4);
    \path[name intersections={of=helper_line3 and section, name=Int12}];

    \draw[black, dashed, line width=0.3pt, dash pattern=on 1.5pt off 1.5pt] (Int11-1) -- (Int12-1);        
    
    \fill[red, opacity=0.15] (E) -- (Int1) -- (D) -- cycle;
    \fill[red, opacity=0.15] (H) -- (A) -- (Int6) -- cycle;
    \fill[red, opacity=0.15] (F) -- (Int2) --(Int9) -- (Int8-2) -- cycle;

    \fill[violet, opacity=0.15] (Int1) -- (D) -- (Int11-1) -- (Int12-1) -- cycle;
    \fill[violet, opacity=0.15] (G) -- (Int10-2) -- (Int5) -- (Int4) -- cycle;

    \fill[blue, opacity = 0.15] (Int12-1) -- (Int3) -- (C) -- (Int11-1) -- cycle;
    \fill[blue, opacity = 0.15] (Int10-2) -- (H) -- (Int6) -- (B) -- cycle;
    \fill[blue, opacity = 0.15](E) -- (Int8-2) -- (Int9) -- cycle;

    \fill[teal, opacity = 0.15] (Int3) -- (Int5) -- (B) -- (C) -- cycle;
    \fill[teal, opacity = 0.15] (Int2) -- (Int4) -- (G) -- (F) -- cycle;
\end{tikzpicture}
        \caption{An IET on 4 intervals induced by the vertical flow on a tilted regular octagon.}
        \label{fig:octagon_IET}
    \end{figure}
    We can see that the Poincaré map $T \colon S \to S$ has $\sum_{\sigma \in \Sigma} \operatorname{deg}(\sigma) = 3$ points of discontinuity. We obtain the IET pictured in Figure \ref{fig:octagon_IET_2}, which is given by the length and combinatorial data
    \begin{equation*}
        \boldsymbol{\pi} = \begin{pmatrix}
            A & B & C & D \\
            C & A & D & B
        \end{pmatrix},
    \end{equation*}
    and
    \begin{equation*}
        \boldsymbol{\lambda} = (\lambda_A, \lambda_B, \lambda_C, \lambda_D)
    \end{equation*}
    for some $\lambda_\alpha > 0, \alpha \in \{A,B,C,D\}$. 
    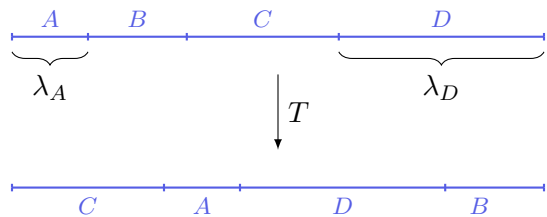
\begin{figure}[ht]
        \centering
        \begin{tikzpicture}
    \draw[color = green!10!blue!95!red!70, thick] (0,2) -- (7,2);
    
    \draw[color =green!10!blue!95!red!70, thick] (0,2) - ++(0,0.05); 
    \draw[color =green!10!blue!95!red!70, thick] (0,2) - ++(0,-0.05); 

    \draw[color =green!10!blue!95!red!70, thick] (1,2) - ++(0,0.05); 
    \draw[color =green!10!blue!95!red!70, thick] (1,2) - ++(0,-0.05); 

    \draw[color =green!10!blue!95!red!70, thick] (2.3,2) - ++(0,0.05); 
    \draw[color =green!10!blue!95!red!70, thick] (2.3,2) - ++(0,-0.05); 

    \draw[color =green!10!blue!95!red!70, thick] (4.3,2) - ++(0,0.05); 
    \draw[color =green!10!blue!95!red!70, thick] (4.3,2) - ++(0,-0.05);

    \draw[color =green!10!blue!95!red!70, thick] (7,2) - ++(0,0.05); 
    \draw[color =green!10!blue!95!red!70, thick] (7,2) - ++(0,-0.05);

    \draw[>=latex, ->] (3.5,1.5) to node[right] {$T$} (3.5,0.5);

    \draw[color = green!10!blue!95!red!70, thick] (0,0) -- (7,0);

    \draw[color =green!10!blue!95!red!70, thick] (0,0) - ++(0,0.05); 
    \draw[color =green!10!blue!95!red!70, thick] (0,0) - ++(0,-0.05); 

    \draw[color =green!10!blue!95!red!70, thick] (2,0) - ++(0,0.05); 
    \draw[color =green!10!blue!95!red!70, thick] (2,0) - ++(0,-0.05); 

    \draw[color =green!10!blue!95!red!70, thick] (3,0) - ++(0,0.05); 
    \draw[color =green!10!blue!95!red!70, thick] (3,0) - ++(0,-0.05); 

    \draw[color =green!10!blue!95!red!70, thick] (5.7,0) - ++(0,0.05); 
    \draw[color =green!10!blue!95!red!70, thick] (5.7,0) - ++(0,-0.05); 

    \draw[color =green!10!blue!95!red!70, thick] (7,0) - ++(0,0.05); 
    \draw[color =green!10!blue!95!red!70, thick] (7,0) - ++(0,-0.05); 

    \node[above, color =green!10!blue!95!red!70] at (1/2,2) {\footnotesize $A$};
    \node[above, color =green!10!blue!95!red!70] at (1.65,2) {\footnotesize $B$};
    \node[above, color =green!10!blue!95!red!70] at (3.3,2) {\footnotesize $C$};
    \node[above, color =green!10!blue!95!red!70] at (5.65,2) {\footnotesize $D$};

    \node[below, color =green!10!blue!95!red!70] at (1,0) {\footnotesize $C$};
    \node[below, color =green!10!blue!95!red!70] at (2.5,0) {\footnotesize $A$};
    \node[below, color =green!10!blue!95!red!70] at (4.35,0) {\footnotesize $D$};
    \node[below, color =green!10!blue!95!red!70] at (6.15,0) {\footnotesize $B$};

    \draw [decorate, decoration={brace, amplitude=5pt}] (1,1.8) -- node[below=5pt] {$\lambda_A$}  (0,1.8);
    \draw [decorate, decoration={brace, amplitude=5pt}] (7,1.8) -- node[below=5pt] {$\lambda_D$}  (4.3,1.8);
\end{tikzpicture}
        \caption{The IET on 4 intervals induced by the vertical flow pictured in Figure \ref{fig:octagon_IET}.}
        \label{fig:octagon_IET_2}
    \end{figure}
\end{example}

The following definition is motivated by the second part of Lemma \ref{lem:self_similarity}.

\begin{definition}[Minimal interval convention -- MIC]\label{def:mic}
    We say that a transverse segment $S$ of a linear flow on a translation surface satisfies the \emph{minimal interval convention}, or \emph{MIC}, if the induced IET $\colon T \to T$ has the minimal possible number
    \begin{equation*}
        n = \sum_{\sigma \in \Sigma} \operatorname{deg}(\sigma) +1
    \end{equation*}
    of subintervals.
\end{definition}
\begin{remark}
    We need to add 1 in the definition above when comparing to the formula in Lemma \ref{lem:self_similarity} because $m$ points of discontinuity correspond to $m+1$ intervals.
\end{remark}

Note that by using the genus formula (Proposition \ref{prop:genus_formula}), we can also relate the minimal number of intervals of the Poincaré map induced by a linear flow on a translation surface to the genus $\mathbf{g}$ by
\begin{equation*}
    n = 2 \mathbf{g} + |\Sigma| - 1.
\end{equation*}
For example, one can check that the regular octagon with the obvious identification of the sides corresponds to a translation surface of genus $\mathbf{g} = 2$ with one conical singularity (of cone angle $6\pi$), so that we indeed obtain $n = 4$ as the minimal number of subintervals.

\subsection{Rauzy--Veech Induction}

The self-similarity discussed in the previous section will be the most important property of IETs to define the Rauzy--Veech procedure, since it allows us to define a map from the \emph{space} of interval exchanges (of a certain number of intervals) to itself by finding a good way of choosing a smaller subinterval. Then, we \emph{renormalize} the new IET induced on the smaller section such that this smaller section becomes of the same length as the section we started with initially. Before we delve into the construction of this algorithm in full generality, we will initially apply the idea to the simple case of the torus, which makes the main idea behind the algorithm highly transparent.

\subsubsection{Finding Continued Fraction Expansion via Renormalization on the Torus} \label{sec:cf_torus}

To illustrate the technique of renormalization, let us show that it can be used to derive the continued fraction expansion of any number using the linear flow on a torus. The strategy explained here is largely taken from \cite{zorich2006flat}, where the same procedure is used to link the linear flow on the torus first to rotations on the circle and then to the Euclidean algorithm. It is well-known that the Euclidean algorithm and the standard algorithm to find continued fraction expansions (which we will define below) are equivalent, see e.g., \cite{anglin1995continued}. Continued fractions are studied extensively in Number Theory and we refer interested readers to the excellently written texts \cite{khinchin1963continued} and \cite{olds1963continued}. 

\begin{figure}[ht]
    \centering
    \begin{tikzpicture}
    \draw (0,0) rectangle (3,3);
    \draw[color = purple, thick] (0,0) -- (2,0); \draw[color = purple, thick] (1,3) -- (3,3);
    \draw[color = green!10!blue!95!red!70, thick] (2,0) -- (3,0); \draw[color = green!10!blue!95!red!70, thick] (0,3) -- (1,3);
    \draw[dashed] (0,0) -- (1,3); \draw [dashed] (2,0) -- (3,3);

    \node at (3.5,1.5) {$\cong$};

    \draw[color = purple, thick] (4,1) -- (6,1);
    \draw[color = purple, thick] (5,2) -- (7,2);
    \draw[color = green!10!blue!95!red!70, thick] (6,1) -- (7,1);
    \draw[color = green!10!blue!95!red!70, thick] (4,2) -- (5,2);

    \draw[>=latex, ->] (5.5,1.8) to node[right] {} (5.5,1.2);

    \draw (6,1) - ++(0,0.05);
    \draw (6,1) - ++(0,-0.05);

    \draw (4,1) - ++(0,0.05);
    \draw (4,1) - ++(0,-0.05);

    \draw (7,1) - ++(0,0.05);
    \draw (7,1) - ++(0,-0.05);

    \draw (4,2) - ++(0,0.05);
    \draw (4,2) - ++(0,-0.05);

    \draw (5,2) - ++(0,0.05);
    \draw (5,2) - ++(0,-0.05);

    \draw (7,2) - ++(0,0.05);
    \draw (7,2) - ++(0,-0.05);

    \node at (7.5, 1.5) {$\cong$};

    \draw[thick, color = purple] (9, 1.5) circle ++(0:1) arc (0:240:1);
    \draw[thick, color = green!10!blue!95!red!70] (9, 1.5) circle ++(0:1) arc (0:-120:1);

    \draw[dashed] (9,1.5) -- ++(0:1.1);
    \draw[dashed] (9,1.5) -- ++(-120:1.1);

    \node at (10.5,1.5) {$\mapsto$};

    \draw[thick, color = purple] (12, 1.5) circle ++(120:1) arc (120:360:1);
    \draw[thick, color = green!10!blue!95!red!70] (12, 1.5) circle ++(0:1) arc (0:120:1);

    \draw[dashed] (12,1.5) -- ++(120:1.1);
    \draw[dashed] (12,1.5) -- ++(0:1.1);

    \draw (9, 1.5) circle ++(45:0.3) arc (45:300:0.3); 
    \draw[>=latex, ->] (9, 1.5) ++(300:0.3) -- ++(45:0.2); 
\end{tikzpicture}
    \caption{An IET on two intervals can be identified with a rotation on $S^1$.}
    \label{fig:torus_IET_rotation}
\end{figure}
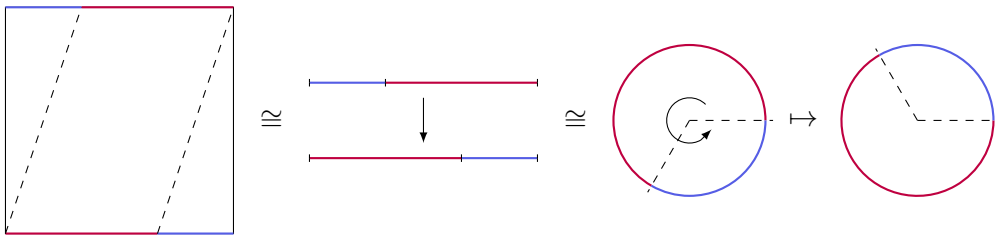

It is easy to see that we can link the linear flow an a torus to a rotation on the circle $S^1$ using an induced IET $T \colon S \to S$. Choosing as transverse section $S$ the whole base of the torus, it is obvious that $S$ satisfies the MIC, so that $T$ is an IET on two intervals. In this case, the combinatorial datum $\boldsymbol{\pi}$ is particularly simple, since the only possibility is to exchange the position of the two intervals. Moreover, if $\boldsymbol{\lambda} = (\lambda_A, \lambda_B)$ then we can write the Poincaré return map $T\colon S \to S$ explicitly as

\begin{equation*}
    T\colon \lambda \mapsto 
    \begin{cases}
        \lambda + \lambda_B \quad &\text{if } \lambda \in I_A,\\
        \lambda + \lambda_B - 1 \quad &\text{if } \lambda \in I_B
    \end{cases},
\end{equation*}
where we use that $\lambda_B - 1 = -\lambda_A$. This is exactly the rotation on $S^1$ by $\lambda_B$. This identification of the IET on two intervals and a rotation on $S^1$ is illustrated in Figure \ref{fig:torus_IET_rotation}.

For the sake of the reader, we also want to very briefly explain the standard procedure to obtain the continued fraction expansion of a number $\lambda \in [0,1)$ which involves the so-called \emph{Gauss map}.

\begin{definition}[Gauss map]\label{def:Gauss_map}
    The \emph{Gauss map} $G \colon [0,1) \to [0,1)$ is defined by
    \begin{equation*}
        G(\lambda) = \begin{cases}
            0 \quad & \text{if } \lambda = 0, \\
            \left\{\frac{1}{\lambda}\right\} & \text{if } \lambda \in (0,1),
        \end{cases}
    \end{equation*}
    where $\{\lambda\}$ denotes the fractional part of $\lambda$.  
\end{definition}

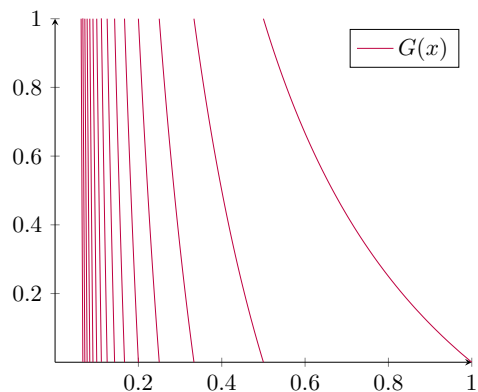
\begin{figure}[ht]
    \centering
    \begin{tikzpicture}[scale = 0.8]
\begin{axis}[
    axis lines=middle,
    xmin=0, xmax=1,
    ymin=0, ymax=1,
    domain=0:1,
    samples=100, 
    legend pos=north east, 
]

\foreach \n in {1,...,15}
    \addplot[
    purple,
    domain = 1/(\n+1):1/\n,
    samples = 100,
    unbounded coords = jump,
    ]{1/x - \n};
    \addlegendentry{$G(x)$}
    \end{axis}
\end{tikzpicture}
    \caption{Part of the graph of the Gauss map $G$.}
    \label{fig:Gauss_map}
\end{figure}

A part of the graph of the Gauss map can be seen in Figure \ref{fig:Gauss_map}. Note that the Gauss map is piecewise smooth on intervals of the form $P_n = \left(\frac{1}{n+1}, \frac{1}{n}\right]$ for $n \in \N_{\geq 1}$. We will call these intervals $P_n$ the \emph{branches of the Gauss map}.

It is convenient to use the notation
\begin{equation*}
    \lambda = a_0 + \cfrac{1}{a_1+\cfrac{1}{a_2+\cfrac{1}{a_3+\cfrac{1}{a_4+\cfrac{1}{\ddots}}}}} \coloneqq [a_0; a_1,a_2,a_3,a_4, \ldots] 
\end{equation*}
for the continued fraction expansion of $\lambda \in \R$. If $\lambda \in [0,1)$, we may omit the first entry and write
\begin{equation*}
    x = [a_1, a_2, a_3, \ldots]. 
\end{equation*}

Clearly, we lose no generality by only considering real numbers in the unit interval. To obtain the continued fraction expansion using the Gauss map, one must simply keep track of the branch that the iterates of $x$ belong to. More precisely, one can show that the digits $a_i$ of the continued fraction expansion are exactly given by the relation
\begin{equation*}
    a_i = k \quad \Leftrightarrow \quad G^i(\lambda) \in P_k. 
\end{equation*}

Now we want to see how one can use linear flows on the torus to obtain the same expansion. We will illustrate the procedure using an example. We start with the standard torus $\T = \R^2 / \Z^2$ and denote by $S$ the segment given by the lower horizontal line as on the left-hand side of Figure \ref{fig:translation_surfaces}.

Say, for instance, we want to compute the continued fraction expansion of $\lambda_1=  \sqrt{6} -2\approx 0.44949$. For this, we will consider the flow in direction $\theta_1 = \arctan \frac{1}{\lambda_1}$. It will be clear shortly why we choose exactly this angle. 

Now we consider trajectories which hit the single conical singularity. We start at the top right corner in the square representation of $\T$ and follow the flow backwards in time until we hit $S$, the base of the torus. 

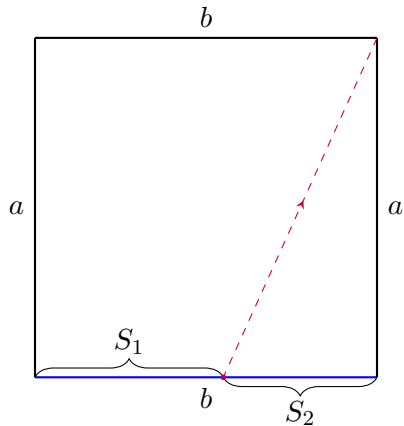
\begin{figure}[ht]
    \centering
\begin{tikzpicture}[scale = 1.5]
    \draw[thick, blue, name path=blue_line] (0,0) -- (3,0) node[midway, below, black] {$b$};
    \draw[thick] (0,0) -- (0,3) node[midway, left] {$a$};
    \draw[thick] (3,0) -- (3,3) node[midway, right] {$a$};
    \draw[thick] (0,3) -- (3,3) node[midway, above] {$b$};
    
    \draw[dashed, purple, postaction={decorate}, decoration={
        markings,
        mark=at position 0.5 with {\arrow{stealth reversed}}
    }, name path=purple_line] (3,3) -- ++(-114.2:3.3);
    
    \fill [purple, name intersections={of=blue_line and purple_line, by=intersection}] (intersection) circle (0.7pt);
    
    \begin{scope}[on background layer]
        \draw [decorate, decoration={brace, amplitude=7pt}] (0,0) -- node[above=5pt] {$S_1$} (intersection);
        \draw [decorate, decoration={brace, amplitude=7pt, mirror}] (intersection) -- node[below=5pt] {$S_2$} (3,0);
    \end{scope}
\end{tikzpicture}
    \caption{Flowing backwards from the conical singularity.}
    \label{fig:torus_renormalization_2.1}
\end{figure}

We obtain a decomposition $S = S_1 \sqcup S_2$ and a corresponding IET $T \colon S \to S$ which exchanges the two intervals $S_1$ and $S_2$. Note that we chose $\theta_1$ in a way such that the interval $S_2$ has length $\lambda_1$. We now consider the subinterval $S' \subseteq S$ which is given by the interval starting at the left endpoint of $S$ of length $\lambda_1$. By construction, this subinterval satisfies the MIC. Indeed, the left endpoint of $S'$ is itself a conical singularity, and the right endpoint of $S'$ hits the singularity under the backwards flow.

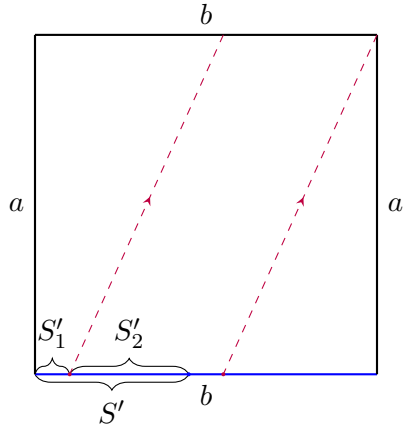
\begin{figure}[ht]
    \centering
    \begin{tikzpicture}[scale = 1.5]
    \draw[thick, blue, name path=blue_line] (0,0) -- (3,0) node[midway, below, black] {$b$};
    \draw[thick] (0,0) -- (0,3) node[midway, left] {$a$};
    \draw[thick] (3,0) -- (3,3) node[midway, right] {$a$};
    \draw[thick] (0,3) -- (3,3) node[midway, above] {$b$};
    
    \draw[dashed, purple, postaction={decorate}, decoration={
        markings,
        mark=at position 0.5 with {\arrow{stealth reversed}}
    }, name path=purple_line] (3,3) -- ++(-114.2:3.3);
    
    \fill [purple, fill opacity = 1, name intersections={of=blue_line and purple_line, by=intersection_1}] (intersection_1) circle (0.5pt);
    
    \path (0,0) -- (3,0) coordinate[pos={(1-0.5505)}] (new_point);
    \fill [blue] (new_point) circle (0.5pt);
    
    \begin{scope}[on background layer]
        \draw [decorate, decoration={brace, amplitude=6pt, mirror}] (0,0) -- node[below=6pt] {$S'$} (new_point);
    \end{scope}
    
    \path let \p1=(intersection_1), \p2=(new_point) in 
        coordinate (start_point_1) at (\x1,3);
    \path[name path global=purple_line_1] (start_point_1) -- ++(-114.2:3.3);
    \path[name intersections={of=blue_line and purple_line_1, by=intersection_1}];
    \draw[dashed, purple, postaction = {decorate}, decoration = {
    markings, mark=at position 0.5 with {\arrow{stealth reversed}}
    }] (start_point_1) -- (intersection_1);

    \fill [purple] (intersection_1) circle (0.5pt);
    \begin{scope}[on background layer]
        \draw[decorate, decoration={brace, amplitude=6pt}] (0,0) -- node[above=6pt] {$S'_1$} (intersection_1);
        \draw[decorate, decoration={brace, amplitude=6pt}] (intersection_1) -- node[above=6pt] {$S'_2$} (new_point);
    \end{scope}

    \end{tikzpicture}
    \caption{Repeating the same strategy for the shorter segment $S'$, which induces an IET $T' \colon S' \to S'$.}
    \label{fig:torus_renormalization_2.2}

\end{figure}

Therefore, the new segment $S'$ induces an IET on 2 intervals $T'\colon S' \to S'$ as well, which exchanges the two intervals $S'_1$ and $S'_2$. As indicated in Figure \ref{fig:torus_renormalization_2.2}, we obtain these two intervals exactly as before by flowing backwards from the conical singularity until we hit $S'$. 

To obtain the first digit of the continued fraction expansion of $\lambda_1$, all we need to do is count the number of times an orbit starting from $S'_1$ visits the segment $S_1$ before returning to $S'$. Indeed, geometrically we can see that this is the same as counting how many times $S_2$ fits into the base of the torus. We will count the starting point, but not the first return. Using Figure \ref{fig:torus_renormalization_2.2}, we can see that the first digit in the continued fraction expansion of $\lambda_1$ is therefore 2. 

To be able to continue this algorithm, we need to \emph{renormalize}. We want to normalize the segment $S'$ and apply exactly the same procedure to the torus, which as a base has a normalized version of $S'$. What we need is the renormalized angle of the flow, or, equivalently, the length of the segment $S'_1$. At first glance, we might think that the length of $S_2'$ is the length of interest, but notice that the question we answered was how often the segment $S_2$ fits into the base of the torus and to continue with the algorithm we want to continue by considering the subinterval which is \enquote{left over} after removing the correct amounts of copies of $S_2$, which is the segment $S_1'$ and find how often it will fit into the segment of length $S_2$ which is of the same length as $S'$. This is then completely analogous to the Euclidean algorithm for integer division.

Nonetheless, we start by computing the length of the segment $S_2'$, which is given by

\begin{equation*}
    \left\{\ceil*{\frac{1}{\lambda_1}} \cdot \lambda_1 \right\} = \ceil*{\frac{1}{\lambda_1}} \cdot \lambda_1 - 1 = \lambda_1 \left(\ceil*{\frac{1}{\lambda_1}} - \frac{1}{\lambda_1} \right) = \lambda_1 \left(1 - \left\{\frac{1}{\lambda_1}\right\}\right),
\end{equation*}
where $\{\cdot\}$ denotes the fractional part and $\ceil{\cdot}$ is the ceiling function. Therefore, the length of the segment $S'_1$ is
\begin{equation*}
    \lambda_1 - |S'_2| = \lambda_1 \cdot \left\{\frac{1}{\lambda_1}\right\},
\end{equation*}
so that the length of the normalized version of $S'_1$ is given by $\lambda_2 = \left\{\frac{1}{\lambda_1}\right\}$. At this point it becomes clear that the algorithm described here indeed produces the digits of the continued fraction expansion of $\lambda_1$, since we recognize the definition of the Gauss map, i.e., Definition \ref{def:Gauss_map} and furthermore the number $\lfloor{\frac{1}{\lambda_1}}\rfloor$ is exactly the number of times we can fit the subinterval $S_2$ into the base of the torus, or equivalently, the number of visits to $S_1$ when starting in $S'_1$ before returning to $S'$. 

\begin{figure}[ht]
    \centering
    \begin{tikzpicture}[scale = 1.5]
    \draw[thick, blue, name path=blue_line] (0,0) -- (3,0) node[midway, below, black] {$b$};
    \draw[thick] (0,0) -- (0,3) node[midway, left] {$a$};
    \draw[thick] (3,0) -- (3,3) node[midway, right] {$a$};
    \draw[thick] (0,3) -- (3,3) node[midway, above] {$b$};
    
    \draw[dashed, purple, postaction={decorate}, decoration={
        markings,
        mark=at position 0.5 with {\arrow{stealth reversed}}
    }, name path=purple_line] (3,3) -- ++(-102.666:3.1);
    
    \fill [purple, name intersections={of=blue_line and purple_line, by=intersection}] (intersection) circle (0.7pt);
    
    \begin{scope}[on background layer]
        \draw [decorate, decoration={brace, amplitude=7pt}] (0,0) -- node[above=5pt] {$S'_2$} (intersection);
        \draw [decorate, decoration={brace, amplitude=7pt, mirror}] (intersection) -- node[below=5pt] {$S'_1$} (3,0);
    \end{scope}
\end{tikzpicture}
    \caption{The torus after renormalizing once.}
    \label{fig:torus_renormalization_2.3}
\end{figure}
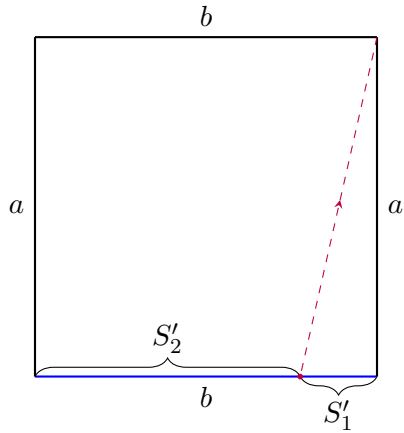

Figure \ref{fig:torus_renormalization_2.3} illustrates the situation after we have renormalized the torus. 

We proceed exactly as above, i.e., we follow the flow starting from the conical singularity in the top right corner backwards in time until we hit the segment $S''$, which is again given by the line whose left endpoint is the bottom left corner of the square representing the torus with length $\lambda_2$.

\begin{figure}[b]
    \centering
    \begin{tikzpicture}[scale = 1.5]
    \draw[thick, blue, name path=blue_line] (0,0) -- (3,0) node[midway, below, black] {$b$};
    \draw[thick] (0,0) -- (0,3) node[midway, left] {$a$};
    \draw[thick] (3,0) -- (3,3) node[midway, right] {$a$};
    \draw[thick] (0,3) -- (3,3) node[midway, above] {$b$};
    
    \draw[dashed, purple, postaction={decorate}, decoration={
        markings,
        mark=at position 0.5 with {\arrow{stealth reversed}}
    }, name path=purple_line] (3,3) -- ++(-102.666:3.08);

    \fill [purple, fill opacity = 1, name intersections={of=blue_line and purple_line, by=intersection_1}] (intersection_1) circle (0.5pt);
    
    \path (0,0) -- (3,0) coordinate[pos={(0.224744)}] (new_point);
    \fill [blue] (new_point) circle (0.5pt);
    
    \begin{scope}[on background layer]
        \draw [decorate, decoration={brace, amplitude=6pt, mirror}] (0,0) -- node[below=6pt] {$S''$} (new_point);
    \end{scope}
    
    \path let \p1=(intersection_1), \p2=(new_point) in 
        coordinate (start_point_1) at (\x1,3);
    \path[name path global=purple_line_1] (start_point_1) -- ++(-102.666:3.1);
    \path[name intersections={of=blue_line and purple_line_1, by=intersection_1}];
    \draw[dashed, purple, postaction = {decorate}, decoration = {
    markings, mark=at position 0.5 with {\arrow{stealth reversed}}
    }] (start_point_1) -- (intersection_1);

    \foreach \i in {1,...,2} {
    \path let \p1 = (intersection_1), \p2 = (new_point) in
        coordinate (start_point_\i) at (\x1, 3);
    \path[name path global= purple_line_\i] (start_point_\i) -- ++(-102.666:3.5);
    \path[name intersections={of=blue_line and purple_line_\i, by=intersection_1}];
    \draw[dashed, purple, postaction = {decorate}, decoration = {
    markings, mark=at position 0.5 with {\arrow{stealth reversed}}
    }] (start_point_\i) -- (intersection_1);
    }

    \fill [purple] (intersection_1) circle (0.5pt);
    \begin{scope}[on background layer]
        \draw[decorate, decoration={brace, amplitude=6pt}] (0,0) -- node[above=6pt] {$S''_1$} (intersection_1);
        \draw[decorate, decoration={brace, amplitude=6pt}] (intersection_1) -- node[above=6pt] {$\quad S''_2$} (new_point);
    \end{scope}

    \end{tikzpicture}
    \caption{We repeat the process with the renormalized torus.}
    \label{fig:torus_renormalization_2.4}
\end{figure}
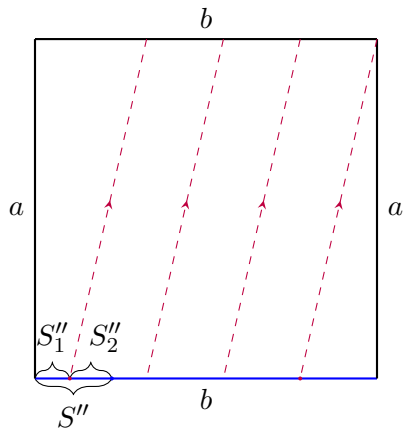

As before, we obtain the digit of the expansion by counting how many times a trajectory omitted from $S''_1$ visits $S'_2$ (note that the order of the intervals changed by the renormalization) before hitting $S''$ again. As can be seen in Figure \ref{fig:torus_renormalization_2.4}, the digit is 4. 

We can proceed with this algorithm indefinitely, since every irrational number has an infinite continued fraction expansion. In this particular case here, one can check that the new length $\lambda_3$ satisfies $\lambda_3 = \lambda_1$, so that we are again in the same situation as we were in the beginning. It follows that the continued fraction expansion of $\sqrt{6}-2$ is periodic, and more precisely that we have
\begin{equation*}
    \sqrt{6} = 2 + \cfrac{1}{2+\cfrac{1}{4+\cfrac{1}{2+\cfrac{1}{4+\cfrac{1}{\ddots}}}}} = [2; 2,4,2,4,2,4, \ldots].
\end{equation*}

It is easily seen that the Gauss map applied directly to $\lambda_1$ produces exactly the same expansion. 

\subsubsection{Surfaces of Higher Genus: Rauzy--Veech Induction}\label{sec:RV_induction}
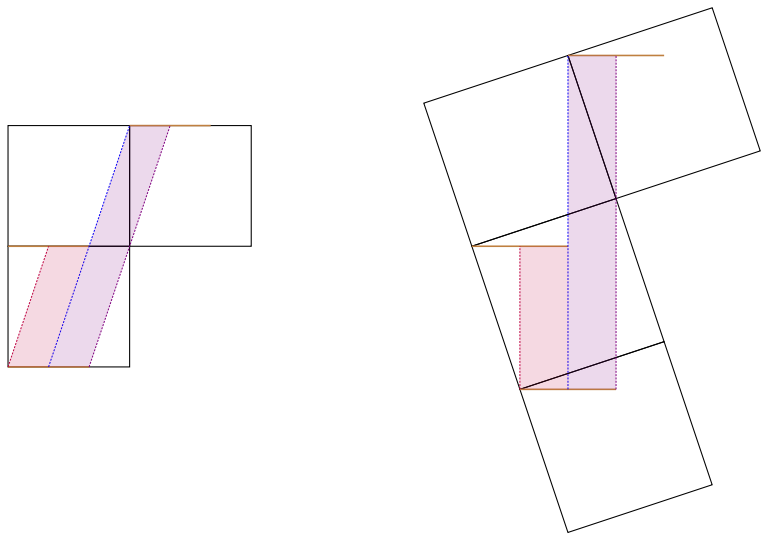
\begin{figure}[ht]
    \centering
    \begin{tikzpicture}[scale = 0.8, baseline={(0,0)}]

    \def\squaresize{2cm}
    
    \draw (0,0) rectangle (\squaresize,\squaresize);
    \draw (\squaresize,0) rectangle ++(\squaresize,\squaresize);
    \draw (0,-\squaresize) rectangle ++(\squaresize,\squaresize);
    
    \draw[brown, line width=0.6pt] (0,-\squaresize) -- ({2*\squaresize/3},-\squaresize);
    \draw[brown, line width = 0.6pt] (0,0) -- (2*\squaresize/3, 0);
    \draw[brown, line width = 0.6pt] (\squaresize, \squaresize) -- (5*\squaresize/3, \squaresize); 
    
    \draw[purple, dashed, line width=0.3pt, dash pattern=on 0.75pt off 0.5pt](0,-\squaresize) -- ++({\squaresize/3},\squaresize);
    \draw[blue, dashed, line width = 0.3pt, dash pattern = on 0.75pt off 0.5pt] (\squaresize/3, -\squaresize) -- ++ (2*\squaresize/3, 2*\squaresize);
    \draw[violet, dashed, line width = 0.3pt, dash pattern = 0n 0.75pt off 0.5pt]
    (2*\squaresize/3, -\squaresize) -- ++ (2*\squaresize/3, 2*\squaresize);
    
    \fill[purple, opacity=0.15] (0,-\squaresize) -- ++({\squaresize/3},0) -- ++(\squaresize/3, \squaresize) -- ++(-\squaresize/3, 0) -- cycle;
    
    \fill[violet, opacity = 0.15] (\squaresize/3, -\squaresize) -- ++ (\squaresize/3, 0) -- ++ (2*\squaresize/3, 2*\squaresize) --++ (-\squaresize/3, 0) -- cycle;
    
    \end{tikzpicture}
    \hspace{2cm}
    \begin{tikzpicture}[rotate around={18.43:(0,0)}, transform shape, scale =2, baseline={(0,0)}]
    
    \def\squaresize{1cm}
    
    \draw (0,0) rectangle (\squaresize,\squaresize);
    \draw (\squaresize,0) rectangle ++(\squaresize,\squaresize);
    \draw (0,-\squaresize) rectangle ++(\squaresize,\squaresize);
    \draw (0,-2*\squaresize) rectangle ++(\squaresize,\squaresize);
    
    \draw[brown, line width = 0.6pt] (0,0) -- (-90+71.57:0.6325*\squaresize) coordinate[pos=0.5] (mid_mid) coordinate[pos = 1] (right_mid);
    \draw[brown, line width = 0.6pt] (\squaresize,\squaresize) --++ (-90+71.57:0.6325*\squaresize)coordinate[pos=0.5] (mid_high);
    \draw[brown, line width = 0.6pt] (0,-\squaresize) --++ (-90+71.57:0.6325*\squaresize) coordinate[pos=0.5] (mid_low) coordinate[pos = 1] (right_low);
    
    \draw[purple, dashed, line width=0.3pt, dash pattern=on 0.75pt off 0.5pt](0,-\squaresize) -- (mid_mid);
    \draw[blue, dashed, line width = 0.3pt, dash pattern = on 0.75pt off 0.5pt] (mid_low) -- (\squaresize, \squaresize);
    \draw[violet, dashed, line width = 0.3pt, dash pattern = 0n 0.75pt off 0.5pt](right_low) -- (mid_high);
    
    \fill[purple, opacity=0.15] (0,-\squaresize) -- (mid_low) -- (right_mid) -- (mid_mid) -- cycle;
        
    \fill[violet, opacity = 0.15] (mid_low) -- (right_low) -- (mid_high) --(\squaresize, \squaresize) -- cycle;
    
    \end{tikzpicture}

    \caption{Equivalence of renormalization procedure for slanted and non-slanted torus. Clearly, the procedure explained in the previous section would not change when applied to the image on the right.}
    \label{fig:torus_straight_vs_slanted}
\end{figure}

Now we want to generalize the renormalization procedure seen in the previous section to IETs of any number of subintervals which will be induced by linear flows on surfaces of higher genus. It will be convenient to consider the \emph{vertical} linear flow as in Example \ref{ex:octagon_IET} and still choose the base $S$ to be horizontal, giving not just a transverse but a perpendicular section. The previous section is readily translated to this setting by considering a \enquote{tilted} torus. The description of the algorithm tremains the same, as can be seen in Figure \ref{fig:torus_straight_vs_slanted}. 

The generalization of the procedure from section \ref{sec:cf_torus} described now is known as \textit{Rauzy--Veech induction}. Let us stress that Rauzy--Veech induction does not give the algorithm from \ref{sec:cf_torus} immediately when restricted to a surface of genus 1, there is a subtle adjustment needed. We will explain precisely how Rauzy--Veech induction is connected to the procedure from section \ref{sec:cf_torus} (or, more generally, to the usual multiplicative algorithm of finding continued fractions using the Gauss map) in the next section. 

Recall the IET on 4 intervals obtained as a discretization of the vertical flow on a regular octagon from Example \ref{ex:octagon_IET}. Note that not only do we obtain a Poincaré map on the section, but we also obtain a new decomposition of the corresponding translation surface into polygons, namely into rectangles that are formed by following each subinterval along the flow until it returns. 

\begin{figure}[!hb]
    \centering
    \begin{tikzpicture}[rotate around={16:(0,0)},transform shape, scale =1.3, baseline={(0,0)}]
    \def\radius{2cm} 
    
    \foreach \angle/\coord in {0/{A}, 45/{B}, 90/{C}, 135/{D}, 180/{E}, 225/{F}, 270/{G}, 315/{H}}
    \coordinate (\coord) at (\angle:\radius);
      
    \draw[name path = octagon] (A) -- (B) -- (C) -- (D) -- (E) -- (F) -- (G) -- (H) -- cycle;
    
    \path[name path=section] (E) --+ (-16:4);
    \path[name path=vert_5] (B) --+(-106:4);

    \draw[line width = 1pt,purple,name intersections={of=section and vert_5,by={Int5}}] (E) -- (Int5);

    \path[name path=vert_1] (D) --+ (-106:4);
    \path[name path=vert_2] (F) --+ (74:4);
    \path[name path=vert_3] (C) --+ (-106:4);
    \path[name path=vert_4] (G) --+ (74:4);
    \path[name path=vert_6] (H) --+ (74:4);

    \path[name intersections={of=vert_1 and section, by=Int1}];
    \path[name intersections={of=vert_2 and section, by=Int2}];
    \path[name intersections={of=vert_3 and section, by=Int3}];
    \path[name intersections={of=vert_4 and section, by=Int4}];
    \path[name intersections={of=vert_6 and octagon,by=Int6}];

    \draw[black, dashed, line width=0.3pt, dash pattern=on 1.5pt off 1.5pt] (Int1) -- (D);
    \draw[black, dashed, line width=0.3pt, dash pattern=on 1.5pt off 1.5pt] (Int2) -- (F);    
    \draw[black, dashed, line width=0.3pt, dash pattern=on 1.5pt off 1.5pt] (Int3) -- (C);
    \draw[black, dashed, line width=0.3pt, dash pattern=on 1.5pt off 1.5pt] (Int4) -- (G);
    \draw[black, dashed, line width=0.3pt, dash pattern=on 1.5pt off 1.5pt] (Int5) -- (B);
    \draw[black, dashed, line width=0.3pt, dash pattern=on 1.5pt off 1.5pt] (Int6) -- (H);

    \path[name path=vert_7] (Int5) --+ (-106:2);
    \path[name intersections={of=vert_7 and octagon, by=Int7}];

    \draw[black, dashed, line width=0.3pt, dash pattern=on 1.5pt off 1.5pt] (Int7) -- (Int5);    
    
    \path[name path = helper_line] (Int6) --+ (-157.5:4);
    \path[name intersections={of=helper_line and octagon, name=Int8}];
    \path[name path = vert_8] (Int8-2) --+ (74:3);
    \path[name intersections={of=vert_8 and section, by=Int9}];

    \draw[black, dashed, line width=0.3pt, dash pattern=on 1.5pt off 1.5pt] (Int9) -- (Int8-2);    

    \path[name intersections={of=vert_5 and octagon, name=Int10}];
    \path[name path = helper_line_2] (Int10-2) --+ (-247.5:4);
    \path[name intersections={of=helper_line_2 and octagon, name = Int11}];
    \path[name path = helper_line3] (Int11-1) --+ (-106:4);
    \path[name intersections={of=helper_line3 and section, name=Int12}];

    \draw[black, dashed, line width=0.3pt, dash pattern=on 1.5pt off 1.5pt] (Int11-1) -- (Int12-1);        
    
    \fill[red, opacity=0.15] (E) -- (Int1) -- (D) -- cycle;
    \fill[red, opacity=0.15] (H) -- (A) -- (Int6) -- cycle;
    \fill[red, opacity=0.15] (F) -- (Int2) --(Int9) -- (Int8-2) -- cycle;

    \fill[violet, opacity=0.15] (Int1) -- (D) -- (Int11-1) -- (Int12-1) -- cycle;
    \fill[violet, opacity=0.15] (G) -- (Int10-2) -- (Int5) -- (Int4) -- cycle;

    \fill[blue, opacity = 0.15] (Int12-1) -- (Int3) -- (C) -- (Int11-1) -- cycle;
    \fill[blue, opacity = 0.15] (Int10-2) -- (H) -- (Int6) -- (B) -- cycle;
    \fill[blue, opacity = 0.15](E) -- (Int8-2) -- (Int9) -- cycle;

    \fill[teal, opacity = 0.15] (Int3) -- (Int5) -- (B) -- (C) -- cycle;
    \fill[teal, opacity = 0.15] (Int2) -- (Int4) -- (G) -- (F) -- cycle;
\end{tikzpicture}
    \hspace{2cm}
    \begin{tikzpicture}[rotate around={16:(0,0)},transform shape, scale =1.3, baseline={(0,0)}]
    \def\radius{2cm} 
    
    \foreach \angle/\coord in {0/{A}, 45/{B}, 90/{C}, 135/{D}, 180/{E}, 225/{F}, 270/{G}, 315/{H}}
    \coordinate (\coord) at (\angle:\radius);
      
    \draw[name path = octagon] (A) -- (B) -- (C) -- (D) -- (E) -- (F) -- (G) -- (H) -- cycle;
    
    \path[name path=section] (E) --+ (-16:4);
    \path[name path=vert_5] (B) --+(-106:4);

    \draw[line width = 1pt,purple,name intersections={of=section and vert_5,by={Int5}}] (E) -- (Int5);

    \path[name path=vert_1] (D) --+ (-106:4);
    \path[name path=vert_2] (F) --+ (74:4);
    \path[name path=vert_3] (C) --+ (-106:4);
    \path[name path=vert_4] (G) --+ (74:4);
    \path[name path=vert_6] (H) --+ (74:4);

    \path[name intersections={of=vert_1 and section, by=Int1}];
    \path[name intersections={of=vert_2 and section, by=Int2}];
    \path[name intersections={of=vert_3 and section, by=Int3}];
    \path[name intersections={of=vert_4 and section, by=Int4}];
    \path[name intersections={of=vert_6 and octagon,by=Int6}];

    \draw[black, dashed, line width=0.3pt, dash pattern=on 1.5pt off 1.5pt] (Int1) -- (D);
    \draw[black, dashed, line width=0.3pt, dash pattern=on 1.5pt off 1.5pt] (Int2) -- (F);    
    \draw[black, dashed, line width=0.3pt, dash pattern=on 1.5pt off 1.5pt] (Int3) -- (C);
    \draw[black, dashed, line width=0.3pt, dash pattern=on 1.5pt off 1.5pt] (Int4) -- (G);
    \draw[black, dashed, line width=0.3pt, dash pattern=on 1.5pt off 1.5pt] (Int5) -- (B);
    \draw[black, dashed, line width=0.3pt, dash pattern=on 1.5pt off 1.5pt] (Int6) -- (H);

    \path[name path=vert_7] (Int5) --+ (-106:2);
    \path[name intersections={of=vert_7 and octagon, by=Int7}];

    \draw[black, dashed, line width=0.3pt, dash pattern=on 1.5pt off 1.5pt] (Int7) -- (Int5);    
    
    \path[name path = helper_line] (Int6) --+ (-157.5:4);
    \path[name intersections={of=helper_line and octagon, name=Int8}];
    \path[name path = vert_8] (Int8-2) --+ (74:3);
    \path[name intersections={of=vert_8 and section, by=Int9}];

    \draw[black, dashed, line width=0.3pt, dash pattern=on 1.5pt off 1.5pt] (Int9) -- (Int8-2);    

    \path[name intersections={of=vert_5 and octagon, name=Int10}];
    \path[name path = helper_line_2] (Int10-2) --+ (-247.5:4);
    \path[name intersections={of=helper_line_2 and octagon, name = Int11}];
    \path[name path = helper_line3] (Int11-1) --+ (-106:4);
    \path[name intersections={of=helper_line3 and section, name=Int12}];

    \draw[black, dashed, line width=0.3pt, dash pattern=on 1.5pt off 1.5pt] (Int11-1) -- (Int12-1);        
    
    \fill[red, opacity=0.15] (E) -- (Int1) -- (D) -- cycle;

    \fill[violet, opacity=0.15] (Int1) -- (D) -- (Int11-1) -- (Int12-1) -- cycle;

    \fill[blue, opacity = 0.15] (Int12-1) -- (Int3) -- (C) -- (Int11-1) -- cycle;

    \fill[teal, opacity = 0.15] (Int3) -- (Int5) -- (B) -- (C) -- cycle;

    \thru[black, dashed, line width=0.3pt, dash pattern=on 1.5pt off 1.5pt]{(D)}{(point_1)}{(F)}{(Int2)}
    \thru[RoyalBlue, line width=0.5pt, decorate,decoration={zigzag,segment length=1mm,amplitude=0.2mm}]{(E)}{(point_2)}{(H)}{(Int6)}

    \draw (point_2) -- (D);

    \thru[black, dashed, line width=0.3pt, dash pattern=on 1.5pt off 1.5pt]{(point_2)}{(point_3)}{(Int8-2)}{(Int9)}

    \draw[line width = 1pt, purple] (point_3) -- (point_1);

    \thru[black, dashed, line width=0.3pt, dash pattern=on 1.5pt off 1.5pt]{(D)}{(point_4)}{(G)}{(Int4)}
    \thru[black, dashed, line width=0.3pt, dash pattern=on 1.5pt off 1.5pt]{(Int11-1)}{(point_5)}{(Int5)}{(Int10-2)}

    \draw[line width = 1pt, purple] (point_4) -- (point_5);

    \thru[black, dashed, line width=0.3pt, dash pattern=on 1.5pt off 1.5pt]{(Int11-1)}{(point_6)}{(Int10-2)}{(B)}
    \thru[RoyalBlue, line width=0.5pt, decorate,decoration={zigzag,segment length=1mm,amplitude=0.2mm}]{(C)}{(point_7)}{(H)}{(Int6)}

    \draw (point_6) -- (point_7);
    \thru[black, dashed, line width=0.3pt, dash pattern=on 1.5pt off 1.5pt]{(point_7)}{(point_8)}{(Int8-2)}{(Int9)}

    \draw[line width = 1pt, purple] (point_6) -- (point_8);
    
    \thru[black, opacity = 0, dashed, line width=0.3pt, dash pattern=on 1.5pt off 1.5pt]{(C)}{(point_9)}{(F)}{(Int2)}    
    \thru[black, dashed, line width=0.3pt, dash pattern=on 1.5pt off 1.5pt]{(B)}{(point_10)}{(G)}{(Int4)}

    \draw[line width = 1pt, purple] (point_9) -- (point_10);

    \fill[red, opacity=0.15] (E) -- (D) -- (point_2)-- cycle;
    \fill[red, opacity=0.15] (point_2) -- (point_3) -- (point_1) -- (D) -- cycle;
    \fill[violet, opacity=0.15] (D) -- (point_4) -- (point_5) -- (Int11-1) -- cycle;
    \fill[blue, opacity = 0.15] (Int11-1) -- (point_6) -- (point_8) -- (C) -- cycle;
    \fill[teal, opacity = 0.15] (C) -- (point_9) -- (point_10) -- (B) -- cycle;
\end{tikzpicture}
    \caption{Zippered rectangle decomposition of a regular octagon.}
    \label{fig:pretzel_zippered_rectangles}
\end{figure}
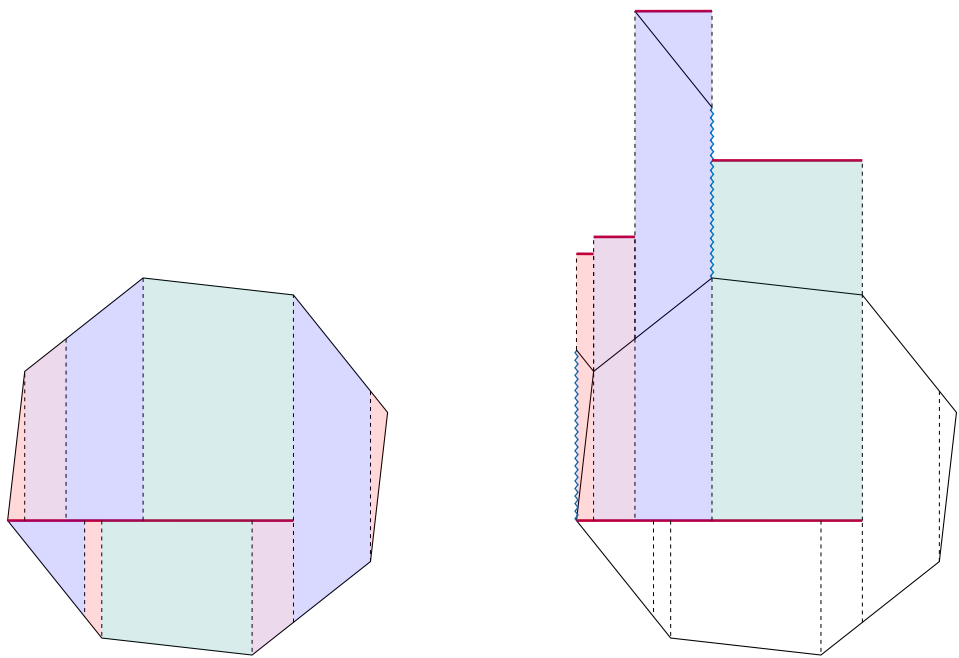

This decomposition into rectangles, seen on the right-hand side in Figure \ref{fig:torus_straight_vs_slanted}, is known as a decomposition into \emph{zippered rectangles}. This decomposition is useful to construct the Rauzy--Veech induction and was first introduced in \cite{veech1982gauss}. Initially, the algorithm was established purely in the language of IETs in \cite{rauzy1979echanges}. The number of rectangles in such a decomposition will be the same as the corresponding number of subintervals exchanged by the Poincaré map, thus it increases with the genus $\mathbf{g}$ of the surface by the formula seen in chapter \ref{ch:IETs_and_Renormalization}. It follows that a regular octagon, which represents a translation surface of genus $\mathbf{g} = 2$ with a single conical singularity of degree $d = 2$, coincides with a decomposition into 4 rectangles. This is illustrated in Figure \ref{fig:pretzel_zippered_rectangles}.

It is important to stress that the identifications of the vertical sides of the rectangles are as pictured on the left-hand side of Figure \ref{fig:pretzel_zippered_rectangles}. This is the reason we refer to the decomposition on the right-hand side as \emph{zippered} rectangles. For example, the part of the right side of the blue rectangle drawn as a blue zigzagged line is \emph{zippered} to the zigzagged segment on the left-hand side of the red rectangle to the far left.

The decomposition into zippered rectangles depicted in Figure \ref{fig:pretzel_zippered_rectangles} induces a permutation $\boldsymbol{\pi}$ by encoding in which way the top horizontal segments of the rectangles are identified to the bottom side of the rectangles, i.e., to the section $S$. It is evident that this permutation contains exactly the same information as the combinatorial datum from Lemma \ref{lem:IET_characterization}, so naturally we will use the same notation both when working with zippered rectangles and when working with IETs.

\example{
If we label the segments of the section $S$ in Figure \ref{fig:pretzel_zippered_rectangles} from left to right by $A,B,C,D$ then the combinatorial datum is given by
\begin{equation*}
    \boldsymbol{\pi} = 
    \begin{pmatrix}
        A & B & C & D \\
        C & A & D & B
    \end{pmatrix},
\end{equation*}
which of course coincides exactly with the datum from Example \ref{ex:octagon_IET}. The interpretation, when working with zippered rectangles, is that say the top side of the rectangle with base $C$ gets identified with the left-most part of the bottom segment and so on. 
}

Let now $(\boldsymbol{\pi}, \boldsymbol{\lambda})$ represent an IET on $d \geq 2$ intervals obtained as a Poincaré map of a vertical flow on a translation surface. Recall that the procedure from section \ref{sec:cf_torus} in essence consisted of two steps. First, we chose a smaller subinterval with which we obtained a new IET. Secondly, we renormalized in order to come back to the same starting point. This is also the overarching structure of the general algorithm. We will refer to the first step as the \emph{induction step} and to the second step as the \emph{renormalization step}.

We begin by comparing the length of the rightmost subintervals before and after the action of the IET. Formally, we denote by $\alpha(\varepsilon)$ the last symbol in the expression of $\pi_\varepsilon$ for $\varepsilon \in \{\operatorname{top}, \operatorname{bot}\}$, i.e., 
\begin{align*}
    \alpha(\mathrm{top}) &= \topp^{-1}(d) = \alpha_d^\mathrm{top},\\
    \alpha(\mathrm{bot}) &= \bott^{-1}(d) = \alpha_d^\mathrm{bot}.
\end{align*}
We will assume that the intervals $I_{\alpha(\mathrm{top})}$ and $I_{\alpha(\mathrm{bot})}$ do not have the same length. Otherwise, the induction step is not defined. Thus, we can differentiate between the following two cases.

\begin{enumerate}[1)]
    \item $\lambda_{\alpha(\mathrm{top})} > \lambda_{\alpha(\mathrm{bot})}$, i.e., the right-most interval before the interval exchange is \emph{longer} than the right-most interval after the exchange. In this case, we will say the IET determined by $(\boldsymbol{\pi}, \boldsymbol{\lambda})$ is of \emph{top type}.
    \item $\lambda_{\alpha(\mathrm{top})} < \lambda_{\alpha(\mathrm{bot})}$, i.e., the right-most interval before the interval exchange is \emph{shorter} than the right-most interval after the exchange. We will say that $(\boldsymbol{\pi}, \boldsymbol{\lambda})$ is of \emph{bot type}.
\end{enumerate}

Regardless of which case we are in, we will always call the larger of the two intervals $I_{\alpha(\mathrm{top})}$ and $I_{\alpha(\mathrm{bot})}$ the \emph{winner} and the shorter one the \emph{loser}. We can now define the subinterval $S' \subseteq S$ by
\begin{equation}\label{eq:subinterval_rv}
    S' = \begin{cases}
        I \setminus T(I_{\alpha(\mathrm{bot})}) \quad &\text{if } (\boldsymbol{\pi}, \boldsymbol{\lambda}) \text{ is of top type},\\
        I \setminus I_{\alpha(\mathrm{top})} \quad &\text{if } (\boldsymbol{\pi}, \boldsymbol{\lambda}) \text{ is of bot type}.
    \end{cases}
\end{equation}
Geometrically, it is very easy to see how exactly we construct $S'$. Let us first consider the case where $(\boldsymbol{\pi}, \boldsymbol{\lambda})$ is of top type, as is the case for the IET induced by the regular octagon from Figure \ref{fig:pretzel_zippered_rectangles}. To obtain $S'$, we move the right endpoint of the section $S$ to the left until either the backwards vertical flow emitted from the endpoint hits a singularity. Note that since we assume the IET is of top type, we will always hit a singularity with the backwards flow before we hit a singularity with the forwards flow emitted from the right endpoint. 

\begin{figure}[ht]
    \centering
    \begin{tikzpicture}[scale = 0.8]
    \draw (0,0) rectangle (0.8,1); 
    \draw (0.8,0) rectangle (2.3,1.3); 
    \draw (2.3,0) rectangle (4.8, 3); 
    \draw (4.8,0) rectangle (8.8,2); 

    \draw[thick, color = teal] (0,0) -- (2.5,0);
    \draw[thick, color = teal] (2.3, 3) -- (4.8, 3);

    \draw[thick, color = purple] (3.3, 0) -- (7.3, 0);
    \draw[thick, color = purple] (4.8,2) -- (8.8, 2);

    \draw[thick, color = green!10!blue!95!red!70] (2.5,0) -- (3.3, 0); 
    \draw[thick, color = green!10!blue!95!red!70] (0,1) -- (0.8,1);

    \draw[thick, color = violet] (7.3,0) -- (8.8, 0);
    \draw[thick, color = violet] (0.8,1.3) -- (2.3,1.3);

    \draw[dashed] (7.3,-0.5) -- (7.3, 2.5);
    \node[rotate =90] at (7.3,2.7) {\Large \Leftscissors};
    \fill[color = gray, opacity = 0.2] (7.3,0) rectangle (8.8, 2);
\end{tikzpicture}\qquad \qquad
\begin{tikzpicture}[scale = 0.8]
    \draw (0,0) rectangle (0.8,1); 
    \draw (0.8,0) rectangle (2.3,3.3); 
    \draw (2.3,0) rectangle (4.8, 3); 
    \draw (4.8,0) rectangle (7.3,2); 

    \draw[thick, color = teal] (0,0) -- (2.5,0);
    \draw[thick, color = teal] (2.3, 3) -- (4.8, 3);

    \draw[thick, color = purple] (3.3, 0) -- (5.8, 0);
    \draw[thick, color = purple] (4.8,2) -- (7.3, 2);

    \draw[thick, color = green!10!blue!95!red!70] (2.5,0) -- (3.3, 0); 
    \draw[thick, color = green!10!blue!95!red!70] (0,1) -- (0.8,1);

    \draw[thick, color = orange] (5.8,0) -- (7.3, 0);
    \draw[dashed, color = violet] (0.8,1.3) -- (2.3,1.3);
    \draw[thick, color = orange] (0.8, 3.3) -- (2.3, 3.3);

    \draw[white] (7.4,-0.5) -- (7.4, 2.5);

    \fill[color = gray, opacity = 0.2] (0.8,1.3) rectangle (2.3, 3.3);   
\end{tikzpicture}
\begin{tikzpicture}[remember picture, overlay]
    \draw[>=latex, ->, bend left=30, thick] (-7.3,2.5) to node[above] {} (-6.2,2.5);
\end{tikzpicture}
    \caption{An induction step for a zippered rectangle decomposition over an interval of top type.}
    \label{fig:RV_top_move_1}
\end{figure}
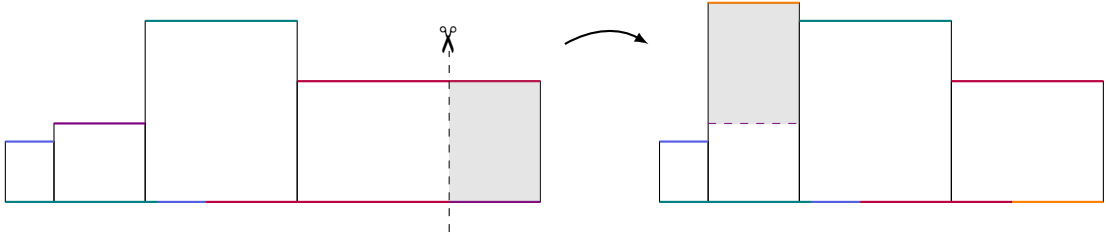

This geometric procedure coincides exactly with the definition of $S'$ in \eqref{eq:subinterval_rv}. Note also, that by construction $S'$ satisfies the MIC, so that the IET induced on $S'$ is again of the same number of intervals by Lemma \ref{lem:self_similarity}. Moreover, it is particularly easy to obtain a representation of the new IET on $S'$. What we do is we cut away a strip of the length of the loser from the winner and glue this strip to the top of the rectangle whose top side is identified with the bottom side of the strip we have cut away. This is pictured in Figure \ref{fig:RV_top_move_1}, which qualitatively is the same picture as in Figure \ref{fig:pretzel_zippered_rectangles} with some rescaling to obtain a clearer picture.

Let us explore how such this induction step acts on the combinatorial and length data $(\boldsymbol{\pi}, \boldsymbol{\lambda})$. The action on the length datum is very easy, since only the last interval $I_{\alpha(\mathrm{top})}$ is affected. So if the initial length data is given by
\begin{equation*}
    \boldsymbol{\lambda} = (\lambda_{\alpha_1^\mathrm{top}}, \lambda_{\alpha_2^\mathrm{top}}, \ldots, \lambda_{\alpha(\mathrm{top)}}),
\end{equation*}
then after the induction step the length data will be
\begin{equation}\label{eq:top_move_length}
    \boldsymbol{\lambda'} = (\lambda_{\alpha_1^\mathrm{top}}, \lambda_{\alpha_2^\mathrm{top}}, \ldots,  \lambda_{\alpha(\mathrm{top)}}-\lambda_{\alpha(\mathrm{bot})}).  
\end{equation}

Now we want to investigate the action on the combinatorial data $\boldsymbol{\pi}$. Note that since we do not fully remove the right-most interval $I_{\alpha(\mathrm{top})}$, the top row of the combinatorial data is unaffected, i.e., we have
\begin{equation*}
    \alpha_j^\mathrm{top} = {\alpha'}_j^\mathrm{top} \quad \text{for all } j \in [d],
\end{equation*}
where we use the notation ${\alpha'}_j^\varepsilon = {\pi'}_\varepsilon^{-1}(j)$ for the combinatorial datum $\boldsymbol{\pi'} = ({\pi'}_\mathrm{top}, \pi'_\mathrm{bot})$ of the IET induced on the subinterval $S'$. We do however need to adjust the bottom row of $\boldsymbol{\pi}$ as follows. There exists some unique $k$ such that $\alpha_k^\mathrm{bot} = \alpha(\mathrm{top})$, which represents the position the top side of the right-most rectangle gets glued. Because of the cutting operation, which cuts the interval $I_{\alpha(\mathrm{top})}$ into two pieces, we need to insert a letter from the alphabet $\mathcal{A}$ here. More precisely, the left part of $I_{\alpha(\mathrm{top})}$ still corresponds to the letter $\alpha_k^\mathrm{bot}$ and the part where the right part of $I_{\alpha(\mathrm{top})}$ gets glued appears next in the bottom row of $\boldsymbol{\pi}$ and is given by $\alpha(\mathrm{bot})$. Indeed, this follows from the construction of the induction step, since we cut away exactly a rectangle of the width of the loser $I_{\alpha(\mathrm{bot})}$, i.e., in the new IET on $S'$ the points that get identified with the right part of the winner $I_{\alpha(\mathrm{top})}$ originate in the interval $I_{\alpha(\mathrm{bot})}$. This description tells us exactly that we can formalize the action of a top type induction move on the combinatorial datum as
\begin{equation*}
    \boldsymbol{\pi'} = 
    \begin{pmatrix}
        \pi'_\mathrm{top} \\
        \pi'_\mathrm{bot}
    \end{pmatrix} = 
    \begin{pmatrix}
        \alpha_1^\mathrm{top} & \alpha_2^\mathrm{top} & \cdots & \alpha_{k-1}^\mathrm{top} & \alpha_k^\mathrm{top} & \alpha_{k+1}^\mathrm{top} & \cdots & \cdots & \alpha(\mathrm{top})\\
        \alpha_1^\mathrm{bot} & \alpha_2^\mathrm{bot} & \cdots & \alpha_{k-1}^\mathrm{bot} & \alpha(\mathrm{top}) & \alpha(\mathrm{bot}) & \alpha_{k+1}^\mathrm{bot} & \cdots & \alpha_{d-1}^\mathrm{bot}
    \end{pmatrix},
\end{equation*}
or, more concisely, as
\begin{align*}
    \alpha_j^\mathrm{top} &= {\alpha'}_j^\mathrm{top} \quad \quad\text{for all } j \in [d], \quad \text{and } \\
    {\alpha'}_j^\mathrm{bot} &= 
    \begin{cases}
        \alpha_j^\mathrm{bot} \quad &\text{if } j \leq k,\\
        \alpha(\mathrm{bot}) \quad &\text{if } j = k+1, \\
        \alpha_{j-1}^\mathrm{bot} &\text{if } j > k+1,
    \end{cases}
\end{align*}
where $k$ is implicitly defined by $\alpha_k^\mathrm{bot} = \alpha(\mathrm{top})$. 

    \begin{figure}[ht]
        \centering
        \begin{tikzpicture}[scale=1.3, transform shape]
  \node (pentagon1) [regular polygon, regular polygon sides=5, minimum size=3cm, draw, thick, rotate = -10] at (0,0) {};
  
  \node (pentagon2) [regular polygon, regular polygon sides=5, minimum size=3cm, draw, thick, rotate=26] at (2.14,-1.14) {};
  
  \clip (pentagon1.corner 1) -- (pentagon1.corner 2) -- (pentagon1.corner 3) -- (pentagon1.corner 4) --  (pentagon2.corner 2) -- (pentagon2.corner 3) -- (pentagon2.corner 4) -- (pentagon2.corner 5) -- (pentagon2.corner 1) -- cycle;

    \draw[dashed, name path = vert] (pentagon2.corner 5) -- ($(pentagon2.corner 5) + (0,-2)$);

    \path[name path = p1] (pentagon1.corner 3) -- ($(pentagon1.corner 3) + (5,0)$);

    \path[name intersections={of=p1 and vert,by=int}];

    \draw[thick, color = green!10!blue!95!red!70, name path = section] (pentagon1.corner 3) -- (int);

    \draw[dashed] (pentagon1.corner 3) -- ($(pentagon1.corner 3) + (0,3)$);

    \path[name path = p2] (pentagon1.corner 1) -- ($(pentagon1.corner 1) + (0,-3)$);
    \path[name intersections={of=p2 and section, by=int2}];
    \draw[dashed] (pentagon1.corner 1) -- (int2);

    \path[name path = p3] (pentagon1.corner 4) -- ($(pentagon1.corner 4) + (0,1)$);
    \path[name intersections={of=p3 and section, by=int3}];
    \draw[dashed] (pentagon1.corner 4) -- (int3);

    \path[name path = p4] (pentagon1.corner 5) -- ($(pentagon1.corner 5) + (0,-2)$);
    \path[name intersections={of=p4 and section, by=int4}];
    \draw[dashed] (pentagon1.corner 5) -- (int4);

    \path[name path = p5] (pentagon2.corner 3) -- ($(pentagon2.corner 3) + (0,2)$);
    \path[name intersections={of=p5 and section, by=int5}];
    \draw[dashed] (pentagon2.corner 3) -- (int5);

    \draw[dashed] ($(int5) + (0.2,0)$) -- ($(int5) + (0.2, -3)$);
    \draw[dashed] ($(int2) + (-0.2,0)$) -- ($(int2) + (-0.2, 3)$);

    \fill[red, opacity=0.15]  (pentagon1.corner 3) rectangle ($(pentagon1.corner 3) + (-2,3)$);
    \fill[red, opacity=0.15] ($(int2) + (-0.2,0)$) rectangle (pentagon1.corner 1);
    \fill[red, opacity=0.15] (int5) rectangle ($(int5) + (0.2,-2)$);
    \fill[red, opacity=0.15] (pentagon2.corner 5) rectangle ($(pentagon2.corner 5) + (1,-3)$);

    \fill[blue, opacity=0.15] (pentagon1.corner 3) rectangle ($(pentagon1.corner 1) - (0.2,0)$);
    \fill[blue, opacity = 0.15] (int) rectangle ($(pentagon2.corner 3) + (0.2,0)$);

    \fill[violet, opacity=0.15] (pentagon1.corner 1) rectangle (int4);
    \fill[violet, opacity=0.15] (pentagon2.corner 3) rectangle (int3);

    \fill[teal, opacity=0.15] (pentagon1.corner 3) rectangle (pentagon1.corner 4);
    \fill[teal, opacity=0.15] (pentagon1.corner 5) rectangle (int);
\end{tikzpicture}
        \caption{Zippered rectangle decomposition of a double pentagon used in Example \ref{ex:double_pentagon_IET_1}.}
        \label{fig:double_pentagon_IET_1}
    \end{figure}
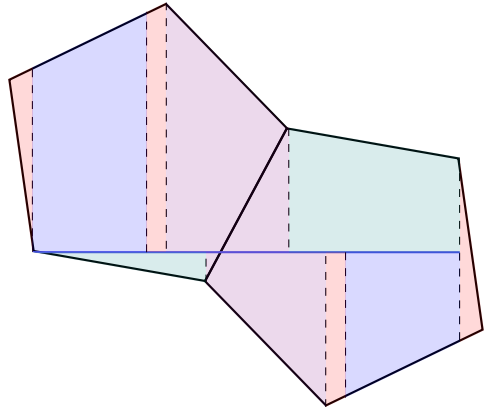
    
\begin{figure}[hb]
        \centering
        \begin{tikzpicture}[scale = 0.5]
    \draw (0,0) rectangle (2,4);
    \draw (2,0) rectangle (10/3,8);
    \draw (10/3,0) rectangle (17/3,6);
    \draw (17/3,0) rectangle (32/3, 3);

    \draw[thick, color = green!10!blue!95!red!70] (0,4) -- (2,4);
    \draw[thick, color = purple] (2,8) -- (10/3, 8);
    \draw[thick, color = violet] (10/3, 6) -- (17/3,6);
    \draw[thick, color = teal] (17/3,3) -- (32/3,3);

    \draw[thick, color = teal] (0,0) -- (5,0);
    \draw[thick, color = violet] (5,0) -- (22/3,0);
    \draw[thick, color = purple] (22/3, 0) -- (26/3,0);
    \draw[thick, color = green!10!blue!95!red!70] (26/3,0) -- (32/3,0);

    \draw[dashed] (26/3,-1) -- (26/3, 4);
    \node[rotate =90] at (26/3,4.2) {\Large \Leftscissors};
    \fill[color = gray, opacity = 0.2] (26/3,0) rectangle (32/3, 3);
\end{tikzpicture}
\qquad \qquad
\begin{tikzpicture}[scale = 0.5]
    \draw (2,0) rectangle (10/3,8);
    \draw (10/3,0) rectangle (17/3,6);
    \draw (17/3,0) rectangle (26/3, 3);

    \fill[color = gray, opacity = 0.2] (0,4) rectangle (2,7);
    \draw (0,0) rectangle (2,7);

    \draw[color = green!10!blue!95!red!70, dashed] (0,4) -- (2,4);
    \draw[thick, color = purple] (2,8) -- (10/3, 8);
    \draw[thick, color = violet] (10/3, 6) -- (17/3,6);
    \draw[thick, color = teal] (17/3,3) -- (26/3,3);
    \draw[thick, color = orange] (0,7) -- (2,7);

    \draw[thick, color = teal] (0,0) -- (3,0);
    \draw[thick, color = violet] (5,0) -- (22/3,0);
    \draw[thick, color = purple] (22/3, 0) -- (26/3,0);
    \draw[thick, color = orange] (3,0) -- (5,0);

    \draw[opacity = 0, dashed] (26/3,-1) -- (26/3,4);
\end{tikzpicture}
\begin{tikzpicture}[remember picture, overlay]
    \draw[>=latex, ->, bend left=30, thick] (-6,3) to node[above] {} (-5,3);
\end{tikzpicture}
        \caption{A top type induction move applied to the zippered rectangle decomposition from Figure \ref{fig:double_pentagon_IET_1}.}
        \label{fig:double_pentagon_IET_2}
    \end{figure}
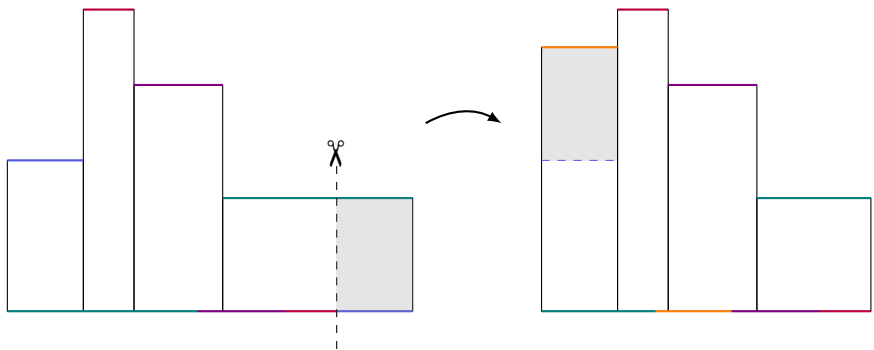

\begin{example}\label{ex:double_pentagon_IET_1}
    Consider the IET induced by the vertical flow on the double pentagon as in Figure \ref{fig:double_pentagon_IET_1}. One checks that the combinatorial datum is given by

    \begin{equation*}
        \boldsymbol{\pi} = 
        \begin{pmatrix}
            A & B & C & D \\
            D & C & B & A
        \end{pmatrix}.
    \end{equation*}
    After applying a top type induction move, the new combinatorial datum is given by
    \begin{equation*}
        \boldsymbol{\pi'} = 
        \begin{pmatrix}
            A & B & C & D \\
            D & A & C & B
        \end{pmatrix},
    \end{equation*}
    as is pictured in Figure \ref{fig:double_pentagon_IET_2}. In the language above, $k = 1$, i.e., the right most interval gets glued to the first position. so the first letter D is left unchanged while the remaining three letters are cyclically permuted. \qedhere
\end{example}

\begin{figure}[ht]
    \centering
    \begin{tikzpicture}[rotate around={30:(0,0)},transform shape, scale =1.3, baseline={(0,0)}]
    \def\radius{2cm} 
    
    \foreach \angle/\coord in {0/{A}, 45/{B}, 90/{C}, 135/{D}, 180/{E}, 225/{F}, 270/{G}, 315/{H}}
    \coordinate (\coord) at (\angle:\radius);
      
    \draw[name path = octagon] (A) -- (B) -- (C) -- (D) -- (E) -- (F) -- (G) -- (H) -- cycle;
    
    \path[name path=section] (E) --+ (-30:4);
    \path[name path=vert_5] (B) --+(-120:4);

    \draw[line width = 1pt,purple,name intersections={of=section and vert_5,by={I1}}] (E) -- (I1);

    \path[name path=vert_1] (C) --+ (-120:4);
    \path[name path=vert_2] (F) --+ (60:4);
    \path[name intersections={of=vert_1 and section, by=I2}];
    \path[name intersections={of=vert_2 and section, by=I3}];

    \path[name path=helpvert_1] (A) --+ (-120:3);
    \path[name intersections={of=helpvert_1 and octagon, name = HI1}];
    \path[name path=help_1] (HI1-2) --+ (112.5:4);
    \path[name intersections={of=help_1 and octagon, name = HI2}];
    \path[name path=vert_3] (HI2-1) --+ (-120:3);
    \path[name intersections={of=vert_3 and section, by=I4}];

    \path[name path=helpvert_2] (G) --+ (60:4);
    \path[name intersections={of=helpvert_2 and octagon, name=HI3}];
    \path[name path=help_2] (HI3-1) --+ (202.5:4);
    \path[name intersections={of=help_2 and octagon, name=HI4}];
    \path[name path=vert_4] (HI4-2) --+ (60:3);
    \path[name intersections={of=vert_4 and section, by=I5}];

    \path[name path = helpvert_3] (I1) --+ (-120:3);
    \path[name intersections={of=helpvert_3 and octagon, by=HI5}];
    \path[name path = help_3] (HI5) --+ (67.5:4);
    \path[name intersections={of=help_3 and octagon, name=HI6}];
    \path[name path = vert_5] (HI6-1) --+(-120:3);
    \path[name intersections={of=vert_5 and section, by=I6}];

    \path[name path = vert_6] (E) --+ (60:3);
    \path[name intersections={of=vert_6 and octagon, name = HI7}];

    \path[name path = help_4] (HI7-1) --+ (-67.5:4);
    \path[name intersections={of=help_4 and octagon, name=HI8}];
    \path[name path = helpvert_4] (HI8-2) --+ (60:4);
    \path[name intersections={of=helpvert_4 and octagon, name = HI9}];
    \path[name path = help_5] (HI9-1) --+ (202.5:4);
    \path[name intersections={of=help_5 and octagon, name = HI10}];
    \path[name path = vert_7] (HI10-2) --+ (60:4);
    \path[name intersections={of=vert_7 and section, by=I7}];
    \path[name intersections={of=vert_7 and octagon, name=HI11}];

    \path[name path = help_6] (HI11-1)--+(-67.5:4);
    \path[name intersections={of=help_6 and octagon, name=HI12}];
    \path[name path = helpvert_5] (HI12-2) --+ (60:4);
    \path[name intersections={of=helpvert_5 and octagon, name=HI13}];
    \path[name path = help_7] (HI13-1) --+ (202.5:4);

    \draw[black, dashed, line width=0.3pt, dash pattern=on 1.5pt off 1.5pt] (I1) -- (B);
    \draw[black, dashed, line width=0.3pt, dash pattern=on 1.5pt off 1.5pt] (I2) -- (C);    
    \draw[black, dashed, line width=0.3pt, dash pattern=on 1.5pt off 1.5pt] (I6) -- (HI6-1);
    \draw[black, dashed, line width=0.3pt, dash pattern=on 1.5pt off 1.5pt] (I1) -- (HI5);
    \draw[black, dashed, line width=0.3pt, dash pattern=on 1.5pt off 1.5pt] (G) -- (HI3-1);
    \draw[black, dashed, line width=0.3pt, dash pattern=on 1.5pt off 1.5pt] (E) -- (HI7-1);
    \draw[black, dashed, line width=0.3pt, dash pattern=on 1.5pt off 1.5pt] (F) -- (I3);
    \draw[black, dashed, line width=0.3pt, dash pattern=on 1.5pt off 1.5pt] (HI4-2) -- (I5);
    \draw[black, dashed, line width=0.3pt, dash pattern=on 1.5pt off 1.5pt] (HI8-2) -- (HI9-1);
    \draw[black, dashed, line width=0.3pt, dash pattern=on 1.5pt off 1.5pt] (HI10-2) -- (I7);
    \draw[black, dashed, line width=0.3pt, dash pattern=on 1.5pt off 1.5pt] (HI2-1) -- (I4);
    \draw[black, dashed, line width=0.3pt, dash pattern=on 1.5pt off 1.5pt] (HI1-2) -- (A);
    
    \fill[red, opacity=0.15] (E) -- (I4) -- (HI2-1) -- (HI7-1) -- cycle;
    \fill[red, opacity=0.15] (HI1-2) -- (HI8-2) -- (HI9-1) -- (A) -- cycle;
    \fill[red, opacity=0.15] (F) -- (I3) -- (I7) -- (HI10-2) -- cycle;

    \fill[violet, opacity=0.15] (HI2-1) -- (I4)  -- (I2) -- (C) -- cycle;
    \fill[violet, opacity=0.15] (G) -- (HI8-2) -- (HI9-1) -- (HI3-1) -- cycle;
    \fill[violet, opacity=0.15] (HI4-2) -- (I5) -- (I7) --(HI11-2) -- cycle;
    \fill[violet, opacity=0.15] (E) -- (D) -- (HI7-1) --cycle;
    \fill[violet, opacity=0.15] (H) -- (A) -- (HI1-2) -- cycle;
    
    \fill[blue, opacity = 0.15] (I2) -- (I6) -- (HI6-1) -- (C) -- cycle;
    \fill[blue, opacity = 0.15] (F) -- (HI5) -- (I1) -- (I3) -- cycle;

    \fill[teal, opacity = 0.15] (I6) -- (I1) -- (B) -- (HI6-1) -- cycle;
    \fill[teal, opacity = 0.15] (G) -- (HI3-1) -- (B) -- (HI5) -- cycle;
    \fill[teal, opacity = 0.15] (E) -- (HI4-2) -- (I5) -- cycle;
\end{tikzpicture}
    \hspace{2cm}
    \begin{tikzpicture}[rotate around={30:(0,0)},transform shape, scale =1.3, baseline={(0,0)}]
    \def\radius{2cm} 
    
    \foreach \angle/\coord in {0/{A}, 45/{B}, 90/{C}, 135/{D}, 180/{E}, 225/{F}, 270/{G}, 315/{H}}
    \coordinate (\coord) at (\angle:\radius);
      
    \draw[name path = octagon] (A) -- (B) -- (C) -- (D) -- (E) -- (F) -- (G) -- (H) -- cycle;
    
    \path[name path=section] (E) --+ (-30:4);
    \path[name path=vert_5] (B) --+(-120:4);

    \draw[line width = 1pt,purple,name intersections={of=section and vert_5,by={I1}}] (E) -- (I1);

    \path[name path=vert_1] (C) --+ (-120:4);
    \path[name path=vert_2] (F) --+ (60:4);
    \path[name intersections={of=vert_1 and section, by=I2}];
    \path[name intersections={of=vert_2 and section, by=I3}];

    \path[name path=helpvert_1] (A) --+ (-120:3);
    \path[name intersections={of=helpvert_1 and octagon, name = HI1}];
    \path[name path=help_1] (HI1-2) --+ (112.5:4);
    \path[name intersections={of=help_1 and octagon, name = HI2}];
    \path[name path=vert_3] (HI2-1) --+ (-120:3);
    \path[name intersections={of=vert_3 and section, by=I4}];

    \path[name path=helpvert_2] (G) --+ (60:4);
    \path[name intersections={of=helpvert_2 and octagon, name=HI3}];
    \path[name path=help_2] (HI3-1) --+ (202.5:4);
    \path[name intersections={of=help_2 and octagon, name=HI4}];
    \path[name path=vert_4] (HI4-2) --+ (60:3);
    \path[name intersections={of=vert_4 and section, by=I5}];

    \path[name path = helpvert_3] (I1) --+ (-120:3);
    \path[name intersections={of=helpvert_3 and octagon, by=HI5}];
    \path[name path = help_3] (HI5) --+ (67.5:4);
    \path[name intersections={of=help_3 and octagon, name=HI6}];
    \path[name path = vert_5] (HI6-1) --+(-120:3);
    \path[name intersections={of=vert_5 and section, by=I6}];

    \path[name path = vert_6] (E) --+ (60:3);
    \path[name intersections={of=vert_6 and octagon, name = HI7}];

    \path[name path = help_4] (HI7-1) --+ (-67.5:4);
    \path[name intersections={of=help_4 and octagon, name=HI8}];
    \path[name path = helpvert_4] (HI8-2) --+ (60:4);
    \path[name intersections={of=helpvert_4 and octagon, name = HI9}];
    \path[name path = help_5] (HI9-1) --+ (202.5:4);
    \path[name intersections={of=help_5 and octagon, name = HI10}];
    \path[name path = vert_7] (HI10-2) --+ (60:4);
    \path[name intersections={of=vert_7 and section, by=I7}];
    \path[name intersections={of=vert_7 and octagon, name=HI11}];

    \path[name path = help_6] (HI11-1)--+(-67.5:4);
    \path[name intersections={of=help_6 and octagon, name=HI12}];
    \path[name path = helpvert_5] (HI12-2) --+ (60:4);
    \path[name intersections={of=helpvert_5 and octagon, name=HI13}];
    \path[name path = help_7] (HI13-1) --+ (202.5:4);

    \draw[black, dashed, line width=0.3pt, dash pattern=on 1.5pt off 1.5pt] (I1) -- (B);
    \draw[black, dashed, line width=0.3pt, dash pattern=on 1.5pt off 1.5pt] (I2) -- (C);    
    \draw[black, dashed, line width=0.3pt, dash pattern=on 1.5pt off 1.5pt] (I6) -- (HI6-1);
    \draw[black, dashed, line width=0.3pt, dash pattern=on 1.5pt off 1.5pt] (I1) -- (HI5);
    \draw[black, dashed, line width=0.3pt, dash pattern=on 1.5pt off 1.5pt] (G) -- (HI3-1);
    \draw[black, dashed, line width=0.3pt, dash pattern=on 1.5pt off 1.5pt] (E) -- (HI7-1);
    \draw[black, dashed, line width=0.3pt, dash pattern=on 1.5pt off 1.5pt] (F) -- (I3);
    \draw[black, dashed, line width=0.3pt, dash pattern=on 1.5pt off 1.5pt] (HI4-2) -- (I5);
    \draw[black, dashed, line width=0.3pt, dash pattern=on 1.5pt off 1.5pt] (HI8-2) -- (HI9-1);
    \draw[black, dashed, line width=0.3pt, dash pattern=on 1.5pt off 1.5pt] (HI10-2) -- (I7);
    \draw[black, dashed, line width=0.3pt, dash pattern=on 1.5pt off 1.5pt] (HI2-1) -- (I4);
    \draw[black, dashed, line width=0.3pt, dash pattern=on 1.5pt off 1.5pt] (HI1-2) -- (A);
    
    \fill[red, opacity=0.15] (E) -- (I4) -- (HI2-1) -- (HI7-1) -- cycle;

    \fill[violet, opacity=0.15] (HI2-1) -- (I4)  -- (I2) -- (C) -- cycle;
    
    \fill[blue, opacity = 0.15] (I2) -- (I6) -- (HI6-1) -- (C) -- cycle;

    \fill[teal, opacity = 0.15] (I6) -- (I1) -- (B) -- (HI6-1) -- cycle;

    \thruu[black, dashed, line width=0.3pt, dash pattern=on 1.5pt off 1.5pt]{(HI7-1)}{(point_1)}{(HI8-2)}{(HI9-1)}
    \thruu[black, dashed, line width=0.3pt, dash pattern=on 1.5pt off 1.5pt]{(HI2-1)}{(point_2)}{(HI1-2)}{(A)}

    \draw (point_1) -- (point_2);

    \thruu[black, dashed, line width=0.3pt, dash pattern=on 1.5pt off 1.5pt]{(point_1)}{(point_3)}{(HI10-2)}{(I7)}
    \thruu[black, dashed, line width=0.3pt, dash pattern=on 1.5pt off 1.5pt]{(point_2)}{(point_4)}{(F)}{(I3)}

    \draw[line width = 1pt, purple] (point_3) -- (point_4);

    \fill[red, opacity=0.15] (HI7-1) -- (HI2-1) -- (point_4) -- (point_3) -- cycle;

    \thruu[black, dashed, line width=0.3pt, dash pattern=on 1.5pt off 1.5pt]{(HI2-1)}{(point_5)}{(HI1-2)}{(A)}
    \draw (point_5) -- (C);
    \thruu[black, dashed, line width=0.3pt, dash pattern=on 1.5pt off 1.5pt]{(C)}{(point_6)}{(E)}{(HI7-1)}
    \draw (point_6) -- (point_5);
    \thruu[black, dashed, line width=0.3pt, dash pattern=on 1.5pt off 1.5pt]{(point_5)}{(point_7)}{(G)}{(HI3-1)}
    \thruu[black, dashed, line width=0.3pt, dash pattern=on 1.5pt off 1.5pt]{(point_6)}{(point_8)}{(HI8-2)}{(HI9-1)}    
    \draw (point_7) -- (point_8);
    \thruu[black, dashed, line width=0.3pt, dash pattern=on 1.5pt off 1.5pt]{(point_7)}{(point_9)}{(HI4-2)}{(I5)} 
    \thruu[black, dashed, line width=0.3pt, dash pattern=on 1.5pt off 1.5pt]{(point_8)}{(point_10)}{(HI10-2)}{(I7)} 

    \draw[line width = 1pt, purple] (point_9) -- (point_10);
    \fill[violet, opacity=0.15] (HI2-1) -- (point_9)  -- (point_10) -- (C) -- cycle;

    \thruu[black, dashed, line width=0.3pt, dash pattern=on 1.5pt off 1.5pt]{(C)}{(point_11)}{(F)}{(I3)}
    \thruu[black, dashed, line width=0.3pt, dash pattern=on 1.5pt off 1.5pt]{(HI6-1)}{(point_12)}{(HI5)}{(I1)}

    \draw[line width = 1pt, purple] (point_11) -- (point_12);
    \fill[blue, opacity = 0.15] (C) -- (point_11) -- (point_12) -- (HI6-1) -- cycle;
    
    \thruu[black, dashed, line width=0.3pt, dash pattern=on 1.5pt off 1.5pt]{(HI6-1)}{(point_13)}{(HI5)}{(B)}    
    \thruu[black, dashed, line width=0.3pt, dash pattern=on 1.5pt off 1.5pt]{(B)}{(point_14)}{(G)}{(HI3-1)}
    \draw (point_13) -- (point_14);

    \thruu[black, dashed, line width=0.3pt, dash pattern=on 1.5pt off 1.5pt]{(point_14)}{(point_15)}{(HI4-2)}{(I5)}
    \draw[line width = 1pt, purple] (point_13) -- (point_15);

    \fill[teal, opacity = 0.15] (HI6-1) -- (point_13) -- (point_15) -- (B) -- cycle;

\end{tikzpicture}
    \caption{Zippered rectangle decomposition of a regular octagon. Here, the induced IET is of bot type.}
    \label{fig:pretzel_zippered_rectangle_2}
\end{figure}
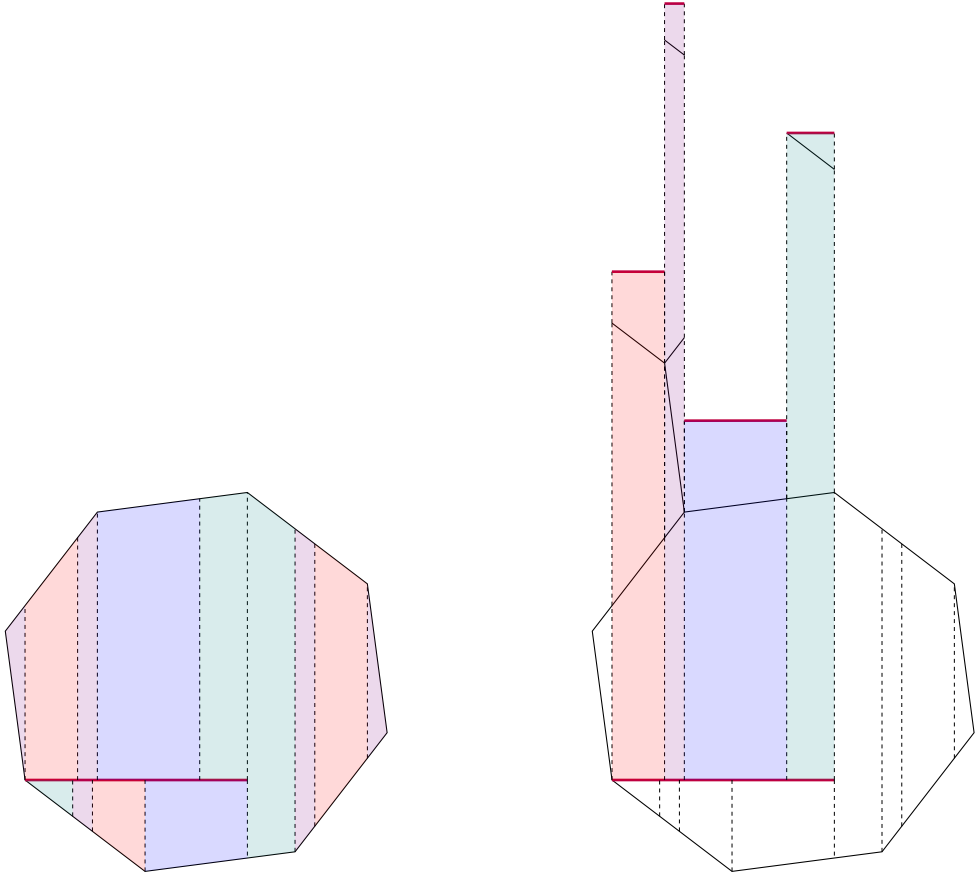

It remains to check how the data gets transformed by a bot type move. Figure \ref{fig:pretzel_zippered_rectangle_2} shows a different translation surface which is also given by a regular octagon, where the vertical flow induces an IET of bot type. 

Again, the action on the length datum $\boldsymbol{\lambda}$ can be expressed in a simple way. Analogously to the case of the IET of top type, it is now only the interval $I_{\alpha(\mathrm{bot})}$ that changes in length. Thus, the new length datum $\boldsymbol{\lambda'}$ is given by
\begin{equation}\label{eq:bot_move_length}
    \lambda'_\alpha = \lambda_\alpha \quad \text{if } \alpha \neq \alpha(\mathrm{bot}), \quad \text{and} \quad
    \lambda'_{\alpha(\mathrm{bot})} = \lambda_{\alpha(\mathrm{bot})} - \lambda_{\alpha(\mathrm{top})}. 
\end{equation}
For the combinatorial datum $\boldsymbol{\pi}$, we can argue as above to conclude that in the case of a bot type move the bottom row of the representation of $\boldsymbol{\pi}$ is unchanged, while we need to insert a symbol in the top row at the correct location. More precisely, the combinatorial datum $\boldsymbol{\pi'}$ after a bot type move is given by
\begin{equation*}
    \boldsymbol{\pi'} = 
    \begin{pmatrix}
        \pi'_\mathrm{top} \\
        \pi'_\mathrm{bot}
    \end{pmatrix} = 
    \begin{pmatrix}
        \alpha_1^\mathrm{top} & \alpha_2^\mathrm{top} & \cdots & \alpha_{k-1}^\mathrm{top} & \alpha(\mathrm{top}) & \alpha(\mathrm{bot}) & \alpha_{k+1}^\mathrm{top} & \cdots & \alpha_{d-1}^{\mathrm{top}}\\
        \alpha_1^\mathrm{bot} & \alpha_2^\mathrm{bot} & \cdots & \alpha_{k-1}^\mathrm{bot} & \alpha_k^\mathrm{bot} & \alpha_{k+1}^\mathrm{bot} & \cdots & \cdots & \alpha(\mathrm{bot})
    \end{pmatrix},
\end{equation*}
or, again written more concisely,

\begin{align*}
    \alpha_j^\mathrm{bot} &= {\alpha'}_j^\mathrm{bot} \quad \quad\text{for all } j \in [d], \quad \text{and } \\
    {\alpha'}_j^\mathrm{top} &= 
    \begin{cases}
        \alpha_j^\mathrm{top} \quad &\text{if } j \leq k,\\
        \alpha(\mathrm{top}) \quad &\text{if } j = k+1, \\
        \alpha_{j-1}^\mathrm{top} &\text{if } j > k+1,
    \end{cases}
\end{align*}
where $k$ is defined implicitly by $\alpha_k^\mathrm{top} = \alpha(\mathrm{bot})$. 

\begin{figure}[ht]
    \centering
    \begin{tikzpicture}[scale=1.3, transform shape]
  \node (pentagon1) [regular polygon, regular polygon sides=5, minimum size=3cm, draw, thick, rotate = 30] at (0,0) {};
  
  \node (pentagon2) [regular polygon, regular polygon sides=5, minimum size=3cm, draw, thick, rotate= 66] at (2.375,0.505) {};
  
  \clip (pentagon1.corner 1) -- (pentagon1.corner 2) -- (pentagon1.corner 3) -- (pentagon1.corner 4) --  (pentagon2.corner 2) -- (pentagon2.corner 3) -- (pentagon2.corner 4) -- (pentagon2.corner 5) -- (pentagon2.corner 1) -- cycle;

    \path[name path = vert] (pentagon2.corner 3) -- ($(pentagon2.corner 3) + (0,4)$);
    \path[name path = horz] (pentagon1.corner 2) -- ($(pentagon1.corner 2) + (6,0)$);
    \path[name intersections={of=horz and vert, by=int1}];

    \draw[thick, color = green!10!blue!95!red!70, name path = section] (pentagon1.corner 2) -- (int1);

    \path[name path = p1] (pentagon1.corner 3) -- ($(pentagon1.corner 3) + (0,3)$);
    \path[name intersections={of=p1 and section, by=int2}];
    \draw[dashed] (pentagon1.corner 3) -- (int2);
    
    \draw[dashed, name path = help] (pentagon1.corner 1) -- ($(pentagon1.corner 1) + (0,-3)$);
    \path[name intersections={of=help and section, by=inter}];
    \path[name path = p3] (pentagon1.corner 4) -- ($(pentagon1.corner 4) + (0,1)$);
    \path[name intersections={of=p3 and section, by=int3}];
    \draw[dashed] (pentagon1.corner 4) -- (int3);

    \path[name path = p4] (pentagon1.corner 5) -- ($(pentagon1.corner 5) + (0,-2)$);
    \path[name intersections={of=p4 and section, by=int4}];
    \draw[dashed] (pentagon1.corner 5) -- (int4);

    \draw[dashed] ($(int1) + (0,-2)$) -- ($(int1) + (0,4)$);

    \path[name path = p5] (pentagon2.corner 5) -- ($(pentagon2.corner 5) + (0,-3)$);
    \path[name intersections={of=p5 and section, by=int5}];
    \draw[dashed] (pentagon2.corner 5) -- (int5);

    \fill[red, opacity=0.15]  (pentagon1.corner 1) rectangle (int4);
    \fill[red, opacity=0.15]  (pentagon2.corner 3) rectangle (int3);

    \fill[blue, opacity=0.15] (pentagon1.corner 2) rectangle (pentagon1.corner 1);
    \fill[blue, opacity = 0.15] (pentagon2.corner 3) rectangle ($(pentagon2.corner 4) + (0,2)$);
    \fill[blue, opacity = 0.15] (pentagon1.corner 3) rectangle (inter);

    \fill[violet, opacity=0.15] (int4) rectangle (pentagon2.corner 5);
    \fill[violet, opacity=0.15] (pentagon1.corner 3) rectangle (int3);

    \fill[teal, opacity=0.15] (pentagon2.corner 5) rectangle (int1);
    \fill[teal, opacity=0.15] ($(pentagon1.corner 2) + (0,-1)$) rectangle (inter);
\end{tikzpicture}
    \caption{Zippered rectangle decomposition of a double pentagon. The induced IET is of bot type.}
    \label{fig:double_pentagon_IET_3}
\end{figure}
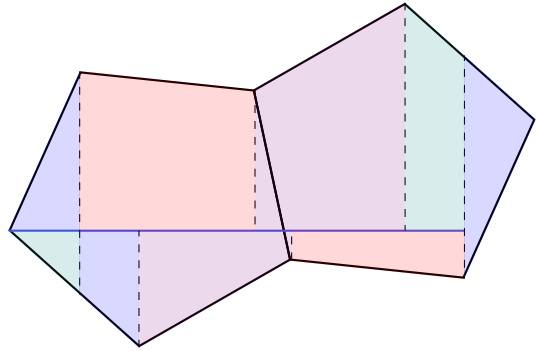

\begin{example}\label{ex:double_pentagon_IET_2}
    Figure \ref{fig:double_pentagon_IET_3} depicts a translation surface induced by a double pentagon, where the IET constructed from the vertical flow is of bot type. The combinatorial datum is given by
    \begin{equation*}
        \boldsymbol{\pi} = 
        \begin{pmatrix}
            A & B & C & D \\
            D & A & C & B
        \end{pmatrix},
    \end{equation*}
    and after applying a bot type induction move it is given by
    \begin{equation*}
        \boldsymbol{\pi'} = 
        \begin{pmatrix}
            A & B & D & C \\
            D & A & C & B
        \end{pmatrix}.
    \end{equation*}
    Using the language above, in this example we have $k = 2$, meaning that the first two letters of the bottom row remain the same while the other two letters are permuted.
\end{example}

\begin{figure}[hb]
    \centering
    \begin{tikzpicture}[scale = 0.4]
    \draw (0,0) rectangle (2,9);
    \draw (2,0) rectangle (6.5,4.5);
    \draw (6.5,0) rectangle (10.5,7.2);
    \draw (10.5,0) rectangle (12.2, 5.7);

    \draw[thick, color = green!10!blue!95!red!70] (0,9) -- (2,9);
    \draw[thick, color = purple] (2,4.5) -- (6.5, 4.5);
    \draw[thick, color = violet] (6.5, 7.2) -- (10.5,7.2);
    \draw[thick, color = teal] (10.5,5.7) -- (12.2,5.7);

    \draw[thick, color = teal] (0,0) -- (1.7,0);
    \draw[thick, color = violet] (3.7,0) -- (7.7,0);
    \draw[thick, color = purple] (7.7, 0) -- (12.2,0);
    \draw[thick, color = green!10!blue!95!red!70] (1.7,0) -- (3.7,0);

    \draw[dashed] (10.5,-1) -- (10.5, 8.2);
    \node[rotate =90] at (10.5,8.6) {\Large \Leftscissors};
    \fill[color = gray, opacity = 0.2] (10.5,0) rectangle (12.2, 5.7);
\end{tikzpicture}
\qquad \qquad
\begin{tikzpicture}[scale = 0.4]
    \draw (0,0) rectangle (2,9);
    \draw (2,0) rectangle (6.5,4.5);
    \draw (6.5,0) rectangle (10.5,7.2);

    \draw[thick, color = green!10!blue!95!red!70] (0,9) -- (2,9);
    \draw[thick, color = purple] (2,4.5) -- (4.8, 4.5);
    \draw[thick, color = violet] (6.5, 7.2) -- (10.5,7.2);
    \draw[thick, color = orange] (4.8,10.7) -- (6.5,10.7);

    \draw[line width = 2pt, color = white] (4.8, 4.5) -- (6.47, 4.5);
    \draw[color = purple, dashed] (4.8, 4.5) -- (6.5, 4.5);

    \draw (4.8, 4.5) -- (4.8,10.7);
    \draw (6.5,7.2) -- (6.5, 10.7);
    
    \fill[color = gray, opacity = 0.2] (4.8,4.5) rectangle (6.5, 10.7);

    \draw[thick, color = orange] (0,0) -- (1.7,0);
    \draw[thick, color = violet] (3.7,0) -- (7.7,0);
    \draw[thick, color = purple] (7.7, 0) -- (10.5,0);
    \draw[thick, color = green!10!blue!95!red!70] (1.7,0) -- (3.7,0);

    \draw[opacity = 0, dashed] (26/3,-1) -- (26/3,4);
\end{tikzpicture}
\begin{tikzpicture}[remember picture, overlay]
    \draw[>=latex, ->, bend left=30, thick] (-5.75,3.5) to node[above] {} (-4.75,3.5);
\end{tikzpicture}
    \caption{A bot type induction move applied to the zippered rectangle decomposition from Figure \ref{fig:double_pentagon_IET_3}.}
    \label{fig:double_pentagon_IET_4}
\end{figure}
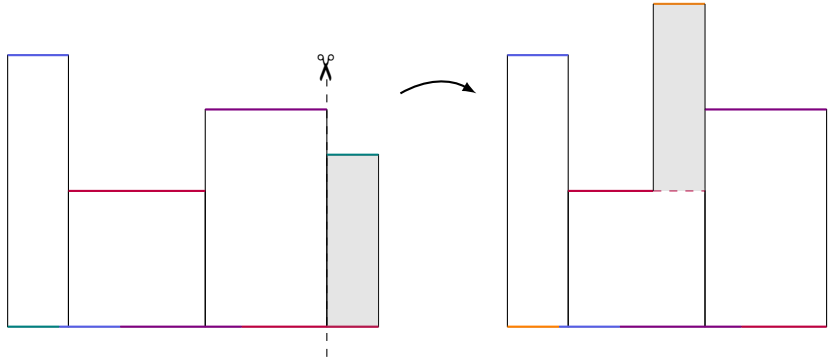

We remark that everything we have described so far can be entirely stated in terms of the IETs $T \colon S \to S$ and $T' \colon S' \to S'$, and we have used the geometric zippered rectangles construction only as a tool to derive the actions of the induction moves on the data $(\boldsymbol{\pi}, \boldsymbol{\lambda})$. This action on the level of IETs is the historical origin of this idea and was proposed by Rauzy in \cite{rauzy1979echanges}. To describe the action on the level of zippered rectangles, we would need more data, namely we need to record the heights of the rectangles and the position of the singularities. We will briefly mention this additional data in section \ref{sec:RV_Teichmuller_connection}, where we will explain how to link the Rauzy--Veech induction procedure to the Teichmüller geodesic flow on the moduli space of translation surfaces. For a detailed treatment we refer the reader to \cite{veech1982gauss}, where the zippered rectangles construction appeared first.

\subsubsection{Renormalization and the Space of IETs}\label{sec:renormalization_space_of_IET}

In section \ref{sec:RV_induction} we have explained in detail how given an IET $T \colon S \to S$ we can construct a new IET $T'\colon S' \to S'$, where the new IET is characterized by the combinatorial and length data $(\boldsymbol{\pi'}, \boldsymbol{\lambda'})$. The second step of the procedure, the \emph{renormalization step}, consists of normalizing the length of each subinterval such that the total length is of the same length as the initial segment $S$. Formally, the renormalization step is the map defined by
\begin{equation*}
    \lambda'_\alpha \mapsto \lambda'_\alpha \frac{|S|}{|S'|} \quad \text{for all }\alpha \in \mathcal{A}.
\end{equation*}

What we gain by this renormalization is that we are again in a comparable situation as in the beginning: We obtain an IET defined on some interval \emph{of the same length} as the initial interval $S$. Let us adopt the convention, that the initial interval was already normalized to have length $|S| = 1$, so that we may from now on assume that all the IETs obtained by the procedure above are defined on unit length intervals. For generic IETs, meaning that they are associated to a translation surface that satisfies the hypotheses of Keane's theorem, we can repeatedly perform the procedure infinitely often, which is why the name Rauzy--Veech \emph{induction} is indeed appropriate. We will say that such a generic IET satisfies \emph{Keane's condition}.

Let us go back to the case of an IET on 2 intervals, or equivalently, a linear flow on the torus satisfying the MIC for a brief moment. Note that for any IET on 2 intervals, there is no choice for the combinatorial datum $\boldsymbol{\pi}$, since the IET will always have to reorder the two intervals, i.e., the only possible combinatorial datum is given by
\begin{equation*}
    \boldsymbol{\pi} = \begin{pmatrix}
        A & B \\
        B & A
    \end{pmatrix}
\end{equation*}
up to a change of the alphabet $\mathcal{A}$. Since we assume that the initial segment $S$ is of unit length, we can thus parametrize an IET on 2 intervals completely by a single parameter $\lambda$ giving, for instance, the length of the right subinterval before the exchange. Applying Rauzy--Veech induction yields a new IET which is again completely characterized by the length of the right subinterval. 

In the case of a surface of higher genus, or an IET of $d \geq 3$ intervals, we obtain a similar parametrization. A priori the length datum $\boldsymbol{\lambda}$ consists of $d$ real numbers, but since we assume that the interval is normalized we lose one degree of freedom, i.e., having fixed the first $d-1$ lengths of the subintervals the length of the last subinterval is already determined. Therefore, the length datum $\boldsymbol{\lambda}$ is parametrized by the $(d-1)$-dimensional standard simplex
\begin{equation*}
    \Delta^{d-1} = \left\{(\lambda_\alpha)_{\alpha \in \mathcal{A}} \mid \sum_{\alpha \in \mathcal{A}} \lambda _\alpha = 1, \lambda_\alpha > 0\right\}.
\end{equation*}

\begin{figure}[ht]
    \centering
    \begin{tikzpicture}
    \draw[>=latex, ->, thick] (0,0) -- (3,0);
    \draw[>=latex, ->, thick] (0,0) -- (0,3);
    \draw[>=latex, ->, thick] (0,0) -- (-2,-1);

    \coordinate (p1) at (2,0);
    \coordinate (p2) at (0,2);
    \coordinate (p3) at (-4/3, -2/3);

    \draw[dashed, teal] (p1) -- (p2) -- (p3) -- cycle;
    \fill[color = teal, opacity = 0.2] (p1) -- (p2) -- (p3) -- cycle;
\end{tikzpicture}
    \caption{The two dimensional standard simplex.}
    \label{fig:standard_simplex_2d}
\end{figure}
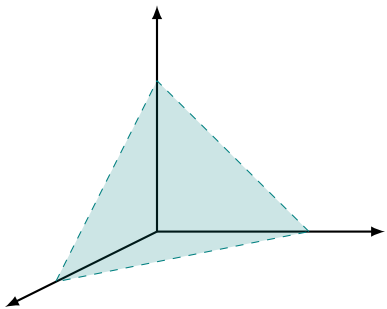

An illustration of the case when $d = 3$ can be seen in Figure \ref{fig:standard_simplex_2d}.

Contrary to the simple case of the torus, with $d \geq 3$ different pairs of bijections $\boldsymbol{\pi}$ may be possible. In particular, it may happen that the IET we are considering consists in fact of two (or more) independent IETs that act simultaneously. This is undesirable, since it is a general principle that if we can decompose our problem into smaller subproblems than we ought to be studying the subproblems instead. The following definition makes this idea more precise.

\begin{definition}[Reducible IET]
    We say that an IET $T\colon S \to S$ on $d \geq 2$ intervals is \emph{reducible}, if there exists $k \in \{1, \ldots, d-1\}$ such that the monodromy invariant $p = \pi_\mathrm{bot}\circ \pi_\mathrm{top}^{-1}$ satisfies
    \begin{equation*}
        p(\{1, \ldots, k\}) = \{1, \ldots, k\}.
    \end{equation*}
    In this case, the subinterval
    \begin{equation*}
        S' = \bigcup_{\pi_\mathrm{top}(\alpha) \leq k} I_\alpha = \bigcup_{\pi_\mathrm{bot}(\alpha) \leq k} I_\alpha
    \end{equation*}
    is invariant under the transformation $T$, and so is the complement of this subinterval. 
\end{definition}

\begin{example}[Reducible IET]
    Consider an IET with combinatorial datum
    \begin{equation*}
        \boldsymbol{\pi} = 
        \begin{pmatrix}
            A & B & C & D \\
            B & A & D & C
        \end{pmatrix}.
    \end{equation*}
    It is easy to see that the monodromy invariant satisfies
    \begin{equation*}
        p(\{1,2\}) = \{1,2\},
    \end{equation*}
    hence the IET is reducible. Indeed, irrespective of the length data $\lambda$, the IET describes two simultaneous IETs on 2 intervals
    \begin{equation*}
        \boldsymbol{\pi}_1 = \begin{pmatrix}
            A & B \\
            B & A
        \end{pmatrix}\quad \text{and} \quad
        \boldsymbol{\pi}_2 = \begin{pmatrix}
            C & D \\
            D & C
        \end{pmatrix}
    \end{equation*}
    that act independently. This is illustrated in Figure \ref{fig:reducible_IET_1}.

    \begin{figure}[ht]
        \centering
        \begin{tikzpicture}
    \draw[thick] (0,0) -- (5,0);
    \draw[thick] (0,-1.5) -- (5,-1.5);

    \foreach \p in {0,1,3,4.5,5}{
        \draw[thick] (\p,0) - ++(0,0.07);
        \draw[thick] (\p,0) - ++(0,-0.07);
    }

    \node[above] at (1/2,0) {\footnotesize $A$};
    \node[above] at (2,0)   {\footnotesize $B$};
    \node[above] at (3.75,0){\footnotesize $C$};
    \node[above] at (4.75,0){\footnotesize $D$};

    \foreach \p in {0,2,3,3.5,5}{
        \draw[thick] (\p,-1.5) - ++(0,0.07);
        \draw[thick] (\p,-1.5) - ++(0,-0.07);
    }

    \node[below] at (1,-1.5)   {\footnotesize $B$};
    \node[below] at (2.5,-1.5) {\footnotesize $A$};
    \node[below] at (3.25,-1.5){\footnotesize $D$};
    \node[below] at (4.25,-1.5){\footnotesize $C$};    

    \draw[>=latex, ->, thick] (2.5,-0.5) -- (2.5,-1);

    \draw[dashed, purple] (3,0.5) -- (3,-2);
    \node[rotate =90, purple] at (3,0.7) {\Large \Leftscissors};    
\end{tikzpicture}
\qquad\qquad
\begin{tikzpicture}
    \draw[thick] (0,0) -- (3,0); \draw[thick] (3.5,0) -- (5.5,0);
    \draw[thick] (0,-1.5) -- (3,-1.5); \draw[thick] (3.5,-1.5) -- (5.5,-1.5);

    \draw[color = white] (3,-1.9) -- (3,-2);

    \foreach \p in {0,1,3,3.5,5,5.5}{
        \draw[thick] (\p,0) - ++(0,0.07);
        \draw[thick] (\p,0) - ++(0,-0.07);
    }
    \foreach \p in {0,2,3,3.5,4,5.5}{
        \draw[thick] (\p,-1.5) - ++(0,0.07);
        \draw[thick] (\p,-1.5) - ++(0,-0.07);
    }

    \node[above] at (1/2,0) {\footnotesize $A$};
    \node[above] at (2,0)   {\footnotesize $B$};
    \node[above] at (4.25,0){\footnotesize $C$};
    \node[above] at (5.25,0){\footnotesize $D$};
    
    \node[below] at (1,-1.5)   {\footnotesize $B$};
    \node[below] at (2.5,-1.5) {\footnotesize $A$};
    \node[below] at (3.75,-1.5){\footnotesize $D$};
    \node[below] at (4.75,-1.5){\footnotesize $C$};    

    \draw[>=latex, ->, thick] (1.5,-0.5) -- (1.5,-1);
    \draw[>=latex, ->, thick] (4.5,-0.5) -- (4.5,-1);
\end{tikzpicture}
\begin{tikzpicture}[remember picture, overlay]
    \draw[>=latex, ->, bend left=30, thick, color =purple] (-6.9,1.2) to node[above] {} (-5.9,1.2);
\end{tikzpicture}
        \caption{An illustration of how a reducible IET splits into two smaller IETs. }
        \label{fig:reducible_IET_1}
    \end{figure}
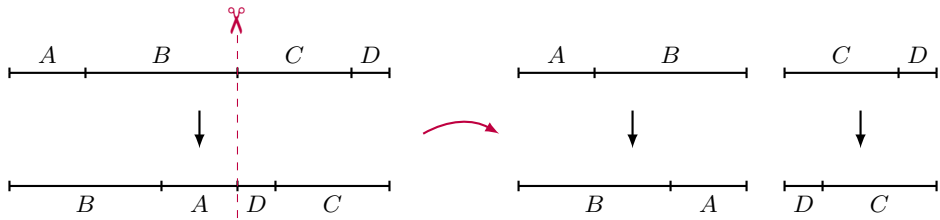    
\end{example}

\begin{remark}
    Abusing the language slightly, we will sometimes also call the combinatorial datum $\boldsymbol{\pi}$ irreducible, if this is true for the associated IET. 
\end{remark}

We see that if we have a \emph{reducible} IET, then we can decompose it into IETs of smaller size with simpler combinatorics. Studying these \enquote{building blocks} is enough, since every result we obtain for irreducible IETs will generalize readily to general IETs by putting together several irreducible IETs. For this reason, from now on we only consider irreducible IETs. 

It turns out that even when restricting to irreducible IETs, it is in general not possible to reach \emph{every} possible pair of bijections by applying Rauzy--Veech induction. Thus we cannot simply parametrize the combinatorial datum by all possible irreducible pairs of bijections. This motivates the following definition.

\begin{definition}[Rauzy class]\label{def:rauzy_class}
    A collection $\mathfrak{R}$ of irreducible pairs of permutations $\boldsymbol{\pi}^j = (\pi^j_\mathrm{top}, \pi^j_\mathrm{bot})$ is called a \emph{Rauzy class}, if for any $\boldsymbol{\pi}^i \neq \boldsymbol{\pi}^j \in \mathfrak{R}$, there exists a finite sequence of Rauzy--Veech induction moves starting from an IET with combinatorial datum $\boldsymbol{\pi}^i$ to an IET with combinatorial datum $\boldsymbol{\pi}^j$. 
\end{definition}

\begin{remark}
    Note that if $\boldsymbol{\pi'}$ is obtained from $\boldsymbol{\pi}$ by any of the two induction moves, then $\boldsymbol{\pi'}$ is irreducible if and only if $\boldsymbol{\pi}$ is irreducible. Moreover, it follows directly from the definition that we can represent any Rauzy class $\mathfrak{R}$ as a (connected) graph $\mathcal{G}$ with the property that every vertex has exactly two outgoing and two incoming edges.
\end{remark}

If we start with some irreducible combinatorial datum $\boldsymbol{\pi}$, we can associate to it a Rauzy class by applying top and bot moves until we cannot reach any new combinatorial datum. It is a fact, that such a class is independent of the starting permutation. 

Rauzy classes have been completely classified in \cite{kontsevich2003connected}. We will give two examples of Rauzy classes in the case $d = 4$, namely the Rauzy classes that correspond to the IETs from Example \ref{ex:double_pentagon_IET_1} and from Figure \ref{fig:irregular_octagon_IET_1}.

\begin{example}\label{ex:rauzy_class_1}
    Suppose we start with the IET given in Example \ref{ex:double_pentagon_IET_1}, which is characterized by the combinatorial datum
    \begin{equation*}
        \boldsymbol{\pi} = 
        \begin{pmatrix}
            A & B & C & D \\
            D & C & B & A
        \end{pmatrix}.
    \end{equation*}
    The following graph represents the associated Rauzy class $\mathfrak{R} = \mathfrak{R}(\boldsymbol{\pi})$.

    \begin{adjustbox}{max width=\textwidth}
    \begin{tikzcd}
	&& {\left(\genfrac{}{}{0pt}{}{A \, D \, B \, C}{D \, C \, A \, B}\right)}\arrow[loop above, distance = 1cm, in=60, out = 120, "\mathrm{bot}"] &&&& {\left(\genfrac{}{}{0pt}{}{A \, B \, D \, C}{D \, A \, C \, B}\right)}\arrow[loop above, distance = 1cm, in=60, out=120, "\mathrm{top}"] \\
	\\
	{\left(\genfrac{}{}{0pt}{}{A \, C \, D \, B}{D \, C \, B \, A}\right)}\arrow[loop left, distance = .9cm, in = 160, out = 200, "\mathrm{top}"] && {\left(\genfrac{}{}{0pt}{}{A \, D \, B \, C}{D \, C \, B \, A}\right)} && \textcolor{MidnightBlue}{\left(\genfrac{}{}{0pt}{}{A \, B \, C \, D}{D \, C \, B \, A}\right)} && {\left(\genfrac{}{}{0pt}{}{A \, B \, C \, D}{D \, A \, C \, B}\right)} && {\left(\genfrac{}{}{0pt}{}{A \, B \, C \, D}{D \, B \, A \, C}\right)}\arrow[loop right, distance = .9cm, in =-20, out = 20, "\mathrm{bot}"]
	\arrow["{\mathrm{bot}}"', curve={height=12pt}, from=3-5, to=3-3]
	\arrow["{\mathrm{bot}}"', curve={height=12pt}, from=3-3, to=3-1]
	\arrow["{\mathrm{bot}}"', curve={height=24pt}, from=3-1, to=3-5]
	\arrow["{\mathrm{top}}", curve={height=-12pt}, from=3-3, to=1-3]
	\arrow["{\mathrm{top}}", curve={height=-12pt}, from=1-3, to=3-3]
	\arrow["{\mathrm{top}}", curve={height=-12pt}, from=3-5, to=3-7]
	\arrow["{\mathrm{bot}}", curve={height=-12pt}, from=3-7, to=1-7]
	\arrow["{\mathrm{bot}}", curve={height=-12pt}, from=1-7, to=3-7]
	\arrow["{\mathrm{top}}", curve={height=-12pt}, from=3-7, to=3-9]
	\arrow["{\mathrm{top}}", curve={height=-24pt}, from=3-9, to=3-5]
\end{tikzcd}
\end{adjustbox}     

\begin{remark}
    We will refer to a graph associated to a Rauzy class as in Example \ref{ex:rauzy_class_1} as a \emph{Rauzy graph}.
\end{remark}
    
\end{example}

\begin{example}\label{ex:rauzy_class_2}

    Consider the IET induced by the vertical flow in the translation surface obtained by identifying parallel sides in the irregular octagon pictured in Figure \ref{fig:irregular_octagon_IET_1}.

    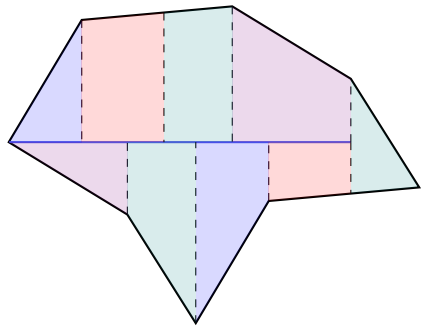
\begin{figure}[ht]
        \centering
        \begin{tikzpicture}[scale = 0.6, transform shape]
    \coordinate (p1) at (0,0);
    \coordinate (p2) at (1.6, 2.7);
    \coordinate (p3) at (4.9, 3);
    \coordinate (p4) at (7.5, 1.4);
    \coordinate (p5) at (9, -1);
    \coordinate (p6) at (5.7, -1.3);
    \coordinate (p7) at (4.1, -4);
    \coordinate (p8) at (2.6, -1.6);

    \draw[thick, name path = octagon] (p1) -- (p2) -- (p3) -- (p4) -- (p5) -- (p6) -- (p7) -- (p8) -- cycle;
    \clip (p1) -- (p2) -- (p3) -- (p4) -- (p5) -- (p6) -- (p7) -- (p8) -- cycle;

    \path[name path = horz] (p1) -- ($(p1) + (10,0)$);
    \path[name path = vert] (p4) -- ($(p4) + (0,-4)$);
    \path[name intersections={of=horz and vert, by=int1}];
    \draw[thick, color = green!10!blue!95!red!70, name path = section] (p1) -- (int1);

    \path[name path = l1] (p2) -- ($(p2) + (0,-3)$);
    \path[name intersections={of=l1 and section, by=int2}];
    \draw[dashed] (p2)-- (int2);

    \path[name path = l2] (p8) -- ($(p8) + (0,3)$);
    \path[name intersections={of=l2 and section, by=int3}];
    \draw[dashed] (p8)-- (int3);

    \path[name path = l3] (p7) -- ($(p7) + (0,5)$);
    \path[name intersections={of=l3 and section, by=int4}];
    \draw[dashed] (p7)-- (int4);

    \path[name path = l4] (p3) -- ($(p3) + (0,-3)$);
    \path[name intersections={of=l4 and section, by=int5}];
    \draw[dashed] (p3)-- (int5);

    \path[name path = l5] (p6) -- ($(p6) + (0,2)$);
    \path[name intersections={of=l5 and section, by=int6}];
    \draw[dashed] (p6)-- (int6);

    \draw[dashed] (int1) -- ++(90:2);
    \draw[dashed, name path = vert2] (int1) -- ++(-90:2);
    \draw[name intersections={of=vert2 and octagon, by=inter}];

    \coordinate (special) at ($(inter) -(p6) + (p2)$);
    \path[name path =l6] (special) --($(special) + (0,-3)$);
    \path[name intersections={of=l6 and section, by=int7}];
    \draw[dashed] (special) -- (int7);


    \fill[color = blue, opacity = 0.15] (p1) rectangle (p2);
    \fill[color = blue, opacity = 0.15] (p7) rectangle (int6);

    \fill[color = red, opacity = 0.15] (special) rectangle (int2);
    \fill[color = red, opacity = 0.15] (p6) rectangle (int1);

    \fill[color = violet, opacity = 0.15] (p1) rectangle (p8);
    \fill[color = violet, opacity = 0.15] (p3) rectangle (int1);

    \fill[color = teal, opacity = 0.15] (p7) rectangle (int3);
    \fill[color = teal, opacity = 0.15] (inter) -- (p4) -- (p5) -- cycle;
    \fill[color = teal, opacity = 0.15] (int5) rectangle ($(special) + (0,1)$);
    
\end{tikzpicture}
        \caption{An IET induced by the vertical flow on a translation surface obtained by identifying parallel sides of an irregular octagon.}
        \label{fig:irregular_octagon_IET_1}
    \end{figure}

    We can see that the associated combinatorial datum is given by
    \begin{equation*}
        \boldsymbol{\pi} = 
        \begin{pmatrix}
            A & B & C & D \\
            D & C & A & B
        \end{pmatrix}.
    \end{equation*}

    The Rauzy class $\mathfrak{R}(\boldsymbol{\pi})$ obtained from this combinatorial datum is represented in the graph in Figure \ref{fig:rauzy_graph}. 

    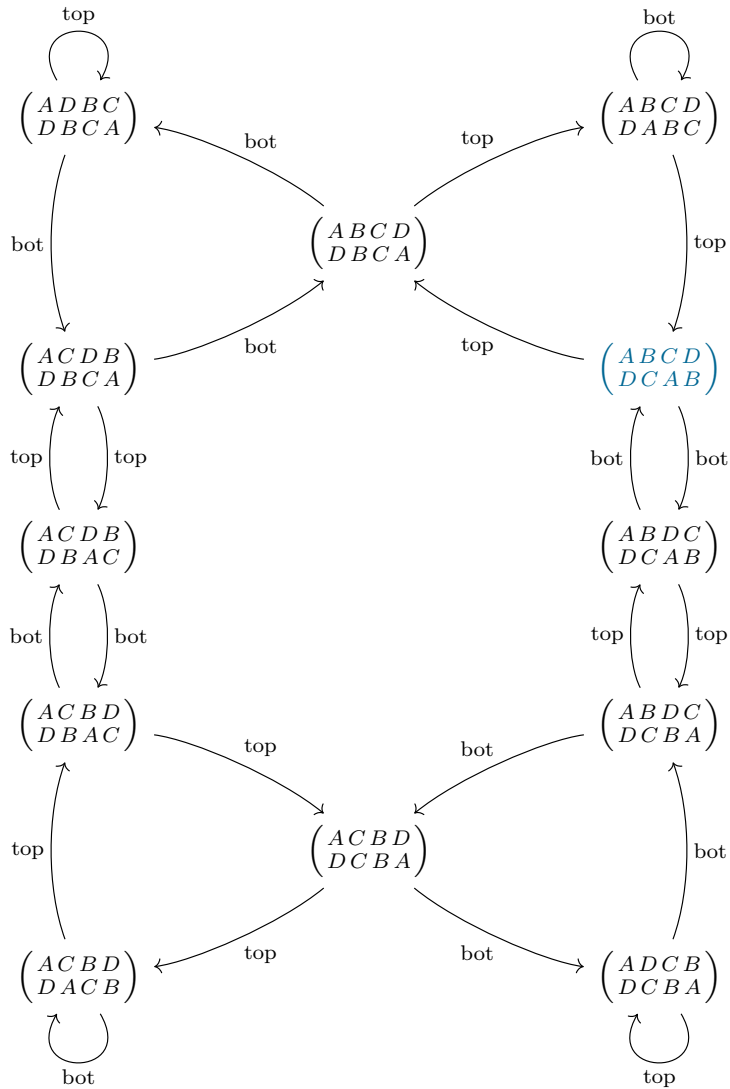
\begin{figure}[ht]
        \centering
\[\begin{tikzcd}
	{\left(\genfrac{}{}{0pt}{}{A \, D \, B \, C}{D \, B \, C \, A}\right)}\arrow[loop above, distance = 1cm, in=60, out = 120, "\mathrm{top}"] &&&& {\left(\genfrac{}{}{0pt}{}{A \, B \, C \, D}{D \, A \, B \, C}\right)}\arrow[loop above, distance = 1cm, in=60, out = 120, "\mathrm{bot}"] \\
	&& {\left(\genfrac{}{}{0pt}{}{A \, B \, C \, D}{D \, B \, C \, A}\right)} \\
	{\left(\genfrac{}{}{0pt}{}{A \, C \, D \, B}{D \, B \, C \, A}\right)} &&&& {\textcolor{MidnightBlue}{\left(\genfrac{}{}{0pt}{}{A \, B \, C \, D}{D \, C \, A \, B}\right)}} \\
	\\
	{\left(\genfrac{}{}{0pt}{}{A \, C \, D \, B}{D \, B \, A \, C}\right)} &&&& {\left(\genfrac{}{}{0pt}{}{A \, B \, D \, C}{D \, C \, A \, B}\right)} \\
	\\
	{\left(\genfrac{}{}{0pt}{}{A \, C \, B \, D}{D \, B \, A \, C}\right)} &&&& {\left(\genfrac{}{}{0pt}{}{A \, B \, D \, C}{D \, C \, B \, A}\right)} \\
	&& {\left(\genfrac{}{}{0pt}{}{A \, C \, B \, D}{D \, C \, B \, A}\right)} \\
	{\left(\genfrac{}{}{0pt}{}{A \, C \, B \, D}{D \, A \, C \, B}\right)}\arrow[loop below, distance = 1cm, in=240, out = 300, "\mathrm{bot}"] &&&& {\left(\genfrac{}{}{0pt}{}{A \, D \, C \, B}{D \, C \, B \, A}\right)}\arrow[loop below, distance = 1cm, in=240, out = 300, "\mathrm{top}"]
	\arrow["{\mathrm{top}}", curve={height=-12pt}, from=3-5, to=2-3]
	\arrow["{\mathrm{top}}", curve={height=-12pt}, from=2-3, to=1-5]
	\arrow["{\mathrm{top}}", curve={height=-12pt}, from=1-5, to=3-5]
	\arrow["{\mathrm{bot}}", curve={height=-12pt}, from=3-5, to=5-5]
	\arrow["{\mathrm{bot}}", curve={height=-12pt}, from=5-5, to=3-5]
	\arrow["{\mathrm{top}}", curve={height=-12pt}, from=5-5, to=7-5]
	\arrow["{\mathrm{top}}", curve={height=-12pt}, from=7-5, to=5-5]
	\arrow["{\mathrm{bot}}"', curve={height=12pt}, from=9-5, to=7-5]
	\arrow["{\mathrm{bot}}"', curve={height=12pt}, from=7-5, to=8-3]
	\arrow["{\mathrm{bot}}"', curve={height=12pt}, from=8-3, to=9-5]
	\arrow["{\mathrm{top}}", curve={height=-12pt}, from=8-3, to=9-1]
	\arrow["{\mathrm{top}}", curve={height=-12pt}, from=9-1, to=7-1]
	\arrow["{\mathrm{top}}", curve={height=-12pt}, from=7-1, to=8-3]
	\arrow["{\mathrm{bot}}", curve={height=-12pt}, from=7-1, to=5-1]
	\arrow["{\mathrm{bot}}", curve={height=-12pt}, from=5-1, to=7-1]
	\arrow["{\mathrm{top}}", curve={height=-12pt}, from=5-1, to=3-1]
	\arrow["{\mathrm{top}}", curve={height=-12pt}, from=3-1, to=5-1]
	\arrow["{\mathrm{bot}}"', curve={height=12pt}, from=3-1, to=2-3]
	\arrow["{\mathrm{bot}}"', curve={height=12pt}, from=2-3, to=1-1]
	\arrow["{\mathrm{bot}}"', curve={height=12pt}, from=1-1, to=3-1]
\end{tikzcd}\] 
        \caption{The Rauzy graph from Example \ref{ex:rauzy_class_2}.}
        \label{fig:rauzy_graph}
    \end{figure}   
\end{example}

Note how both graphs above are symmetric with respect to the vertical axis. This is no coincidence, this symmetry corresponds exactly to interchanging the roles of $\pi_\mathrm{top}$ and $\pi_\mathrm{bot}$, e.g., if we consider the two combinatorial data
\begin{align*}
    \boldsymbol{\pi}_1 &=
    \begin{pmatrix}
        A & D & B & C \\
        D & B & C & A
    \end{pmatrix}\quad \text{and} \\
    \boldsymbol{\pi}_2 &=
    \begin{pmatrix}
        A & B & C & D \\
        D & A & B & C
    \end{pmatrix},
\end{align*}
which can be found on the top left and top right of the Rauzy graph in Example \ref{ex:rauzy_class_2}, then it is easy to verify that $\boldsymbol{\pi}_1$ and $\Tilde{\boldsymbol{\pi}}_2$ have the same monodromy invariant, where 
\begin{equation*}
    \Tilde{\boldsymbol{\pi}}_2 = 
    \begin{pmatrix}
        D & A & B & C \\
        A & B & C & D
    \end{pmatrix}
\end{equation*}
is obtained by exchanging the top and bottom row of $\boldsymbol{\pi}_2$. Moreover, the Rauzy graph from Example \ref{ex:rauzy_class_2} has an additional symmetry: Pairs of combinatorial data that are symmetric with respect to the center have the same monodromy invariant. Since such transformations correspond essentially to the same IET, it is sometimes convenient to identify these pairs.

\begin{definition}[Reduced Rauzy class]
    Given a Rauzy class $\mathfrak{R}$, the associated \emph{reduced Rauzy class} $\Tilde{\mathfrak{R}}$ is obtained by identifying combinatorial data that have the same monodromy invariant. Furthermore, we will call the graph associated to $\Tilde{\mathfrak{R}}$ the \emph{reduced Rauzy graph}. 
\end{definition}

\begin{example}[Reduced Rauzy class]\label{ex:reduced_rauzy_class_1}
    The reduced Rauzy graph associated to the class from Example \ref{ex:rauzy_class_2} is pictured in Figure \ref{fig:reduced_rauzy_graph}.

    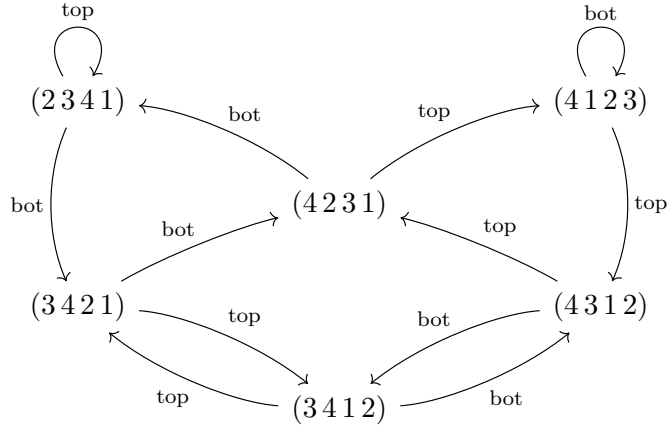
\begin{figure}[ht]
        \centering
\[\begin{tikzcd}
	{(2\, 3 \, 4\, 1)}\arrow[loop above, distance = 1cm, in=60, out = 120, "\mathrm{top}"] &&&& {(4 \, 1 \, 2 \, 3)}\arrow[loop above, distance = 1cm, in=60, out = 120, "\mathrm{bot}"] \\
	&& {(4 \, 2 \, 3 \, 1)} \\
	{(3\, 4 \, 2\, 1)} &&&& {(4 \, 3\, 1\, 2)} \\
	&& {(3\, 4\, 1\, 2)}
	\arrow["{\mathrm{top}}", curve={height=-12pt}, from=2-3, to=1-5]
	\arrow["{\mathrm{top}}", curve={height=-12pt}, from=1-5, to=3-5]
	\arrow["{\mathrm{top}}"', curve={height=6pt}, from=3-5, to=2-3]
	\arrow["{\mathrm{bot}}"', curve={height=12pt}, from=4-3, to=3-5]
	\arrow["{\mathrm{bot}}"', curve={height=12pt}, from=3-5, to=4-3]
	\arrow["{\mathrm{top}}", curve={height=-12pt}, from=3-1, to=4-3]
	\arrow["{\mathrm{top}}", curve={height=-12pt}, from=4-3, to=3-1]
	\arrow["{\mathrm{bot}}"', curve={height=12pt}, from=2-3, to=1-1]
	\arrow["{\mathrm{bot}}"', curve={height=12pt}, from=1-1, to=3-1]
	\arrow["{\mathrm{bot}}", curve={height=-6pt}, from=3-1, to=2-3]
\end{tikzcd}\]
        \caption{The reduced Rauzy graph from Example \ref{ex:reduced_rauzy_class_1}.}
        \label{fig:reduced_rauzy_graph}
    \end{figure}
\end{example}

\begin{remark}
    For IETs on 4 intervals, one can check that the Rauzy classes from Examples \ref{ex:rauzy_class_1} and \ref{ex:rauzy_class_2} are exhaustive: There are 24 possibilities for the monodromy invariant of which only 13 are irreducible. Seven of these monodromy invariants appear in the Rauzy graph in Example \ref{ex:rauzy_class_1}, the other 6 appear in the Rauzy graph from Example \ref{ex:rauzy_class_2} or likewise in the reduced Rauzy graph from Example \ref{ex:reduced_rauzy_class_1}. 
\end{remark}

We have now defined all the objects needed to state very precisely on which space the Rauzy--Veech induction, that is the induction step followed by the renormalization step, acts. Namely, the procedure described in the last two sections \ref{sec:RV_induction} and \ref{sec:renormalization_space_of_IET} applied to an IET on $d$ subintervals of a unit length interval describes a map

\begin{equation*}
    \mathcal{T} \colon  \mathfrak{R} \times \Delta^{d-1} \to \mathfrak{R} \times \Delta^{d-1}
\end{equation*}
from the \emph{space of IETs} (associated to the Rauzy class $\mathfrak{R}$) to itself. In other words, we have found an explicit form of a discrete dynamical system which we now can study, which may not only reveal the behavior of IETs but may also tell us about the linear flow on translation surfaces of higher genus. 

\begin{remark}
    To be more precise, since we want to be able to iterate $\mathcal{T}$ indefinitely we define it on a subset of $\Delta^{d-1}$ consisting of length data $\boldsymbol{\lambda}$ satisfying Keane's condition. This subset has full Lebesgue measure. 
\end{remark}

Having this precise definition at hand, a first natural question concerns the invariant measures of this dynamical system. The following result, proven independently by Masur \cite{masur1982interval} and Veech \cite{veech1982gauss}, gives a very satisfying characterization.

\begin{theorem}\label{thm:rauzy_veech_invariant_measure}
    For each Rauzy class $\mathfrak{R}$, the Rauzy--Veech induction map $\mathcal{T}\colon \mathfrak{R} \times \Delta^{d-1} \to \mathfrak{R} \times \Delta^{d-1}$ admits an invariant measure $\nu$ which is absolutely continuous with respect to $\dd \boldsymbol{\pi} \otimes \operatorname{Leb}$, where $\dd \boldsymbol{\pi}$ denotes the counting measure on $\mathfrak{R}$ and $\operatorname{Leb}$ denotes the Lebesgue measure on $\Delta^{d-1}$. 

    Moreover, this measure $\nu$ is unique up to multiplication by a scalar as well as ergodic. Its density with respect to the Lebesgue measure is given by the restriction to $\mathfrak{R} \times \Delta^{d-1}$ of a function defined on $\mathfrak{R} \times \R^d_+ = \mathfrak{R} \times \{\lambda_1, \ldots, \lambda_d \mid \lambda_j > 0\}$, which is a homogeneous rational function of degree $-d$ and bounded away from zero.
\end{theorem}

We will see an explicit example in section \ref{sec:connecting_cf_to_rv_induction}. 

Note that the measure $\nu$ from Theorem \ref{thm:rauzy_veech_invariant_measure} is infinite, i.e.,
\begin{equation*}
    \nu(\mathfrak{R} \times \Delta^{d-1}) = +\infty.
\end{equation*}
This issue is significant, since the greater part of the tools from ergodic theory are only applicable when working with finite measure spaces. But there is a general method to overcome issues of this type, which is sometimes called the method of \emph{jump transformations}. Morally, the mechanism that allows us to obtain a \emph{finite} measure is based on the following heuristic. Say that the infinitude of the reference measure, in this case the measure $\nu$ from Theorem \ref{thm:rauzy_veech_invariant_measure}, is caused because the system is very slow in some part of its domain: A trajectory tends to spend a lot of time in some small part of the space, which means that any invariant measure needs to assign a lot of measure to this part, possible even an infinite amount of measure. The idea is to \emph{accelerate} the transformations. Instead of keeping track of \emph{all} the iterates of the map, we are only interested in a carefully chosen subset of iterates and we \enquote{jump over} all the rest. Remarkably, choosing this subset in a good way ensures that properties enjoyed by the accelerated sytem, such as ergodicity, translate to the same properteis of the original map. We refer the reader to \cite{dajani2021first} for more details on the theory of jump transformations. 

Following this main idea, let us now define an acceleration of the map $\mathcal{T} \colon \mathfrak{R} \times \Delta^{d-1} \to \mathfrak{R} \times \Delta^{d-1}$ which from now on we will refer to as \emph{slow} Rauzy--Veech induction. 

To obtain an accelerated version of the map $\mathcal{T}$, instead of doing the top (or bot) moves one by one, we apply \emph{at once} as many top (or bot) moves as possible. Formally, let $\mathfrak{R}$ be any Rauzy class and let $(\boldsymbol{\pi}, \boldsymbol{\lambda})$ be such that $\boldsymbol{\lambda}$ satisfies Keane's condition and $\boldsymbol{\pi} \in \mathfrak{R}$. Let us write $\varepsilon \in \{\mathrm{top}, \mathrm{bot}\}$ for the type of $(\boldsymbol{\pi}, \boldsymbol{\lambda})$ and for $j \geq 1$ let $\varepsilon^{(j)}$ be the type of $\mathcal{T}^j(\boldsymbol{\pi}, \boldsymbol{\lambda}) = (\boldsymbol{\pi}^{(j)}, \boldsymbol{\lambda}^{(j)})$.

Then we define
\begin{equation}\label{eq:fast_RV}
    n = n(\boldsymbol{\pi}, \boldsymbol{\lambda}) = \inf\{j \in \N \mid \varepsilon^{(j)} \neq \varepsilon\}.
\end{equation}

\begin{lemma}\label{lem:moves_change}
    The number $n$ defined above is finite. 
\end{lemma}
\begin{proof}
    Let us write $\Tilde{\mathcal{T}}$ for the operator given by just the induction step of $\mathcal{T}$, i.e., without renormalizing the interval, so the map $\Tilde{\mathcal{T}}$ is an endomorphism on $\mathfrak{R} \times \R_{>0}^d$. Towards a contradiction, suppose that $n = +\infty$, meaning that $\varepsilon^{(j)}$ is constantly equal to $\varepsilon$. In that case, the winner of each iterate is constant as well, say it is constantly equal to $\alpha_\mathrm{win} \in \mathcal{A}$. Every application of $\Tilde{\mathcal{T}}$ reduces the length of $I_{\alpha_\mathrm{win}}$ by the length of the loser. Upon application of $\Tilde{\mathcal{T}}$, the only length that is changed is the length of $I_{\alpha_\mathrm{win}}$ itself, meaning that the length of the possible loser intervals is bounded away from zero. But then, the length of $I_{\alpha_{\mathrm{win}}}$ is eventually negative, which is absurd. 
\end{proof}

\begin{remark}
    Note that Lemma \ref{lem:moves_change} immediately implies that when iterating the slow Rauzy--Veech induction $\mathcal{T}$, both top and bot moves appear infinitely many times.
\end{remark}

Let us also remark that one can show something even stronger than Lemma $\ref{lem:moves_change}$. A proof of the following proposition can be found in \cite{viana2006ergodic}. 

\begin{proposition}
    If we denote by $\alpha_\mathrm{win}^n$ the winner and by $\alpha_\mathrm{lose}^n$ the loser of $\mathcal{T}^n = (\boldsymbol{\pi}^n, \boldsymbol{\lambda}^n)$, then both sequences $(\alpha_\mathrm{win}^n)_{n \in \N}$ and $(\alpha_\mathrm{lose}^n)_{n \in \N}$ take every value $\alpha \in \mathcal{A}$ infinitely many times. 
\end{proposition}

We can now define the accelerated version of the (slow) Rauzy--Veech induction which was first introduced in \cite{zorich1996finite}. We will refer to the accelerated version as \emph{fast} Rauzy--Veech induction, but want to remark that in the literature it is also referred to as the \emph{Zorich transformation} or \emph{Zorich induction}.

\begin{definition}[Fast Rauzy--Veech induction]
    The \emph{fast Rauzy--Veech induction} is defined to be the map
    \begin{align*}
        \mathcal{G} \colon \mathfrak{R} \times \Delta^{d-1} &\to \mathfrak{R} \times \Delta^{d-1}, \\
        (\boldsymbol{\pi} , \boldsymbol{\lambda}) & \mapsto \mathcal{T}^n(\boldsymbol{\pi}, \boldsymbol{\lambda}) = (\boldsymbol{\pi}^n, \boldsymbol{\lambda}^n),
    \end{align*}
    where $n$ is given by \eqref{eq:fast_RV}.
\end{definition}

\begin{figure}[ht]
    \centering
    \begin{tikzpicture}
    \draw[thick] (0,0) -- (14,0);
    \draw[thick] (0,-1) -- (14,-1);

    \foreach \p in {0, 42/5, 62/5, 14}{
        \draw[thick] (\p,0) - ++(0,0.07);
        \draw[thick] (\p,0) - ++(0,-0.07);    
    }
    \foreach \p in {0, 8/5, 10, 14}{
        \draw[thick] (\p,-1) - ++(0,0.07);
        \draw[thick] (\p,-1) - ++(0,-0.07);    
    }
    \node[above] at (21/5,0) {\footnotesize $A$};
    \node[above] at (52/5,0)   {\footnotesize $B$};
    \node[above] at (66/5,0){\footnotesize $C$};

    \node[below] at (4/5,-1) {\footnotesize $C$};
    \node[below] at (29/5,-1)   {\footnotesize $A$};
    \node[below] at (12,-1){\footnotesize $B$};

    \draw[>=latex, ->, thick] (5,-1.5) -- (5,-2.5);
    \node[right, xshift = 5] at (5, -2) {$\mathcal{T}$ \footnotesize(bot)};

    \draw[thick] (0,-3) -- (62/5,-3);
    \draw[thick] (0,-4) -- (62/5,-4);

    \foreach \p in {0, 42/5, 54/5 ,62/5}{
        \draw[thick] (\p,-3) - ++(0,0.07);
        \draw[thick] (\p,-3) - ++(0,-0.07);    
    }
    \foreach \p in {0, 8/5, 10, 62/5}{
        \draw[thick] (\p,-4) - ++(0,0.07);
        \draw[thick] (\p,-4) - ++(0,-0.07);    
    }
    \node[above] at (21/5,-3) {\footnotesize $A$};
    \node[above] at (48/5,-3)   {\footnotesize $B$};
    \node[above] at (58/5,-3){\footnotesize $C$};

    \node[below] at (4/5,-4) {\footnotesize $C$};
    \node[below] at (29/5,-4)   {\footnotesize $A$};
    \node[below] at (56/5,-4){\footnotesize $B$};

    \draw[>=latex, ->, thick] (5,-4.5) -- (5,-5.5);
    \node[right, xshift = 5] at (5, -5) {$\mathcal{T}$ \footnotesize(bot)};

    \draw[thick] (0,-6) -- (54/5,-6);
    \draw[thick] (0,-7) -- (54/5,-7);

    \foreach \p in {0, 42/5, 46/5, 54/5}{
        \draw[thick] (\p,-6) - ++(0,0.07);
        \draw[thick] (\p,-6) - ++(0,-0.07);    
    }
    \foreach \p in {0, 8/5, 10, 54/5}{
        \draw[thick] (\p,-7) - ++(0,0.07);
        \draw[thick] (\p,-7) - ++(0,-0.07);    
    }
    \node[above] at (21/5,-6) {\footnotesize $A$};
    \node[above] at (44/5,-6)   {\footnotesize $B$};
    \node[above] at (50/5,-6){\footnotesize $C$};

    \node[below] at (4/5,-7) {\footnotesize $C$};
    \node[below] at (29/5,-7)   {\footnotesize $A$};
    \node[below] at (52/5,-7){\footnotesize $B$};

    \draw[>=latex, ->, thick] (5,-7.5) -- (5,-8.5);
    \node[right, xshift = 5] at (5, -8) {$\mathcal{T}$ \footnotesize\textcolor{purple}{(top)}};

    \draw[thick] (0,-9) -- (10,-9);
    \draw[thick] (0,-10) -- (10,-10);

    \foreach \p in {0, 42/5, 46/5, 10}{
        \draw[thick] (\p,-9) - ++(0,0.07);
        \draw[thick] (\p,-9) - ++(0,-0.07);    
    }
    \foreach \p in {0, 4/5, 8/5, 10}{
        \draw[thick] (\p,-10) - ++(0,0.07);
        \draw[thick] (\p,-10) - ++(0,-0.07);    
    }
    \node[above] at (21/5,-9) {\footnotesize $A$};
    \node[above] at (44/5,-9)   {\footnotesize $B$};
    \node[above] at (48/5,-9){\footnotesize $C$};

    \node[below] at (2/5,-10) {\footnotesize $C$};
    \node[below] at (29/5,-10)   {\footnotesize $A$};
    \node[below] at (6/5,-10){\footnotesize $B$};

    \draw[>=latex, ->, thick, dashed] (5,-10.5) -- (5,-11.5);

    \draw[>=latex, ->, thick, purple, bend left = 30] (13,-1.5) to node[midway, right, xshift = 5] {$\mathcal{G}$ \footnotesize(bot)}(11.5, -6.5);

\end{tikzpicture}
    \caption{An example of how the accelerated map $\Tilde{\mathcal{G}}$ applies multiple moves of the Rauzy--Veech induction step $\Tilde{\mathcal{T}}$ at once. }
    \label{fig:fast_rauzy_veech_1}
\end{figure}
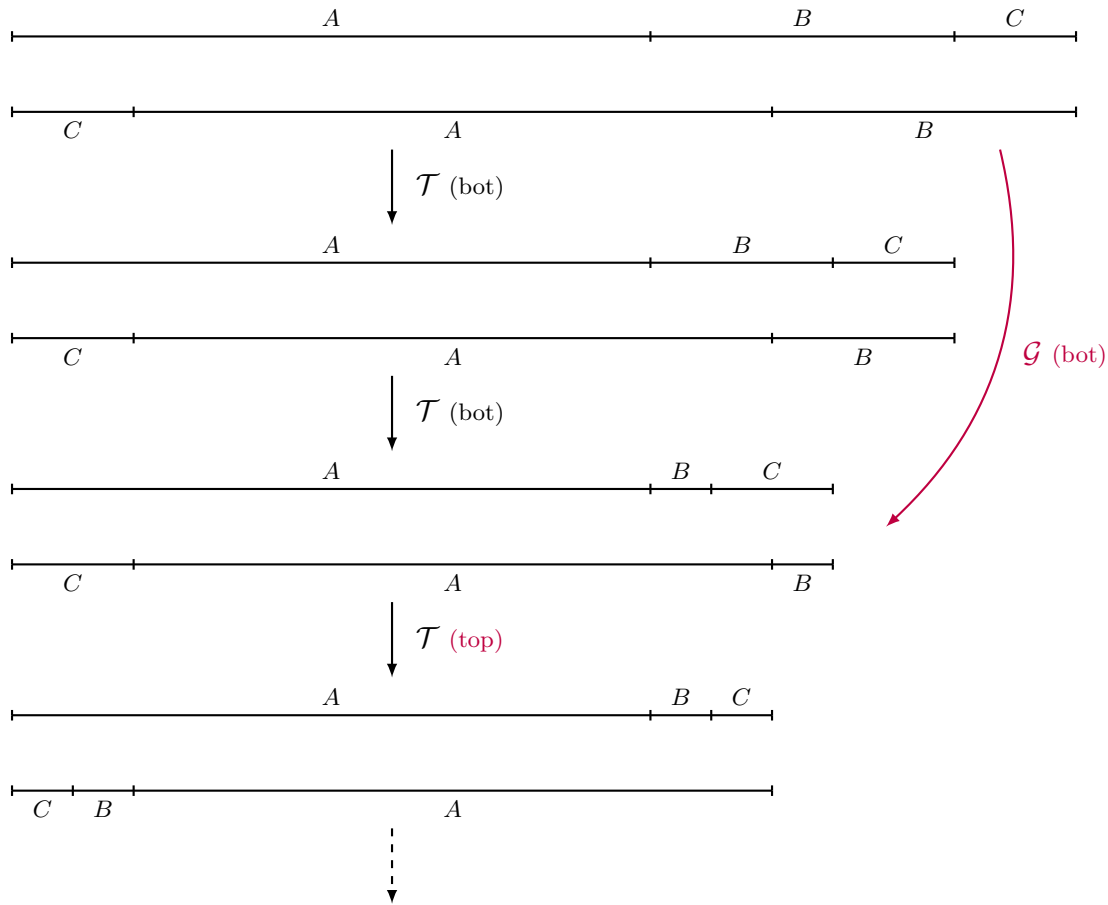

\begin{example}
    Consider the IET on 3 intervals given by the combinatorial and length data
    \begin{align*}
        \boldsymbol{\pi} &=
        \begin{pmatrix}
            A & B & C \\
            C & A & B
        \end{pmatrix}, \\
        \boldsymbol{\lambda} &= \left(\frac{3}{5}, \frac{2}{7}, \frac{4}{35}\right).
    \end{align*}
    For ease of notation we will consider the map $\Tilde{\mathcal{T}}$, i.e., only the induction step without renormalization and analogously its acceleration $\Tilde{\mathcal{G}}$ which is defined analogously to $\mathcal{G}$. Note that the IET characterized by $(\boldsymbol{\pi}, \boldsymbol{\lambda})$ is of bot type. Further, we have
    \begin{alignat*}{5}
        \boldsymbol{\pi}^1 &=
        \begin{pmatrix}
            A & B & C \\
            C & A & B
        \end{pmatrix}, \quad
        &&\boldsymbol{\pi}^2 &&= \begin{pmatrix}
            A & B & C \\
            C & A & B
        \end{pmatrix} \quad \text{and} \quad
        &&\boldsymbol{\pi}^3 &&= 
        \begin{pmatrix}
            A & B & C \\
            C & B & A
        \end{pmatrix},\\   
        \boldsymbol{\lambda}^1 &= \left(\frac{3}{5}, \frac{6}{35}, \frac{4}{35}\right),  
        &&\boldsymbol{\lambda}^2 &&= \left(\frac{3}{5}, \frac{2}{35}, \frac{4}{35}\right) \quad \text{and}  &&\boldsymbol{\lambda}^3 &&= \left(\frac{3}{5}, \frac{2}{35}, \frac{2}{35}\right). 
    \end{alignat*}
    We see that $(\boldsymbol{\pi}^1, \boldsymbol{\lambda}^1)$ is again of bot type and $(\boldsymbol{\pi}^2, \boldsymbol{\lambda}^2)$ is of top type. By definition, the accelerated map $\Tilde{\mathcal{G}}$ applied to $(\boldsymbol{\pi}, \boldsymbol{\lambda})$ is therefore given by
    \begin{equation*}
        \Tilde{\mathcal{G}}(\boldsymbol{\pi}, \boldsymbol{\lambda}) = \Tilde{\mathcal{T}}^2(\boldsymbol{\pi}, \boldsymbol{\lambda}).
    \end{equation*}
    In words, $\Tilde{\mathcal{G}}$ applies two moves of $\Tilde{\mathcal{T}}$ at once. This example is illustrated in Figure \ref{fig:fast_rauzy_veech_1}.
\end{example}

Also in \cite{zorich1996finite}, the following crucial result is proven which shows that the accelerated map accomplishes exactly what we intended.

\begin{theorem}\label{thm:fast_rauzy_veech_invariant_measure}
    For each Rauzy class $\mathfrak{R}$, the fast Rauzy--Veech induction $\mathcal{G} \colon \mathfrak{R} \times \Delta^{d-1} \to \mathfrak{R}\times \Delta^{d-1}$ admits an invariant \emph{probability} measure $\mu$ which is absolutely continuous with respect to $\dd \boldsymbol{\pi} \otimes \operatorname{Leb}$. This probability measure $\mu$ is unique and ergodic. Moreover, its density with respect to the Lebesgue measure is given by the restriction to $\mathfrak{R} \times \Delta^{d-1}$ of a function on $\mathfrak{R}\times \R_+^d$ which is a homogeneous rational function of degree $-d$ and bounded away from zero. 
\end{theorem}

Below in section \ref{sec:connecting_cf_to_rv_induction} we will work out explicitly how the maps $\mathcal{T}$ and $\mathcal{G}$ act in the simple case of the torus, and we will also be able to find the explicit densities $\nu$ from Theorem \ref{thm:rauzy_veech_invariant_measure} and $\mu$ from Theorem \ref{thm:fast_rauzy_veech_invariant_measure}. 

\subsubsection{Connecting Continued Fractions via Renormalization to Rauzy--Veech Induction}\label{sec:connecting_cf_to_rv_induction}

Finally, we will explain in detail the connection between the classical additive and multiplicative continued fraction algorithm and Rauzy--Veech induction performed on a flat surface of genus $\mathbf{g} = 1$. In this case, the space of IETs is parametrized by the unit interval. To see this connection clearly, let us work out an explicit formula for the Rauzy--Veech induction map $\mathcal{T} \colon (0,1) \to (0,1)$. Let us adopt the convention, that $\lambda$ is the length of the right interval (before the intervals are exchanged). 

\begin{figure}[ht]
    \centering
    \begin{tikzpicture}[scale = 0.8]
\begin{axis}[
    axis lines=middle,
    xmin=0, xmax=1,
    ymin=0, ymax=1,
    domain=0:1,
    samples=100, 
    legend pos=north west, 
]

\addplot[
purple,
domain=0:1/2,
samples=200,
unbounded coords=jump,
]{(x/(1-x))};
\addlegendentry{$\mathcal{T}(\lambda)$}

\addplot[
purple,
domain = 1/2:1,
samples = 200,
unbounded coords = jump,
]{2-1/x};
\end{axis}
\end{tikzpicture}
    \caption{The graph of $\mathcal{T}\colon (0,1) \to (0,1)$, the Rauzy--Veech renormalization map induced by the linear flow on the torus.}
    \label{fig:RVTorus_explicit_slow}
\end{figure}
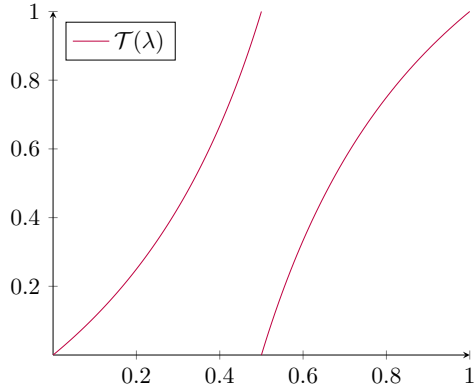

We claim that the map $\mathcal{T}$ is given by
\begin{align}\label{eq:RVTorus_explicit_slow}
    \begin{split}
    \mathcal{T}&\colon (0,1) \to (0,1), \\
    \lambda &\mapsto \begin{cases}
        \frac{\lambda}{1-\lambda} &\text{ if } \lambda < \frac{1}{2},\\
        2- \frac{1}{\lambda} &\text{ if } \lambda > \frac{1}{2}.
    \end{cases}
    \end{split}
\end{align}

A graph of $T$ can be seen in Figure \ref{fig:RVTorus_explicit_slow}. We will briefly explain how to derive this formula. Suppose that $\lambda > \frac{1}{2}$, which means that the IET is of top type. From \eqref{eq:top_move_length} we know that the length of the right interval after the exchange (before renormalization) is given by
\begin{equation*}
    \lambda' = \lambda - (1-\lambda) = 2\lambda - 1.
\end{equation*}
The new interval has length $1-(1-\lambda) = \lambda$ (since we remove a segment of length $1-\lambda$) so that after normalizing by $\frac{1}{\lambda}$ we obtain the corresponding expression in \eqref{eq:RVTorus_explicit_slow}. If $\lambda < \frac{1}{2}$, by \eqref{eq:bot_move_length} we know that the length of the right interval after the exchange is simply given by $\lambda$, which upon multiplication with the normalizing factor $\frac{1}{1-\lambda}$ yields the second case of \eqref{eq:RVTorus_explicit_slow}.

Note that the fact that the invariant measure constructed by Veech in \cite{veech1982gauss} is infinite corresponds to the fact that we have a neutral fixed point at the origin and at 1, meaning that $T'(0) =T'(1) = 1$. Therefore, any trajectory close to any of these points spends a long time there before it can escape, meaning that any invariant measure needs to assign a large (infinite, even) value to neighborhoods of 0. Let us give the density of the invariant measure explicitly.

\begin{proposition}[Invariant Measure of Rauzy--Veech Induction on Torus]\label{prop:slow_RV_induction_invariant_measure}
    Let $\mathcal{T}$ be the the slow Rauzy--Veech induction map from \eqref{eq:RVTorus_explicit_slow}. The invariant ergodic measure $\nu$ from Theorem \ref{thm:rauzy_veech_invariant_measure}, unique up to scaling, has density with respect to the Lebesgue measure which is given by
    \begin{equation*}
        \rho(x) = \frac{1}{2x(1-x)}.
    \end{equation*}
\end{proposition}
\begin{proof}
    Showing that $\frac{1}{2x(1-x)}\dd x$ is $\mathcal{T}$-invariant amounts to a straight-forward computation. Once we know it is invariant, Theorem \ref{thm:rauzy_veech_invariant_measure} implies it is ergodic and unique up to scaling.
\end{proof}
\begin{remark}
    Note that $\rho$ is the restriction of $(x,y) \mapsto \frac{1}{2xy}$ to $ \{(x,y) \in \R^2 \mid x + y = 1\}$, which is indeed a homogeneous rational function of degree -2. 
\end{remark}

Before moving to the acceleration $\mathcal{G}$, let us connect this map to the classical additive continued fraction algorithm, which is based on the \emph{Farey map}.

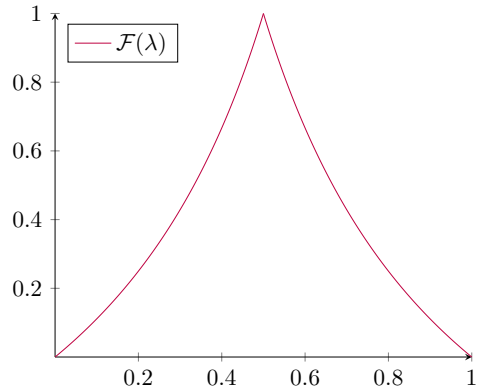
\begin{figure}[ht]
    \centering
    \begin{tikzpicture}[scale =0.8]
\begin{axis}[
    axis lines=middle,
    xmin=0, xmax=1,
    ymin=0, ymax=1,
    domain=0:1,
    samples=100, 
    legend pos=north west, 
]

\addplot[
purple,
domain=0:1/2,
samples=200,
unbounded coords=jump,
]{(x/(1-x))};
\addlegendentry{$\mathcal{F}(\lambda)$}

\addplot[
purple,
domain = 1/2:1,
samples = 200,
unbounded coords = jump,
]{(1-x)/x};
\end{axis}
\end{tikzpicture}
    \caption{The graph of the Farey map $F$ from Definition \ref{def:Farey_map}}
    \label{fig:farey_map}
\end{figure}

\begin{definition}[Farey Map]\label{def:Farey_map}
    The map $F \colon (0,1) \to (0,1)$ given by
    \begin{equation*}
        F(\lambda) = \begin{cases}
            \frac{\lambda}{1-\lambda} \text{ if } \lambda < \frac{1}{2}, \\
            \frac{1-\lambda}{\lambda} \text{ if } \lambda > \frac{1}{2}
        \end{cases}
    \end{equation*}
    is called the \emph{Farey map}. The graph of $F$ can be seen in Figure \ref{fig:farey_map}.
\end{definition}
\begin{remark}
    For our purposes, it is enough to define the Farey map on a set of full Lebesgue measure. Of course, one can easily extend the definition to the full closed interval $[0,1]$.
\end{remark}

\begin{figure}[hb]
    \centering
    \begin{tikzpicture}[scale =0.8]
\begin{axis}[
    width = .5\linewidth,
    axis lines=middle,
    xmin=0, xmax=1,
    ymin=0, ymax=10,
    domain=0:1,
    samples=100, 
    legend style={at={(0.5,1)}, anchor=north, yshift=-10pt}, 
]

\addplot[
purple,
domain=0.05:0.95,
samples=200,
]{1/(2*x*(1-x))};
\addlegendentry{$\rho(x) = \frac{1}{2x(1-x)}$}

\end{axis}
\end{tikzpicture}
\begin{tikzpicture}[scale = 0.8]
\begin{axis}[
    width = .5\linewidth,
    axis lines=middle,
    xmin=0, xmax=1,
    ymin=0, ymax=10,
    domain=0:1,
    samples=100, 
    legend style={at={(0.5,1)}, anchor=north, yshift=-10pt}, 
]

\addplot[
teal,
domain=0.05:1,
samples=200,
]{1/x};
\addlegendentry{$\frac{1}{x}$}

\end{axis}
\end{tikzpicture}

    \caption{The densities of the slow Rauzy--Veech induction map $\mathcal{T}\colon (0,1) \to (0,1)$ and the Farey map $F \colon (0,1) \to (0,1)$.}
    \label{fig:densities_slow_maps}
\end{figure}
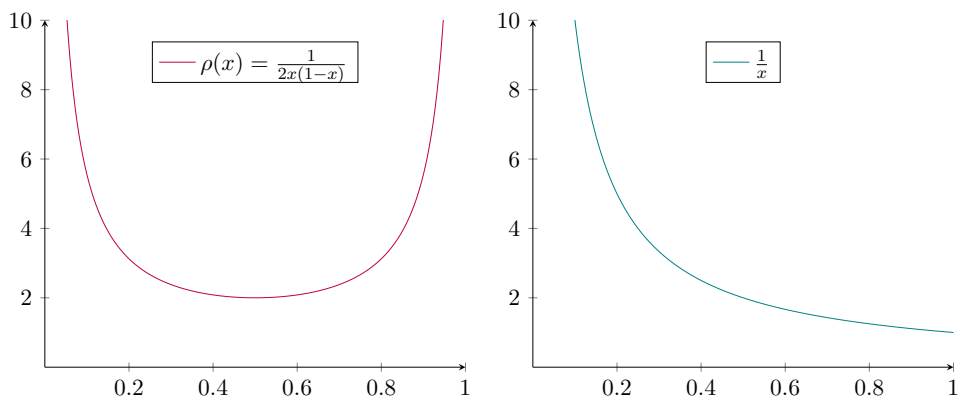

It is well-known that the Farey map preserves the (infinte) measure on $(0,1)$ given by $\frac{1}{x} \dd x$. The densities $\rho$ of the fast Rauzy--Veech induction and $\frac{1}{x}$ of the Farey map are depicted in Figure \ref{fig:densities_slow_maps}.

The connection of the Farey map to continued fractions is most easily seen in the following proposition. A proof can e.g., be found in \cite{dajani2021first}.

\begin{proposition}[Gauss Map is Acceleration of Farey Map]
    The Gauss map (Definition \ref{def:Gauss_map}) is an acceleration of the Farey map, i.e., it is given by the jump transformation of $F$ defined using the interval $(\frac{1}{2}, 1)$. 
\end{proposition}

The slow Rauzy--Veech induction on the torus and the Farey map are not equivalent in the sense of dynamical systems, but they are closely related. The Farey map $F$ is a factor of the slow Rauzy--Veech induction map $\mathcal{T}$ and interestingly, the factor map is given by the Farey map itself.

\begin{proposition}[Farey Map is a Factor of Slow Rauzy--Veech Induction Map]
The diagram
\[\begin{tikzcd}
	{(0,1)} && {(0,1)} \\
	\\
	{(0,1)} && {(0,1)}
	\arrow["{\mathcal{T}}", from=1-1, to=1-3]
	\arrow["F"', from=1-1, to=3-1]
	\arrow["F"', from=3-1, to=3-3]
	\arrow["F", from=1-3, to=3-3]
\end{tikzcd}\]

commutes, and if we write $\mu$ for the invariant measure of $\mathcal{T}$ from Proposition \ref{prop:slow_RV_induction_invariant_measure} and $\nu$ for the invariant measure of $F$, then
\begin{equation*}
    \nu = \mu \circ F^{-1}.
\end{equation*}
\end{proposition}

\begin{proof}
    The fact that $\mathcal{T}$ and $F$ are related in the above way can be checked by a straight-forward computation.
\end{proof}

We will now turn our attention to the accelerated Rauzy--Veech induction map $\mathcal{G}$. Denoting again by $\lambda$ the length of the right interval (before the exchange), the first question in order to find an explicit formula for $\mathcal{G}$ is \enquote{how many times do we apply a move of top (or bot) type, before the type changes?}.

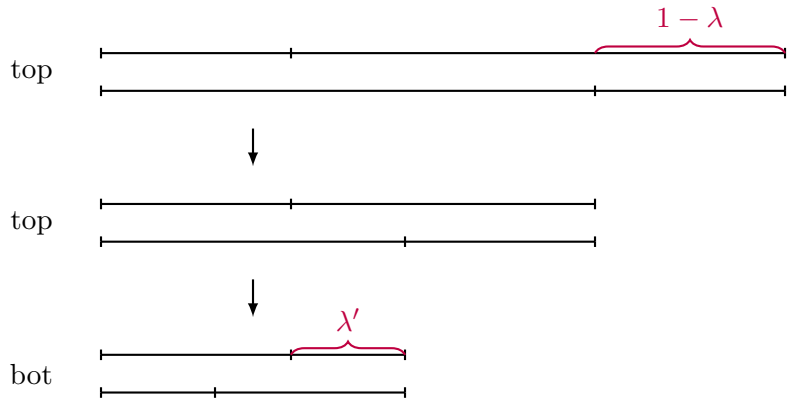
\begin{figure}[ht]
    \centering
    \begin{tikzpicture}[scale = 1, transform shape]
    \draw[thick] (0,0) -- (9,0);
    \draw[thick] (0,-1/2) -- (9,-1/2);

    \node[left] at (-1/2,-0.25) {top};

    \foreach \p in {0, 5/2, 9}{
        \draw[thick] (\p,0) - ++(0,0.07);
        \draw[thick] (\p,0) - ++(0,-0.07);
    }
    \foreach \p in {0, 13/2, 9}{
        \draw[thick] (\p,-1/2) - ++(0,0.07);
        \draw[thick] (\p,-1/2) - ++(0,-0.07);
    }
    \draw[>=latex, ->, thick] (2,-1) -- (2,-1.5);

    \draw [thick, purple, decorate, decoration={brace, amplitude=5pt}] (13/2,0) -- node[above=5pt] {$1-\lambda$}  (9,0);

    \draw[thick] (0,-2) -- (13/2,-2);
    \draw[thick] (0,-2.5) -- (13/2,-2.5);

    \node[left] at (-0.5,-2.25) {top};

    \foreach \p in {0, 5/2, 13/2}{
        \draw[thick] (\p,-2) - ++(0,0.07);
        \draw[thick] (\p,-2) - ++(0,-0.07);
    }
    \foreach \p in {0, 4, 13/2}{
        \draw[thick] (\p,-2.5) - ++(0,0.07);
        \draw[thick] (\p,-2.5) - ++(0,-0.07);
    }
    \draw[>=latex, ->, thick] (2,-3) -- (2,-3.5);

    \draw[thick] (0,-4) -- (4,-4);
    \draw[thick] (0,-4.5) -- (4,-4.5);

    \node[left] at (-0.5,-4.25) {bot};

    \foreach \p in {0, 5/2, 4}{
        \draw[thick] (\p,-4) - ++(0,0.07);
        \draw[thick] (\p,-4) - ++(0,-0.07);
    }
    \foreach \p in {0, 3/2, 4}{
        \draw[thick] (\p,-4.5) - ++(0,0.07);
        \draw[thick] (\p,-4.5) - ++(0,-0.07);
    }
    
    \draw [thick, purple, decorate, decoration={brace, amplitude=5pt}] (5/2,-4) -- node[above=5pt] {$\lambda'$}  (4,-4);
\end{tikzpicture}
    \caption{Action of a top move on an IET on 2 intervals. Each application decreases the length by $1-\lambda$.}
    \label{fig:RVTorus_explicit_fast_1}
\end{figure}

Say the IET is of top type, then we can picture the Rauzy--Veech induction step (without the renormalization) to decrease the length of the interval by $(1-\lambda)$ every time a move of top type is performed. The resulting IET is still of top type, as long as the right interval is longer than the left interval. Formally, this means that we can perform a move of top type $k$ times, where $k$ is the largest number such that
\begin{equation*}
    \lambda - k(1-\lambda) > 0,
\end{equation*}
i.e., 
\begin{equation*}
    k = \left\lfloor \frac{\lambda}{1-\lambda}\right\rfloor. 
\end{equation*}
Consider Figure \ref{fig:RVTorus_explicit_fast_1} for an example. Therefore, the length of the new interval on the right before renormalization is given by
\begin{equation*}
    \lambda' = \lambda - \left\lfloor \frac{\lambda}{1-\lambda}\right\rfloor = (1-\lambda)\left\{\frac{\lambda}{1-\lambda}\right\}.
\end{equation*}

Finally, normalizing by the length of the last interval, which is given by $\lambda' + (1-\lambda)$, we obtain
\begin{equation*}
    \mathcal{G}(\lambda) = 1 - \frac{1}{\left\{\frac{\lambda}{1-\lambda}\right\}+1}.
\end{equation*}
By a similar analysis of moves of bot type, we can conclude that $\mathcal{G}$ is given explicitly by
\begin{align}\label{eq:RVTorus_explicit_fast}
    \begin{split}
    \mathcal{G} &\colon (0,1) \to (0,1), \\
    \lambda &\mapsto \begin{cases}
        \frac{1}{\left\{\frac{1}{\lambda}\right\}+1} \quad &\text{if } \lambda < \frac{1}{2}, \\
        1 - \frac{1}{\left\{\frac{\lambda}{1-\lambda}\right\}+1} &\text{if } \lambda > \frac{1}{2}.
    \end{cases}
    \end{split}
\end{align}

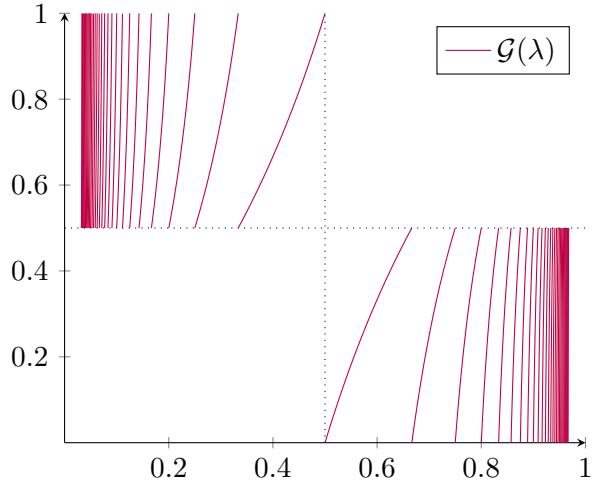
\begin{figure}[ht]
    \centering
    \begin{tikzpicture}
\begin{axis}[
    axis lines=middle,
    xmin=0, xmax=1,
    ymin=0, ymax=1,
    domain=0:1,
    samples=100, 
    legend pos=north east, 
]

\foreach \n in {2,...,30}
    \addplot[
    purple,
    domain = 1/(\n+1):1/\n,
    samples = 150,
    unbounded coords = jump,
    ]{1/(1/x - \n +1)};

\foreach \n in {2,...,30}
    \addplot[
    purple,
    domain = (1-1/\n):(1-1/(\n+1)),
    samples = 150,
    unbounded coords = jump,
    ]{1 - 1/(1/(1-x) - \n + 1))};

\draw[dotted] (axis cs:0.5,\pgfkeysvalueof{/pgfplots/ymin}) -- (axis cs:0.5,\pgfkeysvalueof{/pgfplots/ymax});

\draw[dotted] (axis cs:\pgfkeysvalueof{/pgfplots/xmin},0.5) -- (axis cs:\pgfkeysvalueof{/pgfplots/xmax},0.5);

\addlegendentry{$\mathcal{G}(\lambda)$}

\end{axis}
\end{tikzpicture}
    \caption{The fast Rauzy--Veech induction map $\mathcal{G} \colon (0,1) \to (0,1).$}
    \label{fig:fast_rauzy_veech_2}
\end{figure}

The graph of $\mathcal{G}$ can be seen in Figure \ref{fig:fast_rauzy_veech_2}. 

We are again able to provide an explicit expression for the invariant measure from Theorem \ref{thm:fast_rauzy_veech_invariant_measure} in the case of the torus.

\begin{proposition}[Invariant Measure of fast Rauzy--Veech Induction on the Torus]\label{prop:fast_RV_induction_invariant_measure}
    Let $\mathcal{G}\colon (0,1) \to (0,1)$ be the fast Rauzy--Veech induction map from \eqref{eq:RVTorus_explicit_fast}. The unique invariant ergodic probability measure $\nu$ from Theorem \ref{thm:fast_rauzy_veech_invariant_measure} has density with respect to the Lebesgue measure which is given by
    \begin{equation*}
        \rho(x) = \begin{cases}
            \frac{1}{2\log(2)} \frac{1}{1-x} \quad &\text{if } x < \frac{1}{2}, \\
            \frac{1}{2\log(2)} \frac{1}{x} \quad &\text{if } x > \frac{1}{2}.
        \end{cases}
    \end{equation*}
\end{proposition}
\begin{proof}
    It is easy to check by integrating $\rho$ on the unit interval to see that indeed it defines a probability measure. Moreover, given the formula above it is again a straight-forward computation to check that indeed $\mathcal{G}$ is $\rho(x)\dd x$-invariant, so uniqueness and ergodicity follow by Theorem \ref{thm:fast_rauzy_veech_invariant_measure}.  
\end{proof}

Let us remark that even though it is checking the invariance of $\rho$ is a routine verification once we are given the explicit expression, it is in general not easy to find invariant measures of a transformation. In our case here, Theorem \ref{thm:fast_rauzy_veech_invariant_measure} is constructive, so the proof yields a general formula for the invariant measure. In some special cases, namely when the combinatorial datum $\boldsymbol{\pi}$ has monodromy invariant
\begin{equation*}
    p = (d, d-1, \ldots, 2, 1),
\end{equation*}
the expression becomes particularly nice. More explicitly, it is shown in \cite{viana2006ergodic} that in this case the density of the fast Rauzy--Veech induction map is given by
\begin{equation*}
    \rho(\boldsymbol{\pi}, \boldsymbol{\lambda}) =
    \begin{cases}
        \displaystyle \frac{1}{d!} \prod_{\alpha \neq \alpha(\mathrm{top})} \left(\frac{1}{\lambda_\alpha + \lambda_{\alpha_\mathrm{top}^+}}\right)\cdot \frac{1}{\sum_{\beta \neq \alpha(\mathrm{bot})}\lambda_\beta}, \quad &\text{if }\lambda_{\alpha(\mathrm{top})} > \lambda_{\alpha(\mathrm{bot})},\\
        \displaystyle \frac{1}{d!} \prod_{\alpha \neq \alpha(\mathrm{bot})} \left(\frac{1}{\lambda_\alpha + \lambda_{\alpha_\mathrm{bot}^+}}\right)\cdot \frac{1}{\sum_{\beta \neq \alpha(\mathrm{top})}\lambda_\beta},  &\text{if }\lambda_{\alpha(\mathrm{top})} < \lambda_{\alpha(\mathrm{bot})},
    \end{cases}
\end{equation*}
where $\alpha_\varepsilon^+$ denotes the symbol to the right of $\alpha_\varepsilon$. 
So in the case $d = 2$ this gives
\begin{equation*}
    \rho(\lambda_A, \lambda_B) = 
    \begin{cases}
        \frac{1}{2\lambda_B} \quad &\text{if } \lambda_A < \lambda_B, \\
        \frac{1}{2\lambda_A} & \text{if } \lambda_A > \lambda_B,
    \end{cases}
\end{equation*}
which after multiplying with the normalizing constant $\log(2)$ coincides exactly with the density from Proposition \ref{prop:fast_RV_induction_invariant_measure}.

The Gauss map from Definition \ref{def:Gauss_map} and the fast Rauzy--Veech induction map are related in exactly the same way as the Farey map and the slow Rauzy--Veech induction map, i.e., the former is a factor of the latter. Even more is true, since the factor map is again given by the Farey map $F$. To relate the to maps fully, we also need an invariant measure for the Gauss map. 

\begin{proposition}
    The probability measure $\nu$ defined by
    \begin{equation*}
        \nu(A) = \frac{1}{\log(2)}\int_A \frac{1}{1+x} \, \dd x 
    \end{equation*}
    preserves the Gauss map $G \colon (0,1) \to (0,1), \, x\mapsto  \left\{\frac{1}{x}\right\}$.
\end{proposition}

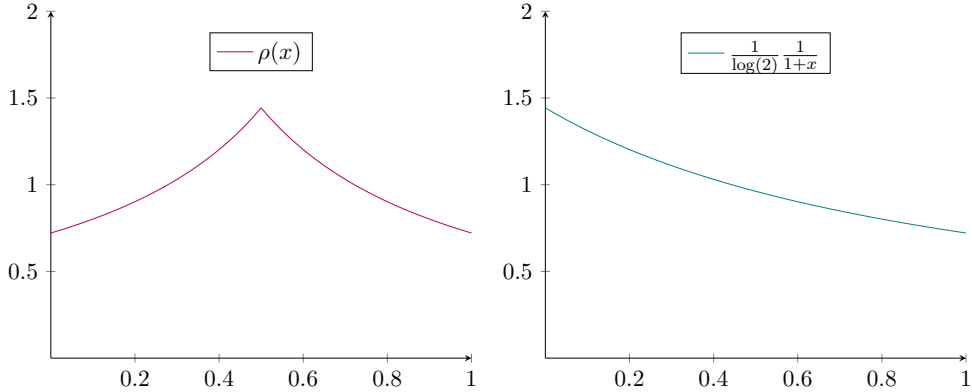
\begin{figure}[ht]
    \centering
    \begin{tikzpicture}[scale = 0.8]
\begin{axis}[
    width = .5\linewidth,
    axis lines=middle,
    xmin=0, xmax=1,
    ymin=0, ymax=2,
    domain=0:1,
    samples=100, 
    legend style={at={(0.5,1)}, anchor=north, yshift=-10pt}, 
]

\addplot[
purple,
domain=0:0.5,
samples=100,
]{1/(2*ln(2)) * 1/(1-x)};
\addlegendentry{$\rho(x)$}
\addplot[
purple,
domain= 0.5:1,
samples=100,
]{1/(2*ln(2)) * 1/x};

\end{axis}
\end{tikzpicture}
\begin{tikzpicture}[scale =0.8]
\begin{axis}[
    width = .5\linewidth,
    axis lines=middle,
    xmin=0, xmax=1,
    ymin=0, ymax=2,
    domain=0:1,
    samples=100, 
    legend style={at={(0.5,1)}, anchor=north, yshift=-10pt}, 
]

\addplot[
teal,
domain=0:1,
samples=200,
]{1/(ln(2)) * 1/(x+1)};
\addlegendentry{$\frac{1}{\log(2)} \frac{1}{1+x}$}

\end{axis}
\end{tikzpicture}

    \caption{The densities of the fast Rauzy--Veech induction map $\mathcal{G}\colon (0,1) \to (0,1)$ and the Gauss map $G \colon (0,1) \to (0,1)$.}
    \label{fig:densities_fast_maps}
\end{figure}

The graphs of the invariant densities of the fast Rauzy--Veech induction map and the Gauss map can be seen in Figure \ref{fig:densities_fast_maps}.

\begin{proposition}[Gauss Map is a Factor of Fast Rauzy--Veech Induction Map]\label{prop:Gauss_factor_of_RV} The diagram

\[\begin{tikzcd}
	{(0,1)} && {(0,1)} \\
	\\
	{(0,1)} && {(0,1)}
	\arrow["{\mathcal{G}}", from=1-1, to=1-3]
	\arrow["F"', from=1-1, to=3-1]
	\arrow["G"', from=3-1, to=3-3]
	\arrow["F", from=1-3, to=3-3]
\end{tikzcd}\]

commutes, and if we write $\mu$ for the invariant measure of $\mathcal{G}$ from Proposition \ref{prop:fast_RV_induction_invariant_measure} and $\nu$ for the Gauss measure, then
\begin{equation*}
    \nu = \mu \circ F^{-1}.
\end{equation*}
\end{proposition}
\begin{proof}
    It is possible to again verify the claims doing the computations. Alternatively, we can also verify that the map $n$ defined in \eqref{eq:fast_RV}, which was used to define the fast Rauzy--Veech induction map, coincides exactly with the first passage time associated to the acceleration of the Farey map yielding the Gauss map. Then, the fact that the above diagram commutes as well as the condition on the measures follows from the general theory of jump transformations.
\end{proof}

It is clear that Proposition \ref{prop:Gauss_factor_of_RV} implies that the fast Rauzy--Veech induction map can be used to find the continued fraction expansion of $\lambda \in (0,1)$. Analogously to the classical case explained in section \ref{sec:cf_torus} we just need to record in which branch of the map $\mathcal{G}$ the iterates of $\lambda$ lie. Because the factor map $F$ from \ref{prop:Gauss_factor_of_RV} is 2-to-1, for any $n \in \N$ we obtain two branches $P_n = \left(\frac{1}{n+2}, \frac{1}{n+1}\right)$ and $Q_n = \left(\frac{n}{n+1}, \frac{n+1}{n+2}\right)$, which is clearly depicted in Figure \ref{fig:fast_rauzy_veech_2}. To compute the continued fraction expansion of $\lambda \in (0,1)$, we must iterate $\mathcal{G}$ on one of the two preimages of $F$, i.e., we compute $\mathcal{G}^n(\Tilde{\lambda})$, where $\Tilde{\lambda} \in F^{-1}(\{\lambda\})$. It is clear that both preimages yield the same result. We will end this section by reexamining the example from section \ref{sec:cf_torus}, i.e., we will use the fast Rauzy--Veech induction map to compute the continued fraction expansion of $\lambda = \sqrt{6}-2$.

\begin{example}[Continued Fractions via fast Rauzy--Veech]
    Let $\lambda = \sqrt{6}-2 \in (0,1)$. To compute the continued fraction expansion of $\lambda$ using (fast) Rauzy--Veech induction, we first choose a preimage $\Tilde{\lambda} \in F^{-1}(\{\lambda\}) = \left\{\frac{1}{5}(4-\sqrt{6}), \frac{1}{5}(1+\sqrt{6})\right\}$. Say we pick $\Tilde{\lambda} = \frac{1}{5}(4-\sqrt{6})$, then
    \begin{equation*}
        \frac{1}{4}<\Tilde{\lambda}<\frac{1}{3} \quad \Rightarrow \Tilde{\lambda} \in P_2.
    \end{equation*}
    Therefore, the first digit in the continued fraction expansion of $\lambda$ is 2. Now we compute
    \begin{equation*}
        \mathcal{G}(\Tilde{\lambda}) = \frac{\sqrt{6}}{3}, 
    \end{equation*}
    for which we have
    \begin{equation*}
        \frac{4}{5}<\frac{\sqrt{6}}{3}<\frac{5}{6},
    \end{equation*}
    hence $\frac{\sqrt{6}}{3} \in Q_4$ and the second digit is 4. Since
    \begin{equation*}
        \mathcal{G}\left(\frac{\sqrt{6}}{3}\right) = \frac{1}{5}(4-\sqrt{6}) = \Tilde{\lambda},
    \end{equation*}
    we again see that the expansion is periodic, hence
    \begin{equation*}
        \lambda = [2,4,2,4,\ldots].
    \end{equation*}
\end{example}

\subsection{Rauzy-Veech Induction as Discretization of Teichmüller Geodesic Flow}\label{sec:RV_Teichmuller_connection}

The final property of the Rauzy--Veech induction map that we want to investigate is the fact that it is given by a Poincaré first return map of the Teichmüller geodesic flow (see Definition \ref{def:teichmueller_geodesic_flow}) to a carefully chosen transverse section. 

There are several ways in which we can parametrize a translation surface $X$ (locally) within its stratum. One such possibility uses the zippered rectangle decomposition introduced in section \ref{sec:RV_induction}. We will give a brief description of these coordinates following \cite{zorich2006flat}. For a more formal treatment, we refer the reader to \cite{veech1982gauss} and to \cite{zorich1996finite}. We obtain local coordinates as follows.

We start by choosing a horizontal segment satisfying the MIC (Definition \ref{def:mic}) which we have seen induces a decomposition into zippered rectangles. 

\begin{figure}[ht]
    \centering
    \begin{tikzpicture}[scale = 0.32]
    \draw (0,0) rectangle (1,6); 
    \fill[red, opacity = 0.15] (0,0) rectangle (1,6);
    \draw (1,0) rectangle (3,7.3); 
    \fill[violet, opacity = 0.15] (1,0) rectangle (3, 7.3);
    \draw (3,0) rectangle (7.5, 12); 
    \fill[blue, opacity = 0.15] (3,0) rectangle (7.5, 12);
    \draw (7.5,0) rectangle (14,9); 
    \fill[teal, opacity = 0.15] (7.5,0) rectangle (14,9);

    \draw[fill] (0,6) circle [radius=2pt] node[left] {$h_1$};
    \draw[fill] (1,7.3) circle [radius=2pt] node[left] {$h_2$};
    \draw[fill] (3,12) circle [radius=2pt] node[left] {\textcolor{BrickRed!90!black!100!}{$a_2=$} $h_3$};
    \draw[fill] (7.5,9) circle [radius=2pt] node[left] {$h_4$};

    \draw[fill, BrickRed!90!black!100!] (1,4) circle [radius=2pt] node[right] {$a_1$};
    \draw[fill, BrickRed!90!black!100!] (7.5,7) circle [radius=2pt] node[right] {$a_3$};
    \draw[fill, BrickRed!90!black!100!] (14,6) circle [radius=2pt] node[right] {$a_4$};

    \draw[dashed, white] (0,-1) -- (0,-2.4);
\end{tikzpicture}
\begin{tikzpicture}[scale = 0.4]
    \draw (0,0) rectangle (2,8);
    \fill[blue, opacity = 0.15] (0,0) rectangle (2,8);
    \draw (2,0) rectangle (7,5);
    \fill[red, opacity = 0.15] (2,0) rectangle (7,5);
    \draw (7,0) rectangle (11,7);
    \fill[violet, opacity = 0.15] (7,0) rectangle (11,7);
    \draw (11,0) rectangle (12.5, 6);
    \fill[teal, opacity = 0.15] (11,0) rectangle (12.5, 6);

    \draw[fill] (0,8) circle [radius=2pt] node[left] {$h_1$};
    \draw[fill] (2,5) circle [radius=2pt] node[left] {$h_2$};
    \draw[fill] (7,7) circle [radius=2pt] node[left] {$h_3$};
    \draw[fill] (11,6) circle [radius=2pt] node[left] {$h_4$};

    \draw[fill, BrickRed!90!black!100!] (2,4) circle [radius=2pt] node[right] {$a_1$};
    \draw[fill, BrickRed!90!black!100!] (7,3.5) circle [radius=2pt] node[right] {$a_2$};
    \draw[fill, BrickRed!90!black!100!] (11,3.75) circle [radius=2pt] node[right] {$a_3$};
    \draw[fill, BrickRed!90!black!100!] (12.5,-1.5) circle [radius=2pt] node[right] {$a_4$};
    \draw[dashed] (12.5,0) -- (12.5, -1.5);
\end{tikzpicture}
    \label{fig:rectangle_coordinates}
    \caption{Zippered rectangle decompositions with coordinates.}
\end{figure}
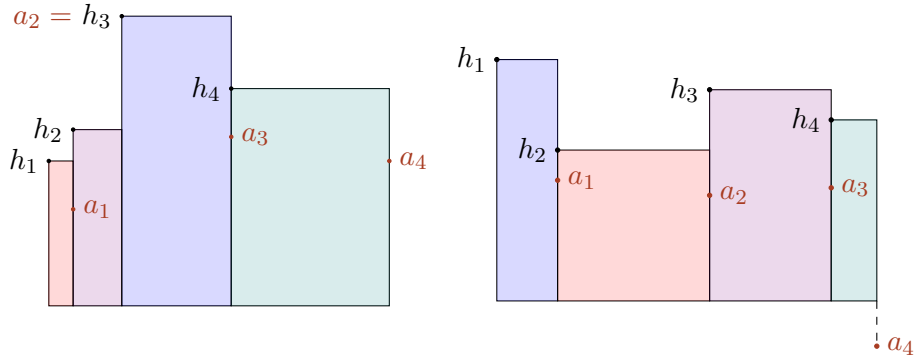

\begin{definition}[Zippered Rectangle Coordinates]
    Given a decomposition of a translation surface $X$ into zippered rectangles, we will call the following data the \emph{rectangle coordinates} of $X$.
    \begin{itemize}
        \item The \emph{widths} of the rectangles are denoted by $\boldsymbol{\lambda} = (\lambda_1, \ldots, \lambda_d)$.
        \item The \emph{heights} of the rectangles are denoted by
        $\boldsymbol{h} = (h_1, \ldots, h_d)$.
        \item The \emph{altitudes responsible for the position of the singularities} are denoted by $\boldsymbol{a} = (a_1, \ldots, a_d)$.
        \item A \emph{pair of permutations} $\boldsymbol{\pi} = (\topp, \bott)$, which is a discrete parameter. 
    \end{itemize}
\end{definition}
In Figure \ref{fig:rectangle_coordinates} we have two zippered rectangles, where the coordinates are marked as well. The left-hand side is the decomposition from Example \ref{ex:octagon_IET}, while the right-hand side comes from Example \ref{ex:double_pentagon_IET_2}. 

Note that the heights of the rectangles and the altitudes of the singularities are not independent. They satisfy some linear constraints which can be found in \cite{veech1982gauss}, where the reader can also find the proof of the following theorem.

\begin{theorem}[\cite{veech1982gauss}]\label{thm:suspensions_existence}
    For any $(\boldsymbol{\pi}, \boldsymbol{\lambda}) \in \mathfrak{R}\times \Delta^{d-1}$ in the space of IETs on $d$ intervals there is an $m$-dimensional open linear cone of solutions $(\boldsymbol{h}, \boldsymbol{a})$, where $m$ is a positive even integer. 
\end{theorem}

\begin{remark}
    The integer $m$ depends on the combinatorial datum $\boldsymbol{\pi}$.
\end{remark}

Put differently, Theorem \ref{thm:suspensions_existence} tells us that starting with any IET $T\colon S \to S$ we can always find a translation surface $X$ and a horizontal segment $S \subseteq X$ inside such that the first return of the vertical flow to $S$ is exactly the IET $T$. Moreover, in general we may find multiple essentially different translation surfaces with this property: There is an $m$-dimensional family of such flat surfaces, meaning that there exists a collection $(\boldsymbol{h}_1, \boldsymbol{a}_1, \ldots, \boldsymbol{h}_m, \boldsymbol{a}_m)$ of height and altitude data such that all translation surfaces with this property are given by rectangular coordinates that are linear combinations of the $(\boldsymbol{h}_i, \boldsymbol{a}_i)$. 

\begin{definition}[Suspensions over an IET]
    The $d$-dimensional family of translation surfaces just described are called \emph{suspensions} over the IET $T \colon S \to S$. 
\end{definition}

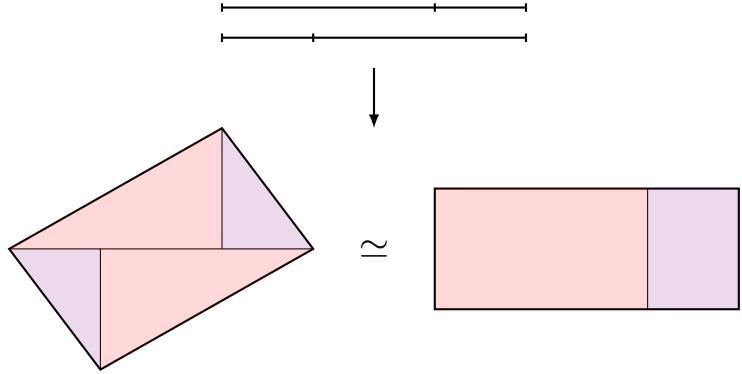
\begin{figure}[ht]
    \centering
    \begin{tikzpicture}[scale = 0.8]
    \draw[thick] (-2.5,4) -- (2.5,4);
    \draw[thick] (-2.5,3.5) -- (2.5,3.5);

    \foreach \p in {-2.5, 1, 2.5}{
        \draw[thick] (\p,4) - ++(0,0.07);
        \draw[thick] (\p,4) - ++(0,-0.07);   
    }
    \foreach \p in {-2.5, -1, 2.5}{
        \draw[thick] (\p,3.5) - ++(0,0.07);
        \draw[thick] (\p,3.5) - ++(0,-0.07);   
    }

    \draw[>=latex, ->, thick] (0,3) -- (0,2);

    \draw[thick] (-6, 0) -- (-2.5,2) -- (-1, 0) -- (-4.5,-2) -- cycle;
    \draw (-6,0) -- (-1,0); \draw (-2.5,2) -- (-2.5, 0); \draw (-4.5, -2) -- (-4.5,0);

    \fill[red, opacity = 0.15] (-6,0) -- (-2.5,2) -- (-2.5,0) -- (-1,0) -- (-4.5,-2) -- (-4.5,0) -- cycle; 

    \fill[violet, opacity = 0.15] (-6,0) -- (-4.5,-2) -- (-4.5,0) -- cycle;
    \fill[violet, opacity = 0.15] (-1,0) -- (-2.5, 2) -- (-2.5, 0) --cycle;

    \node at (0,0) {\Large $\simeq$};

    \draw[thick] (1,-1) rectangle (6,1); \draw (4.5, -1) -- (4.5, 1);
    \fill[red, opacity = 0.15] (1,-1) rectangle (4.5, 1);
    \fill[violet, opacity = 0.15] (4.5, -1) rectangle (6, 1);
\end{tikzpicture}
    \caption{The first type of rectangular coordinates from Example \ref{ex:suspension_torus}.}
    \label{fig:suspension_torus_1}
\end{figure}

\begin{figure}[ht]
    \centering
    \begin{tikzpicture}[scale = 0.8]
    \begin{scope}
    \draw[thick] (-2.5,4) -- (2.5,4);
    \draw[thick] (-2.5,3.5) -- (2.5,3.5);

    \foreach \p in {-2.5, 1, 2.5}{
        \draw[thick] (\p,4) - ++(0,0.07);
        \draw[thick] (\p,4) - ++(0,-0.07);   
    }
    \foreach \p in {-2.5, -1, 2.5}{
        \draw[thick] (\p,3.5) - ++(0,0.07);
        \draw[thick] (\p,3.5) - ++(0,-0.07);   
    }

    \draw[>=latex, ->, thick] (0,3) -- (0,2);

    \draw[thick] (-7.5, 0) -- (-2.5,2) -- (-1, 0) -- (-6,-2) -- cycle;
    \clip (-7.5, 0) -- (-2.5,2) -- (-1, 0) -- (-6,-2) -- cycle;
    \draw (-7.5,0) -- (-2.5,0); \draw (-2.5,2) -- (-2.5, -1); \draw (-6, -2) -- (-6,0);
    \draw (-4,0) -- (-4,2);

    \fill[red, opacity = 0.15] (-7.5,0) rectangle (-4,2);
    \fill[red, opacity = 0.15] (-6,-2) rectangle (-2.5,0);

    \fill[violet,opacity =0.15] (-7.5,0) rectangle (-6,-2);
    \fill[violet, opacity = 0.15] (-4,0) rectangle (-2.5,2);
    \fill[violet, opacity = 0.15] (-2.5,2) rectangle (-1,-1);
    \end{scope}

    \node at (0,0) {\Large $\simeq$};

    \begin{scope}
        \draw[thick] (1,-2) -- (1,0) -- (6,-0) -- (6,2) -- (7.5,2) -- (7.5,-2) -- cycle;
        \draw (6,0) -- (6,-2);
        \fill[red, opacity = 0.15] (1,-2) rectangle (6, 0);
        \fill[violet, opacity = 0.15] (6, -2) rectangle (7.5, 2);
    \end{scope}

\end{tikzpicture}
    \caption{The second type of rectangular coordinates from Example \ref{ex:suspension_torus}.}
    \label{fig:suspension_torus_2}
\end{figure}
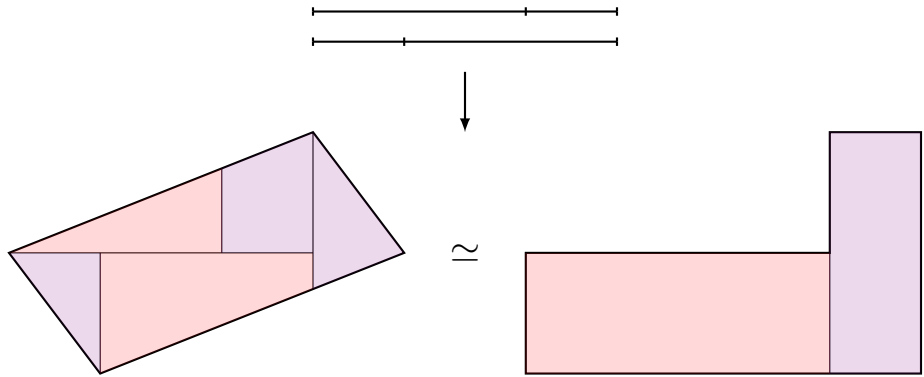

\begin{example}[Suspensions over IET on 2 Intervals]\label{ex:suspension_torus}
    Consider an IET $T \colon S \to S$ on two intervals. Since the length of $S$ can be arbitrary, such IET is parametrized by $\boldsymbol{\lambda} = (\lambda_1, \lambda_2) \in \R_{>0}^2$. An explicit formula for the dimension $m$ of the solution space is given in \cite{veech1982gauss} and it follows that there are 2 essentially different translation surfaces with horizontal segment $S$ satisfying the MIC such that the Poincaré first return map to $S$ is exactly given by $T$. Let us give the rectangle coordinates of each of these translation surfaces.  

    \begin{enumerate}
        \item Let $\xi \in \R_{>0}$. We can choose the following height and altitude data.
        \begin{align*}
            \boldsymbol{h}_1 &= (\xi, \xi), \\
            \boldsymbol{a}_1 &= (\xi, 0).
        \end{align*}
        The corresponding translation surface can be seen in Figure \ref{fig:suspension_torus_1}.

        \item Let again $\xi \in \R_{>0}$. The following choice gives another translation surface for which the rectangle coordinates are not a scalar multiple of the coordinates in the first example. 
        \begin{align*}
            \boldsymbol{h}_2 &= (\xi, 2\xi), \\
            \boldsymbol{a}_2 &= (2\xi, \xi).
        \end{align*}
        Figure \ref{fig:suspension_torus_2} shows an example of such a translation surface. 
    \end{enumerate}
\end{example}
\begin{remark}
    It is evident, that the height and altitude data given above are linearly independent, i.e., that for any choice of $\xi > 0$, these two sets of data constitute a basis of the family of translation surfaces from Theorem \ref{thm:suspensions_existence}.
\end{remark}

The situation becomes even more interesting if we consider IETs of 3 intervals. We give a basis of the cone of solutions in the next example.

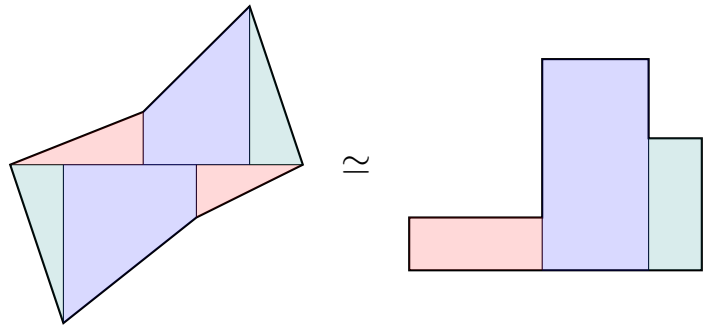
\begin{figure}[ht]
    \centering
    \begin{tikzpicture}[scale = 0.7]
    \begin{scope}
        \draw (0,0) -- (5.5,0);
        \draw[thick] (0,0) -- (2.5,1) -- (4.5,3) -- (5.5,0) -- (3.5,-1) -- (1,-3) -- cycle;

        \draw (2.5,0) -- (2.5,1); \draw (4.5,0) -- (4.5, 3); \draw (1,-3) -- (1,0); \draw (3.5, -1) -- (3.5, 0);

        \clip (0,0) -- (2.5,1) -- (4.5,3) -- (5.5,0) -- (3.5,-1) -- (1,-3) -- cycle;

        \fill[red, opacity = 0.15] (0,0) rectangle (2.5,1);
        \fill[red, opacity = 0.15] (5.5,0) rectangle (3.5,-1);

        \fill[blue, opacity = 0.15] (2.5,0) rectangle (4.5,3);
        \fill[blue, opacity = 0.15] (3.5,0) rectangle (1,-3);
        
        \fill[teal, opacity = 0.15] (0,0) rectangle (1,-3);
        \fill[teal, opacity = 0.15] (5.5,0) rectangle (4.5,3);
    \end{scope}
    
    \node at (6.5,0) {\Large $\simeq$};
    
    \begin{scope}
        \draw[thick] (7.5,-2) -- (7.5,-1) -- (10,-1) -- (10,2) -- (12,2) -- (12,0.5) -- (13,0.5)-- (13,-2) -- cycle;

        \draw (10,-1) -- (10,-2); \draw (12,2) -- (12,-2);

        \fill[red, opacity = 0.15] (7.5, -2) rectangle (10,-1);
        \fill[blue, opacity = 0.15] (10,-2) rectangle (12,2);
        \fill[teal, opacity = 0.15] (12,-2) rectangle (13,0.5);
    \end{scope}
\end{tikzpicture}
    \caption{The first type of rectangular coordinates from Example \ref{ex:suspension_hexagon}.}
    \label{fig:suspension_hexagon_1}
\end{figure}

\begin{figure}[hb]
    \centering
    \begin{tikzpicture}[scale = 0.7]
    \begin{scope}
        \draw (0,0) -- (5.5,0);
        \draw[thick] (0,0) -- (2.5,0.5) -- (4.5,1.5) -- (5.5,0.8) -- (7.5,0) -- (5.5, -1) -- (3,-1.5) -- (1,-0.7) -- cycle;
        \clip (0,0) -- (2.5,0.5) -- (4.5,1.5) -- (5.5,0.8) -- (7.5,0) -- (5.5, -1) -- (3,-1.5) -- (1,-0.7) -- cycle;

        \draw (5.5, 1) -- (5.5, -1);
        \draw (2.5,0) -- (2.5,2);
        \draw (4.5, 0) -- (4.5, 1.5);
        \draw (1,0) -- (1,-2);
        \draw (3,0) -- (3,-5);

        \fill[red, opacity = 0.15] (0,0) rectangle (2.5,0.5);
        \fill[red, opacity = 0.15] (3,-1.5) rectangle (5.5,0);

        \fill[blue, opacity = 0.15] (2.5,0) rectangle (4.5,1.5);
        \fill[blue, opacity = 0.15] (3,-1.5) rectangle (1,0);
        \fill[blue, opacity = 0.15] (5.5,-1) rectangle (7.5,1);

        \fill[teal, opacity = 0.15] (0,0) rectangle (1,-1);
        \fill[teal, opacity = 0.15] (4.5,1.5) rectangle (5.5,0);
        
    \end{scope}

    \node at (8.5,0) {\Large $\simeq$};

    \begin{scope}
        \draw[thick] (9.5,-1.5) -- (9.5,0) -- (12,0) -- (12, 1.5) -- (14,1.5) -- (14, 0) -- (15,0) -- (15,-1.5) -- cycle;
        \clip (9.5,-1.5) -- (9.5,0) -- (12,0) -- (12, 1.5) -- (14,1.5) -- (14, 0) -- (15,0) -- (15,-1.5) -- cycle;
        \draw (12,-1.5) -- (12,0); \draw (14,-1.5) -- (14, 0);
        \fill[red, opacity = 0.15] (9.5,-1.5) rectangle (12,0);
        \fill[blue, opacity= 0.15] (12,-1.5) rectangle (14,1.5);
        \fill[teal, opacity = 0.15] (14,-1.5) rectangle (15,0);
    \end{scope}

\end{tikzpicture}
    \caption{The second type of rectangular coordinates from Example \ref{ex:suspension_hexagon}.}
    \label{fig:suspension_hexagon_2}
\end{figure}
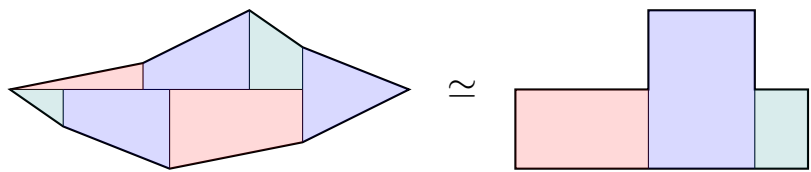

\begin{example}[Suspensions over IET on 3 intervals]\label{ex:suspension_hexagon}
    Consider an IET $T \colon S \to S$ on three intervals, where we parametrize the IET by
    \begin{align*}
        \boldsymbol{\lambda} &= (\lambda_1, \lambda_2, \lambda_3), \\
        \boldsymbol{\pi} &=
        \begin{pmatrix}
            A & B & C \\
            C & B & A
        \end{pmatrix}.
    \end{align*}
    Using again the explicit form of the dimension $m$ of the solution space from \cite{veech1982gauss}, we see that here we have $m = 2$. Let us give two pairs of vectors $(\boldsymbol{h}_i, \boldsymbol{a}_i), i \in \{1,2\}$ which form a basis of the solution space from \ref{thm:suspensions_existence}.
    \begin{enumerate}
        \item The first height and altitude data we choose is given by
        \begin{align*}
            \boldsymbol{h}_1 &= (1, 3, 2), \\
            \boldsymbol{a}_1 &= (1,2,0).
        \end{align*}
        A corresponding translation surface can be seen in Figure \ref{fig:suspension_hexagon_1}.
        \item The second height and altitude data we choose is given by
        \begin{align*}
            \boldsymbol{h}_2 &= (2, 4, 2), \\
            \boldsymbol{a}_2 &= (3,2,1).
        \end{align*}
        A corresponding translation surface can be seen in Figure \ref{fig:suspension_hexagon_2}.
    \end{enumerate}
\end{example}

We have found a way to move from IETs to translation surfaces and back, but we saw that switching between these two concepts requires making a choice. We have discussed the impact of this choice above when moving from an IET to a suspension by witnessing that there is a solution space of some positive explicit dimension. 

The next natural question to ask is that given some translation surface $X$, is there a canonical choice for a zippered rectangle decomposition (which would in particular translate to a canonical choice of IET)? Unfortunately, there seems to be no such canonical choice of decomposition. In \cite{zorich2006flat}, an \emph{almost} canonical choice is proposed which comes from the original construction of zippered rectangles in \cite{veech1982gauss}, leaving only an \enquote{arbitrariness of finite order}, which means that for any translation surface $X$ there are only finitely many zippered rectangle decompositions to choose from. Let us explain how such a choice is made.

Choosing a horizontal segment $S \subseteq X$ completely determines a zippered rectangle decomposition, so we may as well describe how we choose such a segment. The choice is done as follows: We put the left endpoint of $S$ at one of the singularities and choose the length $|S|$ of the segment such that $S$ is the shortest possible interval satisfying the MIC with $|S| \geq 1$. Since each translation surface has a finite number of conical singularities and each singularity emits only a finite number of horizontal segments, it is clear that the number of possible choices is finite. 

\begin{example}[Choice of decomposition for the Torus]
    Note that since a translation surface of genus $\mathbf{g} = 1$, i.e., a torus, has exactly one conical singularity of cone angle $2\pi$ we do get a canonical choice of zippered rectangle in this case, since the choice described above is unique. Depending on the specific coordinates, the decomposition will either as in Figure \ref{fig:suspension_torus_1} or as in Figure \ref{fig:suspension_torus_2}.
\end{example}

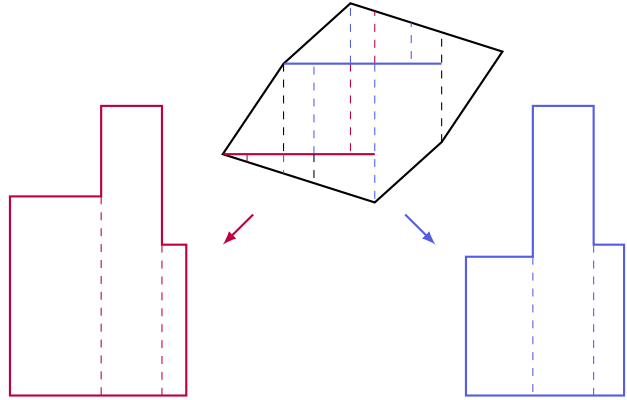
\begin{figure}[ht]
    \centering
    \begin{tikzpicture}[scale = 0.8]
\begin{scope}
    \coordinate (p1) at (0,0);
    \coordinate (p2) at (1,1.5);
    \coordinate (p3) at (2.1,2.5);
    \coordinate (p4) at (4.6,1.7);

    \coordinate (p5) at (3.6,0.2);
    \coordinate (p6) at (2.5, -0.8);

    \draw[thick] (p1) -- (p2) --(p3) -- (p4) -- (p5) -- (p6) -- cycle;
    \clip (p1) -- (p2) --(p3) -- (p4) -- (p5) -- (p6) -- cycle;

    \draw[thick, purple] (p1) -- ($(p1) + (2.5,0)$);
    \draw[dashed, green!10!blue!95!red!70] ($(p1) + (2.5,1.5)$) -- (p6);
    \draw[dashed, purple] ($(p1) + (2.5,1.5)$) -- ($(p1) + (2.5,0) + (0,5)$);
    \draw[dashed, green!10!blue!95!red!70] ($(p3) + (0,-1)$) -- (p3);
    \draw[dashed, purple] ($(p3) + (0,-1)$) --($(p3) + (0,-2.5)$);
    \draw[dashed] (p2) -- ($(p2) + (0,-1.5)$);
    \draw[dashed, green!10!blue!95!red!70] ($(p2) + (0,-1.5)$) -- ($(p2) + (0,-2)$);
    \draw[dashed] (p5) -- ++(90:5);

    \draw[dashed, purple] ($(p1) + (0.4,0)$) -- ++(270:1);

    \draw[dashed, green!10!blue!95!red!70] ($(p3) + (1,0)$) -- ++(270:1);
    \draw[dashed] ($(p1) + (1.5,0)$) -- ++(270:1);
    \draw[dashed, green!10!blue!95!red!70] ($(p1) + (1.5,0)$) -- ++(90:1.5);
    
    \draw[thick, green!10!blue!95!red!70] (p2) -- ($(p2) + (2.6,0)$);
\end{scope}

\begin{scope}
    \draw[>=latex, ->, thick, green!10!blue!95!red!70] (3,-1) -- (3.5,-1.5);
    \draw[>=latex, ->, thick, purple] (0.5,-1) -- (0,-1.5);
\end{scope}

\def\x{-3.5}; \def\y{-4};
\draw[thick, purple] (\x,\y) -- (\x, \y + 3.3) -- (\x + 1.5, \y + 3.3) -- (\x + 1.5, \y + 4.8) -- (\x + 2.5, \y + 4.8) -- (\x + 2.5, \y + 2.5) -- (\x + 2.9, \y + 2.5)-- (\x + 2.9, \y) -- cycle;
\draw[purple, dashed] (\x + 1.5, \y + 3.3) -- (\x + 1.5, \y);
\draw[purple, dashed] (\x + 2.5, \y + 2.5) -- (\x + 2.5, \y);

\def\x{4}; \def\y{-4}
\draw[thick, green!10!blue!95!red!70] (\x,\y) -- (\x, \y + 2.3) -- (\x + 1.1, \y + 2.3) -- (\x + 1.1, \y + 4.8) -- (\x + 2.1, \y + 4.8) -- (\x + 2.1, \y + 2.5) -- (\x + 2.6, \y + 2.5)-- (\x + 2.6, \y) -- cycle;
\draw[green!10!blue!95!red!70, dashed] (\x + 1.1, \y + 2.3) -- (\x + 1.1, \y);
\draw[green!10!blue!95!red!70, dashed] (\x + 2.1, \y + 2.5) -- (\x + 2.1, \y);
\end{tikzpicture}
    \caption{Two possible choices of zippered rectangle decompositions satisfying the given conditions.}
    \label{fig:hexagon_canonical_rectangle}
\end{figure}

\begin{example}[Nonunique choice of decomposition]
    Gluing a translation surface from a hexagon results in a surface of genus $\mathbf{g} = 1$ with two conical singularities, each of cone angle $2\pi$. Therefore, we can choose two distinct sections that satisfy the conditions outlined above. Figure \ref{fig:hexagon_canonical_rectangle} illustrates an instance of this.
\end{example}

\begin{remark}
    We can apply this kind of decomposition if the translation surface satisfies Keane's condition, so it is valid outside of a set of measure zero. Formally, this allows us to view the space of zippered rectangles essentially as a finite branched covering over the corresponding component of the stratum (\cite{zorich2007geodesics}).
\end{remark}

This (almost) canonical choice of coordinates provides us with a fundamental domain in the space of zippered rectangles. The Teichmüller geodesic flow $g_t$ naturally lifts to the space of zippered rectangles, where it acts in a particularly nice way. As we have seen, $g_t$ contracts a translation surface in vertical and expands in horizontal direction exactly by the factor $\e^{\pm t}$. In rectangle coordinates $(\boldsymbol{\pi}, \boldsymbol{\lambda}, \boldsymbol{h}, \boldsymbol{a})$, the action of $g_t$ is thus simply given by
\begin{equation}\label{eq:teichmüller_rectangle_explicit}
    g_t(\boldsymbol{\pi}, \boldsymbol{\lambda}, \boldsymbol{h}, \boldsymbol{a}) = (\boldsymbol{\pi}, \e^{t}\boldsymbol{\lambda}, \e^{-t} \boldsymbol{h}, \e^{-t}\boldsymbol{a}).
\end{equation}
Clearly, for $t$ large enough, $g_t(X)$ will no longer belong to the fundamental domain, i.e., the horizontal will be too long at some point and we could find a shorter horizontal section of length at least 1 satisfying the MIC. We want to have a procedure that, given some zippered rectangle decomposition not belonging to the fundamental domain because the horizontal segment is too long (or one that belongs to its boundary), produces the representative belonging to (the interior of) the fundamental domain. It is easy to see that this is achieved by applying the Rauzy--Veech induction step, which shortens the horizontal segment by $\lambda(\varepsilon)$ for $\varepsilon \in \{\mathrm{top}, \mathrm{bot}\}$ (depending on the type of the underlying IET) and if necessary we can repeat this step until the segment is the shortest possible segment satisfying the MIC. Let us develop this idea in a more formal way.

Instead of working in the space of (all) zippered rectangles, we start by passing to a subspace $\Omega$ of codimension 1 defined by the condition
\begin{equation*}
    \langle \boldsymbol{\lambda}, \boldsymbol{h}\rangle = 1,
\end{equation*}
where $\langle \cdot , \cdot \rangle$ denotes the usual inner product. This condition restricts us to consider only translation surfaces of total area 1. In the same sense as above, $\Omega$ is essentially a branched covering of the space $\mathcal{H}_1^\mathrm{comp}(d_1, \ldots, d_m)$, a component of the stratum restricted to unit area translation surfaces. 

Now we consider a further subspace $\Upsilon \subseteq \Omega$ of codimension 1 defined by
\begin{equation*}
    \|\boldsymbol{\lambda}\|_1 = \lambda_1 + \ldots + \lambda_d = 1,
\end{equation*}
that is, we consider zippered rectangles of area 1 whose base length is 1 as well. 

\begin{remark}
    The space $\Upsilon$ is naturally a fiber bundle over $\mathfrak{R}\times \Delta^{d-1}$, the space of IETs on $d$ intervals. The projection just \enquote{forgets} the height and altitude data.
\end{remark}

As we have mentioned above, we may view $\Omega$ as a fundamental domain in the space of all unit area zippered rectangles. Then, $\Upsilon$ is exactly the boundary of the fundamental domain. Applying the Teichmüller geodesic flow $g_t$ to a zippered rectangle $X \in \Upsilon\subseteq \Omega$, we will eventually leave $\Omega$ since the length of the base of $X$ will become too long, meaning that we can find a shorter possible base satisfying the MIC. In fact, we can give explicitly the exact time when this happens.  

\begin{proposition}[Hitting Time of Boundary]\label{prop:hitting_time_boundary}
    Let $X \in \Upsilon \subseteq \Omega$ be a unit area zippered rectangle with base length 1. Applying the Teichmüller geodesic flow $g_t$ to $X$, the first positive time it hits the boundary $\Upsilon$ of $\Omega$ again is given by
    \begin{equation*}
        t_0 = -\log\big(1-\min(\lambda_{\alpha(\mathrm{top})}, \lambda_{\alpha(\mathrm{bot})})\big).
    \end{equation*}
\end{proposition}
\begin{proof}
    Note that the expansion and contraction resulting from the Teichmüller geodesic flow $g_t$ commutes with the Rauzy--Veech induction step. Applying the move of the correct type, we obtain a zippered rectangle with base of length
    \begin{equation*}
        1 - \min(\lambda_{\alpha(\mathrm{top})}, \lambda_{\alpha(\mathrm{bot})}).
    \end{equation*}
    For $t_0$ as defined above, it follows from $\eqref{eq:teichmüller_rectangle_explicit}$ that the length of the base of $g_{t_0}(X)$ is given by
    \begin{equation*}
        \e^{t_0} (1 - \min(\lambda_{\alpha(\mathrm{top})}, \lambda_{\alpha(\mathrm{bot})})) = \frac{1}{1 - \min(\lambda_{\alpha(\mathrm{top})}, \lambda_{\alpha(\mathrm{bot})})} \cdot (1 - \min(\lambda_{\alpha(\mathrm{top})}, \lambda_{\alpha(\mathrm{bot})})) = 1. 
    \end{equation*}
\end{proof}
The way to think of what we have just described, is that we start at a point $X \in \Upsilon$, flow up to the time $t_0(X)$ defined in Proposition \ref{prop:hitting_time_boundary} and immediately apply the Rauzy--Veech induction step to obtain $g_{t_0}(X) \in \Upsilon$. Of course, this description is nothing but a first return map 
\begin{equation*}
    \mathcal{S} \colon \Upsilon \to \Upsilon
\end{equation*}
of the Teichmüller geodesic flow $g_t$ to the transverse section $\Upsilon$. Moreover, the action of $g_{t_0}$ corresponds exactly to the renormalization step of Rauzy--Veech induction, so we may in fact view the map $\mathcal{S} \colon \Upsilon \to \Upsilon$ as Rauzy--Veech induction performed at the level of zippered rectangles. Put differently, if we denote by $\iota$ the natural projection from the space $\Upsilon(\mathfrak{R})$, which is just a restriction of the space $\Upsilon$ to some Rauzy class $\mathfrak{R}$, to $\mathfrak{R} \times \Delta^{d-1}$, the space of IETs on $d$ intervals, then

\[\begin{tikzcd}
	{\Upsilon(\mathfrak{R})} && {\Upsilon(\mathfrak{R})} \\
	\\
	{\mathfrak{R}\times \Delta^{d-1}} && {\mathfrak{R}\times \Delta^{d-1}}
	\arrow["{\mathcal{S}}", from=1-1, to=1-3]
	\arrow["{\mathcal{T}}"', from=3-1, to=3-3]
	\arrow["\iota"', from=1-1, to=3-1]
	\arrow["\iota", from=1-3, to=3-3]
\end{tikzcd}\]

commutes. In particular, the invariant measure on $\mathfrak{R}\times \Delta^{d-1}$ is a push forward of the natural invariant measure on the space $\Upsilon$ of zippered rectangles, i.e., the restriction of the Lebesgue measure on the noncombinatorial data which is exactly the Masur--Veech measure from Definition \ref{def:masur_veech_measure}.

In the same way, one can find that the fast Rauzy--Veech induction $\mathcal{G} \colon \mathfrak{R}\times \Delta^{d-1} \to \mathfrak{R}\times \Delta^{d-1}$ corresponds to choosing a subsection $\Upsilon' \subseteq \Upsilon$ and consider the first return map of $\Upsilon'$. This subsection has a simple description in terms of the zippered rectangle coordinates $(\boldsymbol{\pi}, \boldsymbol{\lambda}, \boldsymbol{h}, \boldsymbol{a})$ given by an additional condition on the parameter $a_d$, that is the last altitude parameter. 

From the construction of the solution space in Theorem \ref{thm:suspensions_existence}, it follows that
\begin{equation*}
    a_j > 0 \quad \text{for all } 1\leq j \leq d-1,
\end{equation*}
but the last parameter $a_d$ may take both positive and negative values. If the rightmost rectangle $R_d$ contains a singularity on its right side, we have $a_d > 0$, where $a_d$ is just the height of the singularity. But if the right side of $R_d$ does not contain a singularity, we have that $a_d<0 $ which is then the distance to the next singularity \emph{below} the segment $S$. This fact allows us to define the transverse section $\Upsilon' \subseteq \Upsilon$ as follows.

\begin{definition}[Section for fast Rauzy--Veech on Zippered Rectangles]
    We define the subsection $\Upsilon' \subseteq \Upsilon$ as
    \begin{equation*}
        \Upsilon' \coloneqq \{(\boldsymbol{\pi}, \boldsymbol{\lambda}, \boldsymbol{h}, \boldsymbol{a}) \in \Upsilon \mid a_d > 0, \lambda_{\alpha(\mathrm{top})} > \lambda_{\alpha(\mathrm{bot})}\}\cup\{(\boldsymbol{\pi}, \boldsymbol{\lambda}, \boldsymbol{h}, \boldsymbol{a}) \in \Upsilon\mid a_d <0,  \lambda_{\alpha(\mathrm{top})} < \lambda_{\alpha(\mathrm{bot})}\}.
    \end{equation*}
    In words, the subsection $\Upsilon'$ consists of rectangle coordinates that either are of top type with $a_d > 0$ or of bot type with $a_d < 0$. 
\end{definition}
\begin{remark}
    The condition imposed on $\Upsilon$ can conveniently be written as
    \begin{equation*}
        a_d \cdot (\lambda_{\alpha(\mathrm{top})} - \lambda_{\alpha(\mathrm{bot})})  > 0. 
    \end{equation*}
\end{remark}

One checks that $\Upsilon'$ is still a fibre bundle over the space of IETs. Moreover, the corrsponding first return map $\mathcal{S}' \colon \Upsilon'  \to  \Upsilon'$ projects to the fast Rauzy--Veech induction $\mathcal{G} \colon \mathfrak{R} \times \Delta^{d-1} \to \mathfrak{R}\times \Delta^{d-1}$, i.e., the diagram

\[\begin{tikzcd}
	{\Upsilon'(\mathfrak{R})} && {\Upsilon'(\mathfrak{R})} \\
	\\
	{\mathfrak{R}\times \Delta^{d-1}} && {\mathfrak{R}\times \Delta^{d-1}}
	\arrow["{\mathcal{S}'}", from=1-1, to=1-3]
	\arrow["{\mathcal{G}}"', from=3-1, to=3-3]
	\arrow["\iota"', from=1-1, to=3-1]
	\arrow["\iota", from=1-3, to=3-3]
\end{tikzcd}\]

commutes. To see this, we notice that after a move of top type we always have $a_d < 0$, so we continue applying moves of top type until the first point where we need to apply a bot move. This corresponds exactly to the algorithmic description of the acceleration of the Rauzy--Veech induction considered on the space of IETs. 
 \pagebreak

\section{Diagonal Changes}\label{sec:diagonal_changes}
\thispagestyle{plain}

In section \ref{sec:rauzy_veech} we have introduced a discretization of the Teichmüller geodesic flow called the \emph{Rauzy--Veech Induction} and we have mentioned that this algorithm has been successfully employed to give answers to many questions about the dynamics of translation surfaces. In particular, we have seen that the associated zippered rectangle decomposition comes with coordinates on which the Teichmüller geodesic flow acts in a particularly nice way, namely it affects the different coordinates simply by scaling. 

However, for certain applications this decomposition, or more broadly the algorithm given by Rauzy--Veech induction, may not clearly encompass all the properties of the translation surfaces one is interested in. Suppose for instance one is interested in the \emph{systoles} of a family of translation surfaces (see Definition \ref{def:systole}). Given a zippered rectangle decomposition, there is no clear and easy way to identify the systole of the given translation surface. The question arises, whether another discretization might yield coordinates that, among other things, give a simple way to see the systole of any given translation surface. As we will see, the algorithm known as \emph{diagonal changes}, introduced first by Delecroix and Ulcigrai in \cite{delecroix2015diagonal}, fulfills this requirement. Diagonal changes exhibit many more nice properties that sets it apart from Rauzy--Veech induction. In particular, as we will discuss in section \ref{sec:counting}, diagonal changes may be used to count closed Teichmüller geodesics in certain components of strata. These types of results are very hard to obtain using Rauzy--Veech induction, mostly due to the fact that the branched cover given by the zippered rectangle decomposition discussed above tends to be too complicated to recover a closed geodesic in the space of translation surfaces from a closed geodesic in the space of zippered rectangles. 

\subsection{The simple case of the torus}

We again begin by describing the algorithm for the simple case of the torus, i.e., a translation surface of genus $\mathbf{g} = 1$. As was the case in section \ref{sec:rauzy_veech}, on the torus the algorithm is closely tied to the continued fraction expansion.

\begin{figure}[ht]
    \centering
    \begin{tikzpicture}
    \coordinate (xleft) at (-4,0); \coordinate (xright) at (4,0);
    \coordinate (ytop) at (0,5); \coordinate (ybot) at (0,-1);
    
    \draw[>=latex, ->, thick] (xleft) -- (xright);
    \draw[>=latex, ->, thick] (ybot) -- (ytop);

    \coordinate (B1) at ({1.6*pi/3}, {1});
    \coordinate (B2) at ({1/2}, {1*exp(1)});

    \clip (-4,0) rectangle (4,4.5);

    \foreach \a in {-4,...,4}{
        \foreach \b in {0,1, 2, 3, 4}{
            \draw[fill = black] ($\a*(B1) + \b*(B2)$) circle (1pt);
        }
    }

    \draw[>=latex, ->, purple] (0,0) -- (B1) node[midway, below]{$1$};
    \draw[>=latex, ->, purple] (0,0) -- ($(B2) - (B1)$) node[midway, left]{$2$};
    \draw[>=latex, ->, purple] (0,0) -- (B2) node[midway, right]{$3$};
    \draw[>=latex, ->, purple] (0,0) -- ($2*(B2)-(B1)$) node[midway, left]{$4$};
    
\end{tikzpicture}
    \caption{A geometric interpretation of the classical continued fraction algorithm.}
    \label{fig:continued_fraction_approximate_vertical}
\end{figure}
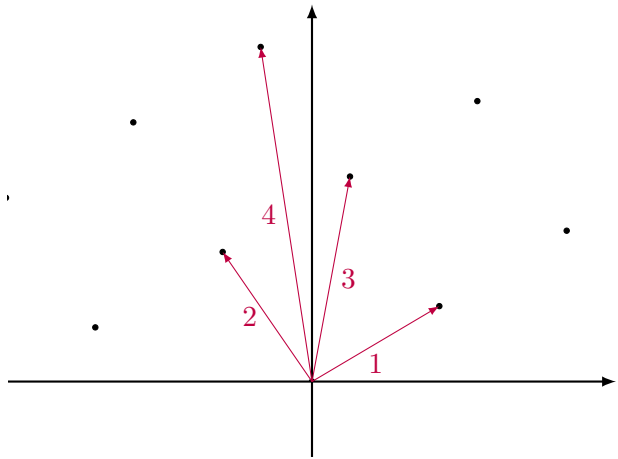

Let $\Lambda \subseteq \C$ be a lattice, i.e., a discrete additive subgroup of $\C$. The classical (additive) continued fraction algorithm given by the Farey map (Definition \ref{def:Farey_map}) can be seen as a way to construct a sequence of vectors in $\Lambda$ that approximate the vertical direction. This is illustrated in Figure \ref{fig:continued_fraction_approximate_vertical}.

\begin{figure}[ht]
    \centering
    \begin{tikzpicture}
    \draw[>=latex, ->, thick, noamblue] (0,0) -- (3,1) node[midway, below, yshift = -5]{$w_r$};
    \draw[>=latex, ->,thick, noamblue] (0,0) -- (-2,2) node[midway, left, yshift = -5] {$w_\ell$};
    \draw[thick] (-2,2) -- (1,3) -- (3,1);
    \draw[dashed] (0,0) -- ++(90:3.5);

    \def \o{7};
    \draw[>=latex, ->, thick, noamblue] (\o,0) -- (\o-3,1) node[midway, below, yshift=-5] {$w_\ell$};
    \draw[>=latex, ->, thick, noamblue] (\o,0) -- (\o+2,2) node[midway, right, yshift=-5] {$w_r$};
    \draw[thick] (\o-3,1) -- (\o-1,3) -- (\o+2,2);
    \draw[dashed] (\o,0) -- ++(90:3.5);
\end{tikzpicture}
    \caption{The two different types of parallelograms induced by the wedge $\{w_\ell, w_r\}$. On the left, we have a left-slanted paralellogram, on the right a right-slanted parallelogram.}
    \label{fig:torus_wedges}
\end{figure}
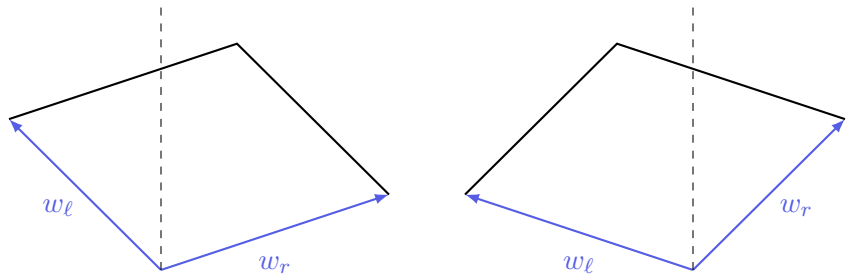

Diagonal changes, when performed on a translation surface of genus $\mathbf{g} = 1$, can be seen as the same algorithm expressed in geometric terms. Let us provide a description of this fact. We start by choosing a basis $\{w_\ell, w_r\}$ of $\Lambda$ such that
\begin{align}\label{eq:parallelogram_normal_form}
    \begin{split}
    \RE(w_\ell) < 0 \quad &\text{and} \quad \RE(w_r) > 0, \\ \IM(w_\ell) > 0 \quad &\text{and} \quad \IM(w_r) > 0.
    \end{split}
\end{align}
This is clearly possible as long as $\Lambda$ does not contain any horizontal or vertical (nonzero) vectors, which we shall assume unless otherwise stated. The basis $\{w_\ell, w_r\}$ forms what we will call a \emph{wedge} which contains (part of) the vertical ray emitted from the origin as pictured in Figure \ref{fig:torus_wedges}. Phrased differently, the vertical direction is contained in the positive cone generated by the basis. The wedge $\{w_\ell, w_r\}$ induces a quadrilateral $Q = Q(w_\ell, w_r)$ which in this case is a parallelogram. Note that $Q$ is a fundamental domain for the action of $\Lambda$ on $\C$. 

We will call the quadrilateral $Q$ \emph{left-slanted}, if the horizontal ray emitted from the origin, that is the set $\{z \in \C \mid \RE(z) = 0, \IM(z) \geq 0\}$, crosses the left top side which is parallel to $w_r$. Analogously, we call the quadrilateral $Q$ \emph{right-slanted} if the vertical ray crosses the top right side parallel to $w_\ell$.

\begin{figure}[ht]
    \centering
    \begin{tikzpicture}[scale = 0.75]
    \draw[thick] (0,0) -- (-3,1) -- (-1,3) -- (2,2) -- cycle;
    \draw[dashed](0,0) -- (-1,3);
    \fill[teal, opacity = 0.15] (0,0) -- (-3,1) -- (-1,3) -- cycle;
    \node[rotate =110, xshift = 8] at (-1,3) {\Large \Leftscissors};

    \draw[>=latex, ->] (2,1.5) -- (3,1.5);

    \def \x1{6.1};

    \draw[thick] (\x1,0) -- (\x1 - 3,1) -- (\x1 - 1,3); \draw[dashed] (\x1 - 1,3) -- (\x1,0);
    \draw[thick] (\x1 + 0.3,0) -- (\x1 + 2.3,2) -- (\x1 - 0.7,3);\draw[dashed] (\x1+0.3,0) -- (\x1-0.7,3);
    \fill[teal, opacity = 0.15] (\x1,0) -- (\x1-3,1) -- (\x1-1,3) -- cycle;

    \draw[>=latex, ->] (8.4,1.5) -- (9.4,1.5);

    \def \x2{10.1};
    \draw[thick] (\x2 + 0.3,0) -- (\x2 + 2.3,2) -- (\x2 - 0.7,3);
    \draw[dashed] (\x2+0.3,0) -- (\x2-0.7,3);
    \draw[thick] (\x2 + 2.3,2.2) -- (\x2 -0.7,3.2) -- (\x2 + 1.3,5.2);
    \draw[dashed] (\x2 + 1.3,5.2) -- (\x2 + 2.3, 2.2);
    \fill[teal, opacity = 0.15] (\x2 + 2.3,2.2) -- (\x2 -0.7,3.2) -- (\x2 + 1.3,5.2) -- cycle;

    \draw[>=latex, ->] (12.8,1.5) -- (13.8,1.5);

    \def \x3{14.5};
    \draw[thick] (\x3,0) -- (\x3 + 2.3,2.2) -- (\x3 - 0.7, 3.2) -- (\x3 + 1.3, 5.2);
    \draw[dashed] (\x3 + 1.3, 5.2) -- (\x3 + 2.3,2.2);
    \draw[dashed] (\x3, 0) -- (\x3 - 0.7, 3.2);
    \fill[teal, opacity = 0.15] (\x3 + 2.3,2.2) -- (\x3 - 0.7, 3.2) -- (\x3 + 1.3, 5.2) -- cycle;
\end{tikzpicture}
    \caption{A right move applied to a right-slanted quadrilateral $Q$.}
    \label{fig:diagonal_change_torus_right}
\end{figure}
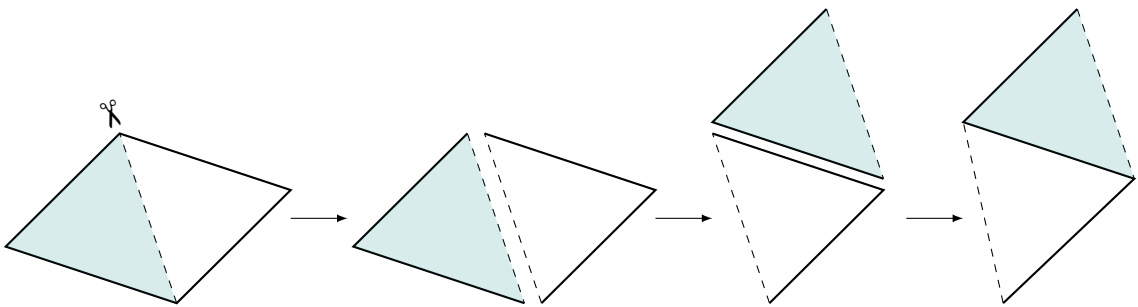

Let us now describe one step of the algorithm. Let us suppose the quadrilateral $Q$ is right-slanted, the left-slanted case works analogously. We obtain a new quadrilateral $Q'$ from the basis $\{w'_\ell = w_\ell + w_r, w'_r = w_r\}$. Geometrically, we cut the old quadrilateral $Q$ along the diagonal coming from the origin and paste the triangle on the left to the top of the triangle on the right-hand side. This geometric interpretation can be seen in Figure \ref{fig:diagonal_change_torus_right}. We call this procedure a \emph{right move}. 

\begin{remark}
    The vertical direction is contained in the positive cone generated by the new basis associated to $Q'$ as well, independent of whether we apply a left or a right move.
\end{remark}

The full algorithm is obtained by successive application of the moves we have just described. Formally, we can write $w_\ell^{(0)} = w_\ell$ and $w_r^{(0)} = w_r$, so that consecutive applications of the moves give a sequence $\left(w_\ell^{(n)}, w_r^{(n)}\right)_{n \in \N}$ of bases of $\Lambda$.

Note that the cut and paste procedure described above is invertible, and the inverse is given by cutting along the diagonal of $Q$ that does not contain the origin and gluing the top triangle to the bottom triangle in the unique way such that the newly obtained quadrilateral $Q'(w_\ell, w_r)$ satisfies the constraints in \eqref{eq:parallelogram_normal_form}. Therefore, we can in fact obtain a bi-infinite sequence $\left(w_\ell^{(n)}, w_r^{(n)}\right)_{n \in \Z}$ of bases of $\Lambda$. The following lemma justifies all these assertions.

\begin{lemma}[Geometric Continued Fraction Algorithm]\label{{lem:geometric_cf_algorithm}}
    The algorithm given by successive applications of right and left moves, which we will call \emph{geometric continued fraction algorithm}, is defined for all $n \in \Z$. Moreover, 
    \begin{equation*}
        \IM\left(w_\varepsilon^{(n)}\right) \xrightarrow{n \to +\infty} +\infty \quad \text{for } \varepsilon \in \{\ell, r\}.
    \end{equation*}
\end{lemma}
\begin{proof}
    The algorithm can only terminate in finite time in the forward direction, if for some $k \in \N$ we have $\RE(w_\ell) + \RE(w_r) = 0$, i.e., if the diagonal containing the origin is exactly vertical so that it is not evident whether a left or a right move is to be applied. But in this case, this vertical diagonal $d^{(k)}$ is written as
    \begin{equation*}
        d^{(k)} = w_\ell^{k} + w_r^{k} = a_1 w_\ell^{(0)} + a_2 w_r^{(0)}
    \end{equation*}
    for some $a_1, a_2 \in \N$ by the definition of the algorithm. But since we assume that the lattice $\Lambda$ does not contain any vertical vectors, this is impossible. Note that this proof also shows the converse statement, i.e., the algorithm terminates if and only if the lattice $\Lambda$ contains a vertical vector. The case for the backwards direction is obtained analogously by replacing verticals with horizontals. 

    For the second statement, it is enough to show that we apply both left and right moves infinitely many times. Indeed, if this is the case, then we can write $c_\varepsilon^{(n)}$ for the function counting the occurrences for left or right moves, so that
    \begin{equation*}
        \IM\left(w_\varepsilon^{(n)}\right) = c_\varepsilon^{(n)} \cdot \IM\left(w_\varepsilon^{(0)}\right) \xrightarrow{n \to \infty} + \infty,
    \end{equation*}
    since $c_\varepsilon^{n} \to +\infty$. To see that this is the case, suppose towards a contradiction that eventually we only apply, say, left moves. In particular, we have
    \begin{equation*}
        \RE(w_r^{(n+1)}) = \RE(w_\ell^{(n)}) + \RE(w_r^{(n)}),
    \end{equation*}
    where we recall that $\RE(w_\ell^{(n)}) < 0$. Inductively, we then obtain
    \begin{equation*}
        \RE(w_r^{(m)}) = (m-n_0)\cdot \RE(w_\ell^{(n_0)}) + \RE(w_r^{n_0)}),
    \end{equation*}
    which would be negative for $m$ large enough. This contradicts the fact that applying the appropriate move ensures that the vertical ray from the origin is always contained in the positive cone generated by the basis elements at any step. Phrased differently, this means that at some point we should have applied a right move. We conclude that both moves appear infinitely often, which finishes the proof.
\end{proof}

Let us briefly explain why the name \emph{geometric continued fraction algorithm} is the appropriate name for this procedure. Let $\Gamma \subseteq \Lambda$ be the set of \emph{primitve} vectors of $\Lambda$, that is, the set of vectors that are not a nontrivial integer multiple of another vector in $\Lambda$. We can decompose $\Gamma$ into
\begin{equation*}
    \Gamma = \Gamma_\ell \sqcup \Gamma_r,
\end{equation*}
where we set $\Gamma_\ell = \{w \in \Gamma \mid \RE(w) < 0\}$ and $\Gamma_r = \{w \in \Gamma \mid \RE(w) > 0\}$. 

\begin{lemma}
    For any $n \in \N$, we have
    \begin{equation*}
        w_\ell^{(n)} \in \Gamma_\ell, \quad w_r^{(n)} \in \Gamma_r. 
    \end{equation*}
\end{lemma}
\begin{proof}
    It is clear from the construction that the condition on the real parts is satisfied, so it suffices to show that $w_\varepsilon^{(n)}$ is a primitive vector for all $n \in \N$. Suppose this is not the case, so that we can write
    \begin{equation*}
        w_\varepsilon^{(n+k)} = a\cdot w_\varepsilon^{(n)}
    \end{equation*}
    for some $k \geq 1$ and some $a \in \Z\setminus\{0, \pm 1\}$. Since the imaginary parts of all wedge vectors is positive, we must have $a \in \N_{>1}$. But applying the algorithm in forward time can only decrease the real part of the wedge vectors in modulus, thus
    \begin{equation*}
        a \cdot |\RE(w_\varepsilon^{(n)})| = |\RE(w_\varepsilon^{(n+k)})| \leq |\RE(w_\varepsilon^{(n)})|, 
    \end{equation*}
    giving a contradiction since $a \geq 2$ and $\RE(w_\varepsilon^{(n)}) \neq 0$. 
\end{proof}

\begin{definition}[Geometric Best Approximation]\label{def:geometric_best_approximation}
    A vector $w \in \Gamma_\ell$ is a \emph{left geometric best approximation}, if
    \begin{equation*}
        \IM(v) < \IM(w) \quad \Rightarrow \quad |\RE(v)| > |\RE(w)| \quad \text{for all } v \in \Gamma_\ell.
    \end{equation*}
    Replacing $\Gamma_\ell$ by $\Gamma_r$, we obtain the definition of \emph{right geometric best approximation}. 
\end{definition}

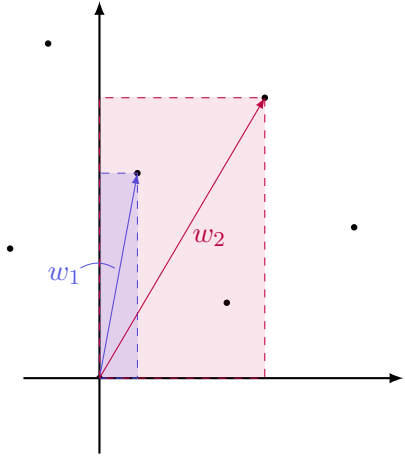
\begin{figure}[ht]
    \centering
    \begin{tikzpicture}
    \coordinate (xleft) at (-1,0); \coordinate (xright) at (4,0);
    \coordinate (ytop) at (0,5); \coordinate (ybot) at (0,-1);
    
    \draw[>=latex, ->, thick] (xleft) -- (xright);
    \draw[>=latex, ->, thick] (ybot) -- (ytop);

    \coordinate (B1) at ({1.6*pi/3}, {1});
    \coordinate (B2) at ({1/2}, {1*exp(1)});

    \clip (-1.5,0) rectangle (4,4.5);

    \foreach \a in {-4,...,4}{
        \foreach \b in {0,1, 2, 3, 4}{
            \draw[fill = black] ($\a*(B1) + \b*(B2)$) circle (1pt);
        }
    }

    \draw[>=latex, ->, noamblue] (0,0) -- (B2);
    \draw[dashed, noamblue] (0,0) rectangle (B2) node[midway, left, xshift = -10] {$w_1$};
    \draw[noamblue] (-1/4,e/2+0.1) to[out=30,in =150] (1/5,e/2+0.1);
    \fill[blue, opacity = 0.1] (0,0) rectangle (B2);

    \draw[>=latex, ->, purple] (0,0)-- ($(B1) + (B2)$) node[midway, right] {$w_2$};
    \draw[dashed, purple] (0,0) rectangle ($(B1) + (B2)$);
    \fill[purple, opacity =0.1] (0,0)rectangle ($(B1) + (B2)$);
    
 \end{tikzpicture}
    \caption{The vector $w_1$ is a right geometric best approximation, since the blue rectangle does not contain any lattice points in its interior. In contrast, $w_2$ is \emph{not} a right geometric best approximations since the red rectangle contains two lattice points in its interior.}
    \label{fig:geometric_best_approximation}
\end{figure}

In words, Definition \ref{def:geometric_best_approximation}
says that if $w$ is a best approximation and $v$ is a complex number closer to the real axis, then it is necessarily the case that the real part of $v$ is larger than the real part of $w$. We can state this also in purely geometric terms, saying that $w$ is a left geometric best approximation if and only if the rectangle $R(w) \coloneqq [\RE(w), 0]\times[0, \IM(w)]$ contains no other points of $\Lambda$ in its interior. This is illustrated in Figure \ref{fig:geometric_best_approximation}.

The geometric continued fraction algorithm generates geometric best approximations in a sense made precise in the following theorem.

\begin{theorem}[\cite{delecroix2015diagonal} Theorem 1]\label{thm:DU15_thm1}
    Let $\Lambda$ be a lattice in $\C$ that contains neither vertical nor horizontal vectors. Then, the sequence of bases 
    \begin{equation*}
        \left(w_\ell^{(n)}, w_r^{(n)}\right)_{n \in \Z}
    \end{equation*}
    built from the geometric continued fraction algorithm is uniquely defined up to a shift in the numbering. Moreover, the vectors $w_\ell^{(n)}$ and $w_r^{(n)}$ are exactly the geometric best approximations.
\end{theorem}

Theorem \ref{thm:DU15_thm1} can be stated and proven fully in terms of Diophantine approximations and corresponds exactly to the fact, that intermediate convergents of a real number $\alpha$ are exactly the approximations of the first kind, see e.g., \cite{khinchin1963continued}. In particular, the sequence of bases of $\Lambda$ contains all wedges that encode the (non-intermediate) convergents of $\alpha$, that is, the partial quotients of the continued fraction expansion of $\alpha$.

Let us give a proof based on the geometric language developed above. The proof here is adapted from the proof of Theorem 10 in \cite{delecroix2015diagonal}, which is the analogous statement for the generalization to structures related to surfaces of higher genus which we will also develop below.

\begin{proof}[Proof of Theorem \ref{thm:DU15_thm1}]
    We start by showing that any saddle connection belonging to a wedge of the sequence generated by the algorithm is a geometric best approximation. So let $(w_\ell^{(n)}, w_r^{(n)})$ be any such wedge. We will show that $w_\ell^{(n)}$ is a left geometric best approximation, the same argument adapts readily to the right part of the wedge. 

    Let $R\left(w_\ell^{(n)}\right) = [\RE(w_\ell^{(n)}), 0] \times [0, \IM(w_\ell^{(n)})]$ be the rectangle with horizontal and vertical sides that has $w_\ell^{(n)}$ as its diagonal. We claim that $R\left(w_\ell^{(n)}\right)$ has no lattice points in its interior. Indeed, if $(w_\ell^{(n)}, w_r^{(n)})$ induces a right-slanted parallelogram, the rectangle is contained in the strip of the lattice as indicated in Figure \ref{fig:lattice_strip}. In case of a left-slanted parallelogram, we obtain a comparable picture.
    \begin{figure}[ht]
        \centering
        \begin{tikzpicture}
    \coordinate (xleft) at (-4,0); \coordinate (xright) at (2,0);
    \coordinate (ytop) at (0,5); \coordinate (ybot) at (0,-1);
    
    \draw[>=latex, ->, thick] (xleft) -- (xright);
    \draw[>=latex, ->, thick] (ybot) -- (ytop);

    \draw[purple] (0,0) -- (-2,2.5);
    \draw[purple, dashed] (-2,0) -- (-2,2.5) -- (0,2.5);
    \fill[purple, opacity =0.15] (0,0) rectangle (-2,2.5);

    \coordinate (p) at (0,0);
    \draw (p) -- ([xshift = 17, yshift = 27]p) -- ($([xshift = 17, yshift = 27]p) + (-2,2.5)$) -- ($(p) + (-2,2.5)$);

    \coordinate (p) at ([xshift=17,yshift=27]0,0);
    \draw (p) -- ([xshift = 17, yshift = 27]p) -- ($([xshift = 17, yshift = 27]p) + (-2,2.5)$) -- ($(p) + (-2,2.5)$);

    \coordinate (p) at ([xshift=34,yshift=54]0,0);
    \draw[dashed] (p) -- ([xshift = 17, yshift = 27]p);
    \draw[dashed] ($([xshift = 17, yshift = 27]p) + (-2,2.5)$) -- ($(p) + (-2,2.5)$);
    
    \coordinate (p) at ([xshift=-17,yshift=-27]0,0); 
    \draw ([xshift = 17, yshift = 27]p) -- (p)  -- ($(p) + (-2,2.5)$) -- ($([xshift = 17, yshift = 27]p) + (-2,2.5)$);

    \coordinate (p) at ([xshift=-34,yshift=-54]0,0); 
    \draw ([xshift = 17, yshift = 27]p) -- (p)  -- ($(p) + (-2,2.5)$) -- ($([xshift = 17, yshift = 27]p) + (-2,2.5)$);

    \coordinate (p) at ([xshift=-51,yshift=-81]0,0); 
    \draw[dashed] ([xshift = 17, yshift = 27]p) -- (p);
    \draw[dashed] ($(p) + (-2,2.5)$) -- ($([xshift = 17, yshift = 27]p) + (-2,2.5)$);
 \end{tikzpicture}
        \caption{The rectangle $R\left(w_\ell^{(n)}\right)$ is contained in the strip induced by a subset of $\Lambda$.}
        \label{fig:lattice_strip}
    \end{figure}
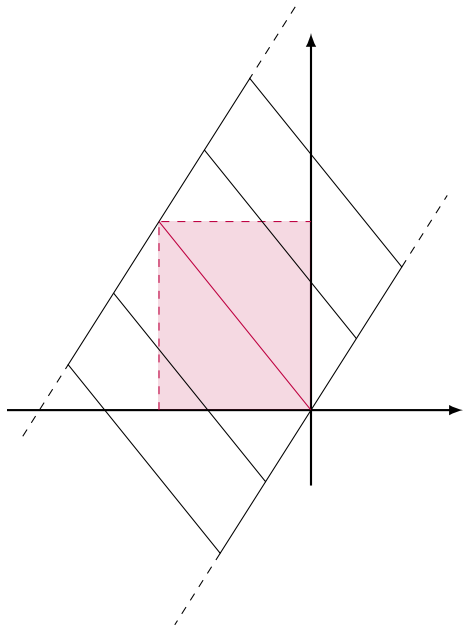

    Notice that any parallelogram induced by any of the wedges in the sequence cannot contain lattice points in their interior. This follows immediately from the definition of the algorithm and the fact that the initial choice of wedges form a basis. This shows that $R\left(w_\ell^{(n)}\right)$ cannot contain lattice points in its interior either, thus $w_\ell^{(n)}$ is a left geometric best approximation.

    It remains to show that all geometric best approximations appear in the sequence in their natural order. We have shown that the sequence $(w_\ell^{(n)})_{n \in \Z}$ consists entirely of best approximations which by the definition of the algorithm are naturally ordered by increasing imaginary part. Thus it suffices to show that if $w_\ell^{(n)}$ and $w_\ell^{n+1}$ are two successive left basis vectors according to this order, there is no geometric best approximation with imaginary part strictly in between the imaginary parts of these two vectors.

    For this, it is enough to verify that the rectangle $R = [\RE(w_\ell^{(n)}),0] \times [0, \IM(w_\ell^{(n+1)})]$ satisfies
    \begin{equation*}
        \Gamma_\ell \cap R = \{w_\ell^{(n)}, w_\ell^{(n+1)}\},
    \end{equation*}
    since if this is the case, then there are no primitive vectors $v \in \Gamma^\ell$ with $\IM w_\ell^{(n)} < \IM v < \IM w_\ell^{(n+1)}$ and $0 < |\RE v| \leq |\RE w_\ell^{(n)}|$. So the only possibility for primitive vectors $v$ satisfying $\IM w_\ell^{(n)} < \IM v < \IM w_\ell^{(n+1)}$ is to also satisfy $|\RE v | > |\RE w_\ell^{(n)}|$, but then it is not a best approximation. 
    
    To see that $\Gamma_\ell \cap R$ does not contain anything else, notice that by construction $w_\ell^{(n+1)}$ is the diagonal $d^{(n)}$ of the parallelogram spanned by $w_\ell^{(n)}$, so that the line segment $w$ joining the two points $w_\ell^{(n)}$ and $w_\ell^{(n+1)}$ is a translation of $w_r^{(n)}$, see Figure \ref{fig:best_approximation_thm}. Now, the rectangle $R$ is given by the (intersecting) union of rectangles $R_1, R_2$ and $R_3$ which have as diagonals $w_\ell^{(n)}, w_\ell^{(n+1)}$ and $w$ respectively. 
    \begin{figure}[ht]
        \centering
        \begin{tikzpicture} [scale = 0.75]
     \draw[>=latex, ->, noamblue] (0,0) -- (2,6) node[midway, right, xshift= -2] {$w_r^{(n)}$};
     \draw[>=latex, ->] (0,0) -- (-4,2) node[midway, below, xshift = -3] {$w_\ell^{(n)}$};
     \draw[>=latex, ->] (0,0) -- (-2,8) node[midway, above, yshift= 50, xshift=3] {$w_\ell^{(n+1)}$};
     \draw[>=latex, ->, noamblue] (-4,2) -- (-2,8) node[midway, left] {$w$};

     \draw[dashed, purple] (0,0) rectangle (-4,8);
     \draw[dashed, noamblue] (0,0) rectangle (2,6);

     \fill[teal, opacity =0.1] (0,0) rectangle (-2,8);
     \fill[red, opacity = 0.1] (0,0) rectangle (-4,2);
     \fill[blue, opacity = 0.1] (-2,2) rectangle (-4,8);

     \draw[fill = purple, color = purple] (1,5) circle (1pt);
     \draw[fill = purple, color = purple] ($(1,5) + (-4,2)$) circle (1pt);
 \end{tikzpicture}
        \caption{The rectangles used in the proof of Theorem \ref{thm:DU15_thm1}}
        \label{fig:best_approximation_thm}
    \end{figure}
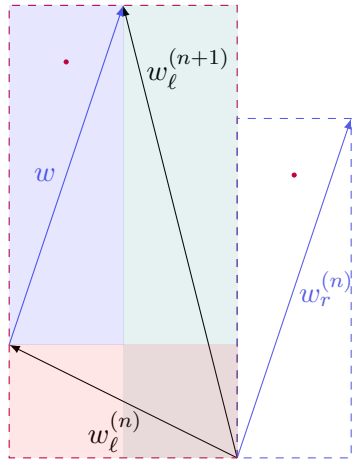
    
    By the first part of this proof, $R_1$ and $R_2$ cannot contain any elements of $\Gamma_\ell$ in their interior, so any element $v \in \Gamma_\ell$ in the interior of $R$ must lie in $R_3$. If by contradiction there was such a point, we would find another lattice point by translating it by $-w_\ell^{(n)}$, so that this translated point now lies in the interior of the rectangle $[0, \RE w_r^{(n)}] \times [0, \IM w_r^{(n)}]$. This contradicts the fact that $w_r^{n}$ is a right geometric best approximation, which we know from the first part of the theorem. 
\end{proof}

Note that the quotient $\T_\Lambda = \C/\Lambda$ is a flat torus where the origin is a marked point, so that all of the above can be interpreted in terms of translation surfaces as well. Moreover, the wedges $w_\varepsilon^{(n)}$ correspond exactly to saddle connections, see Definition \ref{def:saddle_connection}. The algorithm, therefore, generates saddle connections that progressively offer refined approximations of the vertical linear flow.

\subsection{Extend the idea to higher genera}
The natural way to extend the above idea is to move from the flat torus with a marked point to general translation surfaces. We have given most of the relevant definitions in section \ref{sec:introduction}, but it is convenient to introduce some more concepts adapted particularly to the development of the diagonal changes algorithm. 

Throughout this section we denote by $X$ a translation surface with singularity set $\Sigma$ and write $\varphi_t^\theta$ for the linear flow in direction $\theta$. 

\begin{definition}[Holonomy vector]\label{def:holonomy_vector}
    Given a saddle connection $(\varphi_t^\theta (\mathbf{x}))_{t \in I}$ with $\mathbf{x} \in X \setminus \Sigma$ and $I = [a,b] \subseteq \R$, meaning in particular that $\varphi_a^\theta(\mathbf{x}), \varphi_b^\theta(\mathbf{x}) \in \Sigma$, its \emph{holonomy} (or \emph{discplacement}) \emph{vector} is
    \begin{equation*}
        (b-a)\e^{\ii \theta} \in \C. 
    \end{equation*}
\end{definition}

We want to remark that the definition of holonomy vector is a reformulation of the concept of period coordinates introduced in Definition \ref{def:period_coordinates}. Here we can directly associate a complex number to any saddle connection, not necessarily an edge of a polygon. Cutting along such a saddle connection as in Definition \ref{def:translation_surface_equivalence} shows that the two concepts are the same.

\begin{figure}[ht]
    \centering
    \begin{tikzpicture}[scale = 0.75]
     \draw[thick] (0,0) rectangle (-3,3);
     \node at (0,1.5) {\Large $\sim$}; \node at (-3,1.5) {\Large $\sim$};
    \node at (-2,0) {//}; \node at (-2,3) {//};

    \draw[thick, noamblue] (-3,0) -- (-1.5,3);
    \draw[thick, noamblue](-1.5,0) -- (0,3);

    \coordinate (xleft) at (2,-1); \coordinate (xright) at (7,-1);
    \coordinate (ytop) at (3,6); \coordinate (ybot) at (3,-2);
    
    \draw[>=latex, ->, thick] (xleft) -- (xright);
    \draw[>=latex, ->, thick] (ybot) -- (ytop);     

    \draw[>=latex, ->, thick, noamblue] (3,-1) -- (6,5);

    \draw[>=latex, ->, thick] (1,1.5) -- (2,1.5);
 \end{tikzpicture}
    \caption{A simple example of a saddle connection (on the left) and its associated holonomy vector (on the right).}
    \label{fig:holonomy_vector}
\end{figure}
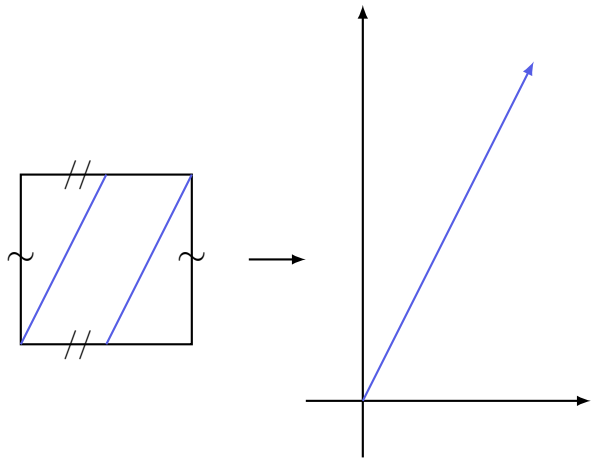

\begin{example}[Holonomy vector]\label{ex:holonomy_vector}
    Consider the standard torus $\T^2 = \R^2 / \Z^2$, with $\mathbf{x} = (1/2, 1), a=-\frac{\sqrt{5}}{2}, b=\frac{\sqrt{5}}{2}$ and $\theta = \frac{3\pi}{8}$. Under the usual assumption that the linear flow is parametrized by arclength, $(\varphi_t^\theta(\mathbf{x}))_{t \in I}$ is exactly the saddle connection depicted in Figure \ref{fig:holonomy_vector}. Its holonomy vector, also depicted in Figure \ref{fig:holonomy_vector}, is given by
    \begin{equation*}
        \sqrt{5}\e^{\ii \frac{3\pi}{8}}.
    \end{equation*}
\end{example}

The following fact is important and implicit in all discussion that concern counting of saddle connections. A sketch of a proof can be found both in \cite{masur2006ergodic} or \cite{hubert2005affine}.

\begin{theorem}\label{thm:holonomy_discrete}
    For any translation surface $X$, the set of displacement vectors of saddle connections on $X$ is a discrete subset of $\C$.
\end{theorem}

Every saddle connection comes with a choice of orientation and we will call the following choice \emph{natural.}

\begin{definition}[Natural orientation of saddle connections]
    We call \emph{natural orientation} of a saddle connection $\gamma$ on a translation surface $X$ the unique orientation of $\gamma$ such that its holonomy vector has nonnegative imaginary part. Furthermore, we say that $\gamma$ starts (or ends) at $p \in \Sigma$ if $p$ is the first (or last) endpoint of the saddle connection with respect to its natural orientation. 

    Lastly, we will say that a saddle connection is \emph{left-slanted} (or \emph{right-slanted}) if the real part of its associated holonomy vector with respect to the natural orientation is negative (or positive, respectively). 
\end{definition}

Let us write $\Gamma = \Gamma(X)$ for the set of all saddle connections on a translation surface $X$ and let $\Gamma_\ell$ be the subset of all left-slanted saddle connections and analogously let $\Gamma_r$ be the subset of all right-slanted saddle connections. We can subdivide these sets of saddle connections further, where we will use the same convention as in \cite{delecroix2015diagonal} which originally was introduced by Marchese in \cite{marchese2012khinchin}. Suppose that $\Sigma = \{p_1, \ldots, p_n\}$ with associated cone angles $2\pi k_1, \ldots, 2\pi k_n$. If the cone angle at $p_i \in \Sigma$ is $2\pi k_i$, then there are $k_i$ different outgoing trajectories of the vertical linear flow from $p_i$ and the same number of outgoing trajectories of the horizontal linear flow. This follows from the fact that $p_i$ has a neighborhood that is isomorphic to $2k_i$ half-disks glued together. For each $p_i \in \Sigma$, we choose a reference horizontal half ray $v_i$ starting from $p_i$ and for any two linear trajectories $\gamma, \gamma'$ starting at $p_i$ we write
\begin{equation*}
    \sphericalangle(\gamma, \gamma') \in [0, 2\pi k_i)
\end{equation*}
for the angle between $\gamma$ and $\gamma'$. Each saddle connection $\gamma$ originating at $p_i$ belongs to one of the $k_i$ \emph{outgoing half-planes}, i.e., the angle $\sphericalangle(\gamma,v_i)$ satisfies
\begin{equation*}
    2\pi j \leq \sphericalangle(\gamma, v_i) < 2\pi j + \pi,
\end{equation*}
for a unique $0\leq j < k_i$. We say that two saddle connections \emph{belong to the same bundle} if they start from the same singularity $p_i$ and belong to the same half-plane. Note that there are $k$ bundles of saddle connections on $X$, where $k = k_1 + \ldots + k_n$ is the total cone angle. We label these \emph{bundles of saddle connections} by the integers $1, \ldots, k$ and denote them by $\Gamma_1, \ldots, \Gamma_k$.

In the case of the torus above, the algorithm produced a sequence of pairs of saddle connections which which are increasingly good approximations of the vertical direction. When extending this idea to surfaces of higher genus, the algorithm will produce a sequence of collections of $k$ pairs of saddle connections, one pair for each of the $k$ vertical rays in $X$ starting at a singularity. We will call such a pair of saddle connections a \emph{wedge}. 

\begin{definition}[Wedge]\label{def:wedge}
    A \emph{wedge} $w$ on a translation surface $X$ is a pair of saddle connections $w = (w_\ell, w_r)$ such that the following properties hold.
    \begin{enumerate}
        \item Both $w_\ell$ and $w_r$ start from the same conical singularity of $X$.
        \item $w_\ell$ is left-slanted and $w_r$ is right-slanted.
        \item The pair $(w_\ell, w_r)$  consists of two edges of an embedded triangle in $X$. 
    \end{enumerate}
\end{definition}
\begin{remark}
    The first and third condition in Definition \ref{def:wedge} together are equivalent to asking that the saddle conections $w_\ell$ and $w_r$ that form the wedge $w$ belong to the same bundle of saddle connections $\Gamma_i$.
\end{remark}

An important property that follows from the definition is that any wedge contains a unique vertical trajectory, that is, there is exactly one trajectory of the vertical flow starting from the conical singularity shared by $w_\ell$ and $w_r$ that intersects the interior of the triangle with edges $w_\ell$ and $w_r$. We have already seen two simple examples of wedges above in Figure \ref{fig:torus_wedges}.

When working with the torus it was clear how to obtain a decomposition into a parallelogram which is completely described by a wedge. We want to obtain an analogous decomposition also in the case of translation surfaces of higher genus, for which we introduce the notions of \emph{quadrilaterals} and \emph{quadrangulations}.

\begin{definition}[Quadrilateral]
    A \emph{quadrilateral} $q$ in a translation surface $X$ is the image of an isometrically embedded quadrilateral in $\C$ so that the vertices of $q$ are singularities of $X$ and there are no other singularities of $X$ in $q$, neither in the interior nor on the edges. 
\end{definition}

Since we want to work with wedges, the following restriction on quadrilaterals is very natural.

\begin{definition}[Admissible quadrilateral]
    A quadrilateral $q$ in $X$ is \emph{admissible}, if there is exactly one trajectory of the vertical linear flow of $X$ starting from one of its vertices and exactly one ending in a vertex. An equivalent way of saying this, is that left-slanted and right-slanted saddle connections alternate when going around the quadrilateral. 
\end{definition}

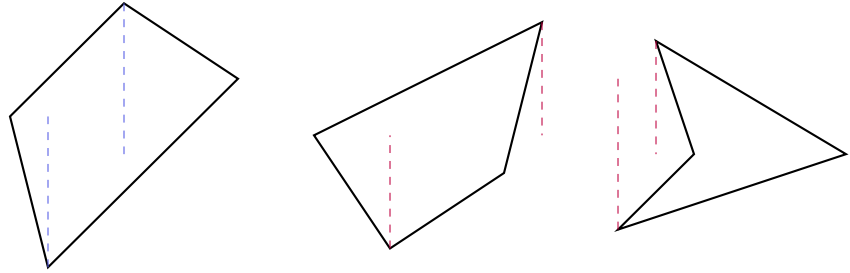
\begin{figure}[ht]
    \centering
    \begin{tikzpicture}[scale = 0.5]
     \draw[thick] (0,0) -- (-1,4) -- (2,7) -- (5,5) -- cycle;
     \draw[thick] (9,1/2) -- (7,3.5) -- (13,6.5) -- (12,2.5) -- cycle;
     \draw[thick] (15,1) -- (17,3) -- (16,6) -- (21,3) -- cycle;

    \draw[dashed, noamblue] (0,0) -- ++(90:4);
    \draw[dashed, noamblue] (2,7) -- ++(-90:4);

    \draw[dashed, purple] (9,1/2) -- ++(90:3);
    \draw[dashed, purple] (13,6.5)-- ++(-90:3);

    \draw[dashed, purple] (15,1) --++(90:4);
    \draw[dashed, purple] (16,6) -- ++(-90:3);
 \end{tikzpicture}
    \caption{An admissible quadrilateral (left) and two non-admissible quadrilaterals (middle and right).}
    \label{fig:quadrilateral_admissible_nonadmissible}
\end{figure}

Figure \ref{fig:quadrilateral_admissible_nonadmissible} shows examples of admissible and non-admissible quadrilaterals.

Let us introduce some terminology that will be useful throughout this section. Let $q$ be an admissible quadrilateral. 
\begin{itemize}
    \item The saddle connections starting from the same singularity are the \emph{bottom sides} of $q$. 
    \item The saddle connections ending in the same singularity are the \emph{top sides} of $q$.
    \item The right-slanted bottom side of $q$ will be called the \emph{bottom right side} of $q$ and likewise we will call \emph{bottom left side} the left-slanted bottom side of $q$.
    \item Analogously, we introduce the terminology \emph{top right side} and \emph{top left side}.
\end{itemize}

It follows directly from the definition that the bottom sides of any admissible quadrilateral $q$ form a wedge, which we will refer to as the \emph{wedge of the quadrilateral} $q$.

Using this terminology, we can now define the decomposition of the translation surface $X$ into quadrilaterals that will be useful for us.

\begin{definition}[Quadrangulation]
    A quadrangulation $Q$ of $X$ is a decomposition of $X$ into a union of admissible quadrilaterals.
\end{definition}

\begin{figure}[ht]
    \centering
    \begin{tikzpicture}[scale = 0.8]
    \pgfmathsetmacro{\s}{2} 
    \pgfmathsetmacro{\rot}{13} 
    
    \draw[rotate=\rot] (0:\s) \foreach \x [count=\i] in {45,90,...,360} {
        -- (\x:\s) coordinate (p\i)
    } -- cycle;

    \draw[dashed] (p3)--(p5);
    \draw[dashed] (p2) -- (p5);
    \draw[dashed] (p1) -- (p7);
    \draw[dashed] (p1) -- (p6);

    \fill[red, opacity =0.15] (p5) -- (p3) -- (p2) -- cycle;
    \fill[red, opacity = 0.15] (p1) -- (p6) -- (p7) -- cycle;

    \fill[violet, opacity = 0.15] (p1) -- (p2) -- (p5) -- (p6) -- cycle;

    \fill[teal, opacity = 0.15] (p3) -- (p4) -- (p5) -- cycle;
    \fill[teal, opacity =0.15] (p7) -- (p8) -- (p1) -- cycle;

    \draw[>=latex, ->, thick] (2.5,0) -- (4.2,0);

    \def \s1{7};
    \draw[thick] ($(p1) + (\s1,0)$) -- ($(p2) + (\s1,0)$) -- ($(p5) + (\s1,0)$) -- ($(p6) + (\s1,0)$) -- cycle;
    \fill[violet, opacity = 0.15] ($(p1) + (\s1,0)$) -- ($(p2) + (\s1,0)$) -- ($(p5) + (\s1,0)$) -- ($(p6) + (\s1,0)$) -- cycle;

    \def \s2{6.5};
    \draw[thick] ($(p5) + (\s2, -1)$) -- ($(p2) + (\s2, -1)$)  -- ($(p2) + (\s2, -1) + (p3) - (p5)$) -- ($(p3) + (\s2, -1)$) -- cycle;
    \fill[red, opacity = 0.15] ($(p5) + (\s2, -1)$) -- ($(p2) + (\s2, -1)$)  -- ($(p2) + (\s2, -1) + (p3) - (p5)$) -- ($(p3) + (\s2, -1)$) -- cycle;

    \def \s3{10.1};
    \draw[thick] ($(p3) + (\s3, 1)$) -- ($(p5) + (\s3, 1)$) -- ($(p5) + (\s3, 1) + (p4) - (p3)$) -- ($(p4) + (\s3, 1)$) -- cycle;
    \fill[teal, opacity = 0.15]($(p3) + (\s3, 1)$) -- ($(p5) + (\s3, 1)$) -- ($(p5) + (\s3, 1) + (p4) - (p3)$) -- ($(p4) + (\s3, 1)$) -- cycle;
\end{tikzpicture}
    \caption{An example of a quadrangulation of the translation surface obtained by identifying parallel sides of a regular octagon. Note that not any decomposition into quadrilaterals automatically gives a quadrangulation $Q$, since we require $q \in Q$ to be admissible.}
    \label{fig:quadrangulation_octagon}
\end{figure}
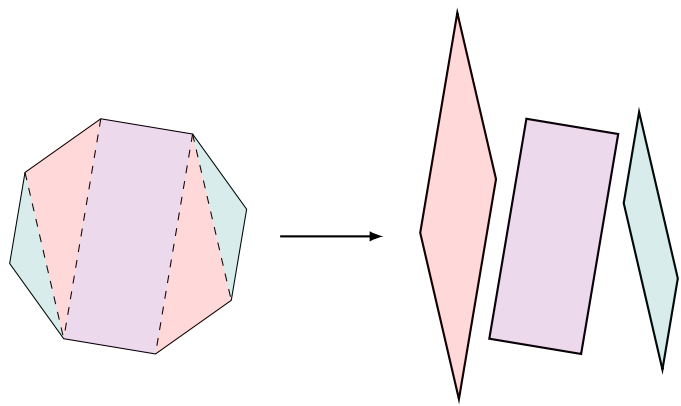

We will write $q \in Q$ if $q$ is a quadrilateral in the decomposition $Q$ and we call the \emph{wedges of} $Q$ the collection of wedges of all quadrilaterals $q \in Q$. See Figure \ref{fig:quadrangulation_octagon} for an example of a quadrangulation.

Note that \emph{every} side of a quadrilateral $q \in Q$ belongs to some wedge of $Q$, since every top side of a quadrilateral $q$ is identified with a bottom side of some (possibly different) quadrilateral $q'$. This establishes the following lemma.

\begin{lemma}
    A quadrangulation $Q$ of $X$ is completely determined by its wedges.
\end{lemma}

Given a quadrangulation $Q$ of a translation surface $X$, the main idea in order to extend the algorithm suggested initially for the torus is to apply right or left moves to the individual quadrilaterals. We will call this process \emph{diagonal change}, i.e., a diagonal change consists of replacing the left or right part of the wedge of a quadrilateral $q \in Q$ by the diagonal of the quadrilateral $q$. What needs to be ensured is that after a diagonal change, we obtain again a quadrangulation of $X$, which has the consequence that we may need to do several diagonal changes simultaneously. The basic object that allows us to do so are so-called \emph{staircases} which we will explain next.

Given a quadrilateral $q \in Q$ represented by the wedge $w = (w_\ell, w_r)$ we denote by $w_d$ the diagonal saddle connection in $q$. We will use the same terminology of \emph{left-slanted} and \emph{right-slanted} quadrilateral as we introduced above in the case of the torus, that is, we call the quadrilateral \emph{left-slanted} if the vertical ray emitted from the bottom singularity crosses the top left side of $q$ and we call it \emph{right-slanted} if the ray crosses the top right side of $q$. On the level of quadrilaterals, the following moves are possible.
\begin{enumerate}
    \item If $q$ is left-slanted, we can either keep the wedge $w = (w_\ell, w_r)$ or we can do a \emph{left diagonal change}, i.e., we replace the wedge by $(w_\ell, w_d)$. Since $q$ is left-slanted, this indeed produces a wedge, i.e., $w_d$ is right-slanted.
    \item Analogously, if $q$ is right-slanted we either keep the wedge $(w_\ell, w_r)$ or we replace it by $(w_d, w_r)$.
\end{enumerate}

\begin{definition}[Staircase]
    Given a quadrangulation $Q$ of $X$, a \emph{left staircase for} $Q$ is a subset $S \subseteq X$, which is the union of quadrilaterals $q_1, \ldots, q_n \in Q$ that are cyclically glued so that the \emph{top left} side of $q_i$ is identified with the \emph{bottom right} side of $q_{i+1}$ for $1 \leq i < n$ and to the bottom right side of $q_1$ for $i = n$.

    We obtain the definition of a \emph{right staircase} by exchanging all instances of \enquote{left} into \enquote{right} and vice versa in the definition above. 
\end{definition}

\begin{definition}[Well-slanted]
    A staircase $S$ is \emph{well-slanted}, if all $q \in S$ are \emph{slanted alike}, i.e., either all $q \in S$ are left-slanted if $S$ is a left staircase or all $q \in S$ are right-slanted if $S$ is a right staircase.
\end{definition} 

\begin{figure}[ht]
    \centering
    \begin{tikzpicture}[scale = 0.75]
    \coordinate (p1) at (0,0);
    \coordinate (p2) at (-2,1);
    \coordinate (p3) at (-0.5,3);
    \coordinate (p4) at (0.1, 5.5);
    \coordinate (p5) at (4,7);
    \coordinate (p6) at (6.3, 6);
    \coordinate (p7) at (5.3, 3.7);
    \coordinate (p8) at (3.2, 2);

    \draw[thick] (p1) -- (p2) -- (p3) -- (p4) -- (p5) -- (p6) -- (p7) -- (p8) -- cycle;
    \draw[thick] (p3) -- (p8);
    \draw[thick] (p4) -- (p7);
    \draw[dotted] (p1) -- (p3);

    \draw[white] (2,-0.5) -- (2.1,-0.5);
\end{tikzpicture}\qquad\qquad
\begin{tikzpicture}[scale = 0.75]
    \coordinate (p1) at (0,0);
    \coordinate (p2) at (-2,1);
    \coordinate (p3) at (-0.5,3);
    \coordinate (p4) at (0.1, 5.5);
    \coordinate (p5) at (4,7);
    \coordinate (p6) at (6.3, 6);
    \coordinate (p7) at (5.3, 3.7);
    \coordinate (p8) at (3.2, 2);

    \coordinate (extra) at ($(p6) + (p3)$);

    \draw[thick] (p1) -- (p3) -- (p4) --  (p5) -- (extra) -- (p6) -- (p7) -- (p8) -- cycle;
    \draw[thick] (p4) -- (p8);
    \draw[thick] (p5) -- (p7);
    \draw[dotted](p6) -- (p5);
\end{tikzpicture}
    \caption{Two right staircases. Each staircase is a union of three quadrilaterals. Performing a staircase move on the staircase on the left yields the staircase on the right.}
    \label{fig:right_staircases}
\end{figure}
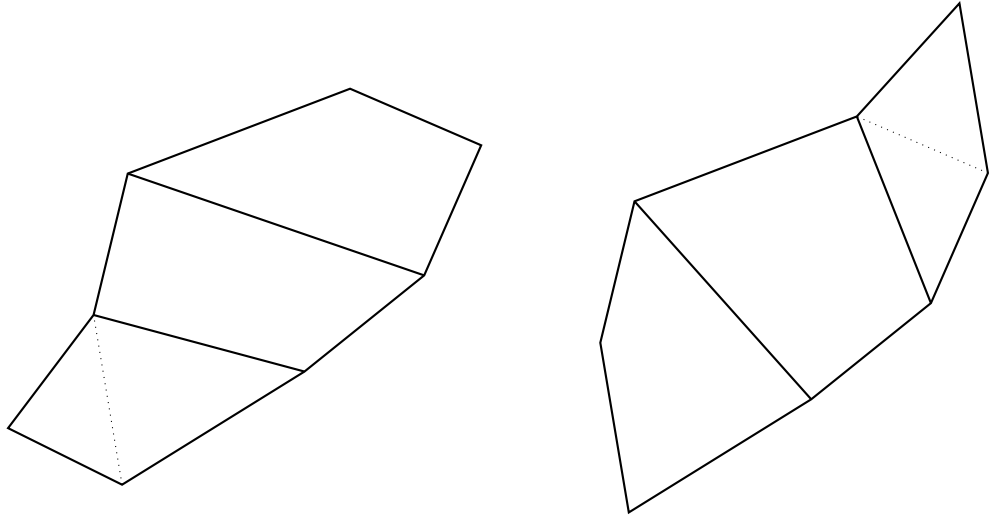

Two right staircases can be seen in Figure \ref{fig:right_staircases}. Both these staircases are given by the union of three quadrilaterals. Note that we can view a staircase topologically as a cylinder whose boundary consists of a union of saddle connections which are all slanted alike, its direction depending on whether the staircase is a left or a right staircase. 

\begin{definition}[Staircase move]\label{def:staircase_move}
    Given a quadrangulation $Q$ and a well-slanted left (or right) staircase $S$, the \emph{staircase move} in $S$ is the operation consisting of performing left (or right) diagonal changes in all quadrilaterals of $S$ simultaneously.
    \end{definition}

Figure \ref{fig:right_staircases} also provides us with an example of a staircase move: Applying such a move to the staircase on the left produces the staircase on the right. We want to stress that given a quadrangulation $Q$, there may exist none, one or several well-slanted staircases. In the first case, no staircase move is possible while when there are multiple staircases there is a choice of moves. We will see that it is still possible in these cases to unambigously define an algorithm. 

The importance of staircases comes from the fact that if $Q$ is a quadrangulation of a surface $X$ and $S$ is a well-slanted staircase in $Q$, then the staircase move in $S$ produces (part of) a new quadrangulation $Q'$ of $X$ (see Lemma \ref{lem:staircase_move_on_data}). Moreover, one can show that staircase moves are in some sense (that we will make precise) the smallest unions of quadrilaterals where one can perform diagonal changes consistently in order to keep a quadrangulation (see Lemma \ref{lem:irreducibility_of_staircases}). 

Before elaborating on the diagonal changes algorithm much in the same way we did on Rauzy--Veech induction in section \ref{sec:rauzy_veech}, we want to mention that contrary to before we will not be able to define the algorithm for any stratum of the moduli space of Abelian differentials, but we will need to restrict ourselves to certain connected components called the \emph{hyperelliptic components of strata}. The surfaces we consider come with extra symmetry which ensures the existence of quadrangulations. We will give more details on hyperelliptic components in section \ref{sec:existence_of_quadrangulations}.

\subsection{The Space of Quadrangulations}\label{sec:space_of_quadrangulations}

In general, a quadrangulation $Q$ of a translation surface $X$ will consist of several quadrilaterals, making it necessary to keep track of which sides of which quadrilaterals are identified. Just as in the case of the Rauzy--Veech induction we will use a \emph{combinatorial datum} which we will also denote by $\boldsymbol{\pi}$ or $\boldsymbol{\pi}_Q$ if we want to emphasize the dependence on the quadrangulation.

\begin{definition}[Combinatorial datum]\label{def:combinatorial_datum_dc}
    Let $Q$ be a quadrangulation with $k$ quadrilaterals, labeled $q_1, \ldots, q_k$. The \emph{combinatorial datum} $\boldsymbol{\pi}$ $(= \boldsymbol{\pi}_Q)$ of $Q$ is a pair $({\pi}_\ell, {\pi}_r)$ of permutations of $[k]$ such that the following hold.
    \begin{enumerate}
        \item For each $1 \leq i \leq k$, the top left side of $q_i$ is identified with the bottom right side of $q_{{\pi}_\ell(i)}$.
        \item For each $1 \leq i \leq k$, the top right side of $q_i$ is identified with the bottom left side of $q_{{\pi}_r(i)}$.
    \end{enumerate}
\end{definition}

\begin{example}[Combinatorial datum]\label{ex:combinatorial_datum_dc}
    Consider the decomposition of the regular octagon into a quadrangulation $Q$ from Figure \ref{fig:quadrangulation_octagon}, where we label the three quadrilaterals $q_1, q_2, q_3$ from left to right as they appear on the right in the figure. We have
    \begin{alignat*}{3}
        {\pi}_\ell(1) &= 3, \quad {\pi}_\ell(2) &&= 1, \quad {\pi}_\ell(3) &&= 2, \\
        {\pi}_r(1) &= 3, \quad
        {\pi}_r(2) &&= 2, \quad
        {\pi}_r(3) &&= 1.
    \end{alignat*}
    Using the usual cyclic notation, we can write the combinatorial datum more succinctly as
    \begin{equation*}
        \boldsymbol{\pi} = ({\pi}_\ell, {\pi}_r) = 
        \big(
        \begin{pmatrix}
            1 & 3 & 2
        \end{pmatrix},
        \begin{pmatrix}
            1 & 3
        \end{pmatrix}
        \big)
    \end{equation*}
\end{example}

\begin{figure}[ht]
    \centering
{\begin{tikzcd}
	& 1 \\
	\\
	2\arrow[loop left, distance = .9cm, in=160, out = 200, "r", color={rgb,255:red,214;green,92;blue,92}] && 3
	\arrow["r", color={rgb,255:red,214;green,92;blue,92}, from=1-2, to=3-3]
	\arrow["r", color={rgb,255:red,214;green,92;blue,92}, curve={height=-12pt}, from=3-3, to=1-2]
	\arrow["\ell", color={rgb,255:red,92;green,92;blue,214}, curve={height=-12pt}, from=1-2, to=3-3]
	\arrow["\ell", color={rgb,255:red,92;green,92;blue,214}, from=3-3, to=3-1]
	\arrow["\ell", color={rgb,255:red,92;green,92;blue,214}, from=3-1, to=1-2]
\end{tikzcd}}
    \caption{The labeled directed graph $G_Q$ representing the combinatorial datum induced by the quadrangulation $Q$ from Example \ref{ex:combinatorial_datum_dc}.}
    \label{fig:combinatorial_datum_graph_dc}
\end{figure}
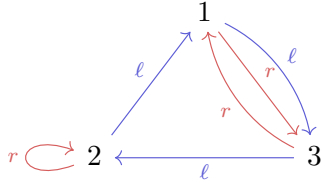

Let us mention that another convenient way to describe the combinatorial datum of a quadrangulation is by using a labeled directed graph $G_Q$. In particular, these graphs are used by Ferenczi and Zamboni in \cite{ferenczi2011eigenvalues} where they describe a related algorithm on the level of a special class of IETs. Each vertex of $G_Q$ corresponds to a quadrilateral $q_i \in Q$ and is denoted by the number $i$. The edges are labeled by $r$ or $\ell$ such that for each $1 \leq i \leq k$ there is an $\ell$-edge from $i$ to ${\pi}_\ell(i)$ and an $r$-edge from $i$ to ${\pi}_r(i)$. The graph $G_Q$ corresponding to the quadrangulation $Q$ from Example \ref{ex:combinatorial_datum_dc} can be seen in Figure \ref{fig:combinatorial_datum_graph_dc}.

For a quadrangulation $Q$ given by quadrilaterals $q_1, \ldots, q_k$, or equivalently given by wedges $w_1, \ldots, w_k$, it follows that there is a one-to-one correspondence between quadrilaterals (or wedges) and bundles of saddle connections. A labeling of the quadrilaterals thus induces a labeling of bundles of saddle connections, i.e., the saddle connections $w_{i, \ell}$ and $w_{i, r}$ belonging to the quadrilateral $q_i$ (or the wedge $w_i$) belong to the bundle $\Gamma_i$. More precisely, we will use the notation that $w_{i, \ell}$ belongs to the bundle $\Gamma_{i, \ell}$ and $w_{i, r}$ belongs to $\Gamma_{i, r}$, subdividing each bundle into its left- and right-slanted saddle connections. 

For any wedge $w_i = (w_{i, \ell}, w_{i, r})$, it is clear from context to which bundle $\Gamma_i = \Gamma_{i, \ell} \cup \Gamma_{i, r}$ this wedge belongs, which allows us to abuse the notation slightly and identify the saddle connections in the wedges with their holonomy vectors. Note that by definition of the combinatorial datum (Definition \ref{def:combinatorial_datum_dc}), the following relation must be satisfied.
\begin{equation}\label{eq:wedges_train_tracks}
    w_{i, \ell} + w_{{\pi}_\ell(i), r} = w_{i, r} + w_{{\pi}_r(i), \ell}, \quad 1 \leq i \leq k.
\end{equation}
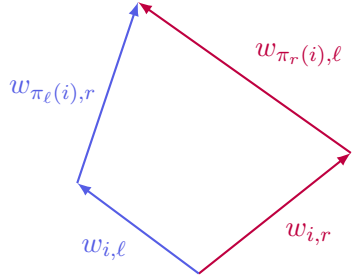
\begin{figure}[hb]
    \centering
    \begin{tikzpicture}[scale = 0.8]
    \draw[thick, purple, >=latex, ->] (0,0) -- (2.5,2) node[midway, below right] {$w_{i,r}$};
    \draw[thick, purple, >=latex, ->] (2.5,2) -- (-1,4.5) node[midway, above right] {$w_{\pi_r(i), \ell}$};
    \draw[thick, noamblue, >=latex, ->] (0,0) -- (-2,1.5) node[midway, below left] {$w_{i, \ell}$};
    \draw[thick, noamblue, >=latex,->] (-2, 1.5) -- (-1, 4.5) node[midway, left] {$w_{\pi_\ell(i), r}$};
\end{tikzpicture}
    \caption{The figure illustrates the train-track relations from \eqref{eq:wedges_train_tracks}. Note that these relations ensure that the saddle connections fit properly together to form quadrilaterals.}
    \label{fig:wedges_train_tracks}
\end{figure}
We will call equations \eqref{eq:wedges_train_tracks} the \emph{train-track relations}. This also motivates the following definition.

\begin{definition}[Length datum]
    Let $Q$ be a quadrangulation characterized by $k$ wedges $w_1, \ldots, w_k$. Identifying the saddle connections of the wedges with their holonomy vectors and writing 
    \begin{equation*}
        w_i = (w_{i, \ell}, w_{i, r}) \in (\R_- \times \R_+) \times (\R_+ \times \R_+),
    \end{equation*}
    where $\R_- = \{x \in \R \mid x < 0\}$ and $\R_+ = \{x \in \R \mid x > 0\}$, the \emph{length datum} $\boldsymbol{w}$ of $Q$ is given by
    \begin{equation*}
        \boldsymbol{w} = \big( (w_{i, \ell}, w_{i, r})\big)_{i = 1}^k  \in \big( (\R_- \times \R_+) \times (\R_+ \times \R_+) \big)^k.
    \end{equation*}
\end{definition}

Any quadrangulation $Q$ comes equipped with a combinatorial and a length datum. Conversely, starting with any combinatorial datum $\boldsymbol{\pi}$, any length datum $\boldsymbol{w}$ satisfying the train-track relations \eqref{eq:wedges_train_tracks} defines a quadrangulation $Q = (\boldsymbol{\pi}, \boldsymbol{w})$. 

\begin{remark}
    When denoting a quadrangulation by $Q = (\boldsymbol{\pi}, \boldsymbol{w})$, it is implicitly assumed that the length datum $\boldsymbol{w}$ satisfies the train-track relations \eqref{eq:wedges_train_tracks}.
\end{remark}

\subsection{Bi-partite IETs}

In section \ref{sec:rauzy_veech}, we saw that IETs and the vertical linear flow on translation surfaces are intimately related. Choosing a vertical section $S$ in a translation surface $X$ induces an IET on this section, and given an IET on some interval we can construct a translation surface $X$ utilizing the zippered rectangle decomposition (see Theorem \ref{thm:suspensions_existence}). The situation for quadrangulations is similar. Since, as we mentioned above, we will need to restrict ourselves to hyperelliptic components of strata, this additional symmetry restriction reduces the number of IETs that are possibly related to the vertical linear flow in a quadrangulation. As we wil see, the union of the wedges of $Q$ provide a convenient section for the vertical linear flow on the associated surface inducing a so-called \emph{bi-partite IET} as the Poincaré first return map to this section. 

The basic idea is the following. Suppose we have an IET $T \colon I \to I$ on some interval $I$ on $d$ intervals, which is characterized by the combinatorial and length data $(\boldsymbol{\pi}, \boldsymbol{\lambda})$. Any such IETs naturally induces \emph{two} partitions of $I$ given by the union of the intervals in the orders specified by the combinatorial datum $\boldsymbol{\pi} = (\pi_\mathrm{top}, \pi_\mathrm{bot})$. If $d$ is even, we can label the intervals as $\ell$ (\enquote{left}) and $r$ (\enquote{right}) in an alternative fashion for both partitions such that the number of intervals labeled $\ell$ and $r$ respectively agree. 

A bipartite IET respects this labeling in the sense that an interval labeled $\ell$ with respect to $\pi_\mathrm{top}$ is mapped to an interval labeled $r$ with respect to $\pi_\mathrm{bot}$ and analogously an interval labeled $r$ with respect to $\pi_\mathrm{top}$ is mapped to an intervall labeled $\ell$ with respect to $\pi_\mathrm{bot}$. Moreover, we require the length of the pairs of intervals to be invariant under the IET, meaning that the length of the pair of intervals in position $i$ with respect to $\pi_\mathrm{top}$ must sum to the same number as the length of the pair of intervals in the same position $i$ with respect to $\pi_\mathrm{bot}$. 

\begin{figure}[ht]
    \centering
    \begin{tikzpicture}[scale = 0.7]
    \draw[thick] (0,0) -- (5,0);
    \draw[thick] (6,0) -- (11.5,0);
    \draw[thick] (12.5,0) -- (17.5,0);

    \foreach \p in {0, 5, 6, 11.5, 12.5, 17.5}{
        \draw[thick] (\p,0) - ++(0,0.1);
        \draw[thick] (\p,0) - ++(0,-0.1);    
    }
    \foreach \p in {4, 10, 15.5}{
        \draw[thick] (\p,0) - ++(0,0.1);
    }
    \foreach \p in {2, 7.5, 13.5}{
        \draw[thick] (\p,0) - ++(0,-0.1);
    }

    \fill[teal, opacity = 0.2] (0,0) rectangle (5,0.1);
    \fill[teal, opacity = 0.2] (7.5,0) rectangle (11.5, -0.1);
    \fill[teal, opacity = 0.2] (12.5,0) rectangle (13.5, -0.1);

    \fill[blue, opacity =0.2] (6,0) rectangle (11.5,0.1);
    \fill[blue, opacity  = 0.2] (13.5,0) rectangle (17.5,-0.1);
    \fill[blue, opacity = 0.2] (6,0) rectangle (7.5,-0.1);

    \fill[red, opacity = 0.2] (12.5,0) rectangle (17.5,0.1);
    \fill[red, opacity = 0.2] (0,0) rectangle (5,-0.1);

    \node[above] at (2,0) {$J_{1,\ell}$};
    \node[above] at (4.5,0) {$J_{1,r}$};
    \node[above] at (8,0) {$J_{2,\ell}$};
    \node[above] at (10.75,0) {$J_{2,r}$};
    \node[above] at (14,0){$J_{3,\ell}$};
    \node[above] at (16.5, 0) {$J_{3, r}$};

    \node[below] at (1,0) {$I_{1,\ell}$};
    \node[below] at (3.5,0) {$I_{1,r}$};
    \node[below] at (6.75,0) {$I_{2,\ell}$};
    \node[below] at (9.5,0) {$I_{2,r}$};
    \node[below] at (13,0){$I_{3,\ell}$};
    \node[below] at (15.5, 0) {$I_{3, r}$};


    \draw[>=latex, ->] (2,0.7) to[out=45, in =90] (5.5,0.5) to [out = -90, in = 90] (5.5,-0.5) to [out = -90, in = 225] (9.5,-0.7);
    \draw[>=latex, ->] (4.5,0.7) to[out=45, in = 90] (5.3,0.5) to [out = -90, in = 90] (5.3, -1) to [out = -90, in = 225] (13,-0.7);
    \draw[>=latex, ->] (8,0.7) to[out=45, in =90] (12,0.5) to[out=-90, in=90] (12, -0.5) to [out=-90, in=225] (15.5, -0.7);
    \draw[>=latex, ->] (10.75,0.7) to[out = 135, in = 90] (5.7,0)  to [out=-90, in =225] (6.75,-0.7);
    \draw[>=latex, ->] (14,0.7) to[out=150, in =20] (0.8,1.4) to[out = 200, in = 130] (-0.5,-0.8) to[out=-50, in=225] (3.5,-0.7);
    \draw[>=latex, ->] (16.5,0.7) to[out=150, in=5] (5.5,3.3) to[out = 185, in=60](-1.3,0.2) to[out=240, in=150] (-1,-1) to[out=-30, in=225] (1,-0.7);
\end{tikzpicture}
    \caption{A bipartite IET on $d = 6$ intervals.}
    \label{fig:bipartite_IET}
\end{figure}
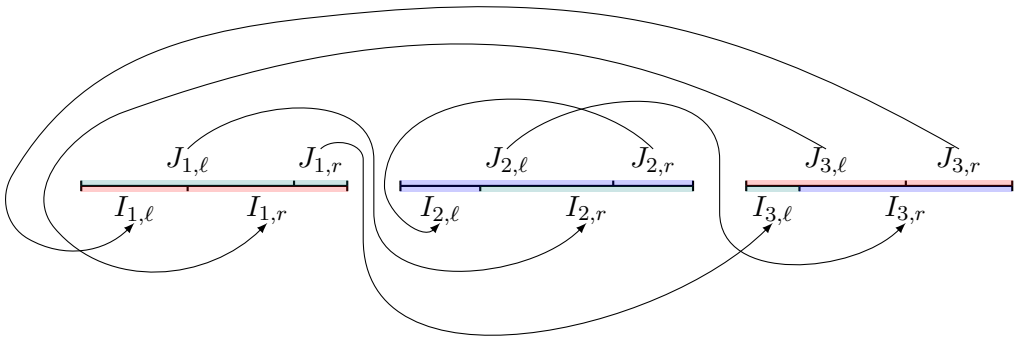

An example corresponding to the combinatorial datum
\begin{equation*}
    \boldsymbol{\pi} = 
    \begin{pmatrix}
        A & B & C & D & E & F \\
        F & E & D & A & B & C
    \end{pmatrix}
\end{equation*}
can be seen in Figure \ref{fig:bipartite_IET}.

Let us formalize this idea. For ease of notation, we will assume that the alphabet $\mathcal{A}$ associated to the IET of $d = 2k$ intervals is always given by
\begin{equation}\label{eq:alphabet_bipartite}
    \mathcal{A} = \{(1, \ell), (1, r), \ldots, (k, \ell), (k, r)\},
\end{equation}
and let us further assume that $\pi_\mathrm{top} = \operatorname{id}$ so that $\pi_\mathrm{bot}$ coincides with the monodromy invariant $p$. 

\begin{definition}[Bi-partite IET]\label{def:bipartite_IET}
    Suppose we are given combinatorial and length data $(\boldsymbol{\pi}, \boldsymbol{\lambda})$ with the alphabet $\mathcal{A}$ from \eqref{eq:alphabet_bipartite} and $\pi_\mathrm{top} = \operatorname{id}$, so that we may write $p = \pi_\mathrm{bot}$ for the monodromy invariant and let the length datum be presented by
    \begin{equation*}
        \boldsymbol{\lambda} = \big( (\lambda_{{1}, \ell}, \lambda_{{1}, r}),\ldots, (\lambda_{{k},\ell}, \lambda_{{k}, r}) \big) \quad \in (\R_+ \times \R_+)^k, 
    \end{equation*}
    such that
    \begin{equation}\label{eq:length_train_tracks}
        \lambda_{i, \ell} + \lambda_{i, r} = \lambda_{p(i,\ell)} + \lambda_{p(i, r)} \quad \text{for all } 1 \leq i \leq k.
    \end{equation}

    For $1\leq i \leq k$, we set $I_i = (\lambda_{i, \ell}, \lambda_{i, r}) \subseteq \R$ and let
    \begin{align*}
        I_{i, \ell} &= (-\lambda_{i, \ell}, 0), &&I_{i, r} = (0, \lambda_{i, r}), \\
        J_{i, \ell} &= (-\lambda_{i, \ell}, -\lambda_{i, \ell} + \lambda_{p(i, r)}), && J_{i, r} = (\lambda_{i, r} - \lambda_{p(i,r)}, \lambda_{i,r}).
    \end{align*}
    The \emph{bipartite IET} with data $(\boldsymbol{\pi}, \boldsymbol{\lambda})$ is the application that maps $J_{i, \ell}$ to $I_{p(i, r)}$ by translation. Note that the map is not defined at the points $\lambda_{i, d} \coloneqq \lambda_{i, \ell} + \lambda_{p(i, \ell)} = \lambda_{i,r} = \lambda_{p(i, r)}$.      
\end{definition}

The notation we have chosen in Definition \ref{def:bipartite_IET} highlights the connection to the general definition of an IET, or more precisely, the connection to Lemma \ref{lem:IET_characterization}. To make the connection to the linear flow on a quadrangulation more apparent, we want to move away from this convention in the following way. Since intervals labeled with $\ell$ get mapped to intervals labeled with $r$ and vice versa, instead of specifying the monodromy invariant $p \colon [d] \to [d]$, we may specify two permutations $\pi_\ell$ and $\pi_r$, both defined on $[k]$, such that $\pi_\ell(i) = j$ such that $p(i, \ell) = (j, r)$ and $\pi_r(i) = j$, where $j$ is such that $p(i, r) = (j, \ell)$. Moreover, we adopt the convention that we multiply the lengths of left intervals, i.e., of intervals labeled $\ell$, by -1 so that $\boldsymbol{\lambda} \in (\R_- \times \R_+)^k$.

Utilizing this adapted notation, equation \eqref{eq:length_train_tracks} transforms to
\begin{equation}\label{eq:length_train_tracks_adapted}
    \lambda_{i, \ell} + \lambda_{\pi_\ell(i), r} = \lambda_{i, r} + \lambda_{\pi_r(i), \ell} \quad \text{for all } 1 \leq i \leq k,
\end{equation}
where we immediately recognize the train-track relations given in equation \eqref{eq:wedges_train_tracks}. For this reason, we will call equation \ref{eq:length_train_tracks_adapted} the \emph{train-track relations for the length datum}. In the new notation, the intervals $I_{i, \varepsilon}$ and $J_{i, \varepsilon}$ are given by
\begin{align*}
    I_{i, \ell} &= (\lambda_{i, \ell}, 0), &&I_{i, r} = (0, \lambda_{i, r}), \\
    J_{i, \ell} &= (\lambda_{i, \ell}, \lambda_{i, \ell} + \lambda_{\pi_\ell(i), r}), && J_{i, r} = (\lambda_{i, r} + \lambda_{\pi_r(i),\ell}, \lambda_{i,r}).
\end{align*}

\begin{figure}[ht]
    \centering
\[\begin{tikzcd}
	1 && 2 && 1 & 2 \\
	3 && 4 && 3 & 4 \\
	5 && 6
	\arrow[from=1-5, to=1-6]
	\arrow[from=1-6, to=2-6]
	\arrow[from=2-5, to=1-5]
	\arrow[from=2-6, to=2-5]
	\arrow[curve={height=-6pt}, from=1-1, to=2-3]
	\arrow[tail reversed, from=1-3, to=3-1]
	\arrow[curve={height=6pt}, from=2-1, to=3-3]
	\arrow[from=2-3, to=2-1]
	\arrow[from=3-3, to=1-1]
\end{tikzcd}\]
    \caption{Representation of combinatorial data as a directed graph. The graph on the left corresponds to the example from Figure \ref{fig:bipartite_IET} and has a bi-partite structure where the partitions are given by the even and odd numbers. The graph on the right corresponds to the IET from Example \ref{ex:octagon_IET} and even though it admits a bi-partite structure, it is not the one where the partitions are given by the even and odd numbers.}
    \label{fig:bipartite_IET_graph}
\end{figure}
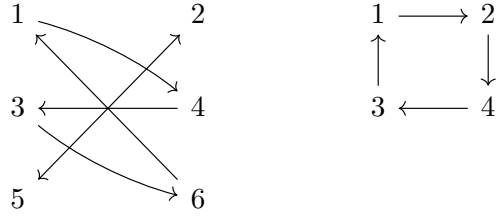

Let us quickly highlight why the name \emph{bi-partite} IET is appropriate. In general we can describe the combinatorial datum of any IET using a directed graph, where the vertices are given by the elements of the alphabet $\mathcal{A}$, in the case here this is just $[d]$, and we put an edge from $i$ to $j$ if and only if $p(i) = j$, where $p$ is the monodromy invariant. For bi-partite IETs, this graph has a specific bi-partite structure, where the partition is given by the odd and the even numbers in $[d]$. This is not the case for IETs in general, see Figure \ref{fig:bipartite_IET_graph}. Let us stress that this particular bi-partite structure on the graph is necessary, but not sufficient, since we also impose the train-track relations for the length datum.

The following result follows at once from the construction of bi-partite IETs. 

\begin{lemma}[Cross sections of quadrangulations]\label{lem:cross_section_quadrangulation}
    Given a quadrangulation $Q = (\boldsymbol{\pi}, \boldsymbol{w})$, the Poincaré first return map $F$ of the vertical flow to the union of the interiors of the wedges of $Q$ is conjugate to the bipartite IET $T$ given by the combinatorial and length data $(\boldsymbol{\pi}, \boldsymbol{\lambda})$, where $\boldsymbol{\lambda}$ consists of the real parts of the wedges $\boldsymbol{w}$.

    More precisely, if we write $\iota$ for the projection that maps a point $z$ of the wedge $w_i$ to the point $\RE(z) \in I_i$, then the following diagram commutes.
\[\begin{tikzcd}
	{\bigsqcup_i w_i} && {\bigsqcup_i w_i} \\
	\\
	{\bigsqcup_i I_i} && {\bigsqcup_i I_i}
	\arrow["F", from=1-1, to=1-3]
	\arrow["T", from=3-1, to=3-3]
	\arrow["\iota"', from=1-1, to=3-1]
	\arrow["\iota", from=1-3, to=3-3]
\end{tikzcd}\]
\end{lemma}

Let us stress the similarity of Lemma \ref{lem:cross_section_quadrangulation} to the way we explained how general IETs arise naturally as Poincaré first return maps in section \ref{sec:rauzy_veech}. Exactly as in the case of the zippered rectangle construction, there is also a way to move in the converse direction, i.e., to obtain a quadrangulation starting with a bi-partite IET characterized by $(\boldsymbol{\pi}, \boldsymbol{\lambda})$. Let us explain how to build such suspensions over bi-partite IETs.

\begin{definition}[Suspension datum]
    A \emph{suspension datum} $\boldsymbol{\tau}$ for the bi-partite IET characterized by $(\boldsymbol{\pi}, \boldsymbol{\lambda})$ is a vector
    \begin{equation*}
        \boldsymbol{\tau} = \big( (\tau_{1, \ell}, \tau_{1, r}), \ldots, (\tau_{k, \ell}, \tau_{k, r})\big) \quad \in (\R_+ \times \R_+)^k
    \end{equation*}
    that satisfies the \emph{train-track relations for suspension data}
    \begin{equation*}
        \tau_{i, \ell} + \tau_{\pi_\ell(i), r} = \tau_{i, r} + \tau_{\pi_r(i), \ell} \quad \text{for } 1 \leq i \leq k.
    \end{equation*}
\end{definition}

From \emph{any} suspension datum $\tau$ satisfying the corresponding train-track relations, we can build a quadrangulation $Q = (\boldsymbol{\pi}, \boldsymbol{\lambda}, \boldsymbol{\tau}) = (\boldsymbol{\pi}, \boldsymbol{w})$ by setting the wedges of $Q$ to be
\begin{equation*}
    w_{i, \ell} = \lambda_{i, \ell} + \ii \tau_{i, \ell}, \quad w_{i, r} = \lambda_{i, r} + \ii \tau_{i, r}.
\end{equation*}

Morally, the suspension data $\boldsymbol{\tau}$ provides the \enquote{height} of the quadrilaterals where the train-track relations ensure that we obtain well-defined quadrilaterals.

\begin{figure}[ht]
    \centering
    \begin{tikzpicture}[scale = 0.5]
    \foreach \y in {-1,7}{
        \draw[thick] (0,\y) -- (7,\y);
        \draw[thick] (8,\y) -- (11,\y);
        \draw[thick] (12,\y) -- (19,\y);
    }

    \draw[thick, noamblue] (0,2) -- (3,0) -- (7,3);
    \draw[thick, noamblue, dashed] (0,2) -- (2,6) -- (7,3);

    \draw[thick, noamblue] (8,2) -- (9,0) -- (11,4);
    \draw[thick, noamblue, dashed] (8,2) -- (10,6) -- (11,4);

    \draw[thick, noamblue] (12,3) -- (17,0) -- (19,4);
    \draw[thick, noamblue, dashed] (12,3) -- (16,6) --  (19,4);

    \foreach \x in {0,2,7,8,10,11,12,16,19}{
        \draw[thick] (\x,7) - ++(0,0.1);
        \draw[thick] (\x,7) - ++(0,-0.1);
    }
    \foreach \x in {0,3,7,8,9,11,12,17,19}{
        \draw[thick] (\x,-1) - ++(0,0.1);
        \draw[thick] (\x,-1) - ++(0,-0.1);
    }

    \node[above] at (1,7) {$J_{1, \ell}$};
    \node[above] at (4,7) {$J_{1, r}$};
    \node[above] at (9,7) {$J_{2, \ell}$};
    \node[above] at (10.5,7) {$J_{2, r}$};
    \node[above] at (14,7) {$J_{3, \ell}$};
    \node[above] at (17.5,7) {$J_{3, r}$};

    \node[below] at (1.5,-1) {$I_{1, \ell}$};
    \node[below] at (5,-1) {$I_{1, r}$};
    \node[below] at (8.5,-1) {$I_{2, \ell}$};
    \node[below] at (10,-1) {$I_{2, r}$};
    \node[below] at (14.5,-1) {$I_{3, \ell}$};
    \node[below] at (18,-1) {$I_{3, r}$};
  
\end{tikzpicture}
    \caption{A suspension corresponding to the bi-partite IET from Example \ref{ex:suspension_bipartite}.}
    \label{fig:bipartite_IET_suspension}
\end{figure}
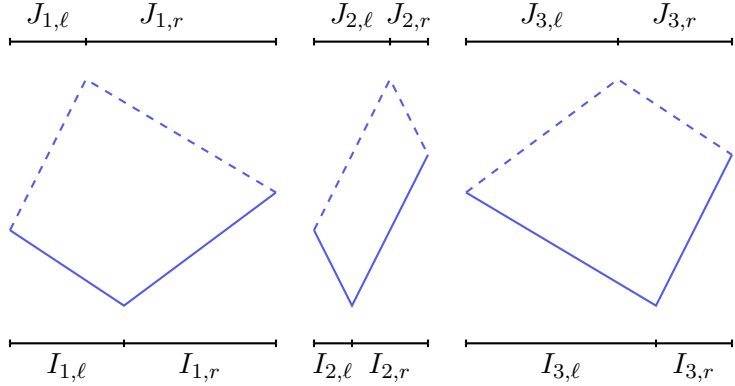

\begin{example}[Suspension over bi-partite IET]\label{ex:suspension_bipartite}
    Consider the bi-partite IET characterized by the combinatorial and length data $(\boldsymbol{\pi}, \boldsymbol{\lambda})$, where $\boldsymbol{\pi} = (\pi_\ell, \pi_r)$ is given by
    \begin{align*}
        \pi_\ell &= \begin{pmatrix}
            1 & 2 & 3
        \end{pmatrix}, \\
        \pi_r &= \begin{pmatrix}
            1 & 3
        \end{pmatrix}
    \end{align*}
    and $\boldsymbol{\lambda}$ is given by
    \begin{equation*}
        \boldsymbol{\lambda} = \big( (-3, 4), (-1, 2), (-5,2)\big).
    \end{equation*}
    It is easily verified that $\boldsymbol{\lambda}$ satisfies the train-track relations for the length datum. An example of a suspension datum verifying the corresponding train-track relations is given by
    \begin{equation*}
        \boldsymbol{\tau} = \big( (2, 3), (2, 4), (3, 4)\big).
    \end{equation*}
    This defines a quadrangulation of a translation surface of genus $\mathbf{g} = 2$, which is illustrated in Figure \ref{fig:bipartite_IET_suspension}.
\end{example}

This construction is done in a way, such that we exactly obtain a \enquote{dual} statement to Lemma \ref{lem:cross_section_quadrangulation}.

\begin{lemma}[Suspension of bi-partite IETs]\label{lem:bipartite_IET_suspension}
    Given a bi-partite IET $T$ characterized by $(\boldsymbol{\pi}, \boldsymbol{\lambda})$ and a suspension datum $\boldsymbol{\tau}$, let $Q = (\boldsymbol{\pi}, \boldsymbol{\lambda}, \boldsymbol{\tau})$ be the associated quadrangulation. Then, the Poincaré first return map of the vertical linear flow on the associated translation surface to the union of the interiors of the wedges of $Q$ is conjugate to $T$. In other words, denoting by $\kappa$ the bijective assignment $\kappa \colon (\lambda_i, \tau_i) \mapsto w_i$, the following diagram commutes.

\[\begin{tikzcd}
	{\bigsqcup_i I_i} && {\bigsqcup_i I_i} \\
	\\
	{\bigsqcup_i w_i} && {\bigsqcup_iw_i}
	\arrow["T", from=1-1, to=1-3]
	\arrow["F", from=3-1, to=3-3]
	\arrow["\kappa"', from=1-1, to=3-1]
	\arrow["\kappa", from=1-3, to=3-3]
\end{tikzcd}\]
\end{lemma}

\begin{remark}
    Lemma \ref{lem:bipartite_IET_suspension} can be seen as analog to Theorem \ref{thm:suspensions_existence}.
\end{remark}

\subsection{Staircase Moves}\label{sec:staircase_moves}
In Definition \ref{def:staircase_move} we have already given a geometric description of a staircase move. Now we want to describe how such a move acts on the data $(\boldsymbol{\pi}, \boldsymbol{w})$ representing a quadrangulation $Q$. This can be seen as an analog to how the Rauzy--Veech induction step introduced in section \ref{sec:rauzy_veech} acts on the combinatorial and length data of an IET.

Given such a quadrangulation, recall that the top right side of $q_i \in Q$ is glued to the bottom left side of $q_{\pi_r(i)}$. So if
\begin{equation*}
    \begin{pmatrix}
        i & \pi_r(i) & \ldots & \pi_r^n(i)
    \end{pmatrix}
\end{equation*}
is a cycle of $\pi_r$, that is $\pi_r^j(i) \neq i$ for $1 \leq j \leq n$ but $\pi_r^{n+1}(i) = i$, the associated quadrilaterals
\begin{equation*}
    \{q_i, q_{\pi_r(i)}, \ldots, q_{\pi_r^n(i)}\}
\end{equation*}
are glued to each other by identifying top right and bottom left sides. Of course, an analogous statement is true if we replace $\pi_r$ by $\pi_\ell$. Thus we have shown the following lemma.

\begin{lemma}[Cycles determine staircases]
    There is a one-to-one correspondence between cycles of $\pi_r$ (or $\pi_\ell)$ and right (or left) staircases in the following sense.
    \begin{enumerate}    
        \item A cycle $c = \begin{pmatrix}
            i & \pi_r(i) & \ldots & \pi_r^n(i)
        \end{pmatrix}$
        of $\pi_r$ corresponds to a right staircase $S_c = \bigcup_{j = 0}^n q_{\pi_r^j(i)}$.
        \item A cycle $c = \begin{pmatrix}
            i & \pi_\ell(i) & \ldots & \pi_\ell^n(i)
        \end{pmatrix}$
        of $\pi_\ell$ corresponds to a left staircase $S_c = \bigcup_{j = 0}^n q_{\pi_\ell^j(i)}$.
    \end{enumerate}
\end{lemma}

Abusing the notation, we will denote by $S = S_c$ both the union of the quadrilaterals as a subset of the translation surface $X$ and the collection of quadrilaterals, meaning that we will both write $S \subseteq X$ and $q \in S$, where $q$ is a quadrilateral of the quadrangulation $Q$ contained in the staircase. 

For each wedge $w_i = (w_{i, \ell}, w_{i, r})$ corresponding to a quadrilateral $q_i \in Q$, we denote by $w_{i,d}$ the (forward) diagonal of the quadrilateral, i.e., we set
\begin{equation*}
    w_{i, d} \coloneqq w_{i, \ell} + w_{\pi_\ell(i), r} = w_{i, r} + w_{\pi_r(i), \ell},
\end{equation*}
where the second equality is justified by the train-track relations on the wedges. The following fact is immediate.

\begin{lemma}
    A right (or left) staircase $S_c$ is well-slanted if and only if $\RE(w_{i, d}) < 0 $ (or $\RE(w_{i,d}) > 0)$ for all $i \in c$.
\end{lemma}

Our goal now is to show that applying a staircase move to a staircase $S_c$ produces a new quadrangulation, and moreover we want to describe its data in relation to the initial quadrangulation.

So to this end, let $c$ be a cycle of $\pi_r$ and assume that $S_c$ is well-slanted. Since in a diagonal change we replace a side of a wedge with its diagonal, the new length data $\boldsymbol{w}'$ is given by
\begin{equation*}
    w_i' =
    \begin{cases}
        (w_{i, d}, w_{i, r}) \quad &\text{if } i \in c,\\
        w_i & \text{else.}
    \end{cases}
\end{equation*}

Since $S_c$ is well-slanted, it follows that $w'_i$ are wedges for all $i$. Indeed, the diagonal $w_{i, d}$ is left-slanted by definition, so that $w'_{i,\ell} \in \R_-\times\R_+$ and $w'_{i, r} \in \R_+\times\R_+$.

Analogously, if $c$ is a cycle of $\pi_\ell$, the new length data $\boldsymbol{w}'$ is given by
\begin{equation*}
    w_i' =
    \begin{cases}
        (w_{i, \ell}, w_{i, d}) \quad &\text{if } i \in c,\\
        w_i & \text{else.}
    \end{cases}
\end{equation*}

The wedges $\boldsymbol{w}'$ determine a new quadrangulation $Q'$ in the following way. 

\begin{enumerate}
    \item If $c$ is a cycle of $\pi_r$ containing $i$, then the wedge $w'_i$ is the wedge of the quadrilateral $q_i'$, which has the following edges:
    \begin{center} 
        \begin{tabular}{ccc}
        \toprule
            & left & right \\
        \midrule  
        top &  $w_{\pi_\ell\pi_r(i), r}$ & $w_{\pi_r(i), d}$ \\
        \midrule
        bottom  & $w_{i,d}$ & $w_{i,r}$ \\
        \bottomrule

    \end{tabular}
    \end{center}
    \vspace{10pt}
     \item If $c$ is a cycle of $\pi_\ell$ containing $i$, then the wedge $w'_i$ is the wedge of the quadrilateral $q_i'$, which has the following edges:
     \begin{center} 
        \begin{tabular}{ccc}
        \toprule
            & left & right \\
        \midrule  
        top &  $w_{\pi_\ell(i), r}$ & $w_{\pi_r\pi_\ell(i), d}$ \\
        \midrule
        bottom  & $w_{i,\ell}$ & $w_{i,d}$ \\
        \bottomrule
    \end{tabular}
    \end{center}
\end{enumerate}
\vspace{10pt}
In particular, this also shows that in the case where $c$ is a cycle of $\pi_r$, the quadrilateral glued to the top right side of $q_i'$ is $q_{\pi_r(i)}'$, while the quadrilateral glued to the top left side of $q_i'$ is $q'_{\pi_\ell\pi_r(i)}$ and the same is true when $c$ is a cycle of $\pi_\ell$ with the roles of $\ell$ and $r$ exchanged. 

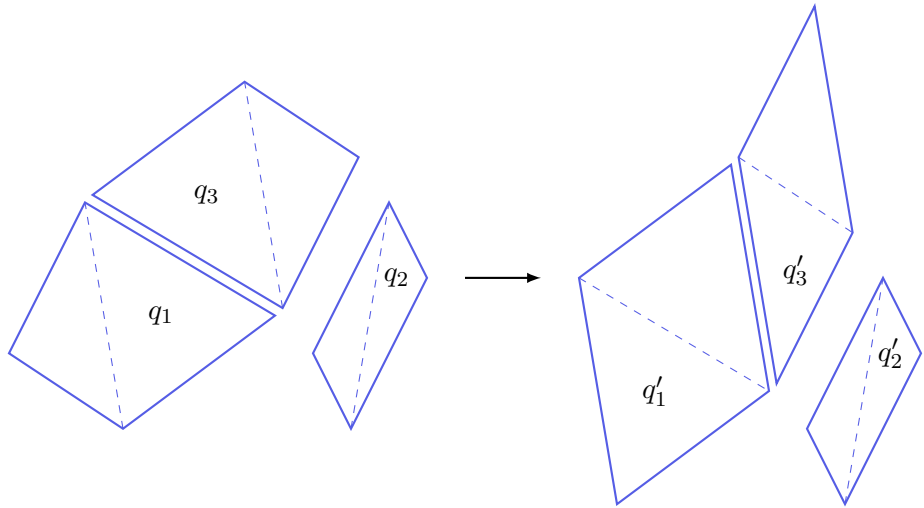
\begin{figure}[ht]
    \centering
    \begin{tikzpicture}[scale = 0.5]
    \draw[thick, noamblue] (0,0) -- (-3,2) -- (-1,6) -- (4,3) -- cycle;
    \draw[noamblue, dashed](0,0) -- (-1,6);
    \draw[noamblue, thick] (-0.8,6.2) -- (4.2,3.2) -- (6.2,7.2) -- (3.2, 9.2) -- cycle;
    \draw[noamblue, dashed] (4.2,3.2) -- (3.2,9.2);
    \draw[thick, noamblue] (5,2) -- (6,0) -- (8,4) -- (7,6) -- cycle;
    \draw[noamblue, dashed] (6,0) -- (7,6);

    \draw[thick, >=latex, ->] (9,4) -- (11,4);

    \def \x{13}; \def \y{-2};

    \draw[thick, noamblue] (\x,\y) -- (\x-1, \y+6) -- (\x+3, \y + 9) -- (\x + 4, \y + 3) -- cycle;
    \draw[dashed, noamblue] (\x-1, \y+6) -- (\x + 4, \y + 3);

    \draw[thick, noamblue] (\x + 4.2, \y + 3.2) -- (\x+3.2, \y + 9.2) -- (\x + 5.2, \y + 13.2) -- (\x+6.2, \y + 7.2) -- cycle;
    \draw[dashed, noamblue] (\x+3.2, \y + 9.2) -- (\x+6.2, \y + 7.2);

    \draw[thick, noamblue] (\x+5,\y+2) -- (\x+6,\y+0) -- (\x+8,\y+4) -- (\x+7,\y+6) -- cycle;
    \draw[noamblue, dashed] (\x+6,\y+0) -- (\x+7,\y+6);

    \node at (1,3) {$q_1$};
    \node at (2.2,6.2) {$q_3$};
    \node at (7.2,4) {$q_2$};

    \node at (\x+1, \y+3) {$q'_1$};
    \node at (\x+4.7, \y+6.2) {$q'_3$};
    \node at (\x+7.2,\y+4) {$q'_2$};
\end{tikzpicture}
    \caption{A right staircase move performed on the quadrangulation from Example \ref{ex:suspension_bipartite}.}
    \label{fig:staircase_move}
\end{figure}

\begin{example}[Action of staircase move]
    Consider the suspension from Example \ref{ex:suspension_bipartite}. The union of the leftmost and rightmost quadrilaterals, which we will denote by $q_1$ and $q_3$, forms a right well-slanted staircase $S$. The associated cycle is given by $c = \pi_r = \begin{pmatrix}
        1 & 3
    \end{pmatrix}$. The action of a staircase move is depicted in Figure \ref{fig:staircase_move}. The new quadrangulation $Q' = (\boldsymbol{\pi}, \boldsymbol{w})$ is given by

    \begin{equation*}
        \pi_\ell = \begin{pmatrix}
            2 & 3
        \end{pmatrix}, \quad
        \pi_r = \begin{pmatrix}
            1 & 3
        \end{pmatrix}
    \end{equation*}
    and
    \begin{align*}
        \boldsymbol{w} &= (w_1', w_2', w_3') \\&= \big((w_{1,d},w_{1, r}), (w_{2, \ell}, w_{2, r}), (w_{3, d}, w_{3, r})\big) \\&=
        \big( (w_{1, \ell} + w_{3, r}, w_{1,r}), (w_{2, \ell}, w_{2, r}), (w_{3,\ell} + w_{3, r}, w_{3,r})\big) \\&=
        \big( (1+6\ii, 4+3\ii), (-1 + 2\ii, 2 + 4\ii), (-1 + 6\ii, 2 + 4\ii).
    \end{align*}
    
\end{example}

Let us summarize the action of a staircase move on the quadrangulation data $Q = (\boldsymbol{\pi}, \boldsymbol{w})$ in the following proposition.

\begin{proposition}[Action of staircase move]\label{prop:action_staircase_move}
    A staircase move acts on the data $(\boldsymbol{\pi}, \boldsymbol{w})$ as follows.

    \begin{enumerate}
        \item $c \in \pi_r$: The new combinatorial datum $\boldsymbol{\pi}' = (\pi'_\ell, \pi'_r)$ is given by
        \begin{align*}
            \pi_\ell'(i) &= 
            \begin{cases}
                \pi_\ell\pi_r(i) \quad &\text{if } i \in c,\\
                \pi_\ell(i) & \text{else},
            \end{cases}\\
            \pi'_r &= \pi_r.
        \end{align*}
        The new length datum $\boldsymbol{w}'$ is given by
        \begin{equation*}
            w'_i =
            \begin{cases}
                (w_{i, d}, w_{i,r}) \quad &\text{if }i \in c,\\
                w_i &\text{else.}
            \end{cases}
        \end{equation*}
        \item $c \in \pi_\ell$: The new combinatorial datum $\boldsymbol{\pi}' = (\pi'_\ell, \pi'_r)$ is given by
        \begin{align*}
            \pi_r'(i) &= 
            \begin{cases}
                \pi_r\pi_\ell(i) \quad &\text{if } i \in c,\\
                \pi_r(i) & \text{else},
            \end{cases}\\
            \pi'_\ell &= \pi_\ell.
        \end{align*}
        The new length datum $\boldsymbol{w}'$ is given by
        \begin{equation*}
            w'_i =
            \begin{cases}
                (w_{i, \ell}, w_{i,d}) \quad &\text{if }i \in c,\\
                w_i &\text{else.}
            \end{cases}
        \end{equation*}
    \end{enumerate}
\end{proposition}

We will write $c \cdot \boldsymbol{\pi}$ for the new combinatorial datum $\boldsymbol{\pi}'$ given by the formulas in Proposition \ref{prop:action_staircase_move}.

\begin{corollary}
    The train-track relations for $\boldsymbol{\pi}'$ are satisfied by $\boldsymbol{w}'$.
\end{corollary}
\begin{proof}
    Let us check the train-track relations if $c$ is a cycle of $\pi_\ell$. If $i \notin c$ the train-track relations remain the same. If $i \in c$  we have
    \begin{align*}
        w'_{i,\ell} + w'_{\pi'_\ell(i), r} &=
        w_{i,\ell} + w'_{\pi_\ell(i), r} \\&=
        w_{i,\ell} + w_{\pi_\ell(i), d} \\&=
        w_{i, \ell} + w_{\pi_\ell(i), r} + w_{\pi_r\pi_\ell(i), \ell} \\&=
        w_{i, r} + w_{\pi_r(i), \ell} + w'_{\pi'_r(i), \ell} \\&=
        w_{i,d} + w'_{\pi'_r(i), \ell} \\&=
        w'_{i,r} + w'_{\pi'_r(i), \ell},
    \end{align*}
    where we use the formulas from Proposition \ref{prop:action_staircase_move}, the definition of the forward diagonal $w_{i,d}$ as well as the train-track relations for $\boldsymbol{w}$.
\end{proof}

It is important to notice that the action on the permutation datum $\boldsymbol{\pi}$ is independent of the length datum $\boldsymbol{w}'$. Furthermore, the action on $\boldsymbol{w}$ is linear. Thus we may collect the information about the action of a staircase move on $\boldsymbol{w}$ into a $2k \times 2k$ matrix $A_{\boldsymbol{\pi}, c}$ in the following way.

We index the rows and columns of the matrix $A_{\boldsymbol{\pi}, c}$ by the indices
\begin{equation*}
    (1, \ell), (1, r), (2, \ell), (2, r), \ldots, (k, \ell), (k, r),
\end{equation*}
which is exactly the alphabet $\mathcal{A}$ from equation \eqref{eq:alphabet_bipartite} we introduced to define bi-partite IETs. We write ${I}_{2k}$ for the $2k \times 2k$ identity matrix and we denote by ${E}_{(i, \varepsilon), (j, \nu)}$ the matrix with a 1 in position $\big((i, \varepsilon), (j, \nu)\big)$ and zeroes everywhere else. 

\begin{definition}[Staircase matrix]\label{def:staircase_matrix}
    We define the matrix $A_{\boldsymbol{\pi}, c}$ to be
    \begin{equation*}
        A_{\boldsymbol{\pi}, c} =
        \begin{cases}
            I_{2k} + \sum_{i \in c} E_{(i, \ell)(\pi_\ell(i), r)} \quad &\text{if } c \text{ is a cycle of } \pi_r,\\
            I_{2k} + \sum_{i \in c} E_{(i, r)(\pi_r(i), \ell)} & \text{if } c \text{ is a cycle of } \pi_\ell.
        \end{cases}
    \end{equation*}
\end{definition}

Viewing $\boldsymbol{w}$ and $\boldsymbol{w'}$ as column vectors, it follows immediately from the definition of matrix-vector multiplication and Proposition \ref{prop:action_staircase_move} that $\boldsymbol{w}' = A_{\boldsymbol{\pi},c} \cdot \boldsymbol{w}$. This establishes the following lemma.

\begin{lemma}[Staircase move on data]\label{lem:staircase_move_on_data}
    Given a labeled quadrangulation $Q = (\boldsymbol{\pi}, \boldsymbol{w})$ and a cycle $c$ of $\boldsymbol{\pi}$, if the staircase $S_c$ is well-slanted, when performing on $Q$ the staircase move in $S_c$ one obtains a new labeled quadrangulation $Q' = (\boldsymbol{\pi}', \boldsymbol{w}')$ with
    \begin{equation*}
        \boldsymbol{\pi}' = c \cdot \boldsymbol{\pi}, \quad \boldsymbol{w}' = A_{\boldsymbol{\pi}, c} \cdot \boldsymbol{w},
    \end{equation*}
    where $c \cdot \boldsymbol{\pi}$ and $A_{\boldsymbol{\pi}, c}$ are given by the formulas in Proposition \ref{prop:action_staircase_move} and Definition \ref{def:staircase_matrix}. 
\end{lemma}

One might ask if it was possible to define such moves on smaller parts of a quadrangulation than on staircases. This is not the case, i.e., one can show that staircases are the smallest unions of quadrilaterals in which one can simultaneously perform diagonal changes to get a new quadrangulation $Q'$ in the sense made precise in the lemma below. Therefore, staircases are the correct \enquote{irreducible} objects to consider.

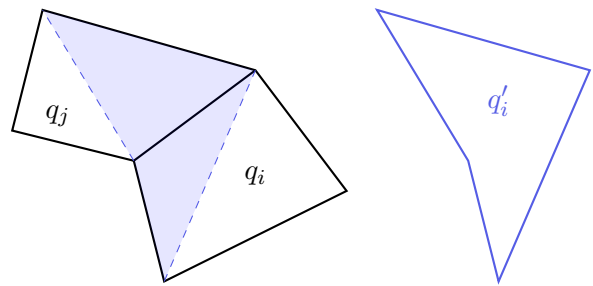
\begin{figure}[ht]
    \centering
    \begin{tikzpicture}[scale = 0.4]
    \draw[thick] (0,0) -- (-1,4) -- (-5,5) -- (-4,9) -- (3,7) -- (6,3) -- cycle;
    \draw[thick] (-1,4) -- (3,7);

    \draw[dashed, noamblue] (0,0) -- (3,7); \draw[dashed, noamblue] (-1,4) -- (-4,9);
    \fill[blue, opacity = 0.1] (0,0) -- (3,7) -- (-4,9) -- (-1,4) -- cycle;

    \node at (3,3.5) {$q_i$};
    \node at (-3.5,5.5) {$q_j$};

    \def \x{11};
    \draw[thick, noamblue] (\x+0,0) -- (\x+3,7) -- (\x-4,9) -- (\x-1,4) -- cycle;
    \node[noamblue] at (\x, 6) {$q'_i$};
\end{tikzpicture}
    \caption{Performing diagonal changes on unions of quadrilaterals that are not staircases does not give a quadrangulation.}
    \label{fig:no_decomposition_not_admissible}
\end{figure}

\begin{lemma}[Irreducibility of staircases]\label{lem:irreducibility_of_staircases}
    Let $Q = (\boldsymbol{\pi}, \boldsymbol{w})$ be a quadrangulation and let $\mathcal{I}_\ell, \mathcal{I}_r \subseteq [d]$ be such that the quadrilaterals $q_i$ with $i \in \mathcal{I}_\ell$ are left-slanted and the quadrilaterals $q_i$ with $i \in \mathcal{I}_r$ are right-slanted.

    The new set of wedges obtained after individual diagonal changes in the quadrilaterals $q_i$ for $i \in \mathcal{I}_\ell \cup \mathcal{I}_r$ is associated to a quadrangulation if and only if the set of indices $\mathcal{I}_\ell$ (or $\mathcal{I}_r$) is a union of cycles of $\pi_\ell$ (or $\pi_r$).
\end{lemma}
\begin{proof}
    If the set of indices of $\mathcal{I}_\ell \cup \mathcal{I}_r$ are unions of cycles, performing staircase moves yields a new quadrangulation as we have shown above. 

    Conversely, if, say, $\mathcal{I}_\ell$ is \emph{not} the union of cycles in $\pi_\ell$, then there exists $i \in \mathcal{I}_\ell$ such that $j = \pi_\ell(i) \notin \mathcal{I}_\ell$, so consequently $j\in \mathcal{I}_r$ hence $q_j$ is right-slanted. The quadrilaterals $q_i$ and $q_j$ are glued as on the left of Figure \ref{fig:no_decomposition_not_admissible}. Performing diagonal changes individually in these two quadrilaterals necessarily yields a quadrilateral $q'_i$ as pictured on the right of Figure \ref{fig:no_decomposition_not_admissible}, which evidently cannot be admissible again so we do not obtain a proper quadrangulation.
\end{proof}

\subsection{Diagonal Changes Algorithms Based on Staircase Moves}

It remains to use the staircase moves described in the previous section to define an algorithm, that is, a procedure that given some starting quadrangulation $Q = Q^{(0)}$ produces a sequence of quadrangulations $Q^{(1)}, Q^{(2)}, \ldots$ such that $Q^{(n+1)}$ is obtained by $Q^{(n)}$ by a certain number of staircase moves. Exactly as in the case of Rauzy--Veech induction, there are several possible ways to define such an algorithm. 

The first remark we want to make is that if $S_1$ and $S_2$ are two \emph{disjoint} well-slanted staircases in $Q$, then staircase moves in $S_1$ and $S_2$ clearly commute, so that the order in which they are performed is irrelevant and the moves can be considered to be performed simultaneously.

\begin{definition}
    If the quadrangulation $Q'$ is obtained from the quadrangulation $Q$ by performing staircase moves (possibly just on a subset of the well-slanted staircases of $Q$), we will say that $Q'$ is obtained from $Q$ by \emph{simultaneous staircase moves}.
\end{definition}

We will describe two algorithms here which correspond exactly to the slow and the fast Rauzy--Veech induction algorithms described in section \ref{sec:rauzy_veech}. The first algorithm corresponds to the algorithm introduced in \cite{ferenczi2010structure} for bipartite IETs.

\begin{definition}[Greedy algorithm]\label{def:greedy_algorithm}
The \emph{greedy diagonal changes algorithm} starting from the quadrangulation $Q = Q^{(0)}$ produces the sequence $\left(Q^{(n)}\right)_{n \in \N}$ of quadrangulations, where $Q^{(n+1)}$ is obtained from $Q^{(n)}$ by performing simultaneous staircase moves in all well-slanted staircases for $Q^{(n)}$. 
\end{definition}

In order to define the second algorithm, the following remark is important. Performing a left (or right) staircase move does not modify $\pi_\ell$ (or $\pi_r$). So even if the quadrilaterals change under the staircase move, the staircases itself, seen as subsets of the surface $X$, remain the same. This justifies the following definition.

\begin{definition}[Multiplicity of staircases]
    The \emph{multiplicity of a left (or right) staircase} $S_c$ is the maximal $n$, such that we can perform $n$ consecutive left (or right) staircase moves in $S_c$. We denote the multiplicity by $\operatorname{mul}(S_c)$.
\end{definition}

\begin{definition}[Left/right algorithm]
The \emph{left/right diagonal changes algorithm} starting from the quadrangulation $Q = Q^{(0)}$ produces a sequence $\left(Q^{(n)}\right)_{n \in \N}$, where, if $n$ is even, $Q^{(n+1)}$ is obtained from $Q^{(n)}$ by performing $\operatorname{mul}(S_c)$ staircase moves in any left-slanted staircase $S_c$ in $Q$ and if $n$ is odd, $Q^{(n+1)}$ is obtained from $Q^{(n)}$ by performing $\operatorname{mul}(S_c)$ staircase moves in any right-slanted staircase $S_c$ in $Q$.
\end{definition}

\begin{remark}
    We may see the left/right algorithm as an acceleration of the greedy algorithm. Also, in the left/right algorithm the succeeding quadrangulation $Q^{(n+1)}$ is in general obtained from $Q^{(n)}$ by several staircase moves that are \emph{not} successive. 
\end{remark}

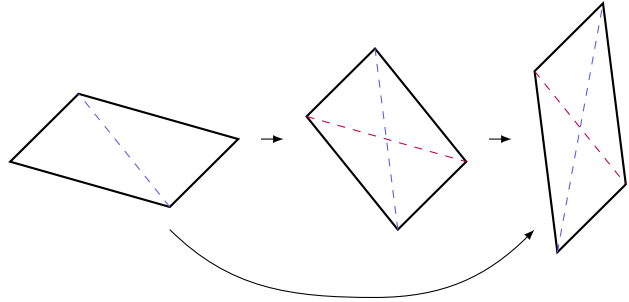
\begin{figure}[ht]
    \centering
    \begin{tikzpicture}[scale = 0.3]
    \draw[thick] (0,0) -- (-7,2) -- (-4,5) -- (3,3) -- cycle;
    \draw[dashed, noamblue] (0,0) -- (-4,5);

    \draw[>=latex, ->] (4,3) -- (5,3);

    \def \x{10}; \def \y{-1};
    \draw[thick] (\x, \y) -- (\x - 4, \y + 5) -- (\x - 1, \y + 8) -- (\x +3, \y+3) -- cycle;
    \draw[dashed, noamblue] (\x, \y) -- (\x - 1, \y + 8);
    \draw[dashed, purple] (\x - 4, \y + 5) -- (\x +3, \y+3);

    \draw[>=latex, ->] (\x+4, 3) -- (\x+5, 3);

    \def \z{17}; \def \w{-2};

    \draw[thick] (\z, \w) -- (\z-1, \w+8) -- (\z+2, \w+11) -- (\z+3, \w+3) -- cycle;
    \draw[dashed, noamblue] (\z, \w) -- (\z+2, \w+11);
    \draw[dashed, purple] (\z-1, \w+8) -- (\z+3, \w+3);

    \draw[>=latex, ->] (0,-1) to[out=-45, in=180] (\x - 1, -4) to[out = 0, in = 225] (\z-1, -1);
    
\end{tikzpicture}
    \caption{Two steps of the greedy algorithm performed on the translation surface given by a single quadrilateral corresponds to one application of the left/right algorithm. The forward diagonal is depicted in blue, the backwards diagonal in red.}
    \label{fig:greedy_left_right_algorithm}
\end{figure}

Figure \ref{fig:greedy_left_right_algorithm} illustrates the difference between these two algorithms in the simplest case, where the surface $X$ is given by a single quadrilateral $q$, for which we have $\operatorname{mul}(q) = 2$. Therefore, one application of the left/right algorithm amounts to two applications of the greedy algorithm.

Not every algorithm based on diagonal changes is induced by \emph{simultaneous} staircase moves. It is possible to define algorithms where the new quadrangulation $Q'$ is obtained by not only applying staircase moves to the initial quadrangulation $Q$, but also to intermediate quadrangulations. We will see one such algorithm in section \ref{sec:section_for_teichmueller_flow}. This algorithm, called \emph{geodesic algorithm} in \cite{delecroix2015diagonal}, is defined as the Poincaré first return map of the Teichmüller geodesic flow $g_t$ to some specified Poincaré section. 

\subsection{Parameter Space and Nice Properties of Diagonal Changes}\label{sec:parameter_space}

Our goal in this section is to introduce the space of (labeled) quadrangulations of surfaces on which the algorithms defined in the end of the previous section act. Moreover, we will show that the staircase moves are invertible and self-dual, meaning that their inverses are again given by staircase moves and we will also explain, that the parameter space exhibits a Markovian structure, i.e., there is a \emph{loss of memory} phenomenon. Lastly, we will extend the definition of best approximations to the setting of translation surfaces of higher genus and show that algorithms performing simultaneous staircase moves are able to detect systoles that are realized along a Teichmüller geodesic.

As we will mention below (see Theorem \ref{thm:geometric_objects_are_the_same}), any algorithm where the new quadrangulation is obtained from the old by simultaneous staircase moves produces the same geometric objects (wedges or quadrilaterals). Thus, the actual choice of algorithm is insubstantial for these kind of questions, as long as we employ only \emph{simultaneous} staircase moves. 

Let us stress that the construction of the parameter space is done for surfaces in a hyperelliptic stratum, a concept we will explain in section \ref{sec:existence_of_quadrangulations}. We start by introducing so-called \emph{diagonal changes classes}, a direct analog to Rauzy classes (see Definition \ref{def:rauzy_class}.

\begin{definition}[Diagonal change classes]\label{def:diagonal_change_classes}
    A collection $\mathfrak{DC}$ of pairs of permutations $\boldsymbol{\pi}^j = (\pi_\ell^j, \pi_r^j)$ is called a \emph{diagonal change class}, or \emph{DC class} for short, if for any $\boldsymbol{\pi}^i \neq \boldsymbol{\pi}^j \in \mathfrak{DC}$, there exists a finite sequence of staircase moves from a quadrangulation with combinatorial datum $\boldsymbol{\pi}^i$ to a surface with combinatorial datum $\pi^j$. 
\end{definition}

Given any combinatorial datum $\boldsymbol{\pi}$, we can represent the $DC$ class to which it belongs as a directed graph $\mathcal{G}(\boldsymbol{\pi}) = \mathcal{G}$, by letting the vertices be given by the set of combinatorial data obtainable from $\boldsymbol{\pi}$ by a sequence of staircase moves and letting there be an edge from $\boldsymbol{\pi}$ to $\boldsymbol{\pi}'$ if and only if $c \cdot \boldsymbol{\pi} = \boldsymbol{\pi}'$. As was the case for Rauzy classes in section \ref{sec:rauzy_veech}, such a graph is independent of the chosen starting datum. We will refer to such a graph $\mathcal{G}$ as a $DC$ graph.

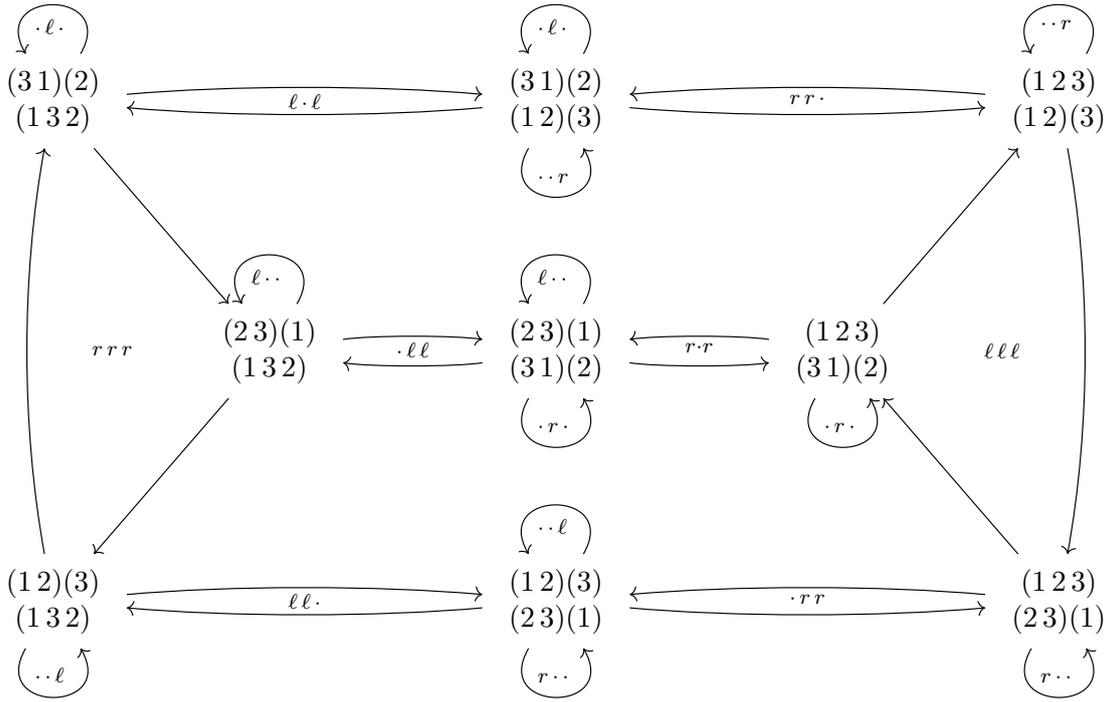
\begin{figure}[ht]
\[\begin{tikzcd}[row sep=large]
    {\begin{array}{c} (3 \,1)(2)\\(1 \, 3 \,2) \end{array}}
    \arrow[loop above, distance = 1cm, in=120, out=60, "\cdot\,\ell\,\cdot" {pos=0.7}] 
    &&& {\begin{array}{c} (3\,1)(2)\\(1\,2)(3) \end{array}}
    \arrow[loop above, distance = 1cm, in=120, out=60, "\cdot\,\ell\,\cdot" {pos=0.7}] 
    \arrow[loop, distance = 1cm, in=300, out=240, "\cdot\,\cdot\,r" {pos=0.3}] 
    &&& {\begin{array}{c}(1\,2\,3)\\(1\,2)(3) \end{array}}
    \arrow[loop above, distance = 1cm, in=120, out=60, "\cdot\,\cdot\,r" {pos=0.7}] 
    \\
    \\
    & {\begin{array}{c}(2 \, 3)(1)\\(1 \, 3 \,2) \end{array}}
    \arrow[loop above, distance = 1cm, in=120, out=60, "\ell\,\cdot\,\cdot" {pos=0.7}] 
    && {\begin{array}{c}(2\,3)(1)\\(3\,1)(2) \end{array}}
    \arrow[loop above, distance = 1cm, in=120, out=60, "\ell\,\cdot\,\cdot" {pos=0.7}]\arrow[loop, distance = 1cm, in=300, out=240, "\cdot\,r\,\cdot" {pos=0.3}] 
    && {\begin{array}{c}(1\,2\,3)\\(3\,1)(2) \end{array}}\arrow[loop, distance = 1cm, in=300, out=240, "\cdot\,r\,\cdot" {pos=0.3}] \\
    \\
    {\begin{array}{c}(1\,2)(3)\\(1\,3\,2) \end{array}}\arrow[loop, distance = 1cm, in=300, out=240, "\cdot\,\cdot\,\ell" {pos=0.3}]
    &&& {\begin{array}{c}(1\,2)(3)\\(2\,3)(1) \end{array}}
    \arrow[loop above, distance = 1cm, in=120, out=60, "\cdot\,\cdot\,\ell" {pos=0.7}]\arrow[loop, distance = 1cm, in=300, out=240, "r\,\cdot\,\cdot" {pos=0.3}] 
    &&& {\begin{array}{c}(1\,2\,3)\\(2\,3)(1) \end{array}}\arrow[loop, distance = 1cm, in=300, out=240, "r\,\cdot\,\cdot" {pos=0.3}]
    \arrow["\ell\,\cdot\,\ell"', curve={height=-6pt}, from=1-1, to=1-4]
    \arrow[curve={height=-6pt}, from=1-4, to=1-1]
    \arrow["r\,r\,\cdot", curve={height=6pt}, from=1-7, to=1-4]
    \arrow[curve={height=6pt}, from=1-4, to=1-7]
    \arrow[from=3-6, to=1-7]
    \arrow[from=1-1, to=3-2]
    \arrow["{\cdot\,\ell\,\ell}"', curve={height=-6pt}, from=3-2, to=3-4]
    \arrow[curve={height=-6pt}, from=3-4, to=3-2]
    \arrow["{r \cdot r}", curve={height=6pt}, from=3-6, to=3-4]
    \arrow[curve={height=6pt}, from=3-4, to=3-6]
    \arrow[from=5-7, to=3-6]
    \arrow[from=3-2, to=5-1]
    \arrow["{\ell\,\ell\,\cdot}"', curve={height=-6pt}, from=5-1, to=5-4]
    \arrow[curve={height=-6pt}, from=5-4, to=5-1]
    \arrow["{\cdot\, r\,r}", curve={height=6pt}, from=5-7, to=5-4]
    \arrow[curve={height=6pt}, from=5-4, to=5-7]
    \arrow["{\qquad r\,r\,r}"', curve={height=-12pt}, from=5-1, to=1-1]
    \arrow["{\ell\,\ell\,\ell\qquad}"', curve={height=-12pt}, from=1-7, to=5-7]
\end{tikzcd}\]

    \caption{The $DC$ graph $\mathcal{G}$ of combinatorial data for quadrangulations in the stratum $\mathcal{H}(2)$ from Example \ref{ex:dc_graph}.}
    \label{fig:dc_graph}
\end{figure}

\begin{example}[DC graph]\label{ex:dc_graph}
    Figure \ref{fig:dc_graph} depicts a $DC$ graph $\mathcal{G}$ of combinatorial data for quadrangulations in the stratum $\mathcal{H}(2)$, i.e., the stratum consisting of translation surfaces with one conical singularity of total cone angle $6\pi$.

    The notation used for the cycles is as follows. If $c$ is a cycle of $\pi_\ell$ (or $\pi_r$), then we write it as a word of length $k$ on the alphabet $\{\cdot, \ell\}$ (or $\{\cdot, r\}$), where the $i\textsuperscript{th}$ letter of the word is $\ell$ (or $r$) if and only if $i \in c$. Here, $k$ denotes again the number such that the associated bi-partite IET has $2k$ intervals. So for example, if $k = 5$ and $c = \begin{pmatrix}
        1 & 3 & 4
    \end{pmatrix} \in \pi_\ell$, then we write $c = \ell\cdot\ell\,\ell\,\cdot$, and if $c = \begin{pmatrix}
        2 & 3
    \end{pmatrix} \in \pi_r$ we write $c = \cdot\, r\, r \cdot \cdot$.
\end{example}

It may happen that, starting with a different combinatorial datum $\boldsymbol{\pi}'$ which is not part of the graph $\mathcal{G}(\boldsymbol{\pi})$, we obtain a distinct graph $\mathcal{G}(\boldsymbol{\pi}')$. However, these two graphs will be isomorphic. In particular, there exists a permutaion $\sigma$ such that the isomorphism is given by $(\pi_\ell, \pi_r) \mapsto (\sigma\pi_\ell\sigma^{-1}, \sigma\pi_r\sigma^{-1})$, see \cite{delecroix2015diagonal}.

For each combinatorial datum $\boldsymbol{\pi}$ let us now introduce the parameter spaces $\Delta_{\boldsymbol{\pi}}$ and $\Theta_{\boldsymbol{\pi}}$, which are all possible length and suspension data satisfying the corresponding train-track relations. That is, we define
\begin{align*}
    \begin{split}
    \Delta_{\boldsymbol{\pi}} &= \left\{\big( (\lambda_{1, \ell}, \lambda_{1, r}), \ldots, (\lambda_{k, \ell}, \lambda_{k, r})\big)\in (\R_-\times\R_+)^k \mid \lambda_{i, \ell} + \lambda_{\pi_\ell(i), r} = \lambda_{i,r} + \lambda_{\pi_r(i), \ell},\, 1 \leq i \leq k\right\},\\
    \Theta_{\boldsymbol{\pi}} &= \left\{\big( (\tau_{1, \ell}, \tau_{1, r}), \ldots, (\tau_{k, \ell}, \tau_{k, r})\big)\in (\R_+\times\R_+)^k \mid \tau_{i, \ell} + \tau_{\pi_\ell(i), r} = \tau_{i,r} + \tau_{\pi_r(i), \ell},\, 1 \leq i \leq k\right\}.
    \end{split}
\end{align*}

\begin{definition}[Space of labeled quadrangulations]\label{def:space_labeled_quadrangulations}
    The \emph{space of labeled quadrangulations} of surfaces corresponding to a bi-partite IET on $2k$ intervals is given by
    \begin{equation*}
        \mathcal{Q}_k = \{(\boldsymbol{\pi}, \boldsymbol{\lambda}, \boldsymbol{\tau}) \mid \boldsymbol{\pi} \in \mathcal{G}, \boldsymbol{\lambda} \in \Delta_{\boldsymbol{\pi}}, \boldsymbol{\tau} \in \Theta_{\boldsymbol{\pi}}\}.
    \end{equation*}
\end{definition}

Given $(\boldsymbol{\pi}, \boldsymbol{\lambda}, \boldsymbol{\tau}) \in \mathcal{Q}_k$ and $c \in \pi_\varepsilon$, notice that the heights $\tau \in \Theta_{\boldsymbol{\pi}}$ play no role in determining whether $S_c$ is well-slanted. Let $\Delta_{\boldsymbol{\pi}, c} \subseteq \Delta_{\boldsymbol{\pi}}$ be the subset of length data for which $S_c$ s well-slanted, i.e., 
\begin{equation*}
    \Delta_{\boldsymbol{\pi},c} = 
    \begin{cases}
        \{\boldsymbol{\lambda} \in \Delta_{\boldsymbol{\pi}} \mid \lambda_{i,d} < 0 \, \text{ for all }i \in c\} \quad \text{ if } c \in \pi_r, \\
        \{\boldsymbol{\lambda} \in \Delta_{\boldsymbol{\pi}}\mid \lambda_{i,d} > 0 \, \text{ for all }i \in c\}\quad \text{ if }c \in \pi_\ell.
    \end{cases}
\end{equation*}
Here we write $\lambda_{i,d} = \RE(w_{i,d})$ for the real part of the forward diagonal of $q_i$. In this new language, one can perform a staircase move in $S_c$ if and only if $\boldsymbol{\lambda} \in \Delta_{\boldsymbol{\pi},c}$. Since the matrix $A_{\boldsymbol{\pi},c}$ (see Definition \ref{def:staircase_matrix}) acts linearly on both the real and imaginary part of the wedges in $\boldsymbol{w}$, we can define a staircase move on the parameter space as follows.

\begin{definition}[Staircase move on parameter space]
    Let $\boldsymbol{\pi} = (\pi_\ell, \pi_r) \in \mathcal{G}$ and let $c$ be a cycle of $\pi_\varepsilon$. The \emph{staircase move} $\hat{m}_{\boldsymbol{\pi},c}$ on $\{\boldsymbol{\pi}\}\times \Delta_{\boldsymbol{\pi},c} \times \Theta_{\boldsymbol{\pi}}\subseteq\mathcal{Q}_k$ is the map
    \begin{align*}
        \hat{m}_{\boldsymbol{\pi},c} \colon \{\boldsymbol{\pi}\}\times \Delta_{\boldsymbol{\pi},c} \times \Theta_{\boldsymbol{\pi}} &\to \mathcal{Q}_k, \\
        (\boldsymbol{\pi},\boldsymbol{\lambda}, \boldsymbol{\tau}) &\mapsto (c\cdot \boldsymbol{\pi}, A_{\boldsymbol{\pi}, c} \cdot \boldsymbol{\lambda}, A_{\boldsymbol{\pi}, c} \cdot \boldsymbol{\tau}).
    \end{align*}
\end{definition}

\subsubsection{Invertibility}
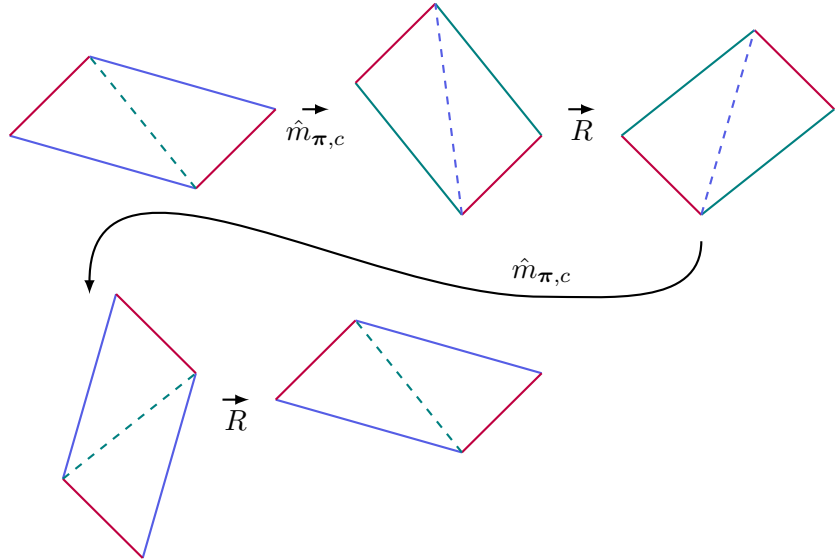
\begin{figure}[ht]
    \centering
    \begin{tikzpicture}[scale =0.35]
    \draw[thick, noamblue] (0,0) -- (-7,2); \draw[thick, noamblue] (3,3) -- (-4,5);
    \draw[thick, purple] (-7,2) -- (-4,5); \draw[thick, purple] (0,0) -- (3,3);
    \draw[thick, teal, dashed] (0,0) -- (-4,5);

    \draw[thick, >=latex, ->] (4,3) -- (5,3) node[midway, below] {$\hat{m}_{\boldsymbol{\pi}, c}$};

    \def \p2{10}; \def \q2{-1};
    \draw[thick, teal] (\p2, \q2) -- (\p2 - 4, \q2 + 5); 
    \draw[thick, teal] (\p2+3, \q2+3) -- (\p2 - 1, \q2 + 8);
    \draw[thick, purple] (\p2, \q2) -- (\p2 + 3, \q2 + 3); 
    \draw[thick, purple] (\p2 -4, \q2+5) -- (\p2 - 1, \q2 + 8);
    \draw[noamblue, thick, dashed] (\p2, \q2) -- (\p2 - 1, \q2 + 8);

    \draw[thick, >=latex, ->] (\p2 + 4, 3) -- (\p2 + 5, 3) node[midway, below] {$R$};

    \def \p3{19}; \def \q3{-1};
    \draw[thick, purple] (\p3, \q3) -- (\p3 - 3, \q3 + 3); 
    \draw[thick, purple] (\p3 + 5, \q3 + 4) -- (\p3 + 2, \q3 + 7);
    \draw[thick, teal] (\p3, \q3) -- (\p3 + 5, \q3 + 4);
    \draw[thick, teal] (\p3 - 3, \q3 + 3) -- (\p3 + 2, \q3 + 7);
    \draw[noamblue, thick, dashed] (\p3, \q3) -- (\p3 + 2, \q3 + 7);

    \draw[thick, >=latex, ->] (\p3, \q3 - 1) to[out=270, in = 0] (13, -4.1) to[out = 180, in = 90] (-4,-4);
    \node[above] at (13,-4.1){$\hat{m}_{\boldsymbol{\pi}, c}$};

    \def \p4{-2}; \def \q4{-14};
    \draw[thick, noamblue] (\p4,\q4) -- (\p4 + 2,\q4 + 7); \draw[thick, noamblue] (\p4-3,\q4+3) -- (\p4 -1,\q4 + 10);
    \draw[thick, purple] (\p4,\q4) --  (\p4-3,\q4+3); \draw[thick, purple] (\p4 + 2,\q4 + 7) -- (\p4 -1,\q4 + 10);
    \draw[thick, teal, dashed] (\p4-3,\q4+3) -- (\p4 + 2,\q4 + 7);    
    \draw[thick, >=latex, ->] (\p4 + 3, \q4+6) -- (\p4 + 4, \q4+6) node[midway, below] {$R$};
    
    \def \p5{10}; \def \q5{-10};
    \draw[thick, noamblue] (\p5+0,\q5+0) -- (\p5-7,\q5+2); \draw[thick, noamblue] (\p5+3,\q5+3) -- (\p5-4,\q5+5);
    \draw[thick, purple] (\p5-7,\q5+2) -- (\p5-4,\q5+5); \draw[thick, purple] (\p5+0,\q5+0) -- (\p5+3,\q5+3);
    \draw[thick, teal, dashed] (\p5+0,\q5+0) -- (\p5-4,\q5+5);
\end{tikzpicture}
    \caption{A staircase move and its inverse applied to a staircase consisting of a single quadrilateral.}
    \label{fig:inverse_example}
\end{figure}

Our next goal is to show that staircase moves are invertible where their inverses are again given by staircase moves. If we consider a staircase consisting of a single quadrilateral, it is easy to visualize how we find the inverse of a staircase move: We rotate the quadrilateral by $\frac{\pi}{2}$, apply the appropriate staircase move to the rotated quadrilateral (which will be opposite of the first staircase move) and rotate back by $-\frac{\pi}{2}$. An example can be seen in Figure \ref{fig:inverse_example}.

Let us formalize this idea. To do this, we introduce the rotation operator $R$ on $\mathcal{Q}_k$ which rotates quadrilaterals by $\frac{\pi}{2}$. First, we will investigate how $R$ acts on on the different data. The only aspect that warrants caution is the fact that we work with \emph{labeled} quadrangulations, so that we need to be careful to label the rotated quadrangulation in a well-defined manner.

Note first that if $q \subseteq \C$ is an admissible quadrilateral, then $\ii q$, i.e., the quadrilateral obtained by multiplying the wedges by the imaginary unit $\ii$ or equivalently, the quadrilateral obtained by rotating by $\frac{\pi}{2}$, is still admissible. It follows that for any quadrangulation $Q$ of $X$, the collection 
\begin{equation*}
    \ii Q = \{\ii q \mid q \in \Q\}
\end{equation*}
is a quadrangulation as well, namely of $\ii X$ which is the translation surface obtained from $X$ by rotating all polygons (as in Definition \ref{def:translation_surface}) by $\frac{\pi}{2}$. 

\begin{definition}[Rotated labeled quadrangulation]\label{def:rotated_labeled_quadrangulation}
    Given a quadrangulation $Q$, we denote by $Q'$ the quadrangulation $\ii Q$ labeled so that the wedge $w_i'$ of the quadrilateral $q_i' \in Q'$ contains the same vertical ray which was contained in $w_i$ of $q_i\in Q$.
\end{definition}

\begin{figure}[ht]
    \centering
    \begin{tikzpicture}[scale = 1]
    \begin{scope}
        \draw[thick, teal, >=latex, ->] (0,0) -- (-2,1.5) node[midway, below left] {$w_{i,\ell}$};
        \draw[thick, teal, >=latex, ->] (0,0) -- (1.5,3) node[pos=0.3, right] {$w_{i,r}$};

        \draw[thick, purple, >=latex, ->] (1,-1.5) -- (0,0) node[midway, left] {$w_{\pi_\ell^{-1}(i), \ell}$};

        \draw[dotted] (1,-1.5) -- (1.5,3) -- (-2, 1.5);
        \draw[dashed] (0,0) -- ++(90:3);

        \node at (-0.75,1.25) {$q_i$};
        \node at (0.68,0) {$q_{\pi_\ell^{-1}(i)}$};
    \end{scope}
    \draw[>=latex, ->, thick] (2,1) -- (3,1) node[midway, above] {rotate} node[midway, below] {by $\nicefrac{\pi}{2}$};
    
    \begin{scope}[rotate = 90, xshift = 30, yshift = -185]
        \draw[thick, teal, >=latex, ->] (0,0) -- (-2,1.5) ;
        \draw[thick, teal, >=latex, ->] (0,0) -- (1.5,3) ;

        \draw[thick, purple, >=latex, ->] (1,-1.5) -- (0,0);

        \draw[dotted] (1,-1.5) -- (1.5,3) -- (-2, 1.5);
        \draw[dashed] (0,0) -- ++(00:1.5);

        \node at (-0.35,1.25) {$q_i$};
        \node at (0.75,0.55) {$q_{\pi_\ell^{-1}(i)}$};
    \end{scope}
    \draw[>=latex, ->, thick] (8,1) -- (9,1) node[midway, above] {relabel};

    \begin{scope}[rotate = 90, xshift = 30, yshift = -350]
        \draw[thick, black, >=latex, ->] (0,0) -- (-2,1.5);
        \draw[thick, noamblue, >=latex, ->] (0,0) -- (1.5,3) node[midway, below] {$w'_{i, \ell}$} ;

        \draw[thick, noamblue, >=latex, ->] (0,0) -- (1,-1.5) node[midway, below] {$w'_{i,r}$};

        \draw[dotted] (1,-1.5) -- (1.5,3) -- (-2, 1.5);
        \draw[dashed] (0,0) -- ++(00:1.5);

        \node at (0.75,0.55) {$q'_i$};
    \end{scope}
\end{tikzpicture}
    \caption{The quadrangulation $Q'$ obtained by rotation and relabeling.}
    \label{fig:correct_relabeling}
\end{figure}
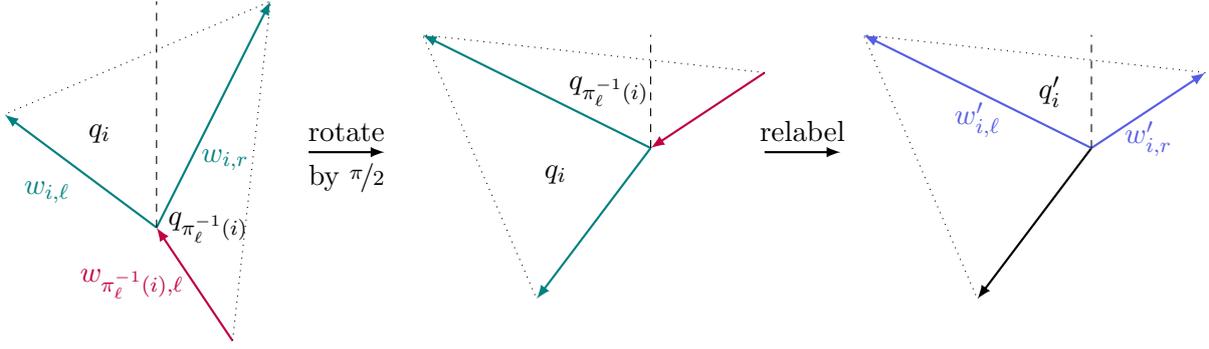

Let us stress that the labeling in Definition \ref{def:rotated_labeled_quadrangulation} is \emph{not} obtained by calling $q'_i$ the quadrilateral $\ii q_i$. An example of a correct relabeling can be seen in Figure \ref{fig:correct_relabeling}, where for clarity the rotation and relabeling are split up into separate steps.

We now work out the explicit formulas for the wedges and combinatorial datum of $q' \in Q'$ obtained from $Q = (\boldsymbol{\pi}, \boldsymbol{w})$. We begin with the action on the combinatorial datum $\boldsymbol{\pi}$.

\begin{enumerate}
    \item The top right side of $q_i'$ is the same up to rotation as the bottom right side of $q_{\pi_\ell}^{-1}(i)$, which is glued to the top left of $q_{\pi_\ell^{-1}\pi_\ell^{-1}}(i) = q'_{\pi_\ell^{-1}(i)}$. Therefore it follows that
    \begin{equation}\label{eq:rotation_1}
        \pi'_r = \pi_\ell^{-1}.
    \end{equation}
    \item The top left side of $q'_i$ is the same up to rotation as the top right side of $q_{\pi_\ell(i)}$, which is glued to the bottom left of $q_{\pi_r\pi_\ell^{-1}(i)} = q'_{\pi_\ell\pi_r\pi_\ell^{-1}(i)}$. From this it follows that
    \begin{equation}\label{eq:rotation_2}
        \pi'_\ell = \pi_\ell\pi_r\pi_\ell^{-1}.
    \end{equation}
\end{enumerate}
The action on the quadrilaterals can be immediately deduced by examining Figure \ref{fig:correct_relabeling} and is as follows.
\begin{align}\label{eq:rotation_3}
    \begin{split}
    q'_i &= \ii \cdot q_{\pi_\ell^{-1}(i)}, \\
    w'_{i, \ell} &= \ii \cdot w_{i,r}, \\
    w'_{i, r} &= \ii \cdot w_{\pi_\ell^{-1}(i), \ell}.
    \end{split}
\end{align}

Thus we obtain the following definition.

\begin{definition}[Rotation operator]\label{def:rotation_operator}
    The rotation operator $R$ sends $Q = (\boldsymbol{\pi}, \boldsymbol{w})\in \mathcal{Q}_k$ to $RQ = (\boldsymbol{\pi}', \boldsymbol{w}')$ given by the formulas \eqref{eq:rotation_1}, \eqref{eq:rotation_2} and \eqref{eq:rotation_3}. 
\end{definition}

In order to check that $R$ is well-defined, we will need to check the following two points.

\begin{enumerate}
    \item Since all $q_i \in Q$ are admissible by definition of quadrangulation, it is clear that $w'$ is a vector of wedges as well, that also satisfies the train-track relations for $\boldsymbol{\pi}'$.

    Hence, if we write $\boldsymbol{w}' = \boldsymbol{\lambda}' + \ii \cdot \boldsymbol{\tau}'$, it follows that
    \begin{equation*}
        \boldsymbol{\lambda}' \in \Delta_{\boldsymbol{\pi}'} \quad \text{and} \quad \boldsymbol{\tau}' \in \Theta_{\boldsymbol{\pi}'}.
    \end{equation*}
    \item Since $\mathcal{Q}_k = \mathcal{G} \times \Delta_{\boldsymbol{\pi}'} \times \Theta_{\boldsymbol{\pi}'}$, it only remains to show that
    \begin{equation*}
        \boldsymbol{\pi}' = (\pi_\ell\pi_r\pi_\ell^{-1}, \pi_\ell^{-1})
    \end{equation*}
    belongs to the same graph $\mathcal{G} = \mathcal{G}(\boldsymbol{\pi})$. This is shown below in Corollary \ref{cor:rotation_well_defined}.
\end{enumerate}

\begin{proposition}[Invertibility of rotation operator]\label{prop:rotation_invertible}
    The operator $R$ is invertible.
\end{proposition}
\begin{proof}
    The inverse $R^{-1}\colon \mathcal{Q}_k \to \mathcal{Q}_k$ is given by $(\boldsymbol{\pi}', \boldsymbol{w}') = R^{-1}(\boldsymbol{\pi}, \boldsymbol{w})$, where 

    \begin{align*}
        \pi'_\ell &= \pi_r^{-1}, \\ \pi'_r &= \pi_r\pi_\ell\pi_r^{-1},
    \end{align*}
    and
    \begin{align*}
        w'_{i,\ell} &= \ii \cdot w_{\pi_r^{-1},r}, \\ w'_{i,r} &= -\ii \cdot w_{i, \ell}.
    \end{align*}
\end{proof}

\begin{remark}
    The rotation operator $R$ exchanges the roles of $\boldsymbol{\lambda}$ and $\boldsymbol{\tau}$, i.e., if $(\boldsymbol{\pi}', \boldsymbol{w}') = R(\boldsymbol{\pi}, \boldsymbol{w})$, then
    \begin{align*}
        w'_{i,\ell} &= -\tau_{i, r} + \ii \cdot \lambda_{i,r}, \quad \text{ and }\\
        w'_{i,r} &= \tau_{\pi_\ell^{-1}(i), \ell} - \ii \cdot \lambda_{\pi_\ell^{-1}(i), \ell}.
    \end{align*}
\end{remark}

Up to now we have only considered the \emph{forward} diagonal $w_{i,d}$ of a quadrilateral $q_i \in Q= (\boldsymbol{\pi}, \boldsymbol{w})$. Let us write $w_{i,d^+} = w_{i,d}$ for clarity, to clearly distinguish the forward diagonal from the \emph{backward} diagonal, which we define next.

\begin{figure}[ht]
    \centering
    \begin{tikzpicture}
    \draw[thick] (0,0) -- (-2,1.5) -- (-1,4) -- (2,3) -- cycle;
    \node at (1,1) {$q_i$};
    \draw[thick, noamblue] (0,0) -- (-1,4) node[pos=0.35, right] {$w_{i,d}^+$};
    \draw[thick, purple] (-2,1.5)-- (2,3) node[pos = 0.6, above] {$w_{i,d}^-$};
\end{tikzpicture}
    \caption{A quadrilateral $q_i$ with its forward diagonal $w_{i,d^+}$ in blue and its backward diagonal $w_{i,d^-}$ in red.}
    \label{fig:backward_diagonal}
\end{figure}
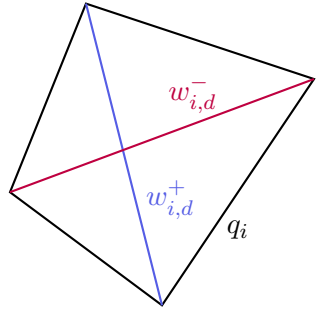

\begin{definition}[Backward diagonal]
    Let $q_i \in Q$ be a quadrilateral. The \emph{backward diagonal} $w_{i,d^-}$ of $q_i$ is the diagonal joining the left vertex to the right vertex of $q_i$, as pictured in Figure \ref{fig:backward_diagonal}.
\end{definition}

The definition of the backward diagonal $w_{i,d^-}$ is given in a way, so that the forward diagonal $w'_{i,d^+}$ of the quadrilateral $q'_i \in Q' = RQ$ (with the new labeling) is obtaind by rotating the backward diagonal of $q_{\pi_\ell^{-1}(i)}$. Formally, we can express this as
\begin{equation*}
    w'_{i,d^+} = \ii \cdot w_{\pi_\ell^{-1}(i), d^-} = \ii \cdot \left(w_{\pi_\ell^{-1}, r} - w_{\pi_\ell^{-1}, \ell}\right).
\end{equation*}

Geometrically, it is easy to see that a left (or right) staircase becomes a right (or left) staircase when being acted upon by $R$. Let us be more precise. 

\begin{figure}[ht]
    \centering
    \begin{tikzpicture}[scale= 0.85]
    \begin{scope}
    \foreach \p in {0,1,2}{
    \draw[thick] (-\p, \p) -- (-\p-1, \p + 1) -- (-\p, \p + 2) -- (-\p+1, \p+1) -- cycle;
    }

    \node at (0,1) {$q_{\pi_\ell^{-1}(i)}$};
    \node at (-1,2) {$q_i$};
    \node at (-2,3) {$q_{\pi_\ell(i)}$};
    
    \draw[thick, dashed] (-3,3) -- (-3.5,3.5);
    \draw[thick, dashed] (-2,4) -- (-2.5, 4.5);
    \draw[thick, dashed] (0,0) -- (0.5,-0.5);
    \draw[thick, dashed] (1,1) -- (1.5, 0.5);
    \end{scope}

    \begin{scope}[rotate around= {90:(-1,2)}, yshift= -140]
    \foreach \p in {0,1,2}{
    \draw[thick] (-\p, \p) -- (-\p-1, \p + 1) -- (-\p, \p + 2) -- (-\p+1, \p+1) -- cycle;
    }

    \node at (0,1) {$q'_{i}$};
    \node at (-1,2) {$q'_{\pi_\ell(i)}$};
    \node at (-2,3) {$q'_{\pi_\ell^2(i)}$};
    
    \draw[thick, dashed] (-3,3) -- (-3.5,3.5);
    \draw[thick, dashed] (-2,4) -- (-2.5, 4.5);
    \draw[thick, dashed] (0,0) -- (0.5,-0.5);
    \draw[thick, dashed] (1,1) -- (1.5, 0.5);
    \end{scope}

    \draw[->, >=latex, thick] (0,3) to[out=30, in=180] (1.5,3.5) to[out=0, in=150] (3,3);

    \begin{scope}[rotate around= {90:(-1,2)}, yshift= -250]
    \foreach \p in {0,1,2}{
    \draw[thick] (-\p, \p) -- (-\p-1, \p + 1) -- (-\p, \p + 2) -- (-\p+1, \p+1) -- cycle;
    }

    \node at (0,1) {$q_{\pi_r(i)}$};
    \node at (-1,2) {$q_i$};
    \node at (-2,3) {$q_{\pi_r^{-1}(i)}$};
    
    \draw[thick, dashed] (-3,3) -- (-3.5,3.5);
    \draw[thick, dashed] (-2,4) -- (-2.5, 4.5);
    \draw[thick, dashed] (0,0) -- (0.5,-0.5);
    \draw[thick, dashed] (1,1) -- (1.5, 0.5);
    \end{scope}

    \begin{scope}[xshift = 390]
    \foreach \p in {0,1,2}{
    \draw[thick] (-\p, \p) -- (-\p-1, \p + 1) -- (-\p, \p + 2) -- (-\p+1, \p+1) -- cycle;
    }

    \node at (0,1) {$q'_{\pi_\ell\pi_r^{-1}(i)}$};
    \node at (-1,2) {$q'_{\pi_\ell(i)}$};
    \node at (-2,3) {$q'_{\pi_\ell\pi_r(i)}$};
    
    \draw[thick, dashed] (-3,3) -- (-3.5,3.5);
    \draw[thick, dashed] (-2,4) -- (-2.5, 4.5);
    \draw[thick, dashed] (0,0) -- (0.5,-0.5);
    \draw[thick, dashed] (1,1) -- (1.5, 0.5);
    \end{scope}

    \draw[->, >=latex, thick, xshift = 253, yshift = -60] (0,3) to[out=-30, in=180] (1.5,2.5) to[out=0, in=210] (3,3);

    \node at (1.5,4) {$R$};
    \node[xshift = 215, yshift=-100] at (1.5,4) {$R$};
\end{tikzpicture}
    \caption{A left (or right) cycle becomes a right (or left) cycle under application of $R$.}
    \label{fig:rotation_cycle}
\end{figure}
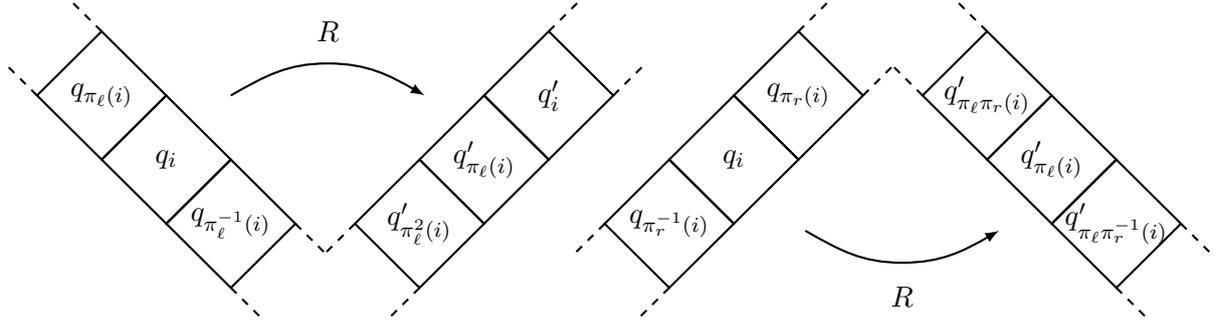

\begin{lemma}[Cycles after rotation]\label{lem:rotated_cycle}
    Let $c \in \pi_\varepsilon$ be a cycle. We define
    \begin{equation*}
        c' \coloneqq
        \begin{cases}
            c^{-1} \quad &\text{if } c \text{ is a cycle of }\pi_\ell,\\
            \pi_\ell c &\text{if } c \text{ is a cycle of }\pi_r.
        \end{cases}
    \end{equation*}
    Then, if $S_c$ is a left (or right) staircase for $Q$, it corresponds to the right (or left) staircase $S_{c'}$ for $Q' = RQ$ under the action of $R$, i.e., $S_{c'}$ is the union of the rotated quadrilaterals $\ii \cdot q$, where $q \in S_c$.
\end{lemma}
\begin{proof}
    Examining Figure \ref{fig:rotation_cycle}, we see that if $c$ is a left cycle of $\boldsymbol{\pi} = (\pi_\ell, \pi_r)$, then its inverse $c^{-1}$ is a right cycle for $\boldsymbol{\pi}' = (\pi_\ell\pi_r\pi_\ell^{-1}, \pi_\ell^{-1})$. If $c = \{i_1, \ldots, i_n\} \in \pi_r$, then $\pi_\ell c = \{\pi_\ell(i_1), \ldots, \pi_\ell(i_n)\}$ is a left cycle of $\boldsymbol{\pi}'$.
\end{proof}

Naturally, we want to look for conditions under which we can perform a staircase move on $S_{c'}$, i.e., conditions under which $S_{c'}$ is well-slanted. Recall that the well-slantedness of $S_c$, the original staircase, was independent of the suspension datum $\boldsymbol{\tau}$. It follows that whether $S_{c'}$ is well-slanted is independent of the (original) length datum $\boldsymbol{\lambda}$ and depends only on $\boldsymbol{\tau}$, since the roles of the length and suspension data are reversed by $R$. Thus, the set of parameters, where a move in $S_{c'}$ is possible, is given by $\{\boldsymbol{\pi}\}\times \Delta_{\boldsymbol{\pi}} \times \Theta_{\boldsymbol{\pi},c}$, where
\begin{equation*}
    \Theta_{\boldsymbol{\pi},c} \coloneqq
    \begin{cases}
        \{\boldsymbol{\tau} \in \Theta_{\boldsymbol{\pi}} \mid \tau_{i,d^-} = \tau_{i,r} - \tau_{i,\ell} < 0, \, i \in c\} \quad &\text{if }c \in \pi_\ell, \\
\{\boldsymbol{\tau} \in \Theta_{\boldsymbol{\pi}} \mid \tau_{i,d^-} = \tau_{i,r} - \tau_{i,\ell} > 0, \, i \in c\} &\text{if }c \in \pi_r.
    \end{cases}
\end{equation*}
This additional information on the action on the cycles and the suspension data gives the following proposition, which is a more refined version of Proposition \ref{prop:rotation_invertible}.

\begin{proposition}[Domain and codomain of rotation operator]\label{prop:domain_codomain_rotation}
    Let $\boldsymbol{\pi} = (\pi_\ell, \pi_r)$ be a combinatorial datum of a quadrangulation $Q \in \mathcal{Q}_k$ and let $\boldsymbol{\pi}' = (\pi_\ell\pi_r\pi_\ell^{-1}, \pi_\ell^{-1})$ be the combinatorial datum of $Q' = RQ$.

    Let $c$ be a cycle of $\boldsymbol{\pi}$ and let $c'$ be the corresponding cycle in $\boldsymbol{\pi}'$ given by the formula in Lemma \ref{lem:rotated_cycle}. Then,
    \begin{enumerate}
        \item $R$ maps $\{\boldsymbol{\pi}\} \times \Delta_{\boldsymbol{\pi}}\times \Theta_{\boldsymbol{\pi},c}$ bijectively onto $\{\boldsymbol{\pi'}\} \times \Delta_{\boldsymbol{\pi}',c'}\times \Theta_{\boldsymbol{\pi}'}$
        \item $R$ maps $\{\boldsymbol{\pi}\} \times \Delta_{\boldsymbol{\pi},c}\times \Theta_{\boldsymbol{\pi}}$ bijectively onto $\{\boldsymbol{\pi'}\} \times \Delta_{\boldsymbol{\pi}'}\times \Theta_{\boldsymbol{\pi}',c'}$
    \end{enumerate}
\end{proposition}

This result now extends to the self-duality we announced in the beginning of the section.

\begin{theorem}[Self-duality]\label{thm:self_duality}
    Let $\boldsymbol{\pi}$ be a permuation datum and $c \in \boldsymbol{\pi}$ be a (left or right) cycle of $\boldsymbol{\pi}$. Then, the map
    \begin{equation}\label{eq:staircase_move}
        \hat{m}_{\boldsymbol{\pi}, c} \colon \{\boldsymbol{\pi}\}\times \Delta_{\boldsymbol{\pi},c}\times \Theta_{\boldsymbol{\pi}} \to \{c\cdot \boldsymbol{\pi}\}\times \Delta_{c\cdot \boldsymbol{\pi}}\times \Theta_{c\cdot \boldsymbol{\pi}, c}
    \end{equation}
    is a bijection. Moreover, the inverse is given by
    \begin{equation}\label{eq:self_duality}
        \hat{m}^{-1}_{\boldsymbol{\pi},c} = R^{-1} \circ \hat{m}_{c\cdot\boldsymbol{\pi},c'} \circ R,
    \end{equation}
    where $R$ is the rotation operator (see \ref{def:rotation_operator}) and $c'$ is given by the formula in Lemma \ref{lem:rotated_cycle}.
\end{theorem}

The result follows essentially from the definitions and Proposition \ref{prop:domain_codomain_rotation}, the latter of which ensures that the composition in \eqref{eq:self_duality} is meaningful. For a full proof we refer to \cite{delecroix2015diagonal}. Note that equation \eqref{eq:self_duality} expresses exactly the notion of self-duality. For a more thorough treatment of (self-)duality, we refer the reader to \cite{schweiger1995ergodic}.

\subsubsection{Markovian Structure}
In some sense we will make precise below, equation \eqref{eq:staircase_move} shows that there is a \emph{loss of memory} phenomenon, or commonly called a \emph{Markov property} tied to staircase moves. The space of quadrangulations $\mathcal{Q}_k$ naturally projects onto the corresopnding space of bi-partite IETs, which is given by
\begin{equation*}
    \{(\boldsymbol{\pi}, \boldsymbol{\lambda)}\mid \boldsymbol{\pi}\in \mathcal{G}, \boldsymbol{\lambda} \in \Delta_{\boldsymbol{\pi}}\}.
\end{equation*}
The projection just \enquote{forgets} the suspension datum. Note that this is completely anologous to the projection from the space $\Upsilon$ constructed in section \ref{sec:RV_Teichmuller_connection}. 

Let us write $m_{\boldsymbol{\pi},c}$ for the restriction of the map $\hat{m}_{\boldsymbol{\pi},c}$ to the space of bi-partite IETs, i.e., we set
\begin{equation*}
    m_{\boldsymbol{\pi},c}(\boldsymbol{\pi}, \boldsymbol{\lambda}) = \left(c\cdot \boldsymbol{\pi}, A_{\boldsymbol{\pi},c} \cdot \boldsymbol{\lambda}\right).
\end{equation*}

The following corollary is immediate.

\begin{corollary}[Staircase move bijective on the level of IETs]\label{cor:staircase_bijective_IET}
    The map
    \begin{equation*}
        m_{\boldsymbol{\pi}, c} \colon \{\boldsymbol{\pi}\}\times \Delta_{\boldsymbol{\pi},c} \to \{c\cdot \boldsymbol{\pi}\} \times \Delta_{c \cdot \boldsymbol{\pi}}
    \end{equation*}
    is a \emph{bijection}.
\end{corollary}

This result corresponds to a Markovian structure as follows. Corollary \ref{cor:staircase_bijective_IET} shows that for any given (oriented) path in the $DC$ graph $\mathcal{G} = \mathcal{G}(\boldsymbol{\pi})$, which is equivalently a sequence of staircase moves, there exists a quadrangulation $Q = (\boldsymbol{\pi}, \boldsymbol{w})$ from which we can apply these moves. This can be seen as follows.

Say we are given a directed path in $\mathcal{G}$ from $(\boldsymbol{\pi}_\mathrm{start})$ to $(\boldsymbol{\pi}_\mathrm{end})$. We may choose any $\boldsymbol{\lambda} \in \Delta_{\boldsymbol{\pi}_{\mathrm{end}}}$ and chase it through the graph in backwards direction, applying the appropriate (inverse) staircase move, until we reach the starting point $\boldsymbol{\pi}_\mathrm{start}$. This is possible by \ref{thm:self_duality}. This provides us with the data $(\boldsymbol{\pi}_\mathrm{start}, \boldsymbol{\lambda}_\mathrm{start})$, and taking any admissible suspension datum $\boldsymbol{\tau}$ gives a quadrangulation $Q$ such applying staircase moves traces out the initial path in $\mathcal{G}$. The \enquote{past} of a quadrangulation is irrelevant, which gives exactly the claimed Markovian property. 

\begin{figure}[ht]
    \centering
\[\begin{tikzcd}
	A && B \\
	\\
	D && C
	\arrow["a", color={rgb,255:red,214;green,92;blue,92}, from=1-1, to=1-3]\arrow[loop, distance = 1cm, in=300, out=360, "c", {pos = 0.1}]
	\arrow["b", color={rgb,255:red,214;green,92;blue,92}, from=1-3, to=3-3]
	\arrow["d", color={rgb,255:red,214;green,92;blue,92}, from=3-3, to=3-1]
\end{tikzcd}\]
    \caption{An example of a graph not exhibiting a Markovian structure, given that in any path, edge $d$ can only follow after edge $c$.}
    \label{fig:markov_counterexample}
\end{figure}
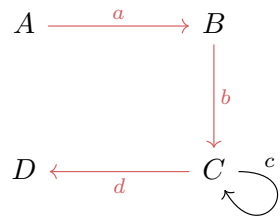

To make things even more clear, let us explain how such a Markovian property could be absent. Consider the graph in Figure \ref{fig:markov_counterexample} and suppose that edge $d$ can only follow after edge $c$. In this case, the past \emph{would} be relevant since, for example, the red path in Figure \ref{fig:markov_counterexample} could not be realized because of the \enquote{past} of the node $C$.

\subsubsection{Best Approximations via Staircase Moves}
The last interesting property of staircase moves we want to highlight is the fact that the algorithm \enquote{sees} all saddle connections of a given surface $X$ which realize the systoles along a Teichmüller geodesic. As we will see, this is the case because these saddle connections are exactly the wedges produced by any algorithm performing simultaneous staircase moves. This sharply contrasts with the situation we observed in the Rauzy--Veech algorithm, where we lack an obvious way to identify saddle connections from the zippered rectangle construction.

We proceed as follows. First, we will define best approximations for translation surfaces of genus $\mathbf{g} \geq 2$ and we will show that all best approximations in each bundle are produced by any diagonal changes algorithm that performs simultaneous staircase moves (Theorem \ref{thm:staircase_produce_best_approximations}). Then, we show that the geometric objects, namely the wedges an well-slanted staircases, produced by any sequence of staircase moves are the same (Theorem \ref{thm:geometric_objects_are_the_same}), before finally showing that the saddle connections which realize the systoles along a Teichmüller geodesic are contained in the set of best approximations (\ref{thm:systoles_along_teichmueller}).

A very natural point of view when trying to extend the notion of geometric best approximation (Definition \ref{def:geometric_best_approximation}) is to compare saddle connections which belong to the same bundle $\Gamma_i$. Let us introduce the following convention. Given a saddle connection on a translation surface $X$, we can associate to it a pair $(i, v)$, where $v \in \C$ is its holonomy vector and $i \in [k]$ is such that the saddle connection belongs to the bundle $\Gamma_i$. Conversely, if we know which bundle the saddle connection belongs to, its holonomy vector $v \in \C$ completely determines the saddle connection. Therefore, we can abuse the notation in the same way we have done above and identify saddle connections and holonomy vectors as long as the bundle is clear from context. For instance, we will use the following notation. For a saddle connection $v$, we write $\RE(v), \IM(v)$ and $|v|$ for the real part, the imaginary part and the length of the associated holonomy vector. Given a bundle $\Gamma_i$, we write $v \,(\in\C)$ for the saddle connection in $\Gamma_i$, i.e., we will write $v \in \Gamma_i$.

\begin{definition}[Geometric best approximation]\label{def:geometric_best_approximation_2}
    A saddle connection $w \in \Gamma_{i, r}$ is a \emph{right (geometric) best approximation}, if
    \begin{equation*}
        \IM(v) < \IM(w) \Rightarrow |\RE(v)| > |\RE(w)| \quad \text{for all } v \in \Gamma_{i,r}. 
    \end{equation*}
    Replacing $\Gamma_{i,r}$ by $\Gamma{i,\ell}$, we obtain the definition of \emph{left geometric best approximation}.
\end{definition}

As for the torus, we have a useful equivalent characterization for geometric best approximations based on immersed rectangles, where we call an \emph{immersed rectangle} $R \subseteq X$ a subset without singularities in its interior, which is obtained by isometrically immersing in $X$ a Euclidean rectangle with horizontal and vertical sies in $\C$. We remark that we do not require the rectangle to be \emph{embedded}, i.e., it is enough that the map is \emph{locally} injective meaning that self-intersections on $X$ may occur. 

\begin{lemma}[Best approximations via immersed rectangles]\label{lem:immersed_rectangles}
    A saddle connection $w$ on $X$ is a geometric best approximation if and only if there exists an immersed rectangle $R(w)$ in $X$ which has $w$ as its diagonal.
\end{lemma}
\begin{proof}
    For clarity, we will differentiate between a saddle connection $\gamma$ on $X$ and its holonomy vector, which we will denote by $\operatorname{hol}(\gamma) \in \C$. Moreover, we will write $\operatorname{hol(\Gamma_i)}$ for the set of holonomy vectors of saddle connections in $\Gamma_i$. For each saddle connection $\gamma$ in $\Gamma_{i, r}$ (or in $\Gamma_{i, \ell}$), let $\Tilde{R}(\gamma)$ be the rectangle given by
    \begin{equation*}
        \Tilde{R}(\gamma) = [0, \RE(\operatorname{hol(\gamma))}] \times [0, \IM(\operatorname{hol}(\gamma))],     
    \end{equation*}
    if $\gamma \in \Gamma_{i,r}$, or
    \begin{equation*}
        \Tilde{R}(\gamma)= [\RE(\operatorname{hol}(\gamma)), 0] \times [0, \IM(\operatorname{hol}(\gamma))],
    \end{equation*}
    otherwise. 

    This allows us to restate Definition \ref{def:geometric_best_approximation_2} as follows. A saddle connection $\gamma \in \Gamma_{i, r}$ (or $\Gamma_{i, \ell}$) is a (geometric) best approximation if and only if the rectangle $\Tilde{R}(\gamma)$ does not contain any element of $\operatorname{hol}(\Gamma_{i,r})$ (or $\operatorname{hol}(\Gamma_{i, \ell})$) in its interior. 

    We will treat the case where $w$ is a right-slanted saddle connection. The case for left-slanted saddle connections is analogous. Suppose that $w$ starts at a point $p_0 \in \Sigma$ and let $\Gamma_{i,r}$ be the bundle to which $w$ belongs. Suppose first that there exists an immersed rectangle $R(w) \subseteq X$ which has $w$ as a diagonal and does not contain any singularities in ints interior. Define the so-called \emph{developing map} by
    \begin{align*}
        \operatorname{dev}_{p_0} \colon R(w) &\to \C, \\
        p \mapsto \int_{p_0}^p \omega,
    \end{align*}
    where $\omega$ is the Abelian differential given by the translation structure on $X$. Morally, the developing map simply gives for any $p\in R(w)$ the distance to the lower left corner $p_0$. Thus it follows that the image $\operatorname{dev}_{p_0}(R(w))$ is exactly the rectangle $\Tilde{R}(w) \subseteq \C$ defined above. Towards a contradiction, let us assume that $w$ is \emph{not} a best approximation, so that by the alternative definition of best approximation developed in the beginning of the proof we know that $\operatorname{hol}(\Gamma_{i,r})$ intersects the interior of $\Tilde{R}(w)$. Therefore, we find a saddle connection $\gamma \in \Gamma_{i,r}$ whose holonomy $\operatorname{hol(\gamma)}$ belongs to the interior of $\Tilde{R}(w)$ so in partiuclar the saddle connection $\gamma$ belongs to the same bundle as $w$. But this means that $\gamma$ is contained in $R(w)$, contradicting our initial assumption.

    Conversely, assume now that $w \in \Gamma_{i,r}$ is a best approximation. We claim that we can immerse the rectangle $\Tilde{R}(w)$ in $X$ such that its image $R(w)$ is an immersed rectangle which has $w$ as diagonal and does not contain any singularities in its interior. Towards this end, we define the immersion $\iota$ by mapping the point $z = \rho \e^{\ii \theta} \in \Tilde{\R}(\gamma)$ to $\iota(z) = \gamma_{\rho}^{\theta}(p_0)$, which has distance $\rho$ from $p_0$ and belongs to the unique linear trajectory $(\gamma_t^\theta(p_0))_{t \geq 0}$ in direction $\theta$ starting at $p_0$ in $\Gamma_{i,r}$, that is $|\sphericalangle(\gamma_t^\theta(p_0), w)|< \frac{\pi}{2}$. This map $\iota$ is well-defined as long as the linear trajectory does not hit any singularities. By showing that this does not occur, we aslo show that the image $\iota(R(w))$ does not intersect $\Sigma$, which finishes the proof. If there was a singularity $q \in \Sigma$ in the interior of $R(w)$, since the saddle connection $v$ connecting $p_0$ to $q$ is inside $R(w)$, it belongs in particular to the same bundle as $w$ and thus satisfies $\operatorname{hol}(v) \in \Tilde{R}(w)$. So in this case, the interior of $\Tilde{R}(w)$ intersects $\operatorname{hol}(\Gamma_{i,r}$, which is impossible again by the equivalent definition of best approximation established in the beginning of this proof. 
\end{proof}

We know state the first result announced in the beginning of the section.

\begin{theorem}[Simultaneous staircase moves produce best approximations]\label{thm:staircase_produce_best_approximations}
    Let $X$ be a translation surface with a quadrangulation $Q$ that corresponds to a bi-partite IET on $2k$ intervals, such that it contains neither hoizontal nor vertical saddle connections. Let $(Q^{(n)})_{n \in \Z}$ be any sequence of labeled quadrangulations $Q^{(n)} = (\boldsymbol{\pi}^{(n)}, \boldsymbol{w}^{(n)})$ of $X$, where $Q^{(n+1)}$ is obtained from $Q^{(n)}$ by simultaneous staircase moves.

    Then, for each $1 \leq i \leq k$, the saddle connections in the sequence $\left(w_{i,\ell}^{(n)}\right)_{n \in \Z}$ (or $\left(w_{i,r}^{(n)}\right)_{n \in \Z})$ are exactly all best geometric approximations in $\Gamma_{i, \ell}$ (or $\Gamma_{i, r}$) ordered by increasing imaginary part. 
\end{theorem}
Thanks to Lemma \ref{lem:immersed_rectangles}, one can prove Theorem \ref{thm:staircase_produce_best_approximations} by slightly adapting the proof of Theorem \ref{thm:DU15_thm1} and we refer the interested reader to \cite{delecroix2015diagonal} for further details. Let us remark that in fact, the proof of \ref{thm:DU15_thm1} in this paper was obtained by simplifying the proof of Theorem \ref{thm:staircase_produce_best_approximations} from \cite{delecroix2015diagonal}, where it appears as Theorem 10.

The following corollary applies in case we are interested in sequences obtained by forward moves only.

\begin{corollary}\label{cor:wedges_are_best_approximations}
    Let $X$ be a translation surface with a quadrangulation $Q$ corresponding to a bi-partite IET on $2k$ intervals without vertical saddle connections and let $\left(Q^{(n)}\right)_{n \in \N}$ be a sequence of labeled quadrangulations of $X$, where $Q^{(n+1)}$ is obtained from $Q^{(n)}$ by simultaneous staircase moves. Then,
    \begin{enumerate}
        \item For each $1\leq i \leq k$, the saddle connections in the sequence $\left(w_{i,r}\right)_{n \in \N}$ and in $\left(w_{i,\ell}\right)_{n \in \N}$ are exactly all best geometric approximations $w$ in $\Gamma_{i,r}$ and in $\Gamma_{i,\ell}$ respectively, which have $\IM(w) \geq \IM(w_{i, \varepsilon})$ or equivalently, $\RE(w) \leq \RE(w_{i, \varepsilon})$, where $\varepsilon\in\{\ell, r\}$.
        \item For each $1 \leq i \leq k$, the set of diagonals $\left(w_{i,d}^{(n)}\right)_{n \in \N}$ coincide with the set of best approximations $v \in \Gamma_i$ such that $\IM(v) > \max\left(\IM(w_{i,\ell}), \IM(w_{i,r})\right)$, or equivalently, they coincide with the set of bottom sides of the quadrilaterals $\left(q_i^{(n)}\right)_{n \in \N}$ different from the one of $q_i^{(0)}$. 
    \end{enumerate}
\end{corollary}

In the proof we will need part of the following lemma, whose proof can be found in \cite{delecroix2015diagonal} as Lemma 14.

\begin{lemma}\label{lem:diagonals_become_sides}
    Let $Q$ be a quadrangulation of a surface $X$ which corresponds to a bi-partite IET on $2k$ intervals. There exists an infinite sequence of staircase moves starting from $Q$ such that the real part of each saddle connection tends to zero if and only if $X$ has no vertical saddle connection. 
    
    Moreover, if $X$ has no vertical saddle connection then for any infinite sequence of staircase moves starting from $Q$, there are infinitely many left and right diagonal changes and the width of each wedge goes to zero. 
\end{lemma}

\begin{remark}
    Another way to state Lemma \ref{lem:diagonals_become_sides} is to say that Keane's condition is exactly the condition needed to ensure that \emph{any} diagonal changes algorithm does not stop. 
    
    Note also, that Lemma \ref{lem:diagonals_become_sides} implies in particular that each diagonal in any quadrilateral eventually becomes part of a wedge under application of staircase moves. 
\end{remark}

\begin{proof}[Proof of Corollary \ref{cor:wedges_are_best_approximations}]
    The first item follows from Theorem \ref{thm:staircase_produce_best_approximations}, since geometric best approximations are produced ordered by increasing imaginary part. 

    By the remark after Lemma \ref{lem:diagonals_become_sides} we know that if $X$ has no vertical saddle connection, then each diagonal of $Q^{(n)}$ eventually becomes a side of a wedge, so that the second item follows from the first one.
\end{proof}

Using the following lemma, which states that diagonals uniquely determine their quadrilaterals, we can obtain the second result announced in the beginning of the section. A proof of this lemma can be found in \cite{delecroix2015diagonal}.

\begin{lemma}[Diagonals determine quadrilaterals]\label{lem:diagonals_determine_quadrilaterals}
    Let $X$ be a translation surface without vertical saddle connections and let $w$ be a saddle connection which is a geometric best approximation. Then there exists a unique admissible quadrilateral $q$ whose sides are all geometric best approximations and that has $w$ as forward diagonal. In particular, $q$ is \emph{embedded} in $X$.

    If moreover $w$ is left-slanted (or right-slanted), then there exists a unique right-slanted (or left-slanted) admissible quadrilateral in $X$ whose sides are best approximations and so that $w$ is its bottom left side (or bottom right side).  
\end{lemma}

\begin{remark}
    If, say, $w$ is left-slanted, then the second part of the lemma does not say anything about \emph{left-slanted} admissible quadrilaterals that have $w$ as bottom left side. I tis indeed the case that there may be none or several such quadrilaterals. 
\end{remark}

\begin{theorem}[Simultaneous staircase moves produce the same geometric objects]\label{thm:geometric_objects_are_the_same}
    Let $X$ be a translation surface corresponding to a bi-partite IET on $2k$ intervals without vertical saddle connections and let $Q$ be a quadrangulation of $X$. Let 
    \begin{equation*}
        \left(Q_1^{(n)}\right)_{n \in \N} \quad \text{and} \quad \left(Q_2^{(n)}\right)_{n \in \N}
    \end{equation*}
    be any two sequences of quadrangulations of the surface $X$ such that $Q_1^{(0)} = Q_2^{(0)} = Q$ and $Q_i^{(n+1)}$ is obtained from $Q_i^{(n)}$ by simultaneous staircase moves for $i \in \{1,2\}$. 

    \begin{enumerate}
        \item The collections of the wedges associated to the quadrangulations of either sequence coincide.
        \item The sets of well-slanted staircases associated to the quadrangulations of either sequence coincide.
    \end{enumerate}
\end{theorem}

\begin{proof}
    For $(Q_1^{(n)})_{n \in \N}$ and $(Q_2^{(n)})_{n\in \N}$ as in the theorem, Lemma \ref{lem:diagonals_become_sides} ensures that each diagonal in $Q_i^{(n)}$ eventually becomes part of a wedge for $i \in \{1,2\}$. By Corollary $\ref{cor:wedges_are_best_approximations}$, the diagonals of the two sequences coincide. Lastly, Lemma \ref{lem:diagonals_determine_quadrilaterals} implies that the quadrilaterals, and thus the wedges, of the two sequences coincide which proves the first statement.

    To prove the second statement, it is enough to show that each quadrilateral uniquely determines the well-slanted staircase to which it belongs, since in the first part we showed that the quadrilaterals of the two sequences are the same. Towards this goal, let $q = q_i$ be a right-slanted quadrilateral in $Q_1^{(n)}$ for some $n \in \N$ and let $w$ be its right top side. It suffices to show that there is a unique quadrilateral $q'$ which is right-slanted and has $w$ as its bottom left side, since such a quadrilateral is necessarily a neighbour of $q$ in a well-slanted right staircase. From Theorem \ref{thm:staircase_produce_best_approximations}, we know that $w$ is a geometric best approximation and Lemma \ref{lem:diagonals_determine_quadrilaterals} guarantees existence and uniqueness. We can continue this argument inductively to conclude that the right well-slanted staircase which contains $q$ is uniquely determined.
\end{proof}

In order to show the last of the announced results, we will establish the following lemma which holds for general translation surfaces, i.e., not just those corresponding to bi-partite IETs.

\begin{lemma}[Saddle connections realizing the systole are best approximations]\label{lem:saddle_connections_realize_systole}
    Let $X$ be a translation surface and let $w$ be a saddle connection on $X$ which realizes the systole for some time $t_0$ along the Teichmüller geodesic $(g_tX)_{t \in \R}.$ Then, $w$ is a geometric best approximation.
\end{lemma}

Before giving a proof of Lemma \ref{lem:saddle_connections_realize_systole}, let us state the next result and show that it follows at once.

\begin{theorem}[Systoles along Teichmüller geodesic are best approximations]\label{thm:systoles_along_teichmueller}
    Let $X$ be a translation surface corresponding to a bi-partite IET of $2k$ intervals with neither horizontal nor vertical saddle connections. Let $\left(Q^{(n)}\right)_{n \in \Z}$ be a sequence of quadrangulations of the surface $X$, where $Q^{(n+1)}$ is obtained from $Q^{(n)}$ by simultaneous staircase moves. Then, the set of saddle connections on $X$ which realize the systoles along the Teichmüller geodesic passing through $X$ is a subset of the sides of the quadrangulations $Q^{(n)}, n \in \Z$. 
\end{theorem}
\begin{proof}
    If $v$ realizes the systole at time $t_0$, by Lemma \ref{lem:saddle_connections_realize_systole} we have that $v$ is a best approximation. It then follows by Theorem \ref{thm:staircase_produce_best_approximations} that $v$ appears in one of the wedges.
\end{proof}

\begin{proof}[Proof of Lemma \ref{lem:saddle_connections_realize_systole}]
    Let $w \in \Gamma_i$ and assume that $w$ realizes the systole for some time $t_0 > 0$. The key property we need is that being a geometric best approximation is invariant under the Teichmüller geodesic flow $(g_t)_{t \in \R}$. Indeed, immersed rectangles with horizontal and vertical sides are mapped to immersed rectangles of the same form, and if the initial rectangle did not contain a singularity in its interior then neither does the rectangle obtained under the flow so the property follows by Lemma \ref{lem:immersed_rectangles}. Thus, we can replace $X$ by $g_{-t}X$ and suppose that $t_0 = 0$. For any saddle connection $u$ in $X$, the length $|u|$ of the associated holonomy vector satisfies
    \begin{equation*}
        |u| \leq |w|,
    \end{equation*}
    by the definition of systole. In particular, the semicirlce in $\C$ of radius $|w|$ centered at the origin does not contain any point of $\Gamma_i$ in its interior which implies that the rectangle $R(w)$, where we use the notation from Lemma \ref{lem:immersed_rectangles}, is an immersed rectangle with diagonal $w$, i.e., the rectangle does not contain any points of $\Gamma_i$ in its interior so that $w$ is indeed a geometric best approximation.
\end{proof}

If we only look forward in time, we need to make sure we don't miss some saddle connection realizing a systole for $t$ too small, meaning that we obtain the following corollary where a small adjustment of the statement is needed.

\begin{corollary}
    Let $X$ be a translation surface corresponding to a bi-partite IET on $2k$ intervals with neither horizontal nor vertcial saddle connections. Let $Q = (\boldsymbol{\pi}, \boldsymbol{w})$ be a quadrangulation of $X$ for which each side $w$ satisfies
    \begin{equation}
        |\RE(w)| > \operatorname{sys}(X).
    \end{equation}
    Let $\left(Q^{(n)}\right)_{n \in \N}$ be any sequence of quadrangulations with $Q^{(0)} = Q$ and such that $Q^{(n+1)}$ is obtained from $Q^{(n)}$ by simultaneous staircase moves. Then, the saddle connections on $X$ which realize the systoly along the Teichmüller geodesic ray $(g_tX)_{t\geq 0}$ are a subset of the sides of the quadrangulations in $\{Q^{(n)}\mid n \in \N\}$. 
\end{corollary}
\begin{proof}
    Since $|\RE(w)| > \operatorname{sys}(X)$ for each side $w$ of $Q$, the first part of Corollary \ref{cor:wedges_are_best_approximations} implies that the set of sides of $(Q^{(n)})_{n \in \N}$ contains all best approximations $w$ that satisfy $\RE(w) \leq \operatorname{sys}(X)$. By Lemma \ref{lem:saddle_connections_realize_systole}, it is enough to show that saddle connections $w$ that realize the systole at a positive time $t_0$ satisfy $|\RE(w)| \leq \operatorname{sys}(X)$. By definition of systole, we always have
    \begin{equation*}
        \operatorname{sys}(g_tX) \leq \e^t \operatorname{sys}(X).
    \end{equation*}
    Therefore, if $w$ is a saddle connection which realizes a systole at time $t_0 > 0$, then we have
    \begin{align*}
        |\RE(g_{t_0}w)| &\leq |g_{t_0}w| \\&= \operatorname{sys(g_{t_0}X)} \\&\leq \e^{t_0}\operatorname{sys}X. 
    \end{align*}
    Since $|\RE(g_{t_0}w)| = \e^{t_0}|\RE(w)|$ the result follows. 
\end{proof}

\subsection{Existence of Quadrangulations and Staircase Moves}\label{sec:existence_of_quadrangulations}

The theory about diagonal changes algorithms based on staircase moves we have developed above relies heavily on the fact that we are always able to find well-slanted staircases. As we have mentioned, this may not always be possible. In the results we have presented above, the condition we needed to impose was that the Poincaré first return map to the union of the wedges corresponds to a \emph{bi-partite} IET which corresponds to an assumption of additional symmetry. Another way to state this hypothesis is to say that we restrict ourselves to translation sufaces in \emph{hyperelliptic components}. In this section, we want to introduce this concept and show that every translation surface in such a hyperelliptic component admits a quadrangulation with at least one well-slanted staircase. As a consequence, we can apply any algorithm which produces new quadrangulations by simultaneous staircase moves in such hyperelliptic components.

\subsubsection{Hyperelliptic Components}

In order to define hyperelliptic components properly, we will need to work with the definition of translation surfaces via holomorphic differentials, i.e., Definition \ref{def:translation_surface_3}. More precisely, we want to extend this definition to \emph{half-}translation surfaces. Our exposition here is largely based on \cite{massart2022short}, with some ideas and notations borrowed from \cite{lanneau2004hyperelliptic}. 

Recall that half-translation surfaces are obtained by gluing polygons by translations \emph{and by rotations by} $\pi$. Equivalently, half-translation surfaces are given by an atlas where the transition maps are of the form

\begin{equation*}
    z \mapsto \pm z +c
\end{equation*}
for some $c \in \C$. Note that the map $z \mapsto -z$ does \emph{not} preserve the linear form $\dd z$. 

Analogously to the definition of a translation structure coming from a holomorphic, or Abelian differential $\omega$ we can define a half-translation structure through a holomorphic \emph{quadratic} differential. 

\begin{definition}[Holomorphic quadratic differential]
    A \emph{holomorphic quadratic differential} $q$ is an application, that assigns to each $x \in X$ a quadratic form $q(x)$ from the tangent space to $X$ at $x$ to $\C$, with the condition that $q(x)$ depends holomorphically on $x$. That is, for each $z \in \C \simeq T_xX$, the map
    \begin{equation*}
        x \mapsto q(x)(z)
    \end{equation*}
    is holomorphic. 
\end{definition}
Viewed as a complex manifold, $X$ has dimension 1 and the only complex-valued quadratic forms from $\C$ to itself are the maps
\begin{equation}\label{eq:quadratic_form}
    z \mapsto \lambda z^2, \quad \lambda \in \C.
\end{equation}
We will write the square map $z \mapsto z^2$ as $\dd z^2$, which by \eqref{eq:quadratic_form} is a basis for all quadratic forms. In charts, a holomorphic quadratic differential reads $f(z)\dd z^2$ for some holomorphic function $f$. We have the usual compatibility requirement, which can be seen as the chain rule, i.e., if two charts overlap and $q$ reads $f(z)\dd z^2$ in the first and $g(z) \dd z^2$ in the second, then the transition map $\phi$ must satisfy
\begin{equation*}
    f(\phi(z))\cdot \phi'(z)^2 = g(z).
\end{equation*}
\begin{remark}
    To be fully correct, we need to allow $f$ to have simple poles, so that we are actually not dealing with \emph{holomorphic} forms but with special kind of \emph{meromorphic} quadratic forms. 
\end{remark}

Given a holomorphic differential $\omega$ that reads $f(z) \dd z$ in charts, we can obtain a quadratic differential by the formula $f(z)^2 \dd z^2$. However, in general this construction only goes one way: Not every holomorphic function is a square. So while not every quadratic differential is a square, it is true that every quadratic differential becomes a square in a suitable double cover of $X$ using the classical Riemann surface construction of the squareroot function $z \mapsto \sqrt{z}$. 

Let us now take a closer look at half-translation surfaces, or, equivalently, a Riemann surface endowed with a quadratic differential $q$ which we will denote by $(X, q)$. As before, we will write $\Sigma \subseteq X$ for the singularities corresponding to the images of the vertices of the polygons in the construction of our half-translation surfaces. Details about the following result can be found in \cite{athreyamasur2023translationsurfaces}.

\begin{proposition}[Conical singularities on half-translation surfaces]
    A quadratic differential induces conical singularities of the form $\pi k, k \geq 1$. 
\end{proposition}
We want to emphasize the difference to translation surfaces, where conical singularities are of the form $2\pi k, k \geq 1$.

On half-translation suface, there is no well-defined notion of directions $\theta \in S^1$, but a quadratic differential defines a well-defined notion of \emph{non-oriented lines} in direction $\theta \in \PP^1\R$, the real projective line. As for translation surfaces given by Abelian differentials, we have a notion of isomorphism.

\begin{definition}[Isomorphic quadratic differentials]
    Two quadratic differentials $(X,q)$ and $(X',q')$ are \emph{isomorphic}, if there exists a homeomorphism $h \colon X \to X'$ such that
    \begin{equation*}
        q = h^*q',
    \end{equation*}
    i.e., we require that $q(x) = q'(f(x))$ for all $x \in X$.
\end{definition}
\begin{remark}
    This notion of isomorphism can also be defined using cut and paste operations on polygons.
\end{remark}

\begin{definition}[Class of quadratic differentials]\label{def:class_quadratic_differentials}
    We write $\mathcal{Q}(k_1 - 2, \ldots, k_n - 2)$ for the equivalence class of quadratic differentials with conical singularities of angles $\pi k_1, \ldots, \pi k_n$.
\end{definition}

The number $k_i - 2$ in Definition \ref{def:class_quadratic_differentials} corresponds to the degree of the quadratic differential $q$, i.e., to the order of the zero at the singularity. Indeed, locally we can write $q$ as
\begin{equation*}
    z^{k_i - 2}\dd z^2
\end{equation*}
around a conical singularity of angle $\pi k_i$.

There is a topological restriction on the cone angles, which is essentially a further consequence of the Gauss--Bonnet Theorem (Theorem \ref{thm:gauss_bonnet}). We have
\begin{equation}\label{eq:topological_restriction_quadratic_differentials}
    \sum_{i = 1}^n k_i = 4\mathbf{g} - 4 + 2n,
\end{equation}
where $\mathbf{g}$ is the genus of the surface. Let us introduce the following notation. If we have multiple singularities of the same cone angle, we denote the multiplicity as an exponent, i.e., if there are $m_i$ singularities with total angle $\pi k_i$ then we will write
\begin{equation*}
    \mathcal{Q}((k_1 - 2)^{m_1}, \ldots, (k_n - 2)^{m_n}). 
\end{equation*}

\begin{example}[Class of quadratic differentials]
    Using the notation introduced above, we write
    \begin{equation*}
        \mathcal{Q}(2,3,3,3,4,4) = \mathcal{Q}(2,3^3,4^2).
    \end{equation*}
\end{example}

Let us provide further details as to how a quadratic differential corresponds to an Abelian differential in a suitable cover, which we will call \emph{orientation cover}.

\begin{definition}[Orientation cover]\label{def:orientation_cover}
    For a quadratic differential $(X,q)$, its canonical \emph{orientation cover} is the unique (up to isomorphism) Abelian differential $(\Tilde{X},\omega)$ such that there is a degree 2 map
    \begin{equation*}
        p \colon \Tilde{X} \to X,
    \end{equation*}
    verifying $p^*q = \omega^2$, where $p^*q$ denotes the pullback of $q$ by $p$.
\end{definition}

The orientation cover satisfies the following properties.

\begin{proposition}[Properties of the orientation cover]\label{prop:properties_orientation_cover}
    Let $(\Tilde{X}, \omega)$ be the orientation cover of the quadratic differential $(X, q)$. 
    \begin{enumerate}
        \item Each singularity of angle $\pi k_i$, with $k_i$ \emph{even}, is not ramified, i.e., it has two preimages in $\Tilde{X}$ of angle $\pi k_i$.
        \item Each singularity of angle $\pi k_i$ with $k_i$ \emph{odd} is ramified, i.e., it only has one preimage under $p$ and it corresponds to a singularity on $\Tilde{X}$ of angle $2\pi k_i$.
    \end{enumerate}
\end{proposition}

\begin{example}
    A quadratic differential in $\mathcal{Q}(2, 3^2)$ corresponds to an Abelian differential in $\mathcal{H}(1^2, 4^2)$, since it follows from $2 = k_1 - 2$ that $k_1$ is even, so it corresponds to two zeros of order 1 in the orientation cover, and similarly each of the two zeros of order $3$ for the quadratic differential corresponds to a zero of order $4$ in the orientation cover, since $k_2 = 5$ is odd.

    The following table provides a convenient overview over the different possibilities.

    \begin{center} 
        \begin{tabular}{cccc}
        \toprule
         $k_i - 2$, &  & $k_i - 1$, &  \\
         order of zero of $q$ & cone angle in $(X,q)$ & order of zero(s) of $\omega$ & cone angle(s) in $(\Tilde{X}, \omega)$\\ 
        \midrule  
        4 & $(4+2)\pi = 6\pi$ & $(2,2)$ & $(2+1)\pi \cdot 2 = 6\pi$ \\
        5 & $(5+2)\pi = 7\pi$ & $6$ & $(6+1) \pi \cdot 2 = 14\pi$ \\
        $2n$ & $(2n+2)\pi$ & $(n,n)$ & $(n+1)\pi \cdot 2 = (2n+2)\pi$ \\
        $2n+1$ & $(2n+3)\pi$ & $2n + 2$ & $(2n+3)\pi \cdot 2$ \\
        $-1$ & $\pi$ & 0 & $2\pi$\\
        \bottomrule
    \end{tabular}
    \end{center}
\end{example}

Since a degree two cover is always normal, an orientation cover always comes with an involution whose quotient is the corresponding quadratic differential. 

When a quadratic differential varies in its stratum $\mathcal{Q}(k_1 - 2, \ldots, k_n - 2)$, the associated orientation cover varies in a connected component of the corresponding stratum of Abelian differentials.

\begin{definition}[Hyperelliptic locus]
    When the stratum of quadratic differentials is a sphere, i.e., the stratum is of the form $\mathcal{Q}(k_1 - 2, \ldots, k_n - 2)$ with $\sum_{i = 1}^n k_i = 2n - 4$ implying that the genus satisfies $\mathbf{g} = 0$, then the corresponding locus of Abelian differentials obtained as a double cover from these quadratic differentials is called a \emph{hyperelliptic locus}. Its elements are \emph{hyperelliptic surfaces} and the associated involution is a \emph{hyperelliptic involution}.

    The points which are fixed by the involution are called \emph{Weierstrass points}. They may be singular or regular. If they are regular, they project down to conical singularities of cone angle $\pi$ on the sphere that are called \emph{poles}, since they correspond to singularities of the form $z^{-1}\dd z$  for the quadratic differential.
\end{definition}

Using a classical result called the \emph{Riemann--Hurwitz formula} from the theory of Riemann surfaces, we are able to give the exact number of Weierstrass points on hyperelliptic surfaces.

\begin{theorem}[Riemann--Hurwitz formula]
    Let $X$ and $X'$ be Riemann surfaces and let $p \colon \Tilde{X} \to X$ be a holomorphic $m$-sheeted covering map with branch points $s_1, \ldots, s_n$ with corresponding ramification indices $e_1, \ldots, e_n$. Write $\mathbf{g}(X')$ for the genus of $X'$ and $\mathbf{g}(X)$ for the genus of $X$. Then, the formula
    \begin{equation*}
        2\mathbf{g}(X') - 2 = m \cdot (2\mathbf{g}(X) - 2) + \sum_{i = 1}^n (e_i - 1)
    \end{equation*}
    holds. 
\end{theorem}

\begin{theorem}[Number of Weierstrass points on hyperelliptic surfaces]
    A hyperelliptic surface $X'$ of genus $\mathbf{g}$ has $2\mathbf{g} + 2$ Weierstrass points. 
\end{theorem}
\begin{proof}
    Since $X'$ is hyperelliptic, it is constructed as the orientation cover of $(X,q)$, which is a sphere. Consequently, $\mathbf{g}(X) = 0$. Moreover, note that the ramification index of each branch point is given by $e_i = 2$, since we have a two-fold covering. Hence, by the Riemann--Hurwitz formula, 
    \begin{equation*}
        2\mathbf{g} - 2 = -4 + \sum_{i = 1}^n 1 \quad \Leftrightarrow \quad 2\mathbf{g}+2 = n,
    \end{equation*}
    which is what we wanted to show.
\end{proof}

In most cases, hyperelliptic loci have positive codimension in the corresponding stratum of Abelian differentials, but an infinite family of hyperelliptic loci have full dimension and form connected components. This is a consequence of the complete classification of strata done by Kontsevich and Zorich in \cite{kontsevich2003connected}.

\begin{theorem}[\cite{kontsevich2003connected}, section 2.1 p.5-7]\label{thm:hyperelliptic_loci}
    In each stratum $\mathcal{H}(2\mathbf{g} - 2)$, the hyperelliptic locus built as the orientation cover of quadratic differentials in $\mathcal{Q}(k-2, -1^{k+2})$ for $k = 2\mathbf{g}-1$ forms a connected component. 

    The same is true for each stratum $\mathcal{H}(\mathbf{g}-1, \mathbf{g}-1)$, i.e., the hyperelliptic locus built as an orientation cover of quadratic differentials in $\mathcal{Q}(k-2, -1^{k+2})$ for $k = 2\mathbf{g}$ forms a connected component.

    Moreover, these are the only hyperelliptic loci that form connected components of strata.
\end{theorem}

\begin{definition}[Hyperelliptic component]
    For each integer $k \geq 1$, we denote by $\mathcal{C}^\mathrm{hyp}(k)$ the unique hyperelliptic component which contains surfaces with total conical angle $2\pi k$.
\end{definition}

\begin{remark}
    If $k$ is odd, then $\mathcal{C}^\mathrm{hyp}(k) \subseteq \mathcal{H}(k-1)$, while if $k$ is even, then $\mathcal{C}^\mathrm{hyp}(k) \subseteq \mathcal{H}(\tfrac{k}{2}-1, \tfrac{k}{2}-1)$. For small $k$ the situation is particularly simple, as we have
    \begin{equation*}
        \mathcal{C}^\mathrm{hyp}(1) = \mathcal{H}(0), \quad \mathcal{C}^{\mathrm{hyp}}(2) = \mathcal{H}(0,0), \quad \mathcal{C}^\mathrm{hyp}(3) = \mathcal{H}(2), \quad \mathcal{C}^\mathrm{hyp}(4) = \mathcal{H}(1,1).
    \end{equation*}
\end{remark}

It is important to note the difference between hyperelliptic \emph{component} and hyperelliptic \emph{loci}. There are surfaces which admit a hyperelliptic involution such that the corresponding space of quadratic differentials is not of the form given in Theorem \ref{thm:hyperelliptic_loci}. This happens exactly because the hyperelliptic locus has positive codimension. For such surfaces, the results on existence of quadrangulations and well-slanted staircases (Theorem \ref{thm:existence_quadrangulation} and Theorem \ref{thm:existence_well_slanted_staircase}) do not hold in general. We will exhibit a specific example in section \ref{sec:outside_hyperelliptic_strata}.

We can see that, since surfaces in $\mathcal{C}^\mathrm{hyp}(k)$ have total cone angle $2\pi k$, any quadrangulation of them consists of exactly $k$ quadrilaterals. Using the definition of surfaces in a hyperelliptic component $\hyp{k}$ as double covers of quadratic differentials in the stratum $\mathcal{Q}(k-2, -1^{k+1})$, one can show the following two important geometric results. 

\begin{lemma}[Action of hyperelliptic involution on staircases]\label{lem:action_involution_staircases}
    Let $Q$ be a quadrangulation of a surface $X$ in a hyperelliptic component $\hyp{k}$. 
    \begin{enumerate}
        \item Each staircase for $Q$ is fixed (as a set) by the hyperelliptic involuation of $X$.
        \item If $q\in Q$ is a quadrilateral, then its image under the hyperelliptic involution is another quadrilateral that belongs to the same left and right staircases for $Q$ to which $q$ belongs. 
    \end{enumerate}
\end{lemma}
A proof of this result can be found in \cite{delecroix2015diagonal}, where it appears as Lemma 7. This result implies in particular, that any quadrilateral of a quadrangulation of a hyperelliptic surface is either a parallelogram $q$, in which case the $q$ is fixed by the hyperelliptic involution, or the quadrilaterals come in pairs $q_i \neq q_j$ such that $q_i \mapsto q_j$ by the involution. In particular, $q_i$ and $q_j$ must have parallel diagonals. This justifies the assertion that hyperellipticity can be seen as a symmetry constraint.

\begin{lemma}[Cut and paste in hyperelliptic surfaces]\label{lem:cut_paste_hyperelliptic}
    Let $X$ be a surface in a hyperelliptic component $\hyp{k}$ and let $s \colon X \to X$ be the hyperelliptic involution. Let $\gamma$ be a saddle connection in $X$ that is not fixed by $s$. Then, $X \setminus (\Sigma \cup \gamma \cup s\gamma)$ has two connected components, both of them having $\gamma$ and $s\gamma$ on their boundaries. Let $X_1, X_2$ be obtained from these components by identifying $\gamma$ and $s \gamma$ by translations.

    Then, $X_1$ and $X_2$ are nonempty translation surfaces in hyperelliptic components. Moreover, if $k_1 \geq 1$ and $k_2 \geq 1$ are such that $X_1 \in \hyp{k_1}$ and $X_2 \in \hyp{k_2}$, then $k = k_1 + k_2$.
\end{lemma}

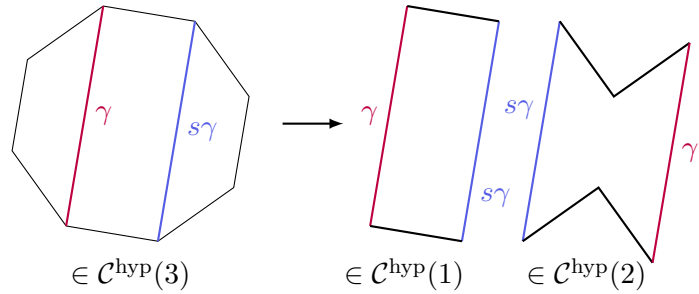
\begin{figure}[ht]
    \centering
    \begin{tikzpicture}[scale = 0.8]
    \pgfmathsetmacro{\s}{2} 
    \pgfmathsetmacro{\rot}{13} 
    
    \draw[rotate=\rot] (0:\s) \foreach \x [count=\i] in {45,90,...,360} {
        -- (\x:\s) coordinate (p\i)
    } -- cycle;

    \draw[thick, noamblue] (p6) -- (p1) node[pos =0.5, right] {$s\gamma$};
    \draw[thick, purple] (p5) -- (p2) node[pos = 0.5,right] {$\gamma$};

    \node at (0,-2.5) {$\in \hyp{3}$};

    \draw[>=latex, ->, thick] (2.5,0) -- (3.5,0);

    \def\x{5}; \def\y{0};

    \draw[thick, purple] ($(p5) + (\x,\y)$) -- ($(p2) + (\x,\y)$) node[pos=0.5, left] {$\gamma$};
    \draw[thick, noamblue] ($(p6) + (\x,\y)$) -- ($(p1) + (\x,\y)$) node[pos=0.2, right] {$s\gamma$}; 
    \draw[thick] ($(p5) + (\x,\y)$) -- ($(p6) + (\x,\y)$);
    \draw[thick] ($(p1) + (\x,\y)$) -- ($(p2) + (\x,\y)$);

    \draw[thick, noamblue] ($(p6) + (\x+1,\y)$) -- ($(p1) + (\x+1,\y)$) node[pos=0.6, left] {$s\gamma$};
    \draw[thick] ($(p6) + (\x+1,\y)$) -- ($(p7) + (\x+1,\y)$) -- ($(p7) + (\x+1,\y) + (p5) - (p4)$);
    \draw[thick] ($(p1) + (\x+1,\y)$) -- ($(p8) + (\x+1,\y)$) -- ($(p8) + (\x+1,\y) + (p2) - (p3)$);
    \draw[thick, purple] ($(p7) + (\x+1,\y) + (p5) - (p4)$) -- ($(p8) + (\x+1,\y) + (p2) - (p3)$) node[pos=0.5, right] {$\gamma$};

    \node at (4.5,-2.5) {$\in \hyp{1}$};
    \node at (7.5,-2.5) {$\in \hyp{2}$};
\end{tikzpicture}
    \caption{The cut and paste procedure described in Lemma \ref{lem:cut_paste_hyperelliptic} applied to a translation surface obtained from a regular octagon.}
    \label{fig:cut_paste_hyperelliptic}
\end{figure}

An example that illustrates Lemma \ref{lem:cut_paste_hyperelliptic} can be seen in Figure \ref{fig:cut_paste_hyperelliptic}.

\begin{proof}
    Let us write $X/s$ for the quotient of the surface by the hyperelliptic involution. The image of $\gamma$ in $X/s$, which coincides with the image of $s\gamma$, is a segment that cannot contain a pole in its interior. This follows from the fact that a saddle connection on $X$ is preserved by $s$ if and only if it contains a Weierstrass point in its interior, so since we assume that $\gamma$ is not fixed by $s$ no interior point of $\gamma$ can be mapped to a pole. We obtain a simple closed curve on the sphere $X/s$, since both ends of the segment are the unique zero of the quadratic differential. Hence, it separates the sphere into two nonempty components whose boundaries consist each of a copy of the image of $\gamma$. Adding a pole each to the interiors of these segments gives two new quadratic differentials. Passing to the orientation covers we obtain the translation surfaces $X_1$ and $X_2$ as in the statement. 

    Lastly, the relation $k = k_1 + k_2$ follows from computing the cone angles. In more detail, separating the sphere along the image of $\gamma$ splits the zero of the quadratic differential, so that in one component we have a zero of order, say, $k_1$ and in the other component we have a zero of order $k-2-k_1$. We conclude by an application of Proposition \ref{prop:properties_orientation_cover}. For example, if both $k$ and $k_1$ are odd, the total cone angles in the orientation covers $X_1$ and $X_2$ are given by $(k_1 + 1)$ and $(k - k_1 -1)$ respectively, which add exactly to $k$. The other cases follow from similar computations.
\end{proof}

A natural question to ask now is how a quadrangulation descends to the sphere via the quotient map given by the hyperelliptic involution. Note that by Lemma \ref{lem:action_involution_staircases}, this question is well-posed since the quadrangulation is fixed by the hyperelliptic involution. Before we answer this question, let us investigate once more the simplest case given by a surface in $\hyp{1} = \mathcal{H}(0)$, i.e., a torus.

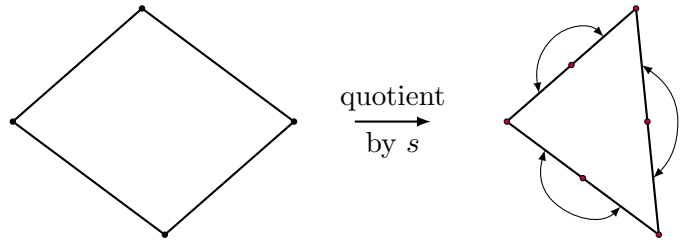
\begin{figure}[ht]
    \centering
    \begin{tikzpicture}
    \draw[thick] (0,0) -- (-2,1.5) -- (-0.3,3) -- (1.7, 1.5) -- cycle;
    \foreach \p in {(0,0), (-2,1.5), (-0.3, 3), (1.7,1.5)} {
        \draw[fill=black] \p circle (1pt);
    }

    \draw[thick, >=latex, ->] (2.5,1.5) -- (3.5,1.5) node[pos=0.5, above] {quotient} node[pos=0.5, below] {by $s$};

    \pgfmathsetmacro{\x}{6.5}

    \draw[thick] (\x,0) -- (\x - 2, 1.5) -- (\x - 0.3,3) -- cycle;

    \foreach \p in {(\x,0), (\x -2, 1.5), (\x - 0.3, 3)}{
        \draw[fill=purple] \p circle (1pt);
    }

    \draw[fill=purple] ($(\x,0)!0.5!(\x-2,1.5)$) circle (1pt);
    \draw[fill=purple] ($(\x,0)!0.5!(\x-0.3,3)$) circle (1pt);
    \draw[fill=purple] ($(\x-0.3,3)!0.5!(\x-2,1.5)$) circle (1pt);

    \draw[>=latex, <->] ($(\x,0)!0.25!(\x-2,1.5)$) to[out=-135, in=-45] ($(\x,0)!0.5!(\x-2,1.5) + (-0.4,-0.3)$) to[out=135, in =-105] ($(\x,0)!0.75!(\x-2,1.5)$);
    \draw[>=latex, <->] ($(\x,0)!0.25!(\x-0.3,3)$) to[out=45, in=-90] ($(\x,0)!0.5!(\x-0.3,3) + (0.4,0)$) to[out=90, in =-35] ($(\x,0)!0.75!(\x-0.3,3)$);
    \draw[>=latex, <->] ($(\x-2,1.5)!0.25!(\x-0.3,3)$) to[out=100, in=225] ($(\x-2,1.5)!0.5!(\x-0.3,3) + (-0.25,0.35)$) to[out=45, in =135] ($(\x-2,1.5)!0.75!(\x-0.3,3)$);
\end{tikzpicture}
    \caption{A simple example of the quotient of a quadrilateral by the hyperelliptic involution $s$.}
    \label{fig:quotient_involution_example}
\end{figure}

\begin{example}[Quotient of quadrilateral by hyperelliptic involution]\label{ex:quotient_involution}
    Suppose we have a hyperelliptic surface $X \in \hyp{1}$ given by a single quadrilateral, as in Figure \ref{fig:quotient_involution_example}. The hyperelliptic involution $s$ can be seen as a rotation by $\pi$. In the quotient, we obtain the triangle on the right-hand side, which belongs to the stratum $\mathcal{Q}(-1^4)$. The singularities are marked in red.
\end{example}

This example motivates the following definitions.

\begin{definition}[Triangle]
    Let $q$ be a quadratic differential on the Riemann sphere $\C\PP^1$ such that $q$ belongs to $\mathcal{Q}(k-1, -1^{k+1})$. Denote by $z_0$ the point on $\C\PP^1$ at which $q$ has the zero of degree $k-2$. We call a \emph{triangle} on $(\C\PP^1, q)$ an open embedded triangle in $(\C\PP^1, q)$ whose boundary consists of saddle connections between $z_0$ and itself that may pass through one pole. 
\end{definition}

\begin{remark}
    Since the cone angle at a pole is $\pi$, an edge passing through a pole consists of two copies of the same segment. In Example \ref{ex:quotient_involution}, identifications are indicated by the double-pointed arrows.
\end{remark}

\begin{definition}[Triangulation]
    A \emph{triangulation} of $(\C\PP^1,q)$ is a set of triangles on $q$ such that their interiors have empty intersection and their union is all of $\C\PP^1$.
\end{definition}

To each triangulation $\tau$ of $(\C\PP^1, q)$ we can associate its dual graph $G_\tau$, which is defined as follows.

\begin{definition}[Dual graph to triangulation]
    Let $\tau$ be a triangulation of $(\C\PP^1, q)$. Its \emph{dual graph} $G_\tau$ is given by the following.
    \begin{enumerate}
        \item The vertices $v_T$ are the triangles $T \in \tau$.
        \item We join two vertices $v_T$ and $v_{T'}$ by an edge if the corresponding triangles $T$ and $T'$ share an edge which has no pole on it.
    \end{enumerate}
\end{definition}

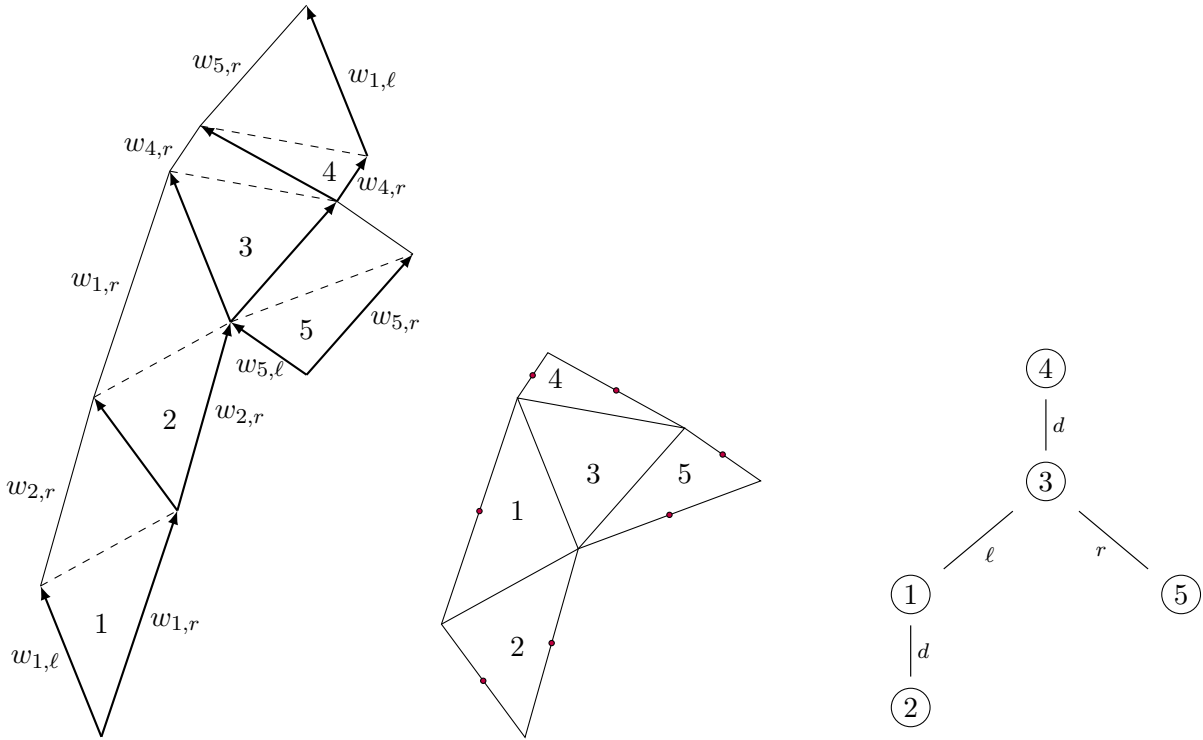
\begin{figure}[ht]
    \centering
    \begin{tikzpicture}
        \coordinate (p1) at (0,0);
        \coordinate (p2) at (-0.8,2);
        \coordinate (p3) at (1,3);
        \coordinate (p4) at (-0.1,4.5);
        \coordinate (p5) at (1.7,5.5);
        \coordinate (p6) at (0.9, 7.5);
        \coordinate (p7) at (3.1, 7.1);
        \coordinate (p8) at (1.3,8.1);
        \coordinate (p9) at (3.5, 7.7);
        \coordinate (p10) at (2.7, 9.7);
        \coordinate (p11) at (2.7, 4.8);
        \coordinate (p12) at (4.1, 6.4);

    \begin{scope}
        \draw[thick, >=latex, ->] (p1) -- (p2) node[midway, left] {$w_{1, \ell}$};
        \draw[thick, >=latex, ->] (p1) -- (p3) node[midway, right] {$w_{1, r}$};
        \draw[thick, >=latex, ->] (p3) -- (p4);
        \draw[thick, >=latex, ->] (p3) -- (p5)node[midway, right] {$w_{2, r}$};
        \draw[thick, >=latex, ->] (p5) -- (p6);
        \draw[thick, >=latex, ->] (p5) -- (p7);
        \draw[thick, >=latex, ->] (p7) -- (p8);
        \draw[thick, >=latex, ->] (p7) -- (p9) node[pos=0.3, right] {$w_{4,r}$};
        \draw[thick, >=latex, ->] (p9) -- (p10)node[midway, right] {$w_{1, \ell}$};
        \draw[thick, >=latex, ->] (p11) -- (p5)node[midway, below, xshift = -3] {$w_{5, \ell}$};
        \draw[thick, >=latex, ->] (p11) -- (p12)node[midway, right, yshift = -3] {$w_{5, r}$};

        \draw (p2) -- node[midway, left] {$w_{2,r}$} (p4) -- node[midway, left] {$w_{1,r}$} (p6) -- node[midway, left] {$w_{4,r}$}(p8) -- node[midway, left] {$w_{5,r}$} (p10);
        \draw (p12) -- (p7);
        
        \draw[dashed] (p2) -- (p3);
        \draw[dashed] (p4) -- (p5);
        \draw[dashed] (p6) -- (p7);
        \draw[dashed] (p8) -- (p9);
        \draw[dashed] (p5) -- (p12);

        \node at ($(p1) + (0,1.5)$) {$1$};
        \node at ($(p3) + (-0.1,1.2)$) {$2$};
        \node at ($(p5) + (0.2,1)$) {$3$};
        \node at ($(p7) + (-0.1,0.4)$) {$4$};
        \node at ($(p11) + (0,0.6)$) {$5$};
    \end{scope} 
\end{tikzpicture}
\begin{tikzpicture}
        \coordinate (p1) at (0,0);
        \coordinate (p2) at (-0.8,2);
        \coordinate (p3) at (1,3);
        \coordinate (p4) at (-0.1,4.5);
        \coordinate (p5) at (1.7,5.5);
        \coordinate (p6) at (0.9, 7.5);
        \coordinate (p7) at (3.1, 7.1);
        \coordinate (p8) at (1.3,8.1);
        \coordinate (p9) at (3.5, 7.7);
        \coordinate (p10) at (2.7, 9.7);
        \coordinate (p11) at (2.7, 4.8);
        \coordinate (p12) at (4.1, 6.4);

        \draw (p3) -- (p4) -- (p6) -- (p8) -- (p7) -- (p12) -- (p5) -- cycle;
        \draw (p4) -- (p5) -- (p6) -- (p7) -- (p5);

        \node at ($(p6) + (0,-1.5)$) {$1$};
        \node at ($(p3) + (-0.1,1.2)$) {$2$};
        \node at ($(p5) + (0.2,1)$) {$3$};
        \node at ($(p8) + (0.1,-0.35)$) {$4$};
        \node at ($(p7) + (0,-0.6)$) {$5$};

        \draw[fill=purple] ($(p3)!0.5!(p4)$) circle (1pt);
        \draw[fill=purple] ($(p4)!0.5!(p6)$) circle (1pt);
        \draw[fill=purple] ($(p6)!0.5!(p8)$) circle (1pt);
        \draw[fill=purple] ($(p8)!0.5!(p7)$) circle (1pt);
        \draw[fill=purple] ($(p7)!0.5!(p12)$) circle (1pt);
        \draw[fill=purple] ($(p12)!0.5!(p5)$) circle (1pt);
        \draw[fill=purple] ($(p3)!0.5!(p5)$) circle (1pt);
\end{tikzpicture}\qquad\qquad
\begin{tikzcd}[baseline=(current bounding box.south)]
	& {\circled{4}} \\
	& {\circled{3}} \\
	{\circled{1}} && {\circled{5}} \\
	{\circled{2}}
	\arrow["d", no head, from=1-2, to=2-2]
	\arrow["\ell", no head, from=2-2, to=3-1]
	\arrow["r"', no head, from=2-2, to=3-3]
	\arrow["d", no head, from=3-1, to=4-1]
\end{tikzcd}
    \caption{A quadrangulation of a surface in $\hyp{5}$, its quotient in $\mathcal{Q}(3,-1^7)$ and the associated dual graph.}
    \label{fig:quadrangulation_triangulation_dual_graph}
\end{figure}

An example of a quadrangulation (of a surface in $\hyp{5})$, its quotient in $\mathcal{Q}(3,-1^7)$ and the associated dual graph can be seen in Figure \ref{fig:quadrangulation_triangulation_dual_graph}.

The dual graph $G_\tau$ to a triangulation $\tau$ to a quadratic differential in a stratum $\mathcal{Q}(k-2, -1^{k+2})$ we used to define hyperelliptic component has a special structure, namely it is a tree.

\begin{lemma}[Dual graph is a tree]\label{lem:dual_graph_is_tree}
    Let $G_\tau$ be the dual graph associated to the triangulation $\tau$ of a quadratic differential $(\C\PP^1, q)$ in a stratum $\mathcal{Q}(k-2, -1^{k+2})$. Then, $G_\tau$ is a tree, i.e., connected and acyclic.
\end{lemma}
\begin{proof}
    The connectedness of $G_\tau$ is a direct consequence of the connectedness of $\C\PP^1$, so it suffices to show that the number of edges of $G_\tau$ is one less than the number of vertices. 

    Since the vertices are by definition given by the triangles of $\tau$, there are $k$ of them. Because it is a triangulation and there are $k+2$ poles, the number of edges is $\frac{1}{2}(3k - (k+2)) = k-1$.
\end{proof}

\begin{figure}[ht]
    \centering
    \begin{tikzpicture}
    \draw[thick] (0,0) -- node[below right]{$d$} (4,-0.5) -- node[right]{$\ell$} (1.8, 3) -- node[left]{$r$} cycle;
    \draw[thick] (5,2.8) -- node[above left] {$d$} (9,2.5) -- node[right] {$r$} (7,-0.5) -- node[left] {$\ell$} cycle;

    \draw[noamblue, dashed] (1.8,3) -- ++(-90:3.5);
    \draw[noamblue, dashed] (7,-0.5) -- ++(90:3.5);
\end{tikzpicture}
    \caption{Labeling triangles obtained as quotients of quadrilaterals under the hyperelliptic involution.}
    \label{fig:triangle_labeling}
\end{figure}
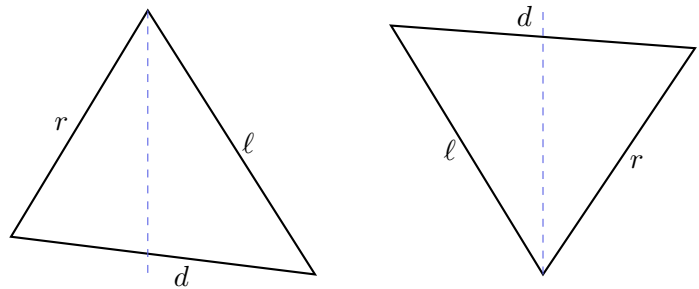

The restriction that quadrilaterals in a quadrangulation need to be admissible translates to a compatibility condition on the triangles as we explain now. Any triangle has exactly one vertex $v$, such that the vertical coming from this vertex is contained in its interior. We assign labels to the edges of the triangle as follows. We assign the label $d$ to the side opposite of the vertex $v$ and we label $\ell$ and $r$ the other two sides of the triangle such that rotating counterclockwise around the vertex, one sees first the side labelled $r$, then the vertical, then the side labelled $\ell$, as shown in Figure \ref{fig:triangle_labeling}.

Admissibility of the quadrilaterals then implies that sides of the triangles that are identified carry the same label, i.e., a side labeled $\ell$ is glued to another side labeled $\ell$ and analogously for the other labels. Because of this, we may define an edge labelling on the associated dual graph $G_T$ by simply labelling an edge of the graph with the corresponding label of the identified edges of the triangle it represents. This labeling can be seen in Figure \ref{fig:quadrangulation_triangulation_dual_graph}. This labeled tree is called the \emph{tree of relations} in \cite{ferenczi2010structure}, where the labels $\hat{+}, \hat{-}$ and $\hat{=}$ are used instead of $\ell, r$ and $d$. In the same paper, Ferenczi and Zamboni also prove the following proposition.

\begin{proposition}[\cite{ferenczi2010structure}, tree of relations encode train-tracks]
    A tree of relations encodes the train-track relations for the length datum $\boldsymbol{w}$ of a quadrangulation $Q$ as follows. If the edge of the tree connecting the vertices $i$ and $j$ carries the label $r$ (or $\ell$), the wedges of the quadrilaterals $q_i$ and $q_j$ obtained by double covers of the triangle $T_i, T_j$ corresponding to the vertices $i$ and $j$ are such that
    \begin{equation*}
        w_{i,r} = w_{j,r} \quad (\text{or } w_{i, \ell} = w_{j,\ell}).
    \end{equation*}
    If the edge connecting $i$ and $j$ carries the label $d$, then
    \begin{equation*}
        w_{i,d} = w_{j,d}.
    \end{equation*}
    These sets of equations are equivalent to the train track relations \eqref{eq:wedges_train_tracks}. 
\end{proposition}

The preimage of a triangle $T$ in $Q$ of the quotient map given by the hyperelliptic involution is a union of a top and a bottom part of a quadrilateral. In order to define a labeling on the triangles, we endow $T$ with the label $i$ if it corresponds to the \emph{bottom} part of $q_i \in Q$. Using this labeling, we can describe the tree of relations fully by specifying a triple $(\sigma_\ell, \sigma_r, \sigma_d)$ of involutions of $[k]$. For $\varepsilon\in\{\ell, r, d\}$, if the nodes $i$ and $j$ share an edge labeled $\varepsilon$, we set $\sigma_\varepsilon(i) = j$ and $\sigma_\varepsilon(j) = i$. For example, the tree in Figure \ref{fig:quadrangulation_triangulation_dual_graph} is encoded by
\begin{align}\label{eq:triangle_combinatorial_datum}
\begin{split}
    \sigma_\ell &= \begin{pmatrix}
        1 & 3
    \end{pmatrix}, \\
    \sigma_r &= \begin{pmatrix}
        3 & 5
    \end{pmatrix}, \\
    \sigma_d &= \begin{pmatrix}
        1 & 2
    \end{pmatrix}\begin{pmatrix}
        3 & 4
    \end{pmatrix},
\end{split}
\end{align}
using the usual cycle notation. The following lemma shows how $\boldsymbol{\sigma} = (\sigma_\ell, \sigma_r, \sigma_d)$ corresponds to the combinatorial datum $\boldsymbol{\pi} = (\pi_\ell, \pi_r)$ or equivalently, it describes the connection between the tree of relations $G_\tau$ and the graph $G_Q$ introduced in section \ref{sec:space_of_quadrangulations}.

\begin{lemma}[Correspondence between combinatorial data of quadrilaterals and triangles]\label{lem:correspondence_combinatorial_data}
    If $\tau$ is a labeled triangulation with combinatorial datum $\boldsymbol{\sigma} = (\sigma_\ell, \sigma_r, \sigma_d)$ induced by a labeled quadrangulation of a surface in $\hyp{k}$ with combinatorial datum $\boldsymbol{\pi} = (\pi_\ell, \pi_r)$, then $\sigma_d = s$ is the action of the hyperelliptic involution on the quadrilaterals of $Q$, i.e., the hyperelliptic involution $s$ maps $q_i$ to $q_{\sigma_d(i)}$ and
    \begin{equation}\label{eq:correspondence_combinatorial_data}
        \pi_\ell = \sigma_r \circ \sigma_d \quad \text{ and } \quad \pi_r = \sigma_\ell \circ \sigma_d.
    \end{equation}
    In particular, $\pi_\ell^{-1} = s\pi_\ell s$ and $\pi_r^{-1} = s \pi_r s$.
\end{lemma}
\begin{proof}
    It follows directly from the construction of the labeling on the tree of relations that $\sigma_d$ corresponds to the action of the hyperelliptic involution $s$. Quadrilaterals $q$ in $Q$ are decomposed into triangles along the backwards diagonal. To reconstruct $\pi_r$ from a bottom triangle, we first need to cross the diagonal, hence apply $\sigma_d$ and then cross the top right side, which according to the labeling convention corresponds to applying $\sigma_\ell$, see Figure \ref{fig:triangle_labeling}. Therefore, $\pi_r = \sigma_r\sigma_d$. The argument for $\pi_\ell$ is completely analogous.

    Finally, the last statement follows from
    \begin{equation*}
        \sigma_d \pi_\ell \sigma_d \pi_\ell = \sigma_d\sigma_r\sigma_d\sigma_d\sigma_r\sigma_d = \operatorname{id},
    \end{equation*}
    where the first equality follows from applying the left equation of \eqref{eq:correspondence_combinatorial_data} twice and the second equality follows from the fact that $\sigma_\varepsilon$ is an involution for $\varepsilon\in\{\ell, r, d\}$. A similar computation leads to the equation for $\pi_r$.
\end{proof}

Let us check that the combinatorial data $\boldsymbol{\pi}$ and $\boldsymbol{\sigma}$ corresponding to the quadrangulation and triangulation in Figure \ref{fig:quadrangulation_triangulation_dual_graph} satisfy the equations in Lemma \ref{lem:correspondence_combinatorial_data}. Note first, that the involution exchanges the quadrilaterals $q_1$ and $q_2$ as well as $q_3$ and $q_4$, while fixing $q_5$. This does indeed correspond exactly to $\sigma_d$ as we can see in \eqref{eq:triangle_combinatorial_datum}. Using the same equations and the fact that
\begin{equation*}
    \boldsymbol{\pi} = (\pi_\ell, \pi_r) = \Big(\begin{pmatrix}
        1 & 2
    \end{pmatrix}\begin{pmatrix}
        3 & 4 & 5
    \end{pmatrix},
    \begin{pmatrix}
        1 & 2 & 3 & 4
    \end{pmatrix}\Big),
\end{equation*}
we can further compute that
\begin{align*}
    \sigma_r \circ \sigma_d = \begin{pmatrix}
        3 & 5
    \end{pmatrix}\circ \begin{pmatrix}
        1 & 2
    \end{pmatrix}\begin{pmatrix}
        3 & 4
    \end{pmatrix}&=
    \begin{pmatrix}
        1 & 2
    \end{pmatrix}\begin{pmatrix}
        3 & 4 & 5
    \end{pmatrix} = \pi_\ell, \\
    \sigma_\ell \circ \sigma_d = \begin{pmatrix}
        1 & 3
    \end{pmatrix}\circ
    \begin{pmatrix}
        1 & 2
    \end{pmatrix}\begin{pmatrix}
        3 & 4
    \end{pmatrix}
    &=
    \begin{pmatrix}
        1 & 2 & 3 & 4
    \end{pmatrix} = \pi_r,
\end{align*}
as well as
\begin{align*}
    \sigma_d \circ \pi_\ell \circ \sigma_d &= \begin{pmatrix}
        5 & 4 & 3
    \end{pmatrix}\begin{pmatrix}
        2 & 1
    \end{pmatrix} = \pi_\ell^{-1}, \\
    \sigma_d \circ \pi_r \circ \sigma_d &= 
    \begin{pmatrix}
        4 & 3 & 2 & 1
    \end{pmatrix} = \pi_r^{-1}.
\end{align*}

One can show that $\sigma_\ell\sigma_r\sigma_d$ is always a $k$-cycle. Geometrically, it corresponds to \enquote{turning around the singularity of angle $\pi k$ on the sphere}. This $k$-cycle is invariant under staircase moves. In fact, it is a \emph{complete} invariant in the sense of the following theorem.

\begin{theorem}\label{thm:complete_invariant_k_cycle}
    Let $Q$ be a quadrangulation of a surface in $\hyp{k}$ and let $\tau_Q$ be the quotient triangulation. Let $\boldsymbol{\pi} = (\pi_\ell, \pi_r)$ and $\boldsymbol{\sigma} = (\sigma_\ell, \sigma_r)$ be the combinatorial data. Then, the permutation
    \begin{equation*}
        \sigma_\ell\sigma_r\sigma_d = \pi_r \sigma_d \pi_\ell
    \end{equation*}
    is a $k$-cycle, which is invariant under the operation of staircase moves $\boldsymbol{\pi} \mapsto c \cdot \boldsymbol{\pi}$.

    Moreover, two combinatorial data $\boldsymbol{\pi}$ and $\boldsymbol{\pi}'$ corresponding to quadrangulations of surfaces in $\hyp{k}$ can be joined by staircase moves, hence they belong to the same $DC$ graph $\mathcal{G} = \mathcal{G}(\boldsymbol{\pi})$ if and only if
    \begin{equation*}
        \pi_r\pi_\ell\sigma_d = \pi'_r\pi'_\ell\sigma_d'.
    \end{equation*}
\end{theorem}
A proof of this result can be found in \cite{cassaigne2011combinatorial}. The following corollary shows that the rotation operator $R \colon \mathcal{Q}_k \to \mathcal{Q}_k$ from Definition \ref{def:rotation_operator} is well-defined.

\begin{corollary}\label{cor:rotation_well_defined}
    Let $\boldsymbol{\pi} = {(\pi_\ell, \pi_r)}$ be a combinatorial datum of a labeled quadrangulation $Q$ in $\mathcal{Q}_k$ and let $\boldsymbol{\pi}' = (\pi_\ell\pi_r\pi_\ell^{-1}, \pi_\ell^{-1})$. Then, $\boldsymbol{\pi}'$ belongs to $\mathcal{G}(\boldsymbol{\pi})$.
\end{corollary}

\begin{proof}
    Consider the quadrangulation $Q' = RQ$, where $R$ is the rotation operator from Definition \ref{def:rotation_operator}. Note that applying $R$ to $Q$ is well-defined, the only thing that is not clear a priori is whether $Q' \in \mathcal{Q}_k$. By definition of $R$, in particular by equation \eqref{eq:rotation_3}, we see that the action $s'$ of the hyperelliptic involution on the labels of $Q'$ is given by
    \begin{equation*}
        s' = \pi_\ell s \pi_\ell^{-1}.
    \end{equation*}
    Using this and the equality $s\pi_\ell^{-1} = \pi_\ell s$ from Lemma \ref{lem:correspondence_combinatorial_data} we compute
    \begin{equation*}
        \pi'_\ell\pi'_rs' = \left(\pi_\ell\pi_r\pi_\ell^{-1}\right) \pi_\ell^{-1}\left(\pi_\ell s \pi_\ell^{-1}\right) = 
        \pi_\ell\pi_r\pi_\ell^{-1}\left(s \pi_\ell^{-1}\right) =
        \pi_\ell\pi_r\pi_\ell^{-1}\pi_\ell s = \pi_\ell\pi_r s.
    \end{equation*}
    By Lemma \ref{lem:correspondence_combinatorial_data} we have $s = \sigma_d$ and $s' = \sigma'_d$, so that by Theorem \ref{thm:complete_invariant_k_cycle} it follows that $\boldsymbol{\pi}' \in \mathcal{G}(\boldsymbol{\pi})$.
\end{proof}

\subsubsection{Proof of Existence of Quadrangulations and Well-Slanted Staircases}

Lastly, we want to outline the proof of the fact that for any surface $X \in \hyp{k}$ admits a quadrangulation and in any quadrangulation of a surface $X \in \hyp{k}$, there exists at least one well-slanted staircase in $Q$.

The following lemma is true for \emph{any} translation surface, not necessarily in a hyperelliptic component. 

\begin{lemma}[Infinitely many best approximations]\label{lem:infinitely_many_approximations}
    Let $X$ be a translation surface which has neither horizontal nor vertical saddle connections. Then, in any bundle of $X$, there are infinitely many left and right best approximations. Moreover, for any bundle $\Gamma_i$ of $X$ we have
    \begin{equation*}
        \min\{\IM(w) \mid w \in \Gamma_i \text{ is a best approximation and } |\RE(w)| < r\} < \frac{\operatorname{Area(X)}}{r},
    \end{equation*}
    where $\operatorname{Area}(X)$ is the area of $X$ induced by the Lebesgue measure on $\R^2$.
\end{lemma}

The first part of Lemma \ref{lem:infinitely_many_approximations} simply tells us that in any outgoing half-plane associated to a singularity, there are infinitely many best approximations. This can be shown by arguments similar to the proof of Keane's Theorem (Theorem \ref{thm:keanes}). If $w$ is a best approximation, the quantity $|\RE(w)|\cdot \IM(w)$ is called the \emph{area of the best approximation} $w$, as this is exactly the area of $R(w)$, the immersed rectangle as defined in Lemma \ref{lem:immersed_rectangles}. The second part of Lemma \ref{lem:infinitely_many_approximations} then states that the area of best approximation is uniformly bounded from above by $\operatorname{Area}(X)$. Indeed, if there was a best approximation $w$ such that $\IM(w) \cdot |\RE(w)| \geq \operatorname{Area}(X)$, we could set $r \coloneqq |\RE(w)|$ and obtain, by the lemma, some best approximation $v$ with $|\RE(v)|<r$ and
\begin{equation*}
    \IM(v) < \frac{\operatorname{Area(X)}}{r}.
\end{equation*}
In particular, since $w$ is a best approximation, the fact that $\RE(v) < \RE(w)$ then implies $\IM(v) > \IM(w)$, so that 
\begin{equation*}
    \operatorname{Area}(X) \leq \IM(w)\cdot |\RE(w)| < \IM(v) \cdot r < \operatorname{Area}(X)
\end{equation*}
gives a contradiction.

We are now ready to state and prove the two main results in this section.

\begin{theorem}[Existence of quadrangulation]\label{thm:existence_quadrangulation}
    Let $X$ be a surface in a hyperelliptic component $\hyp{k}$ that admits neither horizontal nor vertical saddle connections. Then, $X$ admits a quadrangulation.
\end{theorem}
\begin{proof}
    First, pick any singularity of $X$ and choose a saddle connection $w$ which is a best approximation, which exists by Lemma \ref{lem:infinitely_many_approximations}. By Lemma \ref{lem:diagonals_determine_quadrilaterals}, there exists a unique quadrilateral $q$ which has $w$ as diagonal. Let us write $v_\ell$ and $v_r$ for the bottom left and right sides and $v^\ell$ and $v^r$ for the top left and right sides. 

    Consider $s(q)$, the image of $q$ under the hyperelliptic involution $s$. By Lemma \ref{lem:action_involution_staircases}, there are only two possible cases.
    \begin{enumerate}
        \item $q = s(q)$
        \item The interiors of $q$ and $s(q)$ are disjoint.
    \end{enumerate}
    In both cases, we proceed as follows. For each side $v$ of $q$, if $v = s(v)$ we do nothing, while if $v \neq s(v)$ we apply the cut and paste procedure introduced in Lemma \ref{lem:cut_paste_hyperelliptic}. After these operations, we obtain one surface that consists of one (if $q = s(q)$) or two (if $q \neq s(q)$) quadrilaterals and up to four surfaces $X_\ell, X_r, X^\ell$ and $X^r$ that contain respectively $v_\ell, v_r, v^\ell$ and $v^r$. We consider one of these surfaces empty, exactly if $v = s(v)$ for one of the sides. Let us remark that these surfaces may coincide also in the case where they are nonempty. Also by Lemma \ref{lem:cut_paste_hyperelliptic}, each of the (nonempty) surfaces we obtain belong to a hyperelliptic component with strictly smaller total angle. 

    Now we can repeat this construction on each nonempty surface among $X_\ell, X_r, X^\ell$ and $X^r$, i.e., choose a saddle connection that is a best approximation using Lemma \ref{lem:infinitely_many_approximations}, complete it to an admissible quadrilateral by Lemma \ref{lem:diagonals_determine_quadrilaterals} and apply again the cut and paste procedure outlined above. After finitely many steps, this produces $k$ admissible quadrilaterals which give exactly the claimed quadrangulation of $X$.
\end{proof}
\begin{remark}
    This proof does \emph{not} extend to other components of strata. While it is still possible to use a cut and paste construction, the resulting surfaces $X_\ell, X_r, X^\ell$ and $X^r$ in the notation of the proof might be connected to each other. In particular, if two of them are connected then we would obtain a surface in which we would need to complete a set of two saddle connections, i.e., two distinct diagonals, into a quadrangulation.
\end{remark}

Note that Theorem \ref{thm:existence_quadrangulation} is constructive. In fact, Figure \ref{fig:cut_paste_hyperelliptic} depicts one step of the construction in the Theorem, where we are in the case $q = s(q), X_\ell = X^r = \emptyset$ and $X_r = X^\ell$ belongs to $\hyp{2}$. Continuing the construction with this surface one can then, by choosing the appropriate saddle connection, obtain exactly the quadrangulation from Figure \ref{fig:quadrangulation_octagon}.

\begin{theorem}[Existence of well-slanted staircase]\label{thm:existence_well_slanted_staircase}
    Let $Q$ be a quadrangulation of a surface $X$ in a hyperelliptic component $\hyp{k}$ and assume that no quadrilateral in $Q$ has a vertical diagonal. Then, there exists at least one well-slanted staircase in $Q$.
\end{theorem}
\begin{proof}
    The proof of the theorem proceeds by (strong) induction on the number of quadrilaterals $k$, or equivalently on the integer $k$ such that the surface $X$ belongs to $\hyp{k}$. The base case $k = 1$ is trivial since a staircase made of one single quadrilateral is always well-slanted.

    Suppose now that the claim holds for $\{1, 2, \ldots, k-1\}$ and let $Q = (\boldsymbol{\pi}, \boldsymbol{w})$ be a quadrangulation of a surface in $\hyp{k}$. Let us write $s$ both for the hyperelliptic involution itself and for the action on the quadrilaterals, i.e., we write $s(i) = j$ if and only if $q_j = s(q_i)$. Towards a contradiction, assume that there exists no well-slanted staircase in $Q$ so that no staircase move for $Q$ is possible. We claim that there exists a right staircase $S$ which contains both left- and right-slanted quadrilaterals. Indeed, there cannot exist a right staircase consisting of only right-slanted quadrilaterals, since then we would be able to perform a right staircase move. Similarly, there cannot be a right staircase consisting only of left-slanted quadrilaterals, since then we could do a left staircase move. Since there is always at least one right staircase, it follows that there must indeed exist a staircase $S$ with both left- and right-slanted staircases.

    In particular, there exist two consecutive quadrilaterals $q_i$ and $q_{\pi_r(i)}$ such that the first is left-slanted and the second is right-slanted. Note that $s(i) \neq \pi_r(i)$, for if this equality would hold, by Lemma \ref{lem:action_involution_staircases} the diagonals of $q_i$ and $q_{\pi_r(i)} = s(q_i)$ would be parallel, therefore $q_i$ and $q_{\pi_r(i)}$ would be slanted alike. It follows in particular that $w_{\pi_r(i), \ell}$, the common edge of $q_i$ and $q_{\pi_r(i)}$ does not contain a Weierstrass point and $w_{\pi_r(i), \ell} \neq s(w_{\pi_r(i), \ell})$.

    Therefore, we may apply Lemma \ref{lem:cut_paste_hyperelliptic} to the saddle connection $w_{\pi_r(i), \ell}$ and obtain two connected components that give, after identifying on each of the components the appropriate boundary elements, two quadrangulations of surfaces in $\hyp{k_1}$ and $\hyp{k_2}$ respectively, where $k_1, k_2 < k$. More precisely, we cut along the edges $w_{\pi_r(i), \ell}$ and $w_{s(i), \ell}$. Let us write $X'$ for the surface containing $q_i$. By the induction hypothesis, there exists a staircase move in $X'$. Since $q_i$ is left-slanted, the same is true for $q_{s(i)}$ so that applying this staircase move does not change the saddle connection $w_{s(i), \ell}$. Hence, after the move we can glue back the two components into the surface $X$, meaning that we were able to lift the staircase move to $X$, contradicting our assumption that no such move is possible for $Q$. We conclude that there must exist a well-slanted staircase in $Q$.
\end{proof}
Together with Lemma \ref{lem:diagonals_become_sides}, we have shown that under Keane's condition, algorithms given by staircase moves are always well-defined in hyperelliptic components $\hyp{k}$ for $k \geq 1$.

\subsection{Outside of Hyperelliptic Components}\label{sec:outside_hyperelliptic_strata}

Let us finish this section by discussing the situation where the translation surface $X$ does not belong to $\hyp{k}$ for any integer $k$. We will start with the example mentioned in the previous section of a hyperelliptic surface with this property that admits a quadrangulation without any well-slanted staircase. This particular example is taken from \cite{delecroix2015diagonal}. Both examples we will present are constructed using the strategy of finding a quadrangulation consisting of $k$ quadrilaterals such that $\pi_\ell$ and $\pi_r$ both are $k$-cycles. This kind of construction is impossible in $\hyp{k}$ if $k \geq 3$, since the associated dual graph would no longer be a tree, which is impossible due to Lemma \ref{lem:dual_graph_is_tree}.

\begin{example}[No well-slanted staircase]
    The stratum we consider is $\mathcal{H}(0,0,0)$, which is the smallest stratum that does not contain a hyperelliptic component. Let
    \begin{equation*}
        \pi_\ell =\pi_r = \begin{pmatrix}
            1 & 2 & 3
        \end{pmatrix},
    \end{equation*}
    and consider the wedges
    \begin{align*}
        w_{1, \ell} &= (-2, 2),\quad  w_{1,r} = (1,1), \\
        w_{2, \ell} &= w_{3, \ell} = (-1.3, 2), \\
        w_{2, r} &= w_{3,r} = (1.7,1).
    \end{align*}
    A quick computation verifies that the wedges satisfy the train-track relations for $\boldsymbol{\pi} = (\pi_\ell, \pi_r)$ so that indeed we obtain a quadrangulation $Q$. Moreover, we have
    \begin{equation*}
        w_{1,d} = (-0.3,3), \quad w_{2,d} = (0.4,3), \quad w_{3,d} = (-0.3,3).
    \end{equation*}
    In particular, the first component of the diagonals have different signs, meaning that there are both left- and right-slanted quadrilaterals and hence no well-slanted staircase in $Q$. Note that $Q$ does admit a hyperelliptic involution which fixes $q_2$ and exchanges $q_1$ and $q_3$, i.e., the quadrangulation $Q$ is fixed by the hyperelliptic involution. The quotient of the surface by the involution belongs to $\mathcal{Q}(0^2,-1^4)$, which shows that the surface does \emph{not} belong to a hyperelliptic component.
\end{example}

In \cite{delecroix2015diagonal}, there is another example of a surface belonging to $\mathcal{H}(4)$ which admits a quadrangulation into 5 quadrilaterals without any well-slanted staircase. The stratum $\mathcal{H}(4)$ is the smallest one which contains more than one component, one of which is hyperelliptic. Let us provide another example of a surface in $\mathcal{H}(0,0,0,0)$ which is not hyperelliptic.

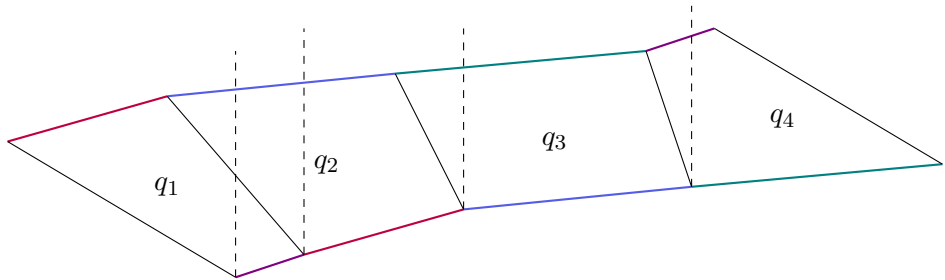
\begin{figure}[ht]
    \centering
    \begin{tikzpicture}[scale = 0.3]
    \coordinate (p1) at (0,0);
    \coordinate (p2) at (-10,6);
    \coordinate (p3) at (3,1);
    \coordinate (p4) at (-3,8);
    \coordinate (p5) at (10,3);
    \coordinate (p6) at (7,9);
    \coordinate (p7) at (20,4);
    \coordinate (p8) at (18,10);
    \coordinate (p9) at (31,5);
    \coordinate (p10) at (21,11);

    \draw[thick, purple] (p2) -- (p4); \draw[thick, purple] (p3) -- (p5);
    \draw[thick, noamblue] (p4) -- (p6); \draw[thick, noamblue] (p5) -- (p7);
    \draw[thick, teal] (p6) -- (p8); \draw[thick, teal] (p7) -- (p9);
    \draw[thick, violet] (p8) -- (p10); \draw[thick, violet] (p1) -- (p3);

    \draw[] (p3) -- (p4); 
    \draw[] (p5) -- (p6);
    \draw[] (p7) -- (p8);
    \draw (p1) -- (p2);
    \draw (p9)-- (p10);

    \node at (-3,4) {$q_1$};
    \node at (4,5) {$q_2$};
    \node at (14,6) {$q_3$};
    \node at (24,7) {$q_4$};

    \draw[dashed] (p1) -- ++(90:10);
    \draw[dashed] (p3) -- ++(90:10);
    \draw[dashed] (p5) -- ++(90:8);
    \draw[dashed] (p7) -- ++(90:8);
    
\end{tikzpicture}
    \caption{A quadrangulation of a surface in $\mathcal{H}(0,0,0,0)$ admitting \emph{no} well-slanted staircase.}
    \label{fig:quadrangulation_no_well_slanted_staircase}
\end{figure}

\begin{example}[Quadrangulation without well-slanted staircase]\label{ex:no_staircase}
    We define a translation surface $X \in \mathcal{H}(0,0,0,0)$ by defining the following quadrangulation $Q$. Let
    \begin{equation*}
        \pi_\ell = \pi_r = \begin{pmatrix}
            1 & 2 & 3 & 4
        \end{pmatrix},
    \end{equation*}
    and further define the wedges
    \begin{align*}
        w_{1, \ell} &= (-10,2), \quad w_{1, r} = (3,1), \\
        w_{2, \ell} &= (-6,2) , \quad w_{2, r} = (7,2), \\
        w_{3, \ell} &= (-3, 2), \quad w_{3, r} = (10,1), \\
        w_{4, \ell} &= (-2, 2), \quad w_{4, r} = (11, 1).
    \end{align*}
    It is again quick to verify that these wedges satsify the train-track relations for $\boldsymbol{\pi} = (\pi_\ell, \pi_r)$ and that $q_1$ is right-slanted and $q_i$ is left-slanted for $i \in \{2,3,4\}$. Therefore, there is no well-slanted staircase in $Q$. The quadrangulation is depicted in Figure \ref{fig:quadrangulation_no_well_slanted_staircase}.
\end{example}

 \pagebreak
\thispagestyle{empty}
\
\pagebreak
\section{Enumerating Closed Geodesics for the Teichmüller Flow}\label{sec:counting}
\thispagestyle{plain}

One question of particular interest in the general theory of dynamical systems is the question of the existence and abundance of closed orbits. More precisely, one is for instance interested in the asymptotic growth of the number of closed orbits, which is known as the \emph{counting problem}. In \cite{Eskin2008CountingCG}, Eskin and Mirzakhani showed that in the case of the Teichmüller geodesic flow $g_t$, the asymptotic growth of the number of closed Teichmüller geodesics with length up to $T$ grows as $\nicefrac{\e^{\mathrm{h}T}}{\mathrm{h}T}$, where $\mathrm{h}$ denotes the entropy of the Teichmüller geodesic flow. Later, analogous results were also proven for closed geodesics in a fixed connected component of a stratum $\stratum$ by Eskin, Mirzakhani and Rafi in \cite{eskin2019counting} and independently by Hamenstäd in \cite{hamenstaedt2013bowen}. Another related problem is to systematically list all closed Teichmüller geodesics in a given connected component of a stratum, for example in the order given by the length of the geodesic. We will refer to this problem as the \emph{enumeration problem}.

The importance of studying closed Teichmüller geodesics can further be seen by establishing a connection between the dynamics of translation surfaces and the seminal work on surfaces by Thurston, see e.g., \cite{fathi1979travaux}. As we will explain below, there is a one-to-one correspondence between closed Teichmüller geodesics and so-called \emph{pseudo-Anosov maps} appearing in the theory of Thurston. The importance of these kind of maps comes from the \emph{Nielsen--Thurston classification.}

\begin{theorem}[Nielsen--Thurston Classification]
    Let $X$ be a compact connected orientable surface. Any homeomorphism $f \colon X \to X$ is isotopic to a homeomorphism $g \colon X \to X$, such that 
    \begin{enumerate}
        \item $g$ is \emph{periodic}, i.e., some iterate $g^n$ is the identity, or
        \item $g$ is \emph{reducible}, i.e., it preserves some finite union of closed simple curves on $X$, or
        \item $g$ is \emph{pseudo-Anosov}.
    \end{enumerate}
    Moreover, if $g$ is pseudo-Anosov, it is neither reducible nor periodic.
\end{theorem}
\begin{remark}
    If $f$ is (isotopic to a map that is) reducible, it is possible to decompose $f$ into irreducible pieces (i.e., into homeomorphisms preserving some subsurface with boundaries) and apply the Nielsen--Thurston Classification to these irreducible maps, see for example Corollary 13.3 in \cite{farb2011primer}. The pseudo-Anosov maps, which are defined below in Definition \ref{def:pseudo_anosov_map}, are therefore the most interesting self-homeomorphisms of a surface. 
\end{remark}    

An important source for this section is Farb and Margalits book \cite{farb2011primer}. A further very approachable exposition of many concepts presented here can be found in \cite{lanneau2017tell}. For the part about the correspondence between closed Teichmüller geodesics and pseudo-Anosov mapping classes, many ideas are taken again from \cite{farb2011primer} as well as from \cite{athreyamasur2023translationsurfaces} and \cite{massart2022short}, where one can find further details. 

\subsection{Connection between Closed Teichmüller Geodesics and Pseudo-Anosov Mapping Classes}

Before we can define pseudo-Anosov maps we need some preliminaries. There are several equivalent ways to define pseudo-Anosov maps, the definition we give here makes use of what are called transverse singular measured foliations, which is a quite general framework do define the concept of pseudo-Anosov. Another approach would be to assume that the surface $X$ comes equipped with a (half-) translation structure, which might seem a lot more restrictive, but turns out to be equivalent to the definition via measured foliations, as we will outline below.  

\subsubsection{Examples of Diffeomorphisms on the Torus induced by Matrices}
Let us present an illustrative example that encapsulates the main idea before giving the general definition of a pseudo-Anosov map. Consider any matrix $A \in \operatorname{SL}(2, \Z)$, which acts linearly on the plane $\R^2$. Since $A(\Z^2) = \Z^2$, the matrix induces a diffeomorphism $\psi$ of the torus $\T^2 = \R^2 / \Z^2$ given by
\begin{equation*}
    \psi(x,y) = A(x,y) \quad \pmod{\Z^2}.
\end{equation*}
The dynamical properties of $\psi$ depend on the eigenvalues $\lambda, \lambda^{-1}$ of $A$, where three possibilities arise:
\begin{enumerate}
    \item $\overline{\lambda} = \lambda^{-1}$ and $\lambda \neq \pm 1$. In this case, we say $\psi$ is of \emph{finite order}. 
    \item $\lambda = \lambda^{-1} = \pm 1$. Here, we say that $\psi$ is \emph{reducible}, meaning that it preserves a closed curve on the torus.
    \item In every other case, when $\lambda$ and $\lambda^{-1}$ are distinct real numbers, we say that $\psi$ is \emph{Anosov}.
\end{enumerate}

Equivalently, we can state the above possibilities in terms of the trace of $A$, i.e., if $|\operatorname{tr}(A)| \in \{0, 1\}$ the map $\psi$ is of finite order, if $|\operatorname{tr}(A)| = 2$ it is reducible and if $|\operatorname{tr}(A)| > 2$ it is Anosov. The terminology is appropriate in the following sense. Since the characteristic polynomial of $A$ is given by $x^2 - \operatorname{tr}(A)x + 1$, if $\tr(A) \in \{0, 1\}$ it follows by the Cayley--Hamilton theorem that $A$ has finite order, i.e., $A^k = \operatorname{Id}$ for some $k \in \N_{\geq 1}$. In fact, it is even the case that $k \in \{2, 3, 4, 6\}$. If $|\tr(A)| = 2$, we can cut the torus along the closed curve and study the map on the simpler structure we obtain this way, which here is an annulus. Let us give explicit examples of these cases.

\begin{figure}[hb]
    \centering
    \begin{tikzpicture}
    \draw[] (0,0) -- (0,2); \draw[] (2,0)-- (2,2);
    \draw[thick, purple] (0,0)-- (2,0); \draw[thick, purple] (0,2) -- (2,2);

    \def \x{3}

    \draw[] (\x,0) -- (\x+2, 2); \draw[] (\x+2, 0) -- (\x+4, 2);
    \draw[thick, purple] (\x+0,0) -- (\x+2,0); \draw[thick, purple] (\x+2,2) -- (\x+4,2);
    \draw[dashed] (\x+2,2) -- (\x+2,0);

    \def \w{8}
    \draw[] (\w,0) -- (\w,2); \draw[] (\w+2,0)-- (\w+2,2);
    \draw[thick, purple] (\w,0)-- (\w+2,0); \draw[thick, purple] (\w,2) -- (\w+2,2);

    \draw[>=latex, ->] (1,2.2) to[out=20, in =180] (\x+2,2.7) to[out = 0, in = 160] (\w+1, 2.2);
    \node[below] at (\x+2, 2.7) {$\psi_B$};
\end{tikzpicture}

\vspace{1cm}
\begin{tikzpicture}[scale=0.65]
    \begin{scope}
    \clip (-3,-1.5) rectangle (3,1.5);
    \draw (0,0) ellipse (3 and 1.5);
    \begin{scope}
        \clip (0,-1.8) ellipse (3 and 2.5);
        \draw (0,2.2) ellipse (3 and 2.5);
    \end{scope}
    \begin{scope}
        \clip (0,2.2) ellipse (3 and 2.5);
        \draw (0,-2.2) ellipse (3 and 2.5);
    \end{scope}
    \draw[purple, thick, dashed, rotate = 5] (0.5,-0.91) ++(270:0.3 and 0.6) arc (270:90:0.3 and 0.6);
    \draw[purple, thick, rotate = 5] (0.5,-0.91) ++(90:0.3 and 0.6) arc (90:-90:0.3 and 0.6);
    \end{scope}
    \draw[>=latex, ->] (3.5,0) -- (4.5,0);

    \draw[thick, purple] (6.25,1) ellipse (1 and 0.3);
    \draw[thick, purple, dashed] (6.25, -1) ++(0:1 and 0.3) arc (0:180:1 and 0.3);
    \draw[thick, purple] (6.25, -1) ++(0:1 and 0.3) arc (0:180:1 and -0.3);

    \draw (5.25,1) -- (5.25,-1);
    \draw (7.25,1) -- (7.25,-1);

    \draw[>=latex, ->] (8,0) -- (9,0);

    \draw[thick, purple] (11,0) circle (0.3cm);
    \draw[thick, purple] (11,0) circle (1.5cm);
\end{tikzpicture}
    \caption{The action of the matrix $B$ from Example \ref{ex:matrix_actions}.}
    \label{fig:reducible_map}
\end{figure}
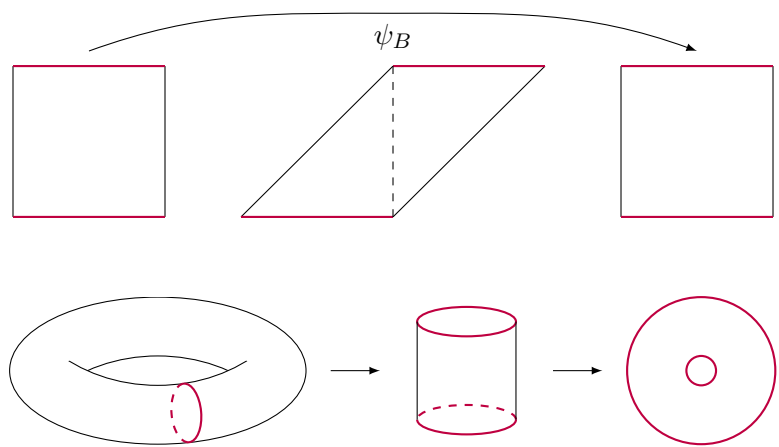

\begin{example}[Finite order and reducible diffeomorphisms on the torus]\label{ex:matrix_actions}

    Consider first the matrix
    \begin{equation*}
         A = \begin{bmatrix}
             -1 & 1 \\
            -1 & 0
        \end{bmatrix},
    \end{equation*}
    which has trace equal to -1 and two conjugate eigenvalues
    \begin{equation*}
        \lambda = \frac{1}{2}(-1 + \ii \sqrt{3}), \quad \overline{\lambda} = \frac{1}{2}(-1 - \ii \sqrt{3}).
    \end{equation*}
    The action of $A$ on $\T^2$ can be visualized as the action on the unit square $[0,1]^2$ followed by a cut \& paste procedure. This is illustrated in Figure \ref{fig:finite_order_map}. A quick computation shows that $A^3 = \operatorname{Id}$, so it follows that for the corresponding diffeomorphism $\psi_A$ we have $\psi_A^3 = \operatorname{Id}$ as well.

    Next we consider the matrix
    \begin{equation*}
        B = 
        \begin{bmatrix}
            1 & 1 \\
            0 & 1
        \end{bmatrix},
    \end{equation*}
    which has trace equal to 2 and a double eigenvalue $\lambda = 1$ with corresponding eigenvector $\mathbf{e}_1 = \begin{bmatrix}
        1\\
        0
    \end{bmatrix}$. Cutting the torus along this eigenvector we see that we may as well study a map on an annulus, which is simpler than the torus, see Figure \ref{fig:reducible_map}.
\end{example}

\begin{figure}[ht]
    \centering
    \begin{tikzpicture}
    \draw[>=latex, ->,thick, noamblue] (0,0) -- (2,0);
    \draw[>=latex, ->, thick, purple] (0,0) -- (0,2);
    \draw (0,2) -- (2,2) -- (2,0);

    \foreach \p in {1/2, 1, 3/2}{
        \draw[dashed] (0,\p) -- ++(0:2);
        \draw[dashed] (\p,0) -- ++(90:2);
    }
    \foreach \p in {0, 1}{
        \fill[red, opacity = 0.15] (\p, 0) rectangle (\p+1/2, 1/2);
        \fill[red, opacity = 0.15] (\p, 1) rectangle (\p+1/2, 3/2);
        \fill[red, opacity = 0.15] (\p+1/2, 1/2) rectangle (\p+1,1);
        \fill[red, opacity = 0.15] (\p+1/2, 3/2) rectangle (\p+1, 2);
    }
    \foreach \p in {1/2, 3/2}{
        \fill[blue, opacity = 0.15] (\p, 0) rectangle (\p+1/2, 1/2);
        \fill[blue, opacity = 0.15] (\p, 1) rectangle (\p+1/2, 3/2);
        \fill[blue, opacity = 0.15] (\p-1/2, 1/2) rectangle (\p,1);
        \fill[blue, opacity = 0.15] (\p-1/2, 3/2) rectangle (\p, 2);
    }

    \def \x{5};

    \draw[>=latex, ->, thick, purple] (\x,2) -- (\x+2,2);
    \draw[->, >=latex, thick, noamblue](\x,2) -- (\x-2,0);
    \draw (\x+2,2) -- (\x,0) -- (\x-2,0);

    \foreach \p in {\x+1/2, \x+1, \x+3/2}{
        \draw[dashed] (\p,2) -- ++(-135:sqrt(8););
        \draw[dashed] (\p, \p-\x) -- ++(180:2);
    }

    \fill[blue, opacity = 0.15] (\x-2,0) -- (\x-3/2,0) -- (\x-1, 1/2) -- (\x-3/2,1/2) -- cycle;
    \fill[blue, opacity = 0.15] (\x-1,0) -- (\x-1/2,0) -- (\x, 1/2) -- (\x-1/2,1/2) -- cycle;
    \fill[blue, opacity = 0.15] (\x-1,1) -- (\x-1/2,1) -- (\x, 3/2) -- (\x-1/2,3/2) -- cycle;
    \fill[blue, opacity = 0.15] (\x,1) -- (\x+1/2,1) -- (\x+1, 3/2) -- (\x+1/2,3/2) -- cycle;
    \fill[blue, opacity = 0.15] (\x-1,1/2) -- (\x-1/2,1/2) -- (\x, 1) -- (\x-1/2,1) -- cycle;
    \fill[blue, opacity = 0.15] (\x,1/2) -- (\x+1/2,1/2) -- (\x+1, 1) -- (\x+1/2,1) -- cycle;    \fill[blue, opacity = 0.15] (\x,3/2) -- (\x+1/2,3/2) -- (\x+1, 2) -- (\x+1/2,2) -- cycle;
    \fill[blue, opacity = 0.15] (\x+1,3/2) -- (\x+3/2,3/2) -- (\x+2, 2) -- (\x+3/2,2) -- cycle;

    \fill[red, opacity = 0.15] (\x-3/2,0) -- (\x-1,0) -- (\x-1/2, 1/2) -- (\x-1,1/2) -- cycle;
    \fill[red, opacity = 0.15] (\x-1/2,0) -- (\x,0) -- (\x+1/2, 1/2) -- (\x,1/2) -- cycle;
    \fill[red, opacity = 0.15] (\x-1/2,1) -- (\x,1) -- (\x+1/2, 3/2) -- (\x,3/2) -- cycle;
    \fill[red, opacity = 0.15] (\x+1/2,1) -- (\x+1,1) -- (\x+3/2, 3/2) -- (\x+1,3/2) -- cycle;
    \fill[red, opacity = 0.15] (\x-3/2,1/2) -- (\x-1,1/2) -- (\x-1/2, 1) -- (\x-1,1) -- cycle;
    \fill[red, opacity = 0.15] (\x-1/2,1/2) -- (\x,1/2) -- (\x+1/2, 1) -- (\x,1) -- cycle;    \fill[red, opacity = 0.15] (\x-1/2,3/2) -- (\x,3/2) -- (\x+1/2, 2) -- (\x,2) -- cycle;
    \fill[red, opacity = 0.15] (\x+1/2,3/2) -- (\x+2/2,3/2) -- (\x+3/2, 2) -- (\x+1,2) -- cycle;

    \def \w{8}
    \draw[] (\w, 0) rectangle (\w+2,2);

    \foreach \p in {1/2, 1, 3/2}{
        \draw[dashed] (\w, \p) -- (\w+2,\p);
    }
    \foreach \p in {0,1/2,1,3/2}{
        \draw[dashed] (\w,\p) -- (\w+2-\p,2);
    }
    \foreach \p in {1/2, 1, 3/2}{
        \draw[dashed] (\w+\p,0) -- (2+\w, 2-\p);
    }

    \foreach \p in {0,1/2,1,3/2}{
        \fill[blue, opacity = 0.15] (\w,\p) -- (\w+1/2,\p) -- (\w+1, \p+1/2) -- (\w+1/2, \p+1/2) -- cycle;
        \fill[blue, opacity = 0.15] (\w+1,\p) -- (\w+3/2,\p) -- (\w+2, \p+1/2) -- (\w+3/2, \p+1/2) -- cycle;
        \fill[purple, opacity = 0.15] (\w+1/2,\p) -- (\w+2/2,\p) -- (\w+3/2, \p+1/2) -- (\w+2/2, \p+1/2) -- cycle;

        \fill[purple, opacity =0.15] (\w,\p) -- (\w+1/2, \p+1/2) -- (\w, \p+1/2) -- cycle;
        \fill[purple, opacity =0.15] (\w+3/2,\p) -- (\w+2, \p) -- (\w+2, \p+1/2) -- cycle;
    }

    \draw[>=latex, ->] (1,2.2) to[out=20, in =180] (\x,2.7) to[out = 0, in = 160] (\w+1, 2.2);
    \node[below] at (\x, 2.7) {$\psi_A$};
\end{tikzpicture}

\begin{tikzpicture} 
    \draw[>=latex, ->,thick, noamblue] (0,0) -- (2,0);
    \draw[>=latex, ->, thick, purple] (0,0) -- (0,2);
    \draw (0,2) -- (2,2) -- (2,0);

    \foreach \p in {1/2, 1, 3/2}{
        \draw[dashed] (0, \p) -- (2,\p);
    }
    \foreach \p in {0,1/2,1,3/2}{
        \draw[dashed] (0,\p) -- (2-\p,2);
    }
    \foreach \p in {1/2, 1, 3/2}{
        \draw[dashed] (\p,0) -- (2, 2-\p);
    }

    \foreach \p in {0,1/2,1,3/2}{
        \fill[blue, opacity = 0.15] (0,\p) -- (1/2,\p) -- (1, \p+1/2) -- (1/2, \p+1/2) -- cycle;
        \fill[blue, opacity = 0.15] (1,\p) -- (3/2,\p) -- (2, \p+1/2) -- (3/2, \p+1/2) -- cycle;
        \fill[purple, opacity = 0.15] (1/2,\p) -- (2/2,\p) -- (3/2, \p+1/2) -- (2/2, \p+1/2) -- cycle;

        \fill[purple, opacity =0.15] (0,\p) -- (1/2, \p+1/2) -- (0, \p+1/2) -- cycle;
        \fill[purple, opacity =0.15] (3/2,\p) -- (2, \p) -- (2, \p+1/2) -- cycle;
    }

    \def \x{5};

    \draw[>=latex, ->, thick, purple] (\x,2) -- (\x+2,2);
    \draw[->, >=latex, thick, noamblue](\x,2) -- (\x-2,0);
    \draw (\x+2,2) -- (\x,0) -- (\x-2,0);

    \foreach \p in {\x+1/2, \x+1, \x+3/2}{
        \draw[dashed] (\p,2) -- ++(-135:sqrt(8););
    }
    \foreach \p in {0,1/2,1,3/2}{
        \draw[dashed] (\x+\p,2) -- (\x+\p, \p);
    }
    \foreach \p in {-1/2, -1, -3/2}{
        \draw[dashed] (\x+\p, 0) -- (\x+\p, 2+\p);
    }

    \foreach \p in {1/2, 1, 3/2, 2}{
        \fill[blue, opacity = 0.15] (\x-2+\p, 1/2) -- (\x-2+\p, 0) -- (\x-3/2+\p, 1/2) -- (\x-3/2+\p, 1) -- cycle; 
        \fill[blue, opacity = 0.15] (\x-1+\p, 3/2) -- (\x-1+\p, 1) -- (\x-1/2+\p, 3/2) -- (\x-1/2+\p, 2) -- cycle; 

        \fill[red, opacity = 0.15] (\x-3/2+\p, 2/2) -- (\x-3/2+\p, 1/2) -- (\x-1+\p, 2/2) -- (\x-1+\p, 3/2) -- cycle;
        \fill[red, opacity = 0.15] (\x-5/2+\p, 0) -- (\x-2+\p,0) -- (\x-2+\p,1/2) -- cycle;
        \fill[red, opacity = 0.15] (\x-1/2+\p, 2) -- (\x-1/2+\p,3/2) -- (\x+\p,2) -- cycle;
    }

    \def \w{8}
    \draw[] (\w, 0) rectangle (\w+2,2);

    \foreach \p in {1/2, 1, 3/2}{
        \draw[dashed] (\w+\p, 0) -- (\w+\p,2);
    }
    \foreach \p in {0,1/2,1,3/2}{
        \draw[dashed] (\w,\p) -- (\w+2-\p,2);
    }
    \foreach \p in {1/2, 1, 3/2}{
        \draw[dashed] (\w+\p,0) -- (2+\w, 2-\p);
    }

    \foreach \p in {0,1/2,1,3/2}{
        \fill[blue, opacity = 0.15] (\w+\p,0) -- (\w+\p, 1/2) -- (\w+\p+1/2, 1) -- (\w+\p+1/2,1/2) -- cycle;
        \fill[blue, opacity = 0.15] (\w+\p,1) -- (\w+\p, 3/2) -- (\w+\p+1/2, 2) -- (\w+\p+1/2,3/2) -- cycle;

        \fill[red, opacity = 0.15] (\w+\p,1/2) -- (\w+\p, 2/2) -- (\w+\p+1/2, 3/2) -- (\w+\p+1/2,2/2) -- cycle;
        \fill[red, opacity = 0.15] (\w+\p,0) -- (\w+\p+1/2, 0)-- (\w+\p+1/2, 1/2) -- cycle;
        \fill[red, opacity = 0.15] (\w+\p, 2) -- (\w + \p, 3/2) -- (\w + \p + 1/2, 2) -- cycle;
    }

    \draw[>=latex, ->] (1,2.2) to[out=20, in =180] (\x,2.7) to[out = 0, in = 160] (\w+1, 2.2);
    \node[below] at (\x, 2.7) {$\psi_A$};
\end{tikzpicture}

\begin{tikzpicture} 
    \draw[>=latex, ->,thick, noamblue] (0,0) -- (2,0);
    \draw[>=latex, ->, thick, purple] (0,0) -- (0,2);
    \draw (0,2) -- (2,2) -- (2,0);

    \foreach \p in {1/2, 1, 3/2}{
        \draw[dashed] (\p, 0) -- (\p,2);
    }
    \foreach \p in {0,1/2,1,3/2}{
        \draw[dashed] (0,\p) -- (2-\p,2);
    }
    \foreach \p in {1/2, 1, 3/2}{
        \draw[dashed] (\p,0) -- (2, 2-\p);
    }

    \foreach \p in {0,1/2,1,3/2}{
        \fill[blue, opacity = 0.15] (\p,0) -- (\p, 1/2) -- (\p+1/2, 1) -- (\p+1/2,1/2) -- cycle;
        \fill[blue, opacity = 0.15] (\p,1) -- (\p, 3/2) -- (\p+1/2, 2) -- (\p+1/2,3/2) -- cycle;

        \fill[red, opacity = 0.15] (\p,1/2) -- (\p, 2/2) -- (\p+1/2, 3/2) -- (\p+1/2,2/2) -- cycle;
        \fill[red, opacity = 0.15] (\p,0) -- (\p+1/2, 0)-- (\p+1/2, 1/2) -- cycle;
        \fill[red, opacity = 0.15] (\p, 2) -- (\p, 3/2) -- (\p + 1/2, 2) -- cycle;
    }

    \def \x{5};

    \draw[>=latex, ->, thick, purple] (\x,2) -- (\x+2,2);
    \draw[->, >=latex, thick, noamblue](\x,2) -- (\x-2,0);
    \draw (\x+2,2) -- (\x,0) -- (\x-2,0);

    \foreach \p in {0,1/2,1,3/2}{
        \draw[dashed] (\x+\p,2) -- (\x+\p, \p);
    }
    \foreach \p in {-1/2, -1, -3/2}{
        \draw[dashed] (\x+\p, 0) -- (\x+\p, 2+\p);
    }
    \foreach \p in {\x+1/2, \x+1, \x+3/2}{
        \draw[dashed] (\p, \p-\x) -- ++(180:2);
    }

    \foreach \p in {0, 1/2, 1, 3/2}{
        \fill[blue, opacity = 0.15] (\x-3/2+\p, \p) rectangle (\x-1+\p, \p+1/2);
        \fill[blue, opacity = 0.15] (\x-1/2+\p, \p) rectangle (\x+\p, \p+1/2);

        \fill[red, opacity = 0.15] (\x-1+\p, \p) rectangle (\x-1/2+\p, \p+1/2);
        \fill[red, opacity = 0.15] (\x-2+\p,\p) -- (\x-3/2+\p,\p) -- (\x-3/2+\p,\p+1/2) -- cycle;
        \fill[red, opacity = 0.15] (\x+\p, \p) -- (\x+\p, \p+1/2) -- (\x+\p+1/2, \p+1/2) -- cycle; 
    }
    \def \w{8}
    \draw[] (\w, 0) rectangle (\w+2,2);

    \foreach \p in {1/2, 1, 3/2}{
        \draw[dashed] (\w,\p) -- ++(0:2);
        \draw[dashed] (\w+\p,0) -- ++(90:2);
    }
    \foreach \p in {0, 1}{
        \fill[red, opacity = 0.15] (\w+\p, 0) rectangle (\w+\p+1/2, 1/2);
        \fill[red, opacity = 0.15] (\w+\p, 1) rectangle (\w+\p+1/2, 3/2);
        \fill[red, opacity = 0.15] (\w+\p+1/2, 1/2) rectangle (\w+\p+1,1);
        \fill[red, opacity = 0.15] (\w+\p+1/2, 3/2) rectangle (\w+\p+1, 2);
    }
    \foreach \p in {1/2, 3/2}{
        \fill[blue, opacity = 0.15] (\w+\p, 0) rectangle (\w+\p+1/2, 1/2);
        \fill[blue, opacity = 0.15] (\w+\p, 1) rectangle (\w+\p+1/2, 3/2);
        \fill[blue, opacity = 0.15] (\w+\p-1/2, 1/2) rectangle (\w+\p,1);
        \fill[blue, opacity = 0.15] (\w+\p-1/2, 3/2) rectangle (\w+\p, 2);
    }

    \draw[>=latex, ->] (1,2.2) to[out=20, in =180] (\x,2.7) to[out = 0, in = 160] (\w+1, 2.2);
    \node[below] at (\x, 2.7) {$\psi_A$};
\end{tikzpicture}
    \caption{The action of the matrix $A$ from Example \ref{ex:matrix_actions}.}
    \label{fig:finite_order_map}
\end{figure}
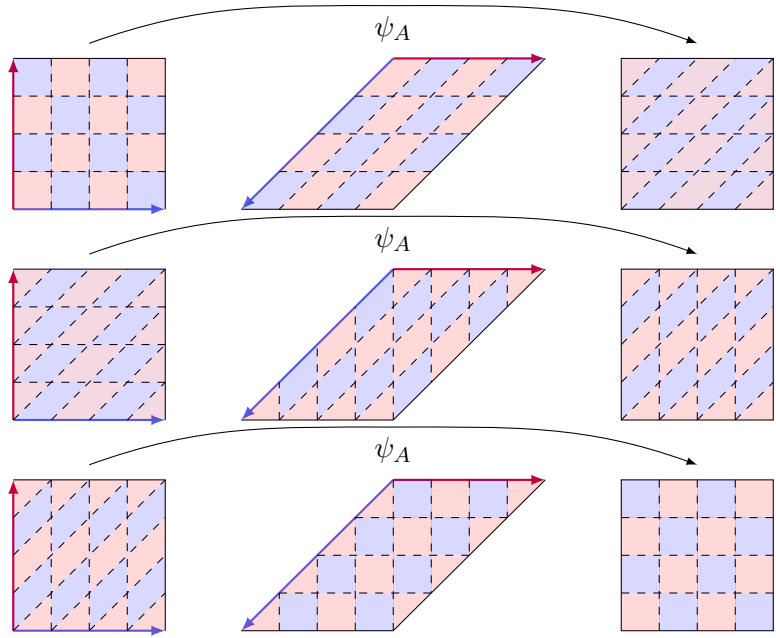

Now, let us focus on the third case which produces the richest dynamics on $\T$, since in this case the corresponding diffeomorphism $\psi$ is chaotic in the sense of Devaney: The set of periodic points is dense on $\T^2$, there is a dense orbit and the map exhibits sensitive dependence on initial conditions. The most well-known example of these types of maps is the so-called CAT map studied by Arnold in \cite{arnold1967problemes}. It is the diffeomorphism $\psi \colon \T^2 \to \T^2$ induced by the matrix
\begin{equation*}
    A = 
    \begin{bmatrix}
        2 & 1 \\
        1 & 1
    \end{bmatrix}.
\end{equation*}

\begin{remark}
    It is sometimes claimed that \enquote{CAT} is an acronym for \enquote{continuous automorphism of the torus}. Another possible explanation for the etymology of the name of this map comes from the fact that Arnold initially used the picture of a cat to illustrate the action of $A$ on $\T^2$. 
\end{remark}

\begin{figure}[ht]
    \centering
    \begin{tikzpicture}
    \draw[thick, >=latex, ->, noamblue] (0,0) -- (2,0);
    \draw[thick, >=latex, ->, purple] (0,0) -- (0,2);
    \draw[] (0,2) -- (2,2) -- (2,0);
    \draw[dashed] (1,0) -- (0,2) -- (2,0) -- (1,2);
    \fill[teal, opacity = 0.15] (0,0) -- (1,0) -- (0,2) -- cycle;
    \fill[violet, opacity = 0.15] (0,2) -- (1,0) -- (2,0) -- cycle;
    \fill[blue, opacity = 0.15] (0,2) -- (2,0) -- (1,2) -- cycle;
    \fill[red, opacity = 0.15] (2,2) -- (2,0) -- (1,2) -- cycle;

    \def\x{3}
    \draw[thick,>=latex, ->,noamblue] (\x,0) -- (\x+4,2);
    \draw[thick,>=latex, ->, purple] (\x,0) -- (\x+2, 2);
    \draw[] (\x+2,2) -- (\x+6,4) -- (\x+4,2);
    \draw[dashed] (\x+2, 1) -- (\x+2,2) -- (\x+4,2) -- (\x+4,3);
    \fill[teal, opacity = 0.15] (\x,0) -- (\x+2,2) -- (\x+2,1) -- cycle;
    \fill[violet, opacity =0.15] (\x+2,1) -- (\x+2,2) -- (\x+4,2) -- cycle;
    \fill[blue, opacity =0.15] (\x+2,2) -- (\x+4,2) -- (\x+4,3) -- cycle;
    \fill[red, opacity = 0.15] (\x + 4, 2) -- (\x+4,3) -- (\x+6,4) --cycle; 

    \def\w{10}

    \draw[] (\w,0) rectangle (\w+2,2);
    \draw[dashed] (\w+2,1) -- (\w,0) -- (\w+2,2) -- (\w,1);
    \fill[teal, opacity = 0.15] (\w,0) -- (\w+2,1) -- (\w+2,2) -- cycle;
    \fill[violet, opacity = 0.15] (\w,1) -- (\w,2) -- (\w+2,2) -- cycle;
    \fill[blue, opacity =0.15] (\w,0) -- (\w+2,0) -- (\w+2,1) -- cycle;
    \fill[red, opacity= 0.15] (\w,0) -- (\w+2,2) -- (\w,1);

    \draw[>=latex, ->] (1,-0.2) to[out=-20, in =180] (\x+3,-1.1) to[out = 0, in = -160] (\w+1, -0.2);
    \node[above] at (\x+3, -1.1) {$\psi_A$};
\end{tikzpicture}
    \caption{The CAT map acting on $\T^2$.}
    \label{fig:cat_map}
\end{figure}
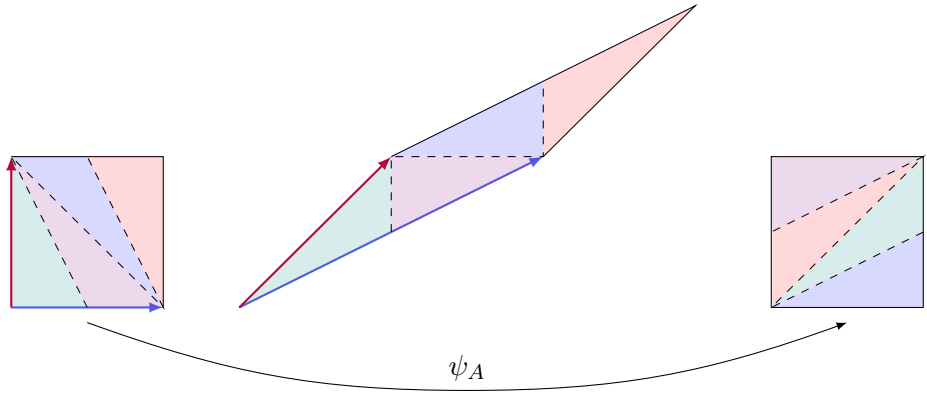

The eigenvalues of $A$ are given by
\begin{equation*}
    \lambda_1 = \frac{1}{2}(3+\sqrt{5}), \quad \lambda_2 = \frac{1}{2}(3 - \sqrt{5}).
\end{equation*}
Note in particular that $\lambda_2 = \lambda_1^{-1}$. The corresponding eigenvectors are
\begin{equation*}
    \mathrm{v}_1 = \begin{bmatrix}
        \frac{1}{2}(1+\sqrt{5}) \\
        1
    \end{bmatrix}, \quad
    \mathrm{v}_2 = \begin{bmatrix}
        \frac{1}{2}(1 - \sqrt{5}) \\
        1
    \end{bmatrix}.
\end{equation*}
The action of the corresponding diffeomorphism $\psi_A \colon \T^2 \to \T^2$ is depicted in figure \ref{fig:cat_map}. Since $\lambda_1 > 1$ and $\lambda_2 < 1$, the corresponding eigendirections get expanded or contracted respectively. As we will see below, this is exactly the defining feature of what we will call a pseudo-Anosov map.

\subsubsection{(Singular) Measured Foliations}

As the next step towards defining pseudo-Anosov maps, we will introduce the concept of singular measured foliations. Once more we will start with the simple case of the torus, where we can even define a \emph{nonsingular} measured foliation.

\begin{figure}[ht]
    \centering
    \begin{tikzpicture}
    \begin{scope}
    \coordinate (xleft) at (-1,0); \coordinate (xright) at (5,0);
    \coordinate (ytop) at (0,5); \coordinate (ybot) at (0,-1);
    
    \draw[>=latex, ->, thick] (xleft) -- (xright);
    \draw[>=latex, ->, thick] (ybot) -- (ytop);

    \draw[dashed] (2,0) -- (2,5);
    \draw[dashed] (4,0) -- (4,5);
    \draw[dashed] (0,2) -- (5,2);
    \draw[dashed] (0,4) -- (5,4);

    \clip (-1,-1) rectangle (5,5);

    \draw[thick, purple] (0,0) -- ++(e*pi/7 r:6) (0,0) -- ++(180+e*pi/7 r:4);
    
    \foreach \p in {5/13, 10/13, 15/13, 20/13, 25/13}{
        \draw (-\p,\p) -- ++(e*pi/7 r:6) (-\p,\p) -- ++(180+e*pi/7 r:4);
    }
    \foreach \p in {10/13, 20/13, 30/13, 40/13, 50/13, 60/13, 70/13, 80/13, 90/13, 100/13}{
        \draw (\p,\p) -- ++(e*pi/7 r:6) (\p,\p) -- ++(180+e*pi/7 r:8);
    }
    \end{scope}
    \node at (1.85,5.25) {\textcolor{purple}{$\ell$}};
    \node at (5.1,2.4) {$\Tilde{\mathcal{F}}_\ell$};
    \begin{scope}
        \def\x{7}; \def\y{1};
        \draw[thick, name path=rectangle] (\x,\y) rectangle (\x+2,\y+2);

        \clip (\x,\y) rectangle (\x+2,\y+2);
        \draw[purple, thick, name path = line] (\x,\y) -- ++(e*pi/7 r:3);
        \path[name intersections={of=line and rectangle, by=int}];

        \draw[thick, purple, name path = line2] ($(int) - (0,2)$) -- ++(e*pi/7 r:5);
        \path[name intersections={of=line2 and rectangle, by = int2}];

        \draw[purple, thick, name path=line3] ($(int2) - (0,2)$) -- ++(e*pi/7 r:5);
        \path[name intersections={of=line3 and rectangle, by = int3}];

        \draw[purple, thick, name path = line4] ($(int3) - (2,0)$) -- ++(e*pi/7 r:5);
        \path[name intersections={of=line4 and rectangle, by = int4}];

        \draw[purple, thick, name path = line5] ($(int4) - (0,2)$) -- ++(e*pi/7 r:2);
        \draw[purple, thick, dashed] ($(int4) - (0,2) + (e*pi/7 r:2)$) -- ++(e*pi/7 r:1);

        \foreach \p in {0,0.5,0.8,1.7,2,2.6,2.9,3.1,3.2,3.8}{
            \draw (\x+\p,\y+2+\p) -- ++(e*pi/7 r:6) (\x+\p,\y+2+\p) -- ++(180+e*pi/7 r:7);
        }
    \end{scope}
    \node at (7+2.3, 2.4) {$\mathcal{F}_\ell$};
\end{tikzpicture}
    \caption{The foliation $\Tilde{\mathcal{F}}_\ell$ of $\R^2$ and the induced foliation $\mathcal{F}_\ell$ of $\T^2$.}
    \label{fig:plane_torus_foliation}
\end{figure}
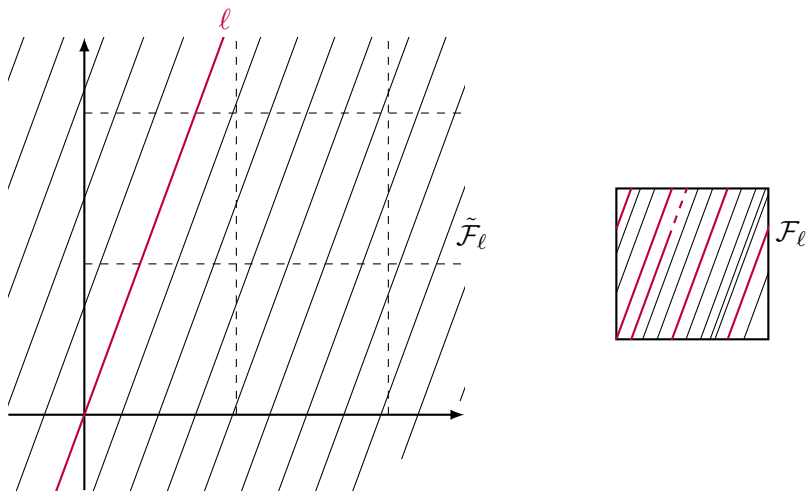

Let $\ell$ be any line through the origin $(0,0)$ in $\R^2$, which determines a {foliation} $\tilde{\mathcal{F}}_\ell$ of $\R^2$ given by all lines parallel to $\ell$. We will call these lines the \emph{leaves} of the foliation. Note that any translation preserves $\tilde{\mathcal{F}}_\ell$, meaning that it sends leaves to leaves. Recall further that the group of deck transformations of a universal covering is isomorphic to the fundamental group of the covered space, so that it follows that deck transformations for the standard universal covering $p \colon \R^2 \to \T^2$ of the torus are translations. Thus, $\Tilde{\mathcal{F}}_\ell$ descends to a foliation $\mathcal{F}_\ell$ on $\T^2$, as is illustrated in Figure \ref{fig:plane_torus_foliation}.

The following Theorem represents a well-known fact about such foliations, or equivalently, about the linear flow on the torus.

\begin{theorem}[Dichotomy of Foliations]\label{thm:dichotomy_foliations}
    Let $\theta$ be the slope of the foliation $\mathcal{F}_\ell$ of $\T^2$.
    \begin{enumerate}
        \item If $\theta$ is rational, then every leaf of $\mathcal{F}_\ell$ is a simple closed geodesic in $\T^2$.
        \item If $\theta$ is irrational, then every leaf of $\mathcal{F}_\ell$ is a dense geodesic in $\T^2$.
    \end{enumerate}
\end{theorem}
To prove Theorem \ref{thm:dichotomy_foliations}, we may show the analogous statement in terms of the linear flow on the torus. One can show that the first return map to the horizontal segment of $[0,1]$ of $\T^2 = [0,1] \times [0,1] / {\sim}$ is a rotation on $S^1$ by $\theta$. The nature of the rotation depends on $\theta$ exactly as stated in the theorem and this behaviour then translates to the linear flow on the torus. 

Given a foliation $\Tilde{\mathcal{F}}_\ell$ (on $\R^2$), we can establish the following definition.

\begin{definition}[Transverse measure]\label{def:transverse_measure_1}
    Let $\nu_\ell\colon \R^2 \to \R$ be the function that records the (signed) distance from $\ell$. For every smooth arc $\alpha$ transverse to the leaves of $\Tilde     {\mathcal{F}}_\ell$, we define a \emph{transverse measure $\mu$ on $\Tilde{\mathcal{F}}_\ell$} by
    \begin{equation*}
        \mu(\alpha) = \int_\alpha \dd \nu_\ell,
    \end{equation*}
    where $\dd \nu_\ell$ is the 1-form given by $\dd \nu_\ell = \frac{\partial \nu_\ell}{\partial x} \dd x + \frac{\partial \nu_\ell}{\partial y} \dd y$.
\end{definition}

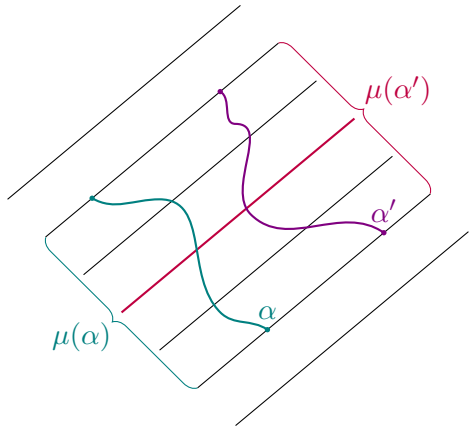
\begin{figure}[ht]
    \centering
    \begin{tikzpicture}
    \foreach \p in {0,0.5,1,2,2.5,3}{
        \draw[] (\p,-\p) -- ++(40:4);
    }
    \draw[thick, purple] (1.5, -1.5) --++(40:4);

    \path (0.5,-0.5) -- ++(40:0.8) coordinate (a1); 
    \path (2.5,-2.5) -- ++(40:1.2) coordinate (a2);

    \draw[teal, thick] (a1) to[out=-30, in=180] (2,0) to[out=0, in=135] (2.8,-1.5) to[out=-45, in = 145] (a2) node[above] {$\alpha$};
    \fill[teal] (a1) circle (1pt);
    \fill[teal] (a2) circle (1pt);

    \path (0.5,-0.5) -- ++(40:3) coordinate (b1); 
    \path (2.5,-2.5) -- ++(40:3.2) coordinate (b2);

    \draw[violet, thick] (b1) to[out=-30, in=180] (3,1) to[out=0, in=135] (3.2,-0.2) to[out=-45, in = 145] (b2) node[above]{$\alpha'$};
    \fill[violet] (b1) circle (1pt);
    \fill[violet] (b2) circle (1pt);

    \path (0.5,-0.5) -- ++(40:4) coordinate (e1);
    \path (2.5, -2.5) --++(40:4) coordinate (e2);

    \draw [purple, decorate, decoration={brace, amplitude=5pt}] (e1) -- node[above right, yshift=1] {$\mu(\alpha')$}  (e2);
    \draw[teal, decorate, decoration={brace, amplitude=5pt}] (2.5, -2.5) -- node[below left] {$\mu(\alpha)$} (0.5, -0.5);
\end{tikzpicture}
    \caption{The transverse measure $\mu$ from Definition \ref{def:transverse_measure_1} is invariant under leaf-preserving isotopies.}
    \label{fig:transverse_measure_isotopy_invariant}
\end{figure}

The transverse measure $\mu$ defined above measures the total variation of $\alpha$ in the direction perpendicular to $\ell$. It is clear that $\mu(\alpha)$ is invariant under isotopies of $\alpha$ such that both ends of $\alpha$ belong to the same leaf of $\Tilde{\mathcal{F}}_\ell$ before and after the isotopy. This is illustrated in Figure \ref{fig:transverse_measure_isotopy_invariant}.

Note that the 1-form $\dd \nu_\ell$ is preserved by translations, so it descends to a 1-form $\omega_\ell$ on $\T^2$ and it induces a transverse measure on the foliation $\mathcal{F}_\ell$.

\begin{definition}[Measured foliation]
    The structure of a foliation $\mathcal{F}_\ell$ on $\T^2$ together with a transverse measure $\mu$ is called a \emph{measured foliation} on $\T^2$.
\end{definition}

We want to extend this idea to closed surfaces $X$ of genus $\mathbf{g} \geq 2$. As we will see below, the Euler--Poincaré formula (Proposition \ref{prop:euler_poincare}) implies that such surfaces cannot support a foliation of the form we have just described. Analogously to the definition of translation surfaces as generalizations of a flat torus, we can adjust the definition by allowing a finite number of singularities of a specific type. 

\begin{figure}[ht]
    \centering
    \begin{tikzpicture}
    \fill[] (0,0) circle (1pt);
    \draw[thick] (0,0) -- ++(90:2) coordinate (e1);\draw[thick] (0,0) -- ++(210:2) coordinate (e2); \draw[thick] (0,0) -- ++ (330:2) coordinate (e3);

    \node[noamblue] at (0,1) {$\vee$};
    \path (0,0) -- ++ (210:1) coordinate (v2);
    \node[rotate = -60, noamblue] at (v2) {$\vee$};
    \path (0,0) -- ++ (330:1) coordinate (v3);
    \node[rotate = 60, noamblue] at (v3) {$\vee$};

    \path (0,0) -- ++ (30:0.4) coordinate (a1);
    \path (0,0) -- ++(150:0.4) coordinate (a2);
    \path (0,0) -- ++(-90:0.4) coordinate (a3);

    \path (0,0) -- ++ (30:0.7) coordinate (b1);
    \path (0,0) -- ++(150:0.7) coordinate (b2);
    \path (0,0) -- ++(-90:0.7) coordinate (b3);

    \draw [] plot [smooth, tension=0.4] coordinates {($(e1)+ (0.2,0)$) (a1) ($(e3)+(0,0.2)$)};
    \node[rotate=34, noamblue] at (a1) {$\vee$};
    \draw [] plot [smooth, tension=0.4] coordinates {($(e1)- (0.2,0)$) (a2) ($(e2)+(0,0.2)$)};
    \node[rotate=-34, noamblue] at (a2) {$\vee$};
    \draw [purple] plot [smooth, tension=0.4] coordinates {($(e2) - (0,0.2)$) (a3) ($(e3)-(0,0.2)$)};

    \draw [] plot [smooth, tension=0.6] coordinates {($(e1)+ (0.4,0)$) (b1) ($(e3)+(0,0.4)$)};
    \draw [] plot [smooth, tension=0.6] coordinates {($(e1)- (0.4,0)$) (b2) ($(e2)+(0,0.4)$)};
    \draw [] plot [smooth, tension=0.6] coordinates {($(e2)- (0,0.4)$) (b3) ($(e3)-(0,0.4)$)};
\end{tikzpicture}
\qquad
\begin{tikzpicture}
    \fill[] (0,0) circle (1pt);
    \draw[thick] (-1.5,0) -- (1.5,0); \draw[thick] (0,-1.5) -- (0,1.5);
    \node[noamblue] at (0,0.75) {$\vee$};
     \node[noamblue] at (0,-0.75) {$\wedge$};

     \node[noamblue, rotate = 90] at (0.75,0) {$\vee$};
     \node[noamblue, rotate = 90] at (-0.75,0) {$\wedge$};

    \begin{scope}
    \draw[] plot [smooth, tension = 0.6] coordinates {(0.2,1.5) (0.4,0.4) (1.5,0.2)};
    \node[noamblue, rotate =52] at (0.4,0.4) {$\vee$};
    \draw[] plot [smooth, tension = 0.6] coordinates {(0.4,1.5) (0.6,0.6) (1.5,0.4)};
    \node[noamblue, rotate =53] at (0.6,0.6) {$\vee$};
    \end{scope}

    \begin{scope}[xscale =-1]
    \draw[] plot [smooth, tension = 0.6] coordinates {(0.2,1.5) (0.4,0.4) (1.5,0.2)};
    \node[noamblue, rotate =-52] at (0.4,0.4) {$\vee$};
    \draw[] plot [smooth, tension = 0.6] coordinates {(0.4,1.5) (0.6,0.6) (1.5,0.4)};
    \node[noamblue, rotate =-53] at (0.6,0.6) {$\vee$};
    \end{scope}

    \begin{scope}[yscale =-1]
    \draw[] plot [smooth, tension = 0.6] coordinates {(0.2,1.5) (0.4,0.4) (1.5,0.2)};
    \node[noamblue, rotate =-52] at (0.4,0.4) {$\wedge$};
    \draw[] plot [smooth, tension = 0.6] coordinates {(0.4,1.5) (0.6,0.6) (1.5,0.4)};
    \node[noamblue, rotate =-53] at (0.6,0.6) {$\wedge$};
    \end{scope}

    \begin{scope}[xscale =-1, yscale = -1]
    \draw[] plot [smooth, tension = 0.6] coordinates {(0.2,1.5) (0.4,0.4) (1.5,0.2)};
    \node[noamblue, rotate =52] at (0.4,0.4) {$\wedge$};
    \draw[] plot [smooth, tension = 0.6] coordinates {(0.4,1.5) (0.6,0.6) (1.5,0.4)};
    \node[noamblue, rotate =53] at (0.6,0.6) {$\wedge$};
    \end{scope}
    
\end{tikzpicture}
    \caption{On the left we have a 3-pronged saddle, which cannot be consistently oriented. In the example here, the red line cannot be given an orientation consistent with the orientation of the prongs. On the right we have a 4-pronged saddle, a possible orientation is indicated by the blue arrows.}
    \label{fig:even_odd_prongs}
\end{figure}
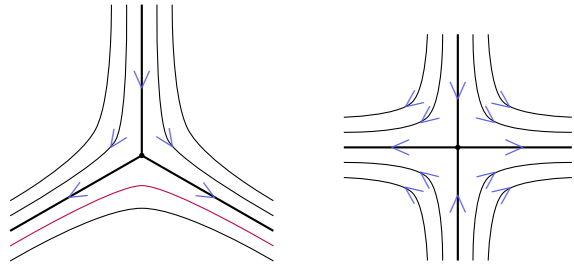

\begin{definition}[Singular foliation]\label{def:singular_foliation}
    A closed surface $X$ is said to have a \emph{singular foliation} $\mathcal{F}$ if it is a disjoint union of subsets, known as the leaves of $\mathcal{F}$, along with a finite set of points on $X$ referred to as the singular points of $\mathcal{F}$, such that the following conditions are satisfied.
    \begin{enumerate}
        \item For a regular point $p \in X$, there is a smooth chart from a neighborhood $U$ of $p$ to $\R^2$ that takes leaves to horizontal line segments. It follows that the transition maps between any two charts around regular points are smooth maps of the form
        \begin{equation*}
            (x,y) \mapsto \big(f(x,y), g(y)\big),
        \end{equation*}
        i.e., the transition maps take horizontal lines to horizontal lines. 
        \item For a singular point $p \in X$, there is a smooth chart from a neighborhood $U$ of $p$ to $\R^2$ that takes the leaves to the level sets of a $k$-pronged saddle, where $k \geq 3$ as in Figure \ref{fig:even_odd_prongs}.
    \end{enumerate}
    Moreover, a singular foliation $\mathcal{F}$ is \emph{orientable}, if the leaves can be consistently oriented. We say $\mathcal{F}$ is \emph{locally orientable}, if every point $p \in X$ has a neighborhood that is orientable. 
\end{definition}

\begin{proposition}[Local orientability of foliations]\label{prop:local_orientation_foliation}
    A singular foliation $\mathcal{F}$ is locally orientable if and only if all the singular points have an even number of prongs.
\end{proposition}
\begin{proof}

    Neighborhoods of regular points are foliated by lines and therefore always orientable. Thus, it suffices to look at singular points.

    If $p \in \Sigma$ corresponds to a $2k$-pronged saddle for some $k \in \N_{>1}$, we can pick any starting prong giving it any orientation. Then we move, say, counterclockwise around the singular points and orient the next prong opposite to the first one. Since the singularity is even-pronged, we are able to continue this alternating pattern until we return to the first prong. Orienting the remaining closeby leaves in the obvious way gives a local orientation.

    Conversely, if $p \in \Sigma$ is $(2k+1)$-pronged for $k \in \N_{\geq 1}$, any choice of orientation of the prongs leads to two adjacent prongs with the same orientation. Leaves between these two prongs cannot be oriented consistently. The situation for both cases is depicted in Figure \ref{fig:even_odd_prongs}.
\end{proof}

\begin{remark}
    Note that there do exist foliations that are locally orientable but not globally orientable. It is however true that as soon as we endow the foliation with a transverse measure as defined below, local and global orientability coincide. We will make further comments concerning this fact below. Even in the case where we do not endow a foliation with a transverse measure, there is a classical argument that may be used to still obtain a global orientation on a related surface. The idea is to move to a so-called \emph{orientation cover}, a double covering of the surface that is constructed in a way such that the local orientations of the foliation lift to a global orientation in the cover. If the foliation admits odd-pronged singularities, we can obtain the same result by considering a branched cover. This argument is used below in the proof of the Euler--Poincaré formula (Proposition \ref{prop:euler_poincare}). For details, we refer the reader to \cite{candel2000foliations}, in particular Proposition 3.5.1. The move to an orientation cover is analogous to the approach in the definition of hyperellipticity (see in particular Definition \ref{def:orientation_cover}). 
\end{remark}

\begin{proposition}[Euler--Poincaré formula]\label{prop:euler_poincare}
    Let $X$ be a surface endowed with a singular foliation $\mathcal{F}$ and let $P_s$ denote the number of prongs at a singular point $s \in \Sigma$. Then,
    \begin{equation*}
        2 \chi(X) = \sum_{s \in \Sigma} (2-P_s).
    \end{equation*}
\end{proposition}

\begin{proof}
    The following argument is taken from \cite{fathi1979travaux}. The Euler--Poincaré formula follows from the Poincaré--Hopf Theorem, which says that for any connected compact and orientable Riemannian surface $X$ with vector field $F$, it holds that
    \begin{equation*}
        \sum_{p \in \operatorname{Sing}} I_p = \chi(X),
    \end{equation*}
    where we sum over the singularities of the vector field $F$ and $I_p$ denotes the index of the singular point $p$.

    To be able to view the foliation as a vector field, we need to be able to orient it globally. This can be ensured by considering a so-called orientation cover if necessary, i.e., a two-to-one covering of $X$ having only even-pronged singularities. Denote by $\Sigma'$ the set of odd-pronged singularities and by $\Sigma''$ the set of even-pronged singularities. There is an orientation homomorphism of the tangent bundle of $\mathcal{F}$,
    \begin{equation*}
        \pi_1(X \setminus \Sigma) \to \Z / 2\Z.
    \end{equation*}
    This map defines a double covering which extends over $\Sigma''$ and is branched over $\Sigma'$. Morally, what we do is we consider both possible orientations except at points in $\Sigma'$, where we double the number of prongs in order to be able to even define a local orientation as well.

    So we have a branched covering $p \colon \Tilde{X} \to X$, which is equipped with a singular \emph{orientable} foliation $\Tilde{\mathcal{F}}$, which now we may think of as coming from a vector field $F$. If $s$ is a singularity of $F$, then $P_s$ is even (by the above) and the index of $F$ at $s$ is given by 
    \begin{equation*}
        -\frac{P_s}{2}+1, 
    \end{equation*}
    which is just the index of a saddle point with $P_s$ separatrices. By the Poincaré--Hopf Theorem mentioned above, it follows that
    \begin{equation*}
        \chi(\Tilde{X}) = \sum_{s \in \operatorname{Sing}(F)} I_s = \sum_{s \in \operatorname{Sing(F)}} \left(-\frac{P_s}{2} + 1\right).
    \end{equation*}
    To connect this to the original surface $X$ we make use of the Riemann--Hurwitz formula, which in the context here states that
    \begin{equation*}
        \chi(\Tilde{X}) = 2\chi(X) - |\Sigma'|.
    \end{equation*}
    Moreover, we also have
    \begin{equation*}
        \sum_{\operatorname{Sing}(F)} 1= 2|\Sigma''| + |\Sigma'|. 
    \end{equation*}
    If $p(s) \in \Sigma''$, then $P_s = P_{p(s)}$ and $s$ has a \enquote{twin}, meaning that on the cover there is another point $s'$ such that $p(s) = p(s')$. If $p(s) \in \Sigma'$, then $P_s = 2P_{p(s)}$. Putting all these equalities together, we obtain
    \begin{equation*}
        2\chi(X) = \sum_{s \in \operatorname{Sing}(F)}\left(-\frac{P_s}{2} + 1\right) +|\Sigma'|= \sum_{s \in \Sigma} -P_s + 2 (|\Sigma'| + |\Sigma''|) = \sum_{s \in \Sigma} (2-P_s),
    \end{equation*}
    which establishes the Euler--Poincaré formula.
    \end{proof}

Since we require $P_s \geq 3$, Proposition \ref{prop:euler_poincare} immediately implies the following.

\begin{corollary}
    If a surface $X$ has $\chi(X) > 0$ it cannot admit a foliation, neither singular nor nonsingular. Any foliation on a surface $X$ with $\chi(X) = 0$ may not have any singular points whereas if the surface has $\chi(X) < 0$ the foliation must admit at least one singular point. 
\end{corollary}

This corollary justifies that we will henceforth not specify whether the foliation under consideration is singular, but will unambigously use the term \emph{foliation} for foliations that have singular points as well as for those who have no singular points. 

As before, in the case of the torus, we want to define a transverse measure that assigns a positive length to arcs transverse to the foliation. To define such a measure formally, the following definition is needed.

\begin{definition}[Transverse arcs, leaf-preserving isotopy]\label{def:leaf_preserving_isotopy}
    Let $(X, \mathcal{F})$ be a foliated surface. A smooth arc $\alpha \colon [0, 1] \to X$ is \emph{transverse to $\mathcal{F}$}, if $\alpha([0,1]) \cap \Sigma = \emptyset$ and it is transverse to each leaf of $\mathcal{F}$ at each point in its interior, i.e., $\alpha'(s)$ and some nonzero vector tangent to $\mathcal{F}$ span $T_{\alpha(s)}X$.

    If $\alpha, \beta \colon [0,1] \to X$ are two smooth arcs transverse to $\mathcal{F}$, a \emph{leaf-preserving isotopy} from $\alpha$ to $\beta$ is a map $H \colon [0,1]^2 \to X$ such that
    \begin{enumerate}
        \item $H(s, 0) = \alpha(s)$ and $H(s, 1) = \beta(s)$ for all $s \in [0,1]$,
        \item $H(s,t)$ is transverse to $\mathcal{F}$ for each fixed $t \in [0,1]$,
        \item $H(0, t)$ and $H(1, t)$ is contained in a single leaf (each) for all $t \in [0,1]$.
    \end{enumerate}
\end{definition}

\begin{definition}[Singular measured foliation]
    A \emph{transverse measure} $\mu$ on a foliation $\mathcal{F}$ is a function that assigns a positive real number to each smooth arc transverse to $\mathcal{F}$ so that $\mu$ is invariant under leaf-preserving isotopy and $\mu$ is absolutely continuous with respect to the Lebesgue measure.

    A \emph{(singular) measured foliation $(\mathcal{F}, \mu)$} on a surface $X$ is a foliation $\mathcal{F}$ endowed with a transverse measure $\mu$.
\end{definition}

Insisting on $\mu$ to be absolutely continuous with respect to the (one-dimensional) Lebesgue measure has the following consequence.

\begin{proposition}\label{prop:transverse_measure_is_pullback}
    Let $X$ be a surface endowed with a measured foliation $(\mathcal{F}, \mu)$. Each regular point $p \in X$ has a neighborhood $U$ and a smooth chart $\varphi \colon U \to \R^2$ so that the measure $\mu$ is given by pulling back the measure $|\dd y|$ on $\R^2$.
\end{proposition}
\begin{proof}
    By the definition of foliation, regular charts map the leaves to horizontal lines. To obtain a chart as in the statement, we reparametrize these charts in vertical direction so that $\mu$ is induced by $|\dd y|$, the Lebesgue measure in vertical direction.

    So for a chart $(U, \varphi)$, fix a horizontal $H$ in $\varphi(U)$, i.e., fix some number $h$ such that
    \begin{equation*}
        H = \{(x,h) \mid x \in \R\} \cap \varphi(U) \neq \emptyset
    \end{equation*}
    By considering a compatible atlas if necessary, we may assume that $\varphi(U)$ is given by a rectangle, so that for any point $(x,y) \in \varphi(U)$ the line perpendicular to $H$ is contained in $\varphi(U)$. This allows us to define the map
    \begin{equation*}
        \Phi(y) = h + \int_y^h \rho(t) \, \dd t,
    \end{equation*}
    where $\rho$ denotes the density of $\mu$ with respect to the Lebesgue measure. We claim that the chart $(U, \Phi\circ(\operatorname{Id},\varphi))$ has the property that the transverse measure $\mu$ is given by pulling back $|\dd y|$ on $(\operatorname{Id},\Phi)(\varphi(U))$. Indeed, we have
    \begin{align*}
        \mu(\alpha) &= \int_{\varphi\left(\alpha\left(0\right)\right)}^{\varphi\left(\alpha(1)\right)} \rho(t) \, \dd t \\&= 
        h + \int_h^{\varphi(\alpha(1))}\rho(t) \, \dd t -\left( h + \int_h^{\varphi(\alpha(0))} \rho(t) \, \dd t\right) \\&=
        \big|\Phi\left(\varphi\left(\alpha\left(1\right)\right)\right) - \Phi\left(\varphi\left(\alpha_2\left(0\right)\right)\right)\big|,
    \end{align*}
    from which the claim follows by the fundamental theorem of calculus. 
\end{proof}

Thus we lose no generality in assuming that we have an atlas consisting of charts such that $\mu$ is induced by $|\dd y|$. This fact has an immediate consequence.

\begin{corollary}\label{cor:natural_charts}
    Given any measured foliation $(\mathcal{F}, \mu)$ on a surface $X$, we can construct an atlas such that the charts are \emph{natural}, i.e., 
    \begin{enumerate}
        \item transition maps between charts are of the form $(x,y) \mapsto (f(x,y), \pm y + c)$ for some constant $c$ depending on the charts, and
        \item $\mu$ is induced by $|\dd y|$.
    \end{enumerate}
\end{corollary}
\begin{proof}
    We claim that the atlas consisting of the charts constructed in Proposition \ref{prop:transverse_measure_is_pullback} consists of natural charts. It remains to show the first assertion. 

    Consider two charts $(U_1, \varphi_1), (U_2, \varphi_2)$ with nontrivial overlap and write $\gamma = \varphi_2 \circ \varphi_1^{-1}$ for the transition map. We know that $\gamma$ maps horizontal lines to horizontal lines, i.e., 
    \begin{equation*}
        \gamma(x,y) = (f(x,y), g(y)),
    \end{equation*}
    for some smooth functions $f$ and $g$. Suppose that $g(y) \neq \pm y +c$ for some constant $c$, so that we are able to find some transverse arc $\alpha$ such that its height, or more precisely, the Lebesgue measure $|\dd y|$ is different. But this contradicts the fact that $\mu$ is induced by $|\dd y|$ in \emph{every} chart.
\end{proof}

In order to define pseudo-Anosov maps, we will need to consider not just one but a pair of measured foliations simultaneously.

\begin{figure}[ht]
    \centering
    \begin{tikzpicture}

    \begin{scope}
    \draw[] (0,0) -- ++(90:2) coordinate (e1);\draw[] (0,0) -- ++(210:2) coordinate (e2); \draw[] (0,0) -- ++ (330:2) coordinate (e3);

    \path (0,0) -- ++ (30:0.4) coordinate (a1);
    \path (0,0) -- ++(150:0.4) coordinate (a2);
    \path (0,0) -- ++(-90:0.4) coordinate (a3);

    \path (0,0) -- ++ (30:0.7) coordinate (b1);
    \path (0,0) -- ++(150:0.7) coordinate (b2);
    \path (0,0) -- ++(-90:0.7) coordinate (b3);

    \draw [] plot [smooth, tension=0.4] coordinates {($(e1)+ (0.2,0)$) (a1) ($(e3)+(0,0.2)$)};
    \draw [] plot [smooth, tension=0.4] coordinates {($(e1)- (0.2,0)$) (a2) ($(e2)+(0,0.2)$)};
    \draw [] plot [smooth, tension=0.4] coordinates {($(e2) - (0,0.2)$) (a3) ($(e3)-(0,0.2)$)};

    \draw [] plot [smooth, tension=0.6] coordinates {($(e1)+ (0.4,0)$) (b1) ($(e3)+(0,0.4)$)};
    \draw [] plot [smooth, tension=0.6] coordinates {($(e1)- (0.4,0)$) (b2) ($(e2)+(0,0.4)$)};
    \draw [] plot [smooth, tension=0.6] coordinates {($(e2)- (0,0.4)$) (b3) ($(e3)-(0,0.4)$)};
    \end{scope}
    \begin{scope}[noamblue, rotate = 180]
    \draw[] (0,0) -- ++(90:2) coordinate (e1);\draw[] (0,0) -- ++(210:2) coordinate (e2); \draw[] (0,0) -- ++ (330:2) coordinate (e3);

    \path (0,0) -- ++ (30:0.4) coordinate (a1);
    \path (0,0) -- ++(150:0.4) coordinate (a2);
    \path (0,0) -- ++(-90:0.4) coordinate (a3);

    \path (0,0) -- ++ (30:0.7) coordinate (b1);
    \path (0,0) -- ++(150:0.7) coordinate (b2);
    \path (0,0) -- ++(-90:0.7) coordinate (b3);

    \draw [] plot [smooth, tension=0.4] coordinates {($(e1)+ (0.2,0)$) (a1) ($(e3)+(0,0.2)$)};
    \draw [] plot [smooth, tension=0.4] coordinates {($(e1)- (0.2,0)$) (a2) ($(e2)+(0,0.2)$)};
    \draw [] plot [smooth, tension=0.4] coordinates {($(e2) - (0,0.2)$) (a3) ($(e3)-(0,0.2)$)};

    \draw [] plot [smooth, tension=0.6] coordinates {($(e1)+ (0.4,0)$) (b1) ($(e3)+(0,0.4)$)};
    \draw [] plot [smooth, tension=0.6] coordinates {($(e1)- (0.4,0)$) (b2) ($(e2)+(0,0.4)$)};
    \draw [] plot [smooth, tension=0.6] coordinates {($(e2)- (0,0.4)$) (b3) ($(e3)-(0,0.4)$)};
    \end{scope}
\end{tikzpicture}
    \caption{Two transverse measured foliations in a neighborhood of a singular point corresponding to a 3-pronged saddle.}
    \label{fig:transverse_measure_singular_point}
\end{figure}
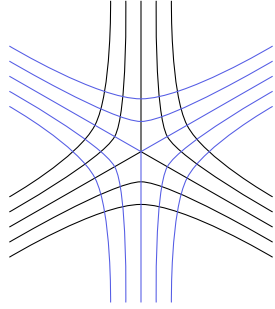

\begin{definition}[Transverse measured foliations]\label{def:transverse_measured_foliations}
    Let $X$ be a surface. We say that two measured foliations $(\mathcal{F}_1, \mu_1)$ and $(\mathcal{F}_2, \mu_2)$ are \emph{transverse}, if their leaves are transverse away from the singularities. 
\end{definition}
Note that it follows immediately from the definition that transverse measured foliations must have the same set of singularities. Indeed, if $p \in X$ was a singular point for just one of the foliations, the other one could not possibly be transverse at $p$. A neighborhood of a (common) singular point $p \in X$ is depicted in Figure \ref{fig:transverse_measure_singular_point}.

The concept of natural charts extends to the case were we endow the surface $X$ with a pair of transverse measured foliations.

\begin{definition}[Natural atlas]
    Let $X$ be a surface endowed with a pair of transverse measured foliations $(\mathcal{F}_1, \mu_1)$ and $(\mathcal{F}_2, \mu_2)$. A \emph{natural atlas} is an atlas, such that, away from the singular points, the foliation $\mathcal{F}_1$ is the pullback of the horizontal foliation of $\R^2$, $\mu_1$ is the pullback of $|\dd y|$, $\mathcal{F}_2$ is the pullback of the vertical foliation of $\R^2$ and $\mu_2$ is the pullback of $|\dd x|$. Each chart of a natural atlas is called \emph{natural with respect to the pair of measured foliations}.
\end{definition}

We will now show how one can construct a natural atlas given \emph{any} two transverse measured foliations. 

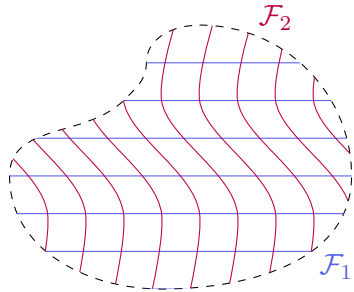
\begin{figure}[ht]
    \centering
    \begin{tikzpicture}
    \node[noamblue] at (1.8,-1.2) {$\mathcal{F}_1$};
    \node[purple] at (1,2.1) {$\mathcal{F}_2$};

    \draw[dashed] (0,2) to [out=0, in = 90] (2,0) to [out=-90, in = 0] (-0.5,-1.5) to [out=180, in = -90] (-2.5,0) to [out = 90, in = 225] (-1, 1) to [out = 45, in = -90] (-0.7,1.5) to [out=90, in = 180] (0,2);
    \clip (0,2) to [out=0, in = 90] (2,0) to [out=-90, in = 0] (-0.5,-1.5) to [out=180, in = -90] (-2.5,0) to [out = 90, in = 225] (-1, 1) to [out = 45, in = -90] (-0.7,1.5) to [out=90, in = 180] (0,2);

    \foreach \p in {-2,-1.5,-1,-0.5,0,0.5,1,1.5}{
        \draw[noamblue] (-3, \p) -- (3, \p);
    }

    \foreach \p in {-2.5,-2,-1.5,-1,-0.5,0,0.5,1,1.5,2}{
        \draw[purple] plot [smooth, tension = 0.5] coordinates {(\p,3) (\p-0.5,1) (\p+0.5,-0.5) (\p,-3)};
    }

\end{tikzpicture}
    \caption{Charts which are natural with respect to the measured foliation $(\mathcal{F}_1, \mu_1)$, which is given by the pullback of the horizontal lines in $\R^2$, before straightening out the verticals.}
    \label{fig:regular_neighborhood_almost_natural}
\end{figure}

If we now have a pair of transverse measured foliations $(\mathcal{F}_1, \mu_1), (\mathcal{F}_2, \mu_2)$, Corollary \ref{cor:natural_charts} allows us to obtain an atlas with charts that are natural with respect to $(\mathcal{F}_1, \mu_1)$. By transversality, neighborhoods of regular points are of the form depicted in Figure \ref{fig:regular_neighborhood_almost_natural}, i.e., the projection to the vertical axis of any leaf is injective. 

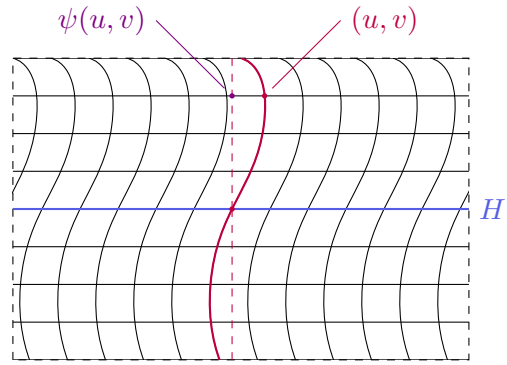
\begin{figure}[ht]
    \centering
    \begin{tikzpicture}
    \begin{scope}
    \draw[dashed] (-3,-2) rectangle (3,2);
    \clip (-3,-2) rectangle (3,2);

    \foreach \p in {-3,-2.5,...,3}{
        \draw (-3,\p) -- (3,\p);
        \draw plot [smooth, tension = 1] coordinates {(\p,2) (\p+0.3,1.1) (\p-0.4,-1) (\p,-2.5)};
    }
    \draw[thick, purple, name path = red] plot [smooth, tension = 1] coordinates {(0,2) (0.3,1.1) (-0.4,-1) (0,-2.5)};
    
    \draw[name path = H, noamblue, thick] (-3, 0) -- (3,0);
    \path[name intersections={of=H and red, by=int}];
    \draw[purple, dashed, name path = horz] (int) -- ++(90:2) (int) -- ++(-90:2);

    \fill[purple] (int) circle (1pt);

    \path[name path = vert1] (-3,1.5) -- (3,1.5);
    \path[name intersections={of=red and vert1, by=int2}];
    \fill[purple] (int2) circle (1pt);

    \path[name intersections={of=horz and vert1, by=int3}];
    \fill[violet] (int3) circle (1pt);
    \end{scope}
    \draw[purple] ($(int2) + (0.1,0.1)$) -- ($(int2) + (1,1)$) node[right]{$(u,v)$};
    \draw[violet] ($(int3) + (-0.1,0.1)$) -- ($(int3) + (-1,1)$) node[left] {$\psi(u,v)$};
    \node[right, noamblue] at (3,0) {$H$}; 
\end{tikzpicture}
    \caption{The map $\Psi$ which straightens out the foliation $\mathcal{F}_2$ in charts.}
    \label{fig:straighten_vertical_flow}
\end{figure}

The idea is now to \enquote{straighten out} the second foliation $(\mathcal{F}_2, \mu_2)$. We apply an analogous construction as in the proof of the existence of so-called tubular neighborhoods for flows. Let $p \in X$ be regular and $(U, \varphi)$ a chart such that $p \in U$. For simplicity, we will assume that the foliation $\mathcal{F}_2$ is globally orientable, so that we may view it as a flow $(f^t)_{t\in I}$ associated to a vector field. If that is not the case, we can move to the orientation double cover (see Definition \ref{def:orientation_cover}) as we have done for example in the proof of Proposition \ref{prop:euler_poincare}. Fix a reference horizontal $H$ in the chart, that is, fix some $h \in \R$ such that $H = \{(x,h)\mid x \in \R\}$ intersects $\varphi(U)$ as in Figure \ref{fig:straighten_vertical_flow}. Considering a compatible atlas if necessary, we may assume that any point $(u,v) \in \varphi(U)$ can now be written as $f^t(h_1, h_2)$, where $(h_1, h_2) \in H$. Next, we simply want to map any point $(u,v)$ to a point at the same height, changing the first coordinate so it is straight above the intersection point of the flow $f^t$ and $H$. Therefore, we define
\begin{align*}
    \Psi \colon \varphi(U) &\to \Psi(\varphi(U)), \\
    (u,v) &\mapsto (h_1, v).
\end{align*}
The map $\Psi$ inherits all regularity properties from the flow $f^t$, so in particular it is a diffeomorphism from $\varphi(U)$ to $\Psi(\varphi(U))$. This construction yields charts $(U, \Psi\circ\varphi)$ such that the foliation $\mathcal{F}_1$ is horizontal and the foliation $\mathcal{F}_2$ is vertical. It remains to ensure that the measure $\mu_2$ associated to the foliation $\mathcal{F}_2$ is the pullback of the measure by $|\dd x|$. This is achieved by a reparametrization that is completely analogous to the one used in the proof of Proposition \ref{prop:transverse_measure_is_pullback} adapted to the horizontal direction.

The construction above proves the following theorem.

\begin{theorem}[Existence of natural atlas]\label{thm:existence_natural_atlas}
    Given a pair of transverse measured foliations $(\mathcal{F}_1, \mu_1)$, $(\mathcal{F}_2, \mu_2)$ on a compact surface $X$, there exists an atlas of $X$ which is natural with respect to the pair of measured foliations.
\end{theorem}

Applying Corollary \ref{cor:natural_charts} to the second foliation $\mathcal{F}_2$ obtained by pulling back vertical lines in $\R^2$, we see that in this natural atlas the transition maps are of the form
\begin{equation}\label{eq:half_translation_transitions}
    (x,y) \mapsto (\pm x +c_1, \pm y + c_2),
\end{equation}
for some constants $c_1, c_2 \in \R$. Since $X$ is assumed to be oriented it follows that the signs of $x$ and $y$ must always agree, i.e., are either both negative or both positive. This means that, away from the singular points, we have recovered the definition of a \emph{half-translation surface}. It remains to verify that the singular points correspond to conical singularities of cone angle $k\pi$, with $k \geq 3$. One way to see this, is to use the proof of the equivalence between the definition of half-translation surfaces as glued polygons and the definition relying on having an atlas with transition maps given by \eqref{eq:half_translation_transitions}. Here, we follow the exposition from \cite{massart2022short}.

As above, we denote the surface by $X$ and we write $\Sigma$ = $\{p_1, \ldots, p_n\}$ for the set of singular points. Maps of the form \eqref{eq:half_translation_transitions} preserve the Euclidean metric, which allows us to equip $X \setminus \Sigma$ with a Riemannian metric which is locally Euclidean. In particular, its geodesics are locally straight lines. Even though the Riemannian metric might not extend to the whole surface $X$, the induced distance function does. It is therefore possible to connect any two singularities by a geodesic. 

\begin{figure}[ht]
    \centering
    \begin{tikzpicture}[scale = 0.8]
    \begin{scope}
    \clip (-3,-2.5) rectangle (3,2.5);
    \draw (0,0) ellipse (3 and 1.5);
    \begin{scope}
        \clip (0,-1.8) ellipse (3 and 2.5);
        \draw (0,2.2) ellipse (3 and 2.5);
    \end{scope}
    \begin{scope}
        \clip (0,2.2) ellipse (3 and 2.5);
        \draw (0,-2.2) ellipse (3 and 2.5);
    \end{scope}
    \draw[purple, thick, dashed, rotate = 5] (0.5,-0.91) ++(270:0.3 and 0.6) arc (270:90:0.3 and 0.6);
    \draw[purple, thick, rotate = 5, name path = purple] (0.5,-0.91) ++(90:0.3 and 0.6) arc (90:-90:0.3 and 0.6);

    \draw[noamblue, thick, name path = blue] (0,0)  ellipse (2.4 and 1.1);

    \path[name intersections={of=blue and purple, by=int1}];
    \fill[] (int1) circle (1.5pt);
    \draw[] (0,1.1) node {$\times$};
    \draw[] (-2.31,-0.3) circle (1.5pt);
    \end{scope}
    \draw[thick, >=latex, ->] (4,0) -- (5,0);

    \def\x{6};\def\y{-1}

    \draw[thick, purple] (\x,\y) -- (\x+2, \y+2); 
    \draw[thick, purple] (\x+6,\y) -- (\x+8, \y+2); 

    \draw[thick, noamblue] (\x+2, \y+2) -- (\x+8, \y+2); 
    \draw[thick, noamblue] (\x,\y) -- (\x+6,\y);

    \fill[] (\x,\y) circle (1.5pt){};
    \fill[] (\x+6,\y) circle (1.5pt){};
    \fill[] (\x+2,\y+2) circle (1.5pt){};
    \fill[] (\x+8,\y+2) circle (1.5pt){};

    \draw[] (\x+2,\y) node {$\times$};
    \draw[] (\x+4,\y+2) node {$\times$};

    \draw[] (\x+4,\y) circle (1.5pt);
    \draw[] (\x+6,\y+2) circle (1.5pt);

    \draw[] (\x+3,\y+2) node {\footnotesize /};
    \draw[] (\x+1,\y) node {\footnotesize /};
    \draw[] (\x+5,\y+2) node {\footnotesize //};
    \draw[] (\x+3,\y) node {\footnotesize //};
    \draw[] (\x+7,\y+2) node {\footnotesize ///};
    \draw[] (\x+5,\y) node {\footnotesize ///};
\end{tikzpicture}
    \caption{Cutting open a surface along Euclidean geodesics between singular points.}
    \label{fig:cutting_translation_surface}
\end{figure}
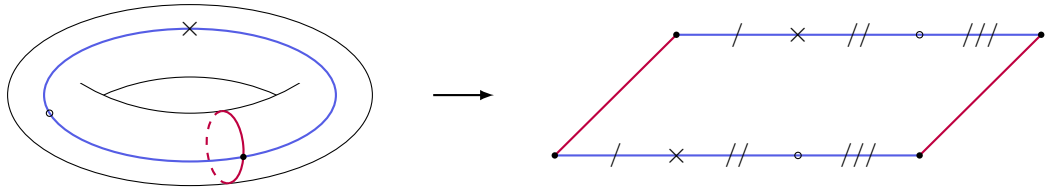

Next, we draw enough geodesics between singularities, such that the connected components of the complement of these geodesics in $X$ are simply connected. These simply connected regions may now be developed to the plane, i.e., we can lay them flat onto the plane. Since the geodesics are straight lines, we obtain a collection of polygons such that $X$ is obtained by identifying the sides of these polygons pairwise by translations or by half-translations. This is illustrated in Figure \ref{fig:cutting_translation_surface}. Thus, the structure on $X$ is indeed a half-translation structure, so in particular the singularities are conical singularities with cone angles $k \pi$.

Lastly, we will be interested in translation surfaces only, so assuming that the given pair of measured foliations is locally orientable also in the neighborhoods of singular points gives the surface $X$ the structure of a translation surface. In short, we have justified that given any pair of transverse measured foliations on a surface $X$, we may assume that the foliations have already been straightened, i.e., we use charts in which the foliations are horizontal and vertical respectively, in the framework of (half-) translation surfaces. 

The arguments at the end of the above discussion suggest a way to construct measured foliations on a closed surface. This particular construction is taken from \cite{farb2011primer}. Given any closed surface $X$, it is possible to realize $X$ as the quotient of a polygon $P$ in $\R^2$, where the sides of the polygon are identified in some suitable way. We emphasize this by writing $X = P / {\sim}$. Additionally, we impose the conditions that any two edges of $P$ that are identified are parallel, and that the total Euclidean angle round each point of $S$ is strictly greater than $\pi$. Notice that the second condition needs to be checked only at the vertices of $P$. Clearly, any foliation of $\R^2$ by parallel lines restricts to give a foliation of the interior of $P$. Let us show that this foliation induces a foliation of $X$.

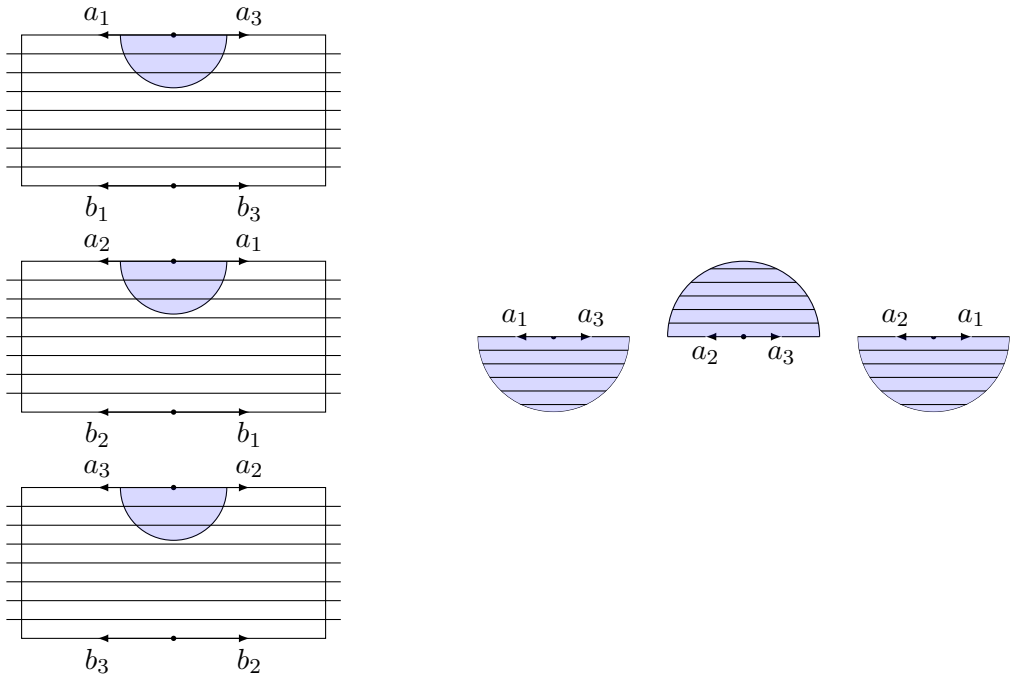
\begin{figure}[ht]
    \centering
    \begin{tikzpicture}
    \draw (0,0) rectangle (4,-2);
    \draw (0,-3) rectangle (4,-5);
    \draw (0,-6) rectangle (4,-8);

    \foreach \p in {0,-2,-3,-5,-6,-8}{
        \fill[] (2,\p) circle (1pt);
    }

    \foreach \p in {0,-3,-6}{
        \draw[] (2,\p) ++(0:0.7) arc (0:-180:0.7);
        \fill[blue, opacity = 0.15] (2,\p) ++(0:0.7) arc (0:-180:0.7) -- (2,\p) ++ (0:0.7) -- cycle;
    }

    \foreach \p in {1/4,1/2,3/4,1,5/4, 3/2, 7/4}{
        \draw[] (-0.2, -\p) -- (4.2,-\p);
        \draw[] (-0.2, -\p-3) -- (4.2,-\p-3);
        \draw[] (-0.2, -\p-6) -- (4.2,-\p-6);
    }

    \draw[->, >=latex, thin] (2,0) -- (3,0) node[above] {$a_3$};
    \draw[->, >= latex,thin] (2,-6) -- (1,-6) node[above] {$a_3$};

    \draw[->, >=latex, thin] (2,0) -- (1,0) node[above] {$a_1$};
    \draw[->, >= latex,thin] (2,-3) -- (3,-3) node[above] {$a_1$};

    \draw[->, >=latex, thin] (2,-3) -- (1,-3) node[above] {$a_2$};
    \draw[->, >= latex,thin] (2,-6) -- (3,-6) node[above] {$a_2$};

    \draw[->, >=latex, thin] (2,-2) -- (3,-2) node[below] {$b_3$};
    \draw[->, >= latex,thin] (2,-8) -- (1,-8) node[below] {$b_3$};

    \draw[->, >=latex, thin] (2,-5) -- (3,-5) node[below] {$b_1$};
    \draw[->, >= latex,thin] (2,-2) -- (1,-2) node[below] {$b_1$};

    \draw[->, >=latex, thin] (2,-8) -- (3,-8) node[below] {$b_2$};
    \draw[->, >= latex,thin] (2,-5) -- (1,-5) node[below] {$b_2$};

    \def\x{7};

    \draw[->, >=latex] (\x,-4) -- (\x-1/2, -4) node[pos=1, above] {$a_1$};
    \draw (\x-1/2, -4) -- (\x-1, -4);
    \draw[->, >=latex] (\x,-4) -- (\x+1/2, -4) node[pos=1, above] {$a_3$};
    \draw (\x+1/2, -4) -- (\x+1, -4);
    \draw[->, >=latex] (\x+2.5,-4) -- (\x+2.5-1/2, -4) node[pos=1, below] {$a_2$};
    \draw (\x+2.5-1/2, -4) -- (\x+2.5-1, -4);
    \draw[->, >=latex] (\x+2.5,-4) -- (\x+2.5+1/2, -4) node[pos=1, below] {$a_3$};
    \draw (\x+2.5+1/2, -4) -- (\x+2.5+1, -4);
    \draw[->, >=latex] (\x+5,-4) -- (\x+5-1/2, -4) node[pos=1, above] {$a_2$};
    \draw (\x+5-1/2, -4) -- (\x+5-1, -4);
    \draw[->, >=latex] (\x+5,-4) -- (\x+5+1/2, -4) node[pos=1, above] {$a_1$};
    \draw (\x+5+1/2, -4) -- (\x+5+1, -4);
    
    \fill (\x+2.5, -4) circle (1pt);
    \draw (\x+2.5,-4) ++ (0:1) arc (0:180:1);
    \fill[blue, opacity = 0.15] (\x+2.5, -4) ++(0:1) arc (0:180:1) -- cycle;

    \begin{scope}
    \clip (\x, -4) ++(0:1) arc (0:-180:1) -- cycle;
    \foreach \p in {0,5}{
        \fill (\p+\x,-4) circle (1pt);
        \draw (\p+\x, -4) ++(0:1) arc (0:-180:1);
        \fill[blue, opacity = 0.15] (\p+\x, -4) ++(0:1) arc (0:-180:1) -- cycle;

        \foreach \q in {0.18,0.36,0.54,0.72, 0.9}{
            \draw (\p+\x - 1.2, -\q-4) -- (\p+\x + 1.2, -\q-4);
        }
    }
    \end{scope}
    \begin{scope}
    \clip (\x+5, -4) ++(0:1) arc (0:-180:1) -- cycle;
    \foreach \p in {0,5}{
        \fill (\p+\x,-4) circle (1pt);
        \draw (\p+\x, -4) ++(0:1) arc (0:-180:1);
        \fill[blue, opacity = 0.15] (\p+\x, -4) ++(0:1) arc (0:-180:1) -- cycle;

        \foreach \q in {0.18,0.36,0.54,0.72, 0.9}{
            \draw (\p+\x - 1.2, -\q-4) -- (\p+\x + 1.2, -\q-4);
        }
    }
    \end{scope}

    \begin{scope}
    \clip (\x+2.5, -4) ++(0:1) arc (0:180:1) -- cycle;
    \foreach \q in {0.18,0.36,0.54,0.72, 0.9}{
        \draw (2.5+\x-1.2, \q-4) -- (2.5+\x+1.2, \q-4); 
    }        
    \end{scope}
\end{tikzpicture}
    \caption{3 half-disks glued together, each foliated by lines parallel to the diameter.}
    \label{fig:half_disks_glued}
\end{figure}

\begin{proposition}
    Any foliation of $\R^2$ induces a foliation on $X = P / {\sim}$. 
\end{proposition}\label{prop:polygon_foliation}
 \begin{proof}
     It is easy to see that any point of $X$ that is not coming from a vertex of $P$ has a regular neighborhood satisfying the definition of regular point from Definition \ref{def:singular_foliation}. At a point $p \in X$ corresponding to a vertex of $P$, since identified sides of $P$ are parallel it follows that the total angle around $p$ is an integer multiple of $\pi$, i.e., it is of the form $k\pi$ for $k \in \N_{\geq 1}$. Moreover, there is some vertex $\Tilde{p}$ of $P$ in the preimage of $p$ under the quotient map, and a vector $\mathrm{v}$ based at $\Tilde{p}$ that is parallel to the foliation of $P$ and points to the interior of $P$, possibly along some edge.

     If we sweep out an angle of $\pi$ starting with $\mathrm{v}$, we obtain a closed Euclidean half-disk in $X$ that is foliated by lines parallel to the diameter. Continuing this way, we see that a neighborhood of $p$ looks like some number of Euclidean half-disks, each foliated by lines parallel to the diameter, and glued along (oriented) radii. See Figure \ref{fig:half_disks_glued} for an illustration of the case $k = 3$. Note that the second assumption we made above ensures that we need at least two half-disks glued at $p$. If there are exactly two half-disks, $p$ is a regular point. If $k \geq 3$, then $p$ is a singular point with $k$ prongs.  
 \end{proof}

We can endow the foliation constructed above with the transverse measure given by the total variation of the Euclidean distance in the direction perpendicular to the foliation of $P$. The charts described in Proposition \ref{prop:polygon_foliation} are then exactly the natural charts for the nonsingular points.

Of course, what we described here is nothing but a construction of a foliation on a half-translation surface. If we now further insist that each edge of $P$ is oriented in a way such that the identifications respect these orientations, meaning that the identifications of the sides are done by translations, we obtain a translation surface and the resulting foliation of $X$ (globally) orientable. Indeed, choosing any of the two orientation of the foliation of $\R^2$ restricted to the interior of $P$ extends to a global orientation of the entire foliation of $X$. Moreover, by Proposition \ref{prop:local_orientation_foliation} and the fact that translation surfaces have conical singularities of cone angle $2\pi k$ for some $k \geq 1$ it follows that in this setting, local and global orientability coincide. The fact that this is the case for any \emph{measured} foliation is due to the fact that \emph{every} measured foliation comes from this polygon construction we have described above. The proof of this fact is in the same spirit as the construction of a natural atlas we gave above and relies on the notion of a \emph{good atlas}. We refer the reader to chapter 14.3 of \cite{farb2011primer} for a detailed discussion.

\subsubsection{The Mapping Class Group}
It is not strictly necessary to introduce the so-called \emph{mapping class group} to give the definition of a pseudo-Anosov map, but it will play an important role when establishing the correspondence to closed Teichmüller geodesics which is why we give the relevant definitions here.

\begin{definition}[Mapping class group]\label{def:mapping_class_group}
    For a connected, closed and orientable surface $X$, we write $\operatorname{Homeo^+(X)}$ for the group of orientation-preserving (or positive) homeomorphisms of $X$. Further we denote by $\operatorname{Homeo_0(X)}$ the connected component of the identity, or equivalently the homeomorphisms of $X$ that are isotopic to the identity. This is a normal subgroup of the positive homeomorphisms of $X$ and the \emph{mapping class group of $X$} is defined as the quotient
    \begin{equation*}
        \operatorname{Mod(X) = \operatorname{Homeo}^+(X) / \operatorname{Homeo}_0}(X).
    \end{equation*}
\end{definition}

The definition of Mod(X) can vary in several ways, for instance we could consider diffeomorphisms instead of homeomorphisms. Moreover, it is also possible to take the quotient over maps homotopic (instead of isotopic) to the identity. These facts can be justfied using the following two results, the first which was proven by Baer in the 1920s (see \cite{baer1928isotopie}) and the second which was proven by Munkres in the 1950s (see Theorem 6.3 in \cite{munkres1960obstructions}). 

\begin{theorem}[\cite{baer1928isotopie}]\label{thm:baer_isotopic}
    Let $X$ be any compact surface and let $f$ and $g$ be homotopic orientation-preserving homeomorphisms of $X$. Then, $f$ and $g$ are isotopic.
\end{theorem}

\begin{theorem}[\cite{munkres1960obstructions}]\label{thm:munkres_homeo_diffeo_isotopic}
    Let $X$ be a compact surface. Then, every homeomorphism of $X$ is isotopic to a diffeomorphism of $X$.
\end{theorem}

\begin{corollary}
    For any compact surface $X$, there are isomorphisms
    \begin{equation*}
        \operatorname{Mod}(X) \simeq \operatorname{Diff^+}(X)/\operatorname{Diff}_0(X) \simeq \operatorname{Homeo}^+(X) / \operatorname{Homeo}'_0(X),
    \end{equation*}
    where $\operatorname{Diff}^+(X)$ is the group of orientation-preserving (or positive) diffeomorphisms on $X$, $\operatorname{Diff}_0(X)$ is the normal subgroup consisting of elements isotopic to the identity and $\operatorname{Homeo}'_0(X)$ is the set of positive homeomorphisms of $X$ homotopic to the identity. 
\end{corollary}

\begin{proof}
    Since every diffeomorphism is a homeomorphism and every smooth isotopy is an isotopy, there is an obvious homomorphism
    \begin{equation*}
        \varphi \colon \operatorname{Diff}^+(X) / \operatorname{Diff}_0(X) \to \operatorname{Mod}(X). 
    \end{equation*}
    Similarly, there is a surjective homomorphism
    \begin{equation*}
        \psi \colon \operatorname{Mod}(X) \to \operatorname{Homeo}^+(X) / \operatorname{Homeo}'_0(X),
    \end{equation*}
    since any isotopy is in particular a homotopy. By Theorem \ref{thm:munkres_homeo_diffeo_isotopic} we know that $\varphi$ is surjective, whereas Theorem \ref{thm:baer_isotopic} implies that the composition $\psi \circ \varphi$ is injective. Then, $\varphi$ must be injective as well, as for any element $[f] \in \operatorname{ker}(\varphi)$ we have $[f] \in \ker (\psi \circ \varphi)$ as well, so that $[f] = [\operatorname{Id}]$ by injectivity of the composition. We have established the claimed isomorphisms.
\end{proof}

There is one last modification to the definition of the mapping class group we wish to make in the situation that we have a set of marked points $\Sigma \subseteq X$. For example, $\Sigma$ could be the set of conical singularities of a translation surface $X$. In this case, $\operatorname{Mod}(X)$ is the group of homeomorphisms (or diffeomorphisms) on $X$ that leave the set of marked points invariant modulo isotopies (or homotopies) that leave the set of marked points invariant. All of the identifications above still hold under this modified definition.

\subsubsection{Pseudo-Anosov Mapping Classes}

We have now all the prerequisites to properly define pseudo-Anosov maps. 

\begin{definition}[Pseudo-Anosov maps and classes]\label{def:pseudo_anosov_map}
    Let $X$ be a closed, connected and orientable surface. An element $[f] \in \operatorname{Mod}(X)$ is called a \emph{pseudo-Anosov mapping class}, if there is a pair of transverse measured foliations $(\mathcal{F}^u, \mu_u)$ and $(\mathcal{F}^s, \mu_s)$ on $X$, a number $\lambda > 1$ and a representative homeomorphism $\phi \in [f]$ so that
    \begin{equation}\label{eq:pA_foliations}
        \phi\cdot (\mathcal{F}^u, \mu_u) = (\mathcal{F}^u, \lambda \mu_u) \quad \text{and} \quad \phi \cdot  (\mathcal{F}^s, \mu_s) = (\mathcal{F}^s, \lambda^{-1}\mu_s),
    \end{equation}
    where the measure obtained by the action of $\phi$ is the pushforward of $\mu_\varepsilon$ by $\phi$ for $\varepsilon \in \{u, s\}$. The measured foliations $(\mathcal{F}^u, \mu_u)$ and $(\mathcal{F}^s, \mu_s)$ are called the \emph{unstable foliation} and the \emph{stable foliation} respectively. We will refer to $\lambda$ as the \emph{dilatation} of $\phi$ (or of $f$). We call the homeomorphism $\phi$ a \emph{pseudo-Anosov map}. 
\end{definition}
The equations \eqref{eq:pA_foliations} show that the basic behavior from the CAT map we presented in the beginning of this section carries over to the general definition of pseudo-Anosov map, i.e., the leaves of the foliation are stretched on one direction and contracted in the other direction. As we mentioned in the definition of the mapping class group (Definition \ref{def:mapping_class_group}), there is some freedom with respect to the required regularity of the maps, i.e., $\phi$ could be a diffeomorphism instead of a homeomorphism. For this reason, one also finds the terminology of pseudo-Anosov \emph{diffeomorphism} in the literature. 

The discussion above, in particular Theorem \ref{thm:existence_natural_atlas}, imply that an equivalent way to formulate the definition of a pseudo-Anosov map is given by the following.

\begin{definition}[Pseudo-Anosov map, equivalent definition]\label{def:pseudo_Anosov_equivalent_definition} A homeomorphism $\phi \colon X \to X$ is a \emph{pseudo-Anosov map}, if $\phi(\Sigma) = \Sigma$ and it is an affine diffeomorphism on $X \setminus \Sigma$ for the Euclidean metric induced by the construction in Theorem \ref{thm:existence_natural_atlas}, where $\Sigma$ are the singular points of the given pair of transverse measured foliations, and its differential
\begin{equation*}
    \mathrm{D}\phi = \begin{bmatrix}
        \pm \lambda & 0 \\
        0 & \pm \lambda^{-1}
    \end{bmatrix}
\end{equation*}
satisfies $|tr(\mathrm{D}\phi)| > 2$. 
\end{definition}

To summarize, given any pseudo-Anosov mapping class $[f]$ we may always choose a representative $\phi\in[f]$ and a natural atlas that endows $X$ with a half-translation structure such that $\phi(\Sigma) = \Sigma$, and there exists $\lambda > 1$ such that for any $x \in X \setminus\Sigma$ one has the equality
\begin{equation*}
    \phi(y) = 
    \begin{bmatrix}
        \pm \lambda & 0 \\
        0 & \pm \lambda^{-1}
    \end{bmatrix}y
\end{equation*}
read in charts around $x$ and $\phi(x)$, where the signs depend on the choice of coordinate charts. Put differently, $\phi$ sends horizontal segments to horizontal segments and vertical segments to vertical segments while expanding by $\lambda$ in the horizontal direction and contracting by $\lambda^{-1}$ in the vertical direction. It follows in particular, that the Lebesgue measure is invariant under $\phi$. 

In the case of translation surfaces, the signs in the matrix above are \emph{global}, indicating whether $\phi$ changes the orientation of the pair of measured foliations. Since in all the constructions above we needed the measures on the foliations to be induced by $|\dd x|$ and $|\dd y|$, see in particular Corollary \ref{cor:natural_charts}, we see that we lose no generality in supposing that the global signs are positive. 


It is clear that the CAT map introduced above gives an example of a pseudo-Anosov map on the torus, where the foliations are non-singular. We want to mention that such maps are known as \emph{Anosov} maps, and were historically studied before the concept was generalized to surfaces of higher genus by Thurston who introduced the \emph{pseudo-Anosov} maps. In our context here, Anosov maps associated to nonsingular foliations play a secondary role appearing only as illustrative examples, so that we will refer to all maps of this form as \emph{pseudo}-Anosov maps. 

Let us proceed by giving two examples of a pseudo-Anosov maps on surfaces of genus $\mathbf{g} = 2$, so that we indeed need to work with a \emph{singular} foliation. The first example is taken from \cite{lanneau2017tell}. The underlying idea, based on the observation that it is often the case that the product of two parabolic elements is hyperbolic, is also known as the \emph{bouillabaisse} construction, introduced initially by Thurston and Veech and popularised by John Hubbard at a talk in Marseille in 2003. 

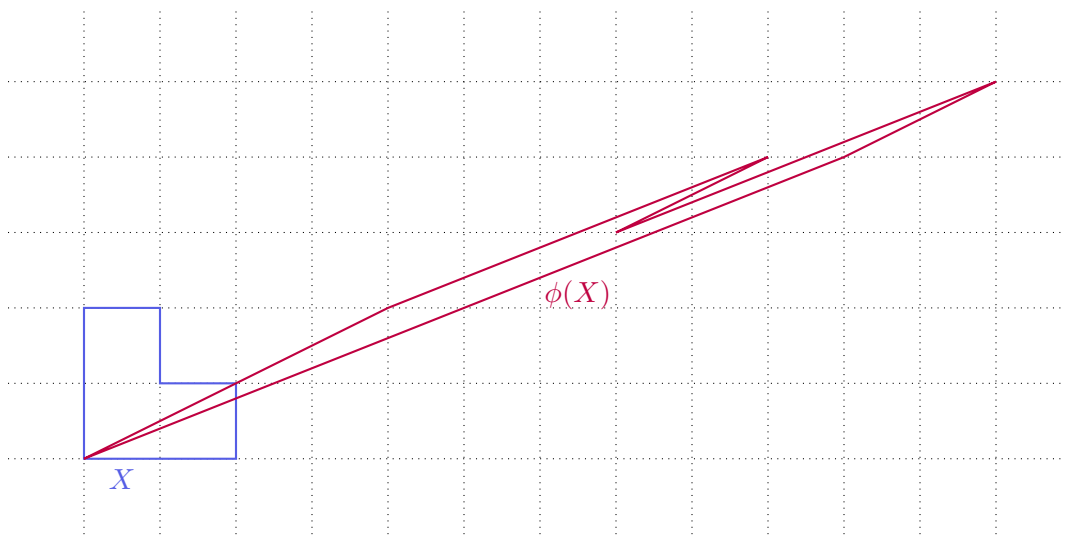
\begin{figure}[ht]
    \centering
    \begin{tikzpicture}[scale = 1]
    \foreach \x in {0,1,...,12}{
        \draw[dotted] (\x,-1) -- (\x,6);
    }    
    \foreach \y in {0,1,...,5}{
        \draw[dotted] (-1, \y) -- (13,\y);
    }
    \draw[thick, noamblue] (0,0) -- (2,0) -- (2,1) -- (1,1) -- (1,2) -- (0,2) -- cycle;
    \node[below, noamblue] at (1/2,0) {$X$};
    \draw[thick, purple] (0,0) -- (10,4) -- (12,5) -- (7,3) -- (9,4) -- (4,2) -- cycle;
    \node[below, purple] at (6.5,2.5) {$\phi(X)$};
\end{tikzpicture}
    \caption{The pseudo-Anosov map $\phi$ from Example \ref{ex:torus_triple_cover}.}
    \label{fig:torus_triple_cover_3}
\end{figure}

\begin{example}[Pseudo-Anosov map on triple cover of the torus]\label{ex:torus_triple_cover}
    Consider the translation surface $X$ given by the triple cover of the torus depicted on the left-hand side in Figure \ref{fig:torus_triple_cover}. We can decompose $X$ into two pieces which are given by a unit square (together with the necessary identifications) and a rectangle of height 1 and width 2. On the top square, the map
    \begin{equation*}
        T_{\alpha_1} \colon \begin{bmatrix}
            x \\ y
        \end{bmatrix} \mapsto \begin{bmatrix}
            1 & 1  \\
            0 & 1
        \end{bmatrix}\cdot  \begin{bmatrix}
            x \\ y
        \end{bmatrix}
    \end{equation*}
    acts as a diffeomorphism which fixes the upper and lower boundary of the square pointwise. Similarly, there is an automorphism on the lower rectangle given by
    \begin{equation*}
        T_{\alpha_2} \colon \begin{bmatrix}
            x \\ y 
        \end{bmatrix}
        \mapsto
        \begin{bmatrix}
            1 & 2 \\
            0 & 1
        \end{bmatrix} \cdot \begin{bmatrix}
            x \\ y
        \end{bmatrix}.
    \end{equation*}

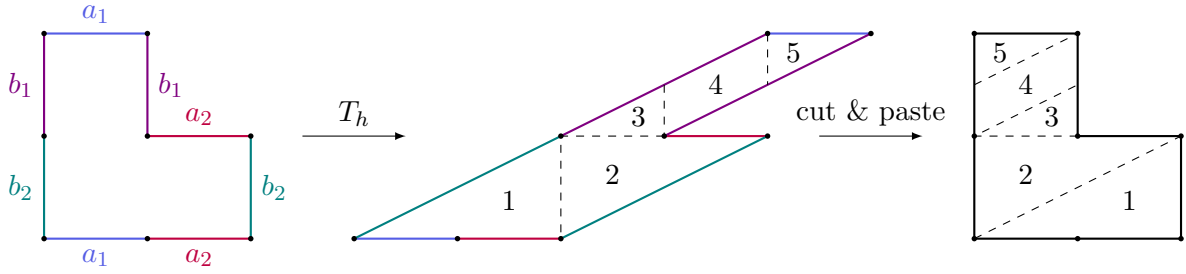
\begin{figure}[ht]
    \centering
    \begin{tikzpicture}[scale =0.68]
    \draw[thick, noamblue] (0,0) -- (2,0) node[midway, below] {$a_1$};
    \draw[thick, noamblue] (0,4) -- (2,4) node[midway, above] {$a_1$};
    \draw[thick, purple] (2,0) -- (4,0) node[midway, below] {$a_2$};
    \draw[thick, purple] (2,2) -- (4,2) node [midway, above] {$a_2$};
    \draw[thick, teal] (0,0) -- (0,2) node [midway, left] {$b_2$};
    \draw[thick, teal] (4,0) -- (4,2)node [midway, right] {$b_2$};
    \draw[thick, violet] (0,2) -- (0,4)node [midway, left] {$b_1$};
    \draw[thick, violet] (2,2) --(2,4)node [midway, right] {$b_1$};

    \fill (0,0) circle (1.5pt);
    \fill (2,0) circle (1.5pt);
    \fill (4,0) circle (1.5pt);
    \fill (0,2) circle (1.5pt);
    \fill (2,2) circle (1.5pt);
    \fill (4,2) circle (1.5pt);
    \fill (0,4) circle (1.5pt);
    \fill (2,4) circle (1.5pt);

    \draw[->, >=latex] (5,2) -- (7,2) node[midway, above] {$T_h$};

    \def\x{6};
    
    \draw[thick, noamblue] (\x,0) -- (\x+2,0);
    \draw[thick, noamblue](\x+8,4) -- (\x+10,4);
    \draw[thick, purple] (\x+2,0) -- (\x+4,0);
    \draw[thick, purple] (\x+6,2) -- (\x+8,2);
    \draw[thick, teal] (\x,0) -- (\x+4,2);
    \draw[thick, teal] (\x+4,0) -- (\x+8,2);
    \draw[thick, violet] (\x+4,2) -- (\x+8,4);
    \draw[thick, violet] (\x+6,2) -- (\x+10,4);

    \fill (\x+0,0) circle (1.5pt);
    \fill (\x+2,0) circle (1.5pt);
    \fill (\x+4,0) circle (1.5pt);
    \fill (\x+4,2) circle (1.5pt);
    \fill (\x+6,2) circle (1.5pt);
    \fill (\x+8,2) circle (1.5pt);
    \fill (\x+8,4) circle (1.5pt);
    \fill (\x+10,4) circle (1.5pt);

    \draw[dashed] (\x+8,4) -- (\x+8,3);
    \draw[dashed] (\x+6,3) -- (\x+6,2);
    \draw[dashed] (\x+6,2) -- (\x+4,2) -- (\x+4,0);

    \node at (\x+3,3/4) {1};
    \node at (\x+5,5/4) {2};
    \node at (\x+5.5,19/8) {3};
    \node at (\x+7, 3) {4};
    \node at (\x+8.5, 29/8) {5};

    \draw[>=latex, ->] (\x+9,2) -- (\x+11,2) node[midway, above] {cut \& paste}; 
    \def\w{18};

    \coordinate (p1) at (\w,0);
    \coordinate (p2) at (\w+2,0);
    \coordinate (p3) at (\w+4,0);
    \coordinate (p4) at (\w+4,2);
    \coordinate (p5) at (\w+2,2);
    \coordinate (p6) at (\w+2,4);
    \coordinate (p7) at (\w,4);
    \coordinate (p8) at (\w,2);

    \draw[thick] (p1) -- (p3) -- (p4)-- (p5) -- (p6) -- (p7) --cycle;

    \foreach \i in {1,2,...,8}{
        \fill (p\i) circle (1.5pt);
    }

    \draw[dashed] (\w,0) -- (\w+4,2);
    \draw[dashed] (\w+2,2) -- (\w,2) -- (\w+2,3);
    \draw[dashed] (\w,3) -- (\w+2,4);

    \node at (\w+3,3/4) {1};
    \node at (\w+1,5/4) {2};
    \node at (\w+1.5,19/8) {3};
    \node at (\w+1, 3) {4};
    \node at (\w+0.5, 29/8) {5};
    
\end{tikzpicture}
    \caption{The action of $T_h$ on the triple cover of the torus.}
    \label{fig:torus_triple_cover}
\end{figure}

    Performing the two maps simultaneously, which is well defined since the boundaries are fixed pointwise, yields a map
    \begin{equation*}
        T_h = T_{\alpha_2} \circ T^2_{\alpha_1},
    \end{equation*}
    which is a diffeomorphism on $X \setminus \Sigma$ with a constant differential given by
    \begin{equation*}
        DT_h = \begin{bmatrix}
            1 & 2 \\
            0 & 1
        \end{bmatrix}.
    \end{equation*}
    The action of $T_h$ is illustrated in Figure \ref{fig:torus_triple_cover}. By symmetry, the same construction works in vertical direction, i.e., we obtain an affine diffeomorphism $T_v \colon X \to X$ with
    \begin{equation*}
        DT_v = \begin{bmatrix}
            1 & 0 \\
            2 & 1
        \end{bmatrix}.
    \end{equation*}
    Finally, we define
    \begin{equation*}
        \phi = T_v \circ T_h,
    \end{equation*}
    which is an affine diffeomorphism away from the singular point and by the chain rule it follows that
    \begin{equation*}
        D\phi = \begin{bmatrix}
            5 & 2 \\
            2 & 1
        \end{bmatrix},
    \end{equation*}
    which satisfies $\tr(D\phi) = 6 > 2$ so that indeed $\phi$ is a pseudo-Anosov map. Computing the eigenpairs of $D\phi$ one further verifies that the dilatation of $\phi$ is given by $\lambda = 3 + 2 \sqrt{2}$ and the two foliations are induced by straight lines in $\R^2$ in the directions of the eigenvectors 
    \begin{equation*}
        \mathbf{v}_1 = \begin{bmatrix}
            1 + \sqrt{2} \\ 1
        \end{bmatrix} \quad \text{and} \quad 
        \mathbf{v}_2 = \begin{bmatrix}
            1 - \sqrt{2} \\1
        \end{bmatrix}.
    \end{equation*}

    Finally, the action of $\phi$ is depicted in Figure \ref{fig:torus_triple_cover_3}. 
\end{example}

An analogous construction can be carried out on the translation surface obtained by identifying the edges of a regular octagon. This particular construction is taken from \cite{smillie2009symbolic}.

\begin{example}[Pseudo-Anosov map on the regular octagon]\label{ex:pA_regular_octagon}
    Let $X$ be the translation surface obtained by identifying parallel edges in a regular octagon. This surface admits a so-called \emph{horizontal cylinder decomposition} which is indicated in Figure \ref{fig:regular_octagon_cylinder_decomposition}.   
\begin{figure}[ht]
    \centering
    \begin{tikzpicture}
    \draw[thick] (22.5:2) \foreach \x [count=\i] in {22.5,67.5,...,337.5} {
        -- (\x:2) coordinate (p\i)
    } -- cycle;

    \draw[dashed] (p1) -- (p4); \draw[dashed] (p8) -- (p5);

    \node[rotate =0] at ($(p4) - (0.2,0)$) {\Large \Rightscissors};

    \fill[gray, opacity = 0.15] (p1) -- (p2) -- (p3) -- (p4) -- cycle;

    \draw[dashed] (p8) -- ++($(p2) - (p3)$) -- ++($(p1)-(p2)$) -- (p7);
    \fill[gray, opacity =0.05] (p8) -- ++($(p2) - (p3)$) coordinate (p9) -- ++($(p1)-(p2)$) coordinate (p10) -- (p7);

    \draw[->, >=latex] plot [smooth, tension = 1] coordinates {(45:2) (20:2.3) (347.5:2.5)};

    \draw[->, >=latex] (4,0) -- (5,0);

    \def\x{8};

    \draw[thick] ($(p1) + (\x,0)$) -- ($(p4) + (\x,0)$) -- ($(p5) + (\x,0)$) -- ($(p6) + (\x,0)$) -- ($(p7) + (\x,0)$) -- ($(p10) +(\x,0)$) -- ($(p9) +(\x,0)$) -- ($(p8) +(\x,0)$) -- cycle;

    \draw[dashed] ($(p8) +(\x,0)$) -- ($(p5) +(\x,0)$);
    \draw[dotted] ($(p8) +(\x,0)$) -- ($(p7) +(\x,0)$);

    \fill[blue, opacity = 0.10] ($(p5) + (\x,0)$) -- ($(p6) + (\x,0)$) -- ($(p7) + (\x,0)$) -- ($(p10) +(\x,0)$) -- ($(p9) +(\x,0)$) -- ($(p8) +(\x,0)$) -- cycle;
    \fill[purple, opacity =0.10] ($(p1) + (\x,0)$) -- ($(p4) + (\x,0)$) -- ($(p5) + (\x,0)$) -- ($(p8) +(\x,0)$) -- cycle;
    
\end{tikzpicture}
    \caption{The translation surface obtained from a regular octagon, decomposed into horizontal cylinders.}
    \label{fig:regular_octagon_cylinder_decomposition}
\end{figure}
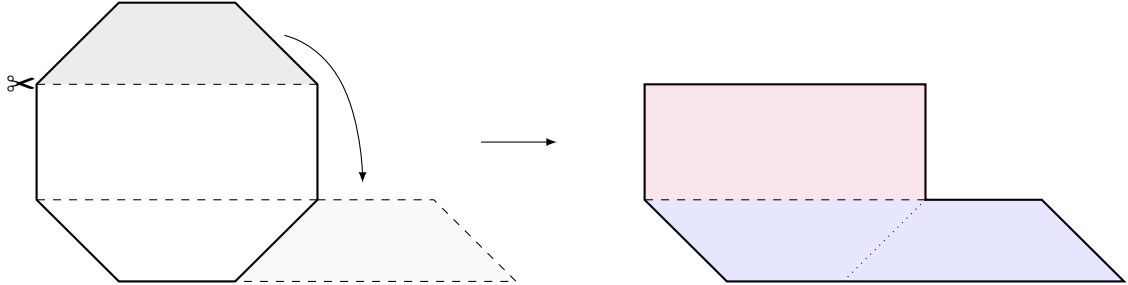 One verifies that the map
    \begin{align*}
        T_h \colon X &\to X, \\
        \begin{bmatrix}
            x \\ y 
        \end{bmatrix} &\mapsto \begin{bmatrix}
            1 & 2+2\sqrt{2} \\
            0 & 1
        \end{bmatrix}\cdot \begin{bmatrix}
            x \\ y
        \end{bmatrix}
    \end{align*}
    is an affine diffeomorphism away from the singular point that fixes the horizontal boundaries of the cylinder decomposition pointwise. The action of $T_h$ is depicted in Figure \ref{fig:regular_octagon_horizontal_shear}. 
    
    \begin{figure}[hb]
    \centering
    \begin{tikzpicture}[scale =0.7]
        \draw[thick] (22.5:2) \foreach \x [count=\i] in {22.5,67.5,...,337.5} {
        -- (\x:2) coordinate (p\i)
    } -- cycle;

    \draw (p7) -- ++($(p1) - (p4)$) coordinate (e1) -- ++ ($(p2) - (p3)$) -- ++($(p3) - (p4)$) coordinate (e2) -- ++($(p1) - (p4)$) coordinate (e3) -- ++($(p1) - (p4)$) -- ++($(p1) - (p8)$) coordinate (e4) -- ++ ($(p2) - (p3)$) coordinate (e5) -- ++($(p1) - (p4)$) -- ++($(p5) - (p6)$) coordinate (s1) -- ++($(p3) - (p2)$) coordinate (e6) -- ++($(p4) - (p1)$) coordinate (e7) -- ++($(p3) - (p2)$) -- ++($(p4) - (p3)$) coordinate (e8) -- ++($(p4) - (p1)$) coordinate (e99) -- ++($(p4) - (p1)$) -- ++($(p8) - (p1)$) coordinate (e9) -- (p8);

    \draw (e1) -- (e9) -- (e2) -- (e99) ;
    \draw (e3) -- (e8) -- (e4) -- (e7);
    \draw (e5) -- (e6);

    \draw[noamblue, dashed] (p6) -- (e9) -- (e8) -- (e6);
    \draw[noamblue, dashed] (p7) -- (e2) -- (e4) -- (s1);

    \draw[] (p5) -- (p8); \draw (p1) -- (p4);
    \draw[noamblue, dashed] ($(p5)!0.5!(p6)$) -- (p8);
    \draw[noamblue, dashed] ($(p1)!0.5!(p8)$) -- (p5);
    \draw[noamblue, dashed] ($(p4)!0.5!(p5)$) -- (p1);
    \draw[noamblue, dashed] ($(p1)!0.5!(p2)$) -- (p4);
    \draw[noamblue, dashed] ($(p3)!0.5!(p4)$) -- (p2);

    \node[shift={(-0.1,0.13)}] at (p7) {\scriptsize 1};
    \node[shift={(0.1,0.25)}] at (p6) {\scriptsize 2};
    \node[shift={(0.45,-0.17)}] at (p5) {\scriptsize 3};
    \node[shift={(-0.3,0.25)}] at (p8) {\scriptsize 4};
    \node[shift = {(0,0)}] at (0,0) {\scriptsize 5};
    \node[shift={(0.3,-0.25)}] at (p4) {\scriptsize 6};
    \node[shift={(-0.45,0.17)}] at (p1) {\scriptsize 7};
    \node[shift={(-0.1,-0.25)}] at (p2) {\scriptsize 8};
    \node[shift={(0.1,-0.13)}] at (p3) {\scriptsize 9};

    \node[shift={(0.1,0.25)}] at (e4) {\scriptsize 2};
    \node[shift={(0.45,-0.17)}] at (e9) {\scriptsize 3};
    \node[shift={(-0.3,0.25)}] at (e2) {\scriptsize 4};
    \node[shift = {(0,0)}] at ($(0,0) + (p1) - (p4) + (p1) - (p4) + (p2) - (p3)$) {\scriptsize 5};
    \node[shift={(0.3,-0.25)}] at (e8) {\scriptsize 6};
    \node[shift={(-0.45,0.17)}] at (e4) {\scriptsize 7};
    \node[shift={(-0.1,-0.25)}] at (e9) {\scriptsize 8};
    \node[shift={(0.1,-0.13)}] at (e6) {\scriptsize 9};
\end{tikzpicture}
    \caption{The map $T_h$ from Example \ref{ex:pA_regular_octagon}.}
    \label{fig:regular_octagon_horizontal_shear}
\end{figure}
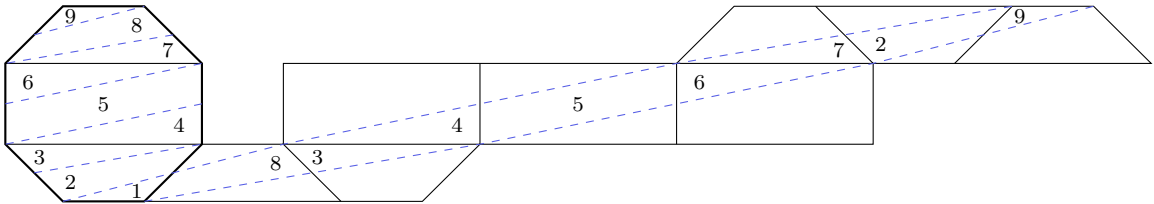

    Using the symmetry of the octagon, we obtain a vertical shear in the same way, i.e., an affine map $T_v$ with differential
    \begin{equation*}
        DT_v = \begin{bmatrix}
            1 & 0 \\
            2 + 2\sqrt{2} & 1
        \end{bmatrix}.
    \end{equation*}
    Defining $\phi = T_v \circ T_h$ yields a pseudo-Anosov map, since using the chain rule we can multiply the two matrices $DT_v$ and $DT_h$ to obtain
    \begin{equation*}
        D\phi = \begin{bmatrix}
            13 + 8\sqrt{2} & 2 + 2 \sqrt{2} \\
            2 + 2\sqrt{2} & 1\\
        \end{bmatrix},
    \end{equation*}
    which satisfies
    \begin{equation*}
        \tr(D\phi) = 14 + 8\sqrt{2} > 2.
    \end{equation*}
    One can again obtain the dilatation of $\phi$ by computing the eigenvalues of $D\phi$. Moreover one verifies that the two eigenvectors $\mathbf{v}_1$ and $\mathbf{v}_2$ are orthogonal, meaning that the associated pair of measured foliations on $X$ are induced by orthogonal lines on the plane $\R^2$.
\end{example}


\subsubsection{Correspondence between Closed Teichmüller Geodesics and Pseudo-Anosov Mapping Classes}

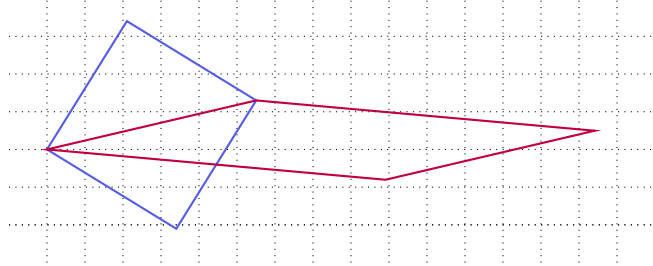
\begin{figure}[ht]
    \centering
    \begin{tikzpicture}
    \begin{scope}
    \foreach \x in {-1, -1,-0.5,0,0.5,...,1.5}{
        \draw[dotted] (-1/2,\x) -- (8,\x);
    }
    \foreach \x in {0,0.5,1,...,7.5}{
        \draw[dotted] (\x, -1.5) -- (\x, 2);
    }
    \end{scope}

    \begin{scope}[rotate = -31.71722467]
        \draw[thick, noamblue] (0,0) rectangle (2,2); 
        \draw[thick, purple] (0,0) -- (4,2) -- (6,4) -- (2,2) -- cycle;
    \end{scope}
\end{tikzpicture}
    \caption{The CAT map acting on a (rotated) torus, which can be seen to be the same as a closed Teichmüller geodesic.}
    \label{fig:cat_map_teichmueller}
\end{figure}

The (constant) differentials of the explicitly constructed pseudo-Anosov maps from Examples \ref{ex:torus_triple_cover} and \ref{ex:pA_regular_octagon} both have orthogonal eigenvectors and by definition, the corresponding pseudo-Anosov map $\phi$ expands in one direction and contracts in the other direction. From these examples it may be already clear, that pseudo-Anosov maps and closed Teichmüller geodesics must be related. This becomes especially transparent if we go back to the CAT map seen above. Rotating the torus such that the two orthogonal eigendirections become horizontal and vertical, the action of the CAT map can now be represented as in Figure \ref{fig:cat_map_teichmueller}.

We can explicitly give the (minimal) period of the corresponding closed Teichmüller geodesic. Recall that the dilataion of the CAT map $\phi \colon \T^2 \to \T^2$ is given by $\lambda = \frac{1}{2}(3 + \sqrt{5})$. If we write $T$ for the minimal period, this imposes
\begin{equation*}
    \e^T = \lambda,
\end{equation*}
since after the rotation the eigenvalue associted to the dilatation is given by the first standard basis vector. Therefore, the CAT map $\phi$ corresponds exactly to a closed Teichmüller geodesic in direction of the eigenvector associated to the dilatation with minimal period
\begin{equation*}
    T = \log \lambda.
\end{equation*}
Conversely, if we are given some $T > 0$ such that
\begin{equation*}
    g_T \cdot \T^2 = \T^2,
\end{equation*}
i.e., $T$ is a period of the Teichmüller geodesic flow $g_t$ on $\T^2$, it is also true that the discrete map
\begin{align*}
    \phi \colon \T^2 &\to \T^2, \\
    \begin{bmatrix}
        x \\
        y
    \end{bmatrix} & \mapsto g_T \cdot \begin{bmatrix}
        x \\
        y
    \end{bmatrix} = \begin{bmatrix}
        \e^T x \\
        \e^{-T}y
    \end{bmatrix}
\end{align*}
defines a pseudo-Anosov diffeomorphism according to the definitions above. Indeed, the pair of transverse measured foliations is simply the pullback of the horizontal and vertical directions and since the matrix $g_T$ is diagonal, it is easy to read off that the eigenvalues are given by $\e^T$ and $\e^{-T}$, as well as $\tr(D\phi) = \e^T + \e^{-T} > 2$. 

Working with the general definition of pseudo-Anosov map using just a pair of transverse measured foliation we gave in the beginning, it may not be clear that this correspondence carries over to the general case. However, by obtaining the equivalent definition making use of a natural atlas, i.e., Definition \ref{def:pseudo_Anosov_equivalent_definition}, the proof of the following theorem becomes very simple.

\begin{theorem}
    Let $X$ be a compact connected orientable surface. There is a one-to-one correspondence between pseudo-Anosov mapping classes $[f]\in \operatorname{Mod}(X)$ corresponding to pairs of transverse measured foliations admitting only even-pronged singularities and closed Teichmüller geodesics through $X$. 
    
    More precisely, given any pseudo-Anosov mapping class $[f]$, we can equip $X$ with a translation structure and associate a minimal period $T = T([f])$ of the Teichmüller geodesic flow to the mapping class, meaning that
    \begin{equation*}
        g_T \cdot X = X.
    \end{equation*}
    Conversely, given any closed Teichmüller geodesic on a translation surface $X$ with minimal period $T$, the map
    \begin{align*}
        \phi \colon X &\to X, \\
        \begin{bmatrix}
            x \\
            y
        \end{bmatrix} &\mapsto g_T \cdot \begin{bmatrix}
            x \\
            y
        \end{bmatrix} = \begin{bmatrix}
            \e^T x \\
            \e^{-T}y
        \end{bmatrix}
    \end{align*}
    defined in charts is a pseudo-Anosov map, where the associated pair of transverse measured foliations is given by the pullback of the horizontal and vertical foliation of $\R^2$ and the dilatation is given by $\lambda = \e^T$. 

    Moreover, the correspondence just described is bijective. 
\end{theorem}

\begin{proof}
    Given any pseudo-Anosov mapping class $[f]$ with only even-pronged singularities, we may choose a representative homeomorphism $\phi$ and a natural atlas that endows $X$ with a translation structure as described in Definition \ref{def:pseudo_Anosov_equivalent_definition}. In particular, we may assume that $X$ is given by a collection of polygons $\mathcal{P}$ in the plane with sides identified by translations and the associated pair of transverse measured foliations is given by the pullback of the horizontal and vertical lines in $\R^2$. Defining $T = \log \lambda$, where $\lambda$ is the dilatation of $\phi$, it follows immediately that
    \begin{align*}
        g_T \cdot X &= \left\{\begin{bmatrix}
            \e^T x \\
            \e^{-T}y
        \end{bmatrix}
        \mid \begin{bmatrix}
            x \\
            y
        \end{bmatrix} \in \mathcal{P}
        \right\} /_\sim \\&= 
        \left\{\begin{bmatrix}
            \lambda x \\
            \lambda^{-1}y
        \end{bmatrix} \mid \begin{bmatrix}
            x \\
            y
        \end{bmatrix} \in \mathcal{P}\right\} /_\sim \\&=
        \left\{\phi\left(\begin{bmatrix}
            x \\ y
        \end{bmatrix}\right)\mid \begin{bmatrix}
            x \\
            y
        \end{bmatrix} \in \mathcal{P}\right\}/_\sim \\&=
        X,
    \end{align*}
    where we denote by $\sim$ the equivalence relation induced by gluing the edges of the polygons and the last equality follows from the fact that $\phi$ is invertible. 

    Conversely, given a closed Teichmüller geodesic through a translation surface $X$ with minimal period $T$, the map
    \begin{align*}
        \phi \colon X &\to X, \\
        \mathbf{x} &\mapsto g_T \cdot \mathbf{x}
    \end{align*}
    is well-defined. Moreover, we have that $\phi(\Sigma) = \Sigma$ since closed Teichmüller geodesics must preserve conical singularities, and defining $\lambda = \e^T$ we have that the map is given by
    \begin{equation*}
        \phi(y) = \begin{bmatrix}
            \lambda & 0 \\
            0 & \lambda^{-1}
        \end{bmatrix} y,
    \end{equation*}
    i.e., it is pseudo-Anosov according to Definition \ref{def:pseudo_Anosov_equivalent_definition} with dilatation $\lambda = \e^T > 1$. It is clear that the correspondence outlined above is bijective, which is essentially just a consequence of the invertibility of the exponential map. 
\end{proof}

\subsection{Counting Pseudo-Anosovs using Diagonal Changes}

The algorithm introduced in section \ref{sec:diagonal_changes} is in particular well-suited to \emph{enumerate} closed geodesics of the Teichmüller geodesic flow $g_t$, a problem very much related to counting closed geodesics. An enumeration of closed geodesics consists of a complete list of closed geodesics, ordered for example by the length of the geodesic. As we have explained above, this is equivalent to a complete list of pseudo-Anosov mapping classes on a surface ordered by increasing dilatation.

There have been attempts to use also Rauzy--Veech induction to obtain counting results for closed Teichmüller geodesics, see in particular \cite{bufetov2005logarithmic}. It is possible to obtain counts on closed geodesics \emph{in the space of zippered recangles}, which can be seen as a finite-to-one cover of the corresponding stratum, however due to the complicated nature of this cover it it appears unfeasible to infer results about closed geodesics in the stratum from these results.

Here, we will explain how diagonal changes can be used in order to enumerate pseudo-Anosov mapping classes. The method we present may be used in all hyperelliptic components $\hyp{k}$, but for technical reasons the method becomes more involved when $k$ is even. For this reason, we will restrict ourselves to the case where $k$ is odd. We refer the reader to \cite{delecroixulcigrai20++enumerating} for further details on the general case.

The main idea is to associate a pseudo-Anosov mapping class to a loop in a graph induced by staircase moves on combinatorial data similar to a $DC$ graph from Definition \ref{def:diagonal_change_classes}. Thus we need an appropriate definition of graph, as well as the correct definition of loop, where we will need to introduce several equivalence classes in order to obtain a one-to-one correspondence between loops and pseudo-Anosov mapping classes. 

\subsubsection{The Unlabeled \texorpdfstring{$DC$}{DC} Graph}
We will start to define an unlabeled counterpart to the $DC$ graph defined above. It will be convenient to fix once and for all a (specific) labeling of bundles, or equivalently a labeling for quadrilaterals, which encodes the action of the hyperelliptic involution $s$. 

\begin{definition}[Cyclical labeling respecting the involution]\label{def:good_labeling}
    We will say that a labeling of quadrilaterals $q_i$ in a Quadrangulation $Q$ is a \emph{cyclical labeling respecting the involution}, if it is constructed as follows.
    \begin{enumerate}
        \item Pick arbitrarily the first bundle $\Gamma_i$, or equivalently the first quadrilateral $q_i$.
        \item Define $\Gamma_{i+1}$ (or $q_{i+1}$) to be the bundle (or quadrilateral) containing the image of the vertical ray in $\Gamma_i$ by $s \circ \rho$, where $s$ denotes the hyperelliptic involution and $\rho$ rotates the ray by $\pi$ around the singularity.
    \end{enumerate}
\end{definition}

We illustrate this idea by an example.

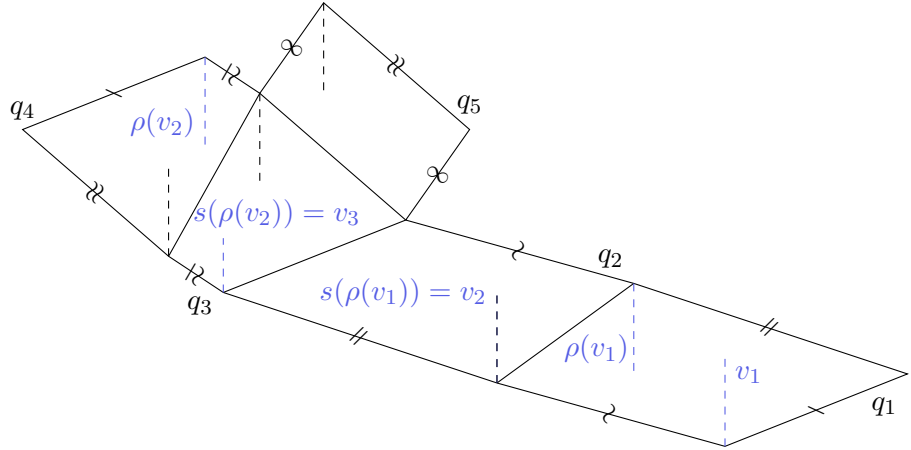
\begin{figure}[ht]
    \centering
    \begin{tikzpicture}[scale = 1.2]
    \begin{scope}[rotate = 90]
    \coordinate (p1) at (0,0);
    \coordinate (p2) at (-0.8,2);
    \coordinate (p3) at (1,3);
    \coordinate (p4) at (-0.1,4.5);
    \coordinate (p5) at (1.7,5.5);
    \coordinate (p6) at (0.9, 7.5);
    \coordinate (p7) at (3.1, 7.1);
    \coordinate (p8) at (1.3,8.1);
    \coordinate (p9) at (3.5, 7.7);
    \coordinate (p10) at (2.7, 9.7);
    \coordinate (p11) at (2.7, 4.8);
    \coordinate (p12) at (4.1, 6.4);   
    \end{scope}

    \draw (p1) -- (p2) -- (p4) -- (p6) -- (p8) -- (p10) --  (p9) -- (p7) -- (p12) -- (p11) -- (p5) -- (p3) -- cycle;
    \draw (p3) -- (p4); \draw (p5) -- (p6); \draw (p7) -- (p8); \draw (p5) -- (p7);

    \node[rotate = 45] at ($(p1)!0.5!(p3)$) {=};
    \node[rotate = 45] at ($(p4)!0.5!(p6)$) {=};

    \node[rotate = 60] at ($(p3)!0.5!(p5)$) {$\sim$};
    \node[rotate = 60] at ($(p2)!0.5!(p4)$) {$\sim$};

    \node[rotate = -25] at ($(p1)!0.5!(p2)$) {$-$};
    \node[rotate = -25] at ($(p9)!0.5!(p10)$) {$-$};

    \node[rotate = 60] at ($(p11)!0.5!(p12)$) {$\approx$};
    \node[rotate = 60] at ($(p10)!0.5!(p8)$) {$\approx$};

    \node[rotate = 60] at ($(p7)!0.5!(p9)$) {$\eqsim$};
    \node[rotate = 60] at ($(p8)!0.5!(p6)$) {$\eqsim$};

    \node[rotate = -30] at ($(p12)!0.5!(p7)$) {\footnotesize $\infty$};
    \node[rotate = -30] at ($(p11)!0.5!(p5)$) {\footnotesize $\infty$};

    \node[below left, yshift = -5] at (p1) {$q_1$};
    \draw[noamblue, dashed] (p2) -- ++(90:1) node[below right] {$v_1$};
    \draw[noamblue, dashed] (p3) -- ++(-90:1) node[above left, xshift= 2] {$\rho(v_1)$};

    \draw[noamblue, dashed] (p4) -- ++(90:1) node[left] {$s(\rho(v_1)) = v_2$};
    \node[above left, yshift = 2] at (p3) {$q_2$};

    \draw[noamblue, dashed] (p9) -- ++(-90:1) node[above left] {$\rho(v_2)$};
    \draw[noamblue, dashed] (p6) -- ++(90:0.6) node [above, xshift = 20] {$s(\rho(v_2)) = v_3$};
    \node[left, yshift= -5] at (p6) {$q_3$};

    \draw[dashed] (p7) -- ++(-90:1);
    \draw[dashed] (p8) -- ++(90:1);
    \draw[dashed] (p4) -- ++(90:1);
    \draw[dashed] (p12) -- ++(-90:1);

    \node[above] at (p10) {$q_4$};
    \node[above, yshift = 2] at (p11) {$q_5$};
\end{tikzpicture}
    \caption{A cyclical labeling respecting the involution of a quadrangulation $Q$.}
    \label{fig:good_labeling}
\end{figure}

\begin{example}[Cyclical labeling respecting the involution]\label{ex:good_labeling}
    Consider the surface $X$ given by the quadrangulation $Q$ in Figure \ref{fig:good_labeling}. Say we pick the bottom right quadrilateral to be $q_1$ and let us denote by $v_1$ the corresponding vertical ray. As we can see in the Figure, the rotatet ray $\rho(v_1)$ belongs to the same quadrilateral, which under the hyperelliptic involution $s$ is mapped to the top left neighbor of $q_1$, which is therefore labeled $q_2$. Continuing this way we obtain the same labeling as for the quadrangulation in Figure \ref{fig:quadrangulation_triangulation_dual_graph}, where the quadrangulation is just rotated by $-\pi$.
\end{example}

\begin{lemma}
    If $k = 2n+1$ is odd, going around the unique singularity of a surface $X \in \hyp{k}$ we see the labels of a labeling respecting the involution in the order $1,3,5, \ldots, k,2,4,\ldots, k-1$.
\end{lemma}
In Example \ref{ex:good_labeling}, it is easy to see that we see the labels in the order $1,3,5,2,4$ when going around the unique singularity.

Suppose now we have a $DC$ graph $\mathcal{G} = \mathcal{G}(k)$ endowed with a labeling according to Definition \ref{def:good_labeling}. We next define the associated \emph{unlabeled} $DC$ graph $\overline{\mathcal{G}}$.

\begin{figure}[ht]
    \centering
\[\begin{tikzcd}
	{\begin{array}{c} (A \,B)(C)\\(A \, C \,B) \end{array}}\arrow[loop above, distance = 1cm, in=120, out=60, "\cdot\,\ell\,\cdot" {pos=0.7}]\arrow[loop, distance = 1cm, in=300, out=240, "r\,r\,r" {pos=0.3}] &&& {\begin{array}{c} (A\,C)(B)\\(A\,B)(C) \end{array}}\arrow[loop above, distance = 1cm, in=120, out=60, "\cdot\,\ell\,\cdot" {pos=0.7}]\arrow[loop, distance = 1cm, in=300, out=240, "\cdot\,\cdot\,r" {pos=0.3}] &&& {\begin{array}{c}(A\,B\,C)\\(A\,B)(C) \end{array}}\arrow[loop above, distance = 1cm, in=120, out=60, "\cdot\,\cdot\,r" {pos=0.7}]\arrow[loop, distance = 1cm, in=300, out=240, "\ell\,\ell\,\ell" {pos=0.3}]
	\arrow["\ell\cdot\ell"', curve={height=-6pt}, from=1-1, to=1-4]
	\arrow[curve={height=-6pt}, from=1-4, to=1-1]
	\arrow["rr\cdot", curve={height=6pt}, from=1-7, to=1-4]
	\arrow[curve={height=6pt}, from=1-4, to=1-7]
\end{tikzcd}\]
    \caption{The unlabeled $DC$ graph $\overline{\mathcal{G}}$ associated to the $DC$ graph $\mathcal{G}$ from Figure \ref{fig:dc_graph}}
    \label{fig:unlabeled_graph}
\end{figure}
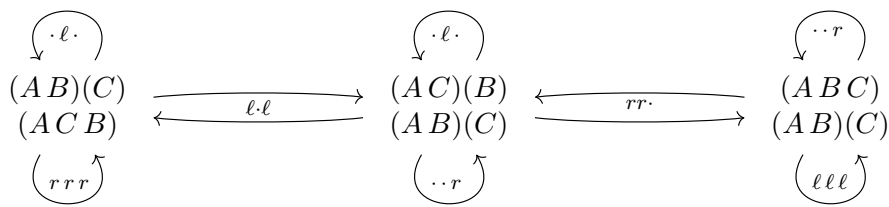

\begin{definition}[Unlabeled $DC$ graph]\label{def:unlabeled_DC_graph}
    We define the \emph{unlabeled $DC$ graph} $\overline{\mathcal{G}}$ associated to the $DC$ graph $\mathcal{G}$ as the quotient of $\mathcal{G}$ by the action of relabeling. More precisely, we define an equivalence relation on the vertices by calling $\boldsymbol{\pi} = (\pi_\ell, \pi_r)$ and $\boldsymbol{\pi}' = (\pi'_\ell, \pi'_r)$ \emph{relabeling equivalent}, if there exists a permutation $\sigma$ of $[k]$ such that 
    \begin{equation}\label{eq:relabeling_equivalence}
        (\pi'_\ell, \pi'_r) = (\sigma\pi_\ell\sigma^{-1}, \sigma\pi_r\sigma^{-1}).
    \end{equation}
    For brevity, we will denote \eqref{eq:relabeling_equivalence} by $\sigma \star \boldsymbol{\pi} = \boldsymbol{\pi}'$. This definition extends to the edges, where $\sigma \star c$ is defined in the obvious way for a (left of right) cycle $c$ representing the edge of the graph.

    Then, we define the unlabeled $DC$ graph $\overline{\mathcal{G}}$ to be the quotient
    \begin{equation*}
        \overline{\mathcal{G}} = \mathcal{G} /_\sim,
    \end{equation*}
    where $\sim$ denotes the equivalence relation we just described. If $\boldsymbol{\pi}_0$ is a vertex in $\mathcal{G}$, we write $[\boldsymbol{\pi}_0]$ for its equivalence class in $\overline{\mathcal{G}}$, i.e., 
    \begin{equation*}
        [\boldsymbol{\pi}_0] = \{\boldsymbol{\pi} \mid \text{ there exists } \sigma \in S_k: \sigma \star \boldsymbol{\pi} = \boldsymbol{\pi}_0\}.
    \end{equation*}
    There exists an oriented edge (or an arrow) $\overline{\gamma}$ in $\overline{\mathcal{G}}$ from $[\boldsymbol{\pi}_1]$ to $[\boldsymbol{\pi_2}]$ if and only if there exists an arrow $\gamma$ in $\mathcal{G}$ from some $\boldsymbol{\pi}_1 \in [\boldsymbol{\pi}_1]$ to some $\boldsymbol{\pi}'_2 = c\cdot \boldsymbol{\pi}_1 \in [\boldsymbol{\pi}_2]$, where $c$ is the cycle represented by the arrow $\gamma$. We will write $[\gamma]$ for the equivalence class of $\gamma$.
\end{definition}
Figure \ref{fig:unlabeled_graph} shows the unlabeled graph $\overline{\mathcal{G}}$ associated to the $DC$ graph $\mathcal{G}$ in Figure \ref{fig:dc_graph}. Note that $\mathcal{G}$ is a finite-to-one cover of $\overline{\mathcal{G}}$, in this case the cover is 3-to-1.

There is an obvious projection map $P \colon \mathcal{G} \to \overline{\mathcal{G}}$, which maps a combinatorial datum or a cycle to its representative. The proof of the following result can be found in the appendix of \cite{delecroixulcigrai20++enumerating}.

\begin{lemma}[Projection has no fixed points if $k$ is odd]
    When $k$ is odd, the projection $P \colon \mathcal{G} \to \overline{\mathcal{G}}$ has no fixed point, hence $\mathcal{G}$ is a \emph{regular cover} of $\overline{\mathcal{G}}$ and $P$ is a \emph{covering map}.
\end{lemma}
\begin{remark}
    This is no longer true if $k$ is even, in which case $\mathcal{G}$ is a \emph{branched} cover of $\overline{\mathcal{G}}$. This is exactly the reason why the treatment of the case when $k$ is odd is simpler. 
\end{remark}

\subsubsection{Enumeration in \texorpdfstring{$\hyp{k}$}{C hyp(k)} when \texorpdfstring{$k$}{k} is odd}
We have given the relevant definition of graph we will need for the enumeration, next we need to define the appropriate notion of loop.

A consequence of the fact that the covering map $P \colon \mathcal{G} \to \overline{\mathcal{G}}$ has no fixed points is, that for every edge $\overline{\gamma}$ in $\overline{\mathcal{G}}$ from $[\boldsymbol{\pi}_0]$ to $[\boldsymbol{\pi}'_0]$ and for every vertex $\boldsymbol{\pi}_1 \in [\boldsymbol{\pi}_0]$, there exists a \emph{unique} edge $\gamma$ in $\mathcal{G}$ starting at $\boldsymbol{\pi}_1$ such that $P(\gamma) = \overline{\gamma}$. If $\gamma$ represents the cycle $c_1$, then in particular
\begin{equation*}
    \boldsymbol{\pi}'_1 = c_1 \cdot \boldsymbol{\pi}_1 \in [\boldsymbol{\pi}_0].
\end{equation*}
This allows us to make the following definitions.

\begin{definition}[Path in $\overline{\mathcal{G}}$]
    A \emph{path} $\overline{\gamma}$ in $\overline{\mathcal{G}}$ is a sequence of arrows $\overline{\gamma}_i$ in $\overline{\mathcal{G}}$ for $0 \leq i \leq n-1$, that can be concatenated. We will write that $\overline{\gamma}_i$ is an arrow from $[\boldsymbol{\pi}_i]$ to $[\boldsymbol{\pi}_{i+1}]$ and furthermore we will use the notation 
    \begin{equation*}
        \overline{\gamma} = \overline{\gamma}_1 \cdots \overline{\gamma}_{n-1}
    \end{equation*}
    for such a path. 

    Further, we say that $\overline{\gamma}$ gives a \emph{presentation of a loop}, if $[\boldsymbol{\pi}_0] = [\boldsymbol{\pi}_{n-1}]$. We say that two presentations of a loop $\overline{\gamma}$ and $\overline{\gamma}'$ are equivalent, or more precisely equivalent up to cyclical reordering, if $\overline{\gamma}'$ is obtained by \emph{cyclical reordering} the arrows of $\overline{\gamma}$, i.e, if $\overline{\gamma} = \overline{\gamma}_0 \cdots \overline{\gamma}_{n-1}$, then $\overline{\gamma}'$ can be written as
    \begin{equation*}
        \overline{\gamma}' = \overline{\beta} \overline{\alpha},
    \end{equation*}
    where $\overline{\alpha} = \overline{\gamma}_0 \cdots \overline{\gamma}_i$ and $\overline{\beta} = \overline{\gamma}_{i+1} \cdots \overline{\gamma}_{n-1}$ for some $i \in [n-1]$.

    A \emph{loop} is then an equivalence class of presentations of loops. We will denote by $\mathcal{Y}_k$ the set of loops in $\overline{\mathcal{G}}$ and abusing the notation slightly, we will write $\overline{\gamma} \in \mathcal{Y}_k$ denoting both a given presentation of a loop $\overline{\gamma}$ and its equivalence class.
\end{definition}

Since $P$ is a covering map for $k$ odd, one can uniquely lift paths to $\mathcal{G}$.

\begin{lemma}[Path lifts]\label{lem:path_lifts}
    Let $k$ be odd. The projection $P \colon \mathcal{G} \to \overline{\mathcal{G}}$ satisfies the \emph{lifting property:} given a path $\overline{\gamma}$ on $\overline{\mathcal{G}}$ and a vertex $\boldsymbol{\pi}'_0\in \mathcal{G}$ which lifts the starting point $[\boldsymbol{\pi}_0]$ of $\overline{\gamma}$, i.e., is such that $P(\boldsymbol{\pi}'_0) = [\boldsymbol{\pi}_0]$, then there exists a \emph{unique} path $\gamma$ starting from $\boldsymbol{\pi}'_0$ such that $P(\gamma) = \overline{\gamma}$, i.e., a unique path $\gamma$ that lifts $\overline{\gamma}$ starting from $\boldsymbol{\pi}'_0$.
\end{lemma}
The proof of this lemma can be found in standard textbooks on covering theory, see e.g., \cite{munkres2019topology}.

We say that the path $\gamma$ in Lemma \ref{lem:path_lifts} is the lift of $\overline{\gamma}$ starting at $\boldsymbol{\pi}$. If $\overline{\gamma}$ is a loop starting and ending at $[\boldsymbol{\pi}]$ and $\gamma$ is a lift of $\boldsymbol{\gamma}$ starting from a lift $\boldsymbol{\pi}_0$ of $[\boldsymbol{\pi}]$, the ending vertex $\boldsymbol{\pi}_{n-1}$ of $\gamma$ must belong to the equivalence class $[\boldsymbol{\pi}]$ as well, so there exists a relabeling $\sigma$, trivial only if $\gamma$ is itself a (proper) loop in $\mathcal{G}$, such that
\begin{equation*}
    \sigma \star \boldsymbol{\pi}_{n-1} = \boldsymbol{\pi}_0.
\end{equation*}
We will say that $\sigma$ \emph{closes the loop $\gamma$} and that the pair $(\gamma, \sigma)$ is a \emph{lifted loop} that \emph{lifts the loop presentation $\overline{\gamma}$ in $\overline{\mathcal{G}}$}, starting at the vertex $\boldsymbol{\pi}_0$. More generally, if $\gamma$ is a path from $\boldsymbol{\pi}_0$ to $\boldsymbol{\pi}_{n-1}$ and $\sigma$ a relabeling such that $\sigma \star \boldsymbol{\pi}_{n-1} = \boldsymbol{\pi}_0$, then $P(\gamma)$ is a presentation of a loop in $\overline{\mathcal{G}}$. This motivates the following definition.

\begin{definition}[Lifted loops]\label{def:lifted_loops}
    We call \emph{lifted loop} a pair $(\gamma, \sigma)$, where $\gamma$ is a path in $\mathcal{G}$ from $\boldsymbol{\pi}_0$ to $\boldsymbol{\pi}_{n-1}$ and $\sigma \in S_n$ is a relabeling such that
    \begin{equation*}
        \sigma \star \boldsymbol{\pi}_{n-1} = \boldsymbol{\pi}_0.
    \end{equation*}
    We will say that the lifted loop $(\gamma, \sigma)$ \emph{lifts the loop} $\overline{\gamma}'$, if it lifts a loop presentation in the equivalence class of the loop presentation $\overline{\gamma}'$, that is, if there exists a cyclical reordering $\overline{\gamma}$ of $\overline{\gamma}'$ such that $P(\gamma) = \overline{\gamma}$. 
\end{definition}

These are \emph{almost} the loops we want to associate to a pseudo-Anosov mapping class. The last notion we need to introduce is that of a \emph{positive} loop. For this definition, we will work with the matrices $A_{\boldsymbol{\pi},c}$ from Definition \ref{def:staircase_matrix}. We will use the following notation for a \emph{sequence} of staircase moves. For a path $\gamma$ in $\mathcal{G}$, which is completely determined by an initial permutation datum $\boldsymbol{\pi}_0$ and a sequence of cycles $c_0, \ldots, c_{n-1}$, we can describe the linear action corresponding to the sequence of staircase moves explicitly by
\begin{equation}\label{eq:path_matrix_1}
    A_\gamma \coloneqq A_{\boldsymbol{\pi}_{n-1}, c_{n-1}} \cdots A_{\boldsymbol{\pi}_1, c_1} \cdot A_{\boldsymbol{\pi}_0, c_0},
\end{equation}
where $\boldsymbol{\pi}_i = c_i \cdot c_{i-1} \cdots c_0 \cdot \boldsymbol{\pi}_0$. In other words, the linear action of the sequence of staircase moves is given simply by the matrix product of the individual moves. 

Let $\overline{\gamma} \in \mathcal{Y}_k$ be a loop in $\overline{\mathcal{G}}$ and let $(\gamma,\sigma)$ be any lifted loop of $\overline{\gamma}$. To the lifted loop $(\gamma, \sigma)$, we now associate the matrix $\Pi_\sigma A_\gamma$, where $\Pi_\sigma$ is the permutation matrix associated to $\sigma$. Choosing a different lifted loop lifting $\overline{\gamma}$ gives a matrix that is cojugated to $\Pi_\sigma A_\gamma$ by a permutation matrix, which follows from the definition of the unlabeled $DC$ graph (Definition \ref{def:unlabeled_DC_graph}), so that we can give the following definitions.

\begin{definition}[Positive loops]\label{def:positive_loops}
    A loop $\overline{\gamma}$ in $\overline{\mathcal{G}}$ is \emph{positive}, if for any lifted loop $(\gamma, \sigma)$ in $\mathcal{G}$ of $\overline{\gamma}$ the associated matrix $\Pi_\sigma A_\gamma$ is \emph{primitive}, i.e., there exists a positive integer $n$ such that $(\Pi_\sigma A_\gamma)^n$ has strictly positive entries.
\end{definition}

\begin{definition}[Set of positive loops]
    We denote by $\mathcal{Y}_k^+ \subseteq \mathcal{Y}_k$ the \emph{space of positive loops} of $\overline{\mathcal{G}}$.
\end{definition}

We will explain below how we can associate to each \emph{positive} loop a pseudo-Anosov map of a surface $X \in \hyp{k}$. It is true that every pseudo-Anosov mapping class can be represented by a positive loop, but this correspondence is not injective, i.e., two different loops can give rise to the same (conjugacy class of) pseudo-Anosov map. To obtain a one-to-one correspondence, we have to consider another equivalence relation on loops in $\mathcal{Y}_k^+$, which fixes the issue that some staircase moves commute with each other, so that it may be the case that applying the moves in a different order creates different loops, but on the level of staircase moves the order is irrelevant. Let us now describe this equivalence relation.

Suppose that $S_c$ and $S_{c'}$ are \emph{disjoint} staircases, that is $c$ and $c'$ are cycles of $\pi_r$ or $\pi_\ell$ which act on disjoint subsets of the indices $[k]$. Then, the staircase moves in $X$ commute, that is, if $S_c$ and $S_{c'}$ are both well-slanted staircases for the quadrangulation $Q = (\boldsymbol{\pi}, \boldsymbol{\lambda}, \boldsymbol{\tau})$, then $c'$ is a cycle of $c\cdot \boldsymbol{\pi}$ and $c$ is a cycle of $c'\cdot \boldsymbol{\pi}$ and we have
\begin{equation*}
    \hat{m}_{c\cdot \boldsymbol{\pi}, c'} \circ \hat{m}_{\boldsymbol{\pi, c}}(\boldsymbol{\pi}, \boldsymbol{\lambda}, \boldsymbol{\tau}) = \hat{m}_{c'\cdot\boldsymbol{\pi},c} \circ \hat{m}_{\boldsymbol{\pi},c'} (\boldsymbol{\pi}, \boldsymbol{\lambda}, \boldsymbol{\tau}).
\end{equation*}
This gives a way of describing the needed conjugacy classes. We will start to define the equivalence relation on paths.

\begin{definition}[Equivalence of paths up to commutation rules]\label{def:equivalent_paths_up_to_commutation}
    We say that two paths $\gamma_1$ and $\gamma_2$  in $\mathcal{G}$ are equivalent up to commutation rules, and we write $\gamma_1 \sim \gamma_2$, if they start at the same combinatorial datum and one path can be obtained from the other by permuting arrows which correspond to staircase moves disjoint staircases.
\end{definition}

\begin{figure}[ht]
    \centering
\[\begin{tikzcd}
	& {\begin{array}{c} (3 \,1)(2)\\(1 \, 3 \,2) \end{array}}\arrow[loop above, distance = 1cm, in=120, out=60, "\cdot\,\ell\,\cdot" {pos=0.7}, color = MidnightBlue] &&& {\begin{array}{c} (3\,1)(2)\\(1\,2)(3) \end{array}}\arrow[loop above, distance = 1cm, in=120, out=60, "\cdot\,\ell\,\cdot" {pos=0.7}, color = Plum]\arrow[loop, distance = 1cm, in=300, out=240, "\cdot\,\cdot\,r" {pos=0.3}, color = ForestGreen] &&& {\begin{array}{c}(1\,2\,3)\\(1\,2)(3) \end{array}}\arrow[loop above, distance = 1cm, in=120, out=60, "\cdot\,\cdot\,r" {pos=0.7}] \\
	\\
	{} && {} &&&& {} && {}
	\arrow["\ell\cdot\ell"', curve={height=-6pt}, from=1-2, to=1-5]
	\arrow[color={rgb,255:red,214;green,92;blue,92}, curve={height=-6pt}, from=1-5, to=1-2]
	\arrow["rr\cdot", curve={height=6pt}, from=1-8, to=1-5]
	\arrow[curve={height=6pt}, from=1-5, to=1-8]
	\arrow[dashed, no head, from=1-8, to=3-7]
	\arrow[dashed, no head, from=1-2, to=3-3]
	\arrow[dashed, from=3-1, to=1-2]
	\arrow[dashed, from=3-9, to=1-8]
\end{tikzcd}\]
    \caption{Commuting staircase moves in (part of) a $DC$ graph.}
    \label{fig:disjoint_staircases}
\end{figure}
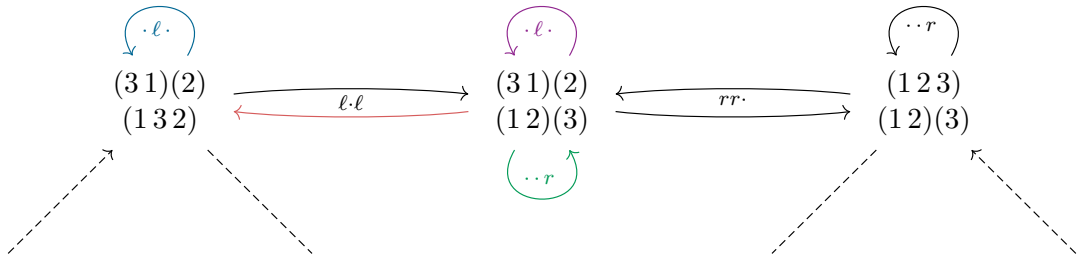

So if $\gamma_1 \sim \gamma_2$ and $Q = (\boldsymbol{\pi}, \boldsymbol{w})$ is a quadrangulation with $\boldsymbol{\pi}$ a common starting vertex of both $\gamma_1$ and $\gamma_2$, then
\begin{equation}\label{eq:commutation_equality}
    Q' = \hat{m}_{\gamma_1}(Q) = \hat{m}_{\gamma_2}(Q),
\end{equation}
meaning that the quadrangulation obtained by following either path is the same. 
\begin{remark}
    In fact, the converse is also true, i.e., if \eqref{eq:commutation_equality} holds, then we necessarily have $\gamma_1 \sim \gamma_2$, see Lemma 6.8 in \cite{delecroixulcigrai20++enumerating}.
\end{remark}
It is easy to check whether two paths $\gamma_1$ and $\gamma_2$ are the same up to commutation rules by using the graph $\mathcal{G}$, since two staircase moves in disjoint staircases are simply discribed by arrows labeled with disjoint cycles. This can for example be seen in Figure \ref{fig:disjoint_staircases}, which is the top row of the $DC$ graph from Figure \ref{fig:dc_graph}. We can see two pairs of commuting staircase moves. First, the move represented by the purple arrow $\cdot \, \ell \,\cdot$ commutes with the move represented by the green arrow $\cdot \, \cdot \, r$. The notation reflects the fact that the two moves are disjoint, as in each position (1, 2 or 3) there is at most one symbol $r$ or $\ell$ in both labels together. Another pair of commuting moves can be seen as follows. Applying first the move represented by the purple arrow $\cdot \, \ell\, \cdot$ and then the move represented by the red arrow $\ell\,\cdot\,\ell$ is the same as doing first the move represented by the red arrow and then the move repressented by the blue arrow $\cdot\,\ell\,\cdot$, since again there is no overlap of the symbols. 

The equivalence relation on paths then induces an equivalence relation on loops.

\begin{definition}[Equivalence of loops up to commutation rules]\label{def:equivalent_loops_up_to_commutation}
    Two \emph{lifted loops} $(\gamma_1, \sigma_1)$ and $(\gamma_2, \sigma_2)$ are \emph{equivalent up to commutation rules}, if $\gamma_1 \sim \gamma_2$ in the sense of Definition \ref{def:equivalent_paths_up_to_commutation} and $\sigma_1 = \sigma_2$. 
    
    We say that two loops $\overline{\gamma}_1$ and $\overline{\gamma}_2$ in $\mathcal{Y}_k^+$ are \emph{equivalent}, and we write $[\overline{\gamma}_1] = [\overline{\gamma}_2]$, if there exist lifted loops $(\gamma'_1, \sigma_1')$ and $(\gamma_2', \sigma_2')$ lifting $\overline{\gamma}_2$ and $\overline{\gamma}_2$ respectively, such that $(\gamma'_1, \sigma_1')$ and $(\gamma_2', \sigma_2')$ are equivalent up to commutation rules. 
\end{definition}

Having defined the correct notion of graph and loop, we can now state the main result on enumeration of pseudo-Anosov diffeomorphisms.

\begin{theorem}[Pseudo-Anosov diffeomorphism enumeration for $k$ odd]
    Let $\hyp{k}$ be a hyperelliptic component of a stratum of translation surfaces, where $k \geq 3$ is odd. There is a one-to-one correspondence between pseudo-Anosov mapping classes $[\phi] \in \operatorname{Mod}(X)$, where $X \in \hyp{k}$ and equivalence classes $[\overline{\gamma}]$ of positive loops $\overline{\gamma}\in \mathcal{Y}_k^+$ on the graph $\overline{\mathcal{G}}$.
\end{theorem}
In \cite{delecroixulcigrai20++enumerating}, the theorem is proven for general $k$ which of course includes the odd case. In the next section, instead of rewriting this proof to only cover the odd case, we will explain how to build this correspondence. 

\subsubsection{Construction of the Correspondence}

Let us know describe how to associate to an equivalence class of positive loops $[\gamma] \in \mathcal{Y}_k^+$, i.e., classes of positive loops of $\overline{\mathcal{G}}$ as defined in Definition \ref{def:equivalent_loops_up_to_commutation}, to a quadrangulation $Q = (\boldsymbol{\pi}, \boldsymbol{\lambda}, \boldsymbol{\tau})$, so that applying the Teichmüller geodesic flow $g_t$ to $Q$ (or more precisely, to the corresponding surface $X(Q)$) results in a closed orbit. Let us first explain the heuristic behind the construction, putting aside some technical issues which we will address further below. 

\begin{enumerate}
    \item First, suppose we are given a loop $\gamma$ in a graph which is associated to the application of staircase moves, i.e., $\gamma = \gamma_0\gamma_1 \cdots \gamma_{n-1}$.
    \item Associated to these staircase moves are matrices
    \begin{equation*}
        A_\gamma = A_{\gamma_0}\cdot A_{\gamma_1} \cdots A_{\gamma_{n-1}}
    \end{equation*}
    as in Definition \ref{def:staircase_matrix} that act on the length datum $\boldsymbol{\lambda}$. For the sake of this heuristic description, we will suppose for now that any length vector $\boldsymbol{\lambda}$ is a positive vector, meaning that all entries are (strictly) positive. 
    \item If $A_\gamma > 0$, i.e., if all entries of the matrix $A_\gamma$ are strictly positive, by the standard Perron--Frobenius theorem we know that there exists a unique positive eigenvector, which we will call $\boldsymbol{\lambda}$, associated to the leading eigenvalue of $A_\gamma$.
    \item Applying $A_\gamma$ to $\boldsymbol{\lambda}$, the length data of the quadrangulations we obtain remain in the infinite ray
    \begin{equation*}
        \{s\boldsymbol{\lambda} \mid s > 0\}.
    \end{equation*}
    Since $\gamma$ is a loop, the combinatorial datum of the newly obtained quadrilateral does not change and we will be able to show that $A_\gamma X(Q) = X(Q)$ as translation surfaces. In other words, applying the correct staircase moves to $Q$ we return back to $Q$ (in the future).
    \item Finally, we will see that applying $A_\gamma$ can be seen as applying the Teichmüller geodesic flow $g_{t_0}$ for some $t_0 > 0$, hence we obtain a closed geodesic.
\end{enumerate}

This description encompasses the main idea. What is still missing in this heuristic description is the height datum $\boldsymbol{\tau}$. Care needs to be taken to establish the existence of data that leads to a well-defined quadrangulation, i.e., it needs to be ensured that the data obtained satisfy the corresponding train-track relations. Moreover, some care needs to be given to the fact that the length datum $\boldsymbol{\lambda}$ is \emph{not} a vector with strictly positive entries but it is a vector in $(\R_-\times \R_+)^k$. The last issue which requires attention is that we will need to work with paths $\gamma$ in $\mathcal{G}$, i.e., with paths in the cover of $\overline{\mathcal{G}}$ which project to a loop but are not necessarily a loop themselves. The equivalence relations we have introduced in the previous section are designed exactly so that the construction outlined above can still be applied. 

Let us now make the above description more formal. As already mentioned, a standard Perron--Frobenius type argument will provide us with an eigenvector whose entries are all strictly positive, which we cannot use as a length datum since $\Delta_{\boldsymbol{\pi}} \not\subseteq \R_+^{2k}$, where $\Delta_{\boldsymbol{\pi}}$ is the set of admissible length data for the permutation datum $\boldsymbol{\pi}$, see section \ref{sec:parameter_space}. We will circumvent this issue by working with the absolute value $|\boldsymbol{\lambda}|$ of $\boldsymbol{\lambda}$, which is defined by
\begin{equation*}
    |\boldsymbol{\lambda}| \coloneqq \big( (|\lambda_{1, \ell}|, |\lambda_{1,r}|), \ldots ,(|\lambda_{k, \ell}|, |\lambda_{k,r}|)\big) = \big((-\lambda_{1, \ell}, \lambda_{1,r},\ldots, (-\lambda_{k, \ell}, \lambda_{k, r})\big).
\end{equation*}

\begin{lemma}\label{lem:JAJA}
    For any matrix $A_{\boldsymbol{\pi},c}$ as in Definition \ref{def:staircase_matrix}, we have
    \begin{equation*}
        A_{\boldsymbol{\pi},c} \cdot \boldsymbol{\lambda} = \boldsymbol{\lambda}'
    \end{equation*}
    if and only if
    \begin{equation*}
        A_{\boldsymbol{\pi},c} \cdot |\boldsymbol{\lambda}'| = |\boldsymbol{\lambda}|.
    \end{equation*}
\end{lemma}
\begin{remark}
    Note that the matrix $A_{\boldsymbol{\pi},c}$ acts on $\boldsymbol{\lambda}$ in the first case, but on $|\boldsymbol{\lambda}'|$ in the second case. 
\end{remark}
\begin{proof}
    Instead of a direct verification, let us introduce an auxiliary matrix $J$ which will be useful beyond this proof. We define the linear application
    \begin{align*}
        J \colon (\R_-\times\R_+)^k &\to \R_+^{2k}\\
        \big((\lambda_{1, \ell}, \lambda_{1,r},\ldots, (\lambda_{k, \ell}, \lambda_{k, r})\big) &\mapsto \big((-\lambda_{1, \ell}, \lambda_{1,r},\ldots, (-\lambda_{k, \ell}, \lambda_{k, r})\big),
    \end{align*}
    and we denote the associated $2k \times 2k$ matrix by $J$ as well. More explicitly, 
    \begin{equation*}
        J = \begin{bmatrix}
            -1& &  & & & &  \\
              &1&  & & & &  \\
              & &-1& & & &  \\
              & &  &1& & &  \\
              & &  & &\ddots& &  \\
              & &  & &      &-1&  \\
              & &  & &      &  & 1
        \end{bmatrix}.
    \end{equation*}
    Note that $J$ is bijective with $J = J^{-1}$, i.e., $J$ is an involution. Also, we have $|\boldsymbol{\lambda}| = J\boldsymbol{\lambda}$. Let us verify that
    \begin{equation}\label{eq:JAJA}
        JA_{\boldsymbol{\pi},c}J= A_{\boldsymbol{\pi},c}^{-1} \quad \text{for all } \boldsymbol{\pi} \in \mathcal{G}, c\in \boldsymbol{\pi}.
    \end{equation}
    Let us write $A \coloneqq A_{\boldsymbol{\pi},c}$ for ease of notation. By definition of matrix multiplication, $J$ and $JA_{\boldsymbol{\pi},c}$ have the same diagonal. If $i \neq j$ and $A_{ij} = 0$, then $(JA)_{ij} = 0$, since
    \begin{equation*}
        (JA)_{ij} = \sum_{m = 1}^{2k} J_{im}\cdot A_{mj} = J_{ii}\cdot A_{ij} = 0.
    \end{equation*}
    If $A_{ij} = a$ and say $c \in \pi_r$, so that $i$ is odd hence $J_{ii} = -1$, the same computation shows that $(JA)_{ij} = -a$. With this, we can verify that $(JA)^2 = \operatorname{Id}$, again using just the definition of matrix multiplication which proves \eqref{eq:JAJA}.

    It now follows, that
    \begin{equation*}
        |\boldsymbol{\lambda}'| = J\boldsymbol{\lambda}' = JA\boldsymbol{\lambda} = JAJJ\boldsymbol{\lambda} = A^{-1}J\boldsymbol{\lambda} = A^{-1}|\boldsymbol{\lambda}|.
    \end{equation*}
    \end{proof}
Analogously to \eqref{eq:path_matrix_1} we can define, for a path $\gamma$ in $\mathcal{G}$,
\begin{equation}\label{eq:dagger_matrix}
    A_\gamma^\dagger \coloneqq A_{\boldsymbol{\pi}_0, c_0}  \cdots A_{\boldsymbol{\pi}_{n-1}, c_{n-1}}.
\end{equation}
Notice that the order of the matrices here is reversed when compared to \eqref{eq:path_matrix_1}. Lemma \ref{lem:JAJA}
implies that
\begin{equation*}
    |\boldsymbol{\lambda}| = A_\gamma^\dagger |\boldsymbol{\lambda}'|,
\end{equation*}
where $\boldsymbol{\lambda}' = A_\gamma \cdot \boldsymbol{\lambda}$. Furthermore, it follows from the proof of the lemma that
\begin{equation}\label{eq:J_conjugation}
    A_\gamma^\dagger = (JA_{\boldsymbol{\pi}_0, c_0}J) \cdots (JA_{\boldsymbol{\pi}_{n-1}, c_{n-1}}J) = JA_\gamma^{-1}J.
\end{equation}
Hence, $\boldsymbol{\lambda}$ is an eigenvector of $A_\gamma^{-1}$, hence also of $A_\gamma$, if and only if $|\boldsymbol{\lambda}| = J\boldsymbol{\lambda}$ is an eigenvector for $A_\gamma^\dagger$. Moreover, the eigenvalue corresponding to the eigenvector $|\boldsymbol{\lambda}|$ of $A_\gamma^\dagger$ is $\mu > 1$ if and only if the eigenvalue associated to the eigenvector $\boldsymbol{\lambda}$ of $A_\gamma$ is $0<\mu^{-1}<1$.

\begin{remark}
    We can think of the matrix $J$ as a coordinate transformation that maps the part of $\R^{2k}$ that contains $\Delta_{\boldsymbol{\pi}}$ to the standard positive cone $\R_+^{2k}$. 
\end{remark}

We have shown that finding an eigenvector for $A_\gamma^\dagger$ with strictly positive entries may be translated to an eigenvector for $A_\gamma$ which we can see as providing a length datum $\boldsymbol{\lambda}$. We characterized positivity of loops $\overline{\gamma}$ through the associated matrix $A_\gamma$ in Definition \ref{def:positive_loops}. The following lemma shows that we may use the same condition on the matrices $A_\gamma^\dagger$ interchangeably. 

\begin{lemma}\label{lem:matrix_positivity}
    Let $\gamma$ be a path on $\mathcal{G}$. Then, $A_\gamma > 0$ if and only if $A_\gamma^\dagger > 0$. 
\end{lemma}
An elementary proof of this fact can be found in \cite{delecroixulcigrai20++enumerating} as Corollary 6.2.

As a consequence of the Markovian structure of staircase moves, we can deduce that any \emph{finite} path $\gamma$ on $\mathcal{G}$ can be realized, i.e., we can find quadrangulations, such that we can apply the staircase moves associated to $\gamma$ to these quadrangulations. The proof is simply given by applying Theorem \ref{thm:self_duality} inductively. If $\gamma = \gamma_0\gamma_1 \cdots \gamma_{n-1}$ is a finite path in $\mathcal{G}$, we will write $\hat{m}_{\gamma_i} = \hat{m}_{\boldsymbol{\pi}_i, c_i}$, where as in \eqref{eq:path_matrix_1} $c_i$ represents the $i\textsuperscript{th}$ edge associated to $\gamma_i$ and $\boldsymbol{\pi}_i = c_i\cdot c_{i-1} \cdots c_0 \cdot \boldsymbol{\pi}_0$ and define
\begin{equation*}
    \hat{m}_\gamma \coloneqq \hat{m}_{\gamma_{n-1}} \circ \cdots \circ \hat{m}_{\gamma_1}\circ \hat{m}_{\gamma_0}.
\end{equation*}

\begin{corollary}[Finite paths can be realized]\label{cor:finite_paths_realized}
    For any finite path $\gamma = \gamma_0\gamma_1 \cdots \gamma_{n-1}$ in $\mathcal{G}$, there exists a \emph{nonempty} subsimplex $\Delta_\gamma \subseteq \Delta_{\boldsymbol{\pi}_0}$, such that
    \begin{equation*}
        \{\boldsymbol{\pi}_0\} \times \Delta_\gamma \times \Theta_{\boldsymbol{\pi}_0} = \hat{m}_\gamma\left(\{\boldsymbol{\pi}_{n-1}\} \times \Delta_{\boldsymbol{\pi}_{n-1}} \times \Theta_{\boldsymbol{\pi}_{n-1}}\right),
    \end{equation*}
    where $\boldsymbol{\pi}_0$ denotes the combinatorial datum at the beginning of $\gamma$ and $\boldsymbol{\pi}_{n-1}$ denotes the combinatorial datum at the end of $\gamma$. In other words, given any quadrangulation $Q = (\boldsymbol{\pi}_0, \boldsymbol{\lambda}, \boldsymbol{\tau}) \in \{\boldsymbol{\pi}_0\}\times \Delta_\gamma \times \Theta_{\boldsymbol{\pi}_0}$, we can apply the sequence of moves $\hat{m}_\gamma$ to $Q$.
\end{corollary}

We want to obtain a result analogous to Corollary \ref{cor:finite_paths_realized} but for \emph{infinite paths. }. As outlined above, this will follow from a Perron--Frobenius type argument. Since these types of arguments rely on the matrices involved to be primitive, a further assumption is needed.

\begin{definition}[Primitive paths]\label{def:primitive_paths}
    An infinite path $\gamma = \gamma_1\gamma_2\cdots\gamma_n\cdots$ on $\mathcal{G}$ is \emph{primitive}, if there exists $n > 0$ such that if we denote by $\gamma(n)$ the initial path of length $n$ given by $\gamma(n) \coloneqq \gamma_0\cdots\gamma_{n-1}$, the matrix $A_{\gamma(n)}^\dagger$ as defined in \eqref{eq:dagger_matrix} is \emph{positive}, i.e., all the entries of $A_{\gamma(n)}$ are strictly positive.
\end{definition}

\begin{remark}
    By Lemma \ref{lem:matrix_positivity}, we could equivalently ask for the matrix $A_{\gamma(n)}$ to be positive.
\end{remark}

The following result shows that an analogous statement to Corollary \ref{cor:finite_paths_realized} holds for \emph{primitive} infinite paths.

\begin{proposition}[Infinite paths can be realized]\label{prop:infinite_path_realized}
    For any \emph{primitive} path $\gamma = \gamma_0\gamma_1 \cdots \gamma_n \cdots$ in $\mathcal{G}$, there exists a quadrangulation $Q = (\boldsymbol{\pi}_0, \boldsymbol{\lambda}, \boldsymbol{\tau})\in \mathcal{Q}_k$, to which we can apply the sequence of moves corresponding to $\gamma$, i.e., if we write $Q^{(i)} = \hat{m}_{\gamma_{i-1}}\circ \cdots \circ \hat{m}_{\gamma_0}(Q)$ as well as $Q^{(i)} = (\boldsymbol{\pi}_i, \boldsymbol{\lambda}_i, \boldsymbol{\tau}_i$), then $\boldsymbol{\lambda}_i \in \Delta_{\boldsymbol{\pi}_i, c_i}$, meaning that $\hat{m}_{\gamma_i}$ can be applied to $Q^{(i)}$ for any $i \in \N$. 
\end{proposition}
\begin{proof}
    For any $n>0$, applying Corollary \ref{cor:finite_paths_realized} to the (finite) truncation $\gamma(n)$ as defined in the statement of the proposition guarantees the existence of a nonempty cone $\Delta_{\gamma(n)}$ in $\R^{2k}$, such that any quadrangulation $Q = (\boldsymbol{\pi_0}, \boldsymbol{\lambda}, \boldsymbol{\tau})$ with $\boldsymbol{\lambda} \in \Delta_{\gamma(n)}$ satisfies the conclusion for $0 \leq i \leq n$. By construction, we have $\Delta_{\gamma(n+1)} \subseteq \Delta_{\gamma(n)}$ for any $n \in \N$. 

    Taking the closures $\operatorname{cl}(\Delta_{\gamma(n)})$ (with respect to the standard topology of $\R^{2k}$ followed by the projection map $\boldsymbol{v} \mapsto \R\boldsymbol{v}$ to the projective space $\PP(\R^{2k})$), we obtain a sequence of nested closed sets in $\PP(\R^{2k})$. Since finite dimensional projective spaces are compact, it follows by the finite intersection property of nested compact sets that the intersection
    \begin{equation*}
        \bigcap_{n > 0} \operatorname{cl}(\Delta_{\gamma(n)})
    \end{equation*}
    is nonempty. 

    We are left to check that any $\boldsymbol{\lambda} \in \bigcap_{n>0}\operatorname{cl}(\Delta_{\gamma(n)})$ belongs to the interior of $\operatorname{cl}(\Delta_{\boldsymbol{\pi}_0})$, or equivalently, that all entries are non-zero. To establish this, we will use the (proof of the) standard Perron--Frobenius theorem making use of the isomorphism $J$ introduced in the proof of Lemma \ref{lem:JAJA}. For any $n \in \N$, since $\boldsymbol{\lambda} \in \Delta_{\gamma(n)}$ there exists in particular some vector $\boldsymbol{v} \in (\R_-\times \R_+)^{k}$, such that $\boldsymbol{\lambda} = A_{\gamma(n)}^{-1}(\boldsymbol{v)}$. Noting that $J^{-1}\boldsymbol{v} \in \R_+^{2k}$, applying the map $J$, we obtain
    \begin{equation*}
        J\boldsymbol{\lambda} \in (JA_{\gamma(n)}^{-1}J)\R_+^{2k} = A_{\gamma_n}^\dagger \R_+^{2k}.
    \end{equation*}
    Since this is true for any $n \in \N$, it follows that
    \begin{equation}\label{eq:PF_intersection}
        |\boldsymbol{\lambda}| = J\boldsymbol{\lambda} \in \bigcap_{n>0} A_{\gamma(n)}^\dagger\R_+^{2k} = \bigcap_{n>0} A_{\boldsymbol{\pi}_0, c_0}\cdots A_{\boldsymbol{\pi}_{n-1}, c_{n-1}}(\R_+^{2k}).
    \end{equation}
    This is an intersection of nested non-negative cones. Since we assume $\gamma$ to be primitive, we know there exists $n_0 \in \N$ such that $A_{\gamma(n_0)}$ has strictly positive entries, which implies that $A_{\gamma(n)}^\dagger(\R_+^{2k})$ consists only of vectors whose entries are all strictly positive. It follows that every entry of $|\boldsymbol{\lambda}|$ must be non-zero. Lastly, since $\boldsymbol{\lambda} = J|\boldsymbol{\lambda}|$ and $J$ preserves the absolute value of the entries, also $\boldsymbol{\lambda}$ has only non-zero entries.
\end{proof}

For our purposes, the special case where the infinite path $\gamma_\infty$ is given by concatenating infinitely many times a given loop $\gamma$ is of particular interest. In this case, if the matrix $A_\gamma^\dagger$ is positive, from which it follows that $\gamma_\infty$ is primitive, we can obtain a stronger conclusion, namely that we obtain a unique (up to scaling) length datum $\boldsymbol{\lambda}$. We will make this precise in the corollary below. Recall first that a positive matrix has a simple largest eigenvalue whose eigenvector has strictly positive entries by the standard Perron--Frobenius theorem. We call \emph{Perron--Frobenius vector} the eigenvector corresponding to this leading eigenvalue of unit norm. The corollary is obtained by an application of the standard Perron--Frobenius theorem using again the coordinate transformation given by the matrix $J$.

\begin{corollary}[Perron--Frobenius uniqueness]\label{cor:perron_frobenius_uniqueness}
    If the infinite path $\gamma_\infty \coloneqq \gamma\gamma\gamma\cdots$ is given by concatenating a loop $\gamma$ on $\mathcal{G}$ infinitely often, such that the matrix $A_\gamma^\dagger$ is positive, then there exists $\boldsymbol{\lambda} \in \Delta_{\boldsymbol{\pi}_0}$, where $\boldsymbol{\pi}_0$ is the initial combinatorial datum of $\gamma$, such that all quadrangulations $Q \in \mathcal{Q}_k$ such that $\hat{m}_\gamma$ can be applied to $Q$ infinitely often belong to 
    \begin{equation*}
        \{Q = (\boldsymbol{\pi}_0, s\boldsymbol{\lambda}, \boldsymbol{\tau}) \mid \boldsymbol{\tau} \in \Theta_{\boldsymbol{\pi}_0}, \, s \in \R_+\}. 
    \end{equation*}
    Moreover, we can choose $\boldsymbol{\lambda}$ to be the unit vector $\boldsymbol{\lambda} = J\boldsymbol{v}$, where $\boldsymbol{v}$ is the unique Perron--Frobenius eigenvector of the matrix $A_\gamma^\dagger$.
\end{corollary}
\begin{proof}
    Let $L$ be the length of $\gamma$ and note that the truncations $\gamma_\infty(nL)$ of length $nL$ of $\gamma_\infty$, where we use the notation introduced in Definition \ref{def:primitive_paths}, are such that
    \begin{equation*}
        A_{\gamma_\infty(nL)}^\dagger = (A_\gamma^\dagger)^n.
    \end{equation*}
    Since $A_\gamma^\dagger > 0$, the infinite path $\gamma_\infty$ is primitive. Arguing as in Proposition \ref{prop:infinite_path_realized}, we can associate to the truncations $\{\gamma_\infty(nL) \mid n > 1\}$ the family of nested cones $\{(A_\gamma^\dagger)^n(\R_+)\mid n > 1\}$ and consider their intersection as in \eqref{eq:PF_intersection}. By positivity of $A_\gamma^\dagger$, these cones are (strictly) positive, so by showing that $A_\gamma^\dagger$ acts projectively as a strict contraction of the Hilbert projective metric as in the standard Perron--Frobenius proof, we can see that the intersection of these cones consists of a \emph{unique} ray of the form $\R_+\boldsymbol{v}$, where $\boldsymbol{v}\in \R_+^{2k}$ is the Perron--Frobenius eigenvector of $A_\gamma^\dagger$. The result now follows by applying $J$ to $\boldsymbol{v}$. 
\end{proof}


Having established these results of Perron--Frobenius type, in particular Corollary \ref{cor:perron_frobenius_uniqueness}, we can now explain how to associate a surface admitting a closed Teichmüller geodesic to a presentation of a positive loop $\overline{\gamma}$ in $\overline{\mathcal{G}}$. 

To this end, let $\overline{\gamma}$ be a presentation of a loop in $\overline{\mathcal{G}}$ starting from the vertex $[\boldsymbol{\pi}_0]$ and let $(\gamma, \sigma)$ be the lifted loop starting at $\boldsymbol{\pi}_0$ that lifts $\overline{\gamma}$ as in Definition \ref{def:lifted_loops}. Starting with $\overline{\gamma}$, we now explain how to construct a quadrangulation $Q$ such that the associated translation surface $X(Q)$ admits a closed Teichmüller geodesic.

First, the path $\gamma$ in $(\gamma, \sigma)$ is not necessarily a loop in $\mathcal{G}$, but it projects to the loop $\overline{\gamma}$ in $\overline{\mathcal{G}}$ and the cover $P \colon \mathcal{G} \to \overline{\mathcal{G}}$ is finite-to-one. If $n$ is the order of $\sigma$, i.e., $n$ satisfies $\sigma^n = \operatorname{id}$ and is the smallest element of $\N_{>0}$ with this property, we can find a loop $\gamma_0$ on $\mathcal{G}$ starting with $\gamma$ such that
\begin{equation*}
    P(\gamma_0) = \underbrace{\overline{\gamma} \cdots \overline{\gamma}}_{n-\text{times}}.
\end{equation*}
In fact, we can give such a loop explicitly, namely one checks that 
\begin{equation*}
    \gamma_0 = \gamma(\sigma^{-1}\star\gamma)(\sigma^{-2}\star\gamma)\cdots (\sigma^{-(n-1)}\star\gamma)
\end{equation*}
is such a loop. Since for any path $\gamma = \gamma_1\cdots\gamma_{m-1}$ and any relabeling $\sigma$ we have that
\begin{equation*}
    A_{\sigma\star\gamma} = \Pi_\sigma \cdot A_\gamma \cdot \Pi_{\sigma^{-1}},
\end{equation*}
it follows that
\begin{align*}
    A_{\sigma^{-1}\star\gamma}^\dagger &=
    A_{\sigma^{-1}\star\gamma_1}\cdots A_{\sigma^{-1}\star \gamma_{m-1}} \\&=
    (\Pi_{\sigma^{-1}}\cdot A_{\gamma_1} \cdot \Pi_\sigma) \cdots (\Pi_{\sigma^{-1}}\cdot A_{\gamma_{m-1}}\cdot \Pi_\sigma) \\&=
    \Pi_\sigma^{-1} \cdot A_\gamma \cdot \Pi_\sigma,
\end{align*}
where we used that $\Pi_{\sigma^{-1}} = \Pi_\sigma^{-1}$ in the last equality. Hence, we also have
\begin{align}\label{eq:staircase_matrix_power}
\begin{split}
    A_{\gamma_0}^\dagger &=
    A_\gamma^\dagger\cdot (\Pi_\sigma^{-1}A_\gamma^\dagger\Pi_\sigma)\cdot(\Pi_\sigma^{-2}A_\gamma^\dagger\Pi_\sigma^2) \cdots (\Pi_\sigma^{-(n-1)}A_\gamma\Pi_\sigma^{n-1}) \\&=
    (A_\gamma^\dagger\Pi_\sigma^{-1})^n,
\end{split}
\end{align}
where we further used that $\Pi_\sigma^{n-1} = \Pi_\sigma^{-1}$. 

We now consider the infinite path $\gamma_\infty = \gamma_0\gamma_0\gamma_0 \cdots$ obtained by concatenating $\gamma_0$ infinitely often. Note that $\gamma_\infty$ starts with $\gamma$. Since $\overline{\gamma}$ is positive as a loop in $\overline{\mathcal{G}}$, we know that $\Pi_\sigma A_\gamma$ is primitive and by a similar computation as the one just above we deduce that $A_\gamma^\dagger \Pi_\sigma^{-1}$ is primitive. By \eqref{eq:staircase_matrix_power} it follows that $A_{\gamma_0}^\dagger$ is primive as well, so that $\gamma_\infty$ is primitive in the sense of Definition \ref{def:primitive_paths}. Therefore, without loss of generality we may assume that $A_{\gamma_0}^\dagger > 0$ by choosing a sufficiently long multiple of $\gamma_0$ if necessary.

By Corollary \ref{cor:perron_frobenius_uniqueness}, we obtain a vector $\boldsymbol{\lambda} \in \Delta_{\boldsymbol{\pi}}$, where $\boldsymbol{\pi}$ is the initial combinatorial datum associated to $\gamma$, which is unique up to positive scaling and such that $\hat{m}_\gamma$ can be applied to all quadrilaterals of the form
\begin{equation*}
    Q_0 = (\boldsymbol{\pi}, s\boldsymbol{\lambda}, \boldsymbol{\tau}), \quad s\in \R_+, \, \tau \in \Theta_{\boldsymbol{\pi}}.
\end{equation*}
Moreover, by normalizing $\boldsymbol{\lambda}$ to have norm 1, we know that $|\boldsymbol{\lambda}| = J\boldsymbol{\lambda}$ is the Perron--Frobenius eigenvector of $A_{\gamma_0}^\dagger$. This $\boldsymbol{\lambda}$ will be the length datum of the quadrangulation we construct. Note that this choice of $\boldsymbol{\lambda}$ ensures exactly, that the path in $\overline{\mathcal{G}}$ corresponding to forwards staircase moves is obtained by repeating the loop $\overline{\gamma}$.

We now need to choose a height datum $\tau$, such that the path in $\overline{\mathcal{G}}$ corresponding to \emph{backwards} staircase moves is obtained by concatenating the loop $\overline{\gamma}$ (in backwards direction). We will do this in essentially the same way that we obtained $\boldsymbol{\lambda}$. 

In particular, we will need to choose $\boldsymbol{\tau}$ such that $\hat{m}_{\gamma_0}^{-1}$, the backwards staircase move associated to $\gamma_0$, can be applied infinitely often. So for any $n \in \N$, we need to find a quadrangulation $Q_{-n}\in \mathcal{Q}_k$ such that
\begin{equation*}
    Q_0 = \hat{m}_{\gamma_0}^n(Q_{-n}).
\end{equation*}
In particular, the height datum $\boldsymbol{\tau}$ of $Q_0$ must be such that
\begin{equation*}
    \boldsymbol{\tau} = A_{\gamma_0}^n \boldsymbol{\tau}_n
\end{equation*}
for some $\boldsymbol{\tau}_n \in \Theta_{\boldsymbol{\pi}}$ for any $n \in \N$. Hence,
\begin{equation}\label{eq:height_intersection}
    \boldsymbol{\tau}\in \bigcap_{n \in \N} A_{\gamma_0}^n(\Theta_{\boldsymbol{\pi}}) \subseteq \bigcap_{n \in \N}A_{\gamma_0}^n(\R_+^{2k}).
\end{equation}
From $A_{\gamma_0}^\dagger > 0$ we may deduce, using Lemma \ref{lem:matrix_positivity} that $A_{\gamma_0}>0$ holds as well. In contrast to the situation of the length datum $\boldsymbol{\lambda}$, the height data $\boldsymbol{\tau} \in \Theta_{\boldsymbol{\pi}}$ are already contained in the standard positive cone $\R_+^{2k}$, hence we may directly apply standard Perron--Frobenius arguments exactly as in Corollary \ref{cor:perron_frobenius_uniqueness} to deduce that there exists a unique (up to positive scaling) height vector $\boldsymbol{\tau}$ satisfying \eqref{eq:height_intersection}, where the one of unit length is exactly the Perron--Frobenius eigenvector of $A_{\gamma_0}$. For our purposes, a different choice will be more convenient, namely let $\boldsymbol{\tau}$ be the vector satisfying \eqref{eq:height_intersection} such that the area of the surface $X(Q_0)$, where $Q_0 = (\boldsymbol{\pi}, \boldsymbol{\lambda}, \boldsymbol{\tau})$ is 1.

We now argue that these choices of length and height data produces a quadrangulation that corresponds to a closed orbit. First, observe that if we let $t_1\in \R$ be such that $\e^{t_1}>0$ is the Perron--Frobenius eigenvalue of $A_{\gamma_0}^\dagger$, then, since $A_{\gamma_0}^\dagger > 0$, we have $\e^{t_1}>1$ and consequently $t_1 > 0$. We already discussed above that \eqref{eq:J_conjugation} implies that $\e^{-t_1}$ is an eigenvalue of $A_{\gamma_0}$ with respect to the eigenvector $\boldsymbol{\lambda}$. A computation similar to \eqref{eq:staircase_matrix_power} reveals that $A_{\gamma_0} = (\Pi^{-1}_\sigma A_\gamma)^n$, from which we deduce that
\begin{equation}\label{eq:closed_orbit_1}
    (\Pi^{-1}_\sigma A_\gamma) \boldsymbol{\lambda} = \e^{-t_0}\boldsymbol{\lambda}, \quad \text{where } t_0 = \frac{t_1}{n}>0
\end{equation}
using the uniqueness part of the Perron--Frobenius theorem. 

Similarly, $\boldsymbol{\tau}$ is an eigenvector of $A_{\gamma_0}$ corresponding to the eigenvalue $\e^{t_2}$ for $t_2 \in \R$. Let us show that $t_2 = t_1$. On one hand, it follows from
\begin{equation*}
    A_{\gamma_0}\boldsymbol{\lambda} = \e^{-t_1} \lambda \quad \text{and}\quad A_{\gamma_0}\boldsymbol{\tau} = \e^{t_2} \boldsymbol{\tau},
\end{equation*}
that the area of $Q_1 = \hat{m}_{\gamma_0}(Q_0)$ is $\e^{t_2 - t_1}$. On the other hand, staircase moves preserve the Lebesgue measure, hence
\begin{equation*}
    \e^{t_2 - t_1} = \operatorname{Area}(Q_1) = \operatorname{Area}(Q_0) = 1,
\end{equation*}
thus $t_2 = t_1$ as claimed. Repeating the argument just above, we deduce that
\begin{equation}\label{eq:closed_orbit_2}
    (\Pi^{-1}_\sigma A_\gamma)\boldsymbol{\tau} = \e^{t_0}\boldsymbol{\tau}, \quad \text{where } t_0 = \frac{t_1}{n} = \frac{t_2}{n}>0.
\end{equation}
Finally, setting $Q_n = \hat{m}_\gamma(Q_0)$ and writing $Q_n = (\boldsymbol{\pi}_n, \boldsymbol{w}_n)$, we can combine equations \eqref{eq:closed_orbit_1} and \eqref{eq:closed_orbit_2} to obtain
\begin{align*}
    \sigma \cdot Q_n &=
    \sigma \cdot \hat{m}_\gamma(Q_0) \\&=
    (\sigma \cdot \boldsymbol{\pi}_n, \Pi_\sigma A_\gamma \cdot \boldsymbol{\lambda}, \Pi_\sigma A_\gamma \cdot \boldsymbol{\tau}) \\&=
    (\boldsymbol{\pi}, \e^{-t_0}\boldsymbol{\lambda}, \e^{t_0}\boldsymbol{\tau}) \\&=
    g_{-t_0}(Q_0),
\end{align*}
where we use that $\sigma$ closes the loop $\gamma$, i.e., $\sigma \cdot \boldsymbol{\pi}_n = \boldsymbol{\pi}$. Equivalently, using the fact that the Teichmüller geodesic flow commutes with staircase moves, we can write
\begin{align}\label{eq:closed_orbit_3}
\begin{split}
    \sigma \cdot \hat{m}_\gamma(g_{t_0}\cdot Q_0) &=
    \sigma \cdot g_{t_0} \cdot Q_n \\&=
    (\sigma \cdot \boldsymbol{\pi}_n, g_{t_0}\cdot \boldsymbol{w}_n) \\&=
    (\boldsymbol{\pi}, g_{t_0}\cdot g_{-t_0} \boldsymbol{w}) \\&=
    (\boldsymbol{\pi}, \boldsymbol{w}) \\&=
    Q_0.
\end{split}
\end{align}
It is now easy to see that $f_\sigma \coloneqq \sigma \circ \hat{m}_\gamma$ defines a translation equivalence, since it is a composition of a cut and paste procedure given by staircase moves and a relabeling action, both of which are equivalences (see in particular Definition \ref{def:relabeling_translation_equivalence} below). Hence, \eqref{eq:closed_orbit_3} implies that
\begin{equation*}
    X(g_{t_0}\cdot Q_0) = X(Q),
\end{equation*}
as translation surfaces, so we did indeed obtain a closed geodesic starting with a positive loop $\overline{\gamma}$ in $\overline{\mathcal{G}}$. 

This construction provides a one-to-one correspondence between positive loops and closed geodesics. To obtain this result, it is necessary to verify the following, which is done in section 6 in \cite{delecroixulcigrai20++enumerating}.

\begin{enumerate}
    \item The mapping from the collection of positive loops $\mathcal{Y}_k^+$ to the collection of closed Teichmüller geodesics is \emph{well-defined}, i.e., it is independent of the representative $\overline{\gamma}$ with respect to the equivalence relation introduced in Definition \ref{def:equivalent_loops_up_to_commutation}. More precisely, if $(\gamma_1, \sigma_1)$ and $(\gamma_2, \sigma_2)$ are two lifted loops on $\mathcal{G}$ lifting  two positive presentations of loops $\overline{\gamma_1}$ and $\overline{\gamma_2}$ respectively, then applying the above construction to either positive loop $\overline{\gamma_1}$ or $\overline{\gamma_2}$ results in the same closed Teichmüller geodesic, or equivalently we obtain two pseudo-Anosov diffeomorphisms which are conjugated.
    \item The assignment is \emph{injective}. This can be shown by using properties of staircase moves $\hat{m}_\gamma$ we have introduced above.
    \item The assignment is \emph{surjective}. The proof of surjectivity in \cite{delecroixulcigrai20++enumerating} uses the language of cutting sequences and how it relates to staircase moves as introduced in \cite{delecroix2015diagonal} as well as the fact that we can see staircase moves as the Poincaré first return map to the Poincaré section $\Upsilon\subseteq \hyp{k}$ of the Teichmüller geodesic flow, which is introduced in the next section, crosses all closed geodesics. 
\end{enumerate}

\subsubsection{Equivalence relations and a Fundamental Domain of Quadrangulations}\label{sec:fundamental_domain_hyperelliptic}

There are two main goals we want to achieve in the next two sections. First, we want to describe an isomorphism between labeled translation surfaces in hyperelliptic components $\hyplab{k}$ and the space of quadrangulations $\mathcal{Q}_k/_\sim$ on which we need to introduce an equivalence relation based on staircase moves. This further allows us to explicitly construct a fundamental domain with respect to this equivalence relation $\sim$ which provides us with a canonical choice of quadrangulations. This isomorphism descends to an isomorphism of $\hyp{k}$ and the space of \emph{unlabeled} quadrangulations, which we will define below.

Second, we will define a section for the Teichmüller geodesic flow $g_t$ such that the Poincaré first return map is a diagonal changes algorithm corresponding to staircase moves. As we will show, this section crosses \emph{all} closed orbits of the Teichmüller flow, which is a fundamental property needed to use diagonal changes when counting or enumerating closed geodesics.

The first issue we are facing is that we need to introduce a labeling (of wedges or of bundles) if we want to talk about the parameter space we have introduced above. We thus cannot use a definition via parameters in order to define \emph{unlabeled} quadrangulations. In order to obtain such a space properly, we will start with the following definition.

\begin{definition}[Moduli space of labeled translation surfaces]
    The \emph{moduli space of labeled translation surfaces} is the collection of translation surfaces $X$ endowed with a labeling of its bundles of saddle connections $\Gamma_i$, where two such surfaces $X$ and $X'$ are \emph{isomorphic}, if they are isomorphic as translation surfaces and the corresponding isomorphism preserves the labels.

    We will write $\hyplab{k}$ for the set of equivalence classes of labeled translation surfaces obtained by labeling the surfaces in $\hyp{k}$ using cyclical labelings respecting the invorlution (see Definition \ref{def:good_labeling}.
\end{definition}

For any $X \in \hyp{k}$, there exist exactly $k$ different cyclical labelings respecting the involution, one for each choice of the starting bundle in Definition \ref{def:good_labeling}. Hence, $\hyplab{k}$ is a $k$-sheeted cover of $\hyp{k}$.

The projection $p \colon \hyplab{k} \to \hyp{k}$ consists of simply forgetting the labels. Taking the quotient with respect to the action by relabeling, we then obtain a faithful representation of connected components of strata, which is key to the purpose of applications such as counting pseudo-Anosov mapping classes in each of these components.

We want to distinguish between \emph{labeled} and \emph{unlabeled} quadrangulations in a similar way.

\begin{definition}[Compatibility of labels]\label{def:compatibility_of_labels}
    Given $X \in \hyp{k}$ with a labeling of the bundles $\Gamma_i, 1\leq i \leq k$, and a labeled quadrangulation $Q$, we say that the labelings of $Q$ and $X$ are \emph{compatible}, if the bottom wedge of $q_i \in Q$ consists of saddle connections in $\Gamma_i$.

    We write $\mathcal{Q}(X)$ for the collection of quadrangulations compatible with $X$.
\end{definition}

The following result, a proof of which can be found in chapter 4 of \cite{delecroixulcigrai20++enumerating}, will be useful below.
\begin{corollary}[Connectedness of $\mathcal{Q}(X)$]\label{cor:connectedness_of_Q(X)}
    For any $Q_1$ and $Q_2$ in $\mathcal{Q}(X)$, there exists a sequence of backward and forward staircase moves from $Q_1$ to $Q_2$. More precisely, there exist two sequences of staircase moves $\gamma_1$ and $\gamma_2$ such that $Q_2 = \hat{m}_{\gamma_2}^{-1}\hat{m}_{\gamma_1}(Q_1)$.
\end{corollary}

Recall from Definition \ref{def:space_labeled_quadrangulations} that we write $\mathcal{Q}_k$ for the space of labeled quadrangulations in $\hyp{k}$. Equivalently, we can write
\begin{equation*}
    \mathcal{Q}_k = \bigcup_{X \in \hyplab{k}} \mathcal{Q}(X),
\end{equation*}
i.e., $\mathcal{Q}_k$ is the union of all compatible quadrangulations $\mathcal{Q}(X)$ as $X$ varies in $\hyplab{k}$. We formally define \emph{unlabeled} quadrangulations as follows.
\begin{definition}[Unlabeled quadrangulations]
    Two quadrangulations $Q = (\boldsymbol{\pi}, \boldsymbol{w})$ and $Q' = (\boldsymbol{\pi}', \boldsymbol{w}')$ are \emph{equivalent as unlabeled quadrangulations}, if there exists a relabeling, that is a permutation $\sigma \in S_k$, such that $w'_i = w_{\sigma(i)}$ for all $1 \leq i \leq k$ and $\boldsymbol{\pi}' = \sigma\boldsymbol{\pi}\sigma^{-1}$. We denote by $\overline{Q}$ the equivalence class of $Q$ as unlabeled quadrangulations.
\end{definition}

We will always refer to an unlabeled quadrangulation as the equivalence class of one of its labeled representatives. 

\begin{definition}[Space of unlabeled quadrangulations]
    We write $\overline{\mathcal{Q}}_k$ for the \emph{space of unlabeled quadrangulations} which correspond to $\mathcal{Q}_k$, i.e., 
    \begin{equation*}
        \overline{\mathcal{Q}}_k = \{\overline{Q} \mid Q \in \mathcal{Q}_k\}.
    \end{equation*}
\end{definition}

Now that we have defined the relevant spaces, let us explain how to obtain an isomorphism
\begin{equation*}
    \hyplab{k} \xrightarrow{\sim} \mathcal{Q}/_\sim,
\end{equation*}
i.e., how to define the correct equivalence relation on $\mathcal{Q}_k$. 

\begin{definition}[Equivalence up to staircase moves]\label{def:equivalence_staircase_moves}
    Two quadrangulations $Q, Q' \in \mathcal{Q}_k$ are \emph{equivalent up to staircase moves} and we write $Q \sim Q'$, if $Q'$ can be obtained from $Q$ by a sequence of forward an backward staircase moves. This equivalence relation induces an equivalence relation $\sim$ on $\overline{\mathcal{Q}}_k$, which more explicitly is defined by setting $\overline{Q} \sim \overline{Q'}$ if and only if there exist representatives $Q'_1$ and $Q'_2$ in $\overline{Q}_1$ and $\overline{Q}_2$ respectively, such that $Q'_1 \sim Q'_2$.
\end{definition}

We obtain the following as a consequence of Corollary \ref{cor:connectedness_of_Q(X)}.

\begin{lemma}\label{lem:iso_labeled}
    The map
    \begin{align*}
        \Psi \colon \mathcal{Q}_k &\to \hyplab{k}, \\
        Q &\mapsto X(Q),
    \end{align*}
    where $X(Q)$ is the surface obtained by gluing the quadrangulation $Q = (\boldsymbol{\pi}, \boldsymbol{w})$ according to $\boldsymbol{\pi}$, is well-defined on $\mathcal{Q}_k/_\sim$ and induces an isomorphism on the subset of $\hyplab{k}$ of surfaces having neither horizontal nor vertical saddle connections.  
\end{lemma}
\begin{proof}
    If $Q \sim Q'$, then it is clear that the surface $X(Q)$ (labeled compatibly with $Q$) is isomorphic to $X(Q')$ (labeled compatibly with $Q'$) as a labeled translation surface, since the sequence of forward and backward staircase moves between $Q$ and $Q'$ provides a sequence of cut and paste operations between $X(Q)$ and $X(Q')$ which preserve the labeling. Therefore, $\Psi$ is well-defined on the quotient $\mathcal{Q}_k/_\sim$.

    In order to show that the map is surjective on the set of surfaces containing neither horizontal nor vertical saddle connections, recall that by Theorem \ref{thm:existence_quadrangulation} any such surface $X \in \hyp{k}$ admits a quadrangulation $Q$, which up to relabeling the bundles of $X$ we can assume to be compatible with $X$.

    For injectivity, assume that $Q$ and $Q'$ are two labeled quadrangulations, such that $X(Q) = X(Q') \eqqcolon X$, where $X$ has neither horizontal nor vertical saddle connections. Then, $Q$ and $Q'$ are both quadrangulations in $\mathcal{Q}(X)$ and by Corollary \ref{cor:connectedness_of_Q(X)} they are connected by a sequence of forward and backward staircase moves, hence $Q \sim Q'$.
\end{proof}

We now give a description of a fundamental domain $\mathcal{F}$ in the space of quadrangulations for the equivalence relation $\sim$ introduced above. We remark that the construction here is very much comparable to the construction we did in Section \ref{sec:RV_Teichmuller_connection} for zippered rectangles and Rauzy--Veech induction. 

\begin{definition}[Width intervals]\label{def:width_intervals}
    Let $Q = (\boldsymbol{\pi}, \boldsymbol{w}) = (\boldsymbol{\pi}, \boldsymbol{\lambda}, \boldsymbol{\tau})$ be a quadrangulation of a surface $X \in \hyp{k}$ and let $q_i \in Q$ be a quadrilateral with (forward) diagonal $w_{i,d}$. We define the intervals $I(q_i)$ and $I'(q_i)$ by
    \begin{align*}
        I(q_i) &= [\lambda_{i, \ell}, \lambda_{i,r}], \\
        I'(q_i) &=
        \begin{cases}
            [\lambda_{i, \ell}, \lambda_{i,d}] \quad &\text{if } q_i \text{ is left-slanted},\\
            [\lambda_{i,d}, \lambda_{i,r}] &\text{if } q_i \text{ is right-slanted}.
        \end{cases}
    \end{align*}
    We write $|I(q_i)|$ for the length of the interval $I(q_i)$, which is the \emph{horizontal width} of the quadrilateral $q_i$. Analogously, $|I'(q_i)|$ is the length of the interval $I'(q_i)$, which can be seen as the horizontal width of the quadrilateral $q'_i$ which is obtained by performing a diagonal change in $q_i$. Using this notation, we have $I'(q_i) = I(q'_i)$.
\end{definition}

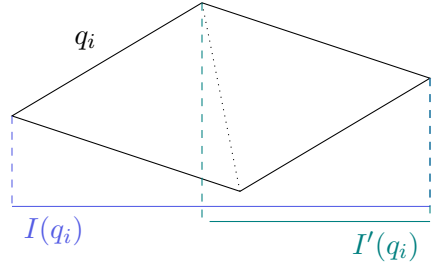
\begin{figure}[ht]
    \centering
    \begin{tikzpicture}
    \draw (0,0) -- (-3,1) -- (-0.5,2.5) -- (2.5,1.5) -- cycle;
    \draw[dotted] (0,0)-- (-0.5,2.5);
    \node[above left] at ($(-3,1)!0.5!(-0.5,2.5)$) {$q_i$};

    \draw[noamblue] (-3,-0.2) -- (2.5,-0.2) node[pos=0.1, below] {$I(q_i)$};
    \draw[teal] (-0.4,-0.4) -- (2.5, -0.4) node[pos=0.8, below] {$I'(q_i)$};

    \draw[noamblue,dashed] (-3,1) -- (-3,-0.2);
    \draw[noamblue, dashed] (2.5,1.5) -- (2.5,-0.2);
    \draw[teal, dashed] (-0.5,2.5) -- (-0.5,-0.4);
    \draw[teal, dashed] (2.5,1.5)-- (2.5, -0.4);
\end{tikzpicture}
    \caption{The intervals $I(q_i)$ and $I'(q_i)$ from Definition \ref{def:width_intervals}.}
    \label{fig:width_intervals}
\end{figure}

Definition \ref{def:width_intervals} is illustrated in Figure \ref{fig:width_intervals}.

\begin{definition}[Fundamental domain]\label{def:fundamental_domain}
    Let $\mathcal{F}$ be the set of all quadrangulations $Q \in \mathcal{Q}_k$ such that the following two conditions hold.
    \begin{enumerate}[start=1,label={($F$\arabic*)}]
        \item For every $q \in Q$ we have $|I(q)| \geq 1$.
        \item In every well-slanted staircase $S$ of $Q$ there exists $q \in S$ such that $|I'(q)|<1$.
    \end{enumerate}
\end{definition}

The following lemma states that $\mathcal{F}$ is indeed a fundamental domain with respect to $\sim$.

\begin{lemma}[Fundamental domain]\label{lem:fundamental_domain}
    Let $X \in \hyplab{k}$ be a labeled translation surface with neither vertical nor horizontal saddle connections. Then, there exists a unique quadrangulation $Q \in \mathcal{Q}(X)$ such that $Q \in \mathcal{F}$.
\end{lemma}

The key idea of the proof is to show that given any quadrangulation of $X$, we can always find another quadrangulation, equivalent up to staircase moves, which belongs to $\mathcal{F}$. For details we refer the reader to \cite{delecroixulcigrai20++enumerating}.

Since for any surface $X \in \hyplab{k}$ there is a \emph{unique} quadrangulation $Q \in \mathcal{F}$, this provides us with a canonical choice of quadrangulation.

\begin{definition}[Canonical quadrangulation]
    Given a labeled translation surface $X \in \hyplab{k}$ with neither vertical nor horizontal saddle connections, the \emph{canonical quadrangulation} associated to $X$, which we will denote by $Q(X)$, is the unique quadrangulation in $\mathcal{Q}_k$ compatible with the labeling of $X$ associated to $X$ by Lemma \ref{lem:fundamental_domain}.
\end{definition}

A key feature in the construction of the fundamental domain and the canonical quadrangulations $Q(X)$ is, that they are independent of the labeling of the bundles. Morally, this is clear since every aspect of the constructions depends only on the geometry of the quadrilaterals and not on the labeling of the bundles. Let us make this more precise. 

\begin{lemma}[Fundamental domain independent of labeling]\label{lem:fundamental_domain_label_independent}
    Let $X \in \hyplab{k}$ be a labeled translation surface and let $\sigma \in S_k$ be a relabeling compatible with the cyclical reordering of the bundles, i.e., the labeling of the relabeld surface $\sigma \cdot X$ is still cyclical respecting the involution. Then, the canonical quadrangulation $Q(\sigma \cdot X)$ of the relabeled surface $\sigma \cdot X$ is the same, up to relabeling, than the canonical quadrangulation of $Q(X)$ of $X$, namely
    \begin{equation}\label{eq:independence_of_labeling}
        Q(\sigma\cdot X) = \sigma \cdot Q(X),
    \end{equation}
    where $\sigma \cdot Q(X)$ denotes the quadrilateral relabeled by $\sigma$.
\end{lemma}
In other words, the lemma states that relabeling first and then obtaining the canonical quadrangulation always gives the same object as first getting the canonical quadrangulation and then applying the relabeling.
\begin{proof}
    The relabeled quadrangulation $\sigma \cdot Q(X)$ still satisfy conditions $(F1)$ and $(F2)$ from Definition \ref{def:fundamental_domain} since these conditions are independent of the labeling of the quadrangulation. By construction, $\sigma \cdot Q(X)$ is a quadrangulation of $\sigma\cdot X$ whose labeling is compatible with the labeling of $\sigma\cdot X$. Hence, $\sigma \cdot Q(X)$ is \emph{the unique} canonical quadrangulation $Q(\sigma\cdot X)$ of the surface $\sigma\cdot X$, i.e., \eqref{eq:independence_of_labeling} holds.
\end{proof}

\begin{corollary}\label{cor:canonical_quadrangulation_descends}
    The map $X \mapsto Q(X)$, which associates the canonical quadrangulation $Q(X)$ to a labeled translation surface $X \in \hyplab{k}$, is well-defined on equivalence classes of translation surfaces up to relabeling, hence it descends to a map on $\hyp{k}$, the hyperelliptic component of \emph{unlabeled} translation surfaces. 
\end{corollary}

In practice, we can use the following algorithmic procedure to explicitly construct the canonical quadrangulation $Q(X)$ for a translation surface $X\in\hyp{k}$ having neither horizontal nor vertical saddle connections. 

\begin{enumerate}
    \item Let $Q$ be any quadrangulation of $X$, which exists by Theorem \ref{thm:existence_quadrangulation}.
    \item Applying a finite number of backwards staircase moves if necessary, we can assume that $|I(q)|\geq 1$ for any $q \in Q$, i.e., $Q$ satisfies condition $(F1)$ from Definition \ref{def:fundamental_domain}.
    \item Define the set
    \begin{equation*}
        \mathcal{N}(Q) \coloneqq \{S_c \mid S_c \text{ is a well-slanted staircase in }Q \text{ for which } (F2) \text{ fails} \}.
    \end{equation*}
    \begin{enumerate}
        \item If $\mathcal{N}(Q) = \emptyset$, the quadrangulation $Q$ satisfies $(F2)$ as well, hence $Q$ is the canonical quadrangulation of $X$.
        \item If $\mathcal{N}(Q) \neq \emptyset$, write $Q'$ for the quadrangulation obtained by a staircase move in $S_c$ for some $S_c \in \mathcal{N}(Q)$. Note that $c$ is also a cycle of $\boldsymbol{\pi}' = c\cdot \boldsymbol{\pi}$. Since
        \begin{align*}
            |I(q'_i)| = |I(q_i)| \geq 1 \forall \, i \notin c, \\
            |I(q'_i)| = |I'(q_i)| \geq 1 \forall \, i \in c,
        \end{align*}
        the quadrangulation $Q'$ still satisfies $(F1)$. 
        \begin{enumerate}
            \item If the newly obtained staircase $S'_c$ is not well-slanted anymore or satisfies $(F2)$, then we have that
            \begin{equation*}
                |\mathcal{N}(Q')| = |\mathcal{N}(Q)| - 1.
            \end{equation*}
            \item If $S'_c$ is still well-slanted and does not satisfy $(F2)$, repeat step b) for $Q'$. By Lemma \ref{lem:diagonals_become_sides}, we know that after finitely many repetitions we must be in case i. Indeed, either the newly obtained staircase $S'_c$ at some point is not well-slanted anymore, or the width $|I'(q)|$ is eventually less than one, since by the lemma it tends to 0 as we continue applying staircase moves. 
        \end{enumerate}
    \end{enumerate}
\end{enumerate}

Lastly, let us state an prove a lemma that shows that translation equivalences act in a particularly nice way on canonical quadrangulations $Q(X)$. 

First, we define a translation equivalence given a quadrangulation $Q$ and a relabeling $\sigma$ such that $\sigma \cdot Q = Q$, meaning that $q_i$ and $q_{\sigma(i)}$ are identical (up to translation) for every $1 \leq i \leq k$ and $\sigma \star \boldsymbol{\pi} = \boldsymbol{\pi}$.

\begin{definition}[Relabeling translation equivalence]\label{def:relabeling_translation_equivalence}
    Let $Q \in \mathcal{Q}(X)$ and $\sigma \in S_k$ such that $\sigma \cdot Q = Q$. Define $f_\sigma \colon Q \to Q$ to be the isometry which sends $q_i$ to $q_{\sigma(i)}$ for each $1\leq i \leq k$. Then, $f_\sigma$ induces a translation equivalence of the surface $X = X(Q)$, which we will call the \emph{relabeling translation equivalence}, that we will also denote by $f_\sigma$. 
\end{definition}

For the canonical quadrangulation, the converse is also true, i.e., any translation equivalence of $X \in \hyp{k}$ is of this form.
\begin{lemma}[Translation equivalence is relabeling of canonical quadrangulation]
    Let $f \colon X \to X$ be a translation equivalence of $X \in \hyp{k}$. Then, there exists a relabeling $\sigma \in S_k$ such that the canonical quadrangulation $Q = Q(X)$ satisfies $\sigma \cdot Q = Q$ and $f$ is induced by the isometry $f_\sigma \colon Q \to Q$.
\end{lemma}
\begin{proof}
    First, fix a cyclical labeling of the bundles of saddle connections respecting the involution of $X$ (see Definition \ref{def:good_labeling}). Since by assumption $f$ is a translation equivalence, $f(X)$ is equal to $X$ as a translation surface in $\hyp{k}$ but possibly the labeling on $f(X)$ is different as the one induced by $f$. Since $f$ is an automorphism of a translation surface, it maps saddle connections to saddle connections and quadrangulations to quadrangulations, so that in particular the image of the canonical quadrangulation $Q$ (of $X$) under $f$ is also a quadrangulation (in $\mathcal{Q}(f(X))$), which we will denote by $Q'$. Since $f$ in particular preserves lengths, $Q'$ still satisfies conditions $(F1)$ and $(F2)$ from Definition \ref{def:fundamental_domain}. Let $\sigma$ be the relabeling such that $f$ maps the bundle $\Gamma_i$ to $\Gamma_{\sigma(i)}$. Clearly, $\sigma \cdot Q'$ still satisfies $(F1)$ and $(F2)$ and by construction, $Q'$ is compatible with $\sigma \cdot X$, i.e., $Q' \in \mathcal{Q}(\sigma\cdot X)$. By the uniqueness part of Lemma \ref{lem:fundamental_domain}, it follows that $Q' = Q(\sigma\cdot X)$ is \emph{the} canonical quadrangulation of $\sigma \cdot X$, which by Lemma \ref{lem:fundamental_domain_label_independent} is equal to $\sigma \cdot Q$. In particular, $Q'$ is isometric to $\sigma \cdot Q$ so that we may act on it by $f_\sigma^{-1}$, where $f_\sigma$ is the relabeling translation equivalence from Definition \ref{def:relabeling_translation_equivalence}. Considering now $f_\sigma^{-1}\colon f$, we obtain
    \begin{equation*}
        f_\sigma^{-1}\circ f (Q) = f_\sigma^{-1}(Q') = f_\sigma^{-1}(Q(\sigma\cdot X)) = f_\sigma\sigma \cdot Q = Q.
    \end{equation*}
    From this computation, it follows that $f_\sigma^{-1}\circ f$ fixes the quadrangulation $Q$ (viewed as a subset of $\R^{2k}$) and, since it is a translation equivalence (as a composition of two translation equivalences) its derivative is equal to the identity and it is continuous when restricted to each quadrilateral. Moreover, by choice of $\sigma$, the map $f_\sigma^{-1}\circ f$ maps each bundle $\Gamma_i$ to itself. It follows that $f_\sigma^{-1}\circ f$ is the identity on the whole domain, hence $f= f_\sigma$. 
\end{proof}

Moreover, the results from this section, in particular Lemma \ref{lem:iso_labeled} and Lemma \ref{lem:fundamental_domain} allow us to establish the following correspondences, the proofs of which are given in \cite{delecroixulcigrai20++enumerating}.

\begin{theorem}
    We have the following isomorphisms.
    \begin{enumerate}
        \item The set $\mathcal{F}$ is a fundamental domain for the quotient $\mathcal{Q}/_\sim$ of the space of (labeled) quadrangulations by the equivalence relation induced by staircase moves and it is isomorphic to the full measure subset of $\hyplab{k}$ consisting of labeled translation surfaces which have neither vertical nor horizontal saddle connections. The isomorphism is given by the map, which to any labeled translation surface $X \in \hyplab{k}$ associates the canonical quadrangulation $Q(X) \in \mathcal{Q}(X)$.
        \item The connected component $\hyp{k}$ is measurably isomorphic to the quotient $\overline{\mathcal{Q}}_k/_\sim$ of \emph{unlabeled}  quadrangulations by the equivalence relation induced by staircase moves. The isomorphism is given by the (almost everywhere defined) map $X \mapsto Q(X)$ from the first part of the Theorem, which is well-defined on $\hyp{k}$ by Corollary \ref{cor:canonical_quadrangulation_descends}, composed with the map which forgets the labeling.
    \end{enumerate}
\end{theorem}

\subsubsection{Section for the Teichmüller Geodesic Flow}\label{sec:section_for_teichmueller_flow}

Let us now define a subset of the fundamental domain which gives a Poincaré section for the Teichmüller geodesic flow $g_t$. 

\begin{definition}[Poincaré section]\label{def:poincare_section}
    Let $\Upsilon_\mathcal{F}$ be the set of all quadrangulations $Q \in \mathcal{F}$ that satisfy the following condition.
    \begin{enumerate}[start=1,label={(S)}]
        \item There exists $q \in Q$ such that $|I(q)| = 1$.
    \end{enumerate}
    In other words, $\Upsilon_\mathcal{F}$ is the set of all quadrangulations $Q \in \mathcal{Q}_k$ such that $Q$ satisfies $(S)$ as well as $(F1)$ and $(F2)$ from Definition \ref{def:fundamental_domain}.

    Further, we denote by $\Upsilon \subseteq \hyp{k}$ the set of all translation surfaces $X \in \hyp{k}$ for which there exists a quadrangulation $Q \in \mathcal{Q}_k$ which belongs to $\Upsilon_\mathcal{F}$.
\end{definition}

Let us remark that, if $X$ has neither vertical nor horizontal saddle connections, Lemma \ref{lem:fundamental_domain} implies that $X \in \Upsilon$ if and only if the canonical quadrangulation $Q(X)$ satisfies $(S)$. Moreover, $\Upsilon$ is a hypersurface in $\mathcal{F}$, i.e., a subspace of codimension 1. 

We will now describe the Poincaré first return map $\mathcal{S} \colon \Upsilon \to \Upsilon$ of the Teichmüller geodesic flow $(g_t)_{t\in \R}\colon \hyp{k} \to \hyp{k}$ to the section $\Upsilon\subseteq \hyp{k}$. Again, the main idea is comparable to the first return map we have described in section \ref{sec:RV_Teichmuller_connection} using the language of zippered rectangles. Since staircase moves commute with the geodesic flow $g_t$, we can describe the first return map as follows.

\begin{enumerate}
    \item View the section $\Upsilon$ as the boundary of $\mathcal{F}$.
    \item Apply (a) staircase move(s) so that the newly obtained quadrangulation $Q'$ is now in the exterior of $\mathcal{F}$, i.e., so that $|I(q)| < 1$ for all $q \in Q'$.
    \item Apply the flow $g_t$, which increases the length of the width intervals $|I(q)|$ for $q \in Q'$, until the first time $t_0$ where $|I(q)| \geq 1$ for all $q \in Q'$ and there exists a $\Tilde{q}\in Q'$ such that $|I(\Tilde{q})| = 1$. This means exactly that $g_{t_0}(X) \in \Upsilon$, i.e., the map indeed describes the first return to the section. 
\end{enumerate}

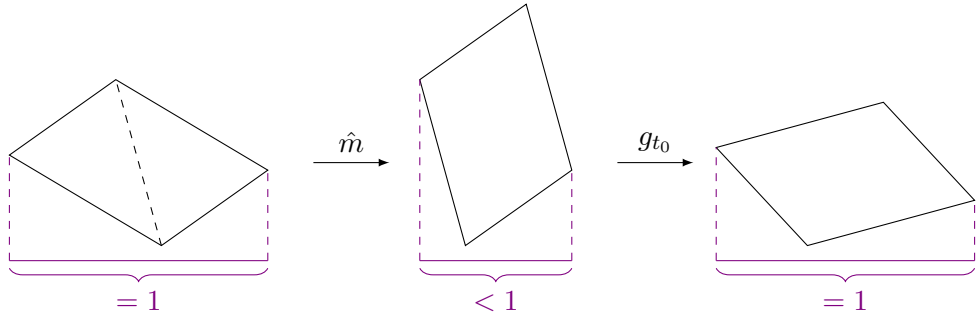
\begin{figure}[ht]
    \centering
    \begin{tikzpicture}
    \draw (0,0) -- (-2,1.2) -- (-0.6,2.2) -- (1.4,1) -- cycle;
    \draw[violet] (-2,-0.2) -- (1.4,-0.2);
    \draw [violet, decorate, decoration={brace, amplitude=5pt}] (1.4,-0.3) -- node[below=5pt] {$=1$}  (-2,-0.3);
    \draw[dashed] (0,0) -- (-0.6,2.2);
    \draw[violet, dashed] (-2,-0.2) -- (-2,1.2);
    \draw[violet, dashed] (1.4,-0.2) -- (1.4,1);

    \draw[->, >=latex] (2,1.1) -- (3,1.1) node[midway, above] {$\hat{m}$};

    \def\x{4};
    \draw (\x,0) -- (\x-0.6, 2.2) -- (\x+0.8, 3.2) -- (\x+1.4,1) -- cycle;
    \draw[violet] (\x-0.6,-0.2) -- (\x+1.4,-0.2);
    \draw[violet, dashed] (\x-0.6,-0.2) -- (\x - 0.6,2.2);
    \draw[violet, dashed] (\x +1.4, -0.2) -- (\x+1.4,1);
    \draw [violet, decorate, decoration={brace, amplitude=5pt}] (\x+1.4,-0.3) -- node[below=5pt] {$<1$}  (\x-0.6,-0.3);

    \draw[->, >=latex] (6,1.1) -- (7,1.1) node[midway, above] {$g_{t_0}$};

    \def\w{8.5}
    \draw (\w-1.2,1.3) -- (\w,0) -- (\w+2.2,0.6) -- (\w+1,1.9) -- cycle;
    \draw[violet] (\w-1.2,-0.2) -- (\w+2.2,-0.2);
    \draw[dashed, violet] (\w-1.2,-0.2) -- (\w-1.2,1.3);
    \draw[dashed, violet] (\w+2.2,-0.2) -- (\w+2.2,0.6);
    \draw [violet, decorate, decoration={brace, amplitude=5pt}] (\w+2.2,-0.3) -- node[below=5pt] {$=1$}  (\w-1.2,-0.3);
\end{tikzpicture}
    \caption{The first return map $\mathcal{S}\colon \Upsilon\to\Upsilon$.}
    \label{fig:first_return_diagonal_changes}
\end{figure}

The first return map $\mathcal{S}$ is illustrated in Figure \ref{fig:first_return_diagonal_changes}, where for simplicity we only pictured a single quadrilateral $q$.

As in section \ref{sec:RV_Teichmuller_connection}, it is possible to work out the exact first return time $t_0$. Let $X \in \hyp{k}$ be a translation surface with neither horizontal nor vertical saddle connections and let $Q(X) = (\boldsymbol{\pi}, \boldsymbol{w})$ be the associated canonical quadrangulation. For each cycle $c$ in $\pi_\ell$ or $\pi_r$ such that the associated staircase $S_c$ is well-slanted, there exists (at least) one quadrilateral $q \in S_c$ such that $|I'(q)| < 1$. Let us denote by $q_c$ the quadrilateral in $S_c$ such that $|I'(q)|$ is minimal, i.e., such that
\begin{equation*}
    |I'(q_c)| = \min\{|I'(q)|\mid  q \in S_c \}.
\end{equation*}
Since by Theorem \ref{thm:existence_well_slanted_staircase} the set of well-slanted staircases in $Q(X)$ is nonempty, we can further define $\Tilde{c}$ to be the cycle of $\pi_\ell$ or $\pi_r$ such that
\begin{equation*}
    |I'(q_{\Tilde{c}})| = \max\{|I'(q_c)| \mid c \text{ a cycle of }\pi\ell \text{ or }\pi_r \text{ such that }S_c \text{ is well-slanted}\}.
\end{equation*}
Finally, we can now define $t_0$ to be the positive real number such that
\begin{equation*}
    \e^{t_0}|I'(q_{\Tilde{c}})| = 1,
\end{equation*}
which is equivalent to
\begin{equation*}
    t_0 = -\log|I'(q_{\Tilde{c}})|.
\end{equation*}
Note that $t_0 > 0$ since by definition of $\Tilde{c}$ we have $|I'(q_{\Tilde{c}})|<1$. A verification that $g_{t_0}$ does indeed correspond to the first return map $\mathcal{S} \colon \Upsilon \to \Upsilon$ can be found as Lemma 5.14 in \cite{delecroixulcigrai20++enumerating}. Moreover, we will give essentially the same argument in section \ref{sec:general_counting} (in the proof of Lemma \ref{lem:fundamental_domain_general}) when extending the ideas to non-hyperelliptic strata. 

Note that by construction of the Poincaré first return map, we can see that it corresponds to first applying staircase moves on the quadrangulation followed by a suitable renormalization, which is exactly the same we saw in the case of Rauzy--Veech induction. To emphasize this idea, we proceed with an explicit example. 

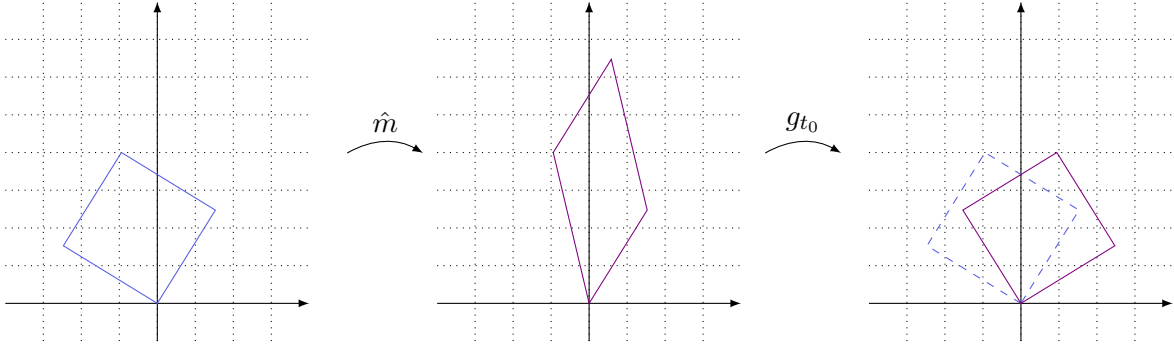
\begin{figure}[ht]
    \centering
    \begin{tikzpicture}[scale = 2]
    \coordinate (xleft) at (-1,0); \coordinate (xright) at (1,0);
    \coordinate (ytop) at (0,2); \coordinate (ybot) at (0,-0.25);

    \foreach \p in {0.25,0.5,...,1.75}{
        \draw[dotted] (-1,\p) -- (1,\p);
    }
    \foreach \p in {-0.75,-0.5,-0.25,0,0.25,0.5,0.75}{
        \draw[dotted] (\p, -0.25) -- (\p, 2);
    }
    
    \draw[>=latex, ->] (xleft) -- (xright);
    \draw[>=latex, ->] (ybot) -- (ytop);

    \coordinate (w1l) at ({-1/2*(sqrt(5)-1)}, {1/2*(3-sqrt(5))});
    \coordinate (w1r) at ({1/2*(3-sqrt(5))}, {1/2*(sqrt(5)-1)});

    \draw[noamblue] (0,0) -- (w1l) -- ($(w1l) + (w1r)$) -- (w1r) -- cycle;
\end{tikzpicture}
\qquad\qquad
\begin{tikzpicture}[scale = 2]
    \coordinate (xleft) at (-1,0); \coordinate (xright) at (1,0);
    \coordinate (ytop) at (0,2); \coordinate (ybot) at (0,-0.25);

    \foreach \p in {0.25,0.5,...,1.75}{
        \draw[dotted] (-1,\p) -- (1,\p);
    }
    \foreach \p in {-0.75,-0.5,-0.25,0,0.25,0.5,0.75}{
        \draw[dotted] (\p, -0.25) -- (\p, 2);
    }
    
    \draw[>=latex, ->] (xleft) -- (xright);
    \draw[>=latex, ->] (ybot) -- (ytop);

    \coordinate (w1l) at ({2-sqrt(5)}, 1);
    \coordinate (w1r) at ({1/2*(3-sqrt(5))}, {1/2*(sqrt(5)-1)});

    \draw[violet] (0,0) -- (w1l) -- ($(w1l) + (w1r)$) -- (w1r) -- cycle;
\end{tikzpicture}
\qquad\qquad
\begin{tikzpicture}[scale = 2]
    \coordinate (xleft) at (-1,0); \coordinate (xright) at (1,0);
    \coordinate (ytop) at (0,2); \coordinate (ybot) at (0,-0.25);

    \foreach \p in {0.25,0.5,...,1.75}{
        \draw[dotted] (-1,\p) -- (1,\p);
    }
    \foreach \p in {-0.75,-0.5,-0.25,0,0.25,0.5,0.75}{
        \draw[dotted] (\p, -0.25) -- (\p, 2);
    }
    
    \draw[>=latex, ->] (xleft) -- (xright);
    \draw[>=latex, ->] (ybot) -- (ytop);

    \coordinate (w1l) at ({-1/2*(sqrt(5)-1)}, {1/2*(3-sqrt(5))});
    \coordinate (w1r) at ({1/2*(3-sqrt(5))}, {1/2*(sqrt(5)-1)});

    \coordinate (v1l) at ({-1/2*(3-sqrt(5))}, {1/2*(sqrt(5)-1)});
    \coordinate (v1r) at ({1/2*(sqrt(5)-1)}, {1/2*(3-sqrt(5))});

    \draw[noamblue, dashed] (0,0) -- (w1l) -- ($(w1l) + (w1r)$) -- (w1r) -- cycle;
    \draw[violet] (0,0) -- (v1l) -- ($(v1l) + (v1r)$) -- (v1r) -- cycle;
\end{tikzpicture}
\begin{tikzpicture}[remember picture, overlay]
    \draw[>=latex, ->, bend left=30] (-11,2.5) to node[above] {$\hat{m}$} (-10,2.5);
    \draw[>=latex, ->, bend left=30] (-5.5,2.5) to node[above] {$g_{t_0}$} (-4.5,2.5);
\end{tikzpicture}
    \caption{The three quadrilaterals $q, q'$ and $q''$ from Example \ref{ex:first_return}.}
    \label{fig:first_return_example}
\end{figure}

\begin{figure}[!ht]
    \centering
    \begin{tikzpicture}[scale = 3]
        \coordinate (xleft) at (-1,0); \coordinate (xright) at (1,0);
    \coordinate (ytop) at (0,3); \coordinate (ybot) at (0,-0.25);

    \foreach \p in {0.25,0.5,...,2.75}{
        \draw[dotted] (-1,\p) -- (1,\p);
    }
    \foreach \p in {-0.75,-0.5,-0.25,0,0.25,0.5,0.75}{
        \draw[dotted] (\p, -0.25) -- (\p, 3);
    }
    
    \draw[>=latex, ->] (xleft) -- (xright);
    \draw[>=latex, ->] (ybot) -- (ytop);

    \coordinate (w1l) at ({-1/2*(sqrt(5)-1)}, {1/2*(3-sqrt(5))});
    \coordinate (w1r) at ({1/2*(3-sqrt(5))}, {1/2*(sqrt(5)-1)});

    \draw[noamblue] (0,0) -- (w1l) -- ($(w1l) + (w1r)$) -- (w1r) -- cycle;

    \coordinate (v1l) at ({2-sqrt(5)}, 1);
    \coordinate (v1r) at ({1/2*(7-3*sqrt(5))}, {1/2*(1+sqrt(5))});

    \draw[purple] (0,0) -- (v1l) -- ($(v1l) + (v1r)$) -- (v1r) -- cycle;

    \foreach \t in {0.2, 0.4, 0.6,0.8}{
        \coordinate (l) at ({exp(\t)*(2-sqrt(5))}, {exp(-\t)});
        \coordinate (r) at ({exp(\t)*1/2*(7-3*sqrt(5))}, {exp(-\t)*1/2*(1+sqrt(5))});
        \draw[dashed, opacity = 0.5, noamblue!\t!purple] (0,0) -- (l) -- ($(l) + (r)$) -- (r) -- cycle;
    }
\end{tikzpicture}
    \caption{The red quadrilateral is obtained by applying two staircase moves to the blue quadrilateral. Applying $g_{t_0}$ for $t_0$ as in Example \ref{ex:first_return} gives back the initial quadrilateral $q$, hence we obtain a \emph{closed} geodesic for the Teichmüller flow.}
    \label{fig:first_return_example_2}
\end{figure}
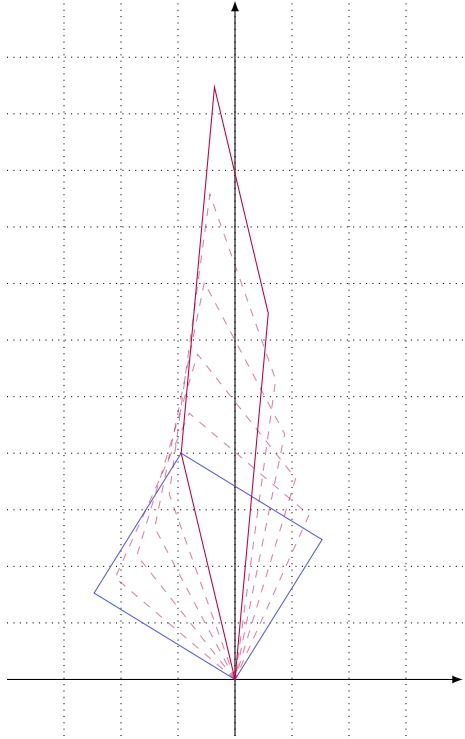

\begin{example}[Diagonal changes as first return of $g_t$]\label{ex:first_return}
    To see the first return map $\mathcal{S} \colon \Upsilon \to \Upsilon$ explicitly, we will apply the construction to the following surface in $\hyp{1}$, which corresponds to a torus.

    The translation surface $X$ is represented by the quadrangulation $Q$ containing a single quadrilateral $q$, which is characterized by the wedge
    \begin{equation*}
        \boldsymbol{w} = (w_\ell, w_r) = \left(\left(-\frac{1}{2}(\sqrt{5}-1), \frac{1}{2}(3 - \sqrt{5})\right), \left(\frac{1}{2}(3-\sqrt{5}), \frac{1}{2}(\sqrt{5}-1)\right)\right).
    \end{equation*}
    The quadrilateral is pictured in Figure \ref{fig:first_return_example} on the left. A direct computation shows that $|I(q)| = 1$ and $|I'(q)| = \frac{1}{2}(\sqrt{5}-1)$, which implies that $X \in \Upsilon$.

    Since staircase moves commute with the Teichmüller geodesic flow $g_t$, we may start by applying a diagonal change. We obtain a new quadrilateral $q'$ which is characterized by the wedge
    \begin{equation*}
        \boldsymbol{w}' = \left((2-\sqrt{5},1),w_r\right),
    \end{equation*}
    where we applied the formulas worked out in section \ref{sec:staircase_moves}. The quadrilateral $q'$ is depicted in the middle of Figure \ref{fig:first_return_example}. Lastly, we apply $g_{t_0}$, where $t_0$ is the first return time given in this case by
    \begin{equation*}
        t_0 = -\log|I'(q)| = -\log\left(\frac{1}{2}(\sqrt{5}-1)\right).
    \end{equation*}
    The resulting quadrilateral $q'' = g_{t_0}(q')$ is depicted on the right in Figure \ref{fig:first_return_example}. One verifies that $|I(q'')| = 1$ and $I'(q'') = I'(q)$, so that indeed we have returned to the section $\Upsilon$.

    We remark that the construction here does \emph{not} characterize a closed geodesic of the Teichmüller flow. The reason is, that the associated matrix $A_{\boldsymbol{\pi},c}$ from Definition \ref{def:staircase_matrix} to the right move applied here is given by
    \begin{equation*}
        A_{\boldsymbol{\pi},c} = \begin{bmatrix}
            1 & 1 \\
            0 & 1
        \end{bmatrix},
    \end{equation*}
    which fails to be \emph{positive}. Note that $q''$ is left-slanted, meaning that if we repeat the above construction with $q''$ and consider the loop $\gamma$ given by concatenating the two staircase moves, we obtain the matrix
    \begin{equation*}
        A_\gamma = A_{\boldsymbol{\pi}'',c''} \cdot A_{\boldsymbol{\pi}, c} = \begin{bmatrix}
            1 & 0 \\
            1 & 1 
        \end{bmatrix}\cdot 
        \begin{bmatrix}
            1 & 1 \\
            0 & 1
        \end{bmatrix} = 
        \begin{bmatrix}
            2 & 1 \\
            1 & 1
        \end{bmatrix},
    \end{equation*}
    which \emph{does} satisfy $A_\gamma > 0$. So applying \emph{two} diagonal changes and flowing for time $2t_0$ (recall that $|I'(q)| = |I'(q'')|$) we obtain a closed Teichmüller geodesic of minimal period $2t_0$, or equivalently, a pseudo-Anosov diffeomorphism with dilatation $2t_0$. This is illustrated in Figure \ref{fig:first_return_example_2}.
\end{example}

Lastly, we present the following fact announced in the beginning of this section. The Poincaré section from Definition \ref{def:poincare_section} intersects \emph{all} closed Teichmüller geodesics. This fact holds significant importance for using diagonal changes in problems concerning enumeration or counting.

\begin{proposition}\label{prop:crosses_all_geodesics}
    The Poincaré section $\Upsilon$ from Definition \ref{def:poincare_section} intersects \emph{all} closed Teichmüller geodesics in $\hyp{k}$.
\end{proposition}
\begin{proof}
    The idea of the proof is simple. We first argue that any surface $X$ belonging to some closed Teichmüller geodesic $\gamma$ admits a canonical quadrangulation $Q(X)$. Then, we apply the \emph{backward} geodesic flow to $Q(X)$, which contracts the widths of the quadrilaterals, until $Q(g_{-t}(X))$ satisfies condition (S) from Definition \ref{def:poincare_section}, showing that the intersection between $\Upsilon$ and $\gamma$ is nonempty. 

    So to this end, let $\gamma$ be a closed Teichmüller geodesic and $X$ any surface in $\gamma$. The surface $X$ is obviously recurrent, both in forward and backward direction, from which it follows immediately that $X$ cannot admit a horizontal nor a vertical saddle connection. Indeed, If $X$ were to have such a saddle connection it would necessarily be divergent under the Teichmüller flow. We may invoke Lemma \ref{lem:fundamental_domain}, by which we deduce that there exists a canonical quadrangulation $Q(X) \in \mathcal{F}$. Since applying the Teichmüller flow in backward direction contracts the widths of the quadrilaterals, it suffices now to to consider $g_{-t_0}(Q(X))$, where $t_0$ is the first time at which condition (S) from Definition \ref{def:poincare_section} is satisfied, hence by the remark after the definition we have that $g_{-t_0}(X) \in \Upsilon$.
\end{proof}
 \pagebreak

\newpage
\section{Diagonal Changes in General Components of Strata}\label{sec:general_counting}
\thispagestyle{plain}

We will now explore further the situation outside of hyperelliptic strata. As discussed in the conclusion of Section \ref{sec:diagonal_changes}, and further exemplified in Example \ref{ex:no_staircase}, it is not possible to generalize the approach using staircases directly to arbitrary strata. However, in \cite{ferenczi2015diagonal}, Ferenczi develops a framework of diagonal changes that work for \emph{every} IET, not just those associated to a surface in a hyperelliptic stratum. In this section, we will show that the geometric description of his algorithm can be applied to \emph{any} translation surface in \emph{any} stratum, thereby extending the theory developed in section \ref{sec:diagonal_changes} and originally proposed in \cite{delecroix2015diagonal}. Analogous to the discussion in \ref{sec:counting}, we are able to construct a fundamental domain in the set of polygonal representations of translation surfaces used to define the algorithm, as well as a Poincaré section such that its associated first return map of the Teichmüller geodesic flow is given exactly by diagonal changes. These constructions primarily involve adapting the concepts from \cite{bell2019coding}, which focuses on algorithms utilizing so-called \emph{Veering triangulations}, to our specific context. 

First, we will introduce the diagonal changes algorithm from \cite{ferenczi2015diagonal} and define the parameter space on which it acts. While Ferenczi does give a formal description of forward moves in \cite{ferenczi2014generalization}, we will extend this description to include backward diagonal changes as well. In the second part of this section, we will show that this algorithm may be applied to any translation surface by showing that it admits an adequate decomposition. Here, the key ideas will be from \cite{bell2019coding}. Lastly, adapting the arguments in \cite{delecroix2015diagonal} for the hyperelliptic case, we explain how to obtain both a fundamental domain and a Poincaré section for the Teichmüller geodesic flow. 

\subsection{Diagonal Changes on Castle Polygons}

The generalized version of the diagonal changes algorithm will be defined on objects called \emph{sets of castle polygons}. We start with some preliminary definitions. Note that the plane comes equipped with an orientation, so it makes sense to talk about left, right, up and down.

\begin{definition}[Base and stack triangles]\label{def:base_stack_triangles}
    Let $T \subseteq \R^2$ be a triangle without horizontal or vertical edges and let $\mathbf{v}$ be its lowest vertex, which is therefore uniquely defined. We say that $T$ is a \emph{base triangle}, if the two edges $w_1, w_2$ from $\mathbf{v}$ form a wedge that is opened upwards, i.e., if
    \begin{equation*}
        (w_i, w_j) \in \big((\R_- \times \R_+) \times (\R_+\times \R_+)\big), \quad i \neq j.
    \end{equation*}

    If $T$ is not a base triangle, we call it a \emph{stack triangle} if its highest vertex lies between the other two vertices with respect to the horizontal. Equivalently, a triangle that is not a base triangle is a stack triangle if the vertical downwards ray from its highest vertex crosses the opposing side.
\end{definition}

\begin{figure}[ht]
    \centering
    \begin{tikzpicture}
        \draw (0,0) node[below] {$\mathbf{v}$} -- node[midway, below right] {$w_2$} (2,1.5) -- (-1,2.2) -- cycle node[midway, left] {$w_1$};

        \def \x{3};
        \draw (\x,0) node[below] {$\mathbf{v}$} -- node[midway, below, yshift = -2] {$w_2$} (\x+3,1) -- (\x+2,2.1) -- cycle node [midway, above left] {$w_1$};
        \draw[dashed, noamblue] (\x+2, 2.1) -- ++(-90:2); 

        \def \w{7.5}
        \draw (\w,0) node[below] {$\mathbf{v}$} -- node[midway, below right] {$w_2$} (\w+3,2.2) -- (\w+1,1.7) -- cycle node[midway, left] {$w_1$};
        \draw[dashed, purple] (\w+3, 2.2) -- ++ (-90:2);

    \end{tikzpicture}
    \caption{A base triangle, a stack triangle and a triangle that is neither a base nor a stack triangle.}
    \label{fig:base_stack_neither}
\end{figure}
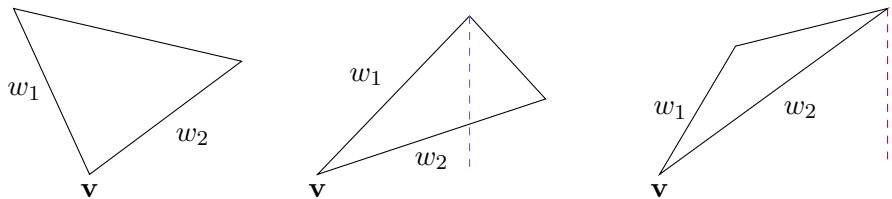

In Figure \ref{fig:base_stack_neither}, we see an example of a base triangle, a stack triangle and a triangle that is neither a base nor a stack triangle. The next definition justifies the terminology of \emph{base} and \emph{stack} triangles.

\begin{definition}[$k$-set of castle polygons]\label{def:castle_polygons}
    Let $k \in \N_{\geq 2}$. A \emph{$k$-set of castle polygons} is a collection of $k$ disjoint polygons equipped with a triangulation, such that the following hold. 
    \begin{enumerate}[i)]
        \item The lowest vertex of polygon $j$ is denoted by $\mathbf{v}_j$ and is the lowest vertex of a base triangle of the triangulation.
        \item If we denote by $(w_{i,\ell}, w_{i, r})$ the wedge defining the base triangle of polygon $i$, then the upper sides of the polygons form a partition of the set
        \begin{equation*}
            \{w_{1, \ell}, w_{1, r}, \ldots, w_{k, \ell}, w_{k, r}\}.
        \end{equation*}
        \item Every triangle in the triangulation that does not contain any lowest vertex $\mathbf{v}_j$ is a stack triangle. 
        \item There is no \emph{strict} subset $J$ of $[k]$ such that all the upper sides of the polygons $i$ for $i \in J$ correspond to the wedges $w_{j, \ell}$ or $w_{j, r}$ for $j \in J$.
    \end{enumerate}
\end{definition}

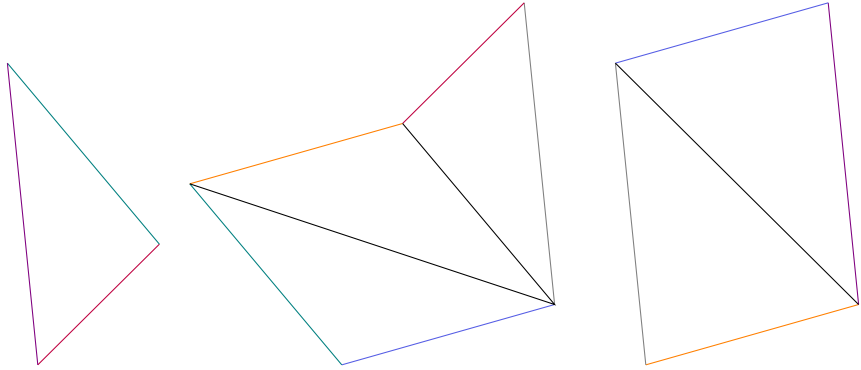
\begin{figure}[ht]
    \centering
    \begin{tikzpicture}[scale = 0.4]
        \coordinate (p1) at (0,0);
        \coordinate (p2) at (4,4);
        \coordinate (p3) at (-1,10);

        \def \x{10}
        \coordinate (p4) at (\x,0);
        \coordinate (p5) at (\x + 7, 2);
        \coordinate (p6) at (\x + 6, 12);
        \coordinate (p7) at (\x + 2, 8);
        \coordinate (p8) at (\x - 5, 6);

        \def \w{20}
        \coordinate (p9) at (\w, 0);
        \coordinate (p10) at (\w+7,2);
        \coordinate (p11) at (\w+6, 12);
        \coordinate (p12) at (\w - 1, 10);

        \draw[purple] (p1) -- (p2);
        \draw[teal] (p2) -- (p3);
        \draw[violet] (p3) -- (p1);

        \draw[noamblue] (p4) -- (p5);
        \draw[gray] (p5) -- (p6);
        \draw[purple] (p6) -- (p7);
        \draw[orange] (p7) -- (p8);
        \draw[teal] (p8) -- (p4);

        \draw[orange] (p9) -- (p10);
        \draw[violet] (p10) -- (p11);
        \draw[noamblue] (p11) -- (p12);
        \draw[gray] (p12) -- (p9);


        \draw (p8) -- (p5) -- (p7);
        \draw (p10) -- (p12);
    \end{tikzpicture}
    \caption{An example of a $3$-set of castle polygons that is not a quadrangulation.}
    \label{fig:castle_polygons_1}
\end{figure}

Note that any quadrangulation gives an example of a $k$-set of castle polygons, but the construction here is more general. For instance, a $3$-set of castle polygons that is not a quadrangulation can be seen in Figure \ref{fig:castle_polygons_1}. The second condition in Definition \ref{def:castle_polygons} ensures that we can identify the upper and lower sides of the polygons by translations to obtain a translation surface. The last condition ensures that we obtain only one such surface, i.e., it ensures that the surface obtained by these identifications is connected. 

\begin{lemma}[Total cone angle of $k$-set of castle polygons]\label{lem:total_cone_angle_castle_polygons}
    Let $X$ be the translation surface obtained by identifying the top and bottom sides of a $k$-set of castle polygons. Then, the total cone angle of $X$ is $2\pi k$, i.e., if $X$ belongs to the stratum $\mathcal{H}(k_1 - 1, \ldots, k_n - 1)$, then 
    \begin{equation*}
        \sum_{j = 1}^n k_j = k. 
    \end{equation*}
\end{lemma}
\begin{proof}
    Note that the number of triangles $T_X$ in \emph{any} triangulation of a translation surface is fixed, depending only on the stratum. Indeed, by Euler's formula we have
    \begin{equation*}
        2 - 2\mathbf{g} = \#\mathcal{F} - \#\mathcal{E}  + \#\mathcal{V} = T_X - \frac{3}{2}T_X + n.
    \end{equation*}
    Solving for the number of triangles, we thus obtain
    \begin{equation*}
        T_X = 2 (2\mathbf{g} - 2 + n).  
    \end{equation*}
    We now claim that the number of base triangles is equal to the number of stack triangles. This follows from the fact that each singularity $j$ of $X$ has exactly $k_j$ incoming and outgoing verticals of the vertical geodesic flow. It is easy to see that each base triangle contains an outgoing ray, whereas each stack triangle contains an incoming ray. Thus, the number of both types of triangles in the triangulation is equal to
    \begin{equation*}
        k = 2\mathbf{g} - 2 + n.
    \end{equation*}
    The claim now follows from Proposition \ref{prop:genus_formula}.
\end{proof}

\begin{corollary}\label{cor:not_only_triangles}
    Not every polygon in a $k$-set of castle polygons is a triangle.
\end{corollary}
\begin{proof}
    The proof of Lemma \ref{lem:total_cone_angle_castle_polygons} shows that there are twice as many triangles as castle polygons.
\end{proof}

In a $k$-set of castle polygons, the vectors representing the sides of the polygons must satisfy the following \emph{train-track relations}. For $1 \leq i \leq k$, the sum of the wedges $w_{i, \ell} + w_{i, r}$ of the base triangle of polygon $i$ must be equal to the sum of all the upper sides of polygon $i$. We can encode the way the stack triangles are organized on top of the base triangle by using \emph{parenthesized} train-track equations, which are obtained by putting a set of parentheses around a sum of vectors that are part of the same stack triangle. We will always omit the outermost parentheses. For example, the parenthesized train-track relations for the $3$-set of castle polygons in Figure \ref{fig:castle_polygons_1} are given by
\begin{align}\label{eq:castle_train_tracks}
    \begin{split}
    w_{1, \ell} + w_{1, r} &= w_{2,\ell}, \\
    w_{2, \ell} + w_{2, r} &= w_{3, r} + (w_{1,r} + w_{3, \ell}), \\
    w_{3, \ell} + w_{3, r} &= w_{2, r} + w_{1, \ell}.
    \end{split}
\end{align}

The following lemma is immediate.
\begin{lemma}[$k$-set of castle polygons determined by wedges and train-tracks]
    A $k$-set of castle polygons is fully determined by its wedges and the parenthesized train-track equalities. 
\end{lemma}

Our goal is to describe the parameter space of $k$-sets of castle polygons in a similar fashion as we did for quadrangulations in section \ref{sec:diagonal_changes}. Recall that the parameter space of quadrangulations essentially consisted of two types of data, namely metric data describing the wedges and a combinatorial datum encoding the necessary identifications. The parameter space of $k$-sets of castle polygons will have the same general structure. Although parenthesized train-track relations are one option for encoding combinatorial data, in \cite{ferenczi2014generalization} Ferenczi proposes using \emph{castle forests} consisting of \emph{castle trees} for this purpose. These structures offer clearer insight into the effects of diagonal changes on the combinatorial data. We will take these trees and forests as our starting point in describing the combinatorics. 

\begin{definition}[Castle trees and castle forests]\label{def:castle_trees}
    The \emph{castle forest} of a $k$-set of castle polygons is a collection of $k$ \emph{castle trees}, which are connected directed acyclic graphs defined as follows.
    \begin{enumerate}[i)]
        \item There is one vertex on each upper side of a triangle in the triangulation, so one vertex for each base triangle and two vertices for each stack triangle. 
        \item The vertices of a stack triangle are labeled $\ell_i$ or $r_i$, if it is also an upper side of a polygon with upper side $w_{i, \ell}$ or $w_{i,r}$, respectively.
        \item For every vertex which is also on some lower side of a stack triangle, there is an upwards oriented edge going to the left or to the right vertex corresponding to the upper left or upper right side of the corresponding stack triangle.
        \item If polygon $i$ is a triangle, there is a vertex at the bottom of the triangle and a single upwards oriented edge from this vertex to the one corresponding to the upper side of the base triangle. 
    \end{enumerate}
\end{definition}

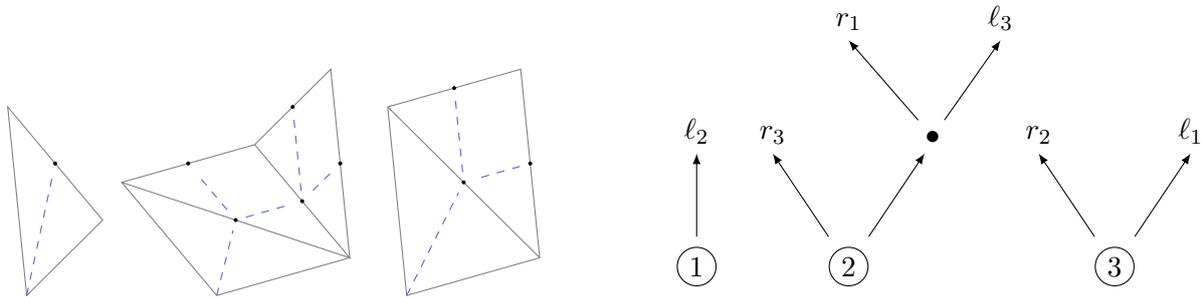
\begin{figure}[ht]
    \centering
    \begin{tikzpicture}[scale = 0.25]
        \coordinate (p1) at (0,0);
        \coordinate (p2) at (4,4);
        \coordinate (p3) at (-1,10);

        \def \x{10}
        \coordinate (p4) at (\x,0);
        \coordinate (p5) at (\x + 7, 2);
        \coordinate (p6) at (\x + 6, 12);
        \coordinate (p7) at (\x + 2, 8);
        \coordinate (p8) at (\x - 5, 6);

        \def \w{20}
        \coordinate (p9) at (\w, 0);
        \coordinate (p10) at (\w+7,2);
        \coordinate (p11) at (\w+6, 12);
        \coordinate (p12) at (\w - 1, 10);

        \draw[gray] (p1) -- (p2) -- node[midway, name = v1] {} (p3) -- cycle;
        
        \draw[gray] (p4) -- (p5) -- node[midway, name = w1] {} (p6) -- node[midway, name = w2] {} (p7) -- node[midway, name = w3] {} (p8) -- cycle;
        \draw[gray] (p8) -- node[midway, name = w4] {} (p5) -- node[midway, name = w5] {} (p7); 

        \draw[gray] (p9) -- (p10) -- node[midway, name = r1] {} (p11) --  node[midway, name=r2] {} (p12) -- cycle;
        \draw[gray] (p10) -- node[midway, name = r3] {} (p12);

        \foreach \p in {v1, w1, w2, w3, w4, w5, r1, r2, r3}{
            \fill (\p) circle (3pt);
        };

        \draw[dashed, noamblue] (p1) -- (v1);
        \draw[dashed, noamblue] (p4) -- (w4) -- (w5); \draw[dashed, noamblue] (w4) -- (w3);  \draw[dashed, noamblue] (w1) -- (w5) -- (w2);
        \draw[dashed, noamblue] (p9) -- (r3); \draw[dashed, noamblue] (r1) -- (r3) -- (r2);
    \end{tikzpicture} \qquad\qquad
    \begin{tikzpicture}
        \node (t1) at (0,0) {\circled{1}};
        \node (t2) at (2,0) {\circled{2}};
        \node (t3) at (5.5,0) {\circled{3}};

        \draw[->, >=latex] (t1) -- (0,1.5) node[above] {$\ell_2$};
        \draw[->, >=latex] (t2) -- (1,1.5) node[above] {$r_3$};
        \draw[->, >=latex] (t2) -- (3, 1.5) node[above, xshift = 3pt, name = i1] {$\bullet$};
        \draw[->, >=latex] (i1) -- (2,3) node[above] {$r_1$};
        \draw[->, >=latex] (i1) -- (4,3) node[above] {$\ell_3$};
        \draw[->, >=latex] (t3) -- (4.5,1.5) node[above] {$r_2$};
        \draw[->, >=latex] (t3) -- (6.5, 1.5) node[above] {$\ell_1$};
    \end{tikzpicture}
    \caption{The castle forest corresponding to the $3$-set of castle polygons from Figure \ref{fig:castle_polygons_1}. The left figure shows how to naturally embed the forest into the polygons.}
    \label{fig:castle_forest_1}
\end{figure}

Given a $k$-set of castle polygons, it is easy to deduce the corresponding castle forest by putting vertices on all of the edges of the stack triangles, putting vertices on the top edges of the base triangles as well as on the bottom of each base triangle. We then obtain the castle forest by adding oriented edges respecting the triangulation and deleting the lowest edge in case the polygon is not a triangle. Figure \ref{fig:castle_forest_1} contains an example of a castle forest corresponding to the 3-set of castle polygons from Figure \ref{fig:castle_polygons_1}.

Every vertex in a castle forest is either a \emph{leaf}, a \emph{node} or a \emph{root}. Instead of representing a castle forest graphically, we may use the following notation. A castle forest is completely determined by partitioning a word of length $2k$ on the letters $\{\varepsilon_i \mid \varepsilon \in \{\ell, r\}, i \in [k]\}$, where every letter appears exactly once, into $k$ parts and then fully parenthesizing these parts with respect to the nodes. Using this notation, the castle forest from Figure \ref{fig:castle_forest_1} is denoted by
\begin{equation}\label{eq:forest_notation}
    \boldsymbol{\pi} = \big(\ell_2\big) \big(r_3(r_1 \ell_3)\big) \big(r_2 \ell_1\big).
\end{equation}
Of course, this notation encodes exactly the train-track relations from \eqref{eq:castle_train_tracks}. It will be more convenient to state the action of a diagonal change on the combinatorial datum in this shorthand notation. But let us mention that the tree structure behind the notation is still valuable, e.g., for computer implementation, since by definition every castle tree is a \emph{full} binary tree which is a data structure which is very efficient to implement. To underline the analogy to the hyperelliptic case, we will write $\boldsymbol{\pi}$ for the representation of the castle forest as in \eqref{eq:castle_train_tracks} and we will refer to $\boldsymbol{\pi}$ either as the combinatorial datum or a \emph{forest word}. 

\subsubsection{Geometric Description}\label{sec:geometric_description}
We will now explain geometrically how to perform a diagonal change on a $k$-set of castle polygons. This description will be somewhat simpler than in the hyperelliptic case, since castle polygons are far more general than quadrilaterals, thus less care needs to be taken with respect to allowed moves. Let $P$ be a $k$-set of castle polygons. In every castle polygon which is not a triangle, there exists (at least) one diagonal, which is a line segment from the lowest vertex to a distinct higher vertex which is not a side of the polygon. We can apply a diagonal change in each of these polygons, which as in the hyperelliptic case we visualize as a cut and paste operation. The diagonal that is cut is always the one corresponding to the lowest vertex. As we will always assume that any polygon in the $k$-set of castle polygons has no vertices on the same horizontal or vertical, i.e., we assume that the corresponding translation surface satisfies Keane's condition, such a diagonal is uniquely defined. If this diagonal has positive slope, we say we apply a \emph{left move}. Analogously, if the diagonal has negative slope, we say we apply a \emph{right move}.

After cutting the polygon along this diagonal, which is not a vertical segment again by assuming Keane's condition is satisfied, we obtain two new polygons of which one is a castle polygon, i.e., it has a base triangle on the bottom with respect to the new naturally induced triangulation. The other part, which necessarily has a stack triangle as the lowest triangle in the new triangulation, is pasted on top of the side $\ell_i$ or $r_i$, depending on whether we applied a right or a left move. To clarify this procedure, we will proceed with a detailed example.

\begin{example}[Diagonal change in castle polygons -- geometric]\label{ex:castle_diagonal_change_geometric}
    Consider the $3$-set of castle polygons given by the following data. The wedges are given by
    \begin{align*}
        w_{1, \ell} &= (-2,2), \\
        w_{1, r} &= (1,1), \\
        w_{2, \ell} = w_{3, \ell} &= (-1.3,2), \\
        w_{2, r} = w_{3,r} &= (1.7, 1).
    \end{align*}
    The castle forest is defined by the forest word
    \begin{equation*}
        \boldsymbol{\pi}^{(0)} = \big(r_2\ell_2\big)\big(r_3\ell_3\big)\big(r_1\ell_1\big),
    \end{equation*}
    where we use the notation introduced in \eqref{eq:forest_notation}. It is straight-forward to verify that the train-track equalities are satisfied. Let us remark that this example comes from \cite{delecroix2015diagonal}, where it is given as an example of a quadrangulation that does not admit any well-slanted staircases, similar as Example \ref{ex:no_staircase}. This 3-set of castle polygons is pictured in Figure \ref{fig:castle_diagonal_change_01}. 

\begin{figure}[ht]
    \centering
    \begin{tikzpicture}[scale = 0.7]
        \draw (0,0) -- node[midway, below right] {$w_{1,r}$} (2,2) -- node[midway, above right] {$w_{2, \ell}$} (-0.5,5) -- node[midway, above] {$w_{2,r}$} (-4,4) -- node[midway, below left] {$w_{1,\ell}$} cycle;
        \draw[] (2,2) -- (-4,4);

        \def \x{5.5}
        \draw (\x, 0) -- node[midway, below, yshift = -2] {$w_{2,r}$} (\x+3.5, 1) -- node[midway, right] {$w_{3, \ell}$} (\x+1, 4) -- node[midway, above] {$w_{3, r}$} (\x-2.5, 3) -- node[midway, below left] {$w_{2, \ell}$} cycle;
        \draw[] (\x+3.5, 1) -- (\x-2.5, 3);

        \def \w{12}
        \draw (\w,0) -- node[midway, below, yshift = -2] {$w_{3, r}$} (\w + 3.5, 1) -- node[midway, above right] {$w_{1, \ell}$} (\w - 0.5, 5) -- node[midway, above left] {$w_{1, r}$} (\w - 2.5, 3) -- node[midway, below left] {$w_{3, \ell}$} cycle; 
        \draw[] (\w+3.5, 1) -- (\w-2.5, 3); 

        \draw[dashed, purple] (0,0) -- (-0.5,5);
        \node[rotate =100, xshift = 7] at (-0.5,5) {\Large \Leftscissors};
        \fill[gray, opacity = 0.1] (0,0) -- (-0.5,5) -- (-4,4) -- cycle;
        
    \end{tikzpicture}
    \caption{The initial 3-set of castle polygons from Example \ref{ex:castle_diagonal_change_geometric}. The first diagonal change move is already indicated.}
    \label{fig:castle_diagonal_change_01}
\end{figure}
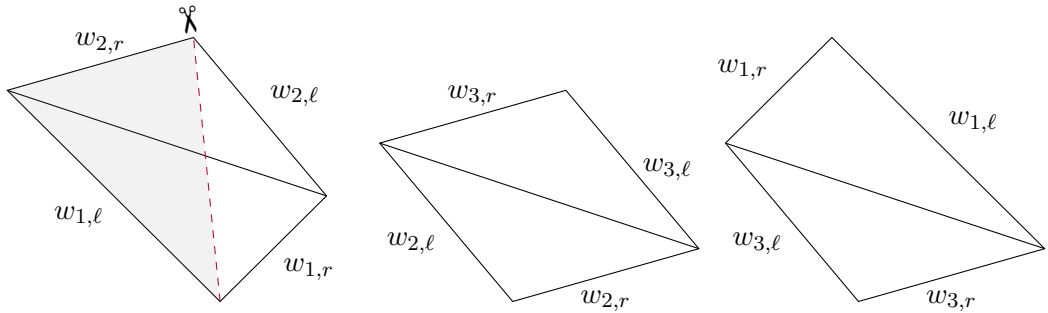

    Since none of the polygons are a triangle, we can apply a diagonal change in each of them. We choose to perform a diagonal change in the first polygon. As described above, we cut this quadrilateral along its diagonal which connects the bottom vertex to the top vertex. This results in two triangles. The right one is a base triangle, whereas the left one is a stack triangle. We glue the stack triangle by its lower side onto the top right side of the right-most polygon. This diagonal change is depicted in Figure \ref{fig:castle_diagonal_change_02}. From the new 3-set of castle polygons, we can read off the corresponding forest word, which is given by
    \begin{equation*}
        \boldsymbol{\pi}^{(1)} = \big(\ell_2\big)\big(r_3\ell_3\big)\big(r_1(r_2\ell_1)\big). 
    \end{equation*}

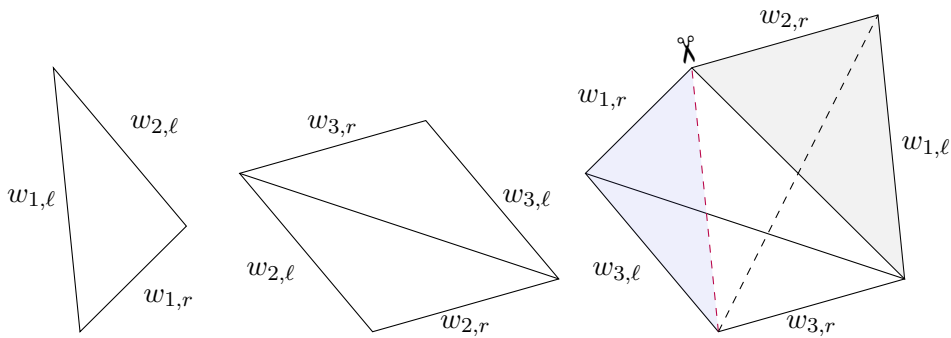
\begin{figure}[hb]
    \centering
    \begin{tikzpicture}[scale = 0.7]
        \draw (0,0) -- node[midway, below right] {$w_{1,r}$} (2,2) -- node[midway, above right] {$w_{2, \ell}$} (-0.5, 5) -- node [midway, left] {$w_{1, \ell}$} cycle;

        \def \x{5.5}
        \draw (\x, 0) -- node[midway, below] {$w_{2,r}$} (\x+3.5, 1) -- node[midway, right] {$w_{3, \ell}$} (\x+1, 4) -- node[midway, above] {$w_{3, r}$} (\x-2.5, 3) -- node[midway, below left] {$w_{2, \ell}$} cycle;
        \draw[] (\x+3.5, 1) -- (\x-2.5, 3);

        \def \w{12}
        \draw (\w,0) -- node[midway, below] {$w_{3, r}$} (\w + 3.5, 1) -- (\w - 0.5, 5) -- node[midway, above left] {$w_{1, r}$} (\w - 2.5, 3) -- node[midway, below left] {$w_{3, \ell}$} cycle; 
        \draw[] (\w-2.5, 3) -- (\w+3.5, 1) -- node[midway, right] {$w_{1, \ell}$}  (\w+3, 6) -- node[midway, above] {$w_{2, r}$} (\w-0.5,5);

        \fill[gray, opacity = 0.1] (\w + 3.5, 1) -- (\w + 3, 6) -- (\w - 0.5, 5) -- cycle;

        \draw[dashed, purple] (\w,0) -- (\w-0.5,5);
        \draw[dashed, black]  (\w,0) -- (\w+3, 6); 
        \node[rotate =100, xshift = 7] at (\w-0.5,5) {\Large \Leftscissors};
        \fill[noamblue, opacity = 0.1] (\w,0) -- (\w-0.5, 5) -- (\w-2.5, 3) -- cycle; 
        
    \end{tikzpicture}
    \caption{The 3-set of castle polygons after applying one diagonal change move. The gray triangle was formerly a part of polygon 1, as in Figure \ref{fig:castle_diagonal_change_01}.}
    \label{fig:castle_diagonal_change_02}
\end{figure}
Continuing, we may only apply a diagonal change to the polygons 2 and 3 since the first one is a triangle. Suppose we choose to apply a diagonal change to the pentagon on the right. There are two diagonals, indicated by the dashed lines in Figure \ref{fig:castle_diagonal_change_02}. By definition, the diagonal change move cuts the polygon along the diagonal whose top vertex is lower, so in our case here it is the left diagonal colored red. This decomposes the pentagon into a quadrilateral, which will become the new polygon 3, and a triangle which is pasted on top of polygon 2 as in Figure \ref{fig:castle_diagonal_change_03}. The new forest word is given by
\begin{equation*}
    \boldsymbol{\pi}^{(2)} = \big(\ell_2\big) \big(r_3(r_1\ell_3)\big) \big(r_2\ell_1\big). 
\end{equation*}
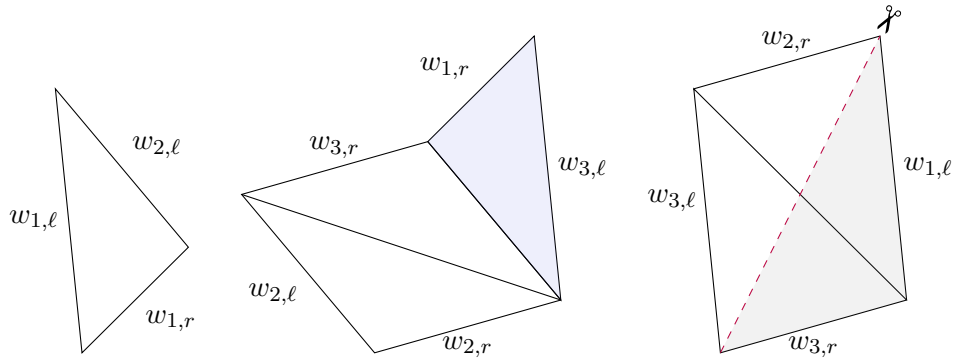
\begin{figure}[!ht]
    \centering
    \begin{tikzpicture}[scale = 0.7]
        \draw (0,0) -- node[midway, below right] {$w_{1,r}$} (2,2) -- node[midway, above right] {$w_{2, \ell}$} (-0.5, 5) -- node [midway, left] {$w_{1, \ell}$} cycle;

        \def \x{5.5}
        \draw (\x, 0) -- node[midway, below] {$w_{2,r}$} (\x+3.5, 1) -- (\x+1, 4) -- node[midway, above] {$w_{3, r}$} (\x-2.5, 3) -- node[midway, below left] {$w_{2, \ell}$} cycle;
        \draw (\x+3.5,1) -- node[midway, right] {$w_{3, \ell}$} (\x+3,6) -- node[midway, above left] {$w_{1, r}$} (\x+1, 4) -- cycle;
        \fill[noamblue, opacity = 0.1] (\x+3.5,1) -- (\x+3,6) -- (\x+1, 4) -- cycle;

        \draw[] (\x+3.5, 1) -- (\x-2.5, 3);

        \def \w{12}
        \draw (\w,0) -- node[midway, below] {$w_{3, r}$} (\w + 3.5, 1) -- (\w - 0.5, 5) -- node[midway, above left] {$w_{3, \ell}$}  cycle; 
        \draw[] (\w+3.5, 1) -- node[midway, right] {$w_{1, \ell}$}  (\w+3, 6) -- node[midway, above] {$w_{2, r}$} (\w-0.5,5);
        \draw[dashed, purple] (\w, 0) -- (\w+3, 6); 
        \node[rotate =60, xshift = 7] at (\w+3, 6) {\Large \Leftscissors};
        \fill[gray, opacity = 0.1] (\w,0) -- (\w+3.5, 1) -- (\w+3, 6) -- cycle;
        
    \end{tikzpicture}
    \caption{The 3-set of castle polygons after applying two diagonal changes. The blue triangle was formerly attached to the polygon to the right as in Figure \ref{fig:castle_diagonal_change_02}.}
    \label{fig:castle_diagonal_change_03}
\end{figure}Note that we are still not able to perform a diagonal change in the first polygon, since we have not pasted anything to it in the last step. Thus, our choice is still restricted to apply a move to polygons 2 or 3. Figure \ref{fig:castle_diagonal_change_04} illustrates the outcome if we choose to apply another move to polygon 3. The corresponding forest word is given by
\begin{equation*}
    \boldsymbol{\pi}^{(3)} = \big(\ell_2\big) \big((r_3\ell_1)(r_1\ell_3\big) \big(r_2\big).
\end{equation*}

In the resulting 3-set of castle polygons, we can now only apply a move to the polygon in the middle. Continuing the algorithm, one can verify that after another move we obtain a 3-set of castle polygons consisting of a triangle to the left, a quadrilateral in the middle and a polygon to the right. 

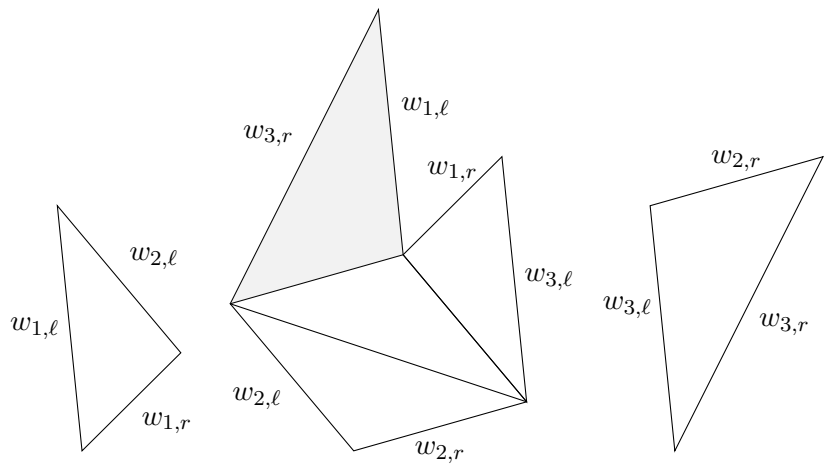
\begin{figure}[!hb]
    \centering
    \begin{tikzpicture}[scale = 0.65]
        \draw (0,0) -- node[midway, below right] {$w_{1,r}$} (2,2) -- node[midway, above right] {$w_{2, \ell}$} (-0.5, 5) -- node [midway, left] {$w_{1, \ell}$} cycle;

        \def \x{5.5}
        \draw (\x, 0) -- node[midway, below, yshift=-2] {$w_{2,r}$} (\x+3.5, 1) -- (\x+1, 4) -- (\x-2.5, 3) -- node[midway, below left] {$w_{2, \ell}$} cycle;
        \draw (\x+3.5,1) -- node[midway, right] {$w_{3, \ell}$} (\x+3,6) -- node[midway, above, yshift = 5] {$w_{1, r}$} (\x+1, 4) -- cycle;
        \draw (\x + 1, 4) -- node[midway, right, yshift = 9] {$w_{1, \ell}$} (\x+0.5,9) -- node[midway, above left] {$w_{3,r}$} (\x-2.5,3); 
        \fill[gray, opacity = 0.1] (\x + 1, 4) -- (\x+0.5, 9) -- (\x-2.5, 3) -- cycle;

        \draw[] (\x+3.5, 1) -- (\x-2.5, 3);

        \def \w{12}
        \draw (\w,0) -- node[midway, below right] {$w_{3, r}$} (\w+3, 6) -- node[midway, above] {$w_{2,r}$}  (\w - 0.5, 5) -- node[midway, above left] {$w_{3, \ell}$}  cycle; 
    \end{tikzpicture}
    \caption{The 3-set of castle polygons after applying three diagonal changes. The gray triangle was formerly attached to the polygon on the right as in Figure \ref{fig:castle_diagonal_change_03}.}
    \label{fig:castle_diagonal_change_04}
\end{figure}
\end{example}

\subsubsection{Action on Data}\label{sec:action_on_data}
As we did in section \ref{sec:staircase_moves} and section \ref{sec:parameter_space}, we will now develop the action of a diagonal change move described geometrically in the previous section on the length and combinatorial data describing a $k$-set of castle polygons. Let $P$ be a $k$-set of castle polygons, which is described by its wedges
\begin{equation*}
    \boldsymbol{w} = \big((w_{1, \ell}, w_{1, r}), \ldots, (w_{k, \ell}, w_{k, r})\big)
\end{equation*}
and its forest word $\boldsymbol{\pi}$ as in \eqref{eq:forest_notation}, so that we write $P = (\boldsymbol{\pi}, \boldsymbol{w})$. 

\begin{definition}[Allowed moves]\label{def:allowed_moves}
    Let $P = (\boldsymbol{\pi}, \boldsymbol{w})$ be a $k$-set of castle polygons. We define $H(P)$ to be the subset of $[k]$ such that $i \in H(P)$ if and only if polygon $i$ is not a triangle. The set $H(P)$ is the \emph{set of allowed moves} in $P$. 
\end{definition}
\begin{remark}
    Corollary \ref{cor:not_only_triangles} ensures that $H(P)$ is nonempty for every $k$-set of castle polygons.
\end{remark}
The set $H(P)$ contains the information in which polygons we may apply a diagonal change. 
\begin{definition}[Choice of diagonal change move]\label{def:choices_of_diagonal_change_moves}
A \emph{choice} of a diagonal change move is an application 
\begin{align*}
    c\colon \mathcal{P}([k])\setminus \{\emptyset\} &\to [k] \times \{\ell, r\}, \\
    H(P) &\mapsto (i, \varepsilon),
\end{align*}
which simply chooses an element $i$ of $H(P)$ and adds the information whether a left or right move is applied. We will use the following notation, which is an adaptation of the notation used for cycles in \ref{sec:parameter_space}. We will write $c$ as a word of length $k$ on the alphabet $\{\cdot, \ell\}$ if the applied move is a left move or the alphabet $\{\cdot, r\}$ if we apply a right move, where the $i$\textsuperscript{th} letter of the word is $\ell$ (or $r$) if and only if $c(H(P)) = (i, \varepsilon)$. 
\end{definition}

\begin{example}[Choices of diagonal change moves]
    In Example \ref{ex:castle_diagonal_change_geometric}, we described three diagonal changes which required us to make three choices. Using the notation from Definition \ref{def:choices_of_diagonal_change_moves}, the choices we made can be written as
    \begin{align*}
        c^{(1)} &= r \cdot \cdot\quad, \\
        c^{(2)} &= \cdot \cdot r\quad, \\
        c^{(3)} &= \cdot \cdot \ell\quad.
    \end{align*}
\end{example}

Our goal now is to show that applying a diagonal change move to a $k$-set of castle polygons produces a new $k$-set of castle polygons, and moreover we want to investigate closer the action on the two types of data. The following lemma is taken from \cite{ferenczi2015diagonal}, where we have adapted the language to fit into our context.

\begin{lemma}[Diagonal change produces castle polygons]\label{lem:diagonal_change produces_castle_polygons}
    The image of a $k$-set of castle polygons by a (forward) diagonal change move given by a decision $c$ is again a $k$-set of castle polygons. 
\end{lemma}
\begin{proof}
    We only need to check that the part that gets cut away from some polygon $p_i$ does not need to get pasted on itself, since in this case the operation as outlined above would not be defined. Towards a contradiction, suppose that this was the case and suppose that $w_{i, \varepsilon}$ is the wedge of polygon $p_i$ which has been cut away. Now, $w_{i, \varepsilon}$ would be an upper side of this polygon too. Write $u_\ell$ for the left end-point of this upper side and $u_r$ for the right endpoint. Since there are no vertical saddle connections, $u_\ell$ and $u_r$ cannot lie directly above the left and right endpoints of the wedge. Suppose that the diagonal change applied was a left move (the situation corresponding to a right move is analogous). Then, the projection of $u_\ell$ to the $x$-axis is to the right of the projection of the lowest vertex of $p_i$. But then, the projection of $u_r$ is to the right of the projection of $w_{i, r}$ as well, which contradicts the fact that $p_i$ is a castle polygon. Indeed, this can only happen if the polygon contains a second base triangle or a triangle that is neither a base nor a stack triangle. 
\end{proof}

For each wedge $w_i = (w_{i, \ell}, w_{i,r})$ corresponding to polygon $p_i \in P$, if $p_i$ is not a triangle we denote by $w_{i, d}$ the (forward) diagonal of $p_i$, which is defined as the diagonal along which we would cut the polygon when applying a diagonal change move. More precisely, if we define
\begin{equation}\label{eq:leaf_sums}
    v_{i, \ell} = \sum_{\substack{\text{left of}\\\text{root } i}} w_{j, r} -  w_{j, \ell}, \qquad v_{i, r} = \sum_{\substack{\text{right of}\\\text{root } i}} w_{j, \ell} - w_{j,r}, 
\end{equation}
then $w_{i, d} = w_{i, r} + v_{i, r}$ if we apply a left move, and $w_{i, d} = w_{i, \ell} + v_{i, \ell}$ if we apply a right move. The sums in \eqref{eq:leaf_sums} are taken over all wedges that correspond to labels in the castle tree that are to the left of the root of tree $i$ or to the right of the root of tree $i$, respectively. Alternatively, we can view these sums as the sums of the upper sides on the left or the right of the forward diagonal. This describes the action on the length datum. 

For the action on the combinatorial datum, let us briefly return to the graphical representations of castle forests and castle trees. This is the language in which Ferenczi describes the action on the combinatorics in \cite{ferenczi2014generalization}. For obvious reasons, the action on the combinatorial datum is referred to as \emph{surgery on trees}.

\begin{proposition}[Surgery on trees]\label{prop:tree_surgery}
    Let $P = (\boldsymbol{\pi}, \boldsymbol{w})$ be a $k$-set of castle polygons and let $c$ be a decision with $c(H(P)) = (i, \varepsilon)$ for $\varepsilon\in\{\ell, r\}$. Applying the associated diagonal change amounts to the action on the castle forest $\boldsymbol{\pi}$ as outlined below.

    Suppose that $c(H(P)) = (i, \ell)$ and write $n$ for the node or leaf at the end of the left edge out of the root of the tree associated to the polygon $p_i$. We distinguish two cases.
    \begin{enumerate}
        \item If $n$ is a \emph{node}, then the part of (former) tree associated to $p_i$ beyond $n$ becomes the new tree associated to $p_i$, including $n$ which becomes its new root.
        
        The rest of the tree, meaning a copy of the node $n$ and everything to the right of the root, is put on top of the leaf labeled $r_i$ which is part of some tree associated to $p_j$. If $j = i$, then we mean the \emph{new} tree associated to $p_i$.

        The copy of $n$ that has not become the root of a new tree becomes a leaf labeled $r_i$. Former leaf $r_i$ becomes a node.
        \item If $n$ is a \emph{leaf} labeled $r_j$, then a new tree associated to $p_i$ is made with a single edge and a leaf labeled $r_j$. The whole former tree associated to $p_i$ is put on top of leaf $r_i$ in its tree which is associated to $p_k$. In this tree, former leaf $r_i$ becomes a node and former leaf $r_j$ (which belongs to the part that was added) is labeled as $r_i$.
    \end{enumerate}
    After any of these two steps, we delete the \enquote{stem} of the tree if necessary, meaning that if the new tree does not correspond to a triangle we make sure that the root has two nodes by removing the root and its single edge otherwise.

    If $c(H(P)) = (i,r)$, we may apply the same procedure by exchanging left and right as well as $\ell$ and $r$. 
\end{proposition}
\begin{proof}
    This follows at once from the description of a diagonal change move given in section \ref{sec:geometric_description}.
\end{proof}

As we have mentioned above, such a surgery on trees can be implemented very efficiently with computer algorithms. Translating this procedure to \emph{forest words} gives an easier description if one wishes to apply diagonal changes by hand.

\begin{corollary}[Surgery on forest words]\label{cor:action_forest_words}
    The procedure outlined in Proposition \ref{prop:tree_surgery} can be translated to the setting of forest words as follows. 
    \begin{itemize}
        \item If $c(H(P)) = (i, \ell)$, i.e., if we apply a \emph{left} move, then we adjust the forest word as follows.
        \begin{enumerate}
            \item Keep the left part of partition $i$, i.e., the part that is associated to the left of the root.
            \item Glue the right part of partition $i$ to the right of $r_i$.
        \end{enumerate}
        \item If $c(H(P)) = (i,r)$, i.e., if we apply a \emph{right} move, then we adjust the forest word as follows.
        \begin{enumerate}
            \item Keep the right part of partition $i$.
            \item Glue the left part of partition $i$ to the left of $\ell_i$.
        \end{enumerate}
    \end{itemize}
\end{corollary}

We will exemplify this procedure by revisiting the diagonal changes outlined in Example \ref{ex:castle_diagonal_change_geometric}.

\begin{example}[Action on forest words]\label{ex:forest_words}
    Recall that the initial combinatorial datum in Example \ref{ex:castle_diagonal_change_geometric} was given by
    \begin{equation}\label{eq:forest_words_01}
        \boldsymbol{\pi}^{(0)} = \combi{\textcolor{BrickRed}{r_2}\textcolor{MidnightBlue}{\ell_2}} \combi{r_3\ell_3} \combi{r_1\ell_1}. 
    \end{equation}
    The first move corresponded to the choice $c(H(P)) = (1, r)$. Following the outline from Corollary \ref{cor:action_forest_words}, we thus keep the right part in partition 1, which we have colored blue in \eqref{eq:forest_words_01}, and we glue the red left part to the left of $\ell_1$. This results in the combinatorial datum
    \begin{equation*}
        \boldsymbol{\pi}^{(1)} = \combi{\textcolor{MidnightBlue}{\ell_2}} \combi{r_3\ell_3} \combi{r_1(\textcolor{BrickRed}{r_2}\ell_1)},
    \end{equation*}
    which coincides with the forest word from Example \ref{ex:castle_diagonal_change_geometric}. After three diagonal changes we arrived at a castle polygon with forest word
    \begin{equation*}
        \boldsymbol{\pi}^{(3)} = \combi{\ell_2}\combi{\textcolor{BrickRed}{(r_3\ell_1)}\textcolor{MidnightBlue}{(r_1\ell_3)}}\combi{r_2},
    \end{equation*}
    and in view of Figure \ref{fig:castle_diagonal_change_04} it is clear that the only possible move is given by the choice $c(H(P)) = (2, \ell)$. Therefore, the forest word for the next iterate would be given by
    \begin{equation*}
        \boldsymbol{\pi}^{(4)} = \combi{\ell_2}\combi{\textcolor{BrickRed}{r_3\ell_1}}\combi{r_2\textcolor{MidnightBlue}{(r_1\ell_3)}}.
    \end{equation*}
\end{example}

Since the action on the combinatorial datum $\boldsymbol{\pi}$ is independent of the length datum $\boldsymbol{w}$ and the action on $\boldsymbol{w}$ is linear, we may again describe the action on the length datum in the form of a matrix $A_{\boldsymbol{\pi},c}$, completely analogous to when we described staircase moves in this way. We will use the same notation, which we recall now. The rows and columns of the matrix $A_{\boldsymbol{\pi}, c}$ are indexed by
\begin{equation*}
    (1, \ell), (1, r), \ldots, (k, \ell), (k, r),
\end{equation*}
and we write $E_{(i,\varepsilon)(j, \nu)}$ for the elementary matrix that has a 1 in position $\combi{(i, \varepsilon),(j,\nu)}$ and zeroes everywhere else. Moreover, we denote by $I_{2k}$ the $2k \times 2k$ identity matrix.

\begin{definition}[Diagonal change matrix]\label{def:diagonal_change_matrix}
    Given a $k$-set of castle polygons $P = (\boldsymbol{\pi}, \boldsymbol{w})$ and a choice $c(H(P)) = (i, \varepsilon)$ of a diagonal change move, we define the matrix $A_{\boldsymbol{\pi},c}$ to be
    \begin{equation*}
        A_{\boldsymbol{\pi},c} \coloneqq
        \begin{cases}
            I_{2k} + \displaystyle\sum_{\substack{\nu_j\text{ left}\\\text{of root}}} \sgn(\nu)\cdot E_{(i, \ell)(j, \nu)} \quad &\text{if } \varepsilon = r, \\
            I_{2k} + \displaystyle\sum_{\substack{\nu_j\text{ right}\\\text{of root}}} -\sgn(\nu)\cdot E_{(i, r)(j, \nu)} &\text{if } \varepsilon = \ell,
        \end{cases}
    \end{equation*}
    where we define $\sgn(\nu)$ as
    \begin{equation*}
        \sgn(\nu) = \begin{cases}
            1 \quad &\text{if } \nu = r, \\
            -1 &\text{if } \nu = \ell.
        \end{cases}
    \end{equation*}
    Writing $\boldsymbol{w}'$ for the length datum after the diagonal change and viewing both $\boldsymbol{w}$ and $\boldsymbol{w}'$ as column vectors, it follows from \eqref{eq:leaf_sums} that
    \begin{equation*}
        \boldsymbol{w}' = A_{\boldsymbol{\pi}, c} \cdot \boldsymbol{w}.
    \end{equation*}
\end{definition}

\begin{example}[Diagonal change matrix]
    We illustrate Definition \ref{def:diagonal_change_matrix} by giving an explicit example relating to Example \ref{ex:castle_diagonal_change_geometric}. We will consider the $3$-set of castle polygons obtained after three diagonal changes seen in Figure \ref{fig:castle_diagonal_change_04}. Under the usual identification $\C \simeq \R^2$, we can write the wedges as a vector in $\C^6$ as
    \begin{equation*}
        \boldsymbol{w} = (-0.3+3\ii, 1+\ii, -1.3+2\ii, 1.7+\ii, -0.3+3\ii, 1.4+4\ii). 
    \end{equation*}
    As we explained in Example \ref{ex:forest_words}, the only choice is given by $c(H(P)) = (2, \ell)$. The corresponding matrix is therefore given by
    \begin{equation*}
        A_{\boldsymbol{\pi}, c} = 
        \begin{bmatrix}
            1 & 0 & 0 & 0 & 0 & 0 \\
            0 & 1 & 0 & 0 & 0 & 0 \\
            0 & 0 & 1 & 0 & 0 & 0 \\
            0 &-1 & 0 & 1 & 1 & 0 \\
            0 & 0 & 0 & 0 & 1 & 0 \\
            0 & 0 & 0 & 0 & 0 & 1
        \end{bmatrix}.
    \end{equation*}
    Therefore, the new length datum $\boldsymbol{w}'$ after applying a further diagonal change in polygon 2 is given by
    \begin{equation*}
        \boldsymbol{w}' = A_{\boldsymbol{\pi}, c} \cdot \boldsymbol{w} = 
        \begin{bmatrix}
            1 & 0 & 0 & 0 & 0 & 0 \\
            0 & 1 & 0 & 0 & 0 & 0 \\
            0 & 0 & 1 & 0 & 0 & 0 \\
            0 &-1 & 0 & 1 & 1 & 0 \\
            0 & 0 & 0 & 0 & 1 & 0 \\
            0 & 0 & 0 & 0 & 0 & 1
        \end{bmatrix}
        \cdot
        \begin{bmatrix}
            -0.3+3\ii \\
            1+\ii \\
            -1.3 + 2\ii \\
            1.7 + \ii \\
            -0.3 + 3\ii \\
            1.4 + 4\ii
        \end{bmatrix} = 
        \begin{bmatrix}
            -0.3 + 3\ii \\
            1 + \ii \\
            -1.3 + 2\ii \\
            0.4 + 3\ii \\
            -0.3 + 3\ii \\
            1.4 + 4\ii
        \end{bmatrix}. 
    \end{equation*}
    One can verify this result geometrically by examining Figure \ref{fig:castle_diagonal_change_04}, where $w_{2, r} = 0.4 + 3\ii$, the fourth coordinate of $\boldsymbol{w}'$, is exactly the forward diagonal of the middle polygon.
\end{example}

Writing $\boldsymbol{\pi}' = c \cdot \boldsymbol{\pi}$ to denote the action on the combinatorial datum from Proposition \ref{prop:tree_surgery} or Corollary \ref{cor:action_forest_words}, we have established the following.

\begin{lemma}[Diagonal change move on data]
    Given a $k$-set of castle polygons $P = (\boldsymbol{\pi}, \boldsymbol{w})$ and a choice $c$, performing the associated staircase move yields a new $k$-set of castle polygons $P' = (\boldsymbol{\pi}', \boldsymbol{w}')$, where
    \begin{equation*}
        \boldsymbol{\pi}' = c \cdot \boldsymbol{\pi}, \quad \boldsymbol{w}' = A_{\boldsymbol{\pi}, c} \cdot \boldsymbol{w}
    \end{equation*}
    are the actions on the data described in Proposition \ref{prop:tree_surgery}, Corollary \ref{cor:action_forest_words} and Definition \ref{def:diagonal_change_matrix}.
\end{lemma}

Before we look closer at the parameter space associated to $k$-sets of castle polygons, let us mention that analogous to the hyperelliptic case it is equally possible to use the description of the action given above to define several different diagonal changes algorithms. In our discussion above, we always assumed that only one diagonal change move is performed at a time. For instance, we could specify our algorithm by requiring to always perform the diagonal change with respect to the left-most polygon where a move is possible. In other words, we apply a diagonal change to the polygon $p_i$ such that $i$ is minimized.

Another choice would be a greedy algorithm, which is defined analogous to the algorithm of the same name in Definition \ref{def:greedy_algorithm}. Namely, at each step we perform all possible diagonal changes (in any order). We will refer to these types of algorithms as \emph{simultaneous diagonal changes}, as opposed to accelerated variants where we may apply several moves in successively constructed $k$-sets of castle polygons. 

Moreover, the following lemma carries over from the hyperelliptic case, see Lemma \ref{lem:diagonals_become_sides}. The proof can be found in \cite{ferenczi2015diagonal}, where it appears as Theorem 1.5.

\begin{lemma}\label{lem:diagonals_become_sides_general}
    Let $P$ be a $k$-set of castle polygons. There exists an infinite sequence of forward diagonal change moves starting from $P$ such that the widths of the polygons tends to zero and the height of the polygons tends to infinity if and only if the translation surface $X(P)$ associated to $P$ has no vertical saddle connections. Moreover, in any such sequence there are infinitely many left and right moves. 
\end{lemma}
\begin{remark}
    An immediate consequence of Lemma \ref{lem:diagonals_become_sides_general} is that any forward diagonal eventually becomes part of a wedge under applications of simultaneous diagonal changes.
\end{remark}

\subsubsection{Parameter Space}

Our goal is to introduce the space of (labeled) $k$-sets of castle polygons on which the algorithms will act. We will follow closely section \ref{sec:parameter_space}, where we developed the same theory in the hyperelliptic case. 

We start by transporting the concept of a diagonal change class to this setting. This classes will depend on the combinatorial datum, which here are given by castle forests or forest words. In Definition \ref{def:castle_trees} we defined these objects as being induced by a $k$-set of castle polygons. We may define castle forests, and by extension forest words as well, independently without mentioning castle polygons as follows.

\begin{definition}[Castle trees and castle forests]\label{def:castle_forest_02}
    A \emph{castle forest} $\boldsymbol{\pi}$ is the disjoint union of $k$ \emph{castle trees}, which are connected directed acyclic graphs numbered from 1 to $k$ satisfying the following.
    \begin{enumerate}[i)]
        \item All vertices are roots, nodes or leaves, where a root is also a node exactly if the tree has more than two vertices. Each tree has exactly one root.
        \item There are $2k$ leaves labeled by $\ell_i$ and $r_i$, where $1\leq i \leq k$.
        \item In a tree with nodes, a leaf $l_i$ (or $r_i$) is at the end of a right (or left) edge. 
        \item There is no strict subset $I$ of $[k]$ such that the trees $i$ for $i \in I$ have only leaves $l_j$ or $r_j$ for $j \in I$.
    \end{enumerate}
\end{definition}

Starting with a castle forest $\boldsymbol{\pi}$ as given by Definition \ref{def:castle_forest_02}, we define
\begin{equation*}
    \Omega_{\boldsymbol{\pi}} \coloneqq \left\{\boldsymbol{w} = \big((w_{1, \ell}, w_{1, r}), \ldots, (w_{k,\ell}, w_{k, r})
    \big) \mid \boldsymbol{w} \text{ satisfies the train-track equalities induced by }F\right\}.
\end{equation*}
In words, $\Omega_{\boldsymbol{\pi}}$ consists of all collections of wedges that yield a $k$-set of castle polygons with castle forest $F$. We have that
\begin{equation*}
    \Omega_{\boldsymbol{\pi}} \subseteq \Big(\big(\R_- \times \R_+\big) \times \big(\R_+ \times \R_+\big)\Big)^k.
\end{equation*}
We say that a castle forest $\boldsymbol{\pi}$ is \emph{admissible}, if $\Omega_{\boldsymbol{\pi}}$ is nonempty. 

\begin{definition}[Diagonal change classes for the general case]
    A collection $\mathfrak{DC}$ of admissible castle forests $\boldsymbol{\pi}^j$ is called a \emph{diagonal change class}, or $DC$ \emph{class} for short, if for any $\boldsymbol{\pi}^i \neq \boldsymbol{\pi}^j$ there exists a finite sequence of diagonal change moves from a $k$-set of castle polygons with combinatorial datum $\boldsymbol{\pi}^i$ to a $k$-set of castle polygons with combinatorial datum $\boldsymbol{\pi}^j$.
\end{definition}

\begin{figure}
    \centering
\begin{tikzcd}[cramped,column sep=tiny,row sep=large]
	&& {\big(r_2\big)\big((r_1\ell_2)\ell_1\big)} &&&& {\big(\ell_2\big)\big(r_1(r_2\ell_1)\big)} \\
	\\
	&& {\big(r_2\big)\big((r_1\ell_1)\ell_2\big)}\arrow[loop above, distance = 1cm, in=120, out=60, "\cdot \,r" {pos=0.7}] &&&& {\big(\ell_2\big)\big(r_2(r_1\ell_1)\big)} \arrow[loop above, distance = 1cm, in=120, out=60, "\cdot \,\ell" {pos=0.7}]\\
	{\big(r_2\ell_1\big)\big(r_1\ell_2\big)} \arrow[loop, distance = 1.5cm, in=240, out=200, "r\,\cdot" {pos=0.15}] \arrow[loop, distance = 1.5cm, in=120, out=160, "r\,\cdot"' {pos=0.15}]&&&& {\big(r_2\ell_2\big)\big(r_1\ell_1\big)} &&&& {\big(r_1\ell_2\big)\big(r_2\ell_1\big)} \arrow[loop, distance = 1.5cm, in=60, out=20, "\ell\,\cdot" {pos=0.15}] \arrow[loop, distance = 1.5cm, in=300, out=340, "\cdot \, \ell"' {pos=0.15}]\\
	&& {\big((r_2\ell_2)\ell_1\big)\big(r_1\big)} \arrow[loop, distance = 1cm, in=300, out=240, "r\,\cdot" {pos=0.3}]&&&& {\big(r_1(r_2\ell_2)\big)\big(\ell_1\big)} \arrow[loop, distance = 1cm, in=300, out=240, "\ell\,\cdot" {pos=0.3}]\\
	\\
	&& {\big((r_2\ell_1)\ell_2\big)\big(r_1\big)} &&&& {\big(r_2(r_1\ell_2\big)\big(\ell_1\big)}
	\arrow["{r \cdot}", from=4-5, to=1-7]
	\arrow["{\cdot r}", from=3-7, to=4-5]
	\arrow["{\cdot r}", from=1-7, to=4-9]
	\arrow["r\cdot", from=4-9, to=3-7]
	\arrow["{\cdot r}"', from=4-5, to=7-7]
	\arrow["{r \cdot}"', from=7-7, to=4-9]
	\arrow["{\cdot r}"', from=4-9, to=5-7]
	\arrow["r\cdot"', from=5-7, to=4-5]
	\arrow["\ell\cdot"', from=4-1, to=3-3]
	\arrow["{\cdot \ell}"', from=3-3, to=4-5]
	\arrow["{\cdot \ell}", from=4-1, to=5-3]
	\arrow["\ell\cdot", from=5-3, to=4-5]
	\arrow["\ell\cdot"', from=4-5, to=1-3]
	\arrow["{\cdot \ell}"', from=1-3, to=4-1]
	\arrow["\ell\cdot", from=7-3, to=4-1]
	\arrow["{\cdot \ell}", from=4-5, to=7-3]
\end{tikzcd}
    \caption{The $DC$ graph $\mathcal{G}$ associated to the castle forest $\boldsymbol{\pi} = \big(r_2\ell_2\big)\big(r_1\ell_1\big)$.}
    \label{fig:dc_graph_general}
\end{figure}
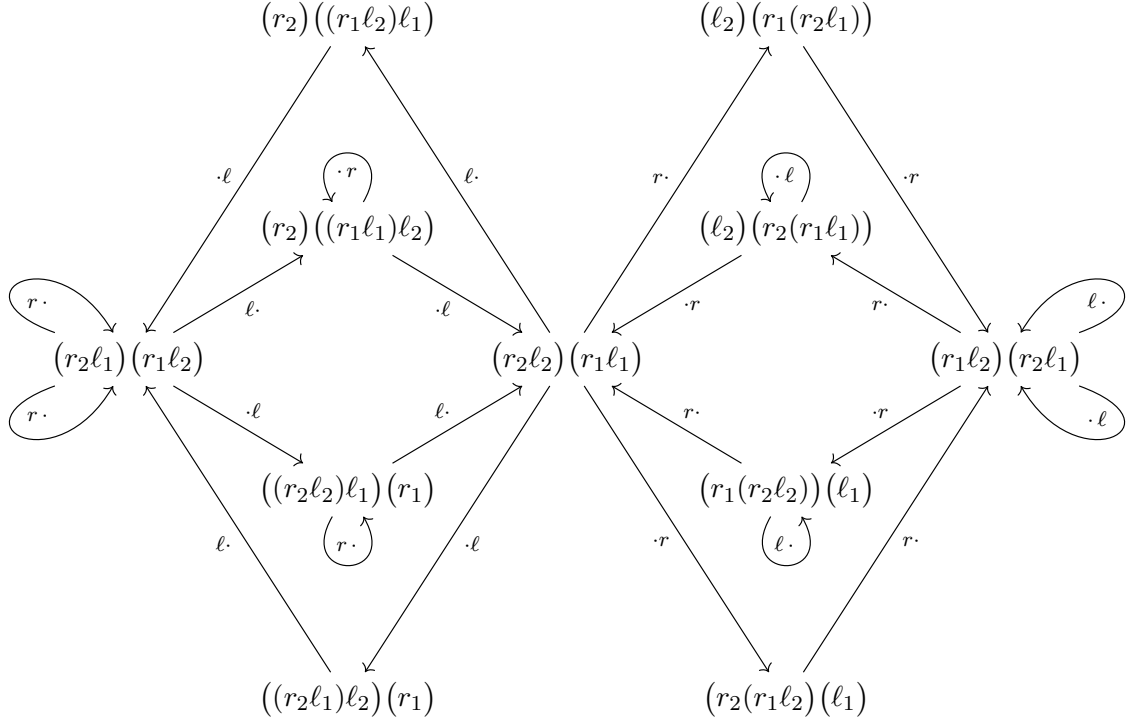

We can again depict any $DC$ class as a directed graph $\mathcal{G} = \mathcal{G}(\boldsymbol{\pi})$, where the graph is obtained by starting with an initial castle forest $\boldsymbol{\pi}$ and putting a directed edge to any admissible castle forest attainable through a diagonal change move. Figure \ref{fig:dc_graph_general} depicts such a graph, where the initial castle forest is given by $\boldsymbol{\pi} = \combi{r_2\ell_2}\combi{r_1\ell_1}$.

The graph depicted in Figure \ref{fig:dc_graph_general} admits a lot of symmetries. For instance, the graph is symmetric along the vertical, where the symmetry is given by exchanging the roles of $r$ and $\ell$. The symmetry along the horizontal axis is given by exchanging the left and right sides in all the castle forests as well as permuting 1 and 2. All these symmetries are typical features of such $DC$ graphs in complete analogy to the graphs arising in the hyperelliptic case. Let us also mention that the complexity of these $DC$ graphs grows extremely fast in $k$. In \cite{ferenczi2014generalization}, Ferenczi depicts a $DC$ graph corresponding to a castle forest consisting of 3 trees, which has 116 vertices, even though he is considering the greedy algorithm which leads to less edges than the single move algorithm we have depicted here. 

We are now equipped to formally define the space of (labeled) $k$-sets of castle polygons, which allows us to give a formal definition of the diagonal change move map. 

\begin{definition}[Space of labeled $k$-sets of castle polygons]
    The \emph{space of labeled $k$-sets of castle polygons} is defined as
    \begin{equation*}
        \mathcal{P}_k \coloneqq \left\{ (\boldsymbol{\pi}, \boldsymbol{w}) \mid \pi \in \mathcal{G}, \boldsymbol{w} \in \Omega_{\boldsymbol{\pi}} \right\}.
    \end{equation*}
\end{definition}

\begin{definition}[Diagonal change move on parameter space]\label{def:move_on_paramter_space}
    Let $\boldsymbol{\pi}$ be a castle forest with associated $DC$ graph $\mathcal{G}$ and let $c$ be a choice of diagonal change move. The \emph{diagonal change move} $\hat{m}_{\boldsymbol{\pi}, c}$ on $\{\boldsymbol{\pi}\} \times \Omega_{\boldsymbol{\pi}} \subseteq \mathcal{P}_k$ is the map
    \begin{align*}
        \hat{m}_{\boldsymbol{\pi}, c} \colon \{\boldsymbol{\pi}\} \times \Omega_{\boldsymbol{\pi}} &\to \mathcal{P}_k, \\
        (\boldsymbol{\pi}, \boldsymbol{w}) &\mapsto (c\cdot \boldsymbol{\pi}, A_{\boldsymbol{\pi},c} \cdot \boldsymbol{w}). 
    \end{align*}
\end{definition}

\subsubsection{Invertibility Through Backward Moves}

Defining the algorithm in the language of translation surfaces instead of IETs has the significant advantage that the moves are invertible. To use more technical terms, the algorithm here provides an invertible natural extension of the corresponding algorithm on the space of IETs. Let us properly define a backward diagonal change move, similar to the exposition in section \ref{sec:action_on_data}, and let us show that the forward and backward moves are mutually inverse.

Recall that in the hyperelliptic case, backward move were defined using the rotation operator $R$. This path was motivated by the abundance of symmetries inherent to hyperelliptic surfaces. In the general case here, this approach is no longer feasible, as is easily seen by the fact that if we rotate a $k$-set of castle polygons we must no longer have a $k$-set of castle polygons. Thus, we must describe backward moves in a different way than before. We begin by defining allowed \emph{backward} moves, which extends Definition \ref{def:allowed_moves}. In the sequel, we will write
\begin{equation}\label{eq:wedge_width_height}
    w_{i, \varepsilon} = (\lambda_{i, \varepsilon}, \tau_{i, \varepsilon}),
\end{equation}
so that $\lambda_{i, \varepsilon}$ represents the width and $\tau_{i, \varepsilon}$ represents the height of the corresponding wedge. Note that $\lambda_{i, \varepsilon} < 0$ if and only if $\varepsilon = \ell$, while all the other parameters are strictly positive. 

\begin{definition}[Allowed backward moves]\label{def:allowed_backward_moves}
    Let $P$ be a $k$-set of castle polygons and define $H_-(P) \subseteq [k]$ to be the set of integers $i$ satisfying one of the following conditions. 
    \begin{enumerate}
        \item $\tau_{i,r} > \tau_{i, \ell}$, the unique upper side $w_{i,r}$ in the $k$-set of castle polygons is not the single upper side of its polygon, meaning that this polygon is not a triangle, and the lower side of its stack triangle has positive slope.
        \item $\tau_{i, r}< \tau_{i, \ell}$, the unique upper side $w_{i, \ell}$ in the $k$-set of castle polygons is not the single upper side of its polygon and the lower side of its stack triangle has negative slope.
    \end{enumerate}
    We call $H_-(P)$ the \emph{set of allowed backward moves} in $P$.
\end{definition}

\begin{lemma}
    For any $k$-set of castle polygons $P$, the set $H_-(P)$ is nonempty.    
\end{lemma}
\begin{proof}
    Suppose towards a contradiction that $H_-(P)$ is empty. Choose some $i \in [k]$ such that $\tau_{i, r} > \tau_{i, \ell}$ and consider the upper side $w_{i,r}$ in $P$. If it is the only upper side in its polygon $p_j$, then we have
    \begin{equation}\label{eq:backward_nonempty_01}
        \tau_{j, r} > \tau_{i, r} > \tau_{i, \ell}.
    \end{equation}
    If $p_j$ is not a trinagle, the lower side of the stack triangle containing $w_{i,r}$ must have negative slope, since we assume $H_-(P)$ to be nonempty. If we write $\tau_{j_1, \ell}$ for the height of the right top side of this triangle, then $\tau_{j_1, \ell} > \tau_{i, r}$. If this right side is not an upper side of $p_j$, we can continue this way until we reach an upper side of $p_j$ after finitely many steps, and we obtain heights
    \begin{equation}\label{eq:backward_nonempty_02}
        \tau_{k, \ell} = \tau_{j_s, \ell} > \ldots > \tau_{j_1, \ell} > \tau_{i, r} > \tau_{i, \ell} 
    \end{equation}
    for some $k$. An analogous argument for the case $\tau_{i,r} < \tau_{i, \ell}$ shows that we find some $k$ such that
    \begin{equation}\label{eq:backward_nonempty_03}
        \tau_{k, r} > \tau_{i, \ell} > \tau_{i, r} \qquad \text{or} \qquad \tau_{k, \ell} > \tau_{i, \ell} > \tau_{i, r}.
    \end{equation}
    Combining equations \eqref{eq:backward_nonempty_01}, \eqref{eq:backward_nonempty_02} and \eqref{eq:backward_nonempty_03}, we have shown that for every $i$ we can find a $k$ such that
    \begin{equation*}
        \sup\{\tau_{k, \ell}, \tau_{k, r}\} > \sup\{\tau_{i, \ell}, \tau_{i, r}\},
    \end{equation*}
    which is impossible since we only have finitely many wedges.
\end{proof}

The definition of allowed backward moves is done exactly so that the procedure we explain next is well-defined. We will again start with a geometric description of backward moves. So let $P$ be a $k$-set of castle polygons and let $i \in H_-(P)$. If we write $w_{i,d}$ for the lower side of the stack triangle appearing in the definition of $H_-(P)$, then a backward move consists of cutting the polygon $p_j$ along $w_{i,d}$, keeping the lower part as the new polygon $p_j$ and pasting the upper part by $w^-_{i, \varepsilon}$ to the bottom of the caslte polygon, which has $w^-_{i, \varepsilon}$ as a lower side. Whether $\varepsilon = r$ or $\varepsilon = \ell$ is also determined by the definition of $H_-(P)$. We will call $w^-_{i,d}$ a \emph{backward diagonal}, which determines if we speak of a left backward move if its slope is positive, or a right backward move if its slope is negative. We now illustrate this process with an example.

\begin{example}\label{ex:backward_move}
    Let $P = (\boldsymbol{\pi}, \boldsymbol{w})$ be the 2-set of castle polygons determined by the data
    \begin{equation*}
        \boldsymbol{\pi} = \combi{(r_2\ell_2)(\ell_1)}\combi{r_1}
    \end{equation*}
    and
    \begin{equation*}
        \boldsymbol{w} = \big((-3.5, 3),(2,1),(-1, 2.5),(1,3.5)\big).
    \end{equation*}
    Applying Definition \ref{def:allowed_backward_moves} to $P$, one obtains that $H_-(P) = \{1,2\}$. This is visualized in Figure \ref{fig:backward_move_01}, where $w^-_{1,d}$ and $w^-_{2,d}$ are highlighted.
    \begin{figure}[ht]
        \centering
        \begin{tikzpicture}[scale = 0.7]
            \coordinate (p1) at (0,0);
            \coordinate (l1) at (-3.5,3);
            \coordinate (r1) at (2,1);
            \coordinate (t1) at (-1.5,4);
            \coordinate (t2) at (-2.5,6.5);

            \def \x{5};
            \coordinate (p2) at (\x,0);
            \coordinate (l2) at (\x - 1,2.5);
            \coordinate (r2) at (\x+1, 3.5);

            \draw (p1) -- node[midway, below right] {$w_{1, r}$} (r1) -- node[midway, above right] {$w_{1, \ell}$} (t1) -- node [midway, right] {$w_{2, \ell}$} (t2) -- node[midway, left]  {$w_{2,r}$}(l1) -- node[midway, below left] {$w_{1, \ell}$}cycle;
            \draw (p2) -- node[midway, right] {$w_{2,r}$}(r2) -- node[midway, above, xshift = -4] {$w_{1,r}$} (l2) -- node [midway, left] {$w_{2, \ell}$} cycle;

            \draw[dashed, noamblue] (t1) -- (l1) node[left] {$w^-_{2,d}$};
            \draw[dashed, purple] (l1) -- (r1) node[right] {$w^-_{1,d}$};
        \end{tikzpicture}
        \caption{The initial 2-set of castle polygons from Example \ref{ex:backward_move}}
        \label{fig:backward_move_01}
    \end{figure}
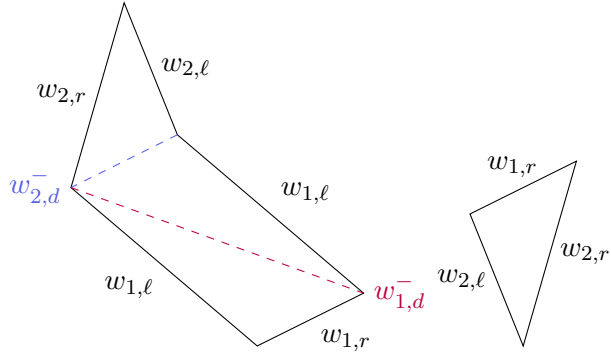
    In Figure \ref{fig:backward_move_02}, we show the outcome of decision $c(H(P)) = (1, r)$ on the left and the outcome of decision $c(H(P))= (2, \ell)$ on the right. The polygons that were pasted by applying the procedure are highlighted in gray.

    \begin{figure}[hb]
        \centering
        \begin{tikzpicture}[scale = 0.6]
            \coordinate (p1) at (0,0);
            \coordinate (l1) at (-5.5,2);
            \coordinate (r1) at (2,1);
            \coordinate (t1) at (-3.5,3);
            \coordinate (t2) at (-4.5,5.5);

            \draw[dotted] (p1) -- (t1);
            \draw[dashed] (l1) -- (r1);
            \draw[dashed] (l1) -- (t1);
            \fill[gray, opacity = 0.1] (p1) -- (l1) -- (t2) -- (t1) -- cycle;

            \def \x{4.5};
            \coordinate (p2) at (\x,0);
            \coordinate (l2) at (\x - 1,2.5);
            \coordinate (r2) at (\x+1, 3.5);

            \draw (p1) -- node[midway, below right] {$w_{1, r}$} (r1) -- node[midway, above right] {$w_{1, \ell}$} (t1) -- node [midway, right] {$w_{2, \ell}$} (t2) -- node[midway, left]  {$w_{2,r}$}(l1) -- node[midway, below left] {$w_{1, \ell}$}cycle;
            \draw (p2) -- node[midway, right] {$w_{2,r}$}(r2) -- node[midway, above, xshift = -4] {$w_{1,r}$} (l2) -- node [midway, left] {$w_{2, \ell}$} cycle; 

            \def\y{12};
            \coordinate (p3) at (\y,0);
            \coordinate (l3) at (\y - 3.5,3);
            \coordinate (r3) at (\y + 2,1);
            \coordinate (t3) at (\y - 1.5, 4);

            \draw (p3) -- node[midway, below right] {$w_{1, r}$} (r3) -- node[midway, above right] {$w_{1, \ell}$} (t3) -- node[midway, above left]{$w_{2, r}$}  (l3) -- node[midway, below left] {$w_{1, \ell}$} cycle;
            \draw[dashed] (r3) -- (l3);

            \def \w{16.5};
            \coordinate (p4) at (\w,0);
            \coordinate (l4) at (\w - 1,2.5);
            \coordinate (r4) at (\w+2, 1);
            \coordinate (t4) at (\w+1,3.5);

            \draw (p4) -- node[midway, below, xshift = 6] {$w_{2,r}$} (r4) -- node[midway, right] {$w_{2, \ell}$} (t4) -- node[midway, above, xshift = -4] {$w_{1, r}$} (l4) -- node[midway, left] {$w_{2,\ell}$} cycle;
            \draw[dashed] (l4) -- (r4);
            \draw[dotted] (p4) -- (t4);
            \fill[opacity = 0.1, gray] (p4) -- (r4) -- (t4) -- cycle;
        \end{tikzpicture}
        \caption{The outcomes of applying the move associated to $c(H_-(P)) = (1, r)$ (on the left) and the move associated to $c(H_-(P)) = (2, \ell)$ (on the right).}
        \label{fig:backward_move_02}
    \end{figure}
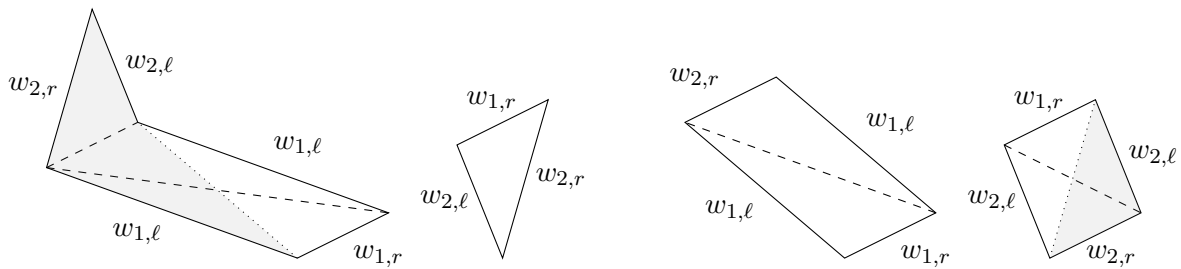
\end{example}

We have the following result which is an analogue of Lemma \ref{lem:diagonal_change produces_castle_polygons}.

\begin{lemma}[Backward diagonal change produces castle polygon]
    The image of a $k$-set of castle polygons by a backward diagonal change move given by a decision $c$ is again a $k$-set of castle polygons.
\end{lemma}
\begin{proof}
    This follows immediately from the definitions.
\end{proof}
\begin{remark}
    In order to prove Lemma \ref{lem:diagonal_change produces_castle_polygons}, it was necessary to check that no cut part has to be pasted to itself. That this is not possible here follows already from the definitions, more precisely it is due to the fact that any triangle is either a stack or base triangle but not both. 
\end{remark}
Equipped with this geometric intuition, define backward moves formally by their action on the data, exactly as we did for (forward) diagonal change moves. Let us start by giving a formula for the backward diagonal $w^-_{i,d}$ appearing in the definition of backward diagonal change move. As in \eqref{eq:leaf_sums}, we define
\begin{equation*}
    u_{i,\ell} = \sum_{\substack{\text{left of}\\\ell_i}} w_{j, \ell} -  w_{j, r}, \qquad u_{i, r} = \sum_{\substack{\text{right of}\\r_i}} w_{i, r} -  w_{i, \ell}.
\end{equation*}
If $\varepsilon = r$, then $w^-_{i, d} = w_{i,\ell} + u_{i, \ell}$, whereas if $\varepsilon = \ell$, then $w^-_{i, d} = w_{i, r} + u_{i, r}$. Applying a backward move corresponding to the choice $c(H_-(P)) = (i, \varepsilon)$ ha no effect on wedges of plygons $p_j$ for $j \neq i$. One part of the wedge of $p_i$ (depending on $\varepsilon$) is exchanged with $w^-_{i,d}$, so these formulas determine the action on the length datum. Written more concisely, for $c(H_-(P)) = (i, \varepsilon)$ we have 
\begin{equation*}
    \begin{cases}
        w'_{j, \eta} = w_{j, \eta}\quad &\text{if } j \neq i \text{ or } \eta \neq \varepsilon, \\
        w'_{i, \varepsilon} = w_{i,d}^-.
    \end{cases}
\end{equation*}

To emphasize the similarity to the action on the combinatorial datum by \emph{forward} diagonal change moves, we will initially describe the action of a backward move in the language of surgery on trees as well. 

\begin{proposition}[Surgery on trees - backward moves]\label{prop:tree_surgery_backward}
    Let $P = (\boldsymbol{\pi}, \boldsymbol{w})$ be a $k$-set of castle polygons and let $c$ be a backward decision with $c(H_-(P)) = (i, \varepsilon)$ for $\varepsilon \in \{\ell, r\}$. Applying the associated backward diagonal change amount sto the action on the castle forest $\boldsymbol{\pi}$ described below.

    Suppose that $c(H_-(P)) = (i, r)$ and write $n$ for the node representing $w^-_{i, d}$, that is, the node below the leaf $\ell_i$ which is part of tree $j$. 

    Cut tree $j$ at node $n$ so that the node gets duplicated and the copy that becomes a leaf is labeled $\ell_i$. Add a stem if the node $n$ is also a root. The other copy of $n$ becomes the root of the new polygon $i$ together with everything above the node. At this point, we distinguish two cases.

    \begin{enumerate}
        \item If former tree $i$, consisting only of the node $n$ and everything below, has a single leaf $\ell_k$, then we rename $\ell_i$ to $\ell_k$ in the new tree and delete former tree $i$.
        \item If former tree $i$, consisting of the node $n$ and everything that is below, has more than one leaf, we put this tree on top of the leaf labeled $\ell_i$. There, the leaf $\ell_i$ becomes a node.
    \end{enumerate}
    If instead we have $c(H_-(P)) = (i, \ell)$, we apply the exact same procedure except for exchanging $\ell$ and $r$ as well as left and right.
\end{proposition}
\begin{proof}
    As for Proposition \ref{prop:tree_surgery}, this follows immediately from the definition of a backward move. 
\end{proof}

Of course, this action admits a description in terms of forest words as well. 

\begin{corollary}[Surgery on forest words - backward moves]\label{cor:action_forest_words_backward}
    Proposition \ref{prop:tree_surgery_backward} can be stated as follows in the language of forest words.
    \begin{itemize}
        \item If $c(H_-(P)) = (i, r)$, i.e., if we apply a right backward move, then we adjust the forest word as follows.
        
        \begin{enumerate}
            \item The subword induced by the parenthesis just to the right of $\ell_i$ is split into two parts, one of which is the single letter $\ell_i$ and the other part is everything else.
            \item The position of letter $\ell_i$ remains unchanged, the other part is glued to the left of the subword containing $\ell_k$, where $k$ is the tree containing $\ell_i$.
        \end{enumerate}
        \item If $c(H_-(P)) = (i, \ell)$, i.e., if we apply a left backward move, then we adjust the forest word as follows.
        \begin{enumerate}
            \item The subword induced by the parenthesis just to the left of $r_i$ is split into two parts, one of which is the single letter $r_i$ and the other part is everything else.
            \item The position of letter $r_i$ remains unchanged, the other part is glued to the right of the subword containing $r_k$, where $k$ is the tree containing $r_i$. 
        \end{enumerate}
    \end{itemize}
\end{corollary}

As an example, we can consider Example \ref{ex:forest_words} read from bottom to top. More precisely, we can check that a backward move corresponding to $c(H_-(P)) = (\ell, 2)$ applied to $\boldsymbol{\pi}^{(4)}$ as explained in Corollary \ref{cor:action_forest_words_backward} results in the forest word $\boldsymbol{\pi}^{(3)}$, and the same action, this time corresponding to $c(H_-(P)) = (1,r)$ on $\boldsymbol{\pi}^{(2)}$, results in the forest word $\boldsymbol{\pi}^{(1)}$. 

Since the action on the combinatorial and length data by a backward moves is defined in complete analogy to a forward move, we can describe the action of a backward move by a matrix as well. Given a choice $c(H(P)) = (i, \varepsilon)$ corresponding to a $k$-set of castle polygons $P = (\boldsymbol{\pi}, \boldsymbol{w})$, we will denote this \emph{backward diagonal change matrix} by $A_{\boldsymbol{\pi, c}}^{-1}$. The following proposition justifies the use of this notation, as it shows that forward and backward moves with respect to the same choice are inverse to each other, which in turn implies that the corresponding matrices are inverse to each other as well.

\begin{proposition}[Inverses of moves]
    Forward and backward moves corresponding to the same choice are mutually inverse. More precisely, if we write $\hat{n}_{\boldsymbol{\pi},c}$ for the map defined as in Definition \ref{def:move_on_paramter_space} but for backward diagonal change moves, then
    \begin{equation*}
        \hat{n}_{c_m\cdot \boldsymbol{\pi}, c_n} \circ \hat{m}_{\boldsymbol{\pi}, c_m} = \hat{m}_{c_n\cdot \boldsymbol{\pi}, c_m} \circ \hat{n}_{\boldsymbol{\pi}, c_n} = \operatorname{Id},
    \end{equation*}
    provided that the choices $c_m$ and $c_n$ for the forward and backward moves coincide. 
\end{proposition}
\begin{proof}
    After a forward diagonal change move, a right move applied to polygon $p_i$, say, the new $\tau_{i, \ell}$ is larger than the old $\tau_{i, r}$ by definition of stack triangle. This ensures that in the image castle, the same decision is always allowed. The fact that applying a backward move corresponding to this decision then yields the original $k$-set of castle polygons follows from the definition of the moves. All the other cases follow by a similar argument. 
\end{proof}

\subsection{Existence of Castle Polygon Decomposition}
We have explained how to apply the general algorithm to a $k$-set of castle polygons. The first question that naturally arises is whether it is always possible to decompose a translation surface ${X \in \mathcal{H}(k_1 - 1, \ldots, k_n - 1)}$ such that it is represented as a $k$-set of castle polygons. By Lemma \ref{lem:total_cone_angle_castle_polygons}, the integer $k$ is necessarily fixed, depending only on the stratum. We will show that this is the case under the usual genericity assumptions. We argue that the existence of a decomposition into $k$-sets of castle polygons is equivalent to the existence of a \emph{veering triangulation}, which allows us to invoke Theorem 2.16 from \cite{bell2019coding} asserting the existence of such a triangulation.

Towards this goal, we will start by introducing the concept of a \emph{veering triangulation}.

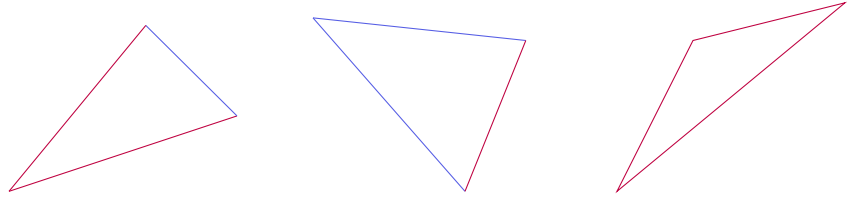
\begin{figure}[ht]
    \centering
    \begin{tikzpicture}
        \draw[purple] (0,0) -- (3,1);
        \draw[noamblue] (3,1) -- (1.8,2.2);
        \draw[purple] (1.8, 2.2) -- (0,0);

        \def \x{4}
        \draw[noamblue] (\x, 2.3) -- (\x+2,0);
        \draw[noamblue] (\x, 2.3) -- (\x+2.8,2);
        \draw[purple] (\x+2.8, 2) -- (\x+2,0);

        \def \w{8}
        \draw[purple] (\w,0) -- (\w+3,2.5) -- (\w + 1, 2) -- cycle;
    \end{tikzpicture}
    \caption{The triangles on the left and in the middle are veering, the triangle on the right is not veering. }
    \label{fig:veering_01}
\end{figure}

\begin{definition}[Veering triangulation]\label{def:veering_triangulation}
    A flat triangulation $\tau = \{T_1, \ldots, T_m\}$ of a translation surface $X \in \mathcal{H}(k_1 - 1, \ldots, k_n - 1)$ (see Definition \ref{def:flat_triangulation}) is called a \emph{veering triangulation}, if every triangle $T_i \in \tau$ contains at least one increasing and one decreasing edge with respect to the natural orientation of the plane. We refer to such triangles $T_i$ as \emph{veering triangles}.
\end{definition}
Following \cite{bell2019coding}, we will give colors to the edges in the following way. Any increasing edge will be colored red, whereas every decreasing edge is colored blue. Using this convention, Definition \ref{def:veering_triangulation} can equivalently be formulated by requiring none of the triangles $T_i \in \tau$ to be monochromatic. Some examples of veering and non-veering triangles can be seen in Figure \ref{fig:veering_01}.

Comparing Figures \ref{fig:base_stack_neither} and \ref{fig:veering_01}, the connection between veering triangulations and $k$-sets of castle polygons becomes apparent. 
\begin{lemma}[Veering triangles are exactly stack or base triangles]\label{lem:veering_base_stack}
    A triangle $T_i$ in a flat triangulation $\tau$ is veering if and only if it is a base triangle or a stack triangle. 
\end{lemma}
\begin{proof}
    If $T_i$ is not veering, it cannot be a base triangle since the wedges of a base triangle are bichromatic by definition. Neither can it be a stack triangle, since the vertical downwards ray from its highest vertex cannot cross the opposing side. If it did, then it would necessarily be the case that the edge to the left of the top vertex is red and the edge to the right is blue, which is impossible in a non-veering triangle. 

    Conversely, if $T_i$ is veering and the two edges adjacent to the lowest vertex are bichromatic, it follows by definition that the triangle is a base triangle. If both edges adjacent to the lowest edge share the same color, the third edge is of a different color meaning that the downwards vertical from the top vertex is contained in the triangle and crosses the opposing side, making $T_i$ a stack triangle. 
\end{proof}

We can now state the main theorem of this section.

\begin{theorem}[Existence of castle polygons]\label{thm:existence_castle_polygons}
    Let $X \in \mathcal{H}(k_1 - 1, \ldots, k_n - 1)$ be a translation surface satisfying Keane's condition. Then, it admits a decomposition into a $k$-set of castle polygons.
\end{theorem}
\begin{proof}
    Let $\tau$ be a veering triangulation of $X$, which by the proof of Lemma \ref{lem:total_cone_angle_castle_polygons} consists of 
    \begin{equation*}
        T_X = 2(2\mathbf{g} - 2 + n)
    \end{equation*}
    triangles. By Lemma \ref{lem:veering_base_stack}, these triangles are either base or stack triangles. In fact, the number of base and stack triangles is the same. To see this, recall from the proof of Lemma \ref{lem:total_cone_angle_castle_polygons} that each singularity $p_j$ of $X$ has $k_j$ outgoing and incoming rays of the vertical flow each. Hence, the singularity $p_j$ corresponds to $k_j$ base and stack triangles each. Repeating this argument finitely often for all singularities, we obtain the claim. 

    To build a $k$-set of castle polygons out of this triangulation, simply put the $k = \sum_{j = 1}^n k_j$ base triangles on top and glue all the stack triangles on top using the identifications induced by the triangulation $\tau$. Note that any bottom side of a stack triangle is identified with a top side of a base or stack triangle, but never with a bottom side. Therefore, we can glue the stack triangles on top of the base triangles in a way such that the number of polygons in the set is always $k$. 
\end{proof}

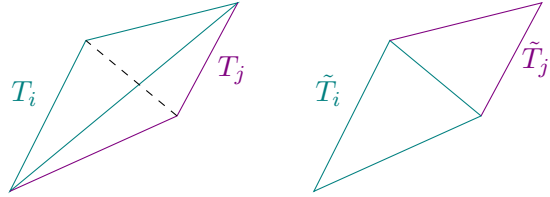
\begin{figure}[ht]
    \centering
    \begin{tikzpicture}
        \def \w{0}
        \draw[teal] (\w,0) -- (\w+3,2.5) -- (\w + 1, 2) -- node[midway, above left] {$T_i$} cycle;
        \draw[violet] (\w,0) -- (\w + 2.2, 1) -- node[midway, below right, yshift = 5] {$T_j$} (\w+3, 2.5);

        \draw[dashed] (\w+2.2,1) -- (\w+1, 2); 

        \def \x{4}
        \draw[teal] (\x,0) -- node[midway, above left] {$\Tilde{T}_i$} (\x + 1, 2) -- (\x+2.2, 1) -- cycle; 
        \draw[violet] (\x+1, 2) -- (\x+3, 2.5) -- node[midway, right] {$\Tilde{T}_j$} (\x+2.2, 1);  
    \end{tikzpicture}
    \caption{Flipping an edge in a triangulation.}
    \label{fig:edge_flipping_01}
\end{figure}

We refer the reader to \cite{bell2019coding} for a full proof on the existence of veering triangulations for translation surfaces. In fact, the same result is shown in the more general setting of half-translation surfaces. Nonetheless, we want to illustrate the main idea how to obtain this result. Say we begin with any flat triangulation $\tau$ of $X$, which as we have mentioned in section \ref{sec:general_principle} always exists. If $\tau$ is not veering, some $T_i \in \tau$ is not a veering triangle, i.e., it is monochromatic. Suppose that all its edges are red. We can consider an adjacent triangle $T_j$ sharing any edge of $T_i$. Clearly, we must have $T_j \neq T_i$ for all adjacent triangles, so that their union forms a quadrilateral $Q_{ij}$ as on the left of Figure \ref{fig:edge_flipping_01}.
We can \emph{flip an edge}, which adjusts the triangulation $\tau$ in the following way. We take the quadrilateral $Q_{ij}$, which has two diagonals $d$, which is the shared edge of $T_i$ and $T_j$, as well as $\Tilde{d}$. We can now obtain two new triangles $\Tilde{T}_i$ and $\Tilde{T}_j$ by removing $d$ and replacing it by $\Tilde{d}$. This is shown on the right of Figure \ref{fig:edge_flipping_01}. This induces a new triangulation $\Tilde{\tau}$ of $X$. We can see that neither $T_i$ nor $T_j$ were veering, but both triangles obtained from flipping an edge \emph{are} veering. \emph{If} the two triangles $T_i$ and $T_j$ were the only triangles preventing the triangulation $\tau$ from being a veering triangulation, we would have obtained a veering triangulation $\Tilde{\tau}$ in this way. But it needs not be the case that we always obtain veering triangles after flipping an edge. However, the proof of Theorem 2.16 in \cite{bell2019coding} on the existence of a veering triangulation reveals, that we always obtain a veering triangulation after finitely many edge flips, given that Keane's condition is satisfied. 

We finish this section by mentioning another possibility to obtain the existence of a decomposition into a $k$-set of castle polygons. Theorem \ref{thm:existence_castle_polygons} fundamentally relies on the existence of a \enquote{good} triangulation, i.e., the existence of a triangulation consisting of an equal number of base and stack triangles. Veering triangulations offer one possibility, another type of triangulations with this property are (weak) $L^\infty$-Delaunay triangulations. I am grateful to Pouya Honaryar for explaining this concept to me. We will now introduce these triangulations and prove that any generic translation surface does in fact admit a $L^\infty$-Delaunay triangulation. In this part, we stay close to the exposition in \cite{zykoski2023polytopal}.
\begin{definition}[Immersed square]\label{def:immersed_square}
    An \emph{immersed square} $S \subseteq X$ is a subset with no singularities in its interior, obtained by isometrically immersing a Euclidean square $[0,h]\times [0,h]$ in $X$, for some $h > 0$. An immersed square is said to be \emph{maximal}, if it is not properly contained in any other immersed square.
\end{definition}
Of course, Definition \ref{def:immersed_square} describes a special case of what we called \emph{immersed rectangles} in section \ref{sec:diagonal_changes}. A maximal immersed square will have singularities on its boundary. If $X$ satisfies Keane's condition, there can be at most 4 such singularities, since no two singularities can lie on the same edge of $S$. Generically, a maximal square will contain 3 singular points on its border. In \cite{zykoski2023polytopal}, the former case with 4 singularities is excluded. In our setting, meaning for translation surfaces satisfying Keane's condition, it may actually happen, that some of the squares we are interested in have 4 singular points on its boundary. To emphasize this difference, we will speak of \emph{weak} $L^\infty$-Delaunay triangles and \emph{weak} $L^\infty$-Delaunay triangulations, in order to highlight the fact that the triangulations we consider are slightly more general. 

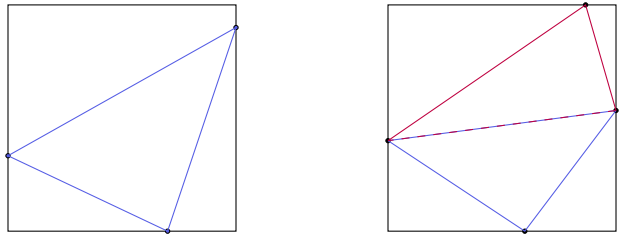
\begin{figure}[ht]
    \centering
    \begin{tikzpicture}
        \draw (0,0) rectangle (3,3);
        \draw[fill=noamblue] (0,1) circle (0.8pt);
        \draw[fill=noamblue] (3,2.7) circle (0.8pt);
        \draw[fill=noamblue] (2.1,0) circle (0.8pt);
        \draw[noamblue] (0,1) -- (3,2.7) -- (2.1,0) -- cycle;

        \def \x{5}
        \draw (\x,0) rectangle (\x+3, 3); 
        \draw[fill = black] (\x,1.2) circle (0.8pt);
        \draw[fill = black] (\x+1.8,0) circle (0.8pt);
        \draw[fill = black] (\x+3, 1.6) circle (0.8pt);
        \draw[fill = black] (\x+2.6,3) circle (0.8pt);

        \draw[noamblue] (\x,1.2) -- (\x+1.8,0) -- (\x+3,1.6);
        \draw[noamblue] (\x,1.2) -- (\x+3, 1.6);
        \draw[purple] (\x,1.2) -- (\x+2.6,3) -- (\x+3,1.6);
        \draw[purple, dashed] (\x+3, 1.6) -- (\x,1.2);
        
    \end{tikzpicture}
    \caption{A circumsquare with three singularities on the boundary containing a unique triangle on the left, and a circumsquare with four singularities on the boundary containing two distinct triangles on the right.}
    \label{fig:circumsquares}
\end{figure}

\begin{definition}[Circumsquare and weak $L^\infty$-Delaunay triangle]
    Let $T$ be a triangle on $X$ whose vertices are singularities and whose edges are geodesics. We say that an immersed square $S$ is a \emph{circumsquare of} $T$, if $T$ is the image of a triangle inscribed in $[0,h] \times [0,h]$ under the immersion from Definition \ref{def:immersed_square}.

    We say that $T$ is a \emph{weak $L^\infty$-Delaunay triangle}, if it has a circumsquare $S$ and the edges of $T$ are neither horizontal nor vertical. 
\end{definition}
Figure \ref{fig:circumsquares} shows to possible situations of circumsquares and weak $L^\infty$-Delaunay triangles. 
\begin{definition}[Weak $L^\infty$-Delaunay triangulation]
    Let $\tau = \{T_1, \ldots, T_m\}$ be a triangulation of a translation surface $X \in \mathcal{H}(k_1 - 1, \ldots, k_n - 1)$. We say that $\tau$ is a \emph{weak $L^\infty$-Delaunay triangulation}, if $T_i$ is a weak $L^\infty$-Delaunay triangle for every $i \in [m]$.  
\end{definition}
\begin{remark}
    One may also consider (weak) $L^p$-Delaunay triangulations, where $1 \leq p < \infty$. In this case, one replaces the squares in all the definitions above by circles with respect to the norm $\|\cdot\|_p$. The most commonly studied Delaunay triangulations are the $L^2$-Delaunay triangulations, i.e., the triangulations based on Euclidean circumcircles. Often in the literature, these triangulations are simply referred to as \emph{Delaunay-triangulations}. 
\end{remark}
The following fact is immediate, which together with Lemma \ref{lem:veering_base_stack} implies that weak $L^\infty$-Delaunay triangles are either base or stack triangles.
\begin{lemma}
    Weak $L^\infty$-Delaunay triangles are veering.
\end{lemma}
Hence, if we can show that any generic translation surface $X$ admits a weak $L^\infty$-Delaunay triangulation, we can argue exactly as above to obtain the existence of a decomposition into a $k$-set of castle polygons. The proof of existence we give here is a slight modification of Lemma II.5 in \cite{zykoski2023polytopal} and centrally uses Proposition 2.1 from \cite{gueritaud2016veering}.
\begin{theorem}[Existence of weak $L^\infty$-Delaunay triangulations]
    Let $X \in \mathcal{H}(k_1 - 1, \ldots, k_n - 1)$ satisfy Keane's condition. Then it admits a (not necessarily unique) weak $L^\infty$-Delaunay triangulation. 
\end{theorem}
\begin{proof}
    Denote by $\Tilde{X}$ the universal cover of $X$. For each maximal immersed square in $\Tilde{X}$, which has either 3 or 4 singularities on its boundary, we can connect these singularities by either 3 or 5 Euclidean geodesics such that these geodesics do not intersect. Note that there is a choice involved in case there are 4 singularities on the boundary, which is the reason that the triangulation we obtain needs not be unique. 

    Write $\Tilde{\tau}$ for the collection of all geodesics obtained this way and let $\tau = p(\Tilde{\tau})$ be the image under the universal covering map. Notice, that $\tau$ is a collection of non-intersecting saddle connections on $X$. By Proposition 2.1 of \cite{gueritaud2016veering}, $\tau$ induces a finite flat triangulation of $X$, and since we assume $X$ to satisfy Keane's condition, the corresponding triangles are weak $L^\infty$-Delaunay triangles and hence $\tau$ induces a weak $L^\infty$-Delaunay triangulation.
\end{proof}
To end this section, we want to remark that it is possible to show that a translation surface $X$ admits a \emph{unique} $L^\infty$-Delaunay triangulation by adjusting the genericity condition to exclude the case where a circumsquare can have 4 singularities on its boundary. In \cite{zykoski2023polytopal}, such surfaces are called $L^\infty$\emph{-generic}, which form an open and dense subset of the whole stratum. 

\subsection{Fundamental Domain and a Section for the Teichmüller Geodesic Flow}
This section develops analogous results to the ones from sections \ref{sec:fundamental_domain_hyperelliptic} and \ref{sec:section_for_teichmueller_flow}. Namely, we explain how to obtain a fundamental domain in the space $\mathcal{P}_k$ of labeled $k$-sets of castle polygons, which descends to \emph{unlabeled} $k$-sets of castle polygons as well, exactly as in the hyperelliptic case. Moreover, we will define a subset of this fundamental domain, giving a Poincaré section for the Teichmüller geodesic flow $g_t$, which \enquote{sees} all closed Teichmüller geodesics.

We start with the discussion on the fundamental domain. The notion of compatibility of labels from Definition \ref{def:compatibility_of_labels} can be adapted to this setting, i.e., we say that given $X \in \stratum$ equipped with a labeling of the bundles $\Gamma_i, i \leq i \leq k$, where $k = \sum_{i = 1}^n k_i$, and a decomposition into a (labeled) $k$-set of castle polygons ${P}$, the labelings of $X$ and $P$ are \emph{compatible} if the lower wedges of the base triangle of $p_i \in P$ consists of saddle connections in $\Gamma_i$. We will denote by $\mathcal{P}_k(X)$ the collection of $k$-sets of castle polygons compatible with $X$. Recall the operation of flipping an edge which we introduced in the description of how to obtain a veering triangulation from any other triangulation. Since the cutting operation of any forward or backward diagonal change move as described previously in this chapter is always applied to a diagonal of a quadrilateral, we can see that restricted to $k$-sets of castle polygons the operation of edge flipping and a diagonal change move is exactly the same. Hence, results about edge flips applied to veering triangulations translate immediately to our setting. In particular, we can leverage this equivalence to show that the space of (labeled) $k$-sets of castle polygons admits a lattice structure in the following way, by following the analogous arguments in \cite{bell2019coding}.

Let $\gamma$ be a vertical half-ray emanating from some singularity $\sigma$ on a translation surface $X$, and suppose $\tau$ is a triangulation associated to a decomposition into a $k$-set of castle polygons $P$. Let us write $T(\gamma, \tau)$ to be the triangle of $\tau$ containing $\gamma$. Developing $\sigma, \gamma$ and $T(\gamma, \tau)$ to the plane, we define $w(\gamma, \tau)$ to be the cone in the plane based at $\sigma$ and spanned by wedge given by the lower two edges if $T(\gamma, \tau)$ is a base triangle, or given by the upper two edges if it is a stack triangle. If now $P'$ is another $k$-set of castle polygons with triangulation $\tau'$ with the property that for every vertical half-ray $\gamma$ emanating from a singularity, we have
\begin{equation*}
    w(\gamma, \tau') \subseteq w(\gamma, \tau),
\end{equation*}
then we write $P \prec P'$ and say that $P$ \emph{precedes} $P'$. Corollary 2.24 and Proposition 2.25 from \cite{bell2019coding} then translate to the following results in our setting.

\begin{proposition}
    If $P \prec P'$, then there exists a sequence of forward diagonal change moves from $P$ to $P'$.
\end{proposition}

\begin{proposition}\label{prop:castle_lattices}
    Given a translation surface $X \in \stratum$, then $(\mathcal{P}_k(X), \prec)$ is a lattice, i.e., a partially ordered sets with unique joins and meets.
\end{proposition}
We will use the notation $P \curlyvee P'$ to denote the join of $P$ and $P'$ and similarly, we will write $P \curlywedge P'$ for their meet. The proof of Proposition \ref{prop:castle_lattices} in \cite{bell2019coding} reveals that the meet and join are distributive, i.e., we have
\begin{equation*}
    (P \curlyvee P') \curlywedge P'' = (P \curlywedge P'') \curlyvee (P' \curlywedge P''),
\end{equation*}
and an analogous equality holds by switching meets and joins. This lattice structure is what allows us to choose a canonical representative of $\mathcal{P}_k(X)$. The following definition is a direct adaptation from Definition \ref{def:width_intervals}.

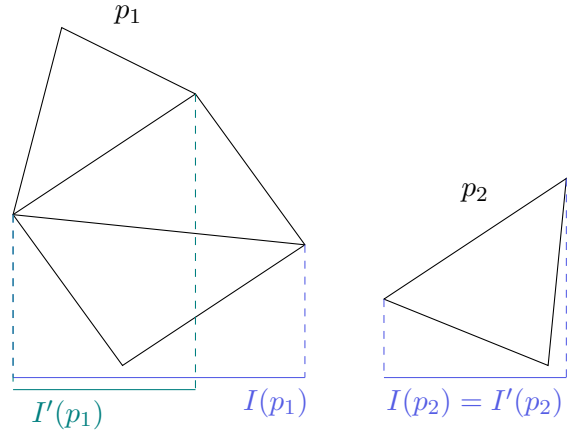
\begin{figure}[ht]
    \centering
    \begin{tikzpicture}[scale = 0.8]
        \draw (0,0) -- (3,2) -- (1.2,4.5) -- (-1.8,2.5) -- cycle;
        \draw (3,2) -- (-1.8,2.5) -- (-1,5.6) -- node[midway, above, yshift = 10] {$p_1$} (1.2,4.5);

        \draw[noamblue] (-1.8, -0.2) -- node[pos = 0.9, below] {$I(p_1)$} (3,-0.2); 
        \draw[teal] (-1.8, -0.4) -- node[pos=0.3, below] {$I'(p_1)$}(1.2, -0.4); 
        \draw[noamblue, dashed] (-1.8, 2.5) -- (-1.8, -0.2); 
        \draw[noamblue, dashed] (3,2) -- (3,-0.2);

        \draw[teal, dashed] (-1.8, 2.5) -- (-1.8, -0.4);
        \draw[teal, dashed] (1.2, 4.5) -- (1.2, -0.4);

        \def \x{7}

        \draw (\x,0) -- (\x+0.3,3.1) -- node[midway, above, yshift = 10] {$p_2$} (\x-2.7, 1.1) -- cycle;

        \draw[noamblue] (\x-2.7, -0.2) -- node[midway, below] {$I(p_2) = I'(p_2)$} (\x+0.3, -0.2);
        \draw[noamblue, dashed] (\x-2.7, 1.1) -- (\x-2.7, -0.2); 
        \draw[noamblue, dashed] (\x+0.3, 3.1) -- (\x+0.3, -0.2);
    \end{tikzpicture}
    \caption{The intervals $I(p_i)$ and $I'(p_i)$ from Definition \ref{def:width_intervals_general}.}
    \label{fig:width_intervals_general}
\end{figure}

\begin{definition}[Width intervals of castle polygons]\label{def:width_intervals_general}
    Let $P$ be a $k$-set of castle polygons and write $\lambda_{i, \varepsilon}$ for the real part of $w_{i, \varepsilon}$ as in \eqref{eq:wedge_width_height}. For any polygon $p_i \in P$, we define
    \begin{align*}
        I(p_i) &= [\lambda_{i, \ell}, \lambda_{i, r}],\\
        I'(p_i) &= \begin{cases}
            I(p_i) & \text{if } p_i \text{ is a triangle}, \\
            [\lambda_{i, \ell}, \lambda_{i, d}] & \text{if we would apply a left move},\\
            [\lambda_{i, d}, \lambda_{i, r}] &\text{if we would apply a right move},
        \end{cases}
    \end{align*}
    where we write $\lambda_{i,d}$ for the real part of the forward diagonal as in \eqref{eq:leaf_sums}. Moreover, we write $|I(p_i)|$ and $|I'(p_i)|$ for the widths of these intervals.
\end{definition}
If we apply a diagonal change move to $p_i$, i.e., if the corresponding choice is given by $c(H(P)) = (i, \varepsilon)$, then $|I(p_i)|$ gives the width of the polygon $p_i$ before the move and $|I'(p_i)|$ gives the with of the polygon after the move. If we write $p'_i$ for the newly obtained polygon, then we have $I'(p_i) = I(p'_i)$. Definition \ref{def:width_intervals_general} is illustrated in Figure \ref{fig:width_intervals_general}.

On $\mathcal{P}_k$, we can define an equivalence relation $\sim$ in complete analogy to Definition \ref{def:equivalence_staircase_moves}. It is for this equivalence relation that we can then define a fundamental domain, a subset of which will give us the required Poincaré section. 
\begin{definition}[Equivalence up to diagonal change moves]\label{def:equivalence_diagonal_changes}
    Two $k$-sets of castle polygons $P, P' \in \mathcal{P}_k$ are \emph{equivalent up to diagonal change moves}, and we write $P \sim P'$, if $P'$ can be obtained from $P$ by a sequence of forward and backward diagonal change moves. 
\end{definition}

\begin{definition}[Fundamental domain]\label{def:fundamental_domain_general}
    Let $\mathcal{F}$ be the set of all $k$-sets of castle polygons $P \in \mathcal{P}_k$ such that the following two conditions hold. 
    \begin{enumerate}[start=1,label={($F$\arabic*)}]
        \item For every $p \in P$ we have $|I(p)| \geq 1$.
        \item For every $p \in P$ that is not a triangle, $|I'(p)|<1$.
    \end{enumerate}
\end{definition}
Adapting the terminology from \cite{bell2019coding}, we will say that a $k$-set of castle polygons $P \in \mathcal{F}$ is \emph{balanced}. In order to show that $\mathcal{F}$ is indeed a fundamental domain with respect to $\sim$, we will show that every translation surface $X \in \stratum$ satisfying Keane's condition admits a \emph{unique} decomposition into a \emph{balanced} $k$-set of castle polygons. We proceed by adapting the proof of Proposition 2.28 in \cite{bell2019coding} to our setting.

\begin{lemma}[Fundamental domain]\label{lem:fundamental_domain_general}
    Let $X \in \stratum$ satisfy Keane's condition. Then it admits a unique decomposition into a balanced $k$-set of castle polygons. 
\end{lemma}

\begin{proof}
    We start by proving the existence of a balanced $k$-set. Let $P$ be any decomposition into a $k$-set of castle polygons, which exists by Theorem \ref{thm:existence_castle_polygons} since $X$ satisfies Keane's condition. Without loss of generality we may assume that $P$ satisfies condition $(F1)$ from Definition \ref{def:fundamental_domain_general}. If not, by Lemma \ref{lem:diagonals_become_sides_general} we can apply finitely many backward diagonal change moves such that the resulting $k$-set does satisfy $(F1)$. 

    Now, let $H_{\mathrm{bad}}^{(0)}(P) \subseteq H(P)$ be the subset consisting only of polygons that do not satisfy $(F2)$. If $H_{\mathrm{bad}}^{(0)}(P) = \emptyset$, we are done as $P$ is already balanced. Otherwise, apply all moves corresponding to indices in $H_{\mathrm{bad}}^{(0)}(P)$. Note that since $(F2)$ was not satisfied for all these polygons, the resulting $k$-set $P^{(1)}$ still satisfies $(F1)$. If $P^{(1)}$ is balanced, we are done, otherwise we proceed with the set $H_{\mathrm{bad}}^{(1)}(P)$ defined analogously as above. Invoking again Lemma \ref{lem:diagonals_become_sides_general}, we know that after repeating this process a finite amount of times, we obtain a $k$-set $P^{(m)}$ such that $H_{\mathrm{bad}}^{(m)}(P) = \emptyset$ and thus $P^{(m)}$ is balanced. 

    For uniqueness, assume towards a contradiction that $X$ admits two distinct decomposition into balanced $k$-sets of castle polygons $P$ and $P'$. Since $\mathcal{P}_k$ is a lattice (see Proposition \ref{prop:castle_lattices}), we can define
    \begin{equation*}
        Q \coloneqq P \curlywedge P',
    \end{equation*}
    the greatest lower bound of $P$ and $P'$ with respect to $\prec$. Denote by $\gamma_P$ and $\gamma_{P'}$ two sequences of moves from $Q$ to $P$ and to $P'$ respectively. Since $P \neq P'$, at least one of these sequences must have strictly positive length, say $\gamma_P$. Let $w_{i,d}$ be the first forward diagonal that gets cut in $Q$ with respect to $\gamma_P$. If at some point in $\gamma_{P'}$ the forward diagonal $w_{i,d}$ appears as well, by using distributivity of the lattice $\mathcal{P}_k$ we may assume that it appears in the beginning as well. But then, applying the corresponding move yields a new lower bound, which is greater than $Q$, contradicting the uniqueness of the meet. Thus, $w_{i,d}$ cannot appear in $\gamma_{P'}$. It follows that $(i,\varepsilon) \in H(P')$, so that since $P'$ is balanced we know that $|I'(p_i)| < 1$ in $P'$. This is also true in $Q$, since we do not make a move in this polygon to arrive at $P'$. On the other hand, we \emph{do} make a move in $p_i$ to arrive at $P$, which implies that $|I'(p_i)| \geq 1$ in $Q$ by the remark after Definition \ref{def:width_intervals_general}. Therefore, we arrive at a contradiction and uniqueness is proven. 
\end{proof}

Note that exactly as in section \ref{sec:fundamental_domain_hyperelliptic}, since every aspect of these constructions here rely solely on the geometry, the unique balanced $k$-set of castle polygons is independent of any labeling. Thus, the same construction descends to the unlabeled counterparts of the objects used here.

\begin{definition}[Canonical $k$-set]\label{def:canonical_k_set}
    Given $X \in \stratum$ satisfying Keane's condition, we call the balanced $k$-set of castle polygons constructed in Lemma \ref{lem:fundamental_domain_general} the \emph{canonical $k$-set of castle polygons} associated to $X$, which we denote by $P(X)$.
\end{definition}

We now follow the ideas from section \ref{sec:section_for_teichmueller_flow} to define a subset of the fundamental domain $\mathcal{F}$ from Definition \ref{def:fundamental_domain_general} that will give us a desirable Poincaré section for the Teichmüller geodesic flow $g_t$.

\begin{definition}[Poincaré section]\label{def:poincare_section_general}
    Let $\Upsilon_\mathcal{F}$ be the collection of all $k$-sets of castle polygons $P \in \mathcal{F}$, that satisfy the following additional condition.
    \begin{enumerate}[start=1,label={($S$)}]
        \item There exists $p \in P$ such that $|I(p)| = 1$.
    \end{enumerate}
    \sloppy We further denote by $\Upsilon \subseteq \stratum$ the set of all translation surfaces ${X \in \stratum}$ for which there exists a decomposition into a $k$-set of castle polygons $P \in \mathcal{P}_k$ that belongs to $\Upsilon_\mathcal{F}$.
\end{definition}
Note that the definition of the Poincaré section is exactly the same as in the hyperelliptic case, i.e., the section defined in Definition \ref{def:poincare_section}. Thus, the conceptual description of the Poincaré first return map $\mathcal{S} \colon \Upsilon \to \Upsilon$ remains exactly the same. The derivation of the first return time follows the same structure as well. Let us redo this computation here.

Let $X \in \stratum$ have neither horizontal nor vertical saddle connections and let $P(X)$ be the associated canonical $k$-set of castle polygons as in Definition \ref{def:canonical_k_set}. By definition, $|I'(p)|< 1$ for every $p \in P(X)$ that is not a triangle. Consider now the set
\begin{equation}\label{eq:support_diagonal_changes}
    A \coloneqq \argmax\{|I'(p)| \mid p \in P(X) \text{ is not a triangle}\}.
\end{equation}
The set $A$ is nonempty, but could contain more than one element. However, for all $p_1, p_2 \in A$, we have $|I'(p_1)| = |I'(p_2)|$. To derive the first return time, only the length and not the polygon itself is important, so that we may take any representative $p_\mathrm{max} \in A$. We now define $t(X)$ to be the unique number such that
\begin{equation*}
    \e^{t(X)}|I'(p_\mathrm{max})| = 1,
\end{equation*}
which is equivalent to
\begin{equation}\label{eq:first_return}
    t(X) = -\log|I'(p_\mathrm{max})|.
\end{equation}

Lastly, it remains to show that the number $t(X)$ defined above is indeed the first return time of $X \in \Upsilon$ to the section $\Upsilon$. We adapt the proof of Lemma 5.14 in \cite{delecroixulcigrai20++enumerating}.

\begin{theorem}[First return time and Poincaré map]
    For any $X \in \Upsilon\subseteq \stratum$ that satisfies Keane's condition, the time $t(X)$ defined in \eqref{eq:first_return} is the first return time of $X$ to the section $\Upsilon$ under the Teichmüller geodeisc flow, i.e., $\mathcal{S}(X) = g_{t(X)}(X)$ and $g_t(X) \notin \Upsilon$ for all $0< t < t(X)$. 
\end{theorem}
\begin{proof}
    The proof procedes in two steps. We first show that $t(X)$ is indeed a return time to $\Upsilon$. To finish, we will then argue that it is the \emph{first} return time, i.e., we will show that $g_t(X) \notin \Upsilon$ for all $0 < t < t(X)$.

    Let $P = P(X) = \{p_1, \ldots, p_k\}$ be the canonical $k$-set of castle polygons associated to $X$, and let $A$ be the set defined as in \eqref{eq:support_diagonal_changes}. We further let $P'$ be the $k$-set of castle polygons obtained by performing a diagonal change on all all $p \in A$, followed by an application of $g_{t(X)}$. For all $p' \in P'$ we have $|I(p')| \geq 1$. Indeed, if $p$ is the corresponding polygon in $P$ and $p \notin A$, then we obtain $p'$ by just applying the flow. Then, $I(p') = \e^{t(X)}I(p)$, so that since $|I(p)| \geq 1$ we even have a strict inequality $|I'(p)|  \e^{t(X)} |I(p)| > 1$. On the other hand, if $p \in A$ then by construction we have $|I(p')| = 1$. This shows that $P'$ satisfies $(F1)$ from Definition \ref{def:fundamental_domain_general} as well as $(S)$ from Definition \ref{def:poincare_section_general}.

    Arguing exactly as in Lemma \ref{lem:fundamental_domain}, we may apply a finite number of forward diagonal changes until $P'$ satisfies $(F2)$ as well. Note that we will never have to apply a move to any polygon which satisfies $|I(p)| = 1$, since these polygons automatically satisfy $(F2)$ already as a consequence of Keane's condition. This shows, that $t(X)$ is a return time to $\Upsilon$.

    For uniqueness, let $t\in(0,t(X))$. If $p \in P$ is not a triangle, by definition of $p_\mathrm{max}$ we have
    \begin{equation}\label{eq:poincare_proof_01}
        \e^{t(X)} = \frac{1}{|p_\mathrm{max}|} \leq \frac{1}{|I'(p)|}.
    \end{equation}
    Let us now write $P'$ for the $k$-set obtained by making a forward move at a non-triangle polygon $p$ followed by an application of $g_t$. As before, we write $p'$ for the polygon obtained from $p$ in this way. Since
    \begin{equation*}
        I(p') = \e^t\cdot|I'(p)| < \e^{t(X)}\cdot |I'(p)| \leq 1,
    \end{equation*}
    where the last inequality follows from \eqref{eq:poincare_proof_01}, we see that $P'$ does not satisfy $(F1)$ from Definition \ref{def:fundamental_domain_general}. Therefore, $P'\notin \Upsilon_\mathcal{F}$. 
    
    Now notice that $g_t(P)$ is too wide to be in $\Upsilon_\mathcal{F}$, since applying the flow for any nonzero time will in particular increase the width of any polygon of width 1 by a nonzero amount. Hence, if there was a $k$-set $P''$ with $g_t(P) \sim P''$, where $\sim$ denotes the equivalence relation introduced in \ref{def:equivalence_diagonal_changes}, then we necessarily have $g_t(P) \prec P''$. But any sequence of forward moves begins with a move as described above, which puts $g_t(P)$ outside of $\Upsilon_\mathcal{F}$. Further forward moves only decrease widths further, so we conclude that there cannot exist any $P'' \in \Upsilon_\mathcal{F}$ equivalent to $g_t(P)$ and hence, $g_t(X) \notin \Upsilon$ for all $t \in (0,t(X))$.
\end{proof}

We will now provide a worked example similar to Example \ref{ex:first_return}. The main purpose of this example is to highlight that the situation is completely analogous to the hyperelliptic case.

\begin{example}[Diagonal changes as first return of $g_t$]\label{ex:3-set_example}
    Consider the translation surface $X \in \mathcal{H}(2)$, decomposed into a 3-set of castle polygons given by the following data. The length datum is given by
    \begin{equation*}
        w_{i, \ell} = \left(-\frac{\sqrt{2}}{2}, \frac{1}{2}(2- \sqrt{2})\right), \quad 
        w_{i, r} = \left(\frac{1}{2}(2 - \sqrt{2}), \frac{\sqrt{2}}{2}\right) \quad \text{ for all } i \in [3], 
    \end{equation*}
    and the combinatorial datum is
    \begin{equation*}
        \boldsymbol{\pi} = \combi{r_1\ell_2}\combi{r_3\ell_1}\combi{r_2\ell_3}.
    \end{equation*}

    This 3-set of castle polygons $P$ is depicted in Figure \ref{fig:3-set_example_01}. Note that $X$ is a rotated version of the translation surface we already saw in \ref{ex:torus_triple_cover}, rotated such that the orthogonal eigenvectors of the pseudo-Anosov map derived in the example become the vertical and the horizontal. 
    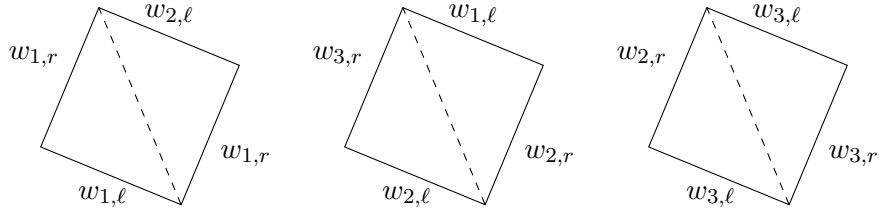
\begin{figure}[ht]
        \centering
        \begin{tikzpicture}[scale = 2]
            \def \wlh{-0.9238};
            \def \wlv{0.3826};
            \def \wrh{0.3826};
            \def \wrv{0.9238};
        
            \def \x{0}
            \draw (\x,0) -- node[midway, below right] {$w_{1, r}$} (\x+\wrh, \wrv) -- node[midway, above] {$w_{2, \ell}$} (\x+\wrh + \wlh, \wrv+\wlv) -- node[midway, above left] {$w_{1, r}$} (\x+\wlh, \wlv) -- node[midway, below, xshift = -3] {$w_{1, \ell}$} cycle; 
            \draw[dashed] (\x,0) -- (\x+\wrh + \wlh, \wrv+\wlv);
            \def \x{2}
            \draw (\x,0) -- node[midway, below right] {$w_{2, r}$}(\x+\wrh, \wrv) -- node[midway, above] {$w_{1, \ell}$}(\x+\wrh + \wlh, \wrv+\wlv) -- node[midway, above left] {$w_{3, r}$}(\x+\wlh, \wlv) -- node[midway, below, xshift = -3] {$w_{2, \ell}$} cycle; 
            \draw[dashed] (\x,0) -- (\x+\wrh + \wlh, \wrv+\wlv);
            \def \x{4}
            \draw (\x,0) -- node[midway, below right] {$w_{3, r}$}(\x+\wrh, \wrv) -- node[midway, above] {$w_{3, \ell}$}(\x+\wrh + \wlh, \wrv+\wlv) -- node[midway, above left] {$w_{2, r}$}(\x+\wlh, \wlv) -- node[midway, below, xshift = -3] {$w_{3, \ell}$} cycle; 
            \draw[dashed] (\x,0) -- (\x+\wrh + \wlh, \wrv+\wlv);
        \end{tikzpicture}
        \caption{The 3-set of castle polygons from example \ref{ex:3-set_example}}
        \label{fig:3-set_example_01}
    \end{figure}
    It is straight-forward to verify that $|I(p_i)| = 1$ for all $i \in [3]$, as well as $|I'(p_i)|= \frac{\sqrt{2}}{2} < 1$, so that $P \in \Upsilon$ and the set $A$ from above coincides with $P$. Moreover, the first return time is given by
    \begin{equation*}
        t(X) = -\log|I'(p_1)| = - \log\frac{\sqrt{2}}{2} = \frac{\log 2}{2}.
    \end{equation*}
    We now apply a diagonal change in all polygons. This gives the 3-set depicted in Figure \ref{fig:3-set_example_02}. In the same figure, we depict in with dashed lines in blue the 3-set obtained by applying $g_{t(X)}$, i.e., applying the Teichmüller geodesic flow up to time $t(X)$. The resulting surface $g_{t(X)}(X)$ is represented by the decomposition $P'$ characterized by the length datum
    \begin{equation*}
        w'_{i, \ell} = \left(\sqrt{2}-2,\frac{\sqrt{2}}{2}\right), \quad w'_{i, r} = \left(\sqrt{2}-1,\frac{1}{2}\right) \quad \text{for all } i \in [3],
    \end{equation*}
    as well as the combinatorial datum
    \begin{equation*}
        \boldsymbol{\pi}' = \boldsymbol{\pi}.
    \end{equation*}
    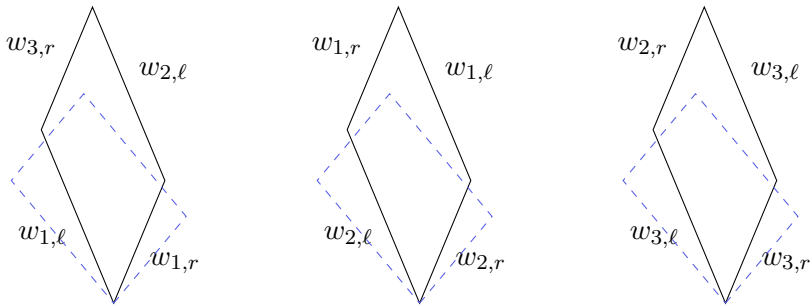
\begin{figure}[ht]
        \centering
        \begin{tikzpicture}[scale = 2.3]
            \def \wlh{-0.4142};
            \def \wlv{1};
            \def \wrh{0.29289};
            \def \wrv{0.7071};            

            \def \x{0}
            \draw (\x,0) -- node[midway, below right] {$w_{1, r}$} (\x+\wrh, \wrv) -- node[midway, above right] {$w_{2, \ell}$} (\x+\wrh + \wlh, \wrv+\wlv) -- node[midway, above left] {$w_{3, r}$} (\x+\wlh, \wlv) -- node[midway, below left] {$w_{1, \ell}$} cycle; 
            \def \x{1.75}
            \draw (\x,0) -- node[midway, below right] {$w_{2, r}$} (\x+\wrh, \wrv) -- node[midway, above right] {$w_{1, \ell}$} (\x+\wrh + \wlh, \wrv+\wlv) -- node[midway, above left] {$w_{1, r}$} (\x+\wlh, \wlv) -- node[midway, below left] {$w_{2, \ell}$} cycle; 
            \def \x{3.5}
            \draw (\x,0) -- node[midway, below right] {$w_{3, r}$} (\x+\wrh, \wrv) -- node[midway, above right] {$w_{3, \ell}$} (\x+\wrh + \wlh, \wrv+\wlv) -- node[midway, above left] {$w_{2, r}$} (\x+\wlh, \wlv) -- node[midway, below left] {$w_{3, \ell}$} cycle;

            \def \wlh{-0.58578};
            \def \wlv{0.7071};
            \def \wrh{0.4142};
            \def \wrv{0.5}; 

            \def \x{0}
            \draw[noamblue, dashed] (\x,0) --  (\x+\wrh, \wrv) --  (\x+\wrh + \wlh, \wrv+\wlv) -- (\x+\wlh, \wlv) -- cycle; 
            \def \x{1.75}
            \draw[noamblue, dashed] (\x,0) --  (\x+\wrh, \wrv) --  (\x+\wrh + \wlh, \wrv+\wlv) -- (\x+\wlh, \wlv) -- cycle; 
            \def \x{3.5}
            \draw[noamblue, dashed] (\x,0) --  (\x+\wrh, \wrv) --  (\x+\wrh + \wlh, \wrv+\wlv) -- (\x+\wlh, \wlv) -- cycle; 
        \end{tikzpicture}
        \caption{Applying 3 diagonal changes and $g_{t(X)}$ to the 3-set in Figure \ref{fig:3-set_example_01}.}
        \label{fig:3-set_example_02}
    \end{figure}
    Similar to what we showed in Example \ref{ex:first_return}, one can show that applying two right and two left moves in each polygon corresponds to a closed Teichmüller geodesic, or equivalently a pseudo-Anosov diffeomorphism, which is exactly the map derived in Example \ref{ex:torus_triple_cover}.
\end{example}

To end this section, let us mention that the section $\Upsilon$ constructed above shares the same central property of the analogous section in the hyperelliptic case: The section $\Upsilon$ intersects all closed Teichmüller geodesics in $\stratum$.

\begin{proposition}
    The Poincaré section $\Upsilon$ from Definition \ref{def:poincare_section_general} intersects \emph{all} closed Teichmüller geodesics in $\stratum$.
\end{proposition}
The proof is obtained by taking the proof of Proposition \ref{prop:crosses_all_geodesics} and exchanging quadrangulations by $k$-sets of castle polygons. 
 \pagebreak

\appendix

\section{A Proof of Keane's Theorem}
\thispagestyle{plain}

We give a proof of Keane's Theorem (Theorem \ref{thm:keanes}), which for the readers convenience we will state again here.
We follow the exposition from \cite{gouezel2006th}.
\begin{theorem}[Keane's Theorem]\label{thm:keanes_appendix}
    If there exists no saddle connection in direction $\theta$, then the linear flow $\varphi_\R^\theta$ is minimal.
\end{theorem}

For the proof, the following two lemmas are essential.

\begin{lemma}\label{lem:keane_lemma_1}
    Let $X$ be a translation surface of genus $\mathbf{g} \geq 2$. If $X$ has a closed vertical geodesic, it also has a vertical saddle connection.
\end{lemma}
\begin{proof}
    Since we assume that $X$ admits a closed vertical geodesic, we can find a point $p \in X$ such that $\varphi_t(p) = p$ for some $t>0$. We can immerse a cylinder around the closed geodesic, small enough so it does not contain any singularities. This is possible since $\Sigma$ is finite. The immersion is of the form
    \begin{equation*}
        f_{r_1,r_2} \colon [-r_1,r_2] \times S^1 \xhookrightarrow{} X,
    \end{equation*}
    for $r_1,r_2 > 0$. An example of such an immersion can be seen in Figure \ref{fig:keane_immersion}.

    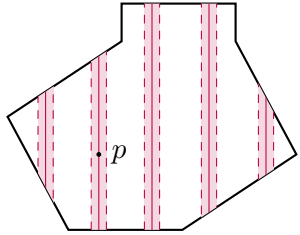
\begin{figure}[ht]
        \centering
        \begin{tikzpicture}
        \coordinate (A) at (0,0);
        \coordinate (B) at (1.5,0);
        \coordinate (C) at (3,1);
        \coordinate (H) at (2.2,2.5);
        \coordinate (D) at (2.2,3);
        \coordinate (E) at (0.7,3);
        \coordinate (F) at (0.7,2.5);
        \coordinate (G) at (-0.8,1.5);

        \draw[thick] (A) -- (B) -- (C) -- (H) -- (D) -- (E) -- (F) -- (G) -- cycle;

        \clip (A) -- (B) -- (C) -- (H) -- (D) -- (E) -- (F) -- (G) -- cycle;

        \fill[purple!100!red!15] (0.3,0) rectangle (0.5,3);
        \draw[line cap = round, purple] (0.4, 0) -- (0.4,3);
        \draw[line cap = round, purple, dashed] (0.3, 0) -- (0.3,3);
        \draw[line cap = round, purple, dashed] (0.5, 0) -- (0.5,3);

        \fill[purple!100!red!15] (2.5,0) rectangle (2.7,3);
        \draw[line cap = round, purple] (2.6, 0) -- (2.6,3);
        \draw[line cap = round, purple, dashed] (2.5, 0) -- (2.5,3);
        \draw[line cap = round, purple, dashed] (2.7, 0) -- (2.7,3);
        
        \fill[purple!100!red!15] (-0.4,0) rectangle (-0.2,3);
        \draw[line cap = round, purple] (-0.3,0) -- (-0.3,3);
        \draw[line cap = round, purple, dashed] (-0.4,0) -- (-0.4,3);
        \draw[line cap = round, purple, dashed] (-0.2,0) -- (-0.2,3);

        \fill[purple!100!red!15] (1.75,0) rectangle (1.95,3);
        \draw[line cap = round, purple] (1.85, 0) -- (1.85, 3);
        \draw[line cap = round, purple, dashed] (1.75, 0) -- (1.75, 3);
        \draw[line cap = round, purple, dashed] (1.95, 0) -- (1.95, 3);

        \fill[purple!100!red!15] (1,0) rectangle (1.2,3);
        \draw[line cap = round, purple] (1.1,0) -- (1.1,3);
        \draw[line cap = round, purple, dashed] (1,0) -- (1,3);
        \draw[line cap = round, purple, dashed] (1.2,0) -- (1.2,3);

        \filldraw[black, draw = black] (0.4,1) circle (0.7pt);
        \node[right, xshift = 1] at (0.4,1) {$p$};
    \end{tikzpicture}
        \caption{An immersion of the form $f_{r_1, r_2}$ as in the proof of Lemma \ref{lem:keane_lemma_1}.}
        \label{fig:keane_immersion}
    \end{figure}

    We can increase the interval $[r_1, r_2]$, but not indefinitely, since the area of the surface is finite. There are exactly two situations where we cannot increase the interval anymore. Either, the two sides of the cylinder coincide at some point, in which case the cylinder would be the whole surface $X$ and consequently, $X$ would be a torus. This is impossible, since we assume $X$ to have a genus $\mathbf{g}$ of at least 2. Therefore, it must be the case that we reach (at least) one singularity with one side of the cylinder. This gives exactly the desired vertical saddle connection.
\end{proof}

\begin{lemma}\label{lem:keane_lemma_2}
    Let $\alpha$ be a non-periodic vertical regular geodesic on a translation surface $X$. Let $p \in \alpha$ and consider a transverse segment $[p,q]$ in $X$. Then,
    \begin{equation*}
        \varphi_{\R_+}(p) \cap (p,q) \neq \emptyset,
    \end{equation*}
    i.e., the vertical trajectory starting at $p$ returns to the interior of the transverse segment.
\end{lemma}

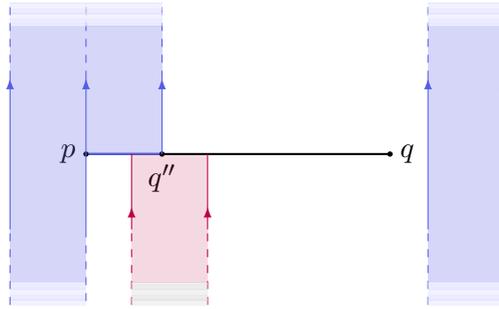
\begin{figure}[ht]
    \centering
    \begin{tikzpicture}
        \coordinate (p) at (0,0);
        \coordinate (q) at (4,0);
        \coordinate (q'') at (1,0);

        \draw[thick] (p) -- (q);
        \draw[thick, green!10!blue!95!red!70] (p) -- (q'');
        \fill[black] (p) circle (1pt); \fill[black] (q) circle (1pt); \fill[black] (q'') circle (1pt);
        \node[left] at (p) {$p$}; \node[right] at (q) {$q$};  \node[below] at (q'') {$q''$};

        \draw[green!10!blue!95!red!70, postaction={
    decorate, decoration={ markings, mark=at position 1 with{\arrow{latex}}}}] (p) -- ($(p) + (0,1)$);
        \draw[dashed, green!10!blue!95!red!70] ($(p) + (0,1)$) -- ($(p) + (0,2)$);

        \draw[green!10!blue!95!red!70, postaction={
    decorate, decoration={ markings, mark=at position 1 with{\arrow{latex}}}}] (q'') -- ($(q'') + (0,1)$);
        \draw[dashed, green!10!blue!95!red!70] ($(q'') + (0,1)$) -- ($(q'') + (0,2)$);

        \fill[pattern = horizontal lines light blue, pattern color = green!10!blue!95!red!70, opacity = 0.5] ($(p) + (0,1.7)$) rectangle ($(q'') + (0,2)$);
        \fill[green!10!blue!95!red!70, opacity = 0.25] (p) rectangle ($(q'') + (0,1.7)$);

        \draw[green!10!blue!95!red!70, postaction={
    decorate, decoration={ markings, mark=at position 1 with{\arrow{latex}}}}] ($(q) + (0.5,-1)$) -- ($(q) + (0.5,1)$);
        \draw[green!10!blue!95!red!70, postaction={
    decorate, decoration={ markings, mark=at position 1 with{\arrow{latex}}}}] ($(q) + (1.5,-1)$) -- ($(q) + (1.5,1)$);

        \draw[green!10!blue!95!red!70, dashed] ($(q) + (0.5,1)$) -- ($(q) + (0.5,2)$);
        \draw[green!10!blue!95!red!70, dashed] ($(q) + (1.5,1)$) -- ($(q) + (1.5,2)$);
        \draw[green!10!blue!95!red!70, dashed] ($(q) + (0.5,-2)$) -- ($(q) + (0.5,-1)$);
        \draw[green!10!blue!95!red!70, dashed] ($(q) + (1.5,-2)$) -- ($(q) + (1.5,-1)$);

        \fill[green!10!blue!95!red!70, opacity = 0.25] ($(q) + (0.5,-1.7)$) rectangle ($(q) + (1.5,1.7)$);
        \fill[pattern = horizontal lines light blue, opacity = 0.5] ($(q) + (0.5,1.7)$) rectangle ($(q) + (1.5,2)$);
        \fill[pattern = horizontal lines light blue, opacity = 0.5] ($(q) + (0.5,-2)$) rectangle ($(q) + (1.5,-1.7)$);

        \fill[green!10!blue!95!red!70, opacity = 0.25] ($(p) + (-1,-1.7)$) rectangle ($(p) + (0,1.7)$);
        \fill[pattern = horizontal lines light blue, opacity = 0.5] ($(p) + (-1,1.7)$) rectangle ($(p) + (0,2)$);
        \fill[pattern = horizontal lines light blue, opacity = 0.5] ($(p) + (-1,-2)$) rectangle ($(p) + (0,-1.7)$);

        \draw[green!10!blue!95!red!70, postaction={
    decorate, decoration={ markings, mark=at position 1 with{\arrow{latex}}}}] ($(p) + (-1,-1)$) -- ($(p) + (-1,1)$);
        \draw[green!10!blue!95!red!70, dashed] ($(p) + (-1,1)$) -- ($(p) + (-1,2)$);  
        \draw[green!10!blue!95!red!70, dashed] ($(p) + (-1,-2)$) -- ($(p) + (-1,-1)$);

        \draw[green!10!blue!95!red!70] (p) -- ($(p) - (0,1)$);
        \draw[green!10!blue!95!red!70, dashed] ($(p) - (0,1)$) -- ($(p) - (0,2)$);

        \draw[purple, postaction = {decorate, decoration = {markings, mark=at position 0.3 with {\arrow{latex}}}}]  ($(p) + (0.6,-1)$) -- ($(p) + (0.6,0)$);
        \draw[purple, postaction = {decorate, decoration = {markings, mark=at position 0.3 with {\arrow{latex}}}}]  ($(q'') + (0.6,-1)$) -- ($(q'') + (0.6,0)$);
        \draw[purple, dashed] ($(p) + (0.6,-2)$) -- ($(p) + (0.6,-1)$);
        \draw[purple, dashed] ($(q'') + (0.6,-2)$) -- ($(q'') + (0.6,-1)$);

        \fill[pattern = horizontal lines light gray, opacity = 0.5] ($(p) + (0.6, -2)$) rectangle ($(q'') + (0.6, -1.7)$);

        \fill[purple, opacity = 0.15] ($(p) + (0.6, -1.7)$) rectangle ($(q'') + (0.6, 0)$);
        
    \end{tikzpicture}
    \caption{The first return map from Lemma \ref{lem:keane_lemma_2}.}
    \label{fig:keane_lemma_2}
\end{figure}

\begin{proof}
    We'll consider the first return map to the segment $[p,q]$. There are only finitely many singularities, hence also only finitely many downward flows starting at these singularities. It follows that the image of $[p,q]$ under the (upward) vertical flow $\varphi_\R$ hits singularities at most finitely many times before it intersects again the segment $[p,q]$. Put differently, the first return map to the segment has finitely many points of discontinuity. We can therefore choose $q' \in (p,q)$ such that the image of $[p,q']$ under $\varphi_\R$ hits no singularities before returning to $[p,q]$. In fact, $[p,q']$ must eventually return to itself. Indeed, if we denote by $R_h$ the rectangle with base $[p,q']$ of height $h$, letting $h \to \infty$ the area of $R_h$ becomes arbitrarily large, so that it is eventually not immersed anymore, since the area of $X$ is finite. 

    Note that the left edge of $R_h$ is included in $\alpha^+ \coloneqq \varphi_{\R_+}(p)$. Let us further denote by $\beta^+$ the right edge of $R_h$. We now have to distinguish two cases. Either $\alpha^+$ intersects $(p,q')$, in which case the lemma is proven, or it is $\beta^+$ that intersects $(p,q')$. In the latter case, let $q''$ be the unique point in $[p,q']$ such that the vertical flow emitted from $q''$ inside $R_h$ hits $p$. We can now repeat the above argument for the segment $[p,q'']$, but we see that now the second case is impossible. Indeed, by construction the rectangle $R'_h$ with base $[p,q'']$ will at some point miss the interval $(p,q)$ to the left as in Figure \ref{fig:keane_lemma_2}, so the first return is of the type that $\alpha^+$ returns to $(p,q)$.
\end{proof}

Having established the lemmas, we can now prove Keane's Theorem.

\begin{proof}[Proof of Theorem \ref{thm:keanes_appendix}]
    Without loss of generality we may assume that the flow is vertical by rotating the translation surface if necessary. Let $\gamma$ be a regular vertical geodesic, which by Lemma \ref{lem:keane_lemma_1} is not closed. Let Y = $\overline{\gamma}$ be the orbit closure. We want to show that $Y = X$.

    If not, take $p \in \partial Y$ and let $\alpha$ be the orbit of the vertical flow starting at $p$, which by invoking Lemma \ref{lem:keane_lemma_1} again we see is not closed. Let $I$ be a transverse segment with endpoint $p$, small enough so its interior is included in $X\setminus Y$. Since $Y$ is invariant under the flow, it follows that $\alpha \subseteq Y$. Lemma \ref{lem:keane_lemma_2} now implies that $\operatorname{Int}(I) \cap \alpha \neq \emptyset$, which is absurd. Therefore, it must be the case that $Y = X$ which shows that all regular geodesics are dense. 
\end{proof}

\section{A Proof of Masur's Criterion}
Here we give a proof of Masur's Criterion (Theorem \ref{thm:masurs_criterion}) which we introduced in section \ref{sec:general_principle}. The proof we present can be found in \cite{monteil2010introduction} and was also communicated in \cite{ulcigrai2022surfacedynamics}.

Note that strata $\mathcal{H}(k_1 - 1, \ldots, k_n - 1)$ are not compact. Indeed, it is easy to construct a sequence converging to a translation outside of the fixed stratum by letting two singularities collide. Preventing this by ensuring that singularities stay away from each other by a fixed amount, we obtain a sequentially compact space. 
\begin{lemma}\label{lem:sequential_compactness}
    Let $\mathcal{H}(k_1 - 1, \ldots, k_n - 1)$ be a stratum and let $\varepsilon>0$. The set
    \begin{equation*}
        K_\varepsilon \coloneqq \{X \in \stratum \mid \operatorname{sys}(X) \geq \varepsilon\} \subseteq \stratum
    \end{equation*}
    is sequentially compact. More explicitly, for any sequence $(X_n)_{n \in \N}$ there exists a subsequence $(X_{n_j})_{j \in \N}$ and some translation surface $X_\infty \in K_\varepsilon$, such that
    \begin{equation*}
        X_{n_j} \xrightarrow{j \to \infty} X_\infty,
    \end{equation*}
    where the convergence is in the sense of Definition \ref{def:convergence_stratum}.
\end{lemma}

In the proof of Lemma \ref{lem:sequential_compactness} we will make use of the following lemma, which can be proven using only Euclidean geometry.

\begin{lemma}\label{lem:triangle_flipping}
    Given $\varepsilon > 0$, there exists a constant $c>0$ such that for all $X \in K_\varepsilon$ there exists a triangulation $\tau = \{T_1, \ldots, T_m\}$ satisfying $\ell(s)\leq c$ for all sides $s$ of the triangles $T_i \in \tau$.
\end{lemma}

\begin{proof}[Proof of Lemma \ref{lem:sequential_compactness}]
    Suppose $(X_n)_{n \in \N}$ is a sequence of translation surfaces in $K_\varepsilon$. For all $n$, let $\tau^n = (T_1^n, \ldots, T_m^n)$ be a flat triangulation of $X_n$ consisting of $m$ triangles. The sides of the triangles $T_i^n$ are saddle connections, hence by assumption and Lemma \ref{lem:triangle_flipping} we have that
    \begin{equation*}
        \varepsilon \leq \ell(s) \leq c,
    \end{equation*}
    for all sides $s$ of the triangles. In particular, the lengths of the sides are contained in a compact subset of $\R$, so that by the Heine-Borel theorem we may extract a subsequence $(\tau^{n_j})_{j \in \N}$, such that
    \begin{equation*}
        T_i^{n_j} \xrightarrow{j \to \infty} T_i^\infty,
    \end{equation*}
    for all $i$, where $T_i^\infty$ is some limit triangle that is non-degenerate. Furthermore, since there are only finitely many gluing patterns for the sides of the triangles, up to taking another subsequence we may assume that the gluing patterns are eventually constant. 

    Let $X_\infty= \{T_1^\infty, \ldots, T_m^\infty\}/_\sim$ be the translation surface obtained by gluing the limit triangles with this gluing pattern. It is clear that $S_\infty\in K_\varepsilon$ and by constructin we have
    \begin{equation*}
        X_{n_j} \xrightarrow{j \to \infty} X_\infty
    \end{equation*}
    in the sense of Definition \ref{def:convergence_stratum}.
\end{proof}

We can now give the proof of Masur's Criterion, which we will state again for convenience.

\begin{theorem}[Masur's Criterion]\label{thm:masurs_criterion_appendix}
    If a translation surface $X \in \mathcal{H}(k_1 - 1, \ldots, k_n - 1)$ is recurrent in direction $\theta$, then the linear flow in direction $\theta$ on $X$ is \emph{uniquely ergodic}. Equivalently, every regular trajectory on $X$ is \emph{equidistributed}.
\end{theorem}
\begin{proof}
    Let $X\in \stratum$ be recurrent in direction $\theta$. Without loss of generality we may assume that $\theta = \frac{\pi}{2}$, i.e., that we are considering the \emph{vertical} linear flow. Otherwise, we may just rotate the translation surface. 

    Towards a contradiction, suppose that the vertical flow $\varphi_\R$ is not uniquely ergodic, so there exist two distinct invariant measures $\nu_1$ and $\nu_2$ as well as an immersed rectangle $R \subseteq X$ with
    \begin{equation*}
        \nu_1(R) \neq \nu_2(R).
    \end{equation*}
    We may assume that both measures are ergodic, since the set of all invariant measures is convex with the ergodic measures as extremal points. Since $X$ is recurrent, there exists a sequence $(t_n)_{n \in \N}$ with $t_n \to \infty$ and $\varepsilon > 0$ such that $\operatorname{sys}(g_{t_n}\cdot X) > \varepsilon$ for all $n \in \N$. By Lemma \ref{lem:sequential_compactness}, passing to a subsequence if necessary, we may assume that $g_{t_n}\cdot X \xrightarrow{n \to \infty} X_\infty$ for some translation surface $S_\infty\in \stratum$.

    Denoting by $\mathbbm{1}_R$ the indicator function of the immersed rectangle $R$, by Birkhoff's ergodic theorem we have
    \begin{equation}\label{eq:birkhoff}
        \lim_{T \to \infty} \frac{1}{T}\int_0^T \mathbbm{1}_R\big(\varphi_t(x)\big) \, \dd t = v_i(R),
    \end{equation}
    for $i \in \{1,2\}$ and $\nu_i$-almost every $x$. For each measure, let $x_i$ be a point such that \eqref{eq:birkhoff} holds. Up to taking a further subsequence, we can assume that $g_{t_n}(x_i) \xrightarrow{n \to \infty} x_i^\infty$ for points $x_i^\infty \in X_\infty$.

    We will first consider the special case where there is an open set disjoint from the singularity set $\Sigma$ that contains a rectangle $R_\infty$ such that $x_1^\infty$ and $x_2^\infty$ are opposite corners of $R_\infty$. For $n$ large enough, we can immerse a rectangle $R_n$ in $g_{t_n}\cdot X$ whose width $w_n$ and height $h_n$ are arbitrarily close to the corresponding dimensions of $R_\infty$ and whose corresponding endpoints are given by $g_{t_n}(x_1)$ and $g_{t_n}(x_2)$.

    We now apply $g_{t_n}^{-1} = g_{-t_n}$ to $R_n$, obtaining a rectangle in $X$ which has height $\e^{t_n}h_n$ and width $\e^{-t_n}w_n$, meaning it is tall and skinny, see the left-hand side of Figure \ref{fig:masurs_criterion_1}. By construction, two opposite endpoints of $g_{t_n}^{-1}R_n$ are exactly $x_1$ and $x_2$.

    \begin{figure}[ht]
        \centering
        \begin{tikzpicture}
            \draw[dashed] (0,0) ellipse (1.5 and 2);
            \draw[thick] (1,0.6) rectangle (-1,-0.6);
            \node[right] at (1,0) {$R$};

            \node[above] at (0,2) {$X$};

            \clip (0,0) ellipse (1.5 and 2);
            \filldraw[fill = green!10!blue!95!red!70, fill opacity = 0.25, draw = black] (0,0) rectangle (-0.1,-3);
            \fill[black] (0,0) circle (1pt);\node[right] at (0,0) {$x_2$};

            \filldraw[fill = green!10!blue!95!red!70, fill opacity = 0.25, draw = black] (-0.7,-1) rectangle (-0.6,2);
            \fill[black] (-0.7,-1) circle (1pt); \node[left] at (-0.7,-1) {$x_1$};

            \filldraw[fill = green!10!blue!95!red!70, fill opacity = 0.25, draw = black] (-1.35,-2) rectangle (-1.25,2);

            \filldraw[fill = green!10!blue!95!red!70, fill opacity = 0.25, draw = black] (-0.45,-2) rectangle (-0.35,2);

            \filldraw[fill = green!10!blue!95!red!70, fill opacity = 0.25, draw = black] (0.8,-2) rectangle (0.9,2);
        \end{tikzpicture}\qquad \qquad
        \begin{tikzpicture}
            \draw[dashed] (0,0) ellipse (1.5 and 2);
            \filldraw[fill = green!10!blue!95!red!70, fill opacity = 0.25, draw = black] (0.8,1) rectangle (-0.5,-0.8);
            
            \fill[black] (-0.5,-0.8) circle (1pt); \node[below] at (-0.5,-0.8) {$g_{t_n}(x_1)$};
            \fill[black] (0.8,1) circle (1pt); \node[above] at (0.45,1) {$g_{t_n}(x_2)$};

            \node[above] at (0,2) {$g_{t_n}\cdot X$};
        \end{tikzpicture}
        \begin{tikzpicture}[remember picture, overlay]
            \draw[>=latex, ->, bend left=30, thick] (-5,4) to node[above] {$g_{t_n}$} (-2.7,4);
        \end{tikzpicture}
        \caption{An illustration of the renormalization procedure done in the proof of Masur's Criterion.}
        \label{fig:masurs_criterion_1}
    \end{figure}
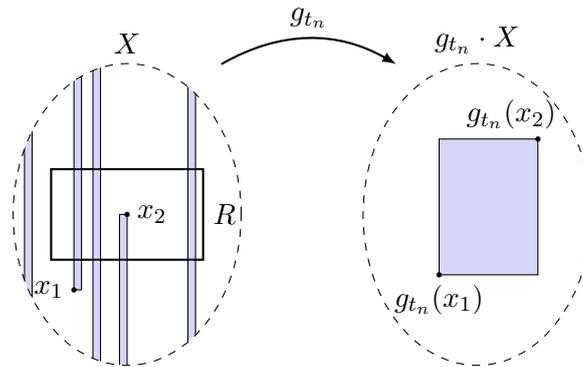

    We set $T_n \coloneqq \e^{t_n}h_n$ and further we write
    \begin{align*}
        I_1^n &\coloneqq \int_0^{T_n} \mathbbm{1}_R\big(\varphi_t(x_1)\big) \, \dd t, \\
        I_2^n &\coloneqq \int_0^{T_n} \mathbbm{1}_R\big(\varphi_t(x_2)\big) \, \dd t,
    \end{align*}
    where $R$ denotes the immersed rectangle from the beginning of the proof. Notice that $I_1^n$ gives the total length of the intersections of the vertical left side of $g_{t_n}^{-1}R_n$, which we will denote by $\ell_n$, and $R$. Formally,
    \begin{equation*}
        I_1^n = \operatorname{Leb}(\ell_n \cap R),
    \end{equation*}
    where $\operatorname{Leb}$ denotes the one-dimensional Lebesgue measure. Similarly, we can write $I_2^n$ as the total lengths of the intersections of the right side of $g_{t_n}^{-1}R_n$, which we will denot by $r_n$, and the rectangle $R$. As mentioned above, $g_{t_n}^{-1}R_n$ is tall and skinny, so it will enter $R$ many times as illustrated in Figure \ref{fig:masurs_criterion_2}.
    \begin{figure}[ht]
        \centering
        \begin{tikzpicture}
            \draw[thick] (-2,-1) rectangle (2,1);
            \node[above] at (1,1) {$R$};

            \foreach \coord in {-1.9,-1.1, -0.1, 0.4, 1.85}{
                \draw[color = purple] (\coord,1.5) -- (\coord, -1.5);
                \draw[color = purple, dashed] (\coord,1.5) -- (\coord, 2);
                \draw[color = purple, dashed] (\coord,-2) -- (\coord, -1.5);
                \draw[color = green!10!blue!95!red!70] (\coord - 0.2, 1.5) -- (\coord - 0.2, -1.5);
                \draw[color = green!10!blue!95!red!70, dashed] (\coord - 0.2, 1.5) -- (\coord - 0.2, 2);
                \draw[color = green!10!blue!95!red!70, dashed] (\coord - 0.2, -2) -- (\coord - 0.2, -1.5);
            }
            \node[below, color = green!10!blue!95!red!70] at (1.65, -2) {$\varphi_t(x_1)$};
            \draw[>=latex, ->, bend left=30, color = green!10!blue!95!red!70] (1.3,-2.1) to (1.6,-1.2);

            \node[below, color = purple] at (-1.1, -2) {$\varphi_t(x_2)$};
            \draw[>=latex, ->, bend right=30, color = purple] (-0.8,-2.1) to (-1.05,-1.2);
            
        \end{tikzpicture}
        \caption{The vertical flows emitted from $x_1$ and $x_2$ crossing the rectangle $R$ in the proof of Masur's Criterion.}
        \label{fig:masurs_criterion_2}
    \end{figure}
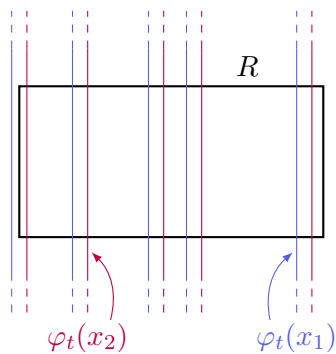
    
    Except for at most two exceptions, either $\ell_n$ and $r_n$ enter the rectangle $R$ together, or neither of them enters $R$. The exceptions may happen on the left and the right side of $R$. This implies that
    \begin{equation*}
        |I_1^n - I_2^n| \leq h(R),
    \end{equation*}
    where $h(R)$ denotes the height of $R$. Since $T_n \xrightarrow{n \to \infty} \infty$, we have
    \begin{equation*}
        \frac{1}{T_n} |I_1^n - I_2^n| \leq \frac{h(R)}{T_n} \conv 0,
    \end{equation*}
    which by \eqref{eq:birkhoff} implies $\nu_1(R) = \nu_2(R)$, a contradiction. We can obtain the general case where the rectangle $R_\infty$ is not of the form assumed above from this first case as follows. Since $x$ and $y$ are not on the same vertical as any of the singularities in $\Sigma$, we can ensure that $g_{t_n}(x_1)$ and $g_{t_n}(x_2)$ are bounded away uniformly from $\Sigma$. Consequently, $x_1^\infty$ and $x_2^\infty$ are not singularities of $X_\infty$. Note that $X_\infty$ is connected and locally path-connected, hence it is path-connected and there exists a path from $x_1^\infty$ to $x_2^\infty$. We can surround this path by an open set that is disjoint from $\Sigma$ as in figure \ref{fig:masurs_criterion_3}.
    \begin{figure}[ht]
        \centering
        \begin{tikzpicture}
            \draw[] (0,0) rectangle (4,0.5);
            \draw[] (4,0.5) rectangle (5,1.2);
            \draw[] (5,1.2) rectangle (5.7, 1.5);
            \draw[] (5.7, 1.5) rectangle (7, 0.7);
            \draw[] (7, 0.7) rectangle (7.7, -1);

            \draw[dashed] (-0.5,-0.5) 
        .. controls (-1,0) and (-1,0.5) ..
        (-0.5,1)
        .. controls (0,1.5) and (2,0.7) ..
        (3,1)
        .. controls (5.3,2) and (6.5,2.5) ..
        (8, 1)
        .. controls (8.5,0.5) and (9,-1) ..
        (7, -1.5)
        .. controls (6.5,-1.6) and (6,0.5) ..
        (5, 0)
        .. controls (4,-1) and (0,-1)..
        (-0.5, -0.5);

        \fill[black] (0,0) circle (1pt); \node[below right] at (-0.3,0) {$x_1^\infty = y_1$};
        \fill[black] (4,0.5) circle (1pt); \node[above left] at (4,0.5) {$y_2$};
        \fill[black] (5,1.2) circle (1pt); \node[below left] at (5,1.2) {$y_3$};

        \fill[black] (5.7,1.5) circle (1pt); \node[above] at (5.7,1.5) {$\ldots$};
        \fill[black] (7,0.7) circle (1pt);
        
        \fill[black] (7.7,-1) circle (1pt); \node[below right] at (7.8,-1) {$y_N = x_2^\infty$};
        
        \end{tikzpicture}
        \caption{The general case in the proof of Masur's Criterion.}
        \label{fig:masurs_criterion_3}
    \end{figure}
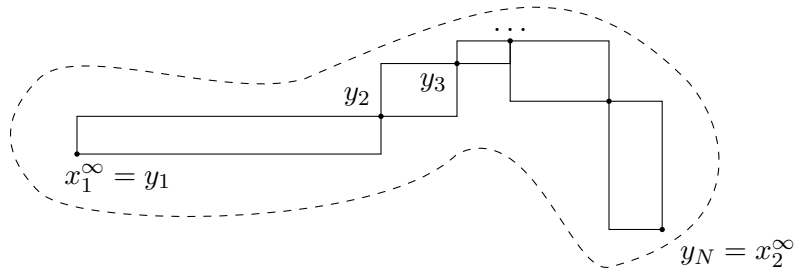

    Thus, there exists a finite sequence $x_1^\infty = y_1, y_2, \ldots, y_{N-1},y_N = x_2^\infty$ of points in $X_\infty$, such that each rectangle with opposite corners $y_i$ and $y_{i+1}$ is contained in this open set. Perturbing the points and possibly passing to subsequences, we may assume that each $y_i$ is the limit point of the trajectory of some $x_i \in X$ under $g_{t_n}$ that is generic for some invariant ergodic measure $\mu_i$. Applying the argument for the basic case finitely many times, we conclude that 
    \begin{equation*}
        \nu_1 = \mu_1 = \mu_2 = \ldots = \mu_{N-1} = \mu_N = \nu_2
    \end{equation*}
    holds in the general case as well.
\end{proof}

Let us mention that Masur's Criterion is one of the two ingredients of the celebrated theorem of Kerkhoff--Masur--Smillie, which states that the behavior of unique ergodicity is generic in a very strong sense.

\begin{theorem}[Kerkhoff--Masur--Smillie]
    For \emph{any} translation surface $X \in \stratum$ and for almost any direction $\theta \in S^1$, the linear flow $\varphi_\R^\theta$ on $X$ is uniquely ergodic. 
\end{theorem}
\begin{remark}
    Note that this result holds for \emph{all} translation surfaces and not just on a set of full measure.
\end{remark}
We can obtain the theorem by combining Masur's Criterion with the following proposition. For details, we refer to \cite{kerckhoff1986ergodicity}.
\begin{proposition}
    For all translation surfaces $X\in \stratum$, the set
    \begin{equation*}
        \{\theta \in S^1 \mid \operatorname{sys}(g_t^\theta \cdot X) \xrightarrow{t \to \infty} 0\}
    \end{equation*}
    has Lebesgue measure zero.
\end{proposition}
 \pagebreak
\thispagestyle{empty}
\
\pagebreak

\thispagestyle{plain}
\listoffigures \pagebreak

\bibliographystyle{alpha}
\bibliography{refs}

\end{document}